\documentclass[12pt,twoside,makeidx]{memo-l}
\usepackage{latexsym}
\usepackage{amsmath}
\usepackage{epsfig}
\usepackage{epic}
\usepackage{amssymb}
\usepackage{enumerate, afterpage, float, tikz, url, mparhack}

\usepackage{hyperref}

\usetikzlibrary{calc,through,patterns}

\colorlet{lines}{green!70!black}
\colorlet{blines}{cyan!70!black}
\colorlet{rlines}{red!70!white}
\colorlet{olines}{brown}
\colorlet{dblines}{blue!65!black}
\colorlet{glines}{green!70!black}
\colorlet{orlines}{orange}
\colorlet{labels}{blue!65!black}
\colorlet{glabels}{green!70!black}
\colorlet{olabels}{brown}
\colorlet{rlabels}{red!70!white}
\colorlet{orlabels}{orange}
\colorlet{rfill}{red}
\colorlet{dblabels}{blue!65!black}

\newtheorem{theorem}{Theorem}[chapter]
\newtheorem{lemma}[theorem]{Lemma}
\newtheorem{corollary}[theorem]{Corollary}
\newtheorem{proposition}[theorem]{Proposition}

\newtheorem{sublemma}{}[theorem]
\newtheorem{conjecture}[theorem]{Conjecture}

\theoremstyle{definition}

\theoremstyle{remark}

\numberwithin{equation}{section}

\newcommand{\ba}{\backslash}

\newcommand{\co}{{\rm co}}
\newcommand{\si}{{\rm si}}
\newcommand{\cl}{{\rm cl}}
\newcommand{\coh}{{\rm coh}}
\newcommand{\fr}{{\rm fr}}
\newcommand{\fcl}{{\rm fcl}}

\newcommand{\clstar}{{\rm \cl^{(*)}}}

\newcommand{\calP}{{\mathcal P}}
\newcommand{\calQ}{{\mathcal Q}}

\newcommand{\subproof}{\begin{proof}[Subproof]}

\newcommand{\FF}{{\bf F}}

\newcommand{\PP}{{\bf P}}
\newcommand{\QQ}{{\bf Q}}
\newcommand{\BB}{{\bf B}}
\newcommand{\CC}{{\bf C}}
\newcommand{\RR}{{\bf R}}

\newcommand{\OO}{{\bf O}}
\newcommand{\AAA}{{\bf A}}

\newcommand{\hQ}{\hat{Q}}
\newcommand{\hP}{\hat{P}}
\newcommand{\hO}{\hat{O}}
\newcommand{\cP}{\check{P}}
\newcommand{\cQ}{\check{Q}}
\newcommand{\less}{\preceq}

\newcommand{\udots}{\cup\cdots\cup }
\newcommand{\smalle}{\preceq}
\newcommand{\more}{\succeq}
\newcommand{\eq}{{\mathcal E}(q)}
\newcommand{\youq}{{\mathcal U}(q)\cap{\mathcal U}^*(q)}

\newcommand{\kstar}{M^*(K_{3,t})}

\makeindex

\begin{document}

\sloppy

\title{Inequivalent Representations 
of Matroids over Prime Fields}

\author{Jim Geelen}
\address{Department of Combinatorics and Optimization\\
University of Waterloo\\
Waterloo, Ontario, Canada}
\email{jfgeelen@math.uwaterloo.ca}
\thanks{The first author was supported by a grant from the
National Sciences and Engineering Research Council of Canada.}


\author{Geoff Whittle}
\address{School of Mathematics, Statistics and Operations Research\\
Victoria University of Wellington,  New Zealand}
\email{geoff.whittle@vuw.ac.nz}
\thanks{The second author was supported by a grant from the Marsden
Fund of New Zealand.}

\subjclass{05B35}
\date{}

\begin{abstract}
It is proved that for each prime field $GF(p)$, there is an integer
$f(p)$ such that a $4$-connected  matroid has
at most $f(p)$ inequivalent representations over $GF(p)$. 
We also prove a stronger theorem that obtains the same conclusion
for matroids satisfying a connectivity condition,
intermediate between $3$-connectivity and $4$-connectivity
that we term ``$k$-coherence''. 

We obtain a variety of other 
results on inequivalent representations
including the following curious one.
For a 
prime power $q$, let ${\mathcal R}(q)$ denote the set of matroids
representable over all fields with at least $q$ elements.
Then there are infinitely many Mersenne primes if and only if,
for each prime power $q$, there is an integer $m_q$ such that a
$3$-connected member of ${\mathcal R}(q)$ has at most 
$m_q$ inequivalent $GF(7)$-representations.

The theorems on inequivalent representations of matroids are
consequences of structural results that do not rely on representability.
The bulk of this paper is devoted to proving such results.
\end{abstract}

\maketitle

\tableofcontents

\mainmatter

\chapter{Introduction}

In this paper
we prove the following theorem.

\begin{theorem}
\label{4-conn}
Let $p$ be a prime number. Then there is an integer $\gamma(p)$
such that a $4$-connected matroid has at most $\gamma(p)$
inequivalent $GF(p)$-representations.
\end{theorem}

We also prove a somewhat stronger theorem that obtains the same
conclusion for a weaker notion of connectivity. Before
discussing this, and other results in this paper, we provide 
some background.

It is easily seen that if a matroid 
$M$ is binary, then $M$ is uniquely
representable over any field for which it is representable. It
is also straightforward to show that ternary matroids are uniquely
representable over $GF(3)$. 
In \cite{ka88}, Kahn proved that 3-connected quaternary matroids are
uniquely 
representable over $GF(4)$. In that paper he made the conjecture
that for any finite field $GF(q)$, there is an integer $\mu(q)$
such that a $3$-connected matroid has at most $\mu(q)$ inequivalent
$GF(q)$ representations. In \cite{oxvewh96} Oxley, Vertigan
and Whittle proved that
Kahn's Conjecture holds for $GF(5)$, but, unfortunately,
examples in that paper show that the conjecture fails for all
fields larger than $GF(5)$. 

One may hope to recover the situation by increasing the connectivity.
What about $4$-connected matroids? For non-prime fields, the situation
is somewhat dire. It is shown in \cite{gegewh09} 
that there are $4$-connected
matroids with an arbitrary number of inequivalent representations over
any non-prime field with at least 9 elements. Indeed, for 
$n\geq m$, there is a vertically $(m+1)$-connected matroid that has
at least $2^{n-1}$ inequivalent representations over any finite field
of non-prime order $q\geq m^m$.

As Theorem~\ref{4-conn} shows, the situation is much 
better for prime fields. Before turning to a more detailed discussion
of the contents of this paper we mention a significant application
of the results of this paper. Seymour \cite{seymour1} showed that
in the worst case it requires exponentially many rank evaluations
to prove that a matroid is binary and this negative result
extends easily to other fields \cite{gewh1}. In contrast to this,
it is proved in \cite[Theorem~1.1]{gewh1} that, for any prime $p$,
an $n$-element matroid can be proved to 
be not representable over $GF(p)$
using only $O(n^2)$ rank evaluations. Results from this paper form an
essential ingredient in the proof of this result.

We also obtain other consequences that we believe are
interesting in their own right. Here is one. Recall
that a {\em Mersenne prime} 
\index{Mersenne prime} is one that has the form
$2^n-1$ for some integer $n$. A very well-known conjecture
is that the number of Mersenne primes is infinite. For a 
prime power $q$, let ${\mathcal R}(q)$ denote the set of matroids
representable over all fields with at least $q$ elements.

\begin{theorem}
\label{mersenne}
There are infinitely many Mersenne primes if and only if,
for each prime power $q$, there is an integer $m_q$ such that a
$3$-connected member of ${\mathcal R}(q)$ has at most 
$m_q$ inequivalent $GF(7)$-representations.
\end{theorem}

The results that we have highlighted so far 
are all to do with representations of matroids over fields.
It is striking that it is not until the short last chapter 
of this paper that
fields play a role. Until then all of our results
are purely structural and are valid independently of any
assumptions about representability. The potential for a matroid
to have inequivalent representations is due to 
certain structural features that are inherent in the matroid.
These are the features that we focus on. The applications to
fields are obtained as consequences of our structural 
results.

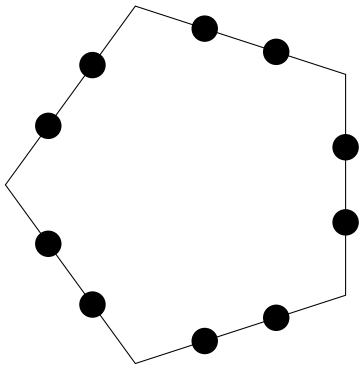
\begin{figure}
\begin{tikzpicture}
	\coordinate (A) at (36:2.5);
	\coordinate (B) at (108:2.5);
	\coordinate (C) at (180:2.5);
	\coordinate (D) at (252:2.5);
	\coordinate (E) at (324:2.5);
	\draw (A) -- (B) -- (C) -- (D) -- (E) -- cycle;
	\coordinate (a1) at ($(A)!0.33!(B)$);
	\coordinate (a2) at ($(A)!0.33!(E)$);
	\coordinate (b1) at ($(B)!0.33!(A)$);
	\coordinate (b2) at ($(B)!0.33!(C)$);
	\coordinate (c1) at ($(C)!0.33!(B)$);
	\coordinate (c2) at ($(C)!0.33!(D)$);
	\coordinate (d1) at ($(D)!0.33!(E)$);
	\coordinate (d2) at ($(D)!0.33!(C)$);
	\coordinate (e1) at ($(E)!0.33!(A)$);
	\coordinate (e2) at ($(E)!0.33!(D)$);
	\foreach \pt in {a1,a2,b1,b2,c1,c2,d1,d2,e1,e2} \fill[black] (\pt) circle (5pt);
\end{tikzpicture}
\caption{A Rank-$5$ Free Swirl}\label{rank-5-free-swirl}
\end{figure}

We now outline our stronger theorem on
inequivalent representations.
Let $M$ be a matroid with a 
basis $B=\{b_1,b_2,\ldots,b_n\}$ such that, for all 
$i\in\{1,2,\ldots,n-1\}$,
there are elements $P_i=\{p_i,q_i\}$ placed freely on the line 
spanned by 
$\{b_i,b_{i+1}\}$ and pairs of 
points $P_n=\{p_n,q_n\}$ placed freely on the 
line spanned by $\{b_n,b_1\}$. Then the matroid $M\ba B$
is the rank-$n$ {\em free swirl}, denoted $\Delta_n$. 
Figure~\ref{rank-5-free-swirl}
illustrates $\Delta_5$. For prime fields, the class of free swirls 
provide the counterexamples to 
Kahn's Conjecture given in \cite{oxvewh96}. It is natural to ask
if all such counterexamples are, in some sense, related to
free swirls. In this paper we provide what is, essentially,
a positive answer to that question.

\begin{figure}
\begin{tikzpicture}[thick,line join=round]
	\coordinate (A) at (11:2.5);
	\coordinate (B) at (83:2.5);
	\coordinate (C) at (155:2.5);
	\coordinate (D) at (227:2.5);
	\coordinate (E) at (299:2.5);
	\coordinate (ab) at ($(A)!0.5!(B)$);
	\coordinate (bc) at ($(B)!0.5!(C)$);
	\coordinate (cd) at ($(C)!0.5!(D)$);
	\coordinate (de) at ($(D)!0.5!(E)$);
	\coordinate (ea) at ($(E)!0.5!(A)$);
	\node[pattern color=lines,draw,circle through=(A),pattern=north east lines] (AB) at (ab) {}; 
	\node[pattern color=lines,draw,circle through=(B),pattern=north west lines] (BC) at (bc) {}; 
	\node[pattern color=lines,draw,circle through=(C),pattern=horizontal lines] (CD) at (cd) {}; 
	\node[pattern color=lines,draw,circle through=(D),pattern=vertical lines] (DE) at (de) {}; 
	\node[pattern color=lines,draw,circle through=(E),pattern=north west lines] (EA) at (ea) {}; 
	\filldraw[fill=white] (A) -- (B) -- (C) -- (D) -- (E) -- cycle;
\end{tikzpicture}
\caption{A Swirl-like Flower}\label{swirl-like-flower}
\end{figure}
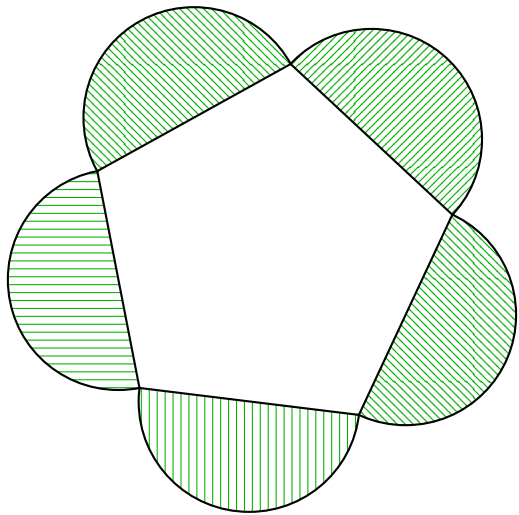

Note that a set $X$ of elements of $\Delta_n$
is non-trivially 3-separating if and only if it is a union of 
members of $(P_1,P_2,\ldots,P_n)$ 
that are consecutive in the cyclic order. 
More generally, a $3$-connected matroid $M$ has a 
{\em swirl-like flower}
with $n$-petals if there is a partition 
$\PP=(P_1,P_2,\ldots,P_n)$ of $E(M)$
into exactly $3$-separating sets called the {\em petals} 
of $\PP$ such that
a union of petals is $3$-separating if and only if it is 
consecutive in
the cyclic order. Figure~\ref{swirl-like-flower} illustrates a 
swirl-like flower with
five petals. It is possible for petals to be degenerate in a way that 
we explain later. If $n\geq 4$, then the {\em order} of a swirl-like 
flower is the number of non-degenerate petals it has. To control
inequivalent representations we control swirl-like flowers.
Let $k\geq 5$ be an integer. Then a matroid is $k$-{\em coherent}
if it 3-connected and has no swirl-like flowers of order $k$. The main
theorem of this paper is really the following.

\begin{theorem}
\label{biggie}
Let $k\geq 5$ be an integer and $p$ be a prime number. Then there
is a function $\gamma(k,p)$ such that a $k$-coherent matroid
has at most $\gamma(k,p)$ inequivalent representations over $GF(p)$.
\end{theorem}

As $4$-connected matroids are $k$-coherent Theorem~\ref{4-conn} is an
immediate corollary of Theorem~\ref{biggie}. 
While Theorem~\ref{biggie} is somewhat more technical than
Theorem~\ref{4-conn}, it is considerably stronger.
We also observe that
it is, in general, easiest to prove a theorem for the weakest version
of connectivity for which the theorem is true. This is because weaker
notions of connectivity are usually easier to keep in 
minors and therefore
facilitate inductive arguments. We know of no way to obtain
a bound on the number of inequivalent representations of 4-connected
matroids other than as a consequence of a stronger theorem using
a weaker connectivity notion.

We now consider the structure of this paper. It was always clear
that this was going to be a long paper, although it never occurred 
to us that it would be {\em this} long. Our original intention was to
partition it into a sequence of papers, but the interconnectivity
of the material was so high that this strategy seemed increasingly
artificial and we eventually abandoned it. In the discussion that
follows we use loosely a number of terms that are defined formally
later in the paper.

Chapters~\ref{background} contains known material,
mainly on connectivity, that is fundamental
to this paper. Flowers are structures in matroids that arise
when 3-separations cross. Flowers come in several different
types and there is a natural equivalence and partial order
on the flowers  in a matroid. Flowers were introduced and 
studied in \cite{flower} and further studied in \cite{flower2}.
It seems that just about every known elementary fact on flowers
is needed at some stage in this paper, and some facts are needed
many times. Chapter~\ref{flower-chapter} is primarily a survey of
basic properties of flowers. 

Let $k\geq 5$ be an integer. As noted above, a 
$k$-coherent matroid is one that
is 3-connected and has no swirl-like flower of order at least $k$. 
This is our basic notion of connectivity. 
Chapter~\ref{k-coherent-chapter} describes properties of this
connectivity notion. If $M$ is a $k$-coherent matroid and
$T$ is a triangle of $M$, then one would typically expect
there to be an element $t\in T$ such that $M\ba t$ is
$k$-coherent. Unfortunately this is not always the case
and triangles that do not have this property are called $k$-{\em wild}.
The structure of $M$ relative to a $k$-wild triangle is described.
If $f$ is an element of the $k$-coherent matroid $M$
and $M\ba f$ and $M/f$ are both $3$-connected, then one
might hope that at least one of these matroids is 
$k$-coherent. Again this is not always the case. If neither
$M\ba f$ nor $M/f$ is $k$-coherent, then we say that $f$ is {\em feral}.
The structure of $M$ relative to a feral element is described.
Feral elements and $k$-wild triangles arise repeatedly in
proofs in later chapters of the paper.  We prove a 
wheels and whirls type theorem for $k$-coherent matroids.

The underlying cause of inequivalent representations in matroids
is that an element may have freedom. A $k$-skeleton is 
a $k$-coherent matroid whose elements have, in some sense, maximum
freedom. It is easily seen, and shown in Chapter~\ref{fields-at-last},
that the number of inequivalent representations of a $k$-coherent
matroid over a finite field is bounded above by the maximum
of the number of inequivalent representations of its $k$-skeleton
minors. In Chapter~\ref{k-skeletons} $k$-skeletons and their
basic properties are described. Again we find that certain
structures arise that are counterexamples to expected behaviour; we
call these structures {\em bogan couples} 
and {\em gangs of three}.
In a way that is entirely analogous to the theorems for
$k$-wild  triangles and feral elements, we 
give theorems that describe the local structure of a matroid
relative to a bogan couple or a gang of three. 

Chapter~\ref{chain-gang} gives a chain theorem for $k$-skeletons.
It is shown that if $M$ is a $k$-skeleton, then unless $M$ is trivially
small, $M$ has a $k$-skeleton minor $N$ such that $|E(M)-E(N)|\leq 4$.
Viewed from a bottom up perspective it gives us a way of building all
$k$-skeletons in a class. A 4-element jump may seem somewhat daunting,
but it is shown that the structure when 3- and 4-element jumps
are required is quite specific. 
Indeed reasonably special structure arises even in the 2-element
case. Apart from a handful of lemmas the results of 
Chapter~\ref{chain-gang} are not needed in the rest of the paper.
We expect that these results will be used to obtain an explicit
bound on the number of inequivalent representations of 
4-connected matroids over $GF(7)$.

Chapter~\ref{paths-of-3-separations} marks a sharp change in the
techniques of this paper. Up to this stage we have focussed on exact
structure. From now on  techniques are extremal. The
rank-$q$ {\em free spike} $\Lambda_q$ is defined in 
Chapter~2.5.
Let $\eq$
\index{$\eq$} 
denote the class of matroids with no $U_{2,q+2}$-,
$U_{q,q+2}$- or $\Lambda_q$-minor. 
The goal is
to eventually show that there are only a finite number of
$k$-skeletons in $\eq$. This gives immediate corollaries
for matroids representable over fields as, if $q$ is prime,
the matroids representable over $GF(q)$ are a subclass of $\eq$.

A {\em path of $3$-separations} 
in a matroid is an ordered partition of the
ground set that induces a nested sequence
of 3-separations in $M$. In Chapter~\ref{paths-of-3-separations}
it is shown that a $k$-skeleton with a sufficiently long  path
of 3-separations cannot be in $\eq$. This begins the process
of taming the structure of a sufficiently large $k$-skeleton
in $\eq$. In Chapter~\ref{taming} this process is 
continued. We show that a sufficiently large $k$-skeleton
must contain, as a minor, a large $4$-connected matroid all
elements of which are neither fixed nor cofixed. This refining 
process is continued and it is shown that a sufficiently
large 4-connected matroid with the above property has a 
large 4-connected minor whose ground set contains many clonal pairs.

The refining process is further continued in 
Chapter~\ref{unavoidable-minors-of-large} where we consider
unavoidable minors of 4-connected matroids whose ground set
contains many clonal pairs. In fact it suffices to focus
on 3-connected matroids whose ground set has a partition into
clonal pairs. It is shown that a sufficiently large
such matroid $M$ in $\eq$ must have a large free-swirl minor. Moreover,
we can also guarantee that all of the clonal pairs in this minor
are clonal pairs in $M$. 

We now have a guaranteed large free-swirl minor. If we delete
a clonal pair from this minor we obtain a matroid with a
natural partition that gives a certain type of 
path of $2$-separations.  This is a minor of 
a matroid in which all of the 2-separations are bridged. 
Moreover the minor has a partition into clonal pairs that
remain clonal pairs in the large matroid.
In Chapter~\ref{strict-2-paths}
it is shown that if the path of  $2$-separations is sufficiently
long, then unless the $2$-separations are bridged in 
a very specific way, the bridging matroid cannot be in $\eq$.

Chapter~\ref{last-rites} returns attention to 
our large free-swirl minor. Such a matroid displays a path of
$3$-separations. These 3-separations are bridged by our larger matroid.
Here we obtain the final win by showing that, if the swirl
is sufficiently large and the clonal pairs of the free
swirl remain clonal in the matroid $M$ that bridges the separations,
then $M$ cannot be in $\eq$. In other words there are a finite
number of $k$-skeletons in $\eq$.

Having achieved the climax of the paper---which 
has no doubt been 
reached in a fever pitch of excitement---the reader can 
relax and enjoy the denouement which consists of the applications
to fields given in the final chapter. Alternatively one could cheat
and begin by reading the last chapter first.
This strategy is recommended.

\subsection*{Acknowledgment} For the last decade we have been
working collaboratively with Bert Gerards, primarily for research 
on another topic. Nonetheless Bert was involved in many discussions
on numerous aspects of the research that led to this paper. 
We thank Bert for many helpful contributions.

\chapter{Background Material}
\label{background}

It is assumed that the reader is familiar with matroid theory
as set forth in Oxley~\cite{ox92}.
Terminology and notation follows \cite{ox92} with the minor
exception that we denote the simple and cosimple matroids canonically
associated with a matroid $M$ by $\si(M)$
\index{si$(M)$} 
and $\co(M)$
\index{co$(M)$}
respectively. This chapter attempts to cover some of the
additional terminology and known results that we will need.
It is largely an accumulation of material scattered in
research papers, although results covered in \cite{ox92}
that are of particular significance here are restated.
Much of it will be needed throughout although some of it is not needed
until later in the paper; for example blocking sequences are not
used until Chapter~\ref{unavoidable-minors-of-large}. 
It just seemed more natural to include blocking sequences
in a chapter primarily
concerned with connectivity techniques. On the other hand
some introductory material is delayed; for example basic facts
on freedom in matroids are not introduced until 
Chapter~\ref{k-skeletons}.

We note that duality plays an integral role throughout this
paper. Almost all of our concepts and results
are either self dual or
have dual formulations that need to be grasped. This applies,
for example, to something as elementary as the closure 
operator of a matroid where we need to understand the
behaviour of the dual operator, the {\em coclosure} 
operator. We often neglect to state dual versions of 
results and in any unexplained context
the phrase ``by Lemma~x'' should always be taken to mean
``by Lemma~x or its dual''.

We  note that, as with all matroid structure theory, in
one way or another it is all about connectivity.

Let $M$ be a matroid on ground set $E$ with rank function $r$.
The {\em connectivity function}
\index{connectivity function}
$\lambda_M$
\index{$\lambda_M$} 
of  $M$ is defined, for all subsets
$A$ of $E$, by $\lambda_M(A)=r(A)+r(E-A)-r(M)$. 
If the matroid is clear from the context then $\lambda_M(A)$
will be denoted by $\lambda(A)$.
We extend the notation to a partition $(A,E-A)$ of $E$, 
by defining $\lambda_M(A,E-A)=\lambda_M(A)$.

The set 
$A$ or the partition $(A,E-A)$ is $k$-{\em separating} 
\index{$k$-separating} if 
$\lambda(A)<k$. The partition $(A,E-A)$ is a $k$-{\em separation}
\index{$k$-separation} if $A$ is $k$-separating
and $|A|,|E-A|\geq k$. Fussing about the distinction
between 3-separating partitions and 3-separations 
leads to constipated prose, but the distinction is 
at times important, so we appear to be stuck with it.
The matroid 
$M$ is $k$-{\em connected} \index{$k$-connected} 
if it has no 
$(k-1)$-separations.
A $k$-separating set $A$, or $k$-separation $(A,E-A)$ is 
{\em exact} if $\lambda(A)=k-1$. A matroid is {\em connected}
if it is 2-connected.

It is immediate from the definition that the connectivity
function of a matroid is {\em symmetric}, that is,
$\lambda_M(X)=\lambda_M(E-X)$ for all subsets $X$ of $E$.
Moreover, one readily checks that if $r^*$ is the rank
function of the dual $M^*$ of the matroid $M$,
then $\lambda_M(X)=r(X)+r^*(X)-|X|$ for all subsets
$X$ of $E$. This establishes the next lemma.

\begin{lemma}
\label{dual-connectivity}
For any matroid $M$ we have $\lambda_M=\lambda_{M^*}$.
\end{lemma}

We freely use Lemma~\ref{dual-connectivity} without mention
throughout the paper. Another elementary
fact about $\lambda$ is that it is monotone under minors.

\begin{lemma}
\label{monotone}
Let $N$ be a minor of the matroid $M$. Then
$\lambda_N(A)\leq \lambda_M(A)$ for any subset $A$ of $E(N)$.
\end{lemma}

An easy rank calculation 
proves that 
$\lambda$ is {\em submodular},  that is 
$\lambda(X)+\lambda(Y)\geq\lambda(X\cup Y)+\lambda(X\cap Y)$ for all
$X,Y\subseteq E$. 
The submodularity of $\lambda$ is frequently used to establish
that certain sets have bounded connectivity. The next lemma
is an instance of this.

\begin{lemma}
\label{uncrossing}
Let $M$ be a matroid and let $X$ and $Y$ be subsets of $E(M)$
such that $\lambda(X)=\lambda(Y)=2$.
If $\lambda(X\cup Y)\geq 2$,
then $\lambda(X\cap Y)\leq 2$. In particular, if $M$ is 
$3$-connected and $\lambda(X)=\lambda(Y)=2$, the following hold.
\begin{itemize}
\item[(i)] If $|X\cap Y|\geq 2$, then $X\cup Y$ is $3$-separating.
\item[(ii)] If $|E(M)-(X\cup Y)|\geq 2$, then $X\cap Y$
is $3$-separating.
\end{itemize}
\end{lemma}

We make frequent use of Lemma~\ref{uncrossing} and 
often write {\em by uncrossing} \index{uncrossing}
to mean ``by an application
of Lemma~\ref{uncrossing}''.

\subsection*{Keeping Connectivity}
A key role that connectivity plays in matroid structure theory
is to eliminate degeneracies caused by low-order separations.
This is precisely the role played by connectivity when 
sufficient connectivity enables
us to bound the number of inequivalent representations of a matroid. 
To make inductive arguments possible it is
necessary to have theorems that enable us to remove elements
keeping a given type of connectivity. 
It seems that Tutte was the first to appreciate the need for
such results. Indeed a number of the results proved by Tutte
have become basic tools. The following is the most elementary.

\begin{lemma}
\label{keep1}
Let $M$ be a connected matroid. Then, for any element $e$
of $M$, either $M\ba e$ or $M/e$ is connected.
\end{lemma}

Let $n\geq 2$ be an integer. Recall that the {\em wheel}
\index{wheel}
$W_n$ 
\index{$W_n$}
is the graph consisting of  a cycle of length $n$
together with another vertex $v$ that is incident with
all of the vertices in the cycle. The {\em rim} 
\index{rim} edges
of $W_n$ are the edges in the cycle. The remaining edges
are the {\em spoke} \index{spoke}
edges. If we lapse and say that a matroid
is a wheel we mean that it is the cycle matroid of a 
wheel. The rim edges of  $W_n$ form a circuit-hyperplane
of $M(W_n)$. The {\em whirl} $W^n$
\index{$W^n$} 
is obtained from 
$M(W_n)$ by declaring this circuit hyperplane to be a basis; see
\cite[Chapter~8.4]{ox92}.
The next theorem of Tutte \cite{tu1} is fundamental.
\index{Wheels and Whirls Theorem}

\begin{theorem}[Tutte's Wheels and Whirls Theorem]
Let $M$ be a $3$-connected matroid. Then, unless $M$
is a whirl or the cycle matroid of a wheel, there is an
element $e$ of $M$ such that either $M\ba e$ or $M/e$
is $3$-connected.
\end{theorem}

A set $X$ of a matroid $M$ is a {\em parallel set} 
\index{parallel set} if
every 2-element subset of $X$ is a circuit. A 
{\em parallel class} 
\index{parallel class}
of $M$ is a maximal parallel set. Dually,
$X$ is a {\em series set} \index{series set}
of $M$ if every 2-element subset
of $M$ is a cocircuit. A {\em series class} 
\index{series class} 
of $M$ is a
maximal series set. A matroid $M$ is {\em $3$-connected up to 
parallel pairs} if whenever $(X,Y)$ is a 2-separation of $M$,
then either $X$ or $Y$ is a parallel pair and is 
{\em $3$-connected up to series pairs} if 
whenever $(X,Y)$ is a 2-separation of $M$,
then either $X$ or $Y$ is a series pair.
Bixby \cite{bi1} proved the next very useful result.

\begin{theorem}[Bixby's Lemma]
Let $e$ be an element of the $3$-connected matroid $M$.
Then either $M\ba e$ is $3$-connected up to series pairs
or $M/e$ is $3$-connected up to parallel pairs.
\end{theorem}

To say that a matroid $M$ has the matroid $N$ {\em as a minor}
means that we may delete or contract elements
from $M$ to obtain
a matroid equal to $N$. To say that $M$ has an $N$-{\em minor} means
that $M$ has a matroid isomorphic to $N$ as a minor.
It is often the case that we would like to keep connectivity 
and keep a minor. Typically isomorphism needs to be invoked.
The version of Seymour's Splitter Theorem \cite{se1}
below is not the
strongest possible. See Oxley \cite[Chapter~11]{ox92}
for a more detailed
discussion of the Splitter Theorem and its consequences.

\begin{theorem}[Seymour's Splitter Theorem]
Let $N$ be a $3$-connected matroid that is not a wheel or a whirl.
If $N$ is a proper minor of the $3$-connected matroid $M$,
then there is an element $e$ of $M$ such that either
$M\ba e$ or $M/e$ is $3$-connected with an $N$-minor.
\end{theorem}

\subsection*{Tutte's Linking Theorem}
\index{Tutte's Linking Theorem}
Let $M$ be a matroid and let $X$ and $Y$ be disjoint 
subsets of $E(M)$. We let
$\kappa_M(X,Y)=\min(\lambda_M(A):X\subseteq A\subseteq E(M)-Y)$.
\index{$\kappa_M(X,Y)$}
If the matroid $M$ is clear we abbreviate $\kappa_M$ to $\kappa$.
Intuitively $\kappa(X,Y)$ measures the connectivity between
$X$ and $Y$ provided by the rest of the matroid.
 
If $N$ is a minor of $M$ and $X,Y\subseteq E(N)$, 
then $\kappa_N(X,Y)\leq \kappa_M(X,Y)$. 
The next
theorem provides a good characterisation for 
$\kappa_M(X,Y)$. This theorem is, in fact, a generalisation
of Menger's Theorem.

\begin{theorem}[Tutte's Linking Theorem \cite{tu2}]
\label{tutte-linking}
Let $M$ be a matroid and let $X$ and $Y$ be disjoint
subsets of $E(M)$. Then there exists a minor 
$N$ on $X\cup Y$ such that 
$\lambda_N(X)=\kappa_M(X,Y)$.
\end{theorem}

The following lemma shows that if we apply Tutte's Linking
Theorem when $\lambda_M(X)=\kappa_M(X,Y)$, the resulting
minor $N$ satisfies $M|X=N|X$.

\begin{lemma}
\label{linking-coroll}
Let $N$ be a minor of a matroid $M$ and let
$X$ be a subset of $E(N)$. If $\lambda_M(X)=\lambda_N(X)$,
then $M|X=N|X$.
\end{lemma}

\subsection*{Local Connectivity} For subsets $X$ and $Y$ in $M$, the
{\em local connectivity} \index{local connectivity}
between $X$ and $Y$, denoted 
$\sqcap_M(X,Y)$, or $\sqcap(X,Y)$
\index{$\sqcap_M(X,Y)$} 
if the matroid is clear from
the context,
is defined by $\sqcap_M(X,Y)=r(X)+r(Y)-r(X\cup Y)$. 
Evidently $\sqcap_M(X,Y)=\lambda_{M|(X\cup Y)}(X,Y)$.
We denote $\sqcap_{M^*}(X,Y)$ by
$\sqcap^*_M(X,Y)$.
\index{$\sqcap^*_M(X,Y)$} 
The next lemma follows from an easy rank calculation.

\begin{lemma}
\label{pi-minor}
Let $x$ be an element of the matroid $M$ and let
$A$ and $B$ be disjoint subsets of $E(M)-\{x\}$. Then
\begin{itemize}
\item[(i)] $\sqcap_{M/x}(A,B)=\sqcap_M(A,B)$ if either
$x\not\in\cl(A\cup B)$ or $x\in \cl(A)$, but $x\not\in\cl(B)$;
\item[(ii)] $\sqcap_{M/x}(A,B)=\sqcap_M(A,B)+1$ if 
$x\in\cl(A\cup B)$ but 
$x\not\in\cl(A)$ and $x\not\in\cl(B)$; and
\item[(iii)] $\sqcap_{M/x}(A,B)=\sqcap_M(A,B)-1$ if
$x\in\cl(A)$ and $x\in\cl(B)$.
\end{itemize}
\end{lemma}

Let $\{A,B,C\}$ be a partition of the ground set $E$ of a matroid $M$.
Then $A$ and $B$ are {\em skew} 
\index{skew}
in $M$ if $r(A\cup B)=r(A)+r(B)$ and
are {\em coskew} 
\index{coskew}
in $M$ if they are skew in $M^*$. Equivalently,
$A$ and $B$ are skew if $\sqcap_M(A,B)=0$ and are coskew if
$\sqcap_M^*(A,B)=0$. 

Recall that sets $X$ and $Y$ in a matroid $M$ form a {\em modular 
pair\ } if $r(X)+r(Y)=r(X\cup Y)+r(X\cap Y)$.
An easy argument from duality proves 

\begin{lemma}
\label{coskew}
Let $M$ be a matroid and let $\{A,B,C\}$ be a partition of $E(M)$.
Then the following are equivalent.
\begin{itemize}
\item[(i)] $A$ and $B$ are coskew in $M$.
\item[(ii)] $r_{M/C}(A)+r_{M/C}(B)=r_{M/C}(A\cup B)$.
\item[(iii)] $r(A\cup C)+r(B\cup C)=r(M)+r(C).$
\item[(iv)] $E-A$ and $E-B$ form a modular pair in $M$.
\end{itemize}
\end{lemma}

The next lemma is proved in \cite{flower}.

\begin{lemma}
\label{lambda-meet}
For disjoint subsets $X$ and $Y$ of $M$,
$$\lambda(X\cup Y)=\lambda(X)+\lambda(Y)-\sqcap(X,Y)-\sqcap^*(X,Y).$$
\end{lemma}

We will say that $A$ and $B$ are {\em fully skew} \index{fully skew}
if they are both skew
and coskew. 
An immediate corollary of Lemma~\ref{lambda-meet} 
and the submodularity of
the connectivity function is

\begin{corollary}
\label{skew-coskew}
$A$ and $B$ are fully skew if and only if 
$\lambda(A\cup B)=\lambda(A)+\lambda(B)$.
\end{corollary}

\section{Structure Related to Connectivity}

\subsection*{Sequential and Equivalent $3$-separations}
Let $A$ be a set in a matroid $M$. The {\em coclosure\ }
\index{coclosure} 
$\cl^*(A)$
\index{cl$^*(X)$} 
of $A$ is the closure of $A$ in $M^*$. If 
$\cl^*(A)=A$, then $A$ is {\em coclosed\ } 
\index{coclosed} 
in $M$. Let 
$(A,\{x\},B)$ be a partition of $E(M)$.
Then $x\in\cl(A)$ if and only if $x$ is a loop of $M/A$. 
So by duality, $x\in\cl^*(A)$ if and only if $x$ is a 
coloop of $M\ba A=M|(A\cup \{x\})$. In other words, we have

\begin{lemma}
\label{cl-clstar}
If $(A,\{x\},B)$ is a partition of $E(M)$, then $x\in\cl^*(A)$
if and only if $x\notin \cl(B)$.
\end{lemma}

Throughout this paper we freely use properties of coclosure that
are obtained by dualising standard properties of closure. For
example, $x\in\cl^*_{M\ba a}(Z)$ if and only if 
$x\in\cl^*_M(Z\cup \{a\})$.

If $A$ is both closed
and coclosed in $M$, then $A$ is {\em fully closed}. 
\index{fully closed} 
The {\em full closure} 
\index{full closure}
of a set $A$, denoted  $\fcl_M(A)$, or $\fcl(A)$
\index{fcl$(X)$} 
if $M$ is clear from the context,
is the intersection of all the 
fully-closed sets containing $A$. Evidently the full closure 
is a closure operator in that it satisfies the
properties that $\fcl(A)\supseteq A$ and $\fcl(\fcl(A))=\fcl(A)$
for all subsets $A$ of $E(M)$. The set $A$ is {\em cohesive}
\index{cohesive}
if $E(M)-A$ is fully closed.

We use the notation $x\in\clstar(A)$
\index{cl$^{(*)}(X)$}
as a shorthand way of saying that
either $x\in\cl(A)$ or $x\in\cl^*(A)$. Note that 
$x\in\clstar(A)$ if and only if 
$\lambda(A\cup\{x\})\leq \lambda(X)$.

\begin{lemma}
\label{fclstuff}
Let $A$ be a set of elements of the matroid $M$. Then the following
are equivalent:
\begin{itemize}
\item[(i)] $B=\fcl(A)$;
\item[(ii)] $B$ is a maximal set containing $A$, 
for which there is an ordering 
$(a_1,\ldots,a_n)$ of $B-A$ such that 
$a_i\in\clstar(A\cup\{a_1,\ldots,a_{i-1}\})$
for all $i\in\{1,\ldots,n\}$.
\item[(iii)] $B$ is a maximal set containing $A$ for which there is
an ordering $(a_1,\ldots,a_n)$ of 
$B-A$ such that 
$\lambda(A\cup\{a_1,\ldots,a_i\})
\leq \lambda(A\cup\{a_1,\ldots,a_{i-1}\})$
for all $i\in\{1,\ldots,n\}$.
\end{itemize}
\end{lemma}

Via the full closure operator we can obtain an equivalence on 
$3$-separating sets of a 3-connected matroid $M$ as follows.
Say $\lambda(A)=\lambda(B)=2$. Then $A$ is {\em equivalent} 
\index{equivalent $3$-separating sets} 
to $B$, denoted $A\cong B$, 
if $\fcl(A)=\fcl(B)$. Say that $(A_1,A_2)$ and $(B_1,B_2)$
are exactly 3-separating partitions in $M$. Then $(A_1,A_2)$ is
{\em equivalent} to $(B_1,B_2)$, 
denoted $(A_1,A_2)\cong (B_1,B_2)$ if, 
for some ordering $(C_1,C_2)$
of $\{B_1,B_2\}$, we have $A_1\cong C_1$ and 
$A_2 \cong C_2$. 

Let $X$ be a 3-separating set of 
the 3-connected matroid $M$. Then $X$ is
{\em sequential\ } 
\index{sequential $3$-separating set}
if it has an ordering $(x_1,\ldots,x_n)$ such
that $\{x_1,\ldots,x_i\}$ is $3$-separating for all 
$i\in\{1,\ldots,n\}$. By the symmetry of the connectivity 
function, $X$ is sequential if and only if $\fcl(E(M)-X)=E(M)$.
Under the equivalence defined earlier one 
can regard sequential $3$-separating partitions
as being equivalent to
trivial 3-separating partitions. If $X$ is not sequential, then we 
say that $X$ is a {\em non-sequential\ } 
\index{non-sequential $3$-separating set}
3-separating set.
A 3-separation $(X,Y)$ of $M$ is {\em non-sequential\ } 
\index{non-sequential $3$-separation} if
neither $X$ nor $Y$ is sequential; otherwise $(X,Y)$
is {\em sequential}. Thus $(X,Y)$ is non-sequential if
$\fcl(X)\neq E(M)$ and $\fcl(Y)\neq E(M)$.

The next lemma gives a test for equivalence of $3$-separations.

\begin{lemma}
\label{non-sequential-equiv}
Let $(A_1,A_2)$ and $(B_1,B_2)$ be   $3$-separating
partitions of the 
$3$-connected matroid $M$. 
\begin{itemize}
\item[(i)] If $(A_1,A_2)$ is sequential, then 
$(A_1,A_2)\cong (B_1,B_2)$ if and only if there is an
ordering $(C_1,C_2)$ of $\{B_1,B_2\}$ such that
$\fcl(A_1)=\fcl(C_1)$ and $\fcl(A_2)=\fcl(C_2)$.
\item[(ii)] If $(A_1,A_2)$ is non-sequential, then 
$(A_1,A_2)\cong (B_1,B_2)$ if and only if there is an
ordering $(C_1,C_2)$ of $\{B_1,B_2\}$ such that
$\fcl(A_1)=\fcl(C_1)$.
\end{itemize}
\end{lemma}

\subsection*{Guts, Coguts, Blocking, Coblocking} 
Let $(A,\{x\},B)$ be a partition of the ground set of
the matroid $M$. 
Assume that $\lambda(A)=\lambda(A\cup\{x\})=k$.
Then either $x\in\cl(A)$, in which case we say that $x$ is in 
the {\em guts\ } 
\index{guts}
of $(A,B\cup\{x\})$, 
or $x\in\cl^*(A)$,
in which case we say that $x$ is in the 
{\em coguts\ } 
\index{coguts}
of $(A,B\cup\{x\})$. Note that, as $\lambda(A)=\lambda(A\cup\{x\})$,
if $x$ is in the guts of $(A,B\cup\{x\})$, then $x\in\cl(B)$
and if $x$ is in the coguts of $(A,B\cup\{x\})$, then 
$x\in\clstar(B)$.

Let $N$ be a minor of $M$ and $(A,B)$ be a partition of
$E(N)$ such that $\lambda_N(A)=k$. Then 
$(A,B)$ is {\em induced\ } 
\index{induced $k$-separation}
in $M$ if there is a partition
$(A',B')$ of $E(M)$ with $A\subseteq A'$ and $B\subseteq B'$
such that $\lambda_M(A')=k$. If $(A,B)$ is not induced
in $M$, then $(A,B)$ is {\em bridged\ } 
\index{bridged $k$-separation}
by $M$.

The case when $N$ is obtained by removing a single element
from $M$ is
of particular interest and has its own terminology.
Assume that $N=M\ba x$ and that $(A,B)$ is 
bridged by $M$. Then we say that $(A,B)$ is {\em blocked\ }
\index{blocked $k$-separation}
by $x$. On the other hand, if $N=M/x$ and $(A,B)$ is
bridged by $M$, then we say that $(A,B)$
is {\em coblocked\ } 
\index{coblocked $k$-separation}
by $x$. The next lemma follows from easy rank calculations.
Recall that $A$ is $k$-separating if $\lambda(A)<k$.

\begin{lemma}
\label{guts-stuff}
Let $M$ be a matroid and let $(A,\{x\},B)$ 
be a partition of $E(M)$ where $\lambda(A)=\lambda(A\cup\{x\})=k$.
Then the following are equivalent.
\begin{itemize}
\item[(i)] $x$ is in the guts of 
$(A,B\cup\{x\})$.
\item[(ii)] $x\in\cl(A)$ and $x\in\cl(B)$.
\item[(iii)] $x\notin\cl^*(A)$ and $x\notin\cl^*(B)$.
\item[(iv)] If $x$ is not a loop of $M$,
then $(A,B)$ is a $k$-separation of $M/x$ that is 
coblocked by $x$.
\end{itemize}
\end{lemma}

Dualising we obtain.

\begin{lemma}
\label{coguts-stuff}
Let $M$ be a matroid and let $(A,\{x\},B)$ 
be a partition of $E(M)$ where $\lambda(A)=\lambda(A\cup\{x\})=k$.
Then the following are equivalent.
\begin{itemize}
\item[(i)] $x$ is in the coguts of 
$(A,B\cup\{x\})$.
\item[(ii)] $x\in\cl^*(A)$ and $x\in\cl^*(B)$.
\item[(iii)] $x\notin\cl(A)$ and $x\notin\cl(B)$.
\item[(iv)] If $x$ is not a coloop of $M$,
then $(A,B)$ is a $k$-separation of $M\ba x$ that is 
blocked by $x$.
\end{itemize}
\end{lemma}

The fact that to coblock the $k$-separation $(A,B)$ of $M/x$ 
we must
have $x\in\cl_M(A)$ and $x\in\cl_M(B)$ is used many times
in this paper. We similarly use the dual observation on blocking,
although rather than observe that $x\in\cl^*_M(A)$
and $x\in\cl^*_M(B)$ we more typically make the 
equivalent observation that $x\notin\cl_M(A)$ and 
$x\notin\cl_M(B)$; no doubt this is because 
we are, after all, more habituated to the
closure operator.

\subsection*{Exposed $3$-separations}
Let $N$ be a $3$-connected minor of the 3-connected matroid
$M$. A $3$-separation $(A,B)$ of $N$ is {\em exposed\ }
\index{exposed $3$-separation}
in $N$ if every $3$-separating partition of $N$ that is equivalent
to $(A,B)$ is bridged in $M$. 
Again the case where $N$ is obtained by deleting or contracting
a single element is of most interest. If $(A,B)$ is 
exposed in $N$ and $N=M\ba x$ or $M/x$ for some element 
$x$, then we say that $(A,B)$ is {\em exposed by $x$}.
The next lemma is immediate.

\begin{lemma}
\label{expose-non-seq}
Let $M$ be a $3$-connected matroid with an element
$x$ such that $M\ba x$ is $3$-connected. If the $3$-separation
$(A,B)$ is exposed by $x$, then $(A,B)$ is non-sequential.
\end{lemma}

It is shown in \cite{upgrade}
that, if  $M$ is a $3$-connected matroid that is not a wheel or 
a whirl, then there is an element $x$ in $E(M)$ such that
either $M\ba x$ or $M/x$ is $3$-connected and does not expose
any $3$-separations. This result has an important application in
this paper and we state it formally as Theorem~\ref{unexpose1}
immediately prior to its use.

If we are looking from a different perspective we alter the
terminology. Let $x$ be an element of the $3$-connected
matroid $M$ such that $M\ba x$ is $3$-connected. 
Let $(A,B)$ be a $3$-separation of $M\ba x$. Then we
say that $(A,B)$ is {\em well-blocked\ } 
\index{well-blocked $3$-separation} 
by $x$ if 
every 3-separating partition of $M\ba x$ that is equivalent to 
$(A,B)$ is blocked by $x$. In other words, $(A,B)$ is 
well-blocked by $x$ if and only if $(A,B)$ is exposed by 
$x$.

\subsection*{Split Sets}
A set $B$ of a matroid $M$ is {\em split\ }
\index{split set}
if some partition $(B',B'')$ of $B$ into nonempty subsets has 
the property that $B'$ and $B''$ are fully skew. An element
$b\in B$ is {\em isolated\ }
\index{isolated element} 
in $B$ if $B-\{b\}$ and $\{b\}$ 
are fully skew.

\begin{lemma}
\label{split}
Let $M$ be a $3$-connected matroid. 
\begin{itemize}
\item[(i)] Let $(X,Y)$ be a $3$-separation of $M$. 
Then neither $X$ nor $Y$
is split.
\item[(ii)] Let $(X,Y)$ be an exact $4$-separation of $M$. 
Then $X$ is split
if and only if there is an element $x\in X$ that is isolated in $X$,
and this holds if and only if $\lambda(X-\{x\})=2$ for some element
$x\in X$.
\end{itemize}
\end{lemma}

We can regard split $4$-separating sets as being in some sense
degenerate. We omit the easy proof of the next lemma.

\begin{lemma}
\label{split.5}
Let $M$ be a $3$-connected matroid and $B$ be a set of elements of
$M$ such that $\lambda(B)=3$. If $b\in B$,
then the following are equivalent.
\begin{itemize}
\item[(i)] $b$ is isolated in $B$.
\item[(ii)] $b\notin\cl(B-\{b\})$ and $b\notin\cl^*(B-\{b\})$.
\item[(iii)] $b\in\cl(E(M)-B)$ and $b\in\cl^*(E(M)-B)$.
\end{itemize}
\end{lemma}

\begin{lemma}
\label{split1}
Let $z$ be an element of the $3$-connected matroid $M$ such that
$M/z$ is $3$-connected. Let $(R,B)$ be a $3$-separation of
$M/z$ that is coblocked by $z$. Assume that $R$ is split in 
$M$ with isolated element $x$. Then, in $M/z$, 
the $3$-separation $(R,B)$ is equivalent
to $(R-\{x\},R\cup\{x\})$ and $x$ is in the guts of $(R,B)$.
\end{lemma}

\begin{proof}
By Lemma~\ref{split.5}, $x\in\cl_M(E(M)-R)$, that is 
$x\in\cl_M(B\cup \{z\})$. Hence $x\in\cl_{M/z}(B)$, so that
$x$ is in the guts of $(R,B)$.
\end{proof}

As a consequence of 
Lemma~\ref{split1}, we obtain the
useful fact that exposed 3-separations in 
$M/z$ correspond to certain unsplit 4-separations
in $M$ and $M\ba z$.

\begin{corollary}
\label{unsplit} 
Let $M$ be a $3$-connected matroid, let $z$ be an
element of $M$ for which $M/z$ is $3$-connected, and
let $(R,B)$ is a $3$-separation of $M/z$
that is exposed by $z$. Then the following 
hold.
\begin{itemize}
\item[(i)] $(R,B\cup\{z\})$ is an unsplit $4$-separation of $M$.
\item[(ii)] If $M\ba z$ is $3$-connected, then 
$(R,B)$ is an unsplit $4$-separation of $M\ba z$.
\end{itemize}
\end{corollary}

\begin{proof}
Assume that $R$ is split in either $M$ or $M\ba z$. Let
$x$ be an isolated element. 
By Lemma~\ref{split1}, $x\in\cl_{M/z}(B)$
so that, in $M/z$, $(R,B)$ is equivalent to $(R-\{z\},B\cup \{z\})$.
But $R-\{z\}$ is not coblocked by $z$, that is
$\lambda_M(R-\{z\})=\lambda_{M\ba z}(R-\{z\})=2$. 
This contradicts the fact that $(R,B)$ is well coblocked by $z$. 
Both parts of the lemma follow from this contradiction.
\end{proof}

\section{Schematic Diagrams} 

Intuition for matroids
can be gained by giving geometric representations. Of
course geometric representations are impossible for matroids
of high rank. If a matroid has 2-separations or 
3-separations some attempt can be made to give geometric
insight by the use of {\em schematic diagrams}.
\index{schematic diagram} 
Such 
diagrams are not infallible and are never a substitute for 
logic, but we find them invaluable as an aid to 
intuition in conducting
research in matroid structure theory, and we believe that
they can aid the reader as well. There are no hard and
fast rules for their precise interpretation and
ambiguity can always threaten. With that warning we
give some examples to illustrate their use. Other examples are
scattered throughout the paper; indeed Figures~\ref{rank-5-free-swirl}
and \ref{swirl-like-flower} of the introduction are schematic
diagrams.

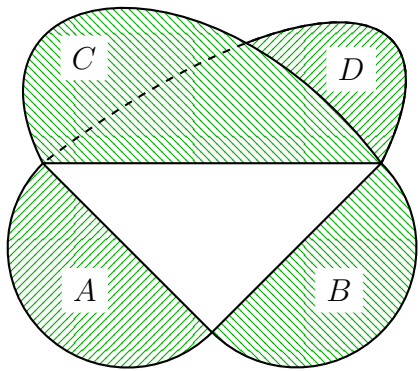
\begin{figure}
\begin{tikzpicture}[thick,line join=round]
	\coordinate (x) at (-2.25,0);
	\coordinate (y) at (2.25,0);
	\coordinate (y') at ($(y) + (-135:2)$);
	\coordinate (x') at ($(x) + (-45:2)$);
	\coordinate (z) at (intersection of x--x' and y--y');
	\filldraw[pattern color=lines,pattern=north east lines] let \p1 = ($(x) - (z)$),
		\n1 = {veclen(\x1,\y1)} in
			(z) arc (180+135:180-45:\n1/2);
	\filldraw[pattern color=lines,pattern=north west lines] let \p1 = ($(x) - (z)$),
		\n1 = {veclen(\x1,\y1)} in
			(y) arc (135+90-180:45-180:\n1/2);
	\filldraw[pattern color=lines,pattern=north east lines] (y) .. controls (3.5,2.5) and (1,2.5) .. (x); 
	\fill[white] (x) .. controls (-3.5,2.5) and (0,3) .. (y); 
	\filldraw[pattern color=lines,pattern=north west lines] (x) .. controls (-3.5,2.5) and (0,3) .. (y); 
	\draw[dashed] (y) .. controls (3.5,2.5) and (1,2.5) .. (x); 
	\draw (x) -- (y) -- (z) -- cycle;
	\node at ($(x)!0.5!(z) + (-135:0.8)$) [rectangle,fill=white,draw=white] {$A$};
	\node at ($(y)!0.5!(z) + (-45:0.8)$) [rectangle,fill=white,draw=white] {$B$};
	\node at (-1.7,1.4) [rectangle,fill=white,draw=white] {$C$};
	\node at (1.85,1.25) [rectangle,fill=white,draw=white] {$D$};
\end{tikzpicture}
\caption{A Schematic Diagram}\label{schematic1}
\end{figure}

Figure~\ref{schematic1} illustrates a matroid where the sets
$A$, $B$, $C$ and $D$ are 3-separating. Note that 
$(A\cup B,C\cup D)$ is also a 3-separation. The sets
$A$, $B$, $C$ and $D$ look like planes, but should
simply be interpreted simply as sets having rank at least
three. Note that $\sqcap(C,D)=2$ and $\sqcap(A,B)=1$,
so that $r(M)=(r(C)+r(D)-2)+(r(A)+r(B)-1)-2=r(A)+r(B)+r(C)+r(D)-5$.
In particular, if $A$, $B$, $C$ and $D$ are in fact planes,
then $r(M)=7$.

\begin{figure}
\begin{tikzpicture}[thick,scale=1.2]
	\coordinate (x) at (0,0);
	\coordinate (y) at (0,2);
	\coordinate (z) at (3,1);
	\coordinate (c) at ($(x)!0.5!(z)$);
	\filldraw[pattern color=lines,pattern=vertical lines] let \p1 = ($(c) - (z)$),
		\n1 = {veclen(\x1,\y1)},
		\n2 = {atan(1/3)},
		\n3 = {atan(3)} in
			(z) arc (-\n2:180-\n2:\n1)
			(x) arc (90-\n3-180:270-\n3-180:\n1);
	\draw (x) -- (y) -- (z) -- cycle;
	\coordinate[label=-135:$a$] (a) at ($(y)!0.5!(z)$);
	\coordinate[label=180:$b$] (b) at ($(x)!0.5!(y)$);
	\foreach \pt in {a,b} \fill[black] (\pt) circle (3pt);
	\node at (1.6,2.2) [rectangle,fill=white,draw=white] {$A$};
	\node at (1.6,2.2) {$A$};
	\node at (1.6,-0.2) [rectangle,fill=white,draw=white] {$B$};
	\node at (1.6,-0.2) {$B$};
\end{tikzpicture}
\caption{Elements in the guts and coguts}
\label{schematic2}
\end{figure}
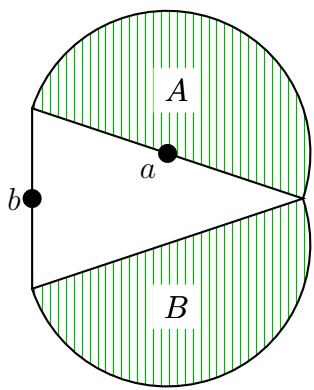

Figure~\ref{schematic2} illustrates a matroid with a 
$3$-separation $(A\cup\{a\},B\cup\{b\})$. Note that 
$a$ is in the guts of this 3-separation and $b$ is in the
coguts. This illustrates a case
of Lemmas~\ref{guts-stuff} and \ref{coguts-stuff}
as deleting $b$ gives a 
schematic diagram of a matroid with a 2-separation
$A\cup\{a\}$ and contracting $a$ corresponds to projection from
$a$ and reveals a 2-separation $(A,B\cup\{b\})$. Note
also that the 3-separating sets $A$, $A\cup\{a\}$,
and $A\cup\{a,b\}$ are equivalent.

\section{Blocking and Bridging Sequences}
\label{bridging-sequences}

Blocking sequences 
\index{blocking sequence}
for matroids 
were introduced by Geelen, Gerards and Kapoor 
in their proof of Rota's Conjecture for $GF(4)$ \cite{ggk}.
They have since proved to be a valuable
technique in matroid theory. In the application 
for \cite{ggk}, matroids are given by standard matrix representations
and hence come naturally with a fixed basis. Thus blocking 
sequences are defined relative to a fixed basis. In this section
we give a basis-free version of blocking sequences that we
call ``bridging sequences''. 

Let $N$ be a matroid with an exact $k$-separation $(X_1,X_2)$,
so that $\lambda_N(X_1,X_2)=k-1$ and let
$M$ be a matroid with an $N$ minor.
Recall that $(X_1,X_2)$ is {\em bridged\ } in $M$ if
\index{bridged}
$\kappa_M(X_1,X_2)\geq k$. 
Let $V=(v_0,v_1,\ldots,v_p)$ be an ordering of the elements of 
$E(M)-E(N)$ and let $S'$ and $T'$ denote the elements of $V$ with
even and odd indices respectively. Then
$V$ is a {\em bridging sequence\ }
\index{bridging sequence} 
for $(X_1,X_2)$ if
there is a permutation $(S,T)$ of $\{S',T'\}$ such that the following 
hold: 
\begin{itemize}
\item[(i)] $S$ is coindependent, $T$ is independent
and $N=M\ba S/T$;
\item[(ii)] if $i\in\{0,1,\ldots,p\}$, then 
$\lambda_M(X_1\cup\{v_0,\ldots,v_i\},X_2\cup\{v_{i+1},\ldots,v_p\})=k$; 
\item[(iii)] if $v_i\in S$, then 
$\lambda_{M\ba v_i}(X_1\cup\{v_0,\ldots,v_{i-1}\},
X_2\cup\{v_{i+1},\ldots,v_p\})=k-1$; and
\item[(iv)] if $v_i\in T$, then 
$\lambda_{M/v_i}(X_1\cup\{v_0,\ldots,v_{i-1}\},
X_2\cup\{v_{i+1},\ldots,v_p\})=k-1$.
\end{itemize}

\begin{lemma}
\label{last-gasp0}
Let $(X_1,X_2)$ be an exact $k$-separation for the 
matroid $N$ and let $M$ be a matroid having $N$ as a minor
such that there is an ordering $V=(v_0,v_1,\ldots,v_p)$
of $E(M)-E(N)$ that is a bridging sequence for 
$(X_1,X_2)$. Then 
$(X_1,X_2)$ is bridged in $M$. 
\end{lemma}

\begin{proof}
Assume that the lemma fails. Then there is a 
subset $W$ of $V$ such that $\lambda_M(X_1\cup W)<k$.
Assume that $W$ is  maximal with this property.
Let $i$ be the first integer such that $v_i\notin W$. 
Such an $i$ certainly exists.
By the definition of bridging sequence
$\lambda(X_1\cup \{v_0,v_1,\ldots,v_{i-1}\})=
\lambda(X_1\cup \{v_0,v_1,\ldots,v_{i}\})=k$.
But then, by uncrossing,
$\lambda(X_1\cup W\cup \{v_{i}\})<k$, contradicting the 
choice of $W$. 
\end{proof}

If $V$ is a bridging sequence for $(X_1,X_2)$
and $(S,T)$ is the partition of $V$ given in the 
definition of bridging sequence, then it is
easily seen that there is a basis $B$ of $M$ that contains
$T$ and avoids $S$. Readers familiar with blocking sequences
will observe that $V$ is a blocking sequence relative to this
basis. Thus properties of bridging sequences can be derived
from properties of blocking sequences. For example we
immediately have

\begin{lemma}
\label{last-gasp1}
Let $(X_1,X_2)$ be an exact $k$-separation in the matroid
$N$ and let $M$ be a matroid in which $(X_1,X_2)$ is bridged.
Then there is a minor $M'$ of $M$ that has a bridging sequence
for $(X_1,X_2)$.
\end{lemma}

If $(X_1,X_2)$ is bridged in $M$ and $M$ has a bridging
sequence $V$, then we say that $M$ is a {\em bridging matroid\ }
for $(X_1,X_2)$. 
\index{bridging matroid}
If no proper minor of $M$ has this property, then we say
that $M$ is  
a {\em minimal bridging matroid} for $(X_1,X_2)$ 
\index{minimal bridging matroid}
and that $V$ is a {\em minimal bridging sequence} for $(X_1,X_2)$.
\index{minimal bridging sequence}

The next lemma is essentially  \cite[Theorem~3.4]{gehlwh}.

\begin{lemma}
\label{bridging-0}
If $M$ is a minimal bridging matroid for the exact $k$-separation
$(X_1,X_2)$ of $N$, then there is a unique partition
$(S,T)$ of $E(M)-E(N)$ such that $N=M\ba S/T$. The set 
$S$ is coindependent and $T$ is independent. Moreover,
there exists an ordering of $E(M)-E(N)$ that is a bridging
sequence for $(X_1,X_2)$. 
\end{lemma}

Let $M$ be a bridging matroid for the exact $k$-separation
$(X_1,X_2)$ of $N$ and let
$V$ be a bridging sequence for $(X_1,X_2)$.
Let $S$ and $T$ be the associated partition of $E(M)-E(N)$
given in the definition of bridging sequence.
We refer to the elements of $S$ and $T$ as the {\em delete}
and {\em contract} elements of $V$ respectively. 
If we are viewing from another perspective, then we 
may refer to them
as {\em extension\ } and {\em coextension\ } elements.
Evidently $M^*$ is a bridging matroid for the exact
$k$-separation $(X_1,X_2)$ of $N^*$ and $V$ is a bridging
sequence for this $k$-separation.
If $Z\subseteq V$, then we denote the matroid
$M\ba (S-Z)/(T-Z)$ by $N[Z]$. In particular we will often
denote the bridging matroid $M$ by $N[V]$.

We will need just a few properties of bridging sequences.

\begin{lemma}
\label{bridge-klonal1}
Assume that $V=(v_0,\ldots,v_p)$ is a bridging sequence for the
$k$-separation $(X_1,X_2)$ of the matroid $M$. 
\begin{itemize}
\item[(i)] If $v_i$ is a delete
element of $V$, then $v_i\notin\cl_{N[v_0,\ldots,v_i]}(X_2)$.
\item[(ii)] If $v_i$ is a contract
element of $V$, then $v_i\notin\cl^*_{N[v_0,\ldots,v_i]}(X_2)$.
\end{itemize}
\end{lemma}

\begin{proof}
Assume otherwise. Then 
$v_i\in\cl_{N[v_1,\ldots,v_p]}(X_2\cup\{v_{i+1},\ldots,v_p\})$,
contradicting the fact that $v_i$ blocks the $k$-separation
$(X_1\cup\{v_1,\ldots,v_{i-1}\},X_2\cup\{v_{i+1},\ldots,v_p\})$
of the matroid $N[v_1,\ldots,v_p]\ba v_i$.
\end{proof}

The proof of the next lemma is even easier and we omit it.

\begin{lemma}
\label{bridge-basic}
Let $V=(v_0,v_1,\ldots,v_p)$ be a bridging sequence for the 
$k$-separation
separation $(X,Y)$ of $N$. Say $i<p$ 
$V$. Then, in $N[v_0,\ldots,v_i]$, we have
$v_i\in\cl(X\cup\{v_0,\ldots,v_{i-1}\})$ and 
$v_i\in\cl^*(X\cup\{v_0,\ldots,v_{i-1}\})$.
\end{lemma}

We will also use the next technical lemma.

\begin{lemma}
\label{last-bridge}
Let $(X_1,X_2)$ be an exact $k$-separation of the matroid
$N$ and let $V$ be a bridging sequence for $(X_1,X_2)$
with associated bridging matroid $N[V]$.
If $(V_1,V_2)$ is a partition of $V$
and there is no $i\in\{0,1,\ldots,p\}$ such that
$V_1=\{v_0,v_1,\ldots,v_i\}$, then 
$\lambda_M(X_1\cup V_1,X_2\cup V_2)>k$.
\end{lemma}

\begin{proof}
Assume that $|V|=2$, so that $V=(v_0,v_1)$.
Consider $\lambda_{N[V]}(X\cup\{v_1\},Y\cup\{v_0\})$.
We may assume that $v_1$ is a delete element, so that 
$v_0$ is a contract element.
We have $\lambda_{N[V]\ba v_1}(X\cup \{v_0\},Y)=k-1$.
But $v_0\in\cl_{N[V]\ba v_1}(X)$ and
$v_0\notin\cl_{N[V]\ba v_1}(Y)$, so that
$\lambda_{N[V]\ba v_1}(X,Y\cup\{v_0\})=k$.
Now $v_1\notin\cl_{N[V]}(X)$, so 
$\lambda_{N[V]}(X\cup\{v_1\},Y\cup\{v_0\})=k+1$.
Thus the lemma holds in this case.

Assume that $|V|>2$. Under the hypotheses of the lemma,
there is an 
$i\in\{0,1,\ldots,p\}$ such that $v_i\in V_1$,
and $v_{i-1}\in V_2$. We may assume that 
$v_{i-1}$ is a delete element of $V$.
Let $M'$ be the minor obtained
by contracting all of the delete elements from
$V-\{v_i,v_{i-1}\}$ and deleting all of the contract elements
from $V-\{v_i,v_{i-1}\}$. It is an easy consequence
of the definition of bridging sequence 
and Tutte's Linking Lemma that 
$\kappa_{M'}(X_1,X_2)=k$. Indeed $(v_{i-1},v_i)$ is a bridging
sequence for the $k$-separation $(X_1,X_2)$ of 
$M'\ba v_{i-1}/v_i$. Thus
$\lambda_{M'}(X_1\cup \{v_i\},X_2\cup \{v_{i-1}\})=k+1$,
and $\lambda_{N[V]}(X_1\cup V_1,X_2\cup V_2) >k$ as required.
\end{proof}

\begin{lemma}
\label{3-bridge}
Let $(P_1,P_2,P_3)$ be a partition of the ground set of a matroid
$N$ such that $\lambda(P_1)=\lambda(P_3)=\kappa(P_1,P_3)=k-1$.
Let $(C,D)$ be a partition of $P_2$, where $C$ is independent,
$D$ is coindependent and $\lambda_{N/C\ba D}(P_1,P_3)=k-1$.
Let $V$ be a bridging sequence for $(P_1,P_2\cup P_3)$
such that $(P_1\cup P_2,P_3)$ is also bridged in $N[V]$.
Then the exact $k$-separation $(P_1,P_3)$
of $M/C\ba D$ is bridged in $N/C\ba D$ and $V$ is a bridging
sequence for this $k$-separation.
\end{lemma}

\begin{proof}
Say $V=(v_0,v_1,\ldots,v_p)$. Let $e$ be an element of $C$.

\begin{sublemma}
\label{3-bridge1}
If $i\in\{0,1,\ldots,p\}$,
then $\lambda_{N[V]/e}(P_1\cup\{v_0,v_1,\ldots,v_i\})=k$
and $\lambda_{N[V]/e}(P_3\cup\{v_{i+1},v_{i+2},\ldots,v_p\})=k$.
\end{sublemma}

\subproof
Assume that $e\in\cl_{N[V]}(P_1\cup V)$. Then
$e\in\cl_N(P_1)$ so that $\lambda_{N/e}(P_1)<\lambda_{N}(P_1)=k-1$,
contradicting the fact that $\lambda_{N/C\ba D}(P_1)=k-1$.
Thus $e\notin\cl_{N[V]}(P_1\cup V)$,
and similarly $e\notin\cl_{N[V]}(P_3\cup V)$.

It follows from the above that
$r_{N[V]/e}(P_1\cup\{v_0,v_1,\ldots,v_i\})
=r_{N[V]}(P_1\cup\{v_0,v_1,\ldots,v_i\})$,
and it follows that
$\lambda_{N[V]/e}(P_1\cup\{v_0,v_1,\ldots,v_i\})
=\lambda_{N[V]}(P_1\cup\{v_0,v_1,\ldots,v_i\})=k$.
Similarly
$\lambda_{N[V]/e}(P_3\cup\{v_{i+1},v_{i+2},\ldots,v_p\})=k$.
\end{proof}

We now show that $V$ is a bridging sequence for 
$(P_1,P_2\cup P_3-\{e\})$
in $N[V]/e$. Let $T$ be the set of contract elements of $V$.
As $C$ is independent in $N$, we see that $C\cup T$ is independent
in $N[V]$, so that $T$ is independent in $N[V]$. Thus
property (i) of bridging sequences holds. Property (ii)
follows from \ref{3-bridge1}. Say $i>1$ and $v_i$
is a delete element of $V$.
Then $k-1=\lambda_{N[V]\ba v_i}(P_1\cup\{v_0,v_1,\ldots,v_{i-1}\})
\geq  \lambda_{N[V]/e\ba v_i}(P_1\cup\{v_0,v_1,\ldots,v_{i-1}\})
\geq \lambda_{N/C\ba D}(P_1)=k-1$. 
Thus (iii), and similarly (iv), also hold and $V$ is indeed a 
bridging sequence for $P_1$ in $N[V]/e$.

It follows from \ref{3-bridge1} and an argument similar to
that of Lemma~\ref{last-gasp0} that $P_3$ is bridged
in $N[V]/e$. The lemma now follows from an obvious induction.
\end{proof}

\section{Special Structures} 

In this section we review properties of certain highly
structured matroids and sets in matroids
that play an important role in this paper. We begin by looking
at 3-separating sets.

\subsection*{Sequential $3$-separators}
Let $M$ be a $3$-connected matroid on $E$ and let 
$A$ be a 3-separating subset of $E$. Recall that
$(A,E-A)$ is {\em sequential} if either $\fcl(A)=E$ or 
$\fcl(E-A)=E$. If the latter case holds we say that 
$A$ is a {\em sequential} $3$-separator.

Evidently any subset of $E$ with at most two elements is a 
sequential $3$-separator. If $A$ is a sequential 3-separator,
then there is an ordering $(a_1,a_2,\ldots,a_n)$ of $A$
such that, for all $i\in \{1,2,\ldots,n\}$, the set
$\{a_1,a_2,\ldots,a_i\}$ is 3-separating. Such an ordering
is said to be a {\em sequential ordering}
\index{sequential ordering} 
of $A$. 
Sequential orderings are typically far from unique. The 
next few lemmas summarise some elementary properties,
most of which
follow immediately from definitions. In  all of 
the lemmas $A$ is a sequential $3$-separator of the $3$-connected
matroid $M$ on $E$ and  $(a_1,a_2,\ldots,a_n)$ is a
sequential ordering of $A$. Recall that a 
{\em triangle}
\index{triangle} 
of a matroid is a $3$-element circuit and a {\em triad}
\index{triad} 
is a 3-element cocircuit.

\begin{lemma}
\label{seq-ord}
\begin{itemize}
\item[(i)] If $n\geq 2$, then $A\subseteq \fcl(\{a_1,a_2\})$.
\item[(ii)] If $n\geq 3$, then $\{a_1,a_2,a_3\}$
is either a triangle or a triad.
\item[(iii)] If $i\in\{1,2,\ldots,n\}$, then 
$a_i\in\clstar(\{a_{i+1},a_{i+2},\ldots,a_n\}\cup (E-A))$.
\item[(iv)] If $i\in\{2,3,\ldots,n-1\}$, then
$a_{i+1}\in\clstar(\{a_1,\ldots,a_i\})$.
\end{itemize}
\end{lemma}

\begin{lemma}
\label{seq-3-con}
If $A$ is fully closed and $|A|\geq 4$, then either
$M\ba a_n$ or $M/a_n$ is $3$-connected. In particular;
\begin{itemize}
\item[(i)] if $a_n\in\cl(E-A)$, then $M\ba a_n$
is $3$-connected, and
\item[(ii)] if $a_n\in\cl^*(E-A)$, then $M/a_n$
is $3$-connected.
\end{itemize}
\end{lemma}  

If $A$ is a sequential $3$-separation set, and 
$B$ is an exactly $3$-separating subset of $A$,
then it is not necessarily the case that $\fcl(B)$ contains
$A$, but the behaviour of such sets is not wild.

\begin{lemma}
\label{sub-seq-3-sep}
If $B$ is a $3$-separating subset of $A$, then $B$ is 
sequential.
\end{lemma}

\begin{proof}
As $A$ is sequential, $\fcl(E(M)-A)=E(M)$. As
$(E(M)-B)\supseteq (E(M)-A)$ we have
$\fcl(E(M)-B)=E(M)$. Thus $B$ is sequential.
\end{proof}

Finally we note

\begin{lemma}
\label{delete-seq}
If $x\in A$ and $M\ba x$ is $3$-connected, then
$A-\{x\}$ is a sequential $3$-separator of $M\ba x$.
\end{lemma}

\subsection*{Fans} 
Let 
$F$ be a set of elements of the $3$-connected matroid $M$.
Then $F$ is a {\em fan} of $M$
\index{fan}
if it has an ordering 
$(f_1,f_2,\ldots,f_n)$ such that
\begin{itemize}
\item[(i)] for all $i\in\{1,2,\ldots,n-2\}$, the triple 
$\{f_i,f_{i+1},f_{i+2}\}$ is either a triangle or a triad, and
\item[(ii)] if $i\in\{1,2\ldots,n-3\}$ then $\{f_i,f_{i+1},f_{i+2}\}$
is a triangle if and only if  $\{f_{i+1},f_{i+2},f_{i+3}\}$ is a triad.
\end{itemize}
An ordering of a fan satisfying (i) and (ii) is a
{\em fan ordering} of $F$ and we will refer to 
such an ordered set as an 
{\em ordered fan}.
\index{ordered fan} 
At times we may blur the distinction between
a fan and an ordered fan.
As triads and triangles are interchanged under duality,
a fan in $M$ is also a fan in $M^*$. 

\begin{lemma}
\label{fanny1}
If $S$ is a fan of a $3$-connected matroid and
$|E(M)-F|\geq 2$, then $F$ is a sequential $3$-separating
set and any fan ordering is a sequential ordering of $F$.
\end{lemma}

\begin{proof}
The lemma clearly holds if $|F|\leq 2$. Assume that
$|F|\geq 3$. Let $(f_1,f_2,\ldots,f_n)$ be a fan ordering of
$F$. Then $F-\{f_n\}$ is a fan and we may assume that the
lemma holds for $F-\{f_n\}$. Now $\{f_{n-2},f_{n-1},f_n\}$
is either a triangle or triad of $M$. In the former
case $f_n\in\cl(F-\{f_n\})$ and in the latter
$f_n\in\cl^*(F-\{f_n\})$. The lemma follows by induction.
\end{proof}

Note that a fan has many sequential orderings; indeed
any triangle or triad of a fan
can initiate such an ordering. Most
are not fan orderings.
If $F=(f_1,f_2,\ldots,f_m)$ is an ordered fan, 
and $i\in\{1,2,\ldots,m\}$,
then we say that $\{f_1,f_2,\ldots,f_i\}$ is an {\em initial section}
of $F$, and $\{f_{i+1},f_{i+2},\ldots,f_m\}$ 
is a {\em terminal section} of $F$.

\begin{lemma}
\label{fan1}
Let $M$ be a $3$-connected matroid 
with at least
four elements. Then $E(M)$ is a fan if and only if $M$ is a 
wheel or a whirl.
\end{lemma}

Thus fans can,
in general,  be thought of as
partial wheels or whirls.
Once we are past degeneracies called by small size, the 
structure of fans becomes quite canonical.
Subsets of size less than two are trivially fans. 
Fans of size three are either triangles or triads and any ordering
is a fan ordering.  
Assume that $F$ is a fan such that $|E(M)-F|\geq 2$.
If $|F|\geq 5$, and $(f_1,f_2,\ldots,f_n)$
is a fan ordering of $F$,
then the only other fan ordering of $F$ is to reverse the order
of the indices.
Fans with four elements are not quite canonical.
If $(f_1,f_2,f_3,f_4)$ is a fan ordering of $F$,
then $(f_1,f_3,f_2,f_4)$ is also a fan ordering of $F$.
This is a fact that is, at times, irritating in the minutiae 
of arguments.

We generalise terminology for elements of wheels and 
whirls to fans as follows.
Let $(f_1,f_2,\ldots,f_n)$ be a fan ordering of a fan $F$
with at least five elements. If $\{f_1,f_2,f_3\}$ is a triangle,
then the elements of $(f_1,f_2,\ldots,f_n)$
with odd indices are {\em spoke} 
\index{spoke element of a fan}
elements
and the elements with even indices are {\em rim} 
\index{rim element of a fan}
elements of
$F$. If $\{f_1,f_2,f_3\}$ is a triad, then the elements
with odd indices are {\em rim} elements and the elements
with even indices are {\em spoke} elements. 
Evidently the above labelling 
of elements is independent of the fan ordering. If
$(f_1,f_2,f_3,f_4)$ is a $4$-element fan, then $f_1$
is a spoke or a rim element according as to whether
$\{f_1,f_2,f_3\}$ is a triangle or a triad respectively.
This gives $f_1$ and $f_4$ labels, but we cannot assign
canonical labels to $f_2$ and $f_3$.  
 
One also sees that any fan with at least 
four elements has well-defined {\em end\ } elements and
{\em internal\ } elements in an obvious way.
We now recall a basic result of Tutte \cite{tu1}. 

\begin{lemma}[Tutte's Triangle Lemma]
\label{triangle-lemma}
Let $T$ be a triangle of a $3$-connected matroid with at
least four elements. If $T$ is not contained in a fan with
at least four elements, then there are elements $t_1,t_2\in T$
such that $M\ba t_1$ and $M\ba t_2$ are both $3$-connected.
\end{lemma}

Thus there is at most one element of a triangle
whose deletion from $M$ destroys 3-connectivity unless
the triangle is in a larger fan.
We can never remove an internal element of a fan to keep
3-connectivity. Nonetheless we can
get close. The next lemma is a straightforward consequence of 
Bixby's Lemma.

\begin{lemma}
\label{bixby-fan}
Let $M$ be a $3$-connected matroid; let 
$F$ be a fan of $M$ with at least five elements; let 
$(f_1,f_2,\ldots,f_n)$ be a fan ordering of $F$; and let
$f_i$ be an internal element of $F$.
\begin{itemize}
\item[(i)] If $f_i$ is a spoke element of $F$, then
$M\ba f_i$ is $3$-connected up to the single series pair
$\{f_{i-1},f_{i+1}\}$.
\item[(ii)] If $f_i$ is a rim element of $F$,
then $M/f_i$ is $3$-connected up to the single
parallel pair $\{f_{i-1},f_{i+1}\}$.
\end{itemize}
\end{lemma}

A fan is {\em maximal} if it is not contained in a larger fan.
If $M$ is not a wheel or a whirl, then any maximal fan $F$ is
exactly 3-separating. The next lemma is not quite a special
case of Lemma~\ref{seq-3-con} as a maximal fan need not be
fully closed.

\begin{lemma}
\label{fan-end}
Let $M$ be a matroid that is not a wheel or a whirl and let
$F$ be a maximal fan of $M$ with at least four elements. Let
$f$ be an end of $F$.
\begin{itemize}
\item[(i)] If $f$ is a spoke element of $F$, then $M\ba f$
is $3$-connected.
\item[(ii)] If $f$ is a rim element of $F$, then 
$M/f$ is $3$-connected.
\end{itemize}
\end{lemma}

\begin{figure}
\begin{tikzpicture}[thick,line join=round]
	\coordinate (x) at (0,0);
	\coordinate (y) at (1,0);
	\coordinate (z) at (4,0);
	\filldraw[pattern color=lines,pattern=north west lines] (0,0) arc (-180:0:2);
	\draw (x) -- (z);
	\coordinate[label=45:{\textcolor{labels}{$f_6$}}] (f6) at (3,0);
	\coordinate[label=135:{\textcolor{labels}{$f_2$}}] (f2) at (1,2);
	\coordinate[label=left:{\textcolor{labels}{$f_1$}}] (f1) at (1,1);
	\coordinate[label=above:{\textcolor{labels}{$f_3$}}] (f3) at (2,2);
	\coordinate[label=45:{\textcolor{labels}{$f_4$}}] (f4) at (3,2);
	\coordinate[label=right:{\textcolor{labels}{$f_5$}}] (f5) at (3,1);
	\draw (y) -- (f2) -- (f4) -- (f6);
	\foreach \pt in {f1,f2,f3,f4,f5,f6} \fill[black] (\pt) circle (3pt);
\end{tikzpicture}
\caption{A $6$-element Fan}
\label{fan-piccy}
\end{figure}
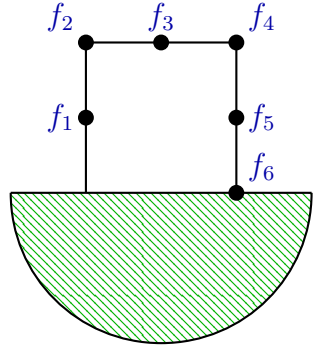

Finally we note that
while fans in matroids generalise the eponymous 
structures in graphs, it can be misleading for
matroidal arguments to visualise them graphically.
Figure~\ref{fan-piccy} is a schematic diagram of a
6-element fan in a matroid. The rim elements are
$\{f_1,f_3,f_5\}$ and the spoke elements 
are $\{f_2,f_4,f_6\}$.

\subsection*{Swirls}
Let $n\geq 3$ be an integer and
$N$ be a simple matroid whose ground set
is the disjoint union of a basis $B=\{b_1,b_2,\ldots,b_n\}$ 
and sets $P_1=\{p_1,q_1\},P_2=\{p_2,q_2\},\ldots,P_n=\{p_n,q_n\}$
such that, for all $i\in\{1,2,\ldots,n-1\}$, we have
$P_i\subseteq \cl(\{b_i,b_{i+1}\})$ and 
$P_n\subseteq \cl(\{b_n,b_1\})$. Then $M=N\ba B$ is a 
rank-$n$ {\em swirl}. 
\index{swirl}
We say
that $\{P_1,P_2,\ldots,P_n\}$ is the
set of {\em legs} 
\index{legs of a swirl}
of the swirl. It is easily seen that swirls are
3-connected.

Evidently $P_i\cup P_{i+1}$
is a circuit of $M$ for all $i$ in the cyclic order on
$\{1,2,\ldots,n\}$. Otherwise, if $C$ is a non-spanning circuit
of a swirl, then $C$ is a transversal of the legs. If $M$
has no such non-spanning circuits, then $M$ is the rank-$n$
{\em free swirl}
\index{free swirl} 
and we denote it by $\Delta_n$.
\index{$\Delta_n$}
As noted in the introduction one can obtain $\Delta_n$
by placing the elements of $P_i$ freely on the
line spanned by $\{b_i,b_{i+1}\}$. Figure~\ref{rank-5-free-swirl}
in the introduction illustrates $\Delta_5$. Evidently
$\Delta_n$ is unique up to isomorphism. Note that
$\Delta_n$ is self dual.

It is shown in \cite{oxvewh96} that if $q$ is a prime
power that exceeds five and is not of the 
form $2^p$, where $2^p-1$ is prime, then
$\Delta_n$ has at least $2^n$ inequivalent representations
over $GF(q)$. In particular, for a prime field
$GF(p)$ of size at least
seven, we can obtain free swirls with an arbitrary number of
inequivalent representations. Free swirls are structures
that necessarily need to be dealt with in understanding
inequivalent representations of matroids over prime fields.

\subsection*{Spikes} 
\begin{figure}
\begin{tikzpicture}[thick,draw=labels]
	\coordinate (a) at (-3,0);
	\coordinate (b) at (2,1);
	\coordinate (c) at (2,-0.7);
	\coordinate (d) at (3,0);
	\coordinate (e) at (3,-1.7);
	\draw[loosely dashed] (b) -- (c) -- (e) -- (d) -- cycle;
	\draw (b) -- (a) -- (c);
	\draw (d) -- (a) -- (e);
	\coordinate (p1) at ($(a)!0.5!(b)$);
	\coordinate (p2) at ($(a)!0.75!(b)$);
	\coordinate (p3) at ($(a)!0.7!(d)$);
	\coordinate (p4) at ($(a)!0.9!(d)$);
	\coordinate (p5) at ($(a)!0.5!(c)$);
	\coordinate (p6) at ($(a)!0.7!(c)$);
	\coordinate (p7) at ($(a)!0.7!(e)$);
	\coordinate (p8) at ($(a)!0.9!(e)$);
	\foreach \pt in {p1,p2,p3,p4,p5,p6,p7,p8} \fill[labels] (\pt) circle (3pt);
\end{tikzpicture}
\caption{Illustration of $\Lambda_4$}
\label{4-spike}
\end{figure}
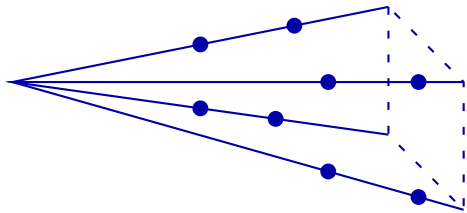

Let $n\geq 3$ be an integer and
let $N$ be a rank-$n$ matroid with ground set
$\{t,p_1,q_1,p_2,q_2,\ldots,p_n,q_n\}$ such that
\begin{itemize}
\item[(i)] $\{t,p_i,q_i\}$ is a triangle for all 
$i\in\{1,2,\ldots,n\}$, and
\item[(ii)] $r(\cup_{j\in J}\{a_j,b_j\})=|J|+1$ for 
every proper subset $J$ of $\{1,2,\ldots,n\}$.
\end{itemize}
Then the matroid $M=N\ba t$ is a rank-$n$ {\em spike}.
\index{spike}
Each pair $\{p_i,q_i\}$ is a {\em leg\ }
\index{leg of a spike} 
of the spike. 
For distinct $i,j\in\{1,2,\ldots,n\}$, the
set $\{p_i,q_i,p_j,q_j\}$ is a circuit of $M$.
As with swirls, any other non-spanning circuit of a 
spike is a transversal of the legs. If all transversals of the
legs are independent then $M$ is the rank-$n$
{\em free spike} 
\index{free spike}
and we denote is by $\Lambda_n$. 
\index{$\Lambda_n$}
Alternatively one can obtain $\Lambda_n$
by taking $M(K_{2,n})$, a matroid that has rank $n+1$,
and truncating it to rank $n$.
Note that
spikes are 3-connected and that $\Lambda_n$ is self dual.
Figure~\ref{4-spike} illustrates $\Lambda_4$.

Spikes in general have turned out to be fundamental 
in matroid theory, frequently as sources of counterexamples to 
superficially reasonable conjectures. 
Free spikes are representable over all
non-prime fields and it is shown in \cite{oxvewh96} that, for
every non-prime field with more than four elements,
$\Lambda_n$ has at least $2^{n-1}$ inequivalent representations.
Fortunately, it is shown in \cite[Lemma~11.6]{totally-free}
that if $p\geq 3$ is prime, then
$\Lambda_p$ is not $GF(p)$-representable so, if our eventual
goal is to understand inequivalent representations over $GF(p)$,
we can exclude $\Lambda_p$ from consideration.

Other properties of spikes and swirls follow from the
fact that the members of the partition of such a matroid into its legs
form the petals of a spike-like or swirl-like flower respectively.
Indeed, these are extremal examples of matroids with
many pairwise-inequivalent mutually-crossing 3-separations
and this is part of the reason for the key role that they
play in matroid structure theory.

\subsection*{Quads}
A {\em quad} of a matroid is a $4$-element set that is both
a circuit and a cocircuit.  
It is easily checked that if
$D$ is a quad, then $\lambda(D)=2$ so that
$D$ is exactly $3$-separating. It is immediate
from the definition that if $D$ is a quad in 
$M$, then $D$ is also a quad in $M^*$. Note that
the complement in $M$ of a quad is fully closed so that
quads are cohesive. Let $M$ be  a 3-connected matroid with at least
five elements. Then any 3-separating set of size three is
sequential as it is either a triangle or a triad. 
If $D$ is a non-sequential 3-separating set with 
four elements, then it is easily checked that $D$ is a quad.
It is no doubt the fact that quads are the unique minimim-sized
non-sequential 3-separating sets that account for their
frequent appearance in proofs in this paper. 

Note also that the union of any 
pair of legs of a spike is a quad, while the union of any consecutive
pair of legs of a swirl is a quad.

\chapter{Flowers}
\label{flower-chapter}

Flowers arise when a matroid
has crossing $3$-separations. They turn out to be fundamental
structures. Their study was initiated in \cite{flower}
and continued in \cite{flower2}. In these papers flowers
were defined only for $3$-connected matroids.
The definition was extended to arbitrary matroids and,
indeed, to structures that arise from higher-order crossing
separations, by Aikin and Oxley \cite{ai-ox}. This chapter
begins by reviewing material from \cite{flower,ai-ox} and then
developing further material that will be needed for
this paper.

\section{Definition and Basic Properties}

We begin by defining a flower in a way that is slightly
less restrictive than that of \cite{flower}, 
and somewhat less general than that of \cite{ai-ox}.
Let $M$ be a connected matroid, and
let $\FF=(P_1,\ldots,P_n)$ be a partition of $E(M)$. Then
$\FF$ is a {\em flower}
\index{flower} 
in $M$ with {\em petals}
\index{petal of a flower}
$P_1,\ldots,P_n$ if  the following hold.
\begin{itemize}
\item[(i)] If $n>1$, then $\lambda(P_i)=2$ for all 
$i\in\{1,2,\ldots,n\}$.
\item[(ii)] If $n>2$,
$\lambda(P_i\cup P_{i+1})=2$
for all $i\in\{1,2,\ldots,n\}$ where all subscripts are interpreted modulo $n$.
\item[(iii)] If $(X,Y)$ is a 2-separation of $M$, then for some 
petal $P_i$, either $X$ or $Y$ is a subset of $P_i$.
\end{itemize}

If $M$ is $3$-connected, then the definition given here is
precisely the definition given in \cite{flower}. This is 
the case we are usually interested in, but, at
times, it will facilitate arguments to allow flowers
in matroids that are not 3-connected. 

The ordering of the petals of a flower is always the 
cyclic order, so that subscripts should always
be interpreted modulo $n$. With this in mind we say that 
a set of petals of $\FF$ is {\em consecutive}
\index{consecutive set of petals} 
if it is of the form
$\{P_i,P_{i+1},\ldots,P_{i+k}\}$, for some $i,k\in\{1,2,\ldots,n\}$.
The flower $\FF$ is a {\em daisy} 
\index{daisy}
if the union of a set of petals
is $3$-separating if and only if the set of petals is consecutive.
The flower $\FF$ is an {\em anemone} 
\index{anemone}
if the union of any set of petals 
is 3-separating. The next theorem is a special case
of \cite[Theorem~1.1]{ai-ox}.

\begin{theorem}
\label{anemone-daisy}
If $\FF$ is a flower in the connected matroid $M$,
then $\FF$ is either an anemone or a daisy.
\end{theorem}

We next introduce a structure related to flowers. Before doing
that we settle some terminology for ordered partitions.
Let $\PP=(P_1,P_2,\ldots,P_n)$ be an ordered partition of 
a set $S$. Then the ordered partition 
$\QQ=(Q_1,Q_2,\ldots,Q_m)$ is a {\em concatenation}
\index{concatenation of a flower} 
of $\PP$ if there are indices $0=k_0<k_1<\cdots<k_m=n$
such that $Q_i=P_{k_{i-1}+1}\udots P_{k_i}$ for 
$i\in \{1,\ldots,m\}$. If $\QQ$ is a concatenation of
$\PP$, then $\PP$ is a {\em refinement}
\index{refinement of a flower} 
of $\QQ$.
Also, in any unexplained context we allow members of 
a partition to be empty 
sets and will 
often abuse notation and denote a singleton set $\{q\}$
by $q$.

Let $\QQ=(Q_1,Q_2,\ldots,Q_n)$
be a partition of the ground set of a connected matroid $M$.  
Then $\QQ$ is a {\em quasi-flower}
\index{quasi-flower} 
of $M$ with {\em petals} $Q_1,Q_2,\ldots,Q_m$
if  
\begin{itemize}
\item[(i)] $\lambda(Q_i\cup Q_{i+1}\udots Q_{i+k})\leq 2$
for all $i,k\in\{1,2,\ldots,n\}$ where subscripts are interpreted
modulo $n$; and 
\item[(ii)] if $(X,Y)$ is a 2-separation of $M$, then for some 
petal $Q_i$, either $X$ or $Y$ is a subset of $Q_i$.
\end{itemize}
The next lemma is clear.

\begin{lemma}
\label{quasi-flower}
Let $\QQ$ be a quasi-flower of the connected matroid
$M$ and let $\QQ'=(Q'_1,Q'_2,\ldots,Q'_m)$ be a concatenation
of $\QQ$. If $\lambda(Q'_i)\geq 2$ for all $i\in\{1,2,\ldots,m\}$,
then $\QQ'$ is a flower in $M$.
\end{lemma}

The concatenations of a quasi-flower that are flowers are
the flowers {\em displayed\ }
\index{displayed flower} 
by the quasi-flower. We also
say that a flower or quasiflower {\em displays} 
\index{displayed $3$-separating set}
a $3$-separating set $X$ or a 3-separation $(X,Y)$ if
$X$ is a union of petals. A $3$-separating set $X$
is {\em contained in a petal}
\index{contained in a petal} 
of a flower or quasiflower if
$X\subseteq \fcl(P)$ for some petal $P$.

Now return attention to flowers.
A {\em trivial\ } 
\index{trivial flower}
flower has just one petal. A flower with
two petals is nothing more than a partition $(P_1,P_2)$
of $E(M)$ for which $\lambda(P_1)=2$. 
If $n=3$, there is no distinction
between an anemone and a daisy. 
For $n\geq 3$, the anemone $(P_1,\ldots,P_n)$ is 
\begin{itemize}
\item[(i)] a {\em paddle}
\index{paddle} 
if $\sqcap(P_i,P_j)=2$ for all distinct
$i,j\in\{1,2,\ldots,n\}$;
\item[(ii)] a {\em copaddle}
\index{copaddle} 
if $\sqcap(P_i,P_j)=0$
for all distinct $i,j\in\{1,2,\ldots,n\}$; and
\item[(iii)] {\em spike-like}
\index{spike-like flower} 
if $n\geq 4$, and 
$\sqcap(P_i,P_j)=1$ for all distinct $i,j\in\{1,2,\ldots,n\}$.
\end{itemize}
A daisy $(P_1,\ldots,P_n)$ is 
\begin{itemize}
\item[(i)] {\em swirl-like}
\index{swirl-like flower} 
if $n\geq 4$ and $\sqcap(P_i,P_j)=1$
for all consecutive $i$ and $j$, and $\sqcap(P_i,P_j)=0$ for
all non-consecutive $i$ and $j$; and
\item[(ii)] {\em V\'amos-like} 
\index{V\'amos-like flower}
if $n=4$ and $\sqcap(P_i,P_j)=1$
for all consecutive $i$ and $j$, while 
$\{\sqcap(P_1,P_3),\sqcap(P_2,P_4)\}=\{0,1\}$.
\end{itemize}
A flower is {\em unresolved} 
\index{unresolved flower}
if $n=3$ and 
$\sqcap(P_i,P_j)=1$ for all distinct $i,j\in\{1,2,3\}$. Due to the
presence of possible additional structure, 
some unresolved flowers are
best regarded as spike-like and others as 
swirl-like.

\begin{theorem}
\label{flower1}
If $\FF$ is a flower, then $\FF$ is either
a paddle, a copaddle, spike-like, V\'amos-like, or swirl-like.
\end{theorem}

Note that $\FF$ is a flower in $M$ if and only if $\FF$
is a flower in $M^*$. Indeed if $\FF$ is spike-like in
$M$, then $\FF$ is spike-like in $M^*$, and if $M$ is swirl-like
in $M$, then $\FF$ is swirl-like in $M^*$. Moreover $\FF$
is a paddle in $M$ if and only if $\FF$ is a copaddle in 
$M^*$.

\section{Equivalent Flowers}

We have defined flowers for connected matroids and there
will be a number of occasions where we will consider flowers
in matroids that are not 3-connected. Having said this the
majority of the time we are interested only in flowers in
3-connected matroids. The more detailed structural descriptions
that we give in the remainder of this chapter apply
only to flowers in 3-connected matroid. It would certainly
be possible to generalise these notions to matroids that are
not 3-connected, but at the cost of additional technicalities. 

Let $M$ be a 3-connected matroid.
Let $\FF_1$ and $\FF_2$ be flowers of $M$, then 
$\FF_1\smalle \FF_2$ if every non-sequential
3-separation displayed by
$\FF_1$ is equivalent to one displayed by $\FF_2$. Clearly
$\smalle$ is a quasi-order on the collection of flowers
of $M$. The flowers $\FF_1$ and $\FF_2$ are {\em equivalent},
\index{equivalent flowers}
denoted $\FF_1\cong \FF_2$, if $\FF_1\smalle \FF_2$ and 
$\FF_2\smalle \FF_1$. Thus equivalent flowers display,
up to equivalence of 3-separations, exactly the same  
non-sequential 3-separations. The flower $\FF_1$ is {\em maximal\ }
\index{maximal flower}
if whenever $\FF_1\smalle \FF_2$, then $\FF_1\cong\FF_2$.
The {\em order}
\index{order of a flower} 
of a flower $\FF$ is the minimum 
number of petals in a 
flower equivalent to $\FF$. 

We give here some examples to illustrate some of the ideas developed
so far and also to motivate some of the future material.
Let $M$ be a rank-$n$ wheel or a whirl, with fan ordering
$(a_1,a_2,\ldots,a_{2n})$ of the elements of $M$. 
Then $(\{a_1,a_2\},\{a_3,a_4\},\ldots,\{a_{2n-1},a_{2n}\})$, 
is a swirl-like flower of $M$. 
However, wheels and 
whirls have no non-sequential 3-separations, so this flower
is equivalent to the trivial flower $(E(M))$. 
Thus there is a sense in
which flowers obtained from wheels and whirls are degenerate.

On the other hand, let 
$\{a_1,b_1\},\{a_2,b_2\},\{a_3,b_3\},\{a_4,b_4\}$
be the set of legs, in order, of a rank-4 swirl. Then 
$(\{a_1,b_1\},\{a_2,b_2\},\{a_3,b_3\},\{a_4,b_4\})$ is a 
swirl-like flower of 
order 4. To see this, note that if $i\in\{1,2,3,4\}$,
then the partition 
$(\{a_i,b_i,a_{i+1},b_{i+1}\},
\{a_{i+2},b_{i+2},a_{i+3},b_{i+3}\})$
is a non-sequential $3$-separation; indeed
both sides of the separation are quads. Certainly no flower with 
fewer petals displays all these $3$-separations. This example 
generalises to larger swirls and to spikes as well. 

Let $\FF$ be a flower of $M$. An element $e$ of $M$ is {\em loose}
\index{loose element}
in $\FF$ if $e\in\fcl(P_i)$ for some petal 
$P_i$. An element that is
not loose is {\em tight}.
\index{tight element} The petal $P_i$ is {\em loose\ }
\index{loose petal} 
if all 
elements of $P_i$ are loose, otherwise it is {\em tight}.
\index{tight petal} A flower
is {\em tight}
\index{tight flower} 
if all of its petals are tight. The next lemma
summarises a number of results from \cite{flower}.

\begin{lemma}
\label{flower2}
Let $\FF$ be a flower of the matroid $M$.
\begin{itemize}
\item[(i)] If $P$ is a tight petal of $\FF$, then $P$ has at least
two tight elements.
\item[(ii)] If $\FF'\cong \FF$, then $\FF$ and $\FF'$ have
the same set of tight elements.
\item[(iii)] If $\FF$ has order at least $3$, then the order of 
$\FF$ is the number of petals in any tight flower equivalent
to $\FF$.
\end{itemize}
\end{lemma}

The condition that $\FF$ has order at least $3$ in
Lemma~\ref{flower2} is important as degeneracies can occur
for flowers of low order.
Consider $U_{3,6}$. Note that
$U_{3,6}\cong \Delta_3\cong \Lambda_3$. Any partition
of the elements of $U_{3,6}$ into
$2$-element subsets gives a tight 3-petal flower.
This flower is equivalent to the trivial flower as it displays
no non-sequential 3-separations. Thus the flower has order 1.

For another example of a similar degenerate situation that we converge
to in case analyses, let
$(C,D)$ be a non-sequential $3$-separation 
of the $3$-connected matroid $M$ where $D$ is a quad.
Let $(D_1,D_2)$ be an arbitrary partition of $D$ into
2-element subsets. Then $(C,D)$ and $(C,D_1,D_2)$ are
equivalent flowers as 
both display just one non-sequential
3-separation, namely $(C,D)$. Evidently this equivalence 
class of flowers has order 2. But $(C,D_1,D_2)$ is a tight
flower, so that (iii) does not hold in this case either.

\section{Structure of Tight Flowers}

We begin by giving more detail from \cite{flower} about 
the structure of equivalence classes of
flowers. Our primary interest is in
swirl-like flowers, but at several places through the
paper we need to understand the structure of other types
of flowers. 
If $(P_1,\ldots,P_n)$ is an
anemone, then any permutation of the petals gives an
equivalent flower, while, if it is a daisy, then
any permutation of the petals that corresponds to a symmetry of
the regular $n$-gon gives an equivalent flower. If one flower
can be obtained from another by such a permutation then
the two flowers are alternative descriptions of the 
same underlying object and we will
say that the flowers are {\em equal up to labels}.
In fact we use the term ``up to labels'' 
\index{up to labels}
to cover a multitude
of sins in this paper and will always mean something like
``by an appropriate relabelling'' or ``by an appropriate reordering
of the indices''.

We now clarify the
status of 3-petal spike-like and swirl-like flowers. 
Recall that we called such flowers
``unresolved''. Let $\FF=(P_1,P_2,P_3)$ be an
unresolved flower. If $\FF$ has no loose elements,
then we can regard it as both spike-like and swirl-like.
Assume that there are loose elements.
If $P_1$ contains loose elements, then there is
an element $x\in P_1$ such that, up to labels, 
$x\in \clstar(P_2)$.
If $x\in \clstar(P_3)$, then we say that the flower 
is {\em spike-like} 
and if $x\notin\clstar(P_3)$, then we say that the 
flower is {\em swirl-like}.
It is shown in \cite{flower} that this definition is consistent
in that, if $\FF$ has loose elements, then $\FF$ is not 
both spike-like and swirl-like.

V\'amos-like flowers are easily described.
The next theorem is \cite[Theorem~6.1]{flower}.

\begin{theorem}
\label{vamoos-vamos}
Let $\FF$ be a V\'amos-like flower of the 
$3$-connected matroid $M$. Then
$\FF$ has no loose elements and any flower equivalent
to $\FF$ is equal to $\FF$ up to labels.
\end{theorem}

\subsection*{Anemone Structure} 
A set $S$ of a elements of a $3$-connected
matroid $M$ is a {\em segment\ }
\index{segment}
if either $|S|\leq 1$ or $|S|\geq 2$
and $M|S\cong U_{2,|S|}$. Equivalently $S$ is a segment if
every 3-element subset of $S$ is a triangle.
The set $S$ is a {\em cosegment\ }
\index{cosegment} 
of $M$ if $S$ is a 
segment of $M^*$. The following is \cite[Theorem~7.1]{flower}.

\begin{theorem}
\label{anemone-structure}
Let $M$ be a $3$-connected matroid and let $\FF$
be a tight flower of $M$ of order $n\geq 3$ that
is a paddle, a copaddle, or is spike-like.
Let $T$ and $L$ denote the sets of tight and loose elements
of $\FF$ respectively. For each petal $P_i$
of $\FF$, let $T_i=P_i\cap T$.
\begin{itemize}
\item[(i)] If $\FF$ is a paddle, then $L$ is a 
segment, $r(T_i)\geq 3$, and $L\subseteq \cl(T_i)$
for all $i\in\{1,2,\ldots,n\}$.
\item[(ii)] If $\FF$ is a copaddle, then $L$ is 
a cosegment, $r^*(T_i)\geq 3$, and $L\subseteq \cl^*(T_i)$
for all $i\in\{1,2,\ldots,n\}$.
\item[(iii)] If $\FF$ is spike-like, then $|L|\leq 2$.
If $L$ contains a single element, then that element is 
either in the closure of $T_i$ for each $i$,
or is in the coclosure of $T_i$ for each $i$. If
$|L|=2$, then one member of $L$ is contained in the closure
of each $T_i$, while the other is contained in the 
coclosure of each $T_i$.
\end{itemize}
Moreover, up to arbitrary permutations of the petals,
the tight flowers equivalent to $\FF$ are precisely the 
partitions of the form 
$(T_1\cup L_1,T_2\cup L_2,\ldots,T_n\cup L_n),$ where 
$(L_1,L_2,\ldots,L_n)$ is a partition of $L$.
\end{theorem}

A loose element of a spike-like flower that is in the 
closure of each petal is called the {\em tip}
\index{tip}
of the spike-like flower and a loose element that is in the
coclosure of each petal is called the {\em cotip}.
\index{cotip}

We illustrate Theorem~\ref{anemone-structure}
with some examples.
Figure~\ref{paddle-pic} illustrates the structure
corresponding to an equivalence class of paddles.
The set of loose elements is $\{l_1,l_2,l_3,l_4\}$.
One flower in the equivalence class is
$(P_1,P_2\cup \{l_3\},P_3\cup \{l_1,l_2,l_4\},P_4)$.

\begin{figure}
\begin{tikzpicture}[thick,line join=round]
	\coordinate (a) at (-2,0);
	\coordinate (b) at (2,0);
	\filldraw[pattern color=lines, pattern=north east lines] (b) .. controls (2,1.5) and (1,1.5) .. (a); 
	\filldraw[pattern color=lines, pattern=north west lines] (b) .. controls (2,-1.5) and (1,-1.5) .. (a); 
	\fill[fill=white] (a) .. controls (-2,-2) and (-1,-2) .. (b); 
	\fill[fill=white] (a) .. controls (-1,2) and (1,2) .. (b); 
	\filldraw[pattern color=lines, pattern=north east lines] (a) .. controls (-2,-2) and (-1,-2) .. (b); 
	\filldraw[pattern color=lines, pattern=north west lines] (a) .. controls (-1,2) and (1,2) .. (b); 
	\draw[dashed] (b) .. controls (2,1.5) and (1,1.5) .. (a); 
	\draw[dashed] (b) .. controls (2,-1.5) and (1,-1.5) .. (a); 
	\draw (a) -- (b);
	\coordinate[label=above:{\textcolor{labels}{$l_1$}}] (l1) at ($(a)!0.275!(b)$);
	\coordinate[label=above:{\textcolor{labels}{$l_2$}}] (l2) at ($(a)!0.425!(b)$);
	\coordinate[label=above:{\textcolor{labels}{$l_3$}}] (l3) at ($(a)!0.575!(b)$);
	\coordinate[label=above:{\textcolor{labels}{$l_4$}}] (l4) at ($(a)!0.725!(b)$);
	\foreach \thg in {l1,l2,l3,l4} \fill[labels] (\thg) circle (3pt);
	\node at (-1.5,1.5) {$P_1$};
	\node at (2.2,0.85) {$P_2$};
	\node at (2.2,-0.75) {$P_3$};
	\node at (-2.2,-1) {$P_4$};
\end{tikzpicture}
\caption{A Paddle}
\label{paddle-pic}
\end{figure}

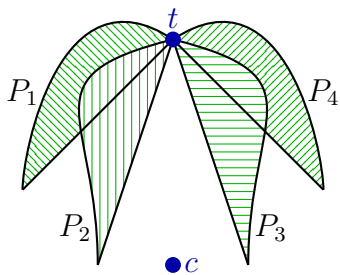
\begin{figure}
\begin{tikzpicture}[thick,line join=round]
	\coordinate (a) at (-2,1);
	\coordinate (b) at (-1,0);
	\coordinate[label=right:{\textcolor{labels}{$c$}}] (c) at (0,0);
	\coordinate (d) at (1,0);
	\coordinate (e) at (2,1);
	\coordinate[label=above:{\textcolor{labels}{$t$}}] (t) at (0,3);
	\filldraw[pattern color=lines,pattern=north west lines] (t) .. controls (-1.5,4) and (-2,1.5) .. (a); 
	\filldraw[pattern color=lines,pattern=north east lines] (t) .. controls (1.5,4) and (2,1.5) .. (e); 
	\filldraw[white] (t) .. controls (-2,2.5) and (-1,2) .. (b); 
	\filldraw[white] (t) .. controls (2,2.5) and (1,2) .. (d); 
	\filldraw[pattern color=lines,pattern=vertical lines] (t) .. controls (-2,2.5) and (-1,2) .. (b); 
	\filldraw[pattern color=lines,pattern=horizontal lines] (t) .. controls (2,2.5) and (1,2) .. (d); 
	\draw (a) -- (t) -- (e);
	\draw (b) -- (t) -- (d);
	\foreach \pt in {c,t} \fill[labels] (\pt) circle (3pt);
	\node at (-2,2.3) {$P_1$};
	\node at (-1.3,0.5) {$P_2$};
	\node at (2,2.3) {$P_4$};
	\node at (1.3,0.5) {$P_3$};
\end{tikzpicture}
\caption{Spike-like Structure}
\label{spike-like-pic}
\end{figure}

Figure~\ref{spike-like-pic} 
illustrates a spike-like flower with tip $t$
and cotip $c$. Note that if $c$ is deleted, then the petals
$P_1$, $P_2$, $P_3$ and $P_4$ are 2-separating. 
If $c$ is contracted we have a spike-like flower without a cotip.

\subsection*{Swirl-like Flower Structure}
The structure of swirl-like flowers is of 
particular interest to us. For a petal $P$ of a flower 
let $\cP$ denote the set of tight elements
of $P$ and let $\hP=\fcl(P)$.
For petals $P_i$ and $P_{i+1}$ of a swirl-like flower,
let $P_i^+=P_{i+1}^-=\hP_i\cap\hP_{i+1}$. 

\begin{theorem}
\label{swirl-like}
Let $M$ be a $3$-connected matroid and let 
$\FF=(P_1,\ldots,P_n)$ be a tight swirl-like flower of $M$
of order at least $3$. Then
$P_1^+\udots P_n^+$ is the set of loose elements of $\FF$.
Moreover, there is a partition
$$\PP=(\cP_1,p_1^i,\ldots,p_1^{k_1},\cP_2,
p_2^1,\ldots,p_2^{k_2},\cP_3,\ldots,\cP_n,p_n^1,\ldots,p_n^{k_n})$$ 
of $E(M)$ having the following properties.
\begin{itemize}
\item[(i)] $\PP$ is a quasi-flower.
\item[(ii)] $(p_i^1,p_i^2,\ldots,p_i^{k_i})$ is an
ordered fan in $M$ and $\{p_i^1,p_i^2,\ldots,p_i^{k_i}\}=P_i^+$
for all $i\in\{1,2,\ldots,n\}$. 
\item[(iii)] The partition $(P'_1,\ldots,P'_n)$ of $E(M)$ is a 
tight flower equivalent to
$\FF$ if and only if $(P'_1,P'_2,\ldots,P'_n)$ is a 
concatenation of $\PP$ such that $\cP_i\subseteq P'_i$
for all $i\in\{1,2,\ldots,n\}$.
\end{itemize}
\end{theorem}

\begin{figure}
\begin{tikzpicture}[thick,line join=round]
	\coordinate (e) at (-22.5:2);
	\coordinate (x) at (22.5:2);
	\coordinate[label=45:{\textcolor{labels}{$a$}}] (a) at (22.5+45:2);
	\coordinate (y) at (90+22.5:2);
	\coordinate (z) at (180-22.5:2);
	\coordinate (w) at (180+22.5:2);
	\coordinate (g) at (270-22.5:2);
	\coordinate (t) at (270+22.5:2);
	\node[pattern color=lines,draw,circle through=(x),pattern=horizontal lines] at ($(x)!0.5!(e)$) {};
	\node[pattern color=lines,draw,circle through=(e),pattern=north west lines] at ($(t)!0.5!(e)$) {};
	\node[pattern color=lines,draw,circle through=(g),pattern=north east lines] at ($(g)!0.5!(w)$) {};
	\node[pattern color=lines,draw,circle through=(z),pattern=north west lines] at ($(z)!0.5!(y)$) {};
	\node[pattern color=lines,draw,circle through=(y),pattern=vertical lines] at ($(y)!0.5!(a)$) {};
	\filldraw[fill=white] (x) -- (e) -- (t) -- (g) -- (w) -- (z) -- (y) -- (a);
	\draw[dashed] (a) -- (x);
	\coordinate[label=45:{\textcolor{labels}{$c$}}] (c) at ($(a) + (1.2,-0.3)$);
	\coordinate[label=45:{\textcolor{labels}{$b$}}] (b) at ($(a)!0.5!(c)$);
	\coordinate[label=right:{\textcolor{labels}{$d$}}] (d) at ($(c)!0.5!(x)$);
	\coordinate[label=135:{\textcolor{labels}{$e$}}] (e') at (e);
	\coordinate[label=45:{\textcolor{labels}{$g$}}] (g') at (g);
	\coordinate[label=above:{\textcolor{labels}{$f$}}] (f) at ($(g)!0.5!(t)$);
	\coordinate[label=right:{\textcolor{labels}{$h$}}] (h) at ($(z)!0.5!(w)$);
	\draw (a) -- (c) -- (x);
	\foreach \pt in {a,b,c,d,e,g,f,h} \fill[black] (\pt) circle (3pt);
	\node at ($(a)!0.5!(y) + (0,1.1)$) {\textcolor{labels}{$\cP_1$}};
	\node at ($(e)!0.5!(x) + (1.1,0)$) {\textcolor{labels}{$\cP_2$}};
	\node at ($(e)!0.5!(t) + (0.8,-0.8)$) {\textcolor{labels}{$\cP_3$}};
	\node at ($(g)!0.5!(w) + (-0.8,-0.8)$) {\textcolor{labels}{$\cP_4$}};
	\node at ($(z)!0.5!(y) + (-0.8,0.8)$) {\textcolor{labels}{$\cP_5$}};
\end{tikzpicture}
\caption{Swirl-like Structure}
\label{swirl-like-pic}
\end{figure}
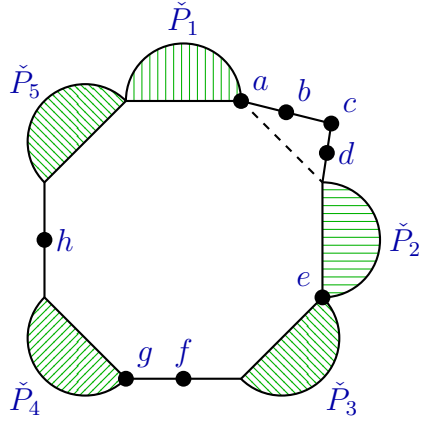

Figure~\ref{swirl-like-pic} illustrates Theorem~\ref{swirl-like}.
The quasi-flower corresponding to the diagram is
$$(\cP_1,a,b,c,d,\cP_2,e,\cP_3,f,g,\cP_4,h,\cP_5).$$
A flower in the equivalence class corresponding to the 
above quasiflower is 
$$(\cP_1\cup\{a\},\{b,c,d\}\cup \cP_2\cup\{e\},
\cP_3,\{f,g\}\cup \cP_4,\{h\}\cup \cP_5).$$

\subsection*{Blooms}
In this subsection we develop further terminology for 
swirl-like flowers
to facilitate
arguments in proofs. We will call the partition of $E(M)$ given by
Theorem~\ref{swirl-like} a {\em bloom}
\index{bloom} 
of $M$.
Note that, associated with a bloom, is a particular ordering
of the fans of loose elements between petals. 
From now on, when we refer
to such a set, we will typically regard it as being endowed
with the ordering induced from a bloom.
To avoid clumsy notation we
will often denote a bloom simply by $(\hP_1,\ldots,\hP_n)$ noting
that it follows from Theorem~\ref{swirl-like} that, if $n\geq 3$,
then the full bloom can be recovered using 
structure in the underlying
matroid. 

For a bloom $(\hP_1,\ldots,\hP_n)$ we use the following
facts and notational conventions freely throughout this paper:
$P_i^+=\hP_i\cap \hP_{i+1}$, $P_i^-=P_{i-1}^+$, and 
$\cP_i=\hP_i-(P_i^-\cup P_i^+)$. 

Recall that blooms are quasi-flowers and that the concatenations
of a quasi-flower that are flowers are the flowers 
{\em displayed\ } by that bloom.
\index{displayed flower}
In particular a $3$-separation is {\em displayed\ }
\index{displayed $3$-separation} 
by a bloom if it is
a concatenation of that bloom. A $3$-separating set
is {\em contained in a petal\ } 
\index{contained in a petal}
of a bloom if 
$A\subseteq \hP_i$ for some
$i$. Of course two blooms are {\em equivalent\ }
\index{equivalent blooms} 
if they display the
same flowers and it is easily seen that equivalent blooms are 
equal up to permutations of the $n$-gon. 

The partial order on swirl-like flowers extends 
easily to include blooms. 
Thus, if 
$\FF$ is a swirl-like flower and $\BB$ is bloom of $M$, 
then $\FF\smalle\BB$ if there is a flower $\FF'$ displayed
by $\BB$ such that $\FF\smalle\FF'$. 

We could easily have broadened the notion of blooms to cover all
types of flowers, but we will have no need for this more general notion
in this paper. The next lemma is easy.

\begin{lemma}
\label{bloom-concat}
If $\BB$ is a maximal bloom, and $\FF\smalle\BB$ is a flower, then
$\FF$ is a concatenation of a flower displayed by $\BB$.
\end{lemma}

\subsection*{Loose elements in Swirl-like Flowers}
The next lemma says that, most of the time, elements in the
guts or coguts of 3-separations displayed by a swirl-like
flower are loose elements. 

\begin{lemma}
\label{fine1}
Let $\FF=(P_1,P_2,\ldots,P_n)$ be a tight swirl-like flower 
of the $3$-connected matroid $M$ and let $i$ be an integer
in $\{1,2,\ldots,n-2\}$. Then the following hold.
\begin{itemize}
\item[(i)] If $x\notin P_1\udots P_i$, then
$x\in\clstar(P_1\udots P_i)$ if and only if $x\in\clstar(P_1)$ or 
$x\in\clstar(P_i)$.
\item[(ii)] $\fcl(P_1\udots P_i)=(P_1\udots P_i)\cup(\hP_1\cup\hP_i)=
\hP_1\cup \hP_2\udots \hP_i$.
\item[(iii)] $\hP_1\udots \hP_i$ is fully closed.
\end{itemize}
\end{lemma}

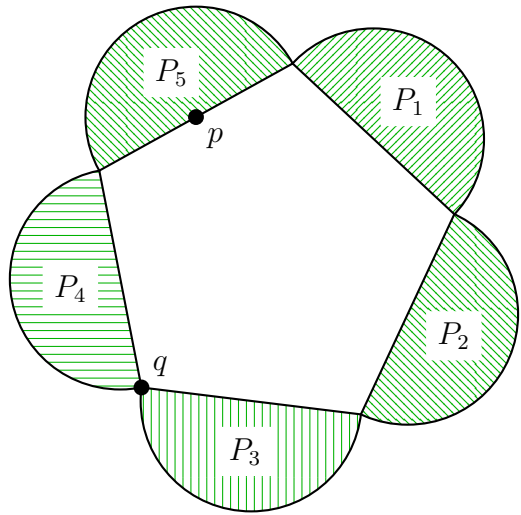
\begin{figure}
\begin{tikzpicture}[thick,line join=round]
	\coordinate (A) at (11:2.5);
	\coordinate (B) at (83:2.5);
	\coordinate (C) at (155:2.5);
	\coordinate (D) at (227:2.5);
	\coordinate (E) at (299:2.5);
	\coordinate (ab) at ($(A)!0.5!(B)$);
	\coordinate (bc) at ($(B)!0.5!(C)$);
	\coordinate (cd) at ($(C)!0.5!(D)$);
	\coordinate (de) at ($(D)!0.5!(E)$);
	\coordinate (ea) at ($(E)!0.5!(A)$);
	\node[pattern color=lines,draw,circle through=(A),pattern=north east lines] (AB) at (ab) {}; 
	\node[pattern color=lines,draw,circle through=(B),pattern=north west lines] (BC) at (bc) {}; 
	\node[pattern color=lines,draw,circle through=(C),pattern=horizontal lines] (CD) at (cd) {}; 
	\node[pattern color=lines,draw,circle through=(D),pattern=vertical lines] (DE) at (de) {}; 
	\node[pattern color=lines,draw,circle through=(E),pattern=north west lines] (EA) at (ea) {}; 
	\filldraw[fill=white] (A) -- (B) -- (C) -- (D) -- (E) -- cycle;
	\coordinate[label=-45:$p$] (p) at ($(B)!0.5!(C)$);
	\coordinate[label=45:$q$] (q) at (227:2.5);
	\foreach \pt in {p,q} \fill[black] (\pt) circle (3pt);
	\node at (47:2.7) [rectangle,fill=white,draw=white] {$P_1$};
	\node at (119:2.7) [rectangle,fill=white,draw=white] {$P_5$};
	\node at (191:2.7) [rectangle,fill=white,draw=white] {$P_4$};
	\node at (263:2.7) [rectangle,fill=white,draw=white] {$P_3$};
	\node at (335:2.7) [rectangle,fill=white,draw=white] {$P_2$};\end{tikzpicture}
\caption{Illustration of Lemma~\ref{fine1}}
\label{guts-not-loose}
\end{figure}

Note that the condition that
$i\leq n-2$ in Lemma~\ref{fine1} is necessary. 
To see this consider the flower illustrated
in Figure~\ref{guts-not-loose}. The element $p$ is in the
closure of $P_1\cup P_2\cup P_3\cup P_4$, but is not a loose
element of the flower. On the other hand the element $q$ is in
the closure of $P_1\cup P_2\cup P_3$ and is guaranteed to be
a loose element.

It follows from Lemma~\ref{fine1} that, 
if $(P_1,\ldots,P_n)$ is a tight swirl-like flower
of $M$, then 
$\hP_1\cup P_2\udots P_{i-1}\cup \hP_i=
\hP_1\cup\hP_2\udots\hP_{i-1}\cup\hP_i$
for any $i\in\{1,2,\ldots,n\}$,
so no harm is done by blurring this distinction.

Let $\FF=(P_1,\ldots, P_n)$ be a swirl-like flower in the 
3-connected matroid $M$. Consider the ordered set
$P_i^+=(p_{1},\ldots,p_{k})$. By Theorem~\ref{swirl-like},
this ordered set is a fan. We say that it is 
the fan of loose elements 
{\em between $\cP_i$ and $\cP_{i+1}$}. We say that the element
$p_i$ is a {\em spoke} element if 
$p_i\in\cl(\cP_i\cup\{p_1,\ldots,p_{i-1}\})$,
and $p_i$ is a {\em rim} element if  
$p_i\in\cl^*(\cP_i\cup\{p_1,\ldots,p_{i-1}\})$. It is easily seen that 
elements of $P_i^+$ alternate between rim and spoke elements. 
Note that, while a fan with less that five elements does not
have a canonical ordering, the structure induced by the bloom
does induce a canonical ordering on the elements of 
fans of loose elements between petals of a swirl-like flower
regardless of the number of elements they have.

\begin{lemma}
\label{fine2}
Let $\FF=(P_1,\ldots,P_n)$ be a swirl-like flower of 
the $3$-connected matroid $M$ and say that
$1\leq i\leq n-2$. Then $x\in\cl(P_1\cup P_2\udots P_i)$ if and only if either
\begin{itemize}
\item[(i)] $x$ is the last element of $P_1^-$ not in $P_1$, and
$x$ is a spoke element of $P_1^-$, or 
\item[(ii)] $x$ is the first element of $P_i^+$ not in $P_1$, and 
$x$ is a spoke element of $P_i^+$.
\end{itemize}
\end{lemma}

Now
consider adjacent petals of a swirl-like flower of order at least
$3$. Label these petals $P_1$ and $P_2$. Assume that $x$ is the first
element of $P_1^+$ not in $P_1$.  Then $x\in\clstar(P_1)$.
Thus $x$ is either in the guts or coguts of the 3-separating set
$P_1$. By Lemma~\ref{fine2}, if $x$ is in the guts of $P_1$,
then $x$ is a spoke element of $P_1^+$ and by the dual of
that lemma, if $x$ is in the coguts of $P_1$, then $x$ is a rim
element of $P_1^+$. For this reason we will frequently refer to the
rim or spoke element $x$  as being in the 
{\em guts} or  {\em coguts} of $(P_1,P_2)$. 
\index{guts of $(P_1,P_2)$}
\index{coguts of $(P_1,P_2)$}
We may also refer to $x$ as being a loose guts or coguts
element {\em between} $P_1$ and $P_2$.
The next lemma is essentially a rephrasing of Lemma~\ref{fine2}.

\begin{lemma}
\label{fine3}
Let $\FF=(P_1,\ldots,P_n)$ be a swirl-like flower of 
the $3$-connected matroid $M$. If $1\leq i\leq n-2$, 
and $x\notin(P_1\udots P_i)$, then $x$ is in
the guts of the $3$-separation $(P_1\udots P_i,P_{i+1}\udots P_n)$
if and only if either 
\begin{itemize}
\item[(i)] $x$ is in the guts 
of $(P_n,P_1)$, and $x$ is the first element of $P_1^-$ not in 
$P_1$, or
\item[(ii)] $x$ is in the
guts of $(P_i,P_{i+1})$, and $x$ is the first element of 
$P_i^+$ not in $P_i$.
\end{itemize}
\end{lemma}

The next lemma describes one situation in which a swirl-like
flower is induced in an extension of a matroid.

\begin{lemma}
\label{page3}
Let $M$ and $M\ba x$ be $3$-connected matroids and 
let $\FF=(P_1,P_2,\ldots,P_n)$
be a swirl-like flower of $M\ba x$. Say $n\geq j>i+1\geq 1$, 
and $x$ is in both $\cl(P_1\cup\cdots\cup P_i)$
and $\cl(P_j\cup\cdots\cup P_n)$. Then 
$(P_1\cup\{x\},P_2,\ldots,P_n)$
is a flower of $M$ and, in this flower, 
$x$ is in the guts of $(P_1\cup \{x\},P_2)$.
\end{lemma}

Ends of fans of loose elements between petals
are good choices for removal without losing 3-connectivity,
or indeed, as we shall see later, other types of connectivity.

\begin{lemma}
\label{red-hill}
Let $\PP$ 
be a swirl-like flower of order at least three in the
$3$-connected matroid $M$. Let $(f_1,f_2,\ldots,f_n)$ 
be a maimal ordered fan of loose elements between consecutive petals
of $\PP$ with ordering induced by a bloom associated with $\PP$.
\begin{itemize}
\item[(i)] If $f_i$ is a spoke element, then $M\ba f_i$ is 
$3$-connected if $i\in\{1,n\}$ and $M\ba f_i$
is $3$-connected up to the single series pair $\{f_{i-1},f_{i+1}\}$
otherwise.
\item[(ii)] If $f_i$ is a rim element, then $M/f_i$ is $3$-connected
if $i\in\{1,n\}$ and $M/f_i$ is $3$-connected
up to the single parallel pair $\{f_{i-1},f_{i+1}\}$
otherwise.
\end{itemize}
\end{lemma} 

\begin{figure}
\begin{tikzpicture}[thick,line join=round]
	\coordinate (a) at (0-18:2);
	\coordinate (b) at (72-18:2);
	\coordinate (c) at (144-18:2);
	\coordinate (d) at (-144-18:2);
	\coordinate (e) at (-72-18:2);
	\node[pattern color=darkgray,draw,circle through=(d),pattern=horizontal lines] at ($(c)!0.5!(d)$) {};
	\node[pattern color=darkgray,draw,circle through=(e),pattern=north east lines] at ($(d)!0.5!(e)$) {};
	\node[pattern color=darkgray,draw,circle through=(a),pattern=north west lines] at ($(e)!0.5!(a)$) {};
	\node[pattern color=darkgray,draw,circle through=(b),pattern=horizontal lines] at ($(a)!0.5!(b)$) {};
	\node[draw,circle through=(c)] at ($(b)!0.5!(c)$) {};	
	\filldraw[fill=white] (a) -- (b) -- (c) -- (d) -- (e) -- cycle;
	\coordinate (x) at ($(b)!0.6!(c)$);
	\coordinate (y) at ($(b)!0.2!(c)$);
	\coordinate (o) at ($(b)!0.9!(c)$);
	\coordinate (z) at ($(b)!0.5!(c) + (0.1,0.9)$);
	\coordinate (t) at ($(z)!0.4!(y)$);
	\coordinate (w) at (intersection of o--t and x--z);
	\draw (x) -- (z) -- (y);
	\draw (o) -- (t);
	\foreach \pt in {x,y,z,t,w} \fill[black] (\pt) circle (3pt);
	\node at ($(a)!0.5!(b) + (0.5,0.15)$) [rectangle,fill=white,draw=white,inner sep=1pt] {$P_2$};
	\node at ($(a)!0.5!(e) + (0.35,-0.35)$) [rectangle,fill=white,draw=white,inner sep=1pt] {$P_3$};
	\node at ($(d)!0.5!(e) + (-0.35,-0.35)$) [rectangle,fill=white,draw=white,inner sep=1pt] {$P_4$};
	\node at ($(d)!0.5!(c) + (-0.5,0.15)$) [rectangle,fill=white,draw=white,inner sep=1pt] {$P_5$};
	\node at ($(b) + (-0.4,1.3)$) {$P_1$};
\end{tikzpicture}
\caption{Not all Fans are Loose}
\label{fan-not-loose}
\end{figure}
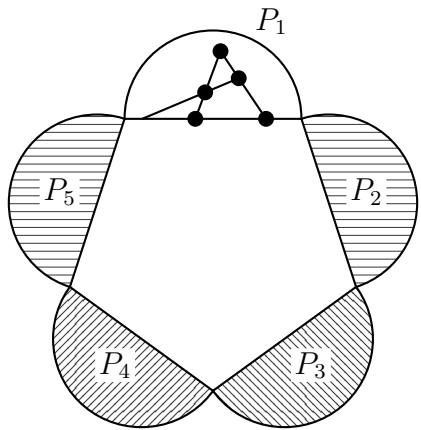

Given the fact that the loose elements of swirl-like flowers form
fans, the reader may be tempted to think that the converse holds
in that a petal whose elements form a fan 
consists entirely of loose elements, that is, 
is a loose petal. This is not 
the case. For example the legs of a swirl are trivially fans
and their elements are not loose. Another example is 
illustrated in Figure~\ref{fan-not-loose}. It all depends on how the 
elements of the fan align with the rest of the matroid.

\section{A Grab Bag of Flower Lemmas}

This section consists of an unordered collection of lemmas on
flowers. The criteria used for inclusion in this section is that
the lemma is needed somewhere in this paper
and that it didn't have a natural
home elsewhere in this chapter.

A $3$-separation $(R,B)$ is {\em well displayed\ }
\index{well-displayed $3$-separation} 
by a bloom $\FF$ if it is
displayed by $\FF$ and neither $R$ nor $B$ is contained in a petal
of $\FF$. The next lemma is important and used many 
times, often without reference.
 
\begin{lemma}
\label{well-displayed}
If $\FF$ is a bloom of $M$ and the $3$-separation $(R,B)$
is well displayed by $\FF$, then $(R,B)$ is non-sequential.
\end{lemma}

\begin{proof}
Say that $(R,B)$ is well-displayed. Then for some
tight 4-petal flower
$(P_1,P_2,P_3,P_4)$ displayed by $\FF$, we have 
$(R,B)=(P_1\cup P_2,P_3\cup P_4)$.
It follows from Lemma~\ref{fine1}, that $\fcl(P_1\cup P_2)\neq E(M)$
and $\fcl(P_3\cup P_4)\neq E(M)$, so that $(R,B)$ is non-sequential.
\end{proof}

Recall that disjoint sets $A$ and $B$ of a matroid $M$ are
{\em fully skew} if $\sqcap(A,B)=\sqcap^*(A,B)=0$.
The next lemma follows from the definition of a 
swirl-like flower and the fact that duals of swirl-like flowers
are swirl-like flowers.

\begin{lemma}
\label{skew-petals}
If $P_i$ and $P_j$ are non-consecutive petals of a 
(not necessarily tight) swirl-like flower,
then $P_i$ and $P_j$ are fully skew.
\end{lemma}

The next lemma follows from results in \cite{flower}.

\begin{lemma}
\label{in-petal}
Let $\PP$ and $\QQ$ be inequivalent maximal flowers of
a $3$-connected matroid $M$ of order at least 
three. Then there is a petal $P$
of $\PP$ and a petal $Q$ of $\QQ$ such that $E(M)-P$ 
is contained in $\fcl(Q)$.
\end{lemma}

As an easy consequence of Lemma~\ref{in-petal} we have

\begin{corollary}
\label{in-petal2}
Let $\PP$ and $\QQ$ be inequivalent maximal swirl-like flowers
of a $3$-connected matroid,
then there is a petal $P\in\PP$ and a petal $Q\in\QQ$ such that
$\cP\subseteq Q$.
\end{corollary}

One might expect that
V\'amos-like flowers are potentially problematic, 
particularly as we 
make no assumptions about representability until the final
chapter. However, it turns out that they 
cause no difficulties for us. One reason for this is that
a V\'amos-like flower is never a concatenation of 
a flower of order greater than four.
Another reason is captured by the following lemma.

\begin{lemma}
\label{keep-spiral}
Let $M$ and $M\ba x$ be $3$-connected matroids,
let $\PP=(P_1,\ldots,P_n)$ be a flower of $M\ba x$.
Assume that, for some $i\in\{1,2,\ldots,n\}$,
we have $x\in\cl(P_i\udots P_n)$ and let 
$\PP'=(P_1,\ldots,P_{i-1},P_i\udots P_n\cup\{x\})$, where $i\geq 3$. 
The $\PP'$ is a flower in $M$ of the same type as $\PP$.
\end{lemma}

Another easy lemma with a similar flavour to Lemma~\ref{keep-spiral}
is

\begin{lemma}
\label{keep-flower-type}
Let $\PP=(P_1,P_2,\ldots,P_n)$ be a flower in the connected
matroid $M$ and let $M'$ be a $3$-connected minor of $M$ on $E'$
with the property that $|P_i\cap E'|\geq 2$ for all
$i\in\{1,2,\ldots,n\}$. Then
$(P_1\cap E',P_2\cap E',\ldots,P_n\cap E')$ is a flower in
$M'$ of the same type as $\PP$.
\end{lemma}

The next lemma is  proved in \cite{flower}.

\begin{lemma}
\label{flower8.2}
Let $\FF$ be a flower of the $3$-connected
matroid $M$, and let $(R,B)$ be a 
$3$-separation of $M$ such that;
\begin{itemize}
\item[(i)] neither $R$ nor $B$ is contained in a petal of $\FF$; and
\item[(ii)] if $(R,B)$ crosses a petal $P$, then 
$|P\cap R|,|P\cap B|\geq 2$.
\end{itemize}
Then there is a flower that refines $\FF$ and displays $(R,B)$.
\end{lemma}

Spikes and swirls show that it is 
possible for a petal of a tight flower to have
two elements, but this does not happen for all flowers.

\begin{lemma}
\label{flower4}
Let $P$ be a petal of a tight flower $\FF$ of a 
$3$-connected matroid. If $|P|=2$, then
$\FF$ is either swirl-like, spike-like or V\'amos-like.
\end{lemma}

We frequently need to test that we have a flower of a certain
type. One way to do this is by local connectivity between
petals. Another way is by the behaviour of loose elements.
These tests are obvious consequences of the structure
theorems for flowers. We note one here.

\begin{lemma}
\label{get-spike}
Let $\FF=(P_1,P_2,\ldots,P_n)$ be a flower in the 
$3$-connected matroid $M$ where $n\geq 3$.
Assume that both $\cl(P_1)\cap P_2$ and 
$\cl^*(P_1)\cap P_2$ are nonempty. Then $\FF$ is 
spike-like.
\end{lemma}

The next fact on loose elements in swirl-like flowers is used 
many times.

\begin{lemma}
\label{loose-swirl-like}
Let $\FF=(P_1,P_2,\ldots,P_n)$ be a swirl-like
flower in the 
$3$-connected matroid $M$, where $n\geq 3$. 
Then $|\cl(P_1)\cap P_2|\leq 1$, and if 
$\cl(P_1)\cap P_2\neq\emptyset$, then 
$\cl^*(P_1)\cap P_2=\emptyset$.
\end{lemma}

Neither of the next two facts is surprising.

\begin{lemma}
\label{display-small}
Let $\FF$ be a flower. Then there is a maximal flower $\FF'$
such that $\FF\less \FF'$ and $\FF'$ displays $\FF$.
\end{lemma}

\begin{lemma}
\label{one-flower}
Let $X$ be a $3$-separating set of elements of the
$3$-connected matroid $M$ and let
$\PP$ and $\QQ$ be  maximal flowers of $M$ that display $X$. Then
$\PP\cong \QQ$.
\end{lemma}

Finally we note that if a loose coguts element of a swirl-like flower
is in a triangle, then that triangle is part of a fan of loose elements.

\begin{lemma}
\label{loosey-loose}
Let $\PP=(P_1,\{b\},P_2,\ldots,P_n)$ be a swirl-like quasi-flower
of the $3$-connected matroid $M$ of order at least four, where
$b$ is a loose element of $\PP$ in the coguts of
$P_1$. If $\{a,b,c\}$ is a
triangle of $M$, then, up to labels, $a\in P_1$, $c\in P_2$,
and both $a$ and $c$ are loose elements in the guts of $P_1$
and $P_2$ respectively.
\end{lemma}

\section{Flowers and Modularity}

Let $A$ and $B$ be sets of elements in a matroid $M$. Recall
that $A$ and $B$ form a {\em modular pair}
\index{modular pair} 
if $r(A)+r(B)=r(A\cup B)+r(A\cap B)$. The next lemma is 
elementary but fundamental.

\begin{lemma}
\label{elementary-modular}
Let $e$ be an element of the matroid $M$ and let $A$ and 
$B$ be a modular pair of sets of $M$. If
$e\in\cl(A)$ and $e\in\cl(B)$, then $e\in\cl(A\cap B)$.
\end{lemma}

The next lemma provides a useful connection between
the connectivity function and modularity.

\begin{lemma}
\label{lambda-modular}
Let $A$ and $B$ be sets of elements of the matroid $M$.
If $\lambda(A)+\lambda(B)=\lambda(A\cup B)+\lambda(A\cap B)$,
then $A$ and $B$ are a modular pair in $M$.
\end{lemma}

\begin{proof}
Let $A'=E(M)-A$ and $B'=E(M)-B$. Since
$\lambda(A)+\lambda(B)=\lambda(A\cup B)+\lambda(A\cap B)$,
we have
$$r(A)+r(A')+r(B)+r(B')=r(A\cup B)+r(A'\cap B')+r(A\cap B)
+r(A'\cup B'),$$
so that
$$r(A)+r(B)-r(A\cup B)-r(A\cap B)=r(A'\cup B')+r(A'\cap B')
-r(A')-r(B').$$
The lemma now follows from the submodularity of the rank
function.
\end{proof}

The next lemma follows from Lemma~\ref{lambda-modular} 
and the definition of flower.

\begin{lemma}
\label{lots-of-modular}
Let $\PP$ be a flower in 
the connected  matroid $M$. 
Let $A$ and $B$ be sets
displayed by $\PP$ such that $A\cap B$ contains at least
one petal of $\PP$ and $A\cup B$ avoids at least
one petal of $\PP$. Then $A$ and $B$ form a modular pair.
\end{lemma}

As an immediate corollary of Lemma~\ref{elementary-modular}
and Lemma~\ref{lots-of-modular}, we have

\begin{corollary}
\label{modular0}
Let $x$ be an element of the matroid $M$ such that
$M\ba x$ is connected with a flower $\PP$. Let 
$A$ and $B$ be sets
displayed by $\PP$ such that $A\cap B$ contains at least
one petal of $\PP$ and $A\cup B$ avoids at least
one petal of $\PP$. If $x\in\cl(A)$ and $x\in\cl(B)$,
then $x\in\cl(A\cap B)$.
\end{corollary}

The next lemma is a strengthening of Corollary~\ref{modular0}
for blooms. We apply it in many places.

\begin{lemma}
\label{modular}
Let $M$ be a matroid with an element $x$ such that
$M\ba x$ is $3$-connected 
with a bloom $\FF$ of order at least $3$. Assume that
the $3$-separating sets $A$ and $B$ are displayed by $\FF$ and that 
$|A\cap B|\geq 2$. Assume further that $x\in\cl(A)$ and
$x\in\cl(B)$. 
\begin{itemize}
\item[(i)] If $|E(M\ba x)-(A\cup B)|\geq 2$, then
$x\in\cl(A\cap B)$.
\item[(ii)] If $|E(M\ba x)-(A\cup B)|=1$ and
this set consists of a loose coguts element of $\FF$,
then $x\in\cl(A\cap B)$.
\end{itemize}
\end{lemma}

\begin{proof}
Let $C=E(M\ba x)-(A\cup B)$. 
Assume that (i) holds. Then $(A-B,A\cap B,B-A,C)$ is a flower in
$M\ba x$ and the result holds by Lemma~\ref{modular0}.

Assume that (ii) holds. Say $C=\{c\}$. As
$(A-B,A\cap B,(B-A)\cup\{c\})$ is a flower in $M\ba x$
and $c$ is a loose coguts element, we see that $c$ is in the coguts of the 3-separation $(A-B,B\cup\{c\})$ of $M\ba x$. Thus
$\sqcap(A-B,B)=\sqcap(A-B,B\cup\{c\})-1=1$. Hence
$$r(A\cup B)=r(B)+r(A-B)+1.$$
Also, as $(A-B,A\cap B,(B-A)\cup\{c\})$ is a 
swirl-like flower
in $M\ba x$, we have $\sqcap(A-B,A\cap B)=1$. Hence
$$r(A)=r(A\cap B)+r(A-B)+1.$$
It follows from the two displayed equations that
$A$ and $B$ form a modular pair in $M\ba x$. Now 
$x\in\cl_M(A\cap B)$ by Lemma~\ref{elementary-modular}.
\end{proof}

\section{Separations Crossing Blooms}

In this section we analyse the precise ways 
in which a 3- or 4-separation
can interact with a maximal bloom. 
The results of this section are lemmas
for results in later sections.

Let $\FF=(\hP_1,\ldots,\hP_n)$ be a maximal bloom of the matroid $M$. 
Then a set $S$ is {\em consecutive}
\index{consecutive set} 
in $\FF$ 
if it is a union of consecutive petals in the
quasi-flower $\FF$. In other words, $S$ is 
consecutive if, for some $i,j\in\{1,\ldots,n\}$,  
$S=L_i^-\cup\cP_i\cup \hP_{i+1}\udots\hP_{j-1}\cup\cP_j\cup L_j^+$, where 
$L_i^-$ is a terminal section of $P_i^-$, 
and $L_j^+$ is an initial section
of $P_j^+$. Note that a $3$-separation $(A,B)$ is
displayed by $\FF$ if and only if $A$ is consecutive. 
Recall that when we say that a set $S$ is contained in
a petal of $\FF$, we mean that $S\subseteq \hP_i$ for some $i$. 

\begin{lemma}
\label{cross1}
Let $M$ be a $3$-connected matroid, 
let $\FF$ be a maximal bloom of order at least $4$, and 
let $(R,B)$ be a $3$-separation
of $M$. If $(R,B)$ is not displayed by $\FF$, 
then either $R$ or $B$ is
contained in a petal of $\FF$.
\end{lemma}

\begin{proof}
If $(R,B)$ is non-sequential, then the result follows from 
\cite[Theorem~8.1]{flower}. Thus we may assume that either $R$ or 
$B$ is sequential. Say that $R$ is sequential. Then
there is a triangle or triad $T=\{a,b,c\}$ of $R$ such that 
$R\subseteq \fcl(T)$.
Up to labels, there is an $i\in\{1,2,\ldots,n-2\}$ such that
$a\in\hP_1$, $b\in\hP_i$, and $c\in\hP_j$ for some $j\geq i$. As 
$i\leq n-2$, by Lemma~\ref{fine1}, $\hP_1\udots\hP_i$  is fully closed
and either $c\in\hP_1$ or $c\in\hP_i$. 
Up to labels we may assume the former holds.
Then $\{a,c\}\subseteq \hP_1$ and $\hP_1$ is fully closed, so 
$T\subseteq \hP_1$ and indeed, $\fcl(T)\subseteq \hP_1$. 
Thus $R$ is contained
in a petal.
\end{proof}

The next lemma gives a more precise description of the the 
way that a $3$-separation can cross a bloom.

\begin{lemma}
\label{3-cross}
Let $M$ be a $3$-connected matroid, 
let $\FF=(\hP_1,\ldots,\hP_n)$ be a maximal bloom of $M$
of order at least $4$, and let $(R,B)$ be a $3$-separation of 
$M$ that is not displayed by $\FF$. Then the following hold.
\begin{itemize}
\item[(i)] Up to labels, $R\subseteq \hP_1$.
\item[(ii)]  $|R\cap P_1^-|,|R\cap P_1^+|\leq 1$,
and if $R\cap P_1^-\neq \emptyset$, then $R\cap P_1^-$ is the last element of
$P_1^-$, and if $R\cap P_1^+\neq \emptyset$, then $R\cap P_1^+$ is the
first element of $P_1^+$.
\item[(iii)] If $|R\cap \cP_1|=1$, then $R$ is a triangle or a triad.
\item[(iv)] If $|B\cap \cP_1|=1$, then $R$ is equivalent to the displayed
$3$-separation $R\cup \cP_1$.
\end{itemize}
\end{lemma} 

\begin{proof}
Part (i) is just a restatement of Lemma~\ref{cross1}.
Consider (ii). Assume that $R\cap \cP_1=\emptyset$.
If $R\subseteq P_1^-$, or 
$R\subseteq P_1^+$,
then $R$ is clearly a consecutive set of loose elements
and is displayed by $\FF$, although a flower 
displayed by $\FF$ that
displayed $R$ would not be tight. Assume that 
$R\cap P_1^-\neq \emptyset$ and $R\cap P_1^+\neq \emptyset$.
In this case an appropriate concatenation of $\FF$ and
Lemma~\ref{skew-petals} show that $R\cap P_1^-$ and
$R\cap P_1^+$ are fully skew contradicting Lemma~\ref{split}.

Thus we may assume that
$R\cap \cP_1\neq\emptyset$. Say that $|R\cap L_1^-|\geq 2$. 
By uncrossing
$R$ with $\hP_n$, we see that $\lambda(\hP_n\cup R)=2$. But 
$R\cap \cP_1\neq\emptyset$, so $\hP_n\cup R$ is not displayed by
$\FF$ and neither it nor its complement is contained in a petal, contradicting
Lemma~\ref{cross1}. Hence $|R\cap P_1^-|\leq 1$. 
Assume that
$R\cap L_1^-=\{f\}$, where $\{f\}$ is not the last element 
of $P_1^-$. 
By Lemma~\ref{fine2},
$f\not\in\cl^{(*)}(R-\{f\})$, so $R$ is split 
and by Lemma~\ref{split},
$(R,B)$ is not a $3$--separation. 
Thus, if $R\cap P_1^-$ is non-empty,
then $R\cap P_1^-$ is the last element of $P_1^-$. 
The same argument 
holds for $R\cap P_1^+$. Thus (ii) holds. 
Parts (iii) and (iv) are easy.
\end{proof}



We now consider how 4-separations cross swirl-like flowers and blooms.
Observe that an unsplit exactly 4-separating set in a matroid has
at least four elements.

\begin{lemma}
\label{4-cross0}
Let $\FF=(P_1,\ldots,P_n)$ be a maximal swirl-like flower of 
the $3$-connected matroid $M$ such that $\FF$ has  order at least
$5$, and $|E(M)|\geq 11$. 
Let $(R,B)$ be an unsplit exact $4$-separation of $M$.
Assume that, amongst all equivalent flowers, $\FF$ has a minimum
number of petals crossed by $(R,B)$. Then, up to labels,
either
\begin{itemize}
\item[(i)] for some $i\in\{1,2,\ldots,n\}$ we have 
$P_1\cup P_2\udots P_{i-1}\subseteq B$; the partition
$(R,B)$ crosses $P_i$; and
$P_{i+1}\cup P_{i+2}\udots P_n\subseteq R$; or
\item[(ii)] $B$ crosses $P_1$ and $P_2$, and $B\subseteq P_1\cup P_2$.
\end{itemize}
\end{lemma}

\begin{proof}
Assume that the lemma fails and that $(R,B)$ satisfies
neither (i) nor (ii). We first show

\begin{sublemma}
\label{4-cross0.1}
There is a petal of $\FF$ that is not crossed by $(R,B)$.
\end{sublemma}

\begin{proof}[Subproof.]
Assume that all petals of $\FF$ are crossed by $(R,B)$.
By uncrossing, either $R\cap(P_1\cup P_2\cup P_3)$,
or $B\cap(P_4\cup P_5\udots P_n)$
is $3$-separating. Assume the former case holds.
Then, by Lemma~\ref{3-cross}, there is an equivalent
flower in which this set crosses at most one petal. Such a flower
has fewer petals crossed by $(R,B)$, contradicting the choice of $\FF$.
Thus the latter case holds. But we may apply the previous argument to
$B\cap(P_4\cup P_5\udots P_n)$ unless this set 
has at most two elements. This means
that $n=5$, and $|B\cap P_4|=|B\cap P_5|=1$. 
This argument extends to show
that  $|R\cap P_i|=|B\cap P_i|=1$
for all $i\in\{1,\ldots,5\}$, so that 
$|E(M)|=10$ contradicting the fact that $M$ has at least 11 elements.
\end{proof}

Thus we may assume that there is at least one petal
contained in $R$. 
Up to labels, we may assume that $\{P_2,P_3,\ldots,P_{i-1}\}$ is 
a maximal consecutive
set of petals each member of which 
is contained in $R$. Also we may assume that
there are at least two crossed petals as otherwise,
under the assumption that the lemma fails we can 
routinely deduce from
Lemma~\ref{skew-petals}, that either $R$ or $B$ is split.

\begin{sublemma}
\label{4-cross0.2}
Both $P_1$ and $P_i$ are crossed
by $(R,B)$.
\end{sublemma}

\begin{proof}[Subproof.]
If neither $P_1$ nor $P_i$ is crossed, then, by Lemma~\ref{skew-petals},
$R$ is split. Thus we may assume that $P_1$ is crossed. Assume that
$P_i\subseteq B$ and let $\{P_i,P_{i+1},\ldots,P_j\}$ be a maximal 
set of consecutive petals each member of which is contained in $B$.
If $P_{j+1}\subseteq R$, then by Lemma~\ref{skew-petals},
$B$ is split. Thus $P_{j+1}$ is crossed. A similar argument shows
that  $P_l$ is crossed for all $l\in \{j+1,j+2,\ldots,n\}$.
In particular, $P_n$ is crossed.
If $P_1=P_{j+1}$, then the lemma holds. Thus $j+1\leq n$.
By uncrossing, either $\lambda(P_1\cup P_2\udots P_{i-1})=2$
or $\lambda(P_i\cup P_{i+1}\udots P_n)=2$. 
In either case, by Lemma~\ref{3-cross},
there is a flower equivalent to $\FF$ that 
displays that $3$-separation.
But again by Lemma~\ref{3-cross}, such a flower has fewer petals
crossed by $(R,B)$. Note that for this to be true, we do need to 
know that $P_n$ is crossed.
\end{proof}

We now know that $P_1$ and $P_i$ both cross $(R,B)$. If $i=n$, then
the lemma holds, so $i<n$. If $l\in\{i+1,i+2,\ldots,n\}$, then 
$P_l\cap B\neq\emptyset$, otherwise, by Lemma~\ref{skew-petals}, 
$B$ is split. As $B$ is not split, $|B|\geq 4$.
Thus, either $|(P_i\cup P_{i+1}\udots P_n)\cap B|\geq 3$,
or $|(P_{i+1}\udots P_n\cup P_1)\cap B|\geq 3$. 
Up to labels we may assume the latter holds. Clearly 
$|(P_2\cup P_3\udots P_i)\cap R|\geq 3$. By uncrossing, one of these
sets is 3-separating. But again, in either case, it follows
straightforwardly from Lemma~\ref{3-cross} that a flower
equivalent to $\FF$ crosses $(R,B)$ in fewer petals.
This final contradiction to the choice of $\FF$ completes the proof.
\end{proof}

We will say that the unsplit exact 4-separation $(R,B)$ is
1-{\em crossing}
\index{$1$-crossing $4$-separation} 
or 2-{\em crossing}
\index{$2$-crossing $4$-separation} 
for $\FF$ according as
to whether Lemma~\ref{4-cross0}(i) or (ii) holds.
We now refine Lemma~\ref{4-cross0} to give a more precise 
description of the way that an unsplit $4$-separation crosses
a bloom.

\begin{lemma}
\label{1-crossing}
Let $\FF=(\hP_1,\hP_2,\ldots,\hP_n)$ be a bloom of the
$3$-connected matroid $M$, where 
$|E(M)|\geq 11$ and $n\geq 5$, and let $(R,B)$ be an
exact unsplit $4$-separation of $M$. If $(R,B)$ is $1$--crossing,
then, up to labels, $|R\cap\hP_1|,|B\cap \hP_1|\geq 2$, and  
one of the following holds.
\begin{itemize}
\item[(i)] $R\subseteq \hP_1$.
\item[(ii)] For some $i\in\{2,3,\ldots,n-1\}$,
$R=(\hP_1\cap R)\cup\cP_2\cup\hP_3\udots\hP_{i-1}\cup\cP_i\cup L_i^+$,
where $L_i^+$ is an initial section of $P_i^+$. Moreover,
$R\cap P_1^-$ is either empty or the last element of $P_1^-$,
and $B\cap P_1^+$ is either empty or the first element of 
$P_1^+$.
\end{itemize}  
\end{lemma}

\begin{proof}
By Lemma~\ref{4-cross0}, there is a flower $(P_1,P_2,\ldots,P_n)$
displayed by $\FF$ such that $(R,B)$ crosses
$P_1$ and $R=(P_1\cap R)\cup P_2\udots P_i$
for some $i\in\{1,2,\ldots,n\}$. Thus, $(R,B)$ crosses $P_1$.
Assume  that $|R\cap P_1|=1$, say $R\cap P_1=\{r\}$. Then
the set $R-\{r\}$ is displayed in the flower
$(P_1,P_2,\ldots,P_n)$ so that $\lambda(R-\{r\})=2$.
Hence $r\notin\clstar(R-\{r\})$ implying that $R$ is split.
Therefore $|R\cap\hP_1|\geq 2$ and similarly $|B\cap \hP_1|\geq 2$.

Assume that (i) does not hold. Then, as
$(R,B)$ is 1-crossing, we may assume up to labels that
$R=(\hP_1\cap R)\cup\cP_2\cup\hP_3\udots\hP_{i-1}\cup\cP_i\cup L_i^+$,
where $L_i^+$ is an initial section of $P_i^+$, and 
$i\in\{2,3,\ldots,n-1\}$. In particular $\cP_2\subseteq R$ and
$\cP_n\subseteq B$. Let $P_1^+=(f_1,f_2,\ldots,f_m)$.
Assume for a contradiction
that $P_1^+\cap B\notin \{\emptyset,\{f_1\}\}$.
Let $k$ be the greatest integer such that $f_k\in B$. By assumption,
$k>1$. If $f_{k-1}\in B$, then $R$ is split by Lemma~\ref{skew-petals}.
Thus $f_{k-1}\in R$. By Lemma~\ref{fine2},
$f_k\notin\clstar(\hP_3\udots\hP_n\cup\cP_1\cup
\{f_1,f_2,\ldots,f_{k-2}\})$
and this set contains $B-\{f_k\}$. Thus $f_k$ is isolated in
$B-\{f_k\}$, contradicting the fact that $B$ is unsplit.
Hence $P_1^+\cap B\in \{\emptyset,\{f_1\}\}$ and symmetrically
$P_1^-\cap R$ is either empty or the last element of $P_1^-$ so 
that (ii) holds.
\end{proof}

\begin{lemma}
\label{2-crossing}
Let ${\bf F}=(P_1,P_2,\ldots,P_n)$ be a maximal tight 
swirl-like flower of order at least $5$ in the $3$-connected
matroid $M$ and let  $(R,B)$ be
an unsplit exact $4$-separation of $M$ that is $2$-crossing for
$\FF$. Assume that $R\subseteq P_1\cup P_2$.
Then the following hold.
\begin{itemize}
\item[(i)] $\lambda(P_1\cap R)=\lambda(P_2\cap R)=2$.
\item[(ii)] $\{\sqcap(P_1\cap R,P_2\cap R),
\sqcap^*(P_1\cap R,P_2\cap R)\}=\{0,1\}$.
\end{itemize}
\end{lemma}

\begin{proof}
Assume that $|P_1\cap R|=1$; say 
$P_1\cap R=\{x\}$. As $(R,B)$ is not split,
$x\in\cl^{(*)}(P_2\cap R)$, and hence $x\in\cl(P_2)$, so that
$(P_1-\{x\},P_2\cup\{x\},P_3,\ldots,P_n)$ is a flower equivalent
to ${\bf F}$, and $(R,B)$ crosses only one petal of this flower,
contradicting the fact that $(R,B)$ is 2--crossing. 
Thus $|P_1\cap R|\geq 2$,
and hence $\lambda(P_1\cap R)\geq 2$.
Assume that $\lambda(P_1\cap R)>2$. Then, by uncrossing,
$\lambda(P_2\cup P_3 \cdots\cup P_n\cup(P_1\cap B))=2$. But, then,
by Lemma~\ref{cross1}, we can obtain a flower equivalent to 
${\bf F}$ that crosses $(R,B)$ in fewer petals. Hence 
$\lambda(P_1\cap R)=2$. Similarly $\lambda(P_2\cap R)=2$, so that
(i) holds.


By Lemma~\ref{lambda-meet} and (i), we have
\begin{align*}
3=\lambda(R)&=
\lambda(R\cap P_1)+\lambda(R\cap P_2)-\sqcap(R\cap P_1,R\cap P_2)
-\sqcap^*(R\cap P_1,R\cap P_2)\\
&=4-\sqcap(R\cap P_1,R\cap P_2)
-\sqcap^*(R\cap P_1,R\cap P_2).
\end{align*}
Part (ii) follows from this equation.
\end{proof}

\chapter{$k$-coherent Matroids}
\label{k-coherent-chapter}

At last we can introduce the key connectivity notion for
this paper.
Let $k\geq 5$ be an integer. Note that if a matroid
has a flower of order $n>k$, then it also has a flower of order
$k$. A matroid $M$ is $k$-{\em coherent\ }
\index{$k$-coherent}
if it is 3-connected and has no swirl-like flowers of order $k$.
The 3-connected matroid $M$ is $k$-{\em fractured\ }
\index{$k$-fracture} 
if it has a swirl-like flower of order $k$. 
A $k$-{\em fracture\ } of $M$ is a swirl-like flower in $M$ 
of order at least $k$.
This chapter 
is devoted to developing properties of $k$-coherent matroids,
focusing particularly on issues related to preserving 
$k$-coherence in
minors. In any unexplained context, if  we say that $M$ is a 
$k$-coherent matroid, then 
it will always be implicit that $M$ is 3-connected
and that $k$ is an integer greater than $4$.

An earlier draft of this paper had a number of purpose
built proofs that focused solely on $k$-coherence,
but it soon became apparent that they were special cases
of more general theorems related to the broader
issue of when it is possible to remove elements without
exposing 3-separations. This study of this question
was undertaken in 
\cite{wild,expose,upgrade}. We now
derive a number of our results as
corollaries of theorems in \cite{wild,expose,upgrade}.

We frequently illustrate situations for 
$k$-coherent matroids with schematic 
diagrams. In such diagrams it is always assumed that $k=5$.

\section{Exposure and $k$-coherence}

We begin by recalling some terminology.
Let $x$ be an element of the  $3$-connected matroid $M$. 
Assume that $M\ba x$ is $3$-connected. Then a 
$3$-separation $(A,B)$ of $M\ba x$ is 
{\em well-blocked\ } by $x$ if every $3$-separation of 
$M\ba x$ equivalent to $(A,B)$ is blocked by $x$. 
In the case that $(A,B)$ is well-blocked by $x$
we say that $(A,B)$ is  {\em exposed\ }
by $x$ or that $x$ {\em exposes $(A,B)$ in} $M/x$. 
On the other hand, if $M/x$ is 3-connected and 
$(A,B)$ is a $3$-separation in this matroid. Then
$(A,B)$ is {\em well-coblocked\ } by $x$ if every 
$3$-separation of $M/x$ equivalent to $(A,B)$ is 
coblocked by $x$. In this case we  say that
$x$ {\em exposes $(A,B)$ in $M/x$.} Recall also that,
if $(R,B)$ is a $3$-separation displayed by a flower $\FF$
of a matroid, then $(R,B)$ is {\em well-displayed\ } if 
both $R$ and $B$ contain at least two tight petals of $\FF.$
The first task of this chapter
is to prove that if $M$ is $k$-coherent, but 
$M\ba x$ is 3-connected and $k$-fractured, then
$x$ exposes a 3-separation in $M\ba x$. 
The next lemma is clear.

\begin{lemma}
\label{not-in-petal}
Let $M$ and $M\ba x$ be $3$-connected matroids. If
$x$ is in the closure of a petal of a $k$-fracture of $M\ba x$,
then $M$ is not $k$-coherent.
\end{lemma}

The next lemma shows that we cannot lose $k$-coherence without
exposing a 3-separation.

\begin{lemma}
\label{x-blocks}
Let $x$ be an element of the 
$k$-coherent matroid $M$ such that  
$M\ba x$ is $3$-connected and $k$-fractured.
Then $x$ exposes a $3$-separation in $M\ba x$.
Moreover, if $\FF$ is a $k$-fracture of $M\ba x$, 
then there is  a $3$-separation of $M\ba x$
that is well displayed by $\FF$ and is
exposed by $x$.
\end{lemma}

\begin{proof}
Assume that the lemma fails. 
Say $\FF=(\hP_1,\ldots,\hP_m)$ is a maximal
bloom that $k$-fractures $M\ba x$. If $(R,B)$ is a well-displayed 
$3$-separation in $\FF$, then either $x\in\cl_M(\fcl_{M\ba x}(R))$
or $x\in\cl_M(\fcl_{M\ba x}(B))$, as otherwise
this is a $3$-separation that is well blocked by $x$. 
Assume that, amongst all such 
$3$-separations $(R,B)$ is chosen 
so that $x\in\cl_M(\fcl_{M\ba x}(R))$
and that $\fcl_{M\ba x}(R)$ contains a minimum number of petals
of $\FF$. By Lemma~\ref{fine1}, we may assume, up to labels, that
$\fcl_{M\ba x}(R)=\hP_1\udots\hP_i$, for some $i\in\{2,\ldots,m-2\}$.
But now, either $x\in\cl(\hP_m\cup \hP_1)$, or 
$x\in\cl(\hP_2\udots \hP_{m-1})$, as otherwise we again
have a $3$-separation exposed by $x$. By Lemma~\ref{modular}
either $x\in\cl(\hP_1)$, or $x\in\cl(\hP_2)$, contradicting
the choice of $(R,B)$.
\end{proof}

The following theorem is \cite[Theorem~1.1]{upgrade}.

\begin{theorem}
\label{unexpose1}
Let $M$ be a $3$-connected matroid that is not a wheel or a 
whirl. Then $M$ has an element $x$ such that either 
$M\ba x$ or $M/x$ is $3$-connected and does not expose
any $3$-separations.
\end{theorem}

The next corollary follows immediately from Theorem~\ref{unexpose1}
and Lemma~\ref{x-blocks}.

\begin{corollary}
\label{unexpose2}
Let $M$ be a $k$-coherent matroid. If $M$ is neither a 
wheel nor a whirl, then $M$ has an element $x$ such
that either $M\ba x$ or $M/x$ is $k$-coherent.
\end{corollary}

\section{Some Local Wins}

While it is good to know that there is 
almost always an element somewhere
that can be removed to keep $k$-coherence,
we also need to identify specific locations where we
may remove elements and keep $k$-coherence.
In this section we describe a number of such situations.
Most often we focus on finding situations in which
removing an element does not expose any 
3-separations. We may omit the obvious corollary
for preserving $k$-coherence.
We begin by describing some straightforward cases. 
The next lemma is \cite[Lemma~2.10]{upgrade}. 
The proof is short,
so we give it here.

\begin{lemma}
\label{unexpose-sequential}
Let $M$ be a $3$-connected matroid and let $A$ be 
a sequential $3$-separating set of $M$ with at least
four elements. If $x\in A$ and $M\ba x$ is $3$-connected,
then $x$ does not expose any $3$-separations in 
$M\ba x$.
\end{lemma}

\begin{proof}
By Lemma~\ref{delete-seq}, $A-\{x\}$ is a sequential 
$3$-separator of $M\ba x$. By Lemma~\ref{seq-ord} 
$A-\{x\}$ has a triangle
or triad $T$ such that $\fcl_{M\ba x}(T)\supseteq A-\{x\}$.
Let $(R,B)$ be a
$3$-separation in $M\ba x$. 
Up to labels we may assume that $|R\cap T| \geq 2$. 
Then $T\subseteq \clstar_{M\ba x}(R\cap T)$, so that 
$A-\{x\}\subseteq \fcl_{M\ba x}(R)$. In $M\ba x$, the
3-separation $(R,B)$
is equivalent to $(\fcl_{M\ba x}(R),B-\fcl_{M\ba x}(R))$. But 
$A\subseteq \fcl_{M\ba x}(R)$, and $x\in\cl_M(A-\{x\})$, so that 
$(\fcl_{M\ba x}(R),B-\fcl_{M\ba x}(R))$ is induced in $M$ and 
$(R,B)$ is not exposed by $x$.
\end{proof}

For $k$-coherent matroids we have

\begin{corollary}
\label{remove-sequential}
Let $M$ be a $k$-coherent matroid and let 
$A$ be  a sequential $3$-separating set in
$M$ such that $|A|\geq 4$.
\begin{itemize}
\item[(i)] If $x$ is an element
of $A$ such that $M\ba x$ is $3$-connected,
then $M\ba x$ is $k$-coherent.
\item[(ii)] If $(a_1,a_2,\ldots,a_n)$ is a sequential
ordering of $A$, then either $M\ba a_n$ or $M/a_n$ is 
$k$-coherent.
\end{itemize}
\end{corollary}

\begin{proof}
Part (i) follows from Lemmas~\ref{x-blocks} and
\ref{unexpose-sequential}. Part (ii) follows from (i)
and Lemma~\ref{seq-3-con}.
\end{proof}

Another easy case, also proved in \cite{upgrade}, 
is provided by quads. The proof is similar
to that of Lemma~\ref{unexpose-sequential}, but easier.

\begin{lemma}
\label{dag-remove}
Let $D$ be a quad of the $3$-connected matroid $M$ and 
let $d$ be an element
of $D$ that is not in a triad. Then
$M\ba d$ is $3$-connected and $d$ does not expose any
$3$-separations in $M\ba d$. 
\end{lemma}

Thus an element of a quad can be removed one way or
the other to keep 3-connectivity and not expose 
3-separations unless it
is in both a triangle and triad, that is, unless it is in a 
4-element fan. One may expect it to be
difficult for every element of a quad to be in a 
4-element fan, but it can happen. Indeed, let $M$ be a spike
with a tip $t$ and cotip $c$ and let $\{p_1,q_1\}$ 
and $\{p_2,q_2\}$ be legs of the spike. Then $\{p_1,q_1,p_2,q_2\}$
is a quad and $(t,p_1,q_1,c)$ is a maximal fan with
ends $t$ and $c$.

\begin{lemma}
\label{daggy-dag}
Let $D$ be a quad of the $3$-connected matroid $M$. If
$D$ has an element that is in both a triangle and a triad,
then every element of $D$ is in a $4$-element fan and 
there is a partition $(D_1,D_2,C)$ of $E(M)$ with 
$D_1\cup D_2=D$ such that the following hold.
\begin{itemize}
\item[(i)] $(D_1,D_2,C)$ is a spike like flower with tip $t\in C$
and cotip $c\in C$. 
\item[(ii)] $D_1\cup\{c,t\}$ and $D_2\cup \{c,t\}$ are $4$-element
fans of $M$.
\end{itemize}
\end{lemma}

\begin{proof}
Say $d\in D$, and that $d$ belongs to a
4-element fan $F$. As $d$ is in both a triangle and a triad,
$d$ is an internal element of $F$. Let $F=(t,d,x,c)$, where
$\{t,d,x\}$ is a triangle and $\{d,x,c\}$ is a triad.
As $D$ is both a circuit and a cocircuit we see that both
$\{t,d,x\}$ and $\{d,x,c\}$ contain exactly 
two elements of $D$ so that $x\in D$ and 
$\{t,c\}\subseteq E(M)-D$. Let $D_1=\{d,x\}$ and 
$D_2=D-D_1$. Consider the flower $\FF=(C,D_1,D_2)$.
Note that $t$ and $c$ are loose elements of this flower,
indeed, $t\in\cl(D_1)$, and $c\in\cl^*(D_1)$. By 
Lemma~\ref{get-spike}, $\FF$ is spike like, so that
$t\in\cl(D_2)$ and $c\in\cl^*(D_2)$. The lemma
now follows routinely.
\end{proof}

The next theorem is \cite[Theorem~7.1]{upgrade}. Its proof
is quite substantial.

\begin{theorem}
\label{remove-a}
Let $(P,\{a,b\},Q)$ be a tight flower of a $3$-connected
matroid $M$ where $\{a,b\}$ is fully closed and both 
$P$ and $Q$ have at least three elements.  
Then the following hold. 
\begin{itemize}
\item[(i)] If $a$ is in a triangle, then 
$M\ba a$ is $3$-connected and has no 
$3$-separations exposed by $a$. 
\item[(ii)] If $a$ is in a triad, then $M/ a$ is $3$-connected 
and has no $3$-separations exposed by $a$. 
\item[(iii)] If $a$ is in neither a triangle nor a triad, 
then both $M\ba a$ and $M/a$ are $3$-connected.
\end{itemize}
Moreover, if $a$ is in neither a triangle nor a triad and  
both $M\ba a$ and $M/a$ have $3$-separations exposed by $a$, 
then $|P|=|Q|=4$, both $M\ba b$ and $M/b$ are $3$-connected, 
and neither $M\ba b$ nor $M/b$ has a 
$3$-separation exposed by $b$. 
\end{theorem}

Combining Theorem~\ref{remove-a} with other facts we
obtain the following consequence.

\begin{corollary}
\label{2-elt-win}
Let $(P,\{a,b\},Q)$ be a tight flower of the $k$-coherent
matroid $M$ where $\{a,b\}$ is fully closed. Then exactly
one of the following holds.
\begin{itemize}
\item[(i)] Either $M\ba a$ or $M/a$ is $k$-coherent.
\item[(ii)] Up to labels, $P\cup\{a,b\}$ is a quad.
Moreover, there is a labelling $p_1,p_2$ of the elements
of $P$ such that $(\{p_1,a\},\{p_2,b\},Q)$
is a spike-like flower with a tip and a cotip.
\end{itemize}
\end{corollary}

\begin{proof}
Note that a $3$-connected matroid with at most nine elements
is $k$-coherent. By
Theorem~\ref{remove-a} and Lemma~\ref{x-blocks} part (i) holds 
unless either $|P|=2$ or
$|Q|=2$. Assume that $|P|=2$. Then, as the flower is
tight, $P\cup \{a,b\}$ is a quad in $M$. If $a$ is not in both a 
triangle and a triad, then it follows from
Lemma~\ref{dag-remove} that (i) holds.
Otherwise it follows from 
Lemma~\ref{daggy-dag} that (ii) holds.
\end{proof}

\begin{figure}
\begin{tikzpicture}[thick,line join=round]
	\coordinate (x) at (0:2);
	\coordinate (y) at (-43:2);
	\coordinate (z) at (-87:2);
	\coordinate (w) at (-130:2);
	\coordinate (t) at (130:2);
	\coordinate (v) at (65:2);
	\node[pattern color=lines,draw,circle through=(z),pattern=north east lines] at ($(w)!0.5!(z)$) {};
	\node[pattern color=lines,draw,circle through=(y),pattern=north east lines] at ($(z)!0.5!(y)$) {};
	\node[pattern color=lines,draw,circle through=(x),pattern=north east lines] at ($(y)!0.5!(x)$) {};
	\node[pattern color=lines,draw,circle through=(v),pattern=north east lines] at ($(x)!0.5!(v)$) {};
	\node[pattern color=lines,draw,circle through=(t),pattern=north east lines] at ($(v)!0.5!(t)$) {};
	\filldraw[fill=white] (w) -- (z) -- (y) -- (x) -- (v) -- (t) -- cycle;
	\draw[dashed] (t) -- (x) -- (w);
	\coordinate[label=left:$b$] (b) at ($(t)!0.5!(w)$);
	\draw[dashed] (b) -- (x);
	\coordinate[label=above:$a$] (a) at ($(b)!0.4!(x)$);
	\foreach \pt in {a,b} \fill[black] (\pt) circle (3pt);
	\node at ($(t)!0.5!(v) + (0,0.5)$) [rectangle,fill=white,draw=white,inner sep=1pt] {$P_1$};
	\node at ($(x)!0.5!(v) + (0.4,0.3)$) [rectangle,fill=white,draw=white,inner sep=1pt] {$P_2$};
	\node at ($(x)!0.5!(y) + (0.3,-0.15)$) [rectangle,fill=white,draw=white,inner sep=1pt] {$P_3$};
	\node at ($(z)!0.5!(y) + (0.05,-0.35)$) [rectangle,fill=white,draw=white,inner sep=1pt] {$P_4$};
	\node at ($(z)!0.5!(w) + (-0.05,-0.35)$) [rectangle,fill=white,draw=white,inner sep=1pt] {$P_5$};
\end{tikzpicture}
\caption{Outcome (i) of Corollary~\ref{2-elt-win}}
\label{2-elt-win1}
\end{figure}
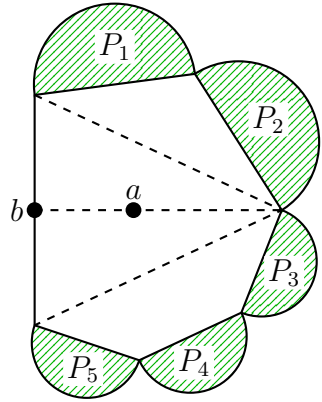

We illustrate Corollary~\ref{2-elt-win} with some examples.
Consider the matroid $M$ illustrated in
Figure~\ref{2-elt-win1}. Let $P=P_1\cup P_2$
and $Q=P_3\cup P_4\cup P_5$. Then the flower $(P,\{a,b\},Q)$
satisfies the hypotheses of Corollary~\ref{2-elt-win}.
Assume that $M$ is 5-coherent. Note that 
$(P_1\cup\{b\},P_2,P_3,P_4,P_5)$ is a 5-fracture of $M\ba a$,
so that $M\ba a$ is not 5-coherent. Nonetheless, $M/a$ is
5-coherent. Due to the way that $a$ and $b$ are aligned
in this particular case, $(P_1,P_2\cup\{b\},P_3,P_4,P_5)$
is a 5-fracture of $M/b$. In this case $M\ba b$ is 5-coherent.

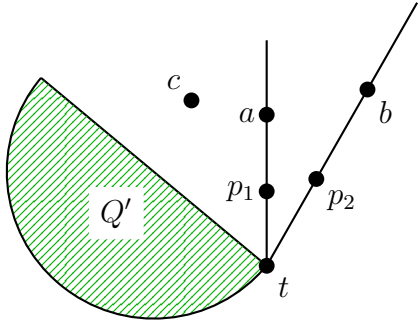
\begin{figure}
\begin{tikzpicture}[thick,line join=round]
	\coordinate[label=-45:$t$] (t) at (0,0);
	\coordinate (x) at (-3,2.5);
	\coordinate (y) at (0,3);
	\coordinate (z) at (2,3.5);
	\coordinate (xt) at ($(x)!0.5!(t)$);
	\coordinate[label=left:$a$] (a) at ($(y)!0.33!(t)$);
	\coordinate[label=-45:$b$] (b) at ($(z)!0.33!(t)$);
	\coordinate[label=left:$p_1$] (p1) at ($(t)!0.33!(y)$);
	\coordinate[label=-45:$p_2$] (p2) at ($(t)!0.33!(z)$);
	\coordinate[label=135:$c$] (c) at (-1,2.2);
	\draw (z) -- (t) -- (y);
	\draw (x) -- (t);
	\filldraw[pattern color=lines,pattern=north east lines] let \p1 = ($(t) - (xt)$),
		\n1 = {veclen(\x1,\y1)},
		\n2 = {atan(-3/2.5)} in
			(x) arc (-\n2+90:180-\n2+90:\n1);
	\foreach \pt in {t,a,b,p1,p2,c} \fill[black] (\pt) circle (3pt);
	\node at (-2,0.7) [rectangle,fill=white,draw=white] {$Q'$};
\end{tikzpicture}
\caption{Outcome (ii) of Corollary~\ref{2-elt-win}}
\label{2-elt-win2}
\end{figure}

The next example illustrates the need for outcome (ii) in
Corollary~\ref{2-elt-win}.
Consider the matroid illustrated in Figure~\ref{2-elt-win2}.
Let $Q=Q'\cup\{c,t\}$. Note that $(\{b,p_2\},\{a,p_1\},Q)$
is a respectable spike-like flower with tip $t$ and cotip $c$.
Unfortunately, we were given the perverse flower
$(\{p_1,p_2\},\{a,b\},Q)$ in which $\{a,b\}$ is tight and
fully closed.

\begin{lemma}
\label{2-3-4-remove}
Let $M$ be a $k$-coherent matroid, 
and let $\RR=(R_1,R_2,\ldots,R_m)$ be a flower of 
$M$, where $m\geq 4$ and $|R_1|\in\{2,3,4\}$. If 
$p\in R_1$ and $M\ba p$ is $3$-connected, then 
$M\ba p$ is $k$-coherent.
\end{lemma}

\begin{proof}
Assume that $|R_1|=2$, say $R_1=\{p,p'\}$. 
Let $\PP=(\{p,p'\},P_1,P_2,P_3)$ be a 4-petal concatenation
of $\RR$. As $M\ba p$ is $3$-connected,
$p\in\cl_M(P_1\cup\{p'\})$ and $p\in\cl_M(P_3\cup\{p'\})$.
Clearly $(\{p'\}\cup P_1,P_2,P_3)$ and
$(P_1,P_2,P_3\cup\{p'\})$ are flowers in $M\ba p$.

Assume that the lemma fails and let 
$\QQ=(\hQ_1,\hQ_2,\ldots,\hQ_n)$ be a bloom of 
$M\ba p$ that $k$-fractures $M\ba p$. Assume that
$P_1\cup\{p'\}$ is not displayed by $\FF$. Then, by 
Lemma~\ref{3-cross}, either 
$P_1\cup\{p'\}$ or $P_2\cup P_3$ is contained in 
$\hQ_i$ for some $i\in\{1,2,\ldots,n\}$. The former
case implies that $p\in\cl_M(\hQ_i)$, so that, 
by Lemma~\ref{not-in-petal},  
$M$ is $k$-fractured. Consider the latter case.
Evidently $p'\in\clstar_{M\ba p}(P_3)$ so that 
$p'\in\clstar_{M\ba p}(\hQ_i)$, a fully closed set.
Thus $P_3\cup\{p'\}\subseteq \hQ_i$. Again we 
see that $p\in\cl_M(\hQ_i)$ and again we obtain 
a contradiction to the fact that $M$ is $k$-coherent.

Thus $P_1\cup\{p'\}$ and, similarly $P_3\cup\{p'\}$
are both displayed in $\QQ$ and neither is 
contained in a petal of $\QQ$. It is now clear that
we have a quasi-flower $(\{p'\},Q'_1,Q'_2,Q'_3,Q'_4)$
displayed by $\QQ$, where $P_1=Q'_1\cup Q'_2$, $P_2=Q'_3$
and $P_3=Q'_4$, and having the further property
that, for some $i\in\{1,2,\ldots,n\}$ we have
$\cQ_i\subseteq Q'_1\subseteq Q'_1\cup\{p'\}\subseteq \hQ_i$.
Now $p\in\cl_M(\{p'\}\cup Q'_1\cup Q'_2)$
and $p\in\cl_M(Q'_4\cup Q'_1\cup\{p'\})$,
so, by Lemma~\ref{modular}, $p\in\cl_M(\{p'\}\cup Q'_1)$.
Thus $p\in\cl(\hP_i)$. Again we obtain the contradiction that
$M$ is $k$-fractured. Thus the lemma holds in the case
that $|R_1|=2$.

A similar, but easier, analysis establishes the lemma
in the case when $R_1$ has three elements.
The case when $|R_1|=4$ is too easy to resist. 
In this  
case, if $R_1$ is  sequential, then it follows
from Corollary~\ref{remove-sequential} that $M\ba x$
is $k$-coherent. On the other hand, if $R_1$ is 
non-sequential, then $R_1$ is a quad, and it follows
from Lemma~\ref{dag-remove} that $M\ba x$
is $k$-coherent.
\end{proof}

If a triangle $T$ is not in a $4$-element fan then,
as we shall see in the next section, it is not always
possible to find an element of $T$ that can be 
deleted to preserve $k$-coherence. 
Nonetheless things often work out well.
The next lemma is 
\cite[Theorem~4.2]{wild}.

\begin{lemma}
\label{a-triangle-win}
Let $\{a,b,c\}$ be a triangle in a $3$-connected matroid M
having at least four element.
Assume that $\{a,b,c\}$ is not contained in a $4$-element fan,
and that $M\ba b$ is not $3$-connected. If 
$z\in\{a,b\}$ then $M\ba z$ is $3$-connected and $z$ does
not expose a $3$-separation in $M\ba z$.
\end{lemma}

\begin{figure}
\begin{tikzpicture}[thick,line join=round]
	\coordinate (A) at (0,0);
	\coordinate (B) at ($1.5*(-2,1.25)$);
	\coordinate (C) at ($1.5*(2,1.25)$);
	\coordinate (D) at ($1.5*(0,3.5)$);
	\coordinate (ab) at ($(A)!0.5!(B)$);
	\coordinate (ac) at ($(A)!0.5!(C)$);
	\filldraw[pattern color=rlines,pattern=north east lines] let \p1 = ($(ab) - (A)$),
		\n1 = {veclen(\x1,\y1)},
		\n2 = {atan(1.25/2)} in
			(B) arc (-\n2-180:180-\n2-180:\n1);
	\filldraw[pattern color=blines,pattern=north west lines] let \p3 = ($(ab) - (A)$),
		\n3 = {veclen(\x3,\y3)},
		\n4 = {atan(2/1.25)} in
			(A) arc (-\n4-90:180-\n4-90:\n3);
	\draw (A) -- (B) -- (D) -- (C) -- cycle;
	\coordinate[label=180:$p_1$] (p1) at ($(B)!0.25!(D)$);
	\coordinate[label=180:$p_2$] (p2) at ($(D)!0.25!(B)$);
	\coordinate[label=right:$p'_2$] (p2') at ($(C)!0.25!(D)$);
	\coordinate[label=right:$p'_1$] (p1') at ($(D)!0.25!(C)$);
	\draw (p1) -- (p1');
	\draw (p2) -- (p2');
	\coordinate[label=below:$l$] (l) at (intersection of p1--p1' and p2--p2');
	\fill[black] (l) circle (3pt);
	\foreach \pt in {p1,p2} \fill[rlines] (\pt) circle (3pt);
	\foreach \pt in {p1',p2'} \fill[blines] (\pt) circle (3pt);
	\node at ($1.5*(1.2,0)$) [rectangle,fill=white,draw=white] {$P_1{''}$};
	\node at ($1.5*(-1.2,0)$) [rectangle,fill=white,draw=white] {$P'_1$};
\end{tikzpicture}
\caption{Illustration for Lemma~\ref{loose-removable}}
\label{lose-removable}
\end{figure}
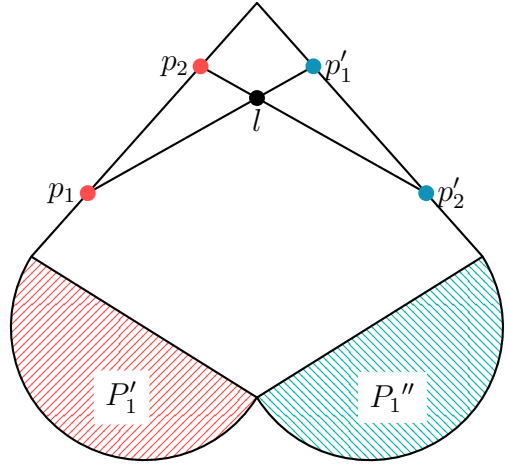

We now consider loose elements in swirl-like or spike-like
flowers. It turns out 
that, apart from some almost degenerate situations
we can remove loose elements of such flowers and 
preserve $k$-coherence.
The degeneracy is again related to the fact that
we may partition quads arbitrarily into 2-element 
subsets to obtain 3-petal flowers.
Figure~\ref{lose-removable} illustrates the needs for the
constraints in part (iii) of Lemma~\ref{loose-removable}.
Here $P_1=\{p_1,p'_1\}$, $P_2=\{p_2,p'_2\}$
and $P_3=P'_3\cup P''_3$. Deleting $l$ exposes the
$3$-separation $(\{p_1,p_2\}\cup P'_3, \{p'_1,p'_2\}\cup P''_3)$.

\begin{lemma}
\label{loose-removable}
Let $M$ be a $3$-connected  matroid  and let 
$\FF=(P_1\cup \{l\},P_2,\ldots,P_n)$ be a 
spike-like or swirl-like flower of 
$M$ with at least three petals and $|P_1|\geq 2$. 
Assume that $l$ is a loose element between 
$P_1\cup\{l\}$ and $P_2$ and that $M\ba l$ is 
$3$-connected.  
\begin{itemize}
\item[(i)] If either $P_1\cup\{l\}$ or $P_2$ is a
loose petal, then $l$ does not expose any $3$-separations
in $M\ba l$.
\item[(ii)] If $n\geq 4$, then $l$ does not expose any
$3$-separations in $M\ba l$.
\item[(iii)] If $n=3$ and at most one member of 
$\{P_1,P_2,P_3\}$ has two elements, then $l$ does not
expose any $3$-separations in $M\ba l$.
\end{itemize}
\end{lemma}

\begin{proof}
We first show that (i) holds. Assume that $P_1\cup\{l\}$
is a loose petal. Then, as $|P_1\cup\{l\}|\geq 3$,
it must be that case that $\FF$ is swirl-like.
By Theorem~\ref{swirl-like}(ii), $P_1\cup\{l\}$ is a fan.
If $|P_1\cup\{l\}|\geq 4$, then it follows from 
Corollary~\ref{remove-sequential} that $l$ does not expose 
a 3-separation in $M\ba l$. Thus we may assume that
$|P_1\cup\{l\}|=3$. Say that, written as
an ordered fan, $P_1\cup\{l\}=(f_1,f_2,l)$. 
In this case $\{f_1,f_2,l\}$ is a triangle
and is not contained in a 4-element fan as otherwise we may
again apply Corollary~\ref{remove-sequential}. 
As $l$ is loose, $l\in\clstar(P_2)$. If $l\in\cl^*(P_2)$,
then $M\ba l$ is not $3$-connected. Thus $l\in\cl(P_2)$.
It now follows from Lemma~\ref{fine2} that 
$f_2\in\cl^*(P_2\cup\{l\})$, so that $M\ba f_2$
is not 3-connected. In this case
$l$ does not expose any $3$-separations in 
$M\ba l$ by Lemma~\ref{a-triangle-win}. The same argument applies
in the case that $P_2$ is loose and (i) holds.

Consider (ii) Assume that $\FF$ has at least four petals.
Assume that $(R,B)$ is a 3-separation of $M\ba l$
that is exposed by $l$. 
By (i) both 
$P_1\cup \{l\}$ and $P_2$ are tight petals. 
Let $\FF''=(P_1,P_2,P'_3,P'_4)$ be a concatenation of 
$(P_1,P_2,\ldots,P_n)$ and let $\FF'$ be a 
maximal flower such that $\FF''\less \FF'$ and such that
$\FF'$ displays $\FF''$. Such a flower exists by 
Lemma~\ref{display-small}.
Say 
$$\FF'=(P_{11},\ldots,P_{1i},P_{21},
\ldots,P_{2j},P_{31},\ldots,P_{3k},P_{41},
\ldots,P_{4m}),$$
where $P_1=P_{11}\udots P_{1i}$ and 
$P_2=P_{21}\udots P_{2j}$, $P'_3=P_{31}\udots P_{3k}$,
and $P'_4=P_{41}\udots P_{4m}$.
Then $l\in\cl(P_1)$ and $l\in\cl(P_2)$, so, by Lemma~\ref{page3},
$l\in\cl(P_{1i})$ and $l\in\cl(P_{21})$. As $\FF'$ is maximal,
$(R,B)$ conforms with $\FF'$. Thus there is a $3$--separation
$(R',B')$, equivalent to $(R,B)$ such that either $(R',B')$
is displayed by $\FF'$, or one of $R'$ or $B'$ is contained in 
a petal of $\FF'$. In either case we deduce 
that either $P_{1i}$ or $P_{21}$
does not cross $(R',B')$. But then, as $l\in\cl(P_{1i})$ and
$l\in\cl(P_{21})$, we see that either $l\in\cl(R')$  
or $l\in\cl(B')$,
contradicting the assumption that $(R,B)$ is well-blocked by $l$.
Thus (ii) holds.

Consider (iii). Here $\FF=(P_1\cup\{l\},P_2,P_3)$.
Again assume that $(R,B)$ is a 3-separation of 
$M\ba l$ exposed by $l$.
Consider the flower $(P_1,P_2,P_3)$ of $M\ba l$.
As $l\in\cl(P_1)$, and $l\in\cl(P_2)$, $(R,B)$, and any 
$3$-separation of $M\ba l$ equivalent to $(R,B)$, 
crosses both $P_1$ and $P_2$.
As at most one member of $\{P_1,P_2,P_3\}$ has two elements,
we may assume up to 
labels that $|P_1|>2$ and that $|P_1\cap R|>2$.

\begin{sublemma}
\label{loose-removable2}
$|P_1\cap B|\geq 2$.
\end{sublemma}

\subproof
Assume that $|B\cap P_1|=1$; let $\{b\}=B\cap P_1$. We have
$|B-P_1|\geq 2$, so that $R\cup P_1$ avoids at least two elements of
$E(M\ba l)$. As $|R\cap P_1|\geq 2$, by uncrossing we see that
$R\cup P_1$ is $3$--separating. But $R\cup P_1=R\cup\{b\}$. Hence
$(R,B)\cong (R\cup \{b\},B-\{b\})$, but this latter $3$--separation
does not cross $P_1$.
\end{proof}

A similar uncrossing argument to that of 
\ref{loose-removable2} shows  that $|P_2\cap R|=1$ if and only
if $|P_2\cap B|=1$. Using an uncrossing argument and possibly
moving to a $3$--separation equivalent to $(R,B)$ we may also
assume that $|P_3\cap R|=1$ if and only if $|P_3\cap B|=1$.

\begin{sublemma}
\label{loose-removable3}
$|P_i\cap R|=|P_i\cap B|=1$ for some $i\in\{2,3\}$.
\end{sublemma}

\begin{proof}[Subproof.]
If neither $P_2$ nor $P_3$ has the property of the sublemma, 
then we may apply 
Lemma~\ref{flower8.2} and obtain a flower $\FF'$ that refines $(P_1,P_2,P_3)$
and displays a 3-separation equivalent to $(R,B)$. 
As $l\in\cl(P_1)$ and $l\in\cl(P_2)$,
it follows from Lemma~\ref{page3}, that, for some petal 
$P$ of $\FF'$, we have $l\in\cl(P)$. But, as $(R,B)$ is displayed,
either $P\subseteq R$ or $P\subseteq B$. In either case we contradict
the assumption that $(R,B)$ is well-blocked by $l$.
\end{proof}

\begin{sublemma}
\label{loose-removable4}
$|P_2\cap R|=|P_2\cap B|=|P_3\cap R|=|P_3\cap B|=1$.
\end{sublemma}

\begin{proof}[Subproof.]
Let $(Q_2,Q_3)$ be a permutation of $(P_2,P_3)$ and assume for a 
contradiction that 
$|R\cap Q_2|\geq 2$. By \ref{loose-removable3},
$|Q_3\cap R|=|Q_3\cap B|=1$. Let $(r,b)=(Q_3\cap R,Q_3\cap B)$. 

Assume that $Q_2\cap B=\emptyset$. Then, by uncrossing,
$(R\cup\{b\},B-\{b\})$ is a $3$-separation in $M\ba l$.
But $(R\cup\{b\},B-\{b\})\cong (R,B)$ and $B-\{b\}\subseteq P_1$,
so $(R\cup\{b\},B-\{b\})$ is not blocked by $l$, contradicting the assumption
that $(R,B)$ is well-blocked by $l$. Thus $Q_2\cap B\neq \emptyset$ and,
as $|Q_2\cap R|\geq 2$, we deduce that $|Q_2\cap B|\geq 2$.

Now 
$|E(M\ba l)-(B\cup(P_1\cup Q_3))|\geq 2$ and 
$|B\cap (P_1\cup Q_3)|\geq 2$.
So by uncrossing, $\lambda_{M\ba l}(B\cap (P_1\cup Q_3))=2$. Note that
$B\cap (P_1\cup Q_3)=(B\cap P_1)\cup\{b\}$. 
A similar uncrossing argument shows that 
$\lambda_{M\ba l}(B\cap P_1)=2$. 
As $B\cap P_1$ and $(B\cap P_1)\cup\{b\}$
are both exactly $3$--separating in $M\ba l$, it follows that
$b\in\clstar_{M\ba l}(B\cap P_1)$ and hence, in $M$, we have
$b\in\clstar_M(P_1\cup\{l\})$. By symmetry we also deduce
that $r\in\clstar_M(P_1\cup\{l\})$. Thus $r$ and $b$ are loose elements
of the flower $\FF$,
and by the structure of loose elements in 
swirl-like or spike-like flowers, we deduce that,
up to labels, we have $b\in\cl_M(P_1\cup\{l\})$. 
As $l\in\cl_M(P_1)$,
we have $\{b,l\}\subseteq \cl_M(P_1)$. But, in the flower 
$(P_1,P_2\cup\{l\},P_3)$, the pair
$\{b,l\}$ is a subset of the closure of a single
petal and, by the structure of swirl-like and spike-like flowers,
there is at most one element
in the petal $P_2\cup\{l\}$ contained in the closure of $P_1$.
This contradiction establishes the sublemma.
\end{proof}

By \ref{loose-removable4}, $|P_2|=|P_3|=2$ contradicting the
assumption that at most one of these sets has two element and (iii)
follows.
\end{proof}

\section{$k$-wild Triangles}
Let $\{a,b,c\}$ be a triangle of the $k$-coherent matroid
$M$ that is not in a 4-element fan. We have seen in 
Lemma~\ref{a-triangle-win} that if $M\ba b$ is not 
$3$-connected, then both $M\ba a$ and $M\ba c$ are 
$k$-coherent. 
Unfortunately the case arises where 
$M\ba a$, $M\ba b$, and $M\ba c$ are 3-connected and $k$-fractured. 
The remainder of this section is devoted to 
describing the structure of
a matroid relative to such a triangle. 
The results of this section are
essentially corollaries of results in \cite{wild}. 

We begin by recalling the $\Delta-Y$ 
\index{$\Delta-Y$ exchange}
operation for matroids. 
Let $\Delta=\{a,b,c\}$ be a triangle of the matroid $M$ and take
a copy of $M(K_4)$ having $\Delta$ as a triangle and
$\{a',b',c'\}$ as the complementary triad, labelled such that
$\{a,b',c'\}$, $\{a',b,c'\}$ and $\{a',b',c\}$ are triangles.
Let $P_\Delta(M(K_4),M)$ denote the generalised parallel
connection of $M(K_4)$ and $M$. 
We write $\Delta M$ for $P_\Delta(M(K_4),M)\ba \Delta$ and say that
$\Delta M$ is obtained from $M$ by a 
$\Delta-Y$ {\em exchange\ } on $\Delta$.
As is common practice, we relabel $a'$, $b'$ and $c'$
as $a$, $b$ and $c$ so that $M$ and $\Delta M$ have the
same ground set.
The matroid $N$ is obtained from $M$ by performing a 
$Y-\Delta$ exchange on a triad $\{a,b,c\}$ if 
$N^*$ is obtained from $M^*$ by performing a 
$\Delta-Y$ exchange on the triangle $\{a,b,c\}$ of $M^*$.

A triangle $T$ of the 3-connected matroid $M$ is {\em wild\ }
\index{wild triangle}
if, for all $t\in T$, either $M\ba t$ is 
not 3-connected, or $M\ba t$ is 3-connected
and exposes a 3-separation in $M$. 
The structure of wild triangles is described
in \cite{wild}, and we will later outline 
the results from there that we
need. First we describe the particular types of wild triangle that
is problematic from the perspective of $k$-coherence. 

\begin{figure}
\begin{tikzpicture}[thick,line join=round]
	\coordinate (x) at (90:2);
	\coordinate (y) at (210:2);
	\coordinate (z) at (330:2);
	\coordinate (q) at ($(z)+(-120:1.8)$);
	\coordinate (r) at ($(y)+(-60:1.8)$);
	\coordinate (u) at ($(x)+(1.8,0)$);
	\coordinate (v) at ($(z)+(60:1.8)$);
	\coordinate (s) at ($(x)+(-1.8,0)$);
	\coordinate (t) at ($(y)+(120:1.8)$);
	\node[pattern color=lines,draw,circle through=(x),pattern=vertical lines] at ($(u)!0.5!(x)$) {};
	\node[pattern color=lines,draw,circle through=(s),pattern=vertical lines] at ($(x)!0.5!(s)$) {};
	\node[pattern color=lines,draw,circle through=(t),pattern=north west lines] at ($(s)!0.5!(t)$) {};
	\node[pattern color=lines,draw,circle through=(y),pattern=north east lines] at ($(t)!0.5!(y)$) {};
	\node[pattern color=lines,draw,circle through=(r),pattern=north east lines] at ($(y)!0.5!(r)$) {};
	\node[pattern color=lines,draw,circle through=(q),pattern=vertical lines] at ($(r)!0.5!(q)$) {};
	\node[pattern color=lines,draw,circle through=(z),pattern=north west lines] at ($(q)!0.5!(z)$) {};
	\node[pattern color=lines,draw,circle through=(v),pattern=north west lines] at ($(z)!0.5!(v)$) {};
	\node[pattern color=lines,draw,circle through=(u),pattern=north east lines] at ($(v)!0.5!(u)$) {};
	\filldraw[fill=white] (t) -- (s) -- (u) -- (v) -- (q) -- (r) -- cycle;
	\draw (x) -- (y) -- (z) -- cycle;
	\coordinate[label=45:{\textcolor{labels}{$a$}}] (a) at ($(x)!0.5!(z)$);
	\coordinate[label=below:{\textcolor{labels}{$b$}}] (b) at ($(y)!0.5!(z)$);
	\coordinate[label=135:{\textcolor{labels}{$c$}}] (c) at ($(x)!0.5!(y)$);
	\draw (a) parabola bend (b) (c);
	\foreach \pt in {a,b,c} \fill[labels] (\pt) circle (3pt);
	\node at ($(x)!0.5!(u) + (0,0.35)$) [rectangle,fill=white,draw=white,inner sep=1pt] {\textcolor{labels}{$A_1$}};
	\node at ($(u)!0.5!(v) + (0.35,0.15)$) [rectangle,fill=white,draw=white,inner sep=1pt] {\textcolor{labels}{$A_2$}};
	\node at ($(v)!0.5!(z) + (0.35,-0.15)$) [rectangle,fill=white,draw=white,inner sep=1pt] {\textcolor{labels}{$A_3$}};
	\node at ($(z)!0.5!(q) + (0.35,-0.15)$) [rectangle,fill=white,draw=white,inner sep=1pt] {\textcolor{labels}{$B_1$}};
	\node at ($(q)!0.5!(r) + (0,-0.35)$) [rectangle,fill=white,draw=white,inner sep=1pt] {\textcolor{labels}{$B_2$}};
	\node at ($(r)!0.5!(y) + (-0.35,-0.15)$) [rectangle,fill=white,draw=white,inner sep=1pt] {\textcolor{labels}{$B_3$}};
	\node at ($(y)!0.5!(t) + (-0.35,-0.15)$) [rectangle,fill=white,draw=white,inner sep=1pt] {\textcolor{labels}{$C_1$}};
	\node at ($(t)!0.5!(s) + (-0.35,0.15)$) [rectangle,fill=white,draw=white,inner sep=1pt] {\textcolor{labels}{$C_2$}};
	\node at ($(s)!0.5!(x) + (0,0.35)$) [rectangle,fill=white,draw=white,inner sep=1pt] {\textcolor{labels}{$C_3$}};
\end{tikzpicture}
\caption{A Standard $5$-wild Triangle}
\label{standard}
\end{figure}
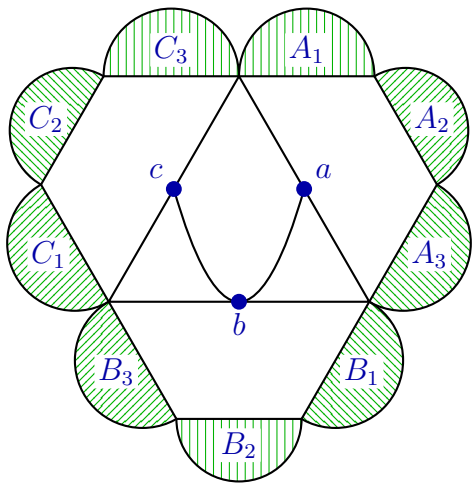

Let $M$ be a $k$-coherent matroid and $T=\{a,b,c\}$ 
be a triangle of
$M$. Then $T$ is $k$-{\em wild\ } 
\index{$k$-wild triangle}
if $M\ba t$ is 3-connected and $k$-fractured
for all $t$ in $T$.
If $T$ is a $k$-wild triangle, then a 
{\em $k$-wild display\ }
\index{$k$-wild display} 
for $T$ is a partition
$$(A_1,A_2,\ldots,A_{k-2},B_1,B_2,\ldots,B_{k-2},
C_1,C_2,\ldots,C_{k-2})$$
of $E(M)-T$ such that the following hold, 
where $A=A_1,A_2,\udots A_{k-2}$,
$B=B_1,B_2,\udots B_{k-2}$ and $C=C_1,C_2,\udots C_{k-2}$.
\begin{itemize}
\item[(i)] $(A_1,A_2,\ldots,A_{k-2},B\cup C\cup T)$, 
$(B_1,B_2,\ldots,B_{k-2},A\cup C\cup T)$,
and $(C_1,C_2,\ldots,C_{k-2},A\cup B\cup T)$ 
are tight swirl-like flowers of $M$.
\item[(ii)] $(A_1,A_2,\ldots,A_{k-2},B\cup\{b\},C\cup \{c\})$,
$(A\cup\{a\},B_1,B_2,\ldots,B_{k-2},C\cup\{c\})$, and 
$(A\cup \{a\},B\cup\{b\},C_1,C_2,\ldots,C_{k-2})$ are $k$-fractures of 
$M\ba a$, $M\ba b$ and $M\ba c$ respectively.
\end{itemize}
Moreover, $T$ is a {\em standard} $k$-wild triangle
\index{standard $k$-wild triangle} 
if it has
a $k$-wild display such that 
$(A\cup\{a\},B\cup\{b\},C\cup\{c\})$ is a swirl-like flower
of $M$.
Figure~\ref{standard} illustrates a standard $5$-wild triangle.

We now describe another type of $k$-wild triangle obtained 
from a $\Delta-Y$ exchange. We first note an elementary lemma.
We omit the easy proof.

\begin{lemma}
\label{delta-conn}
Let $T$ be a triangle of the matroid $M$; let
$\Delta M$ be the matroid obtained by performing a
$\Delta-Y$ exchange on $T$; and let $A$ be a set of elements of 
$M$. Then the following hold.
\begin{itemize}
\item[(i)] If $T\subseteq A$, then $r_{\Delta M}(A)=r_M(A)+1$.
\item[(ii)] If $T\cap A=\emptyset$, then 
$r_{\Delta M}(A)=r_M(A)$.
\item[(iii)] If $T\subseteq A$, then 
$\lambda_{\Delta M}(A)=\lambda_M(A)$.
\item[(iv)] If $t\in T$, then $M\ba t=\Delta M/t$.
\end{itemize}
\end{lemma}

\begin{lemma}
\label{delta-wild}
Let $T$ be a triangle of the 
matroid $M$ that is not contained in a $4$-element fan
and let $M'$ be the
matroid obtained by first performing a $\Delta-Y$ exchange
on $T$ and then taking the dual.
\begin{itemize}
\item[(i)] $M'$ is $3$-connected if and only if $M$ is.
\item[(ii)] Assume that $M$ is $3$-connected,
that $\PP=(P_1,P_2,\ldots,P_m)$ is a swirl-like
flower of $M$ of order at least $3$, 
and $T$ is contained in a petal of 
$\PP$. Then $\PP$ is tight in $M$ if and only if 
$\PP$ is tight in  $M'$.
\item[(iii)] $M'$ is $k$-coherent if and only if $M$ is.
\item[(iv)] Assume that $M$ is $k$-coherent. Then $T$ is 
$k$-wild in $M$ if and only if $T$ is $k$-wild in $M'$.
\end{itemize}
\end{lemma}

\begin{proof}
First note that, by definition, $M$ is obtained from
$M'$ by first taking the dual and then performing a 
$Y-\Delta$ exchange. Again by definition this means that we
first perform a $\Delta-Y$ exchange on $M'$ and then take the
dual. Thus the operation we are considering 
is an involution and we need only
prove the parts of the lemma in one direction.
The straightforward proof of (i) is given in 
\cite[Lemma~8.2]{wild} and we omit it here.

Consider (ii). Assume that $T=\{a,b,c\}$ and $T\subseteq P_2$. 
It is easily seen that $\PP$ is a swirl-like flower of 
$M'$. The only problem that could happen is that $P_2$
either ceases to be tight or becomes tight. 
(Indeed this could happen if $\PP$
were a copaddle, so the issue needs addressing.)
By the observation at the start of the proof we may assume
that $P_2$ is a loose petal of $\PP$ in $M$. 
Say $P_2\neq T$. Then $T$ belongs to a $4$-element
fan in $M'$ and performing a $\Delta-Y$ exchange on a triangle
in a 4-element fan destroys 3-connectivity. Thus 
 we may assume that $T=P_2$.
We now prove that $P_2$ is a loose petal of $\PP$
in $M'$.
Up to labels, we have
$a\in\cl_M(P_1)$, $b\in\cl^*_M(P_1\cup \{a\})$
and $c\in\cl_M(P_3)$. Note that, $\{a',b',c\}$ is a triangle
in $P_\Delta(M(K_4),M)$. Thus 
$r_{\Delta M}(E(M)-(P_1\cup \{c\}))=r_{M}(E(M)-P_1)$.
Evidently $r_{\Delta M}(P_1\cup\{c\})=r_M(P_1)+1$.
Hence $\lambda_{\Delta M}(P_1\cup \{c\})=2$
so that $c\in\cl^*_{\Delta M}(P_1)$. It follows easily
that $b\in\cl_{\Delta M}(P_1\cup\{c\})$. By symmetry
$a\in\cl^*_{\Delta M}(P_3)$ and part (ii) follows.

Consider (iii). Assume that $M$ is not $k$-coherent.
Let $(P_1,P_2,\ldots,P_m)$ be a tight flower that
$k$-fractures $M$. By Lemma~\ref{3-cross} we may assume
up to labels that
$T\subseteq P_1$. By (iii), $(P_1,P_2,\ldots,P_m)$
is also a tight flower of $M'$ so that $M'$ is not
$k$-coherent. 

Assume that $M$ is $k$-coherent and $T$ is $k$-wild. 
Then it has a 
$k$-wild display. By the earlier parts of this lemma,
$M'$ is also $k$-coherent and, indeed, the $k$-wild display
is also a $k$-wild display for $T$ in $M'$. Thus $T$ is 
$k$-wild in $M'$.
\end{proof}

The $k$-wild triangle $T$ of the $k$-coherent matroid $M$ is 
a {\em costandard} $k$-wild triangle
\index{costandard $k$-wild triangle} 
if $T$ is a standard wild triangle 
in the matroid $M'$ constructed in Lemma~\ref{delta-wild}. 
It does not seem that easy to produce a schematic diagram
for a costandard $k$-wild triangle that is at all insightful.

Another type of triangle that has no element that can
be deleted to preserve $k$-coherence is one that
is an internal triangle
in a fan. Note that such a triangle has two 
elements $a$ and $c$ such that
$\co(M\ba a)$ and $\co(M\ba c)$ are both $k$-coherent. 
We can now state the
main result of this section.

\begin{theorem}
\label{wild1}
Let $T$ be a triangle of the $k$-coherent matroid $M$.
Assume that $M\ba t$ is not $k$-coherent for all $t\in T$. 
Then
$T$ is either an internal triangle in a fan of $M$, 
a standard $k$-wild
triangle, or a costandard $k$-wild triangle.
\end{theorem}

To prepare for the proof of Theorem~\ref{wild1} 
we recall material from \cite{wild}. 
Evidently a $k$-wild triangle
is wild. If $T$ is a wild triangle of the 
$3$-connected matroid $M$,
then $\{a,b,c\}$ is a {\em standard} wild triangle if there is a partition 
$\PP=(P_1,P_2,\ldots,P_6)$ of $E(M)-\{a,b,c\}$ 
such that $|P_i|\geq 2$
for all $i$ and the following hold.
\begin{itemize}
\item[(i)] $M\ba a$, $M\ba b$, and $M\ba c$ are 3-connected, 
$M\ba a,b,c$ is connected and $\co(M\ba a,b,c)$ is $3$-connected.
\item[(ii)] $(P_1\cup P_2\cup\{a\},P_3\cup P_4\cup\{b\},
P_5\cup P_6\cup\{c\})$ is a flower in $M$.
\item[(iii)] $(P_2\cup P_3\cup P_4\cup\{b\},
P_5\cup P_6\cup P_1\cup\{c\})$,
$(P_4\cup P_5\cup P_6\cup\{c\},P_1\cup P_2\cup P_3\cup\{a\})$
and $(P_6\cup P_1\cup P_2\cup \{a\},P_3\cup P_4\cup P_5\cup \{b\})$
are $3$-separations exposed in $M$ by $a$, $b$, and $c$ respectively.
\end{itemize}
A partition $\PP$ satisfying these conditions is a partition 
{\em associated} to $\{a,b,c\}$.

If the wild triangle $T$ is obtained from a 
standard wild triangle in a
3-connected matroid by first doing a $\Delta-Y$ exchange and 
dualising in the manner of Lemma~\ref{delta-wild}, then $T$ is
a {\em costandard\ } wild triangle. 
One other type of wild triangle needs to be
described. 

Let $R$ be a 3-separating set $\{a,b,c,s,t,u,v\}$ in a $3$-connected
matroid $M$, where $\{a,b,c\}$ is a triangle. 
Then $R$ is a {\em trident}
\index{trident}
with wild triangle $\{a,b,c\}$ if $\{t,s,u,b\}$, $\{t,u,v,c\}$, and
$\{t,s,v,a\}$ are exposed quads in $M\ba a$, $M\ba b$ and $M\ba c$
respectively. The main theorem of \cite{wild} is

\begin{theorem}
\label{wild2}
Let $T$ be a wild triangle of a $3$-connected 
matroid $M$ with at least
twelve elements. Then $T$ is either a standard 
wild triangle, a costandard
wild triangle, a triangle in a trident of $M$, 
or an internal triangle of a fan of $M$.
\end{theorem}

We will also need the next theorem from \cite{wild}. 

\begin{theorem}
\label{wild3}
Let $\{a,b,c\}$ be a standard wild triangle in a $3$-connected matroid
$M$ where $|E(M)|\geq 12$ and let $(X_1,X_2)$, $(Y_1,Y_2)$, and $(Z_1,Z_2)$
be $3$-separations exposed by $a$, $b$ and $c$, respectively, with 
$a\in Y_2\cap Z_1$, $b\in Z_2\cap X_1$, and $c\in X_2\cap Y_1$. 
Then $(X_1,X_2)$, $(Y_1,Y_2)$ and $(Z_1,Z_2)$ can be replaced by equivalent
$3$-separations such that
$(X_2\cap Y_2,Z_1\cap X_1,Y_2\cap Z_2,X_1\cap Y_1,Z_2\cap X_2,Y_1\cap Z_1)$
is a partition associated to $\{a,b,c\}$.
\end{theorem}

\begin{lemma}
\label{wild4}
Let $T=\{a,b,c\}$ be a $k$-wild triangle 
of the $k$-coherent matroid $M$.
If $T$ is a standard wild triangle of $M$, then $T$ is a standard
$k$-wild triangle of $M$.
\end{lemma}

\begin{proof}
As $T$ is $k$-wild, there are $3$-separations $(X_1,X_2)$, $(Y_1,Y_2)$
and $(Z_1,Z_2)$ exposed by $a$, $b$ and $c$ respectively that are 
displayed in $k$-fractures of $M\ba a$, $M\ba b$ and $M\ba c$
respectively. We may replace these by any equivalent $3$-separation.
So by Theorem~\ref{wild3}, we may assume that  
$(X_2\cap Y_2,Z_1\cap X_1,Y_2\cap Z_2,X_1\cap 
Y_1,X_2\cap Z_2,Y_1\cap Z_1)$
is a partition associated to $\{a,b,c\}$.

Let $P_1=X_2\cap Y_2$, $P_2=X_1\cap Z_1$, $P_3=Y_2\cap Z_2$,
$P_4=X_1\cap Y_1$, $P_5=X_2\cap Z_2$ and $P_6=Y_1\cap Z_1$. Then,
by the definition of an associated partition, 
$(P_1\cup P_2\cup\{a\},P_3\cup P_4\cup\{b\},P_5\cup P_6\cup\{c\})$
is a flower in $M$, so that 
$(P_1\cup P_2,P_3\cup P_4\cup\{b\},P_5\cup P_6\cup\{c\})$ is a flower
in $M\ba a$. The $3$-separation 
$(X_1,X_2)=(P_2\cup P_3\cup P_4\cup\{b\},P_5\cup P_6\cup P_1\cup\{c\})$
crosses this flower so by Lemma~\ref{flower8.2}
$(P_1,P_2,P_3\cup P_4\cup\{b\},P_5\cup P_6\cup\{c\})$ is a flower in 
$M\ba a$. As $(X_1,X_2)$ is displayed in a $k$-fracture of 
$M\ba a$, this flower refines to a $k$-fracture 
${\bf A}$ of $M\ba a$.
This shows that the above flower is swirl-like. 

We now show that ${\bf A}$ is 
obtained by refining $P_1$ and $P_2$. 
Assume otherwise. Then we may assume
that there is a partition $(P',P'')$ of $P_3\cup P_4\cup \{b\}$
such that $(P_1,P_2,P',P'',P_5\cup P_6\cup\{c\})$ is a flower in 
$M\ba a$. Certainly $a\in\cl(P_1\cup P_2)$, so that
$(P_1\cup P_2\cup\{a\},P',P'',P_5\cup P_6\cup\{c\})$ is a flower in 
$M$. We may assume that $b$ is in $P''$. Then, as $b\in\cl(\{a,c\})$,
we see by Lemma~\ref{fine2} that, up to labels,
$(P_1\cup P_2\cup\{a\},P',P''-\{b\},P_5\cup P_6\cup\{b,c\})$ 
is a flower in $M$ and, indeed,
$(P_1\cup P_2,P',P''-\{b\},P_5\cup P_6\cup\{a,b,c\})$ is a flower in
$M$. This contradicts the fact that $a$ exposes the 3-separation
$(X_1,X_2)$.

We deduce that there is a refinement $(A_1,\ldots,A_m)$ of 
$P_1\cup P_2$ such that 
$(A_1,\ldots,A_m,P_3\cup P_4\cup\{b\},P_5\cup P_6\cup\{c\})$
is a $k$-fracture of $M\ba a$. As $M$ is $k$-coherent, and $a$ does not
block $A_1\udots A_m$, we see that $m=k-2$. The lemma follows by repeating the
above argument in $M\ba b$ and $M\ba c$.
\end{proof}

While triangles in tridents are a problem in the general case,
they cause no difficulties in the $k$-coherent case.

\begin{lemma}
\label{wild5}
Let $T$ be a triangle in a trident of the $k$-coherent matroid $M$.
Then $T$ is not $k$-wild.
\end{lemma}

\begin{proof}
Assume that $\{a,b,c\}$ is in the trident
$X=\{a,b,c,t,s,u,v\}$ where the labelling accords with that 
given in the
definition of a trident. Assume that $\{a,b,c\}$ is $k$-wild.
Let $\AAA$ be a maximal $k$-fracture of $M\ba a$. Certainly
$X-\{a\}$ is not contained in a petal of $\AAA$, as otherwise $M$ is 
$k$-fractured. It is now readily verified that the quad
$\{t,s,u,b\}$ of $M\ba a$ is displayed in $\AAA$. 
This shows that, for some $m\geq 2$,
there is a partition $(P_1,\ldots,P_m)$ of $E(M)$ such that
$(X,P_1,\ldots,P_m)$ is a maximal swirl-like flower of $M$, and 
$(\{c,v\},\{t,s,u,b\},P_1,\ldots,P_m)$ 
is a swirl-like flower in $M\ba a$.
(It is conceivable that this flower 
is not maximal in that $\{t,s,u,b\}$
could refine to a pair of petals. 
In fact this does not happen, but we
don't need to establish this fact.) 
We now know that
$$\sqcap(\{c,v\},P_m)=1~~\text{and}~~\sqcap(\{s,t,u,b\},P_m)=0.$$
Repeating the above argument for $M\ba c$, and using 
Lemma~\ref{one-flower} we deduce that either 
$(\{b,u\},\{s,t,a,v\},P_1,\ldots,P_m)$ or 
$(\{s,t,a,v\},\{b,u\},P_1,\ldots,P_m)$ is a swirl-like flower in
$M\ba c$. As $\sqcap(\{s,t,u,b\},P_m)=0$, it must be the former,
and we deduce that $\sqcap(\{s,t,a,v\},P_m)=0$ so that
$$\sqcap(\{a,s\},P_m)=0.$$
Repeating for $M\ba b$ and using the above fact, we deduce that
$(\{a,s\},\{v,u,t,c\},P_1,\ldots,P_m)$ is a swirl-like flower
in $M\ba b$. But this means that $\sqcap(\{v,u,t,c\},P_m)=0$
so that $\sqcap(\{c,v\},P_m)=0$ contradicting the fact that
$\sqcap(\{c,v\},P_m)=1$.
\end{proof}

The proof of Theorem~\ref{wild1} is now just a matter of summing up.

\begin{proof}[Proof of Theorem~\ref{wild1}]
Let $T$ be a $k$-wild triangle of $M$. Then certainly $T$ is wild.
By Lemma~\ref{wild5}, $T$ is not a triangle in a trident. If $T$ is 
a standard wild triangle of $M$, then $T$ is a standard $k$-wild
triangle of $M$ by Lemma~\ref{wild4}. Assume that $T$ is a costandard
wild triangle. Let $M'$ be the matroid obtained by doing a 
$\Delta-Y$ exchange on $T$ and then taking the dual. By 
Lemma~\ref{delta-wild}, $T$ is a $k$-wild triangle of $M'$. By
definition $T$ is a standard wild triangle of $M'$ so that $T$
is a standard $k$-wild triangle of $M'$. Now, by the definition of 
a costandard $k$-wild triangle we deduce that $T$ is indeed a 
costandard $k$-wild triangle of $M$. The theorem now follows
from Theorem~\ref{wild2}.
\end{proof}

If $T$ is a 
triad of the $k$-coherent matroid $M$, then $T$ is 
$k$-{\em wild\ }
\index{$k$-wild triad} 
if $T$ is a $k$-wild triangle of $M^*$. If $T$
is a $k$-wild triad, then $T$ is a {\em standard\ }
\index{standard $k$-wild triad} 
(respectively {\em costandard\ }) $k$-wild triad
\index{costandard $k$-wild triad} 
if $T$ is a standard (respectively
costandard) $k$-wild triangle of $M^*$. A partition of
$E(M)$ is a {\em $k$-wild display\ } for the triad $T$ of $M$
if it is a $k$-wild display for $T$ in $M^*$.

We conclude this section by
giving some elementary properties of $k$-wild 
triangles. We will use the following lemma on
$\Delta-Y$ exchanges.

\begin{lemma}
\label{d-y}
Let $(X,Y)$ be a $3$-separation of the $3$-connected matroid
$M$, let $\{a,b,c\}$ be a triangle of $M$ 
contained in $Y$ and let $N$ be the matroid obtained by
performing a $\Delta-Y$ exchange on
$\{a,b,c\}$. Assume that $M\ba a$, $M\ba b$, and $M\ba c$ are 
$3$-connected and that $\{a,b,c\}\cap \cl_M(X)\neq \emptyset$.
Then $\clstar_N(X)\cap\{a,b,c\}=\emptyset$.
\end{lemma}

\begin{proof}
Assume that $a\in\cl_M(X)$. 
By Lemma~\ref{delta-conn}(ii),  $r_N(X)=r_M(X)$, $r_N(Y)=r_M(Y)+1$,
and $r(N)=r(M)+1$. Thus, $\lambda_N(X)=2$.
Assume that $\clstar_N(X)$ contains the element $t$ in $\{a,b,c\}$.
As $\{a,b,c\}$ is a cocircuit of $N$, $t\notin\cl_N(X)$.
Hence $t\in\cl^*_N(X)$, so that $\lambda_{N\ba t}(X)=1$. 
Say $t=a$. Then $\lambda_{N\ba a/b}(X)=1$. 
By Lemma~\ref{delta-conn}(iv),
$N\ba a/b\cong M\ba a,b$. So $\lambda_{M\ba a,b}(X)=1$. 
But $a\in\cl(X)$
so that $\lambda_{M\ba b}(X)=1$ contradicting the fact that $M\ba b$
is $3$-connected. If $b\in\cl^*_N(X)$, then we see that
$\lambda_{N\ba b/a}(X)=1$. But $N\ba b/a\cong M\ba a,b$ and
we arrive at the same contradiction.
\end{proof}

The next lemma highlights subtle differences in 
the behaviour of standard and costandard $k$-wild triangles.
To avoid a cumbersome statement we omit obvious symmetric
statements. 

\begin{lemma}
\label{k-wild-different}
Let $\{a,b,c\}$ be a $k$-wild triangle of the
$k$-coherent matroid $M$ with $k$-wild display
$(A_1,A_2,\ldots,A_{k-2},B_1,B_2,\ldots,B_{k-2},C_1,C_2,
\ldots,C_{k-2})$.
Let $A=A_1\cup A_2\udots A_{k-2}$. 
Then the following hold.
\begin{itemize}
\item[(i)] If $T$ is standard, then $a\in\cl(A)$.
\item[(ii)] If $T$ is costandard, then $a\notin \cl(A)$.
\item[(iii)] If $T$ is standard, then
$\si(M/a)$ is not $3$-connected.
\item[(iv)] If $T$ is costandard, then $\si(M/a)$ is 
$3$-connected.
\end{itemize}
\end{lemma}

\begin{proof} Part (i) is clear.
Part  (ii) follows from Lemma~\ref{d-y}. 
As $a\in\cl(A)$ we see that $a$ is in the guts of a vertical
$3$-separation of $M$ and (iii) follows.
While routine to prove, we simply note here that
(iv) follows from \cite[Corollary~3.3(iii)]{wild}.
\end{proof}

The fact that any element of a standard $k$-wild triangle is on the
guts of a vertical $3$-separation is helpful for proving that
certain triangles are not standard $k$-wild. The next lemma
is useful for certifying that a triangle is not a 
costandard $k$-wild triangle. For ease of proof we state
the dual form.

\begin{lemma}
\label{cowild-win}
Let $\{a,b,c\}$ be a costandard $k$-wild triad of the 
$k$-coherent matroid $M$. Then $\si(M/a,b)$ is not 
$3$-connected.
\end{lemma}

\begin{proof}
Let $M'$ be the matroid obtained by performing a 
$Y-\Delta$ exchange on $M$. Observe that $T$ is a standard
$k$-wild triangle of $M$. By Lemma~\ref{delta-conn}(iv)
$M/a,b\cong M'\ba a/b$ and it is 
readily verified that $\si(M'\ba a/b)$ is not $3$-connected.
\end{proof}

This next lemma gives some elementary properties that are common
to both types of $k$-wild triangle.

\begin{lemma}
\label{k-wild-common}
Let $\{a,b,c\}$ be a $k$-wild triangle of the
$k$-coherent matroid $M$ with $k$-wild display
$(A_1,A_2,\ldots,A_{k-2},B_1,B_2,\ldots,B_{k-2},
C_1,C_2,\ldots,C_{k-2})$.
Let $A=A_1\cup A_2\udots A_{k-2}$. 
Then the following hold.
\begin{itemize}
\item[(i)] $\lambda(A\cup\{a,b\})>2$
and $\lambda(A\cup\{a,b,c\})>2$.
\item[(ii)] If $i\in\{1,2,\ldots,k-2\}$, then neither 
$a$ nor $b$ is in the full
closure of $A_i$.
\end{itemize}
\end{lemma}

\begin{proof}
Consider (i). Assume that $T$ is standard.
If  $\lambda(A\cup\{a,b,c\})=2$, then,
as $c\in\cl(C_1\cup C_2\udots C_{k-2})$, 
we have $\lambda(A\cup\{a,b\})=2$ so it
suffices to prove that $\lambda(A\cup\{a,b\})>2$.
Assume otherwise. Then $b\in\clstar(A\cup\{a\})$.
But then
$(A\cup\{a,b\},B_1,B_2,\ldots,B_{k-2},C_1,C_2,\udots C_{k-2}\cup\{c\})$
is a $k$-fracture of $M$ contradicting the fact that
$M$ is $k$-coherent.  Thus (i) holds in the case that
$T$ is standard.

Assume that $T$ is costandard. It follows from
(i) and Lemma~\ref{delta-conn} that $\lambda(A\cup\{a,b,c\})>2$.
If $\lambda(A\cup\{a,b\})=2$, then, as $c\in\cl(\{a,b\})$
we have $\lambda(A\cup\{a,b,c\})=2$. Thus (i) also holds
in the case that $T$ is costandard.

If (ii) fails then we easily see that either $a$ or $b$ is in the
full closure of a petal of a $k$-fracture of $M\ba a$
or $M\ba b$ respectively. This contradicts the fact that
$M$ is $k$-coherent.
\end{proof}

\section{Feral Elements}
\label{feral-elements}

Throughout 
this section we assume that $M$ is a $k$-coherent matroid.
An element $f$ of $M$ is {\em feral\ }
\index{feral element} 
if both
$M\ba f$ and $M/f$ are $3$-connected and $k$-fractured. 
The goal of this section is to gain insight into the structure
of a matroid relative to a feral element. 

Let $(P_1,P_2,\ldots,P_m)$ and $(Q_1,Q_2,\ldots,Q_k)$ be partitions of
$E(M)-\{f\}$. Then these partitions form a {\em feral display\ }
\index{feral display} 
for
$f$ if there is an $i\in\{2,3,\ldots,m-1\}$ 
such that the following hold.
\begin{itemize}
\item[(i)] $(P_1,P_2,\ldots,P_m)$ and $(Q_1,Q_2,\ldots,Q_k)$ are 
$k$-fractures of $M\ba f$ and $M/f$ respectively.
\item[(ii)] $\{P_2,P_3,\ldots,P_m,Q_3,Q_4,\ldots,Q_k,Z_1=Q_1\cap P_1,
Z_2=Q_2\cap P_1\}$ partitions $E(N)-\{f\}$ into nonempty sets, 
with the exception that one of $Z_1$ or $Z_2$ may be empty.
\item[(iii)] $P_1=Q_3\cup Q_4\udots Q_k\cup Z_1\cup Z_2$.
\item[(iv)] $Q_1=P_{i+1}\cup P_{i+2}\udots P_m\cup Z_1$.
\item[(v)] $Q_2=Z_2\cup P_2\cup P_3\udots P_i$.
\item[(vi)] $(Q_1\cup Q_2\cup\{f\},Q_3,\dots,Q_{k})$ is a 
swirl-like flower of order $(k-1)$ in $M$.
\item[(vii)] $(P_2,P_3,\ldots,P_i,P_{i+1}\cup P_{i+2}
\udots P_m\cup P_1\cup\{f\})$ 
is a swirl-like flower of order $i$ in $M$.
\item[(viii)] $(P_{i+1},P_{i+2},\ldots,P_m,P_1\udots P_i\cup\{f\})$ 
is a swirl-like flower of order $(m-i+1)$ in $M$.
\item[(ix)] Either $Z_1\neq\emptyset$ or 
$Z_2\neq \emptyset$ and 
$\lambda_M(Z_1),\lambda_M(Z_2)\leq 3$.
\item[(x)] $f$ blocks $P_1$ and $f$ coblocks $Q_1$.
\end{itemize}

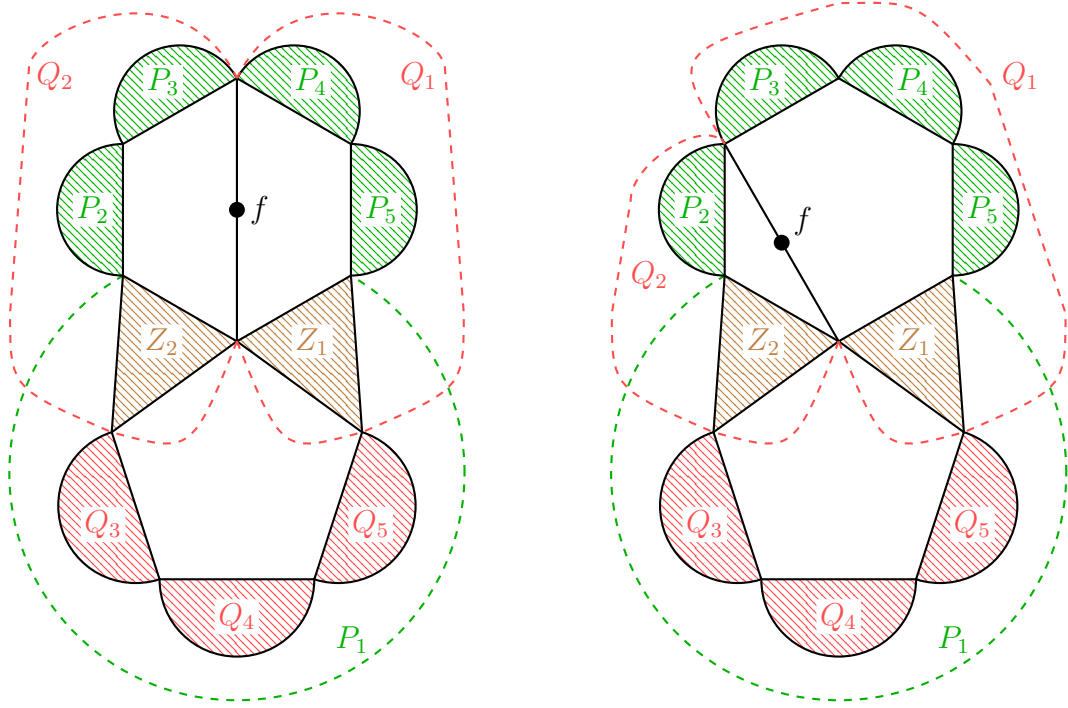
\begin{figure}
\begin{tikzpicture}[thick,line join=round]
	\coordinate (Ao) at (-4,7);
	\coordinate (Bo) at (4,7);
	\coordinate (Af) at ($(Ao) + (0,1.75)$);
	\coordinate (Ap) at ($(Ao) + (0,-1.75)$);
	\coordinate (Aa) at ($(Af) + (30:1.75)$);
	\coordinate (Ab) at ($(Af) + (90:1.75)$);
	\coordinate (Ac) at ($(Af) + (150:1.75)$);
	\coordinate (Ad) at ($(Af) + (210:1.75)$);
	\coordinate (Ae) at ($(Af) + (-30:1.75)$);
	\coordinate (Ag) at ($(Ap) + (18:1.75)$);
	\coordinate (Ah) at ($(Ap) + (162:1.75)$);
	\coordinate (Ai) at ($(Ap) + (234:1.75)$);
	\coordinate (Aj) at ($(Ap) + (306:1.75)$);
	\coordinate (Acd) at ($(Ac)!0.5!(Ad)$);
	\coordinate (Acb) at ($(Ac)!0.5!(Ab)$);
	\coordinate (Aba) at ($(Ab)!0.5!(Aa)$);
	\coordinate (Aae) at ($(Aa)!0.5!(Ae)$);
	\coordinate (Ahi) at ($(Ah)!0.5!(Ai)$);
	\coordinate (Aij) at ($(Ai)!0.5!(Aj)$);
	\coordinate (Ajg) at ($(Aj)!0.5!(Ag)$);
	\node[draw,circle through=(Ae),dashed,glines] at (Ap) {};
	\node[pattern color=lines,draw,circle through=(Ac),pattern=north west lines] at (Acd) {};
	\node[pattern color=lines,draw,circle through=(Ab),pattern=north west lines] at (Acb) {};
	\node[pattern color=lines,draw,circle through=(Aa),pattern=north west lines] at (Aba) {};
	\node[pattern color=lines,draw,circle through=(Ae),pattern=north west lines] at (Aae) {};
	\fill[white] (Ad) -- (Ac) -- (Ab) -- (Aa) -- (Ae) -- (Ao) -- cycle;
	\fill[pattern color=olines,pattern=north west lines] (Ad) -- (Ao) -- (Ag) -- (Ae) -- (Ao) -- (Ah) -- cycle;
	\node[pattern color=rlines,draw,circle through=(Ah),pattern=north west lines] at (Ahi) {};
	\node[pattern color=rlines,draw,circle through=(Ai),pattern=north west lines] at (Aij) {};
	\node[pattern color=rlines,draw,circle through=(Aj),pattern=north west lines] at (Ajg) {};
	\fill[white] (Ao) -- (Ah) -- (Ai) -- (Aj) -- (Ag) -- cycle;
	\draw (Ao) -- (Ad) -- (Ac) -- (Ab) -- (Aa) -- (Ae) -- (Ao) -- (Ah) -- (Ai) -- (Aj) -- (Ag) -- cycle;
	\draw (Ad) -- (Ah);
	\draw (Ae) -- (Ag);
	\draw (Ao) -- (Ab);
	\node at ($(Acd) + (180:0.4)$) [rectangle,fill=white,draw=white,inner sep=1] {\textcolor{glabels}{$P_2$}};
	\node at ($(Acb) + (120:0.45)$) [rectangle,fill=white,draw=white,inner sep=1] {\textcolor{glabels}{$P_3$}};
	\node at ($(Aba) + (60:0.45)$) [rectangle,fill=white,draw=white,inner sep=1] {\textcolor{glabels}{$P_4$}};
	\node at ($(Aae) + (0:0.4)$) [rectangle,fill=white,draw=white,inner sep=1] {\textcolor{glabels}{$P_5$}};
	\node at ($(Ao) + (1,0)$) [rectangle,fill=white,draw=white,inner sep=1] {\textcolor{olabels}{$Z_1$}};
	\node at ($(Ao) - (1,0)$) [rectangle,fill=white,draw=white,inner sep=1] {\textcolor{olabels}{$Z_2$}};
	\node at ($(Ahi) + (-150:0.5)$) [rectangle,fill=white,draw=white,inner sep=1] {\textcolor{rlabels}{$Q_3$}};
	\node at ($(Aij) + (-90:0.5)$) [rectangle,fill=white,draw=white,inner sep=1] {\textcolor{rlabels}{$Q_4$}};
	\node at ($(Ajg) + (-30:0.5)$) [rectangle,fill=white,draw=white,inner sep=1] {\textcolor{rlabels}{$Q_5$}};
	\node at ($(Af) + (0.3,0)$) {$f$};
	\fill[black] (Af) circle (3pt);
	\draw[dashed,rlines] (Ao) .. controls ($(Ao) + (0.6,-1.5)$) .. (Ag) .. controls ($(Ae) + (1.5,-1.5)$) .. ($(Ae) + (1.5,-0.5)$) -- ($(Aa) + (1.25,1)$) .. controls ($(Aa) + (1.25,1.5)$) and ($(Aba) + (0,2)$) .. (Ab);
	\draw[dashed,rlines] (Ao) .. controls ($(Ao) + (-0.6,-1.5)$) .. (Ah) .. controls ($(Ad) + (-1.5,-1.5)$) .. ($(Ad) + (-1.5,-0.5)$) -- ($(Ac) + (-1.25,1)$) .. controls ($(Ac) + (-1.25,1.5)$) and ($(Acb) + (0,2)$) .. (Ab);
	\node at ($(Aa) + (0.9,0.9)$) {\textcolor{rlabels}{$Q_1$}};
	\node at ($(Ac) + (-0.9,0.9)$) {\textcolor{rlabels}{$Q_2$}};
	\node at ($(Aj) + (0.5,-0.8)$) {\textcolor{glabels}{$P_1$}};
	\coordinate (Bq) at ($(Bo) + (0,1.75)$);
	\coordinate (Bp) at ($(Bo) + (0,-1.75)$);
	\coordinate (Ba) at ($(Bq) + (30:1.75)$);
	\coordinate (Bb) at ($(Bq) + (90:1.75)$);
	\coordinate (Bc) at ($(Bq) + (150:1.75)$);
	\coordinate (Bd) at ($(Bq) + (210:1.75)$);
	\coordinate (Be) at ($(Bq) + (-30:1.75)$);
	\coordinate (Bg) at ($(Bp) + (18:1.75)$);
	\coordinate (Bh) at ($(Bp) + (162:1.75)$);
	\coordinate (Bi) at ($(Bp) + (234:1.75)$);
	\coordinate (Bj) at ($(Bp) + (306:1.75)$);
	\coordinate (Bf) at ($(Bo)!0.5!(Bc)$);
	\coordinate (Bcd) at ($(Bc)!0.5!(Bd)$);
	\coordinate (Bcb) at ($(Bb)!0.5!(Bc)$);
	\coordinate (Bba) at ($(Ba)!0.5!(Bb)$);
	\coordinate (Bae) at ($(Be)!0.5!(Ba)$);
	\coordinate (Bhi) at ($(Bh)!0.5!(Bi)$);
	\coordinate (Bij) at ($(Bi)!0.5!(Bj)$);
	\coordinate (Bjg) at ($(Bj)!0.5!(Bg)$);
	\node[draw,circle through=(Be),dashed,glines] at (Bp) {};
	\node[pattern color=lines,draw,circle through=(Bc),pattern=north west lines] at (Bcd) {};
	\node[pattern color=lines,draw,circle through=(Bb),pattern=north west lines] at (Bcb) {};
	\node[pattern color=lines,draw,circle through=(Ba),pattern=north west lines] at (Bba) {};
	\node[pattern color=lines,draw,circle through=(Be),pattern=north west lines] at (Bae) {};
	\fill[white] (Bd) -- (Bc) -- (Bb) -- (Ba) -- (Be) -- (Bo) -- cycle;
	\fill[pattern color=olines,pattern=north west lines] (Bd) -- (Bo) -- (Bg) -- (Be) -- (Bo) -- (Bh) -- cycle;
	\node[pattern color=rlines,draw,circle through=(Bh),pattern=north west lines] at (Bhi) {};
	\node[pattern color=rlines,draw,circle through=(Bi),pattern=north west lines] at (Bij) {};
	\node[pattern color=rlines,draw,circle through=(Bj),pattern=north west lines] at (Bjg) {};
	\fill[white] (Bo) -- (Bh) -- (Bi) -- (Bj) -- (Bg) -- cycle;
	\draw (Bo) -- (Bd) -- (Bc) -- (Bb) -- (Ba) -- (Be) -- (Bo) -- (Bh) -- (Bi) -- (Bj) -- (Bg) -- cycle;
	\draw (Bd) -- (Bh);
	\draw (Be) -- (Bg);
	\draw (Bo) -- (Bc);
	\node at ($(Bcd) + (180:0.4)$) [rectangle,fill=white,draw=white,inner sep=1] {\textcolor{glabels}{$P_2$}};
	\node at ($(Bcb) + (120:0.45)$) [rectangle,fill=white,draw=white,inner sep=1] {\textcolor{glabels}{$P_3$}};
	\node at ($(Bba) + (60:0.45)$) [rectangle,fill=white,draw=white,inner sep=1] {\textcolor{glabels}{$P_4$}};
	\node at ($(Bae) + (0:0.4)$) [rectangle,fill=white,draw=white,inner sep=1] {\textcolor{glabels}{$P_5$}};
	\node at ($(Bo) + (1,0)$) [rectangle,fill=white,draw=white,inner sep=1] {\textcolor{olabels}{$Z_1$}};
	\node at ($(Bo) - (1,0)$) [rectangle,fill=white,draw=white,inner sep=1] {\textcolor{olabels}{$Z_2$}};
	\node at ($(Bhi) + (-150:0.5)$) [rectangle,fill=white,draw=white,inner sep=1] {\textcolor{rlabels}{$Q_3$}};
	\node at ($(Bij) + (-90:0.5)$) [rectangle,fill=white,draw=white,inner sep=1] {\textcolor{rlabels}{$Q_4$}};
	\node at ($(Bjg) + (-30:0.5)$) [rectangle,fill=white,draw=white,inner sep=1] {\textcolor{rlabels}{$Q_5$}};
	\node at ($(Bf) + (45:0.4)$) {$f$};
	\fill[black] (Bf) circle (3pt);
	\draw[dashed,rlines] (Bo) .. controls ($(Bo) + (-0.6,-1.5)$) .. (Bh) .. controls ($(Bd) + (-1.5,-1.5)$) .. ($(Bd) + (-1.5,-0.5)$) -- ($(Bd) + (-1.25,1)$) .. controls ($(Bd) + (-1.25,1.5)$) and ($(Bcd) + (-0.25,1.25)$) .. (Bc);
	\draw[dashed,rlines] (Bo) .. controls ($(Bo) + (0.6,-1.5)$) .. (Bg) .. controls ($(Be) + (1.5,-1.5)$) .. ($(Be) + (1.5,-0.5)$) -- ($(Ba) + (0.5,0.8)$) .. controls ($(Ba) + (-0.5,2)$) and ($(Bb) + (1,1)$) .. ($(Bb) + (0,1)$) .. controls ($(Bcb) + (-1.5,0.8)$) .. (Bc);
	\node at ($(Ba) + (0.9,0.9)$) {\textcolor{rlabels}{$Q_1$}};
	\node at ($(Bd) + (-1,0)$) {\textcolor{rlabels}{$Q_2$}};
	\node at ($(Bj) + (0.5,-0.8)$) {\textcolor{glabels}{$P_1$}};
\end{tikzpicture}
\caption{Some Feral Elements}
\label{some-feral-elements}
\end{figure}

\begin{figure}
\begin{tikzpicture}[thick,line join=round]
	\coordinate (Co) at (-3,2);
	\coordinate (Do) at (3,2);
	\coordinate (Cq) at ($(Co) + (0,1.75)$);
	\coordinate (Cp) at ($(Co) + (0,-1.75)$);
	\coordinate (Ca) at ($(Cq) + (30:1.75)$);
	\coordinate (Cb) at ($(Cq) + (90:1.75)$);
	\coordinate (Cc) at ($(Cq) + (150:1.75)$);
	\coordinate (Cd) at ($(Cq) + (210:1.75)$);
	\coordinate (Ce) at ($(Cq) + (-30:1.75)$);
	\coordinate (Cg) at ($(Cp) + (18:1.75)$);
	\coordinate (Ch) at ($(Cp) + (162:1.75)$);
	\coordinate (Ci) at ($(Cp) + (234:1.75)$);
	\coordinate (Cj) at ($(Cp) + (306:1.75)$);
	\coordinate (Cf) at ($(Co)!0.5!(Cc)$);
	\coordinate (Cx) at ($(Cd)!0.33!(Cc)$);
	\coordinate (Cy) at ($(Cd)!0.67!(Cc)$);
	\coordinate (Ccd) at ($(Cc)!0.5!(Cd)$);
	\coordinate (Ccb) at ($(Cb)!0.5!(Cc)$);
	\coordinate (Cba) at ($(Ca)!0.5!(Cb)$);
	\coordinate (Cae) at ($(Ce)!0.5!(Ca)$);
	\coordinate (Chi) at ($(Ch)!0.5!(Ci)$);
	\coordinate (Cij) at ($(Ci)!0.5!(Cj)$);
	\coordinate (Cjg) at ($(Cj)!0.5!(Cg)$);
	\node[pattern color=lines,draw,circle through=(Cb),pattern=north west lines] at (Ccb) {};
	\node[pattern color=lines,draw,circle through=(Ca),pattern=north west lines] at (Cba) {};
	\node[pattern color=lines,draw,circle through=(Ce),pattern=north west lines] at (Cae) {};
	\fill[white] (Cd) -- (Cc) -- (Cb) -- (Ca) -- (Ce) -- (Co) -- cycle;
	\fill[pattern color=olines,pattern=north west lines] (Cd) -- (Co) -- (Cg) -- (Ce) -- (Co) -- (Ch) -- cycle;
	\node[pattern color=rlines,draw,circle through=(Ch),pattern=north west lines] at (Chi) {};
	\node[pattern color=rlines,draw,circle through=(Ci),pattern=north west lines] at (Cij) {};
	\node[pattern color=rlines,draw,circle through=(Cj),pattern=north west lines] at (Cjg) {};
	\fill[white] (Co) -- (Ch) -- (Ci) -- (Cj) -- (Cg) -- cycle;
	\draw (Co) -- (Cd) -- (Cc) -- (Cb) -- (Ca) -- (Ce) -- (Co) -- (Ch) -- (Ci) -- (Cj) -- (Cg) -- cycle;
	\draw (Cd) -- (Ch);
	\draw (Ce) -- (Cg);
	\draw (Co) -- (Cc);
	\node at ($(Ccd) + (180:0.4)$) [rectangle,fill=white,draw=white,inner sep=1] {\textcolor{glabels}{$P_2$}};
	\node at ($(Ccb) + (120:0.45)$) [rectangle,fill=white,draw=white,inner sep=1] {\textcolor{glabels}{$P_3$}};
	\node at ($(Cba) + (60:0.45)$) [rectangle,fill=white,draw=white,inner sep=1] {\textcolor{glabels}{$P_4$}};
	\node at ($(Cae) + (0:0.4)$) [rectangle,fill=white,draw=white,inner sep=1] {\textcolor{glabels}{$P_5$}};
	\node at ($(Co) + (1,0)$) [rectangle,fill=white,draw=white,inner sep=1] {\textcolor{olabels}{$Z_1$}};
	\node at ($(Co) - (1,0)$) [rectangle,fill=white,draw=white,inner sep=1] {\textcolor{olabels}{$Z_2$}};
	\node at ($(Chi) + (-150:0.5)$) [rectangle,fill=white,draw=white,inner sep=1] {\textcolor{rlabels}{$Q_3$}};
	\node at ($(Cij) + (-90:0.5)$) [rectangle,fill=white,draw=white,inner sep=1] {\textcolor{rlabels}{$Q_4$}};
	\node at ($(Cjg) + (-30:0.5)$) [rectangle,fill=white,draw=white,inner sep=1] {\textcolor{rlabels}{$Q_5$}};
	\node at ($(Cf) + (45:0.4)$) {$f$};
	\foreach \pt in {Cf,Cx,Cy} \fill[black] (\pt) circle (3pt);
	\coordinate (Da) at ($(Do) + (22.5:2.5)$);
	\coordinate (Db) at ($(Do) + (67.5:2.5)$);
	\coordinate (Dc) at ($(Do) + (112.5:2.5)$);
	\coordinate (Dd) at ($(Do) + (157.5:2.5)$);
	\coordinate (De) at ($(Do) + (202.5:2.5)$);
	\coordinate (Dh) at ($(Do) + (247.5:2.5)$);
	\coordinate (Di) at ($(Do) + (292.5:2.5)$);
	\coordinate (Dj) at ($(Do) + (337.5:2.5)$);
	\coordinate (Dab) at ($(Da)!0.5!(Db)$);
	\coordinate (Dbc) at ($(Db)!0.5!(Dc)$);
	\coordinate (Dcd) at ($(Dc)!0.5!(Dd)$);
	\coordinate (Deh) at ($(De)!0.5!(Dh)$);
	\coordinate (Dhi) at ($(Dh)!0.5!(Di)$);
	\coordinate (Dij) at ($(Di)!0.5!(Dj)$);
	\coordinate (Dja) at ($(Dj)!0.5!(Da)$);
	\node[pattern color=glines,pattern=north west lines,draw,circle through=(Da)] at (Dja) {};
	\node[pattern color=glines,pattern=north west lines,draw,circle through=(Db)] at (Dab) {};
	\node[pattern color=glines,pattern=north west lines,draw,circle through=(Dc)] at (Dbc) {};
	\node[pattern color=glines,pattern=north west lines,draw,circle through=(Dd)] at (Dcd) {};
	\node[pattern color=rlines,pattern=north west lines,draw,circle through=(Dh)] at (Deh) {};
	\node[pattern color=rlines,pattern=north west lines,draw,circle through=(Di)] at (Dhi) {};
	\node[pattern color=rlines,pattern=north west lines,draw,circle through=(Dj)] at (Dij) {};
	\filldraw[fill=white] (Da) -- (Db) -- (Dc) -- (Dd) -- (De) -- (Dh) -- (Di) -- (Dj) -- cycle;
	\draw (De) -- (Da);
	\draw (Dd) -- (Dj);
	\fill[pattern=horizontal lines,pattern color=olines] (Do) -- (De) -- (Dd) -- cycle;
	\node at ($(Do) - (1.8,0)$) [rectangle,fill=white,draw=white,inner sep=1] {\textcolor{olabels}{$Z$}};
	\node at ($(Dja) + (0:0.5)$) [rectangle,fill=white,draw=white,inner sep=1] {\textcolor{glabels}{$P_2$}};
	\node at ($(Dab) + (45:0.5)$) [rectangle,fill=white,draw=white,inner sep=1] {\textcolor{glabels}{$P_3$}};
	\node at ($(Dbc) + (90:0.5)$) [rectangle,fill=white,draw=white,inner sep=1] {\textcolor{glabels}{$P_4$}};
	\node at ($(Dcd) + (135:0.5)$) [rectangle,fill=white,draw=white,inner sep=1] {\textcolor{glabels}{$P_5$}};
	\node at ($(Deh) + (-135:0.5)$) [rectangle,fill=white,draw=white,inner sep=1] {\textcolor{rlabels}{$Q_3$}};
	\node at ($(Dhi) + (270:0.5)$) [rectangle,fill=white,draw=white,inner sep=1] {\textcolor{rlabels}{$Q_4$}};
	\node at ($(Dij) + (-45:0.5)$) [rectangle,fill=white,draw=white,inner sep=1] {\textcolor{rlabels}{$Q_5$}};
	\coordinate[label=135:$f$] (Df) at ($(Da)!0.5!(Do)$);
	\coordinate[label=-135:$f'$] (Df') at ($(Dj)!0.5!(Do)$);
	\foreach \pt in {Df,Df'} \fill[black] (\pt) circle (3pt);
\end{tikzpicture}
\caption{Some More Feral Elements}
\label{some-more-feral-elements}
\end{figure}
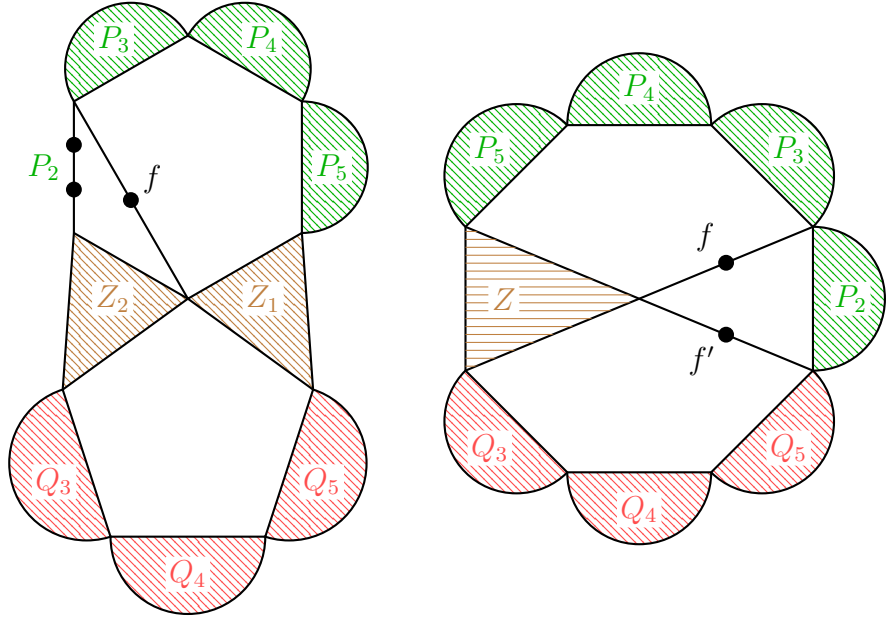

Figures~\ref{some-feral-elements} and 
\ref{some-more-feral-elements} illustrate 
some different cases of feral displays
associated with a feral element $f$. In each of the cases of this 
figure, both $Z_1$ and $Z_2$ are nonempty. But the case that 
one of these sets is empty does arise. Such a case is illustrated by 
``bogan couples'' which are defined and discussed in 
Chapter~\ref{k-skeletons}.

The primary goal is to prove that feral elements are characterised
by the existence of a feral display. In other words we prove

\begin{theorem}
\label{feral}
Let $f$ be a feral element of the matroid $M$ and let
$\FF_d$ and $\FF_c$ be blooms that fracture $M\ba f$
and $M/f$ respectively. Then
there are flowers displayed by $\FF_d$ and $\FF_c$ that form 
a feral display for $f$ in at least one of $M$ or $M^*$.
\end{theorem}

In what follows we assume that 
$f$ is indeed a feral element of $M$ and that 
$\FF_d=(\hP_1,\ldots,\hP_m)$ and $\FF_c=(\hQ_1,\ldots,\hQ_n)$
are maximal blooms that fracture $M\ba f$ and $M/f$ respectively.
We say that $f$ is 1-{\em blocking\ }
\index{$1$-blocking} 
for $\FF_d$ if there is a 
flower $(P_1,P_2,\ldots,P_m)$ displayed by $\FF_d$ such that, 
$f$ blocks $P_i$ and no other petal for some
$i\in\{1,2,\ldots,m\}$. Moreover, $f$ is 
2-{\em spanned\ }
\index{$2$-spanned} 
by $\FF_d$ if  
$f\in\cl(\hP_i\cup\hP_{i+1})$
for some
$i\in\{1,2,\ldots,m\}$. Dually, we say that $f$
is 1-{\em coblocking\ }
\index{$1$-coblocking} 
for $\FF_c$ if $f$ is 1-blocking
for $\FF_c$ in the matroid $M^*$, and $f$ is 2-{\em cospanned\ }
\index{$2$-cospanned}
by $\FF_c$ if $f$ is 2-spanned by $\FF_c$ in $M^*$.

Let $(R,B)$ be a {\em well-coblocked\ }
\index{well coblocked} 
$3$-separation of $M/f$.
Recall that this means that $f$ coblocks every $3$-separation
equivalent to $(R,B)$. By Lemma~\ref{unsplit},
$(R,B)$ is an unsplit $4$-separation in $M\ba f$.
By Lemma~\ref{4-cross0}, $(R,B)$ crosses either
one or two members of $(\check P_1,\ldots,\check P_m)$. Recall that
$(R,B)$ is 1-{\em crossing\ } or 2-{\em crossing\ } for $\FF_d$
according as to which case holds. Note that, if $(R,B)$ is
$1$-crossing, there is a flower $(P_1,\ldots,P_m)$, displayed
by $\FF_d$ such that 
$R=(P_1\cap R)\cup P_2\udots P_i$
for some $i\in\{1,\ldots,m\}$. 
As $f\in\cl_M(R)$ and $f\in\cl_M(B)$,
$i\in\{2,\ldots,m-1\}$, as otherwise $M$ is $k$-fractured.

\begin{lemma}
\label{feral1}
\begin{itemize}
\item[(i)] Up to labels $\hP_1$ is blocked by $f$.
\item[(ii)] If some well-coblocked $3$-separation displayed by
$\FF_c$ is $1$-crossing for $\FF_d$, then $f$ is $1$-blocking
for $\FF_d$.
\item[(iii)] If some well-coblocked $3$-separation displayed
by $\FF_c$ is $2$-crossing for $\FF_d$, then $f$ is 
$2$-spanned by $\FF_d$.
\item[(iv)] $f$ is either $1$-blocking in $\FF_d$ 
or is $2$-spanned by $\FF_d$.
\end{itemize}
\end{lemma}

\begin{proof} By Lemma~\ref{x-blocks}, there is a well-coblocked 
$3$-separation
$(R,B)$ displayed by $\FF_c$. By the observations made prior to the
lemma, 
$(R,B)$ is either 1-crossing or
2-crossing for $\FF_d$. 

\begin{sublemma}
\label{feral1.0}
If $(R,B)$ is $2$-crossing, then (i) and (iii)
hold.
\end{sublemma}

\subproof
Assume that $(R,B)$ is 2-crossing. By Lemma~\ref{4-cross0}, 
up to labels, 
$R\subseteq \hP_1\cup\hP_2$.
As $f\in\cl(R)$, the element $f$ is $2$-spanned by $\FF_d$. 
Thus (iii) holds. Let 
$(\hP_1,P'_2,P'_m)=(\hP_1,\hP_2-P_2^-,E(M\ba f)-(\hP_1\cup\hP_2))$. 
Assume that $f$ blocks $P_2'$. If $f$ does not block $\hP_2$,
then $f\in\cl(\hP_2)$ and we deduce that $M$ is $k$-fractured.
Thus
$f$ blocks $\hP_2$ and 
the result holds by an appropriate relabelling of
$\FF_d$. Thus we may assume that $f$ does not block $P'_2$. 
Also $f$ does not
block $P'_m$. If $f$ does not block $\hP_1$, then we
again contradict the fact that $M$ is $k$-coherent.
Hence (i) holds.
\end{proof}

\begin{sublemma}
\label{feral1.5}
If $(R,B)$ is $1$-crossing, then (i) and (ii)
hold.
\end{sublemma}

\begin{proof}
Assume that $(R,B)$ is 1-crossing. 
As observed prior to the lemma, there is a 
flower $(P_1,\ldots,P_m)$ displayed by $\FF_d$ such that  
$R=(R\cap P_1)\cup P_2\udots P_i$ for some
$i\in\{2,\ldots,m-1\}$.
As $f\in\cl(R)$ and $f\in\cl(B)$ we see that 
$f\in\cl(P_1\cup P_2\udots P_i)$ and 
$f\in\cl(P_{i+1}\udots P_m\cup P_1)$.
Thus, if $i\neq 1$, the element $f$ does not block $P_i$. 
Now let
$(\hP_1,P'_2,P'_m)=
(\hP_1,(P_2\cup P_3\udots P_i)-\hP_1,
(P_{i+1}\cup P_{i+2}\udots P_m)-\hP_1)$. As
$f\in\cl(\hP_1\cup P'_m)$ and $f\in\cl(\hP_1\cup P'_2)$, 
the element $f$ does not block $P'_2$ or $P'_m$. 

Assume that $f$ does not block $\hP_1$. It is easily checked that 
$f\not\in\cl(\hP_1),\cl(P'_2),\cl(P'_m)$ and that 
$f\in\cl(\hP_1\cup P'_2),\cl(\hP_1\cup P'_m),\cl(P'_2\cup P'_m)$. As
$(\hP_1,P'_2,P'_m)$ is a swirl-like flower in $M\ba f$, we have 
$\sqcap(\hP_1,P'_2)=\sqcap(\hP_1,P'_m)=\sqcap(P'_2,P'_m)=1$. Thus, by 
Lemma~\ref{pi-minor}, 
$\sqcap_{M/f}(\hP_1,P'_2)=\sqcap_{M/f}(\hP_1,P'_m)=
\sqcap_{M/f}(P'_2,P'_m)=2$,
so that $(\hP_1,P'_2,P'_m)$ is a paddle in $M/f$. 
But, by Lemmas~\ref{1-crossing} and \ref{2-crossing},
$|R\cap \hP_1|,|B\cap \hP_1|\geq 2$. Thus, by Lemma~\ref{flower8.2}, 
there is a 
paddle in $M/f$ that refines $(\hP_1,P'_2,P'_m)$ 
and displays $(R,B)$. 
But, $(R,B)$ is also displayed in the swirl-like flower
$\FF_c$ of this matroid.
However $(R,B)$
cannot be displayed in both a paddle and a swirl-like flower. 
This contradiction shows that
$f$ blocks $\hP_1$ and the sublemma follows.
\end{proof}

Parts (i), (ii) and (iii) of the lemma follow from
\ref{feral1.0} and \ref{feral1.5}. 
Part (iv) follows
from parts (ii) and (iii) and Lemma~\ref{4-cross0}.
\end{proof}

\begin{lemma}
\label{feral1a}
Assume that $f$ is $1$-blocking but not $2$-spanned. Then there is
a flower $(P_1,P_2,\ldots,P_m)$ displayed by $\FF_d$, and an 
$i\in\{3,4,\ldots,m-2\}$ such that the following hold.
\begin{itemize}
\item[(i)] $P_1=\hP_1$ and $f$ blocks $P_1$.
\item[(ii)]  $f\in\cl_M(P_1\cup P_2\udots P_i)$ and 
$f\in\cl_M(P_{i+1}\cup P_{i+2}\udots P_m\cup P_1)$.
\item[(iii)] If $3\leq s\leq i<t\leq m-1$, then 
$P_s\cup P_{s+1}\udots P_t$ is well blocked by $f$.
\item[(iv)] If $(S,T)$ is a well-coblocked $3$-separation
of $M/f$,
displayed by $\FF_c$, then there is a 
$3$-separation $(S',T')$ of $M/f$, equivalent to $(S,T)$
such that
$\{S'-P_1,T'-P_1\}=\{P_2\cup P_3\udots P_i,
P_{i+1}\cup P_{i+2}\udots P_m\}$.
\item[(v)] If $(S,T)$ is a well-coblocked $3$-separation
of $M/f$,
displayed by $\FF_c$, then there is a flower 
$(P'_1,P'_2\ldots,P'_m)$, displayed by $\FF_d$,
such that 
$\{S-P'_1,T-P'_1\}=\{P'_2\cup P'_3\udots P'_i,
P'_{i+1}\cup P'_{i+2}\udots P'_m\}$.
\end{itemize}
\end{lemma}

\begin{proof}
Let $(R,B)$ be a 3-separation  of $M/f$
displayed by $\FF_c$ that is 
well coblocked by $f$. Then, as $f$ is not $2$-spanned, it
follows from Lemma~\ref{feral1} that $(R,B)$ is 1-crossing.
Thus, for some $i\in\{2,3,\ldots,m-1\}$,
we have, up to labels, that 
$R=(\hP_1\cap R)\cup \cP_2\cup\hP_3\udots 
\hP_{i-1}\cup \cP_i\cup L^+_i$,
where $L_i^+$ is an initial segment of $P_i^+$.
As $f$ is not 2-spanned, $i\in\{3,4,\ldots,m-2\}$.
By Lemma~\ref{feral1}, $f$ blocks $\hP_1$. It is now clear
that there is a flower $(P_1,P_2,\ldots,P_m)$
displayed by $\FF_d$ such that $P_1=\hP_1$, 
$P_1\cup P_2\udots P_i =R\cup P_i$, and 
$P_{i+1}\cup P_{i+2}\udots P_m\cup P_1=B\cup P_1$. 
One consequence of this is that (i) holds.
As $f$ coblocks $(R,B)$ we see that $f\in\cl_M(R)$
and $f\in\cl_M(B)$. Hence
$f\in\cl_M(P_1\udots P_i)$ and $f\in\cl_M(P_{i+1}\udots P_m\cup P_1)$,
so that (ii) holds.

Consider (iii). Assume that $3\leq s\leq i<t\leq m-1$. Let $P$
be a $3$-separating set that is equivalent to 
$P_s\cup P_{s+1}\udots P_t$.
Assume that $f$ does not block $P$. Then either $f\in\cl_M(P)$,
or $f\in\cl_M(E(M\ba f)-P)$. Assume that the former holds. Note that
$\cP_m\cap P=\emptyset$ and $\cP_i\subseteq P$. Thus by 
Lemma~\ref{modular}, $f\in\cl(P\cap (P_1\udots P_i))$. But 
$P\cap(P_1\cup P_2\udots P_i)\subseteq \hP_s\cup \hP_{s+1}\udots \hP_i$.
Hence $f\in\cl(\hP_s\cup\hP_{s+1}\udots \hP_i)$. But we also have
$f\in\cl(P_i\udots P_m\cup P_1)$; 
so again using Lemma~\ref{modular},
we see that $f\in\cl(\hP_i)$ contradicting the fact that
$M$ is $k$-coherent. Assume that the latter case holds,
so that $f\in\cl(E(M\ba f)-P)$. Arguing as above, we deduce
that $f\in\cl(\hP_1\cup \hP_2\udots \hP_{s-1})$ and, as
$f\in\cl(P_{i+1}\udots P_m\cup P_1)$, 
we conclude that $f\in\cl(\hP_1)$,
contradicting the fact that $x$ blocks $\hP_1$. Thus (iii)
holds.

Consider (iv). Note that the 3-separation $(R,B)$ at the
start of the proof
determines a labelling of a flower displayed by the bloom
$\FF_d$. Of course the same conclusions also hold for $(S,T)$, 
but we need to
reconcile the labellings. 
As $f$ is not $2$-spanned, it follows from
Lemma~\ref{feral1}, that, up to the choice of labels $S$ and $T$,
that
$S=(\hP_s\cap S)\cup\cP_{s+1}\cup \hP_{s+2}
\udots \hP_{t-1}\cup \cP_t\cup K_t^+$,
where $K_t^+$ is an initial segment of $P_t^+$, and 
$3\leq t-s\leq m-3$.
Note that $\hP_s$ is blocked by $f$.

\begin{sublemma}
\label{feral1a.1}
$s\in\{m,1,2\}$.
\end{sublemma}

\begin{proof}[Subproof.]
Say $s\not\in\{m,1,2\}$.
Up to symmetry, we may assume that $s\leq i$. 
If $s<i$, then $(E(M\ba f)-\hP_s)\supseteq B$,
so $\hP_s$ is not blocked by $f$. Thus we may assume that 
$s=i$. Say that $\cP_m\subseteq T$. Then, as 
$f\in\cl(\hP_i\cup \hP_{i+1}\udots \hP_t)$ 
and $f\in\cl(\hP_1\cup \hP_2\udots\hP_i)$, we 
see that $f\in\cl(\hP_i)$ contradicting the fact that 
$M$ is $k$-coherent. Thus $\cP_m\subseteq S$.
A similar argument shows that $\cP_2\subseteq T$.
Assume that $\cP_1\subseteq T$.
In this case 
$S\cup B=(\hP_i\cap S)\cup\cP_{i+1}\cup \hP_{i+2}\udots 
\hP_{m-1}\cup\cP_m\cup(\hP_1\cap B)$.
But by uncrossing we have $\lambda_{M/f}(S\cup B)=2$,
so $\lambda_{M\ba f}(S\cup B)\leq 3$. 
Observe that in the case that 
$\lambda_{M\ba f}(S\cup B)= 3$, the 4-separation
$(S\cup B,E(M\ba f)-(S\cup B))$ is unsplit.
However
$S\cup B$ crosses both $\cP_1$ and $\cP_i$, 
the set $S\cup B$ does not have
the form of Lemma~\ref{4-cross0} or \ref{3-cross}.  
Essentially the same argument holds 
in the case that $\cP_1\subseteq T$.
\end{proof}

Up to symmetry we may assume that $t\leq i$.

\begin{sublemma}
\label{feral1a.2}
$t=i$.
\end{sublemma}

\begin{proof}[Subproof.]
Assume that $t<i$. As $f\in\cl(S)$, we have
$f\in\cl(\hP_m\udots\hP_t)$. But then 
$\hP_m\cup \hP_1\udots \hP_{i-1}$ is not blocked by $f$,
contradicting (iii). Thus $t=i$.
\end{proof}

We may now relabel 
$K^+_t$ by $K^+_i$.

\begin{sublemma}
\label{feral1a.3}
Up to separations of $M/f$ equivalent to $(S,T)$,
we may assume that $K^+_i=L_i^+$.
\end{sublemma}

\begin{proof}[Subproof.]
Assume that the claim fails. Then, up to symmetry, we may
assume that $K^+_i$ is a proper subset of $L^+_i$. Let
$y$ be the first element of $L_i^+-K_i^+$. Then either 
$y\in\cl(\cP_i\cup K_i^+)$ or $y\in\cl^*(\cP_i\cup K_i^+)$,
so that either $y\in\cl_{M\ba f}(S)$ or $y\in\cl^*_{M\ba f}(S)$.
In  the former case it is clear that $y\in\cl_{M/f}(S)$. 
Consider the
latter case. 
Note that $B\subseteq E(M\ba f)-(\cP_i\cup K^+_i)$, so $y$
is a coloop of $M\ba (\cP_i\cup K_i^+)$, and thus $y$ is a coloop
$M/f\ba (\cP_i\cup K_i^+)$. Thus $y\in\cl^*_{M/f}(\cP_i\cup K_i^+)$
and hence $y\in\cl^*_{M/f}(S)$.

We now have that in $M/f$ the 3-separation  $(S,T)$
is equivalent to $(S\cup\{y\},T-\{y\})$. 
The sublemma  follows by an obvious induction
\end{proof}

We may now assume that 
$P_3\cup P_4\udots P_i\subseteq S$ and 
$P_{i+1}\cup P_{i+2}\udots P_{m-1}\subseteq T$.

\begin{sublemma}
\label{feral1a.4}
$s=1$.
\end{sublemma}

\begin{proof}[Subproof.]
Assume that the sublemma fails. Up to symmetry we may assume that
$s=2$, so that $S\subseteq \hP_2\cup P_3\udots P_i$.
Recall that $B\subseteq P_{i+1}\cup P_{i+2}\udots P_m\cup\hP_1$. 
Say that $B\cap P_1^+=\emptyset$. 
Then $B\subseteq P_{i+1}\cup P_{i+2}\udots P_m\cup(\hP_1-P_1^+)$.
And, as $f\in\cl_M(B)$ and $f\in\cl_M(S)$, we 
see that $f$ is in the guts of a $3$-separation
of $M$ contradicting the fact that $M/f$ is 3-connected. Thus
$B\cap P_1^+\neq\emptyset$. By Lemma~\ref{4-cross0},
$B\cap P_i^+$ consists of a single element $b$, and this 
element is the first element of $P_i^+$. 
If $b\in\cl_{M\ba f}(\cP_1)$, then $f\in\cl_M(B-\{b\})$
and we again obtain the contradiction that $M/f$
is not $3$-connected. Thus $b\in\cl^*_{M\ba f}(B-\{b\})$.
By symmetry there is a single element $s\in S\cap P_1^+$ and 
$s$ is the last element of $P_1^+$. 
If $b\neq s$, then we again obtain 
a $3$-separation $(X,Y)$ of $M\ba f$ with 
$B\subseteq X$ and $S\subseteq Y$, again contradicting the fact
that $M/f$ is 3-connected. Thus $B_1^+=\{b\}$.

Now $f\in\cl_M(S)$, that is $f\in\cl_M((S-\{b\})\cup \{b\})$. 
Hence $b\in\cl_M((S-\{b\})\cup \{f\})$. But $S-\{b\}\subseteq R$,
so that $b\in\cl_M(R\cup\{f\})$ and $b\in\cl_{M/f}(R)$
Thus, in $M/f$ we have $(R,B)\cong (R\cup\{b\},B-\{b\})$.
As $B$ is well-blocked, $f\in\cl_M(B-\{b\})$.
But $B-\{b\}\subseteq P_{i+1}\cup P_{i+2}\udots P_m\cup\cP_1$,
and $S\subseteq \hP_2\udots P_{i-1}\cup P_i$, and
again we contradict the fact the $M/f$ is 3--connected.
\end{proof}

Part (iv) of the lemma now follows immediately. The proof of (v)  
follows from the fact that the equivalence moves performed
in the proof of (iv) could also have 
been thought of as equivalence moves
in the flower.
\end{proof} 

We now consider the case when $f$ is $2$-spanned. Note that if 
$f\in\cl(\hP_1\cup\hP_2)$, then 
$(\hP_1\cup\hP_2\cup\{f\},\hP_3,\ldots,\hP_m)$
induces a swirl-like flower of $M$ of order $m-1$. Thus $m=k$. 

\begin{lemma}
\label{feral2}
Assume that $f$ is $2$-spanned by $\FF_d$, say 
$f\in\cl(\hP_1\cup\hP_2)$.
If $(P_1,P_2,\ldots,P_k)$ is any flower displayed by $\FF_d$ with 
$P_1\cup P_2=\hP_1\cup\hP_2$, then the following hold.
\begin{itemize}
\item[(i)] If $i\in\{3,4,\ldots,k-1\}$, 
then $P_2\cup P_3\udots P_i$ is a well-blocked
$3$-separation of $M\ba f$.
\item[(ii)] If $(R,B)$ is a well-coblocked $3$-separation of 
$M/f$, then, up to labels, $R\subseteq P_1\cup P_2$.
\end{itemize}
\end{lemma}

\begin{proof}
Consider (i). Say $i\in\{3,4,\ldots,k-1\}$. 
Assume that $P_2\cup P_3\udots P_i$ is not well blocked by $f$.
Then either $f\in\cl(\hP_2\cup \hP_3\udots\hP_i)$, or 
$f\in\cl(\hP_{i+1}\cup\hP_{i+2}\udots\hP_k)$.
Up to symmetry we may assume that  $f\in\cl(\hP_2\udots\hP_i)$.
But $f\in\cl(\hP_1\cup\hP_2)$ and
$(\hP_1\cup\hP_2)\cup(\hP_2\udots \hP_i)$ avoids $\cP_m$.
So, by Lemma~\ref{modular}, $f\in\cl(\hP_2)$ contradicting the
fact that $M$ is $k$-coherent. Thus (i) holds.

Consider (ii). Say $(R,B)$ is a well-coblocked 
3-separation of $M/f$. 
If $(R,B)$ is 2-crossing, then, up to labels, 
$R\subseteq \hP_i\cup\hP_{i+1}$
for some $i$. If $i\neq 1$, then we obtain a 
contradiction by an application
of Lemma~\ref{modular}. 
A similar easy argument shows that (ii) holds in the 
case that $(R,B)$ is 1-crossing
with one exceptional case that we focus on now. 
In this case, up to labels,  we have
$R=(\hP_k\cap R)\cup \cP_1\cup P_1^+\cup \cP_2\cup L_2^+$, where
$L_2^+$ is an initial segment of $P_2^+$. 
Let  $P'_2=\cP_2\cup L_2^+$. 
Assume we are in this case.

\begin{sublemma}
\label{feral2.1}
$f\in\cl(\hP_1\cup P'_2)$.
\end{sublemma}

\begin{proof}[Subproof.]
As $f\in\cl(R)$, we have $f\in\cl(\hP_k\cup\hP_1\cup P'_2)$.
Also $f\in\cl(\hP_1\cup\hP_2)$, so by Lemma~\ref{modular},
$f\in\cl(\hP_1\cup P'_2)$. 
\end{proof}

\begin{sublemma}
\label{feral2.2}
$\hP_1\subseteq R$.
\end{sublemma}

\begin{proof}[Subproof.]
Assume the sublemma fails. Then, by Lemma~\ref{1-crossing}, 
there is a single element $b\in B\cap\hP_1$. Moreover, 
$(\hP_1\cup P'_2)\cap B=\{b\}$. By \ref{feral2.1},
$f\in\cl(\hP_1\cup P'_2)$. Also $f\notin\cl((\hP_1\cup P'_2)-\{b\})$,
otherwise $M/f$ is not $3$-connected. Hence
$b\in\cl(((\hP_1\cup P'_2)-\{b\})\cup\{f\})$,
so that $b\in\cl_M(R\cup\{f\})$ and $b\in\cl_{M/f}(R)$.
Thus, in $M/f$, the 3-separation $(R,B)
is equivalent to (R\cup\{b\},B-\{b\})$.
As $(R,B)$ is well coblocked, this means that $f\in\cl_M(B-\{b\})$.
But, as $f\in\cl(\hP_1\cup P'_2)$, and 
$(B-\{b\})\cap (\hP_1\cup P'_2)=\emptyset$, we have again contradicted
the fact that $M/f$ is $3$-connected.
\end{proof}

Now $B\subseteq E(M\ba f)-(\hP_1\cup P'_2)$, $f\in\cl(B)$, and
$f\in\cl(\hP_1\cup P'_2)$, and again we contradict the fact that
$M/f$ is $3$-connected. Part (ii) of the lemma follows from
this final contradicton.
\end{proof}

\begin{lemma}
\label{feral3}
Either $\FF_d$ is $1$-blocked, or $\FF_c$ is $1$-coblocked.
\end{lemma}

\begin{proof}
Assume the lemma fails. Then, by Lemma~\ref{feral1}(iv), 
we obtain the following up to
labels. For $\FF_d$ we have $f\in\cl(\hP_1\cup \hP_2)$, 
$f$ blocks $\hP_1$, and
$f$ blocks $\hP_2$. For $\FF_c$ we have
$f\in\cl^*(\hQ_1\cup\hQ_2)$, $f$ coblocks $\hQ_1$ 
and $f$ coblocks $\hQ_2$.
Let $(P_1,P_2,\ldots,P_m)$ be a flower displayed by $\FF_d$, where
$P_1\cup P_2=\hP_1\cup\hP_2$ and let
$(Q_1,Q_2,\ldots,Q_n)$ be a flower displayed by $\FF_c$ where
$Q_1\cup Q_2=\hQ_1\cup\hQ_2$. As 
$(P_1\cup P_2\cup\{f\},P_3,\ldots,P_m)$ and
$(Q_1\cup Q_2\cup\{f\},Q_3,\ldots,Q_n)$ are swirl-like flowers in $M$, 
we have $m=n=k$.
Note that $f$ blocks $P_1$ and $P_2$ and $f$ coblocks 
$Q_1$ and $Q_2$,
otherwise $\FF_d$ is $1$-blocked, or $\FF_c$ is $1$-coblocked.

By Lemma~\ref{feral2}(i), $Q_2\cup Q_3$ is a well-coblocked 
$3$-separation of 
$M/f$. Thus, by Lemma~\ref{feral2}(ii), either 
$Q_2\cup Q_3\subseteq P_1\cup P_2$ or 
$Q_4\cup Q_5\udots Q_k\cup Q_1\subseteq P_1\cup P_2$. 
In the latter case $Q_k\cup Q_1\subseteq P_1\cup P_2$
so that up to labels we
may assume that $Q_2\cup Q_3\subseteq P_1\cup P_2$. As
$f\in\cl(Q_2)$, the set
$Q_2$ is not contained in either $P_1$ or $P_2$,
so $Q_2\cap P_1\neq \emptyset$ and 
$Q_2\cap P_2\neq \emptyset$. As
$(Q_1\cup Q_2\cup\{f\},Q_3,\ldots,Q_n)$ are swirl-like flowers in $M$
we see that $(Q_1\cup Q_2,Q_3\cup Q_4\udots Q_k)$ is a 
3-separation of $M\ba f$.
As either $P_1$ or $P_2$ crosses this $3$-separation, it is not
displayed in $(P_1,P_2,\ldots,P_k)$. As $Q_1$ meets both
$P_1$ and $P_2$, the set $Q_1\cup Q_2$ is not contained in a petal
of $(P_1,P_2,\ldots,P_k)$. Thus
$Q_3\cup Q_4\udots Q_k$ is contained in a petal of 
$(P_1,P_2,\ldots,P_k)$. 
As $Q_3\subseteq P_1\cup P_2$,
either $Q_3\cup Q_4\udots Q_k\subseteq \hP_1$ or 
$Q_3\cup Q_4\udots Q_k\subseteq \hP_2$.
Up to labels we may assume that the latter case holds.
Evidently we may also assume that $P_2=\hP_2$. Thus
$Q_3\cup Q_4\udots Q_k\subseteq P_2$.

Altogether we have $Q_2\subseteq P_1\cup P_2$, 
$Q_3\cup Q_4\udots Q_k\subseteq P_2$ and $Q_1\supseteq P_3\udots P_k$.
If $Q_1\cap P_2=\emptyset$, then, as $f\in\cl(Q_1)$, 
the element $f$ does 
not block $P_2$. Thus $Q_1\cap P_2\neq \emptyset$ and 
similarly $Q_1\cap P_1\neq \emptyset$.

Let $R=Q_3\cup Q_4\udots Q_k\cup (Q_2\cap P_2)$. 
By Lemma~\ref{2-crossing},
$\lambda_{M\ba f}(R)=2$. But $f\in\cl_M(E(M)-(\{f\}\cup R))$,
so $\lambda_{M/f}(R)=2$.

\begin{sublemma} 
\label{feral3.1}
In $M/f$ we have $R\ncong R\cup Q_2$ and
$R\ncong R-Q_2$.
\end{sublemma}

\begin{proof}[Subproof.]
If $R\cong R\cup Q_2$, then
$Q_1\cong Q_1\cup(Q_2\cap P_1)$ and, as $f\not\in\cl_M(P_2)$,
we see that $f$ is 1-blocking for $\FF_d$. Assume that 
$R\cong R-Q_2$. In this
case, as $f$ is not 1-blocking for $\FF_d$, we see that either
$f\in\cl(Q_2\cap P_1)$, 
contradicting the fact that $f$ blocks $P_1$,
 or $f\in\cl_M(Q_3\cup Q_4\udots Q_k\cup(P_2\cup Q_2))$, that is, 
$f\in\cl_M(P_2)$, contradicting the fact that $f$ blocks $P_2$.
\end{proof}

By \ref{feral3.1} $R$ is not displayed in $\FF_c$ and neither 
$R$ nor its complement are contained in a petal of $\FF_c$.
This contradicts the fact that $\FF_c$ is a maximal bloom
of $M/f$ and the lemma follows.
\end{proof}

\begin{lemma}
\label{feral5}
Either $f$ is $2$-spanned by $\FF_d$ or $f$ is
$2$-cospanned by $\FF_c$.
\end{lemma}

\begin{proof}
Assume that the lemma fails so that $f$ is not 2-spanned by 
$\FF_d$ and is not $2$-cospanned by $\FF_c$. 
Then there is a flower $(Q_1,Q_2,\ldots,Q_n)$,
displayed by $\FF_c$, and an $i\in\{3,4,\ldots,n-2\}$
such that the dual of Lemma~\ref{feral1a} holds.
By that lemma 
$f\in\cl^*(Q_1\cup Q_2\udots Q_i)$
and $f\in\cl^*(Q_{i+1}\cup Q_{i+2}\udots Q_m\cup Q_1)$. 
Also 
$Q_i\cup Q_{i+1}$ is well coblocked by
$f$. It is easily seen that we may assume that 
$Q_i\cup Q_{i+1}=\hQ_i\cup\hQ_{i+1}$.

Also by Lemma~\ref{feral1a}, there is 
flower $(P_1,P_2,\ldots,P_m)$, displayed 
by $\FF_d$, and a $j\in\{3,4,\ldots,m-2\}$, such that 
$$Q_i\cup Q_{i+1}=(P_1\cap(Q_i\cup Q_{i+1}))\cup P_2\udots P_j.$$
Thus $P_2\cup P_3\udots P_j\subseteq Q_i\cup Q_{i+1}$.

As $Q_i\cup Q_{i+1}$ is coblocked by $f$, we have
$f\in\cl(P_{j+1}\cup P_{j+2}\udots P_m\cup P_1)$. Thus, 
$\lambda_M(P_2\cup P_3\udots P_j)=2$ and hence
$\lambda_{M/f}(P_2\cup P_3\udots P_j)=2$. 
Clearly $P_{j+1}\udots P_m\cup P_1$ is not contained in a petal
of $\FF_c$. Thus $P_2\cup P_3\udots P_j$ 
is either contained in a petal of $\FF_c$
or it is displayed by more than one petal of $\FF_c$. As 
$P_2\cup P_3\udots P_j\subseteq Q_i\cup Q_{i+1}$, the latter case
implies that $Q_i\cup Q_{i+1}\cong P_2\cup P_3\udots P_j$, 
contradicting
the fact that $Q_i\cup Q_{i+1}$ is well coblocked. 
Hence $P_2\cup P_3\udots P_j$ is contained in a petal of $\FF_c$. 
Assume without
loss of generality that $P_2\cup P_3\udots P_j\subseteq \hQ_i$.

Consider $P_j\cup P_{j+1}$. This is a well-blocked 
$3$-separation of $M\ba f$. We may assume that 
$P_j\cup P_{j+1}=\hP_j\cup \hP_{j+1}$. By Lemma~\ref{feral1a},
either $P_j\cup P_{j+1}\supseteq Q_2\cup Q_3\udots Q_i$, 
contradicting the
fact that $P_2\subseteq Q_i$, or 
$P_j\cup P_{j+1}\supseteq Q_{i+1}\cup Q_{i+2}\udots Q_n$ 
contradicting the fact that $P_j\subseteq Q_i$
and the lemma follows. 
\end{proof}

\begin{lemma}
\label{feral5a}
Up to duality, $f$ is $1$-blocking for $\FF_d$ and
$2$-cospanned by $\FF_c$.
\end{lemma}

\begin{proof}
By Lemma~\ref{feral3}, we may assume that $f$ is 1-blocking
for $\FF_d$. Assume that $f$ is not $2$-cospanned by $\FF_c$.
By Lemma~\ref{feral5}, $f$ is 2-spanned by $\FF_d$, and
by Lemma~\ref{feral1}(iii), $f$ is 1-coblocking for 
$\FF_c$. The lemma now follows by taking the dual.
\end{proof}

\begin{lemma}
\label{feral6}
If $f$ is $1$-blocking for $\FF_d$ and
$f$ is $2$-cospanned by $\FF_c$, 
then there is a feral display for $f$ obtained from flowers
displayed by $\FF_d$ and $\FF_c$.
\end{lemma}

\begin{proof}
As $f$ is $2$-cospanned by $\FF_c$, we may assume that 
$f\in\cl(\hQ_1\cup\hQ_2)$. By Lemma~\ref{feral1}(i), 
we may assume that
$\hQ_1$ is coblocked by $f$. Let $(Q_1,Q_2,\ldots,Q_k)$
be a flower displayed by $\FF_c$ where $Q_1\cup Q_2=\hQ_1\cup\hQ_2$.

We split the proof into two cases. 
For the first case assume that $f$ is not 
2-spanned by $\FF_d$. Then Lemma~\ref{feral1a} applies. Let
$(P_1,P_2,\ldots,P_m)$ be a flower displayed by $\FF_d$ satisfying
Lemma~\ref{feral1a}. By Lemma~\ref{feral2}, $Q_2\cup Q_3$
is well coblocked by $f$. Thus, up to labels, there is an
$i\in\{3,4,\ldots,m-2\}$ and a 3-separating set 
$R_3$ of $M/f$, equivalent to $Q_2\cup Q_3$, such that 
$R_3= P_{i+1}\cup P_{i+2}\udots P_m\cup(R_3\cap P_1)$. 
Indeed, by replacing $(Q_1,Q_2,\ldots,Q_k)$ by an 
equivalent flower, we may assume that $R_3=Q_2\cup Q_3$.
Also $Q_2\cup Q_3\udots Q_{k-1}$ is well-coblocked, 
so there is a 
3-separation $R_4$ equivalent to 
$Q_2\udots Q_{k-1}$ such that 
$R_4=P_{i+1}\cup P_{i+2}\udots P_m\cup (P_1\cap R_4)$.
Again we may replace $(Q_1,Q_2,\ldots,Q_k)$ by an equivalent
flower so that $R_4=Q_2\cup Q_3\udots Q_{k-1}$. It is easily
seen that the only petal affected by such equivalence moves is 
$Q_{k-1}$, so that we have both $R_3=Q_2\cup Q_3$
and $R_4=Q_2\cup Q_3\udots Q_{k-1}$. It now follows that
$Q_4\cup Q_5\udots Q_{k-1}\subseteq P_1$.

As $Q_4\cup Q_5\udots Q_{k-1} \subseteq P_1$, 
and $P_1$ is fully closed in $M\ba f$,
it follows from Lemma~\ref{in-petal} that $P_1$ contains all but one
petal of $(Q_1\cup Q_2,Q_3,\ldots,Q_k)$. Assume that $Q_3$ is not
a subset of $P_1$, then 
$Q_3\supseteq P_2\cup P_3\udots P_m$ contradicting the
fact that $\cP_2\subseteq P_1\cup P_2\udots P_i$. Thus $Q_3$, and similarly
$Q_k$, is contained in $P_1$.

Summing up we have $Q_3\cup Q_4\udots Q_k\subseteq P_1$. 
Moreover, by choosing $(Q_1,Q_2,\ldots Q_k)$ to satisfy the 
constraints imposed by $R_3$ and $R_4$, we have 
$P_1\cup P_2\udots P_i\subseteq Q_1$ and 
$P_{i+1}\cup P_{i+2}\udots P_k\subseteq Q_2$. 
Let $Z_1=Q_1\cap P_1$ and $Z_2=Q_2\cap P_1$. As $Q_1$ is 
coblocked by $f$ and $\lambda_M(P_{i+1}\cup P_{i+2}\udots P_m)=2$,
we see that
$Q_1\neq P_{i+1}\cup P_{i+2}\udots P_m$. Thus $Z_1\neq \emptyset$.
Also $\lambda_{M\ba f}(P_1)=2$,
$\lambda_{M\ba f}(Q_1)=3$, 
$\lambda_{M\ba f}(P_1\cup Q_1)=2$. So by uncrossing we have
$\lambda_{M\ba f}(Z_1)\leq 3$. As $f\in\cl(Q_2)$, 
$\lambda_M(Z_1)\leq 3$ and similarly $\lambda_M(Z_2)\leq 3$.
It is now a matter of routine bookkeeping to verify that 
$(P_1,P_2,\ldots,P_m)$ and $(Q_1,Q_2\ldots,Q_k)$ form a feral 
display for $f$.

Consider the second case. Assume that $f$ is 2-spanned by 
$\FF_d$. Note that this means that $m=k$. The bloom $\FF_d$ is
both 1-blocked and 2-spanned and $\FF_c$ is 2-spanned.
By Lemma~\ref{feral1}, $\FF_c$ is also 1-coblocked. Thus we have
a flower $(P_1,P_2,\ldots,P_k)$ displayed by $\FF_d$ such that 
$P_1\cup P_2=\hP_1\cup\hP_2$, the element 
$f$ blocks $P_1$ and no other petal,
and $f\in\cl(P_1\cup P_2)$. We also have a flower
$(Q_1,Q_2,\ldots Q_k)$ displayed by $\FF_c$ such that 
$Q_1\cup Q_2=\hQ_1\cup\hQ_2$, the element $f$ coblocks $Q_1$ 
and no other petal, and
$f\in\cl^*(Q_1\cup Q_2)$. Let $P=P_3\cup P_4\udots P_k$ and 
$Q=Q_3\cup Q_4\udots Q_k$.
Note that 
$\lambda_M(P)=\lambda_M(Q)=\lambda_M(P_2)=\lambda_M(Q_2)$=2.

By Lemma~\ref{in-petal}, either $Q$ or 
$Q_1\cup Q_2$ is contained in 
a petal of $\FF_d$. But $f\in\cl_M(Q_1\cup Q_2)$ and no petal
of $\FF_d$ has that property. Hence $Q\subseteq \hP_j$ for some 
$j\in\{1,2,\ldots,k\}$. By Lemma~\ref{feral2}, 
$(Q_2\cup Q_3,Q_4\udots Q_k\cup Q_1)$ is well-coblocked by
$f$. Thus, by Lemma~\ref{feral2}, either 
$Q_2\cup Q_3\subseteq P_1\cup P_2$,
or $Q_4\cup Q_5\udots Q_k\cup Q_1\subseteq P_1\cup P_2$.

From the above information we draw two conclusions. First we have either
$Q_3$ or $Q_4$ is contained in $P_1\cup P_2$. Thus, as 
$Q_3\cup Q_4\subseteq Q$,
and $Q\subseteq\hP_j$ for some $j\in\{1,2,\ldots,k\}$, we have 
$Q\subseteq \hP_s$ for some $s\in\{1,2\}$. 
The second conclusion is that
$P\subseteq \hQ_t$ for some $t\in\{1,2\}$. 
Up to duality, there are 
three cases: $(s,t)=(1,2)$, $(s,t)=(1,1)$, and $(s,t)=(2,2)$.

\begin{sublemma}
\label{feral6.1}
The lemma holds if $(s,t)=(1,2)$.
\end{sublemma}

\begin{proof}[Subproof.]
In this case $Q\subseteq \hP_1$ and $P\subseteq \hQ_2$. 
We may assume that
$P_1=\hP_1$ and $Q_1=\hQ_1$. As $Q_2$ is not coblocked,
$\lambda_{M\ba f}(Q_2)=2$. 

We now show that we may assume that $Q_2\cap P_2=\emptyset$. By 
Lemma~\ref{3-cross} either $Q_1\cap P_2$ consists of loose elements
of $P_2$ contained in $\fcl_{M\ba f}(P_1)$ or $P_2\cap Q_2$
consists of loose elements of $P_2$ contained
in $\fcl_{M\ba f}(P_3)$. The former
case contradicts the fact that $P_1$ is fully closed. 
Hence the latter
case holds and, by moving to a flower equivalent to 
$(P_1,P_2,\ldots,P_k)$, we may assume that $Q_2\cap P_2=\emptyset$,
so that $P_2\subseteq Q_1$.

Next we show that we may assume that $P_1\cap Q_2=\emptyset$.
Either $Q_2\cap P_1$ or
$P_1-Q_2$ is a set of loose elements of $P_1$
in $M\ba f$. The latter case implies that 
$Q\subseteq \fcl_{M\ba f}(Q_2)$,
so that $Q\subseteq \fcl_M(Q_2\cup \{f\})$ 
contradicting the fact that
$(Q_1\cup Q_2\cup\{f\},Q_3,\ldots,Q_k)$ is a 
swirl-like flower of order $k-1$ in $M$.
Hence $Q_2\cap P_1$ is a set of loose elements of $P_1$, so that
$(P_1,P_2,\ldots,P_k)\cong(P_1-Q_2,P_2,\ldots,P_k\cup(Q_2-P_1))$. By
Lemma~\ref{feral1}, some petal of this flower is blocked by
$f$. Certainly $P_2$ is not blocked.
Moreover, $Q_1\subseteq (P_1-Q_1)\cup P_2$
and $f\in\cl_M(Q_1)$, so no petal in
$\{P_3,P_4,\ldots,P_{k-1},P_k\cup (Q_2-P_1)\}$ is blocked by $f$.
Therefore $P_1-Q_2$ is blocked by $f$. We may now relabel
the above flower to $(P_1,P_2,\ldots,P_k)$. 
In this flower $P_1$ is blocked and
$P_1\cap Q_2=\emptyset$.

Let $Z_1=Q_1\cap P_1$.
We have established the following: 
$Q_2=P_3\cup P_4\udots P_k$, $P_2\subseteq Q_1$, and
$P_1=Q_3\udots Q_k\cup(Q_1\cap P_1)$. 
Let $(P_1',P_2',\ldots,P_k')=(P_1,P_k,P_{k-1},\ldots,P_2)$.
Then we have 
$Q_2=P'_2\cup P'_3\udots P'_{k-1}$, $Q_1=P'_k\cup Z_1$, and 
$P'_1=Q_3\cup Q_4\udots Q_k\cup Z_1$. 
It is now a matter of routing bookkeeping that 
$(P'_1,P'_2,\ldots,P'_k)$ and $(Q_1,Q_2,\ldots,Q_k)$
form a feral display for $f$.
\end{proof}

\begin{sublemma}
\label{feral6.2}
The lemma holds if $(s,t)=(1,1)$.
\end{sublemma}

\begin{proof}[Subproof.]
In this case $Q\subseteq \hP_1$ and $P\subseteq \hQ_1$. Clearly we may
assume that $P_1=\hP_1$ and $Q_1=\hQ_1$. Consider $Q_2$. Note that
$Q_2\cap P_2\neq\emptyset$, otherwise $Q\cup Q_2\subseteq P_1$ so that
$f\notin\cl(Q\cup Q_2)$ implying that $Q_1$ is not coblocked.
As $\lambda_{M\ba f}(Q_2)=2$, and $\FF_d$ is a maximal bloom of 
$M\ba f$, either $P_1\cong P_1\cap Q_2$ or $P_2\cong P_2\cup Q_2$. 
The former case contradicts the fact that $P_1$ is fully closed.
Thus $P_2\cong P_2\cup Q_2$. Thus, by moving to the equivalent flower,
$(P_1-Q_2,P_2\cup Q_2,P_3,\ldots,P_k)$
we have $Q_2\subseteq P_2$. It may be that $P_1-Q_2$ 
is no longer blocked.
But then $P_2\cup Q_2$ is blocked and we are 
in the case covered by 
\ref{feral6.1}.
Therefore we may assume that $P_1-Q_2$ is blocked by $f$ and, 
after an appropriate 
relabelling, we may assume that $Q_2\subseteq P_2$.

We now have $Q_2\subseteq P_2\subseteq Q_2\cup Q_1$. As
$\FF_c$ is a maximal bloom of $M/f$ and $\lambda_{M/f}(P_2)=2$,
either $P_2-Q_2$ is a set of loose elements of $Q_2$, or
$Q_1-P_1$ is a set of loose elements of $Q_1$. The latter
case implies that $P=P_3\cup P_4\udots P_k$ is contained in 
$\fcl_M(P_1\cup P_2\cup\{f\})$, contradicting the fact that
$(P_1\cup P_2\cup\{f\},P_3,\ldots,P_k)$ is a 
swirl-like flower of order
$k-1$ in $M$. Hence $P_2-Q_2$ is a set of loose elements of $P_2$
and, by moving to an equivalent flower we may assume that $P_2=Q_2$.

Summing up, we have the following: $Q_2=P_2$, 
$Q_3\udots Q_k\subseteq P_1$, and $P_3\udots P_k\subseteq Q_1$. Let
$Z_1=P_1\cap Q_1$, $Z_2=P_1\cap Q_2$. Then 
$Z_2=\emptyset$. Moreover, it is easily checked that
$Z_1\neq \emptyset$ and  $\lambda_M(Z_1)\leq 3$.
Thus $(P_1,P_2,\ldots,P_k)$ and $(Q_1,Q_2,\ldots,Q_k)$
form a feral display for $f$.
\end{proof}

\begin{sublemma}
\label{feral6.3}
The lemma holds if $(s,t)=(2,2)$.
\end{sublemma}

\begin{proof}[Subproof.]
In this case $Q\subseteq \hP_2$ and $P\subseteq \hQ_2$. If 
$f$ does not block $(\hP_1\cup\hP_2)-\hP_2$, then
$f$ blocks $\hP_2$ and we are in a case that is equivalent to that
of \ref{feral6.1}. So we may assume that 
$f$ blocks $(\hP_1\cup\hP_2)-\hP_2$, and similarly,
$f$ coblocks $(\hQ_1\cup\hQ_2)-\hQ_2$. Thus there are flowers
$(P_1,P_2,\ldots,P_k)$ and $(Q_1,Q_2,\ldots,Q_k)$ displayed by
$\FF_d$ and $\FF_c$ respectively such that $f$ blocks $P_1$,
$f$ coblocks $Q_1$, $P_1\cup P_2=\hP_1\cup\hP_2$,
$Q_1\cup Q_2=\hQ_1\cup \hQ_2$, $P\subseteq Q_2$, and $Q\subseteq P_2$.

As $f$ does not block $Q_2$, we may argue, just as in the previous case,
that, by moving to an equivalent flower, we may assume that $Q_2=P$.
We now have $Q\subseteq P_2\subseteq Q\cup Q_1$ and $\lambda_{M/f}(P_2)=2$.
Therefore in $M/f$, either $Q\cong P_2$, or
$P_2\cong P_2\cup Q_1$, that is, $P_2\cong P_2\cup P_1$.
The former case implies that there is a flower displayed by
$\FF_c$ such that no petal is coblocked by $f$, so that case does not
occur. Consider the latter case. In this case, $P_1$ is a set of
loose elements of $Q_1$, so we may move to an equivalent flower where
$Q_1\subseteq P_2$. We now have a flower
$(Q_1',Q'_2,\ldots,Q_k')$, where $Q_1'\subseteq P_2$, so that 
$Q_1'$ is not coblocked by $f$. But this means that $f$ blocks
$\hQ_2$ and we are again in the case covered by \ref{feral6.1}.
\end{proof}
All cases have been covered and the lemma follows.
\end{proof}

Theorem~\ref{feral} is an immediate consequence of 
Lemmas~\ref{feral6}
and \ref{feral5a}.

\section{$3$-trees and $k$-coherence}
\label{3tree}

Flowers provide a way of representing certain 
3-separations in a matroid.
It was shown in \cite{flower} that, by using a certain
type of tree, one can simultaneously display a representative of
each equivalence class of non-sequential 3-separations of 
a 3-connected matroid $M$. In this section we describe these
trees and their interaction with $k$-coherence.

Let $\pi$ be a partition of a finite set $E$, where some members of
$\pi$ may be empty, and let $T$ be a tree such that every member of
$\pi$ labels a vertex of $T$. 
We say that $T$ is a $\pi$-{\em labelled tree};
\index{$\pi$-labelled tree} 
labelled vertices are called {\em bag vertices}
\index{bag vertex} 
and members
of $\pi$ are called {\em bags}.
\index{bag}

Let $G$ be a subgraph of $T$ 
with  components $G_1,G_2,\ldots,G_m$. Let $X_i$ be the union of those bags that label vertices of $G_i$.    
Then  {\em the subsets of $E$ displayed by $G$} are 
$X_1,X_2,\ldots,X_m$. In particular, if $V(G) = V(T)$, then 
$\{X_1,X_2,\ldots,X_m\}$ is the {\em partition of $E$  
displayed by $G$}. 
Let $e$ be an edge
of $T$. The {\em partition of $E$ displayed by $e$} 
is the partition
displayed by $T\ba e$. If $e=v_1v_2$ for  vertices 
$v_1$ and $v_2$, 
then $(Y_1,Y_2)$ is the {\it (ordered) 
partition of $E(M)$ displayed by $v_1v_2$} if $Y_1$ 
is the union of the bags in the component of $T\ba v_1v_2$ 
containing $v_1$. 
Let $v$ be a vertex of $T$.
The {\em partition of $E$ displayed by $v$} is the partition
displayed by  $T-v$. The edges incident with $v$  correspond to the components of $T-v$, and
hence to the members of the partition displayed by $v$. Note
that, if $v$ is not a bag vertex, the the partition
displayed by $v$ is a partition of $E$, while, if $v$ is a 
bag vertex, then the partition is a partition of $E-B$ where
$B$ is the bag labelling $v$.
In what follows,
if a cyclic ordering $(e_1,e_2,\ldots,e_n)$ is imposed on the edges
incident with $v$, this cyclic ordering is taken to represent the
corresponding cyclic ordering on the 
members of the partition displayed
by $v$.

Let $M$ be a  $3$-connected matroid with ground set $E$.
An {\em almost partial $3$-tree} $T$ for $M$ is a 
$\pi$-labelled tree, where
$\pi$ is a partition of $E$ such that:
\begin{itemize}
\item[(i)] For each edge $e$ of $T$, the partition 
$(X,Y)$ of $E$ displayed by
$e$ is $3$-separating, and, if $e$ is incident 
with two bag vertices,
then $(X,Y)$ is a non-sequential $3$-separation.
\item[(ii)] Every non-bag vertex $v$ is 
labelled either $D$ or $A$;
if $v$ is labelled $D$, then there is a cyclic
ordering on the edges incident with $v$.
\item[(iii)] If a vertex $v$ is labelled $A$, 
then the partition of
$E$ displayed by $v$ is a tight maximal 
anemone of order at least 3.
\item[(iv)] If a vertex $v$ is labelled $D$, 
then the partition of $E$
displayed by $v$, with the cyclic order 
induced by the cyclic ordering on
the edges incident with $v$, is a tight maximal 
daisy of order at least 3.
\end{itemize}
By conditions (iii) and (iv), a vertex $v$  labelled $D$ or $A$ corresponds to
a flower of $M$. The $3$-separations displayed by this flower 
are the $3$-separations {\em displayed by} $v$. 
A vertex of a partial $3$-tree is referred to as a {\em daisy
vertex\ }
\index{daisy vertex} 
or an {\em anemone vertex\ }
\index{anemone vertex} 
if it is labelled $D$ or $A$,
respectively.  A vertex labelled either $D$ or $A$ is a 
{\em flower vertex}.
\index{flower vertex}  
A 3-separation is 
{\it displayed by an almost partial 
$3$-tree} $T$ if it is displayed by some edge or some 
flower vertex of $T$. 

A $3$-separation $(R,G)$ of $M$
{\em conforms\ } with an almost partial $3$-tree $T$ 
if either $(R,G)$ is equivalent to a
$3$-separation that is displayed by a 
flower vertex or an edge of $T$,
or $(R,G)$ is equivalent to a $3$-separation $(R',G')$ with the
property that either $R'$ or $G'$ is contained in a bag of $T$.

An almost partial $3$-tree for $M$ is a 
{\em partial $3$-tree\ }
\index{partial $3$-tree} 
if every non-sequential $3$-separation of 
$M$ conforms with $T$. 
We now define a quasi order on the set of 
partial $3$-trees for $M$.
Let $T_1$ and $T_2$ be two partial $3$-trees for $M$. 
Then $T_1\preceq T_2$ if all of the 
non-sequential $3$-separations displayed by $T_1$
are displayed by $T_2$. If $T_1\preceq T_2$ and 
$T_2\preceq T_1$, then
$T_1$ is {\em equivalent\ }
\index{equivalent partial $3$-trees} 
to $T_2$. A partial
$3$-tree is {\em maximal\ } if it is maximal 
with respect to this quasi order.

Note  that while flower vertices need to be labelled 
$D$ or $A$, we may suppress these labels when they are
clear from context.
The following theorem is the main result of 
\cite[Theorem 9.1]{flower}.

\begin{theorem}
\label{oswthm}
Let $M$ be a $3$-connected matroid with $|E(M)|\ge 9$,
and let $T$ be a maximal partial $3$-tree for $M$. 
Then every non-sequential $3$-separation of $M$ is 
equivalent to a $3$-separation displayed by $T$.
\end{theorem}

Maximal partial $3$-trees are by no means unique.
Consider the following
situation. Let
$(P_1,P_2,P_3)$ be a maximal flower of order three in a 3-connected
matroid $M$. Then it may be the case that $P_2$ is
sequential, in which case $(P_1,P_2\cup P_3)$ 
and $(P_1\cup P_2,P_3)$ are inequivalent non-sequential  
$3$-separations. Assume that these are the only
inequivalent 
non-sequential 3-separations of $M$.
Then given such a flower, one can obtain distinct
maximal partial 3-trees for $M$ as follows. 
Let $T_1$ be a tree consisting
of a path with vertices $(v_1,v_2,v_3)$ such that $v_i$
labels the bag $P_1$. Let $T_2$ be a star with a flower vertex
of degree $3$ and leaf vertices $v_1$, $v_2$ and $v_3$ labelling
the bags $P_1$, $P_2$ and $P_3$. Then both $T_1$
and $T_2$ are maximal partial 3-trees for $M$. Indeed the situation
can be further complicated by splitting the elements of 
$P_2$ into smaller bags along a path. To get a more canonical
structure we follow \cite{flower2} and say that
a maximal partial 3-tree for $M$ is a {\em $3$-tree} if 
\index{$3$-tree}

\begin{itemize}
\item[(I)] for every tight maximal flower of $M$ of order 
three, there is an equivalent flower that is displayed by
a vertex of $T$; and
\item[(II)] if a vertex is incident with two edges $e$ and 
$f$ that display equivalent 3-separations, then the other
ends of $e$ and $f$ are flower vertices, $v$ has degree two,
and $v$ labels a non-empty bag.
\end{itemize}

The next theorem summarises results from \cite{flower2}.

\begin{theorem}
\label{3-tree-display}
Let $M$ be a $3$-connected matroid. Then $M$ has a $3$-tree $T$.
Moreover, if $\PP$ is a flower of $M$ of order at least
three, then $T$ has a flower vertex $v$ that displays
a flower equivalent to $\PP$.
\end{theorem}

Let $T$ be a $3$-tree for $M$. A bag is a {\em leaf bag\ }
\index{leaf bag} 
if it 
labels a leaf of $T$. We omit the easy proof of the next lemma.

\begin{lemma}
\label{wheel2.5}
Let $M$ be a $3$-connected matroid.
If $L$ is a leaf bag of a $3$-tree for $M$, and 
$x\in\clstar(L)$, then $L\cup\{x\}$
is a leaf bag of a $3$-tree for $M$. Thus $\fcl(L)$ is a leaf bag
for a $3$-tree for $M$.
\end{lemma}

Note that converse of Lemma~\ref{wheel2.5} does not hold 
in that it is not always the case that a
3-separating set equivalent to one in a leaf bag can be 
displayed in
a 3-tree for M.


A subset $L$ of $E(M)$ is a {\em peripheral set\ }
\index{peripheral set} 
of $M$ if 
it is a leaf bag for some $3$-tree $T$ for $M$.
By Lemma~\ref{wheel2.5} the full closure of a peripheral set
is a peripheral set.
An element $x$ of $M$ is 
{\em peripheral\ }
\index{peripheral element} 
if $x$ is an element of a peripheral set.
The theme in what follows is that peripheral elements
of $k$-coherent matroids are well behaved. Indeed this
is true in a broader setting. The next theorem is 
\cite[Theorem~4.2]{upgrade}. 

\begin{theorem}
\label{upgrade}
Let $M$ be a $3$-connected matroid other than a wheel or a 
whirl. Suppose $|E(M)|\geq 9$ and let $S$ be a
peripheral set of $M$. Then $\fcl(S)$ contains an element
$e$ such that either $M\ba e$ or $M/e$ is 
$3$-connected and $e$ does not expose any $3$-separations.
\end{theorem}

This gives the following as an immediate corollary.

\begin{corollary}
\label{upcor}
Let $M$ be a $k$-coherent matroid other than a wheel or a whirl
and let $S$ be a peripheral set of $M$. Then 
$\fcl(S)$ contains an element $e$ such that either
$M\ba e$ or $M/e$ is $k$-coherent.
\end{corollary}

The remainder of this section examines life in peripheral sets
in more detail. We next show that elements of $k$-wild 
triangles are not peripheral.

\begin{lemma}
\label{wheel-wild}
Let $W$ be a $k$-wild triangle of the  $k$-coherent matroid $M$. 
Then no element of $W$ is peripheral.
\end{lemma}

\begin{proof}
Say $W=\{a,b,c\}$ and let $T$ be a 3-tree for $M$. 
By Theorem~\ref{wild1}, $W$ is either a standard
or costandard $k$-wild triangle. 
In either case $\{a,b,c\}$ has a $k$-wild display
$(A_1,A_2,\ldots,A_{k-2},B_1,B_2,\ldots,B_{k-2},
C_1,C_2,\ldots,C_{k-2})$.  
Let
$A=A_1\cup A_2\udots A_{k-2}$, $B=B_1\cup B_2\udots B_{k-2}$ and 
$C=C_1\cup C_2\udots C_{k-2}$.
Then, $\AAA=(A_1,A_2,\ldots,A_{k-2},B\cup C\cup W)$, 
$\BB=(B_1,B_2,\ldots,B_{k-2},A\cup C\cup W)$ and
$\CC=(C_1,C_2,\ldots,C_{k-2},A\cup B\cup W)$ are tight maximal
flowers of $M$ of order at least four. 
Assume that $a$ is peripheral. Then we may assume that
$a\in L$ for some leaf bag of $T$.  
By Theorem~\ref{3-tree-display}, there are flower
vertices $v_A$, $v_B$ and $v_C$ of $T$ that display flowers 
equivalent to $\AAA$, $\BB$ and $\CC$. 
By Lemma~\ref{k-wild-common}, $\fcl(A_i)\cap W=\emptyset$
for all $i\in\{1,2,\ldots,k-2\}$, and symmetric conclusions
hold for the analogous petals of $\BB$ and $\CC$.
It follows that, by moving
to equivalent flowers, we may assume that
it is precisely $\AAA$, $\BB$ and $\CC$ that are displayed
by the vertices $v_A$, $v_B$ and $v_C$. We conclude that
$L\subseteq W$. As $|L|\geq 2$, and $W$ is a triangle,
Lemma~\ref{wheel2.5} implies that we may assume that
$L=\{a,b,c\}$.
As $L$ is sequential, the vertex adjacent to the vertex of
$T$ labelled by $L$ is a flower vertex $v$, 
and as such, displays a tight flower of order at least three.
Indeed, one readily checks that
the partition displayed by $v$ is $(A,B,C,W)$. But then,
by Lemma~\ref{2-3-4-remove}, $M\ba a$ is $k$-coherent,
contradicting the fact that $W$ is $k$-wild.
\end{proof}

Next we show that feral elements are not peripheral.
Recall that a set $A\subseteq E(M)$ is 
{\em cohesive} if $E(M)-A$ is 
fully closed. We denote by $\coh(X)$ 
the set $X-\fcl(E(M)-X)$. Evidently
$\coh(X)$ is the unique maximal 
cohesive set contained in $X$. Note
that the $3$-separating set $X$ is non-sequential 
if and only if $\coh(X)\neq \emptyset$.

\begin{lemma}
\label{feral-internal}
If $f$ is a feral element of the $k$-coherent matroid $M$,
the $f$ is not peripheral.
\end{lemma}

\begin{proof}
Let $f$ be a feral element of $M$. 
Then up to duality there is a feral display 
$(P_1,P_2,\ldots,P_m)$ and $(Q_1,Q_2,\ldots,Q_k)$
for $f$. For this pair of partitions there is an
$i\in\{2,3,\ldots,m-1\}$ such that properties 
(i)---(x) of a feral display hold. By possibly
reversing the order of indices of $(P_1,P_2,\ldots,P_m)$,
we may assume that $i\geq 3$. 

Assume that $f$ is peripheral.
Then there is a $3$-tree $T$ for $M$ with a leaf bag $L$
such that $f\in L$.

\begin{sublemma}
\label{feral-internal-1}
$L$ is non-sequential and $f\in\coh(L)$.
\end{sublemma}

\subproof
Assume that $L$ is sequential. If $|L|\geq 4$,
then $M\ba f$ is $k$-coherent by Corollary~\ref{remove-sequential}.
Certainly $f$ is not in a triangle or triad. Thus
$|L|=2$. As $L$ is sequential, it follows from the definition
of 3-tree that $L$ is a petal of a tight flower of order at least
three in $M$. Only swirl-like and spike-like flowers can have
2-element tight petals, so the flower must be swirl-like or
spike-like. By Corollary~\ref{2-elt-win} and the fact that
$f$ is not in a triangle or a triad we again deduce that
$M\ba f$ or $M/f$ is $k$-coherent. Thus $L$ is non-sequential.

Assume that $f\notin\coh(L)$. By considering the full 
closure of $E(M)-L$, we see that there is a non-sequential
3-separating set $L'\subseteq L$ such that $f$ is in
either the guts or coguts of $L'$ meaning that either 
$M\ba f$ or $M/f$ is not $3$-connected. Hence $f\in\coh(L)$.
\end{proof}

By property (vi) of feral display,
$(Q_1\cup Q_2\cup\{f\},Q_3,\ldots,Q_k)$ is a swirl-like
flower of $M$. By Theorem~\ref{3-tree-display}, there
is a vertex $v$ of $T$ that displays an equivalent 
flower. If $f\in\fcl(Q_i)$ for some $i\in\{3,4,\ldots,k\}$,
then we contradict the fact that $f\in\coh(L)$.
Thus $\coh(L)\subseteq Q_1\cup Q_2\cup\{f\}$.
By property (vii) of feral display,
$(P_2,P_3,\ldots,P_i,P_{i+1}\udots P_m\cup P_1\cup\{f\})$
is a swirl-like flower of order $i$ in $M$.
Arguing as above we deduce that 
$\coh(L)\subseteq P_{i+1}\udots P_m\cup Z_1\cup Z_2$.

Let $K=\coh(L)-\{f\}$. Then $\lambda_{M\ba f}(K)=2$
and $f\in\cl_M(K)$. We now consider the location of 
$K$ relative to the flower $(P_1,P_2,\ldots,P_m)$
of $M\ba f$. If $K\subseteq \fcl_{M\ba f}(P_i)$
for some $i\in\{1,2,\ldots,m\}$, then $M$
is $k$-fractured. Hence $K$ is equivalent to a 3-separation
displayed by $(P_1,P_2,\ldots,P_m)$.
As $Q_3\cup Q_4\udots Q_k\subseteq P_1$,
we see that $K$ is equivalent to the union of a 
consecutive subset of petals in $(P_{i+1},P_{i+2},\ldots,P_m)$.
If there is more than one petal in this set, we
contradict the fact that $L$ is a peripheral 
set of $M$ as in this case the flower 
$(P_1\udots P_i\cup\{f\},P_{i+1},P_{i+2},\ldots,P_m)$
of $M$ has order at least three and 
needs to be displayed by $T$. On the other hand,
we have already shown that there cannot be only one petal in the set.
The lemma follows from this contradiction.
\end{proof}

By Corollary~\ref{upcor} there is always an element
of the full closure of a peripheral set that can be
removed to preserve $k$-coherence. The next lemma
strengthens that outcome in the non-sequential case.

\begin{lemma}
\label{wheel3.5}
Let $L$ be a fully-closed peripheral set of the $k$-coherent matroid 
$M$ such that either $L=E(M)$ or $L$ is a non-sequential
$3$-separating set.  Let $x$ be an element of $\coh(L)$.
\begin{itemize}
\item[(i)] If $M$ is not a wheel or a whirl, and  
$x$ is in a triangle $T$, then $\fcl(T)\subseteq L$
and there is an element $z\in \fcl(T)$ such that either $M\ba z$
or $M/z$ is $k$-coherent.
\item[(ii)] If $x$ is not in a triangle or a triad, then  
either $M\ba x$ or $M/x$ is $k$-coherent.
\end{itemize}
\end{lemma}

\begin{proof}
Say that $x$ is in a triangle $T$. Then $|T\cap \coh(L)|\geq 2$ and 
$\fcl(T)\subseteq L$. Part (i) now follows straightforwardly from
Corollary~\ref{remove-sequential}, 
Lemma~\ref{wheel-wild} and the fact that 
$M$ is not a wheel or a whirl. 

Assume that $x$ is not in a triangle or a triad.
Say $M/x$
is not $3$-connected. Then there is a 3-separation $(A\cup\{x\},B)$
of $M$ with $x\in\cl(A),\cl(B)$. Assume that 
$(A\cup\{x\},B)$ is sequential. Then we may
assume that $A\cup\{x\}$ is sequential.
If $|A|=2$, then $x$ is in a triangle. Hence $|A\cup\{x\}|\geq 4$.
By Bixby's Lemma and the fact that $x$ is in no triads, $M\ba x$ is
$3$-connected. Then, by Corollary~\ref{remove-sequential},
$M\ba x$ is $k$-coherent. It follows that $(A\cup\{x\},B)$
is non-sequential.

If $L=E(M)$, then a
3-tree for $M$ consists of a single vertex and $M$ has no
non-sequential $3$-separations. Thus $M/x$, and similarly $M\ba x$
are $3$-connected. By Lemma~\ref{feral-internal}, one of $M\ba x$
or $M/x$ is $k$-coherent as required.

We may now assume that $\lambda(L)=2$.
Let $K=E(M)-\coh(L)$. Then $K$ is fully closed and $(K,\coh(L))$
is a non-sequential 3-separation.
If $M/x$ is not $3$-connected, then $x$ is in the guts of 
a non-sequential $3$-separation $(A\cup\{x\},B)$. But 
$x\in\fcl(B)$ and $x\notin\fcl(K)$ so $(A\cup\{x\},B)$ is not
equivalent to $(K,\coh(L))$. A $3$-separation 
equivalent to $(A\cup\{x\},B)$
must be displayed in a 3-tree for $M$. As $\coh(L)$ 
is a peripheral set, and 
$K$ is fully closed, it is easily seen that we may assume that
$B\subseteq K$. But now $x\in \cl(B)$ so $x\in\cl(K)$ contradicting
the fact that $K$ is fully closed. 
Therefore both $M/x$ and $M\ba x$
are $3$-connected. Again by Lemma~\ref{feral-internal}, we see that
one of $M/x$ or $M\ba x$ is $k$-coherent.
\end{proof}

\section{Extending a $k$-coherent Matroid}

Let $x$ be an element of the matroid $M$ such that
$M\ba x$ is $3$-connected. If $M$ is not 3-connected,
then it is easily seen that $x$ is either a loop, a coloop,
or is in a parallel pair in $M$. 
The situation when $M\ba x$
is $k$-coherent, but $M$ is not is
a little more complicated. This section gives some 
straightforward lemmas describing the structures that
arise in this and related situations.



We omit the routine proof of the next lemma.

\begin{lemma}
\label{loose1}
Let $M$ be a $3$-connected matroid and let $(l_1,\ldots,l_n)$
be a maximal fan of loose elements between a pair of 
petals in a swirl-like
flower of $M$ of order at least $3$. Assume that $l_1$ is a guts element.
\begin{itemize}
\item[(i)] $M\ba l_1$ is $k$-coherent if and only if $M$ is.
\item[(ii)] If $i$ is odd and $i>1$, then $M\ba l_i/l_{i-1}$ is 
$3$-connected. Moreover $M\ba l_i/l_{i-1}$ is $k$-coherent if and
only if $M$ is.
\end{itemize}
\end{lemma}

Let $M$ be a 3-connected matroid with a maximal $k$-fracture
$\PP$. Then $M$ is {\em uniquely} $k$-fractured {\em by\ } $\PP$
if every $k$-fracture $\QQ$ of $M$ has the property that
$\QQ\less \PP$. Note that if $M$ is uniquely fractured
by $\PP$ and $\PP$ has order $k$, then every $k$-fracture of
$M$ is equivalent to $\PP$. 

Recall the definition of 
quasi-flower. Note that if a quasi-flower has exactly
one $1$-element petal, then this petal is contained in either
the closure or coclosure of either of its adjacent petals,
so it is quite properly regarded as a loose petal.

\begin{lemma}
\label{gain-coherence}
Let $M$ be a $3$-connected matroid with an element
$e$ such that $M\ba e$ is $k$-coherent. If $M$ is not
$k$-coherent, then the following hold.
\begin{itemize}
\item[(i)] $M$ is  uniquely $k$-fractured by a maximal
swirl-like flower $\PP$ of order $k$.
\item[(ii)] If $(P_1,P_2,\ldots,P_k)$ is a $k$-fracture
of $M$ with $e\in P_1$, then $P_1-\{e\}$ is a loose 
petal in the swirl-like quasiflower
$(P_1\ba \{e\},P_2,\ldots,P_k)$ of $M$.
This quasiflower has order $k-1$.
\end{itemize}
\end{lemma}

\begin{proof}
Up to labels $M$ has a $k$-fracture 
$(P_1,\ldots,P_k)$, where $e\in P_1$.
But $(P_1-\{e\},P_2,\ldots,P_k)$ cannot 
be a $k$-fracture of $M\ba e$.
The latter part of the lemma follows easily from this observation.
Say that 
$(Q_1,Q_2,\ldots,Q_k)$ is another $k$-fracture of
$M$. By Lemma~\ref{in-petal}, up to labels,
equivalence and symmetry, we may assume that
$e\in Q_1$, that $Q_1$ is fully closed and that
$Q_1$ contains all but one petal of $(P_1,P_2,\ldots,P_k)$.
By the above, $Q_1-\{e\}$ is a set of loose elements of
the flower $((Q_1\cup Q_2)-\{e\},Q_3,\ldots,Q_k)$ of $M\ba e$.
But it is now easily seen that this is not possible.
Thus $M$ is uniquely fractured by $\PP$ and (i) holds.
\end{proof}

Viewed from the perspective of moving from $M\ba e$ to $M$, we
have

\begin{lemma}
\label{in-closure}
Assume that $M\ba e$ is $k$-coherent and that $M$ is
$3$-connected and $k$-fractured. If $\PP$ is a flower
of $M\ba e$ of order at least three, then, 
$e\in\cl(\hP)$ for some petal $P$ of $\PP$.
\end{lemma} 

\begin{proof}
Let $(P_1,\ldots,P_k)$ be a $k$-fracture of $M$, where $e\in P_1$.
By Lemma~\ref{gain-coherence}, $((P_1\cup P_2)-\{e\},P_3,\ldots,P_k)$
is a maximal flower in $M\ba e$ and it follows that 
$e\in\cl((P_1\cup P_2)-\{e\})$ and the lemma holds in this case.

Say that $\QQ=(Q_1,Q_2\ldots,Q_m)$ is 
another flower in $M$ of order at least
three. By Lemma~\ref{in-petal}, we may assume up to labels 
in $\QQ$ that either $P_k\subseteq \hQ_1$ or $P_2\subseteq \hQ_1$.
Assume without loss of generality that the latter case holds.
Then $\hP_2\subseteq \hQ_1$ and $e\in \cl(\hP_2)$ so that
$e\in\cl(\hQ_1)$ as required.
\end{proof}

We omit the routine proof of the next lemma.

\begin{lemma}
\label{lose-coherence}
Assume that $M\ba e$ is $k$-coherent and that $M$ is $3$-connected
and $k$-fractured. Then there is a swirl-like quasi-flower
$(P_1,L,P_2,\ldots,P_{k-1})$ of $M\ba e$, where 
$L$ is a nonempty loose petal,
such that 
the following hold.
\begin{itemize}
\item[(i)] $(P_1,L\cup \{e\},P_2,\ldots,P_{k-1})$ is a 
$k$-fracture of $M$.
\item[(ii)] $e\in\cl(P_1\cup L)$ and $e\in\cl(P_2\cup L)$.
\item[(iii)] If $|L|>1$, then $e\in\cl(L)$.
\end{itemize}
\end{lemma}

The next lemma is unsurprising. If we block every
$k$-fracture we would expect to become $k$-coherent.

\begin{lemma}
\label{keep-fracture}
Assume that $M$ and $M\ba e$ are $3$-connected and $k$-fractured.
Then there is a $k$-fracture $(P_1,P_2,\ldots,P_n)$ of
$M\ba e$ such that 
$e\in\cl_M(\hP_i)$ for some $i\in\{1,2,\ldots,n\}$.
\end{lemma}

\begin{proof}
Let $(Q_1\cup\{e\},Q_2,\ldots,Q_n)$ be a maximal $k$-fracture
of $M$. If the lemma fails, it must be the case that $n=k$
and that $Q_1$ is a set of loose elements of the maximal
flower $(Q_1\cup Q_2,Q_3,\ldots,Q_k)$ of $M\ba e$.
Note that $e\in\cl_M(\hQ_2)$ and $e\in\cl_M(\hQ_k)$.
As $M\ba e$ is $k$-fractured, it has a maximal $k$-fracture
$(P_1,P_2,\ldots,P_n)$. By Lemma~\ref{in-petal}, there is an
$i\in\{1,2,\ldots,n\}$ such that either $Q_2\subseteq \hP_i$ or 
$Q_k\subseteq \hP_i$. In either case it follows that
$e\in\cl_M(\hP_i)$ as required.
\end{proof}

Note that Lemma~\ref{keep-fracture} would fail if we were to
insist in the statement that the fracture $(P_1,\ldots,P_n)$ was  
{\em maximal}. As an easy corollary of 
Lemma~\ref{keep-fracture}, we have

\begin{corollary}
\label{unique-fracture}
Assume that $M$ and $M\ba e$ are $3$-connected and that, 
$M\ba e$ has a unique maximal $k$-fracture
$(P_1,\ldots,P_k)$. Then $M$ is $k$-fractured if and only if
$e\in\cl(\hP_i)$ for some $i\in\{1,\ldots,k\}$.
\end{corollary}

\begin{lemma}
\label{2-element-fracture}
Assume that $M$ is $3$-connected and that
$(\{e,f\},P_2,\ldots,P_k)$ is a maximal flower that
uniquely $k$-fractures 
$M$, where 
$\{e,f\}$ is fully closed. Then $M\ba e$ is $k$-coherent.
\end{lemma}

\begin{proof}
Assume that $M\ba e$ is not $k$-coherent. 
Then, by Lemma~\ref{keep-fracture},
there is a $k$-fracture $(Q_1,Q_2,\ldots,Q_n)$ of $M\ba e$ such that
$e\in\cl(Q_1)$. But then $(Q_1\cup\{e\},Q_2,\ldots,Q_n)$ 
is a $k$-fracture of
$M$. As $M$ has a unique $k$-fracture and $\{e,f\}$ is fully closed,
$Q_1\cup\{e\}=\{e,f\}$ so that $Q_1=\{f\}$ 
contradicting the assumption that 
$Q_1$ is a petal of a flower of $M\ba e$.
\end{proof}

We omit the easy proof of the next lemma.

\begin{lemma}
\label{split-petal}
Let $(A,B)$ be a $3$-separation of the $k$-coherent matroid
$M$, where $a\in A$, and $|A|\geq 4$. If $M\ba a$ is 
$3$-connected and $k$-fractured, then there is a $k$-fracture
$(P_1,P_2,\ldots,P_n)$ of $M\ba a$ such that, 
for some $i\in \{2,3,\ldots,n-1\}$,
$A=P_1\cup P_2\udots P_i$, and $(A,P_{i+1},P_{i+2}\ldots,P_n)$ is a 
swirl-like flower of $M$.
\end{lemma}

Note that the flower $(A,P_{i+1},P_{i+2}\ldots,P_n)$ of 
Lemma~\ref{split-petal}
may be trivial in the sense that is has only two petals.

\chapter{$k$-skeletons}
\label{k-skeletons}

It was noted in the introduction that the underlying
cause of inequivalent representations of a matroid is
that elements may have ``freedom''. Loosely speaking
 $k$-skeletons  are $k$-coherent matroids whose elements 
have maximal freedom. The next section of this chapter
makes these notions precise. The remaining sections
of the chapter
examine freedom in matroids and the connection with 
$k$-coherence in more detail. The last two sections 
describe certain
structures that turn out to be of importance. Life would have
been easier if we could have avoided worrying about these
structures, but, sadly, this seems not to be the case.

\section{Clones, Fixed Elements, and $k$-skeletons}

Let $M$ be a matroid. Elements $e$ and $f$ of $M$ are
{\em clones\ } 
\index{clones}
if swapping the labels of $e$ and $f$ is an
automorphism of $M$. A {\em clonal set\ }
\index{clonal set} 
of $M$ is a set of
elements every pair of which are clones. A 
{\em clonal class\ }
\index{clonal class} 
of $M$ is a maximal clonal set.
An element
$z$ of $M$ is {\em fixed\ } 
\index{fixed element}
in $M$ if there is no single-element
extension of $M$ by an element $z'$ in which 
$z$ and $z'$ are independent clones. Dually, an element
$z$ in $M$ is {\em cofixed\ } 
\index{cofixed element} 
in $M$ if it is fixed in $M^*$.
Note that if $z$ already has a clone, say $x$, and 
$\{x,z\}$ is independent, then $z$ is not fixed since
we can add a new element $z'$ freely on the line through
$x$ and $z$.

Let  $k\geq 5$ be an integer, and $M$ be a
$k$-coherent matroid. Then $M$ is a 
{\em $k$-skeleton} 
\index{$k$-skeleton}
if the following hold.
\begin{itemize}
\item[(i)] $M$ is not a wheel or a whirl of rank at least 3.
\item[(ii)] If $x$ is fixed in $M$, then $M\ba x$
is not $k$-coherent.
\item[(iii)] If $x$ is cofixed in $M$ then $M/x$ is not
$k$-coherent.
\end{itemize}

As with $k$-coherence, in any unexplained context, when referring
to a $k$-skeleton, we always assume that $k$ is an integer
greater than four. Note that condition (i) is required
in the definition
simply because wheels and whirls 
vacuously satisfy (ii) and (iii). 

It is shown in Lemma~\ref{field1} that the number of inequivalent
representations of a $k$-coherent matroid over a finite field
$\mathbb F$ is bounded above by the maximum of the number of
inequivalent $\mathbb F$-representations of its $k$-skeleton
minors.
The importance of $k$-skeletons in studying inequivalent
representations of $k$-coherent matroids is due to this fact.

\section{Freedom and Cofreedom}

In this section we develop further material related to
freedom in matroids. Most of the results here are either
straightforward or are proved in \cite{short-gf5,totally-free}.

A flat $F$ of the matroid $M$ is {\em cyclic\ }
\index{cyclic flat} 
if, for each
element $e\in F$, there is a circuit $C$ such that 
$e\in C\subseteq F$. It is easily seen that $F$ is a
cyclic flat of $M$ if and only if $E(M)-F$ is a cyclic
flat of $M^*$. The next result is straightforward.

\begin{proposition}
\label{cyclic-clones}
Elements $e$ and $f$ of a matroid $M$ are clones if and
only if $e$ and $f$ are contained in the same cyclic flats.
\end{proposition}

Let $e$ and $f$ be elements of the matroid $M$. Then
$e$ is {\em freer\ }
\index{freer} 
than $f$, denoted 
$e\more f$,   if every cyclic flat containing
$e$ also contains $f$. The {\em freedom\ } 
\index{freedom of an element} 
of an element
$e$ of $M$ is the maximum size of an 
independent clonal class containing $e$ amongst all 
matroids containing $M$ as a restriction, that os, amongst all
extensions of $M$. 
We denote the freedom of $e$ in $M$ by $\fr_M(e)$, or 
$\fr(e)$ if the matroid $M$ is clear.
This maximum does not exist if and only 
if $e$ is a coloop of $M$ in which case $\fr(e)$ is infinity.
A loop has freedom $0$ and an element is fixed if and only if
it has freedom at most 1.

The notion of freedom in matroids was introduced by 
Duke \cite{duke}.
His definition was different from that given above, but it is shown
in \cite[Lemma~2.8]{short-gf5} that Duke's definition is equivalent
to ours. The next two lemmas are also proved in \cite{short-gf5}.

\begin{lemma}
\label{clear}
Let $a$ and $b$ be elements of the matroid $M$ such that
$a\more b$. Then $\fr(a)\geq \fr(b)$. Moreover,
either $a$ and $b$ are clones or $\fr(a)>\fr(b)$.
\end{lemma}

\begin{lemma}
\label{easy-closure}
Let $a$ and $b$ be elements of the matroid $M$.
Then $a\more b$ if and only if the following holds:
for all $X\subseteq E(M)-\{a,b\}$,
if $a\in\cl(X)$, then $b\in\cl(X)$.
\end{lemma}

The next lemma is elementary.

\begin{lemma}
\label{free-minor}
Let $a$ and $b$ be elements of the matroid $M$,
and let $N$ be a minor of $M$ whose ground set contains
$a$ and $b$. If $a$ is freer than $b$ in $M$, then
$a$ is freer than $b$ in $N$.
\end{lemma}

\begin{lemma}
\label{free-at-last}
Let $a$ and $b$ be elements of the matroid $M$. Then
\begin{itemize}
\item[(i)] $\fr_{M\ba a}(b)\geq \fr_M(a)$.
\item[(ii)] $\fr_{M/a}(b)\geq \fr_M(a)-1$, and, if 
$\fr_{M/a}(b)=\fr_M(a)-1$, then $a\less b$ in $M$.
\end{itemize}
\end{lemma}

\begin{proof}
The lemma is trivial if $\fr(b)=0$.
Assume that $b$ has freedom $k\geq 1$ and let
$N$ be a matroid that extends $M$ in which
$b$ belongs to an independent clonal set $B$ of size
$k$. It is easily seen that such an extension exists
with the property that $a\notin B$.
As $B$ is an independent clonal set in $N\ba a$, part (i)
holds.

Consider $N/a$. In this matroid $B$ is a clonal set and contains
an independent subset containing $a$ of size $k-1$,
so $\fr_{M/a}(b)\geq \fr_M(b)-1$. Assume that
$\fr_{M/a}(b)<\fr_M(b)$. Then $B$ is not independent
in $N/a$. Thus $a\in\cl_N(B)$. But every cyclic flat
of $N$ containing $b$ also contains $B$, so every cyclic 
flat of $N$ containing $b$ contains $a$. Hence 
$b$ is freer than $a$ in $N$, and by Lemma~\ref{free-minor}
$b\more a$ in $M$.
\end{proof}

The {\em cofreedom\ }
\index{cofreedom of an element} 
of an element $e$ of $M$, denoted 
$\fr^*(e)$, is the freedom of
$e$ in $M^*$. Note that $e\more f$ in $M^*$ if and only if
$f\more e$ in $M$. This is a consequence of the fact that the
cyclic flats of $M^*$ are the complements of the 
cyclic flats of $M$. The following lemma is the dual of
Lemma~\ref{free-at-last}.

\begin{lemma}
\label{cofree-at-last}
Let $a$ and $b$ be elements of the matroid $M$. Then
\begin{itemize}
\item[(i)] $\fr^*_{M/a}(b)\geq \fr^*_M(a)$.
\item[(ii)] $\fr^*_{M\ba a}(b)\geq \fr^*_M(a)-1$, and, if 
$\fr^*_{M\ba a}(b)=\fr^*_M(a)-1$, then $a\more b$ in $M$.
\end{itemize}
\end{lemma}

The next two results are immediate corollaries of 
Lemmas~\ref{free-at-last} and \ref{cofree-at-last}.
We apply them frequently.

\begin{corollary}
\label{born-free}
Let $a$ and $b$ be elements of the matroid $M$.
\begin{itemize}
\item[(i)] If $b$ is fixed in $M/a$, but not in $M$, then
$a\less b$ in $M$ and $\fr_M(b)=2$. Moreover, either $a$ and $b$
are clones in $M$, or $a$ is fixed in $M$.
\item[(ii)] If $b$ is cofixed in $M\ba a$ but not in $M$,
then $a\more b$ in $M$ and $\fr^*_M(b)=2$. Moreover, either $a$
and $b$ are clones in $M$, or $a$ is cofixed in $M$.
\end{itemize}
\end{corollary}

Elements $a$ and $b$ of $M$ are {\em comparable\ }
\index{comparable elements}
if either $a\more b$ or $b\more a$; otherwise they
are {\em incomparable}.
\index{incomparable elements}

\begin{corollary}
\label{im-free}
Let $a$ and $b$ be incomparable elements of the matroid $M$.
\begin{itemize}
\item[(i)] If $b$ is fixed in $M\ba a$ or $M/a$,
then $b$ is fixed in $M$. 
\item[(ii)] If $b$ is cofixed in $M\ba a$ or $M/a$,
then $b$ is cofixed in $M$.
\end{itemize}
\end{corollary}

The next lemma is  \cite[Lemma~2.13]{short-gf5}.

\begin{lemma}
\label{distinguish}
Let $a,b$ and $e$ be elements of the matroid $M$ such
that $a$ and $b$ are clones and have freedom $2$ in $M\ba e$.
If $a$ and $b$ are not clones in $M$, then either 
$a$ or $b$ is fixed in $M$.
\end{lemma}

At times local connectivity provides a useful way of 
bounding the freedom of an element in a matroid.

\begin{lemma}
\label{freedom1}
Let $X$ and $Y$ be disjoint sets of elements of a matroid
$M$ and suppose that 
$a\in E(M)-(X\cup Y)$. If
$a\in\cl(X),\cl(Y)$, then $\fr(a)\leq\sqcap(X,Y)$.
\end{lemma}

\begin{proof}
Let $N$ be an extension of $M$ in which $a$ belongs to 
an independent clonal set $A$ of size $\fr(a)$.
By Lemma~\ref{easy-closure}, $A\subseteq \cl(X)$
and $A\subseteq \cl(Y)$. Moreover,
$r(\cl_N(X)\cap \cl_N(Y))\leq r_N(X)+r_N(Y)
-r_N(\cl_N(X)\cup\cl_N(Y))=r_M(X)+r_M(Y)-r_M(X\cup Y)=
\sqcap_M(X,Y)$.
We now have $\fr_M(a)=|A|\leq r(\cl_N(X)\cap \cl_N(Y))$,
and the lemma follows.
\end{proof}

The next result follows easily from Lemma~\ref{freedom1}.

\begin{corollary}
\label{fix-cofix}
Let $A,B$ be disjoint subsets of a matroid $M$.
\begin{itemize}
\item[(i)] If $x\in\cl(A),\cl(B)$, and $A,B$ are skew in 
$M/x$, then $x$ is fixed in $M$.
\item[(ii)] If $x\in\cl^*(A),\cl^*(B)$, and $A,B$ are coskew
in $M\ba x$, then $x$ is cofixed in $M$.
\end{itemize}
\end{corollary}

A useful consequence of Lemma~\ref{fix-cofix}(ii) is

\begin{lemma}
\label{swirl-cofix}
Let $a$ be an element of the $3$-connected matroid $M$. If $a$
blocks non-adjacent petals of a swirl-like flower of $M\ba a$,
then $a$ is cofixed in $M$.
\end{lemma}

There are connections between freedom and connectivity. Loosely
speaking, the freer an element $b$ of $M$ is, 
the less likely it is that connectivity will be damaged when 
$b$ is contracted from $M$. If contracting $b$ does damage
connectivity, and $a\less b$, then contracting 
$a$ should also damage connectivity. The next lemma
has easy generalisations. We focus on the case that is 
useful for us.

\begin{lemma}
\label{free-for-all}
Let $b$ be an element of the $3$-connected matroid $M$. Assume
that $M/b$ is not $3$-connected. Then the following hold.
\begin{itemize}
\item[(i)] $\fr(b)\leq 2$.
\item[(ii)] If $a\in E(M)$, and $a\less b$, then
$M/a$ is not $3$-connected.
\end{itemize}
\end{lemma}

\begin{proof}
Let $(P,Q)$ be a 2-separation of $M/b$.
As $M$ is $3$-connected, $b$ coblocks this 2-separation so 
that $b\in\cl_M(P)$ and $b\in\cl_M(Q)$.
But $\sqcap(P,Q)=2$, so by Lemma~\ref{freedom1}, 
$\fr_M(b)\leq 2$, and (i) holds. 
Say $a\less b$. Assume that $a\in P$.
Observe that $a\in\cl((P-\{a\})\cup\{b\})$.
As $a\less b$, we have $a\in\cl(Q)$.
Hence $((P-\{a\})\cup \{b\},Q)$ is a 2-separation
of $M/a$, and (ii) holds.
\end{proof}

We write $a\prec b$ if 
$a\less b$ in $M$, but $a$ and $b$ are not clones in $M$.
By Lemma~\ref{clear}, $a\prec b$ if $a\less b$ 
and $\fr(a)<\fr(b)$. The next result is a routine corollary
of Lemma~\ref{free-for-all} and we omit the proof.

\begin{corollary}
\label{freedom2}
Let $a$ and $b$ be elements of the $3$-connected matroid
$M$ such that $a\prec b$. 
\begin{itemize}
\item[(i)] If $b$ is in a triangle, then there is a
triangle containing $a$ and $b$, and $a$ is fixed in $M$.
\item[(ii)] If $\si(M/b)$ is not
$3$-connected, then $a$ is fixed in $M$ 
and $\si(M/a)$ is not $3$-connected.
\end{itemize}
\end{corollary}

As elements in triangles have freedom at most 2, we do not
need much extra information to fix them.

\begin{lemma}
\label{freedom3}
If $a$ is in a triangle $T$ of the matroid 
$M$ and there is a cyclic flat
of $M$ that contains $a$ but not $T$, then $a$ is fixed in $M$.
\end{lemma}

\begin{proof}
Say $F$ is a cyclic flat of $M$ containing $a$ but not $T$.
Observe that $a\in\cl(F-\{a\})$, $a\in\cl(T-\{a\})$,
and $\sqcap(F-\{a\},T-\{a\})=1$, so that the lemma
holds by Lemma~\ref{freedom1}.
\end{proof}

The next lemma sees a single application in Chapter~\ref{taming}.

\begin{lemma}
\label{freedom21}
Let $\{a,b,c\}$ be a clonal triangle of the matroid $M$ and let
$e$ be an element of $E(M)-\{a,b,c\}$.
\begin{itemize}
\item[(i)] If $e$ is fixed in $M$, then $e$ is fixed in $M\ba a$.
\item[(ii)] If $e$ is cofixed in $M$, then $e$ is cofixed in $M\ba a$.
\end{itemize}
\end{lemma}

\begin{proof}
Consider (i). Let $N=M\ba a$ and assume that $e$ is not fixed
in $M\ba a$. Let $N'$ be a matroid obtained by independently
cloning $e$ by $e'$. Let $M'$ be a  matroid obtained from
$N'$ by freely
placing the point $a$ in the span of $\{b,c\}$. Consider $M'\ba e'$.
We have the following: $M\ba a=M'\ba e',a$, the set
$\{a,b,c\}$ is a triangle in $M'\ba e'$, and $a$ is freely placed
on the line $\{a,b,c\}$. It follows straightforwardly that
 $M'\ba e'=M$. 

Assume that $e$ is fixed in $M$. Then $e$ is fixed in $M'$.
It follows that there is a cyclic flat $F$ of $M'$ containing $e$
but not $e'$. If $F$ does not contain $a$, then we contradict the
fact that $e$ and $e'$ are clones in $N'$. Hence $F$ contains $a$.
But as $a$ is freely placed in $\{a,b,c\}$, we see that
$a$ is freer than both $b$ and $c$ in $M'$. Therefore
$b$ and $c$ are contained in $F$. Let $C$ be a circuit
such that $e\in C$ and $C\subseteq F$. 
If $a\in C$, then we may apply circuit elimination to $C$ 
and $\{a,b,c\}$ to 
obtain a circuit $C'$ containing $e$ that does not
contain $a$. Thus we may assume that $C$ does not contain 
$a$. But then $\cl_{N'}(C)$ is a cyclic flat of $N'$
that contains $e$ but not $e'$ contradicting the fact that
$\{e,e'\}$ is a clonal pair in $N'$. 

Part (ii) follows by a similar argument, but this time by 
coindependently cocloning $e$ by $e'$ in $M\ba a$.
\end{proof} 

As always the connection with flowers is 
important to us.

\begin{lemma}
\label{freedom4}
Let $M$ be a $3$-connected matroid and 
let $\{a,b\}$ be a $2$-element tight petal of 
a swirl-like or spike-like flower with at
least four petals. If $a$ is not fixed in $M$, then the 
following hold.
\begin{itemize}
\item[(i)] $\fr(a)=2$.
\item[(ii)] $b\less a$.
\item[(iii)] If $M'$ is obtained by independently cloning 
$a$ by $a'$, then $\{a,a',b\}$ is a triangle in $M'$.
\end{itemize}
\end{lemma}

\begin{proof}
Say that flower is $(P_1,\{a,b\},P_3,\ldots,P_m)$. 
Let $M'$ be a 
matroid obtained by independently cloning $a$ by $a'$. As $\{a,b\}$
is a tight petal, $a\in\cl(P_1\cup\{b\})$ so that 
$a'\in\cl(P_1\cup\{b\})$, that is, $a'\in\cl(P_1\cup\{a,b\})$
and, similarly, $a'\in\cl(P_3\cup\{a,b\})$. Thus, by
Lemma~\ref{modular},  $\{a,a',b\}$ is a triangle in $M'$. 
Hence (iii) holds. Parts (i) and (ii) follow easily.
\end{proof}

Recall that  the loose elements in a swirl-like
flower can be partitioned
into guts and coguts elements in a canonical way. 
The next lemma is perhaps
a terminological surprise, in that
it is quite easy for a loose element to be fixed.

\begin{lemma}
\label{loose2}
Let $M$ be a $3$-connected matroid and $l$ be a loose element
of a swirl-like or spike-like flower $\FF$ of $M$ 
of order at least $3$. Then the following hold.
\begin{itemize}
\item[(i)] If $\FF$ is spike-like and $l$
is the tip, then $l$ is fixed in $M$.
\item[(ii)] If $\FF$ is spike-like and $l$
is the cotip, then $l$ is cofixed in $M$.
\item[(iii)] If $\FF$ is swirl-like and $l$ is a loose
guts element, then $l$ is fixed in $M$.
\item[(iv)] If $\FF$ is swirl-like and $l$ is a
loose coguts element, then $l$ is cofixed in $M$.
\end{itemize}
\end{lemma}

\begin{proof} 
Assume that either (i) or (iii) holds. Then it is clear
that $M$ has a swirl-like or spike-like
quasi-flower $(P_1,\{l\},P_2,P_3)$
where $l\in \cl(P_1)$, $l\in \cl(P_2)$, and $|P_3|\geq 2$. 
In this case $\sqcap(P_1,P_2)=1$ and the lemma follows by
Lemma~\ref{freedom1}.
\end{proof}

\section{Freedom and $k$-coherence}

In this section we collect some lemmas that deal with the
relationship between the relative freedom of elements in a 
matroid and $k$-coherence.

\begin{lemma}
\label{delete-cofix}
Let $a$ and $b$ be elements of the $k$-coherent matroid $M$, where
$a\less b$. If $M\ba a$ is not $k$-coherent, then $b$ is
cofixed in $M\ba a$.
\end{lemma}

\begin{proof}
Assume that $M\ba a$ is not $k$-coherent.
Say that $M\ba a$ is not $3$-connected. Then $M\ba a$
has a $2$-separation $(X,Y)$. Assume that
$b\notin\cl^*_{M\ba a}(X-\{b\})$.
Then $b\in\cl_{M\ba a}(Y-\{b\})$, that is, $b\in\cl_M(Y-\{b\})$.
As $a\less b$, we have $a\in\cl_M(Y)$, contradicting the fact
that $a$ coblocks $(X,Y)$. Therefore 
$b\in\cl^*_{M\ba a}(X-\{b\})$ and $b\in\cl^*_{M\ba a}(Y-\{b\})$
so that $b$ is cofixed in $M\ba a$ by the dual of 
Lemma~\ref{freedom1}.

Assume that $M\ba a$ is $3$-connected. 
Let  $(P_1,P_2,P_3,\ldots,P_k)$ be a $k$-fracture of 
$M\ba a$, where $b\in P_2$. If $b$ is a coloop of neither
$M|(P_1\cup P_2)$ nor $M|(P_1\cup P_2)$, then,
by Lemma~\ref{easy-closure},
$a\in\cl_M(P_1\cup P_2)$ and $a\in\cl(P_2\cup P_3)$,
so that, by Lemma~\ref{modular}, $a\in\cl_M(P_2)$,
contradicting the fact that $M$ is $k$-coherent.
Hence we may assume that $b$ is a coloop of 
$M|(P_1\cup P_2)$. But this means that $b$ is in the coguts
of $(P_1\cup P_2,P_3\udots P_k)$. In this case,
by Lemma~\ref{fine2}, $b$ is a loose coguts element
of $P_2$. By Lemma~\ref{loose2} $b$ is cofixed in 
$M\ba a$.
\end{proof}

\begin{corollary}
\label{3-comparable} 
Let Let $a$ and $b$ be elements of the 
$k$-coherent matroid $M$, where
$a\less b$. 
\begin{itemize}
\item[(i)] If $\fr(a)\geq 3$, then $M/b$ is $k$-coherent.
\item[(ii)] If $\fr^*(b)\geq 3$, then $M\ba a$ is 
$k$-coherent.
\end{itemize}
\end{corollary}

\begin{proof}
Consider (ii). Say that $\fr^*(b)\geq 3$. Then,
by Lemma~\ref{cofree-at-last}, $\fr^*_{M\ba a}(b)\geq 2$
so that $b$ is not cofixed in $M\ba a$. So by 
Lemma~\ref{delete-cofix} $M\ba a$ is $k$-coherent.
Part (i) is the dual of (ii).
\end{proof}

The next lemma is easily seen to hold for 
arbitrary flowers, but we only need it in the swirl-like case.

\begin{lemma}
\label{comparable-in-petal}
Let $\FF$ be a swirl-like flower
of the $3$-connected matroid $M$ of order at least $4$. 
If $a,b\in E(M)$,
and $a\less b$, then there is a flower equivalent to $\FF$ in
which $a$ and $b$ are in the same petal.
\end{lemma}

\begin{proof}
Say that $(P_1,P_2,\ldots,P_n)$ is a tight 
swirl-like flower equivalent
to $\FF$ that is chosen so that $a\in P_1$, $b\in P_i$, and, 
amongst all equivalent flowers, $i$ is minimal. Assume that 
$i>1$.  By the minimality
assumption $b$ is not a loose coguts element 
between $P_{i-1}$ and $P_i$.
Thus, $b\in\cl((P_i\cup P_{i+1})-\{b\})$. By Lemma~\ref{easy-closure},
$a\in\cl(P_i\cup P_{i+1})$. By Lemma~\ref{fine1},
either $a\in\cl(P_i)$ or $a\in\cl(P_{i+1})$. 
The former case
implies a contradiction.  Thus we 
may assume that $a\in\cl(P_{i+1})$.
Now $a\in P_{i+2}$ and $n\geq 4$, 
so that $i>2$ mod $n$, and $a$ and $b$ are not in
adjacent petals. However putting $a$ in $P_{i+1}$ gives an
equivalent flower in which $a$ and $b$ are in adjacent petals
and we have contradicted the minimality of the choice of $\FF$
and $i$.
\end{proof}

Knowing that elements of a $k$-coherent matroid are comparable
gives us valuable information that enables us to deduce
that certain minors are $k$-coherent. 
This point is exemplified in the next lemma. Let $(P_1,P_2,\ldots,P_n)$
be a swirl-like flower of a matroid. Recall that, when we say
that $x$ is a coguts element of $(P_i,P_{i+1})$ we mean that
$x$ is a rim element of the fan of loose elements between
$P_i$ and $P_{i+1}$.

\begin{lemma}
\label{comparable}
Let $x$ and $y$ be elements of the $k$-coherent matroid $M$ where 
$x\less y$. Assume that  $M\ba x$ is $3$-connected and $k$-fractured. Then 
\begin{itemize}
\item[(i)] $M/y$ is $k$-coherent, and
\item[(ii)] either $x$ and $y$ are clones in 
$M$ or $y$ is cofixed in $M$.
\end{itemize}
\end{lemma}

\begin{proof}
Let ${\bf F}=(P_1\cup\{y\},P_2,\ldots,P_n)$ be a maximal 
$k$-fracture of $M\ba x$. 

\begin{sublemma}
\label{comparable-1}
Up to labels $y$ is in the coguts of $(P_1\cup\{y\},P_2)$. 
Moreover $n=k$.
\end{sublemma}

\subproof
Assume the sublemma does not hold. Then we may assume that
$y$ is not in the coguts of $(P_n,P_1\cup\{y\})$ 
nor in the coguts of
$(P_1\cup\{y\},P_2)$. Thus
$y\in\cl(P_n\cup P_1)$ and $y\in\cl(P_1\cup P_2)$. Then,
as $x\less y$, we have, by Lemma~\ref{easy-closure}, 
$x\in\cl(P_n\cup P_1)$ and 
$x\in\cl(P_1\cup P_2)$ so that $x\in\cl(P_n\cup (P_1\cup\{y\}))$
and $x\in\cl((P_1\cup\{y\})\cup P_2)$. But then, by Lemma~\ref{modular},
$x\in\cl(P_1\cup\{y\})$ implying that $M$ is $k$-fractured.
Thus we may assume that $y$ is in the coguts of $(P_1\cup\{y\},P_2)$.

Now $y\in\cl(P_1\cup P_2)$ so that 
$x\in\cl(P_1\cup P_2)$ and hence
$(P_1\cup P_2\cup \{x,y\},P_3,\ldots, P_n)$
is a tight swirl-like flower in $M$. 
As $M$ is $k$-coherent $n=k$.
\end{proof}

As $y$ is in the coguts of $(P_1,P_2)$ we have
$y\in\cl^*_{M\ba x}(P_1)$ and  $y\in\cl^*_{M\ba x}(P_2)$. 
Also $y\in\cl(P_1\cup P_2)$ and,
as $y\more x$, we have $x\in\cl(P_1\cup P_2)$.

\begin{sublemma}
\label{comparable-2}
$M/y$ and $M\ba x/y$ are $3$-connected.
\end{sublemma}

\subproof
Note that $M\ba x,y$ is not $3$-connected up to series pairs (for
example, $(P_n\cup P_1,P_2\cup P_3\cup \cdots\cup P_{n-1})$ is a 
$2$-separation of $M\ba x,y$).
So, by Bixby's Lemma, 
$M\ba x/y$ is $3$-connected up to parallel pairs.
Hence $M/y$ is $3$-connected up to parallel classes.
To prove the sublemma it suffices to show that $M/y$ has no non-trivial
parallel classes. Assume otherwise.
Then $y$ is in a triangle of $M$. As $y$ is freer than $x$,
there is an element $t$ such that $\{x,y,t\}$ is a triangle of $M$.
As $\{x,y\}\subseteq\cl(P_1\cup P_2)$, we have 
$t\in\cl(P_1\cup P_2)$. By taking an
equivalent flower if necessary,
we may assume that $P_1\cup P_2\cup\{y\}$ is closed
in $M\ba x$. Up to symmetry we may assume that
$t\in P_1$. But now 
$x\in\cl(P_1\cup\{y\})$ and $(P_1\cup\{x,y\},P_2,\ldots,P_n)$ is a
$k$-fracture of $M$, contradicting the fact that $M$ is $k$-coherent. 
\end{proof}

As $M\ba x/y$ is $3$-connected, and $y$ is in the coguts of $(P_1,P_2)$
it follows from Lemma~\ref{loose1} that

\begin{sublemma}
\label{comparable-4}
$(P_1,P_2,\ldots, P_k)$ is a swirl-like 
flower of order $k$ in $M\ba x/y$ so that
$M\ba x/y$ is $k$-fractured.
\end{sublemma}

We now consider the reverse and show that all $k$-fractures of
$M\ba x/y$ have the same form. Assume that $\QQ$ is a $k$-fracture
of $M\ba x/y$. 

\begin{sublemma}
\label{comparable-5}
$\QQ$ has order $k$, and there is a swirl-like flower
$(Q_1,Q_2,\ldots,Q_k)\cong \QQ$
such that $(Q_1\cup\{y\},Q_2,\ldots,Q_k)$ is a $k$-fracture of
$M\ba x$ and $x$ is in the coguts of $(Q_1,Q_2)$.
\end{sublemma}

\subproof
By Lemma~\ref{in-petal}, there is a petal $Q$ of $\QQ$ such that either
$P_1\subseteq \hQ$ or $P_2\subseteq \hQ$. 
But $y\in\cl^*_{M\ba x}(P_1)$ and
$y\in\cl^*_{M\ba x}(P_2)$, so $y\in\cl^*_{M\ba x}(\hQ)$.
Now choose a 
flower $(Q_1,\ldots, Q_m)$ equivalent to 
$\QQ$ such that $Q_1=\hQ$. Then
$(Q_1\cup\{y\},Q_2,\ldots,Q_m)$ is a 
maximal swirl-like flower 
of $M\ba x$. By 
\ref{comparable-1}, $m=k$ and 
we may assume that $y$ is in the coguts of $(Q_1,Q_2)$.
\end{proof}  

Assume that part (i) of the lemma fails. 
Then, as $M/y$ is $3$-connected,
it must be that $M/y$ is $k$-fractured. As $M\ba x/y$ is
$k$-fractured, it follows from Lemma~\ref{keep-fracture} 
that there is a tight swirl-like flower $(R_1,R_2,\ldots,R_k)$ 
of $M\ba x/y$ such that, for some
$i\in\{1,2,\ldots, k\}$, we have $x\in\cl_{M/y}(R_i)$.
By \ref{comparable-5}, we may assume that $y\in\cl_M(R_1\cup R_2)$,
and, as $y\more x$, we have $x\in\cl_M(R_1\cup R_2)$, so that
$x\in\cl_{M/y}(R_1\cup R_2)$. Thus we may assume that $i\in\{1,2\}$
and, up to labels, we have $x\in\cl_{M/y}(R_1)$. But then,
$x\in\cl_M(R_1\cup\{y\})$ so that $y\in\cl_M(R_1\cup\{x\})$,
and $(R\cup\{x,y\},R_2,\ldots,R_k)$ 
is a $k$-fracture of $M$,  contradicting the
fact that $M$ is $k$-coherent. Thus (i) holds.

Consider part (ii).
As $y$ is in the coguts of $P_1$ and $P_2$ in $M\ba x$, 
by Lemma~\ref{delete-cofix}, $y$ is cofixed
in $M\ba x$. As $x\less y$, 
it follows from Lemma~\ref{born-free}(ii),
that either $x$ and $y$ are clones or $y$ is cofixed in $M$.
\end{proof}

Lemma~\ref{comparable} has some useful consequences.
The condition that $\{a,b,c\}$ is not in a 4-element fan
required by the next lemma is necessary.
Let $(p,q,r,s)$ be a 4-element fan
of a $3$-connected matroid $M$. Then it is quite possible
for $q$ and $r$ to be comparable---indeed they can even
be clones. Yet neither $M\ba q$ nor $M/q$ is $3$-connected.

\begin{lemma}
\label{sleazy-triangle}
Let $\{a,b,c\}$ be a triangle of the $3$-connected matroid
$M$ that is not contained in a $4$-element fan.
If $a$ is not fixed, then $M\ba b$ is $3$-connected.
\end{lemma}

\begin{proof}
Assume that $M\ba b$ is not 3-connected. Let
$(A,C)$ be a $2$-separation of $M\ba b$.
As $M$ is $3$-connected we cannot have $\{a,c\}\subseteq A$,
so we may assume that $a\in A$ and $c\in C$. As $\{a,b,c\}$
is not contained in a 4-element fan, $a$ does not 
belong to a triad so that $|A|\geq 3$. If $a\notin \cl(A-\{a\})$,
then $a\in\cl^*_{M\ba b}(C)$, so that  
$(A-\{a\},C\cup\{a\})$
is also a 2-separation of $M\ba b$. But this implies that
$(A-\{a\},C\cup \{a,b\})$ is a 2-separation of $M$.
Thus $a\in\cl(A-\{a\})$.
As $M$ is 3-connected, $b\notin\cl(A-\{a\})$.
Hence, by Lemma~\ref{freedom3}, $a$ is fixed in $M$.
\end{proof}

Lemma~\ref{sleazy-triangle} extends to $k$-coherence.

\begin{corollary}
\label{free-triangle}
Let $\{a,b,c\}$ be a triangle of the $k$-coherent matroid $M$
that is not contained in a $4$-element fan.
If $a$ is not fixed in $M$, then $M\ba b$ is $k$-coherent.
\end{corollary}

\begin{proof}
By Lemma~\ref{sleazy-triangle} $M\ba b$ is 
$3$-connected.
By Lemma~\ref{freedom3}, $b\less a$ and, as $b$ is in a triangle,
$M/b$ is not 3-connected and hence not $k$-coherent. Thus,
by Lemma~\ref{comparable}, $M\ba b$ is $k$-coherent.
\end{proof}

\begin{corollary}
\label{clone-coherent}
Let $x$ and $y$ be clones of the $k$-coherent matroid $M$.
If $\{x,y\}$ is not contained in a $4$-element fan,
then either $M\ba x$ or $M/x$ is $k$-coherent.
\end{corollary}

\begin{proof}
Say that $x$ is in a triangle. Then there is a triangle containing
both $x$ and $y$. By Corollary~\ref{free-triangle}, 
$M\ba x$ is $k$-coherent.
Dually, if $x$ is in a triad, then $M/x$ is $k$-coherent.

Assume that $x$ is not in a triangle or a triad. 
Then, by Bixby's Lemma,
either $M\ba x$ or $M/x$ is 3-connected. Assume the former.
If $M\ba x$ is not $k$-coherent, then, by Lemma~\ref{comparable},
$M/y$ is $k$-coherent. 
But $M/y$ is isomorphic to $M/x$ and the corollary follows.
\end{proof}

\section{Elementary Properties of $k$-skeletons}
\label{skeleton1}

In this section 
we give
some basic properties of $k$-skeletons. 
The first is a more-or-less
immediate consequence of the definition.

\begin{lemma}
\label{skeleton-dual}
The matroid $M$ is a $k$-skeleton if and only if $M^*$ is.
\end{lemma}

A pleasant property of $k$-skeletons is that they do not have
large fans.

\begin{lemma}
\label{no-fan}
If $M$ is a $k$-skeleton, then $M$ has no $4$-element fans.
\end{lemma}

\begin{proof}
Say that $M$ is a skeleton and that $M$ has a 
4-element fan. Then there is a fan $F=(f_1,f_2,f_3,f_4)$
such that $f_1$ is an initial element of a 
maximal fan of $M$. 
As $M$ is not a wheel or a whirl, $F$ is exactly 
$3$-separating in $M$. Up to duality
we may assume that $M\ba f_1$ is $3$-connected, so that,
by Corollary~\ref{remove-sequential}, $M\ba f_1$ is $k$-coherent.
Now $\{f_1,f_2,f_3\}$ is a triangle of $M$ so that
$f_1\in\cl(\{f_2,f_3\})$. Also $f_1\in\cl(E(M)-F)$, and
no other element of $\{f_1,f_2,f_3\}$ is in the closure of this
set, so, by Lemma~\ref{freedom3},
$f_1$ is fixed in $M$, contradicting the definition of a 
$k$-skeleton.
\end{proof}

We can also be quite specific about elements in triangles or triads
in $k$-skeletons.

\begin{lemma}
\label{skeleton-triangle}
Let $T$ be a triangle of the $k$-skeleton $M$. Then the
following hold.
\begin{itemize}
\item[(i)] Either $T$ is a clonal triple, 
or $T$ is a standard or costandard
$k$-wild triangle. 
\item[(ii)] $M\ba t$ is $3$-connected for all $t\in T$.
\end{itemize}
\end{lemma}

\begin{proof}
Say $T=\{t_1,t_2,t_3\}$.
Assume that $T$ is not a standard or costandard $k$-wild triangle.
By Lemma~\ref{no-fan}, $T$ is not contained in a $4$-element fan,
so by Theorem~\ref{wild1}, there is an element $t\in T$ such that
$M\ba t$ is $k$-coherent. Assume that $t=t_1$.
By the definition
of a $k$-skeleton, $t_1$ is not fixed, 
so, by Corollary~\ref{free-triangle},
both $M\ba t_2$ and $M\ba t_3$ are $k$-coherent.
Thus neither $t_1$, $t_2$ nor $t_3$ is fixed. It now follows
easily that $\{t_1,t_2,t_3\}$ is a clonal triple. Thus (i) holds.
Part (ii) follows immediately.
\end{proof}

We had to make a special case of wheels and whirls and
exclude them from the definition of $k$-skeleton. 
The reason for this is that wheels and whirls trivially
satisfy properties (ii) and (iii) of a $k$-skeleton as they do not
have an element $x$ such that either $M\ba x$ or $M/x$
is $k$-coherent. There exists a danger that 
we may fall into this case in
minors hence complicating analyses in proofs. 
The next lemma shows that this will not be a problem for us.

\begin{lemma}
\label{no-wheel}
Let $e$ be an element of the $k$-skeleton $M$ of rank at least $3$.
Then $M\ba e$ is not a wheel or a whirl.
\end{lemma}

\begin{proof}
Assume that $M\ba e$ is a wheel or a whirl. If $T$ is a triangle
of $M\ba e$, then at least two elements of $T$ are fixed in $M\ba e$
and hence fixed in $M$. By Lemma~\ref{no-fan}, $M$ has no
$4$-element fans. This means that $e$ must block every triad of
$M\ba e$. It follows easily from this that 
if $(A,B)$ is a 3-separation of $M$, then either $A$ or $B$ is a
triangle.
Thus $M$ has no non-sequential $3$-separations so that
a $3$-tree for $M$ consists of a single vertex. This means that
all elements of $M$ are peripheral. By Lemma~\ref{wheel-wild}
$T$ has no $k$-wild triangles. By Lemma~\ref{skeleton-triangle},
$T$ is a clonal triple contradicting the fact that $T$
contains fixed elements.
\end{proof}

Quads in $k$-skeletons are also well behaved.

\begin{lemma}
\label{dag-skeleton}
Let $M$ be a $k$-skeleton and let 
$D$  be a quad of $M$.
If $d\in D$, then $d$ is in neither a triangle nor a triad and
both $M\ba d$ and $M/d$ are $k$-coherent.
\end{lemma}

\begin{proof}
Consider a 3-tree $T$ for $M$. This displays
a 3-separation equivalent to $D$, and it follows easily that
the members of $D$ are peripheral elements of $T$. 
By Lemma~\ref{wheel-wild},
no element of $D$ is in a $k$-wild triangle or triad. 
Consider $d\in D$. Assume that $d$ is in a triangle $T$.
Then, as $D$ is both a circuit and a cocircuit, $T$ must have
exactly one element $t$ that is not in $D$. It is easily
seem that Lemma~\ref{freedom3} applies
so that $t$ is fixed in $M$. However, by
Lemma~\ref{skeleton-triangle}, $T$ is a clonal 
triple so that no element of $T$ is fixed. This contradiction shows
that $d$ is not in a triangle, and dually, $d$ is not in a triad.
By Lemma~\ref{dag-remove} both $M\ba d$ and $M/d$ are $k$-coherent.
\end{proof}

A flower is {\em canonical\ }
\index{canonical}
 if it has no loose elements.

\begin{lemma}
\label{clean}
If $\FF$ is a tight swirl-like flower of a 
$k$-skeleton with at least three petals
then $\FF$ is canonical.
\end{lemma}

\begin{proof}
Say $\FF$ is not canonical. Then, by taking an appropriate
concatenation, we see that $M$ has a tight 
swirl-like flower
$(P_1,P_2,P_3)$ with loose elements between $P_1$ and $P_2$.
Let $l$ be the initial element of the fan of loose elements
between $P_1$ and $P_2$. Up to equivalence of flowers and 
duality, we may assume that $l\in P_1$ and $l\in\cl(P_2)$.
By Lemma~\ref{red-hill}, $M\ba l$ is $3$-connected.
As $(P_1,P_2,P_3)$ is tight, $|P_1|>2$. If $|P_2|=|P_3|=2$,
then $P_2\cup P_3$ is a quad and $l$ is in a triangle
with $P_2$, contradicting Lemma~\ref{dag-skeleton}.
Thus the hypotheses of Lemma~\ref{loose-removable}(iii) hold,
and, by that lemma, we deduce that $l$ does not expose any
$3$-separations in $M\ba l$. Hence $M\ba l$ is $k$-coherent.
But $l\in\cl(P_2)$ and $l\in\cl(P_1-\{l\})$. Moreover,
$\sqcap(P_1-\{l\},P_2)=1$, so, by Lemma~\ref{freedom1},
$x$ is fixed in $M$ and we have contradicted the definition of
$k$-skeleton.
\end{proof}

Ideally we would like to find an element $e$ in a $k$-skeleton 
such that
either $M\ba e$ or $M/e$ is a $k$-skeleton. 
This would give us a neat  chain
theorem for skeletons. Unfortunately, life is not that
simple in the world of $k$-skeletons.
In the remaining sections of this chapter
we prove theorems that identify
situations when it is possible to remove elements to keep
a $k$-skeleton and identify
certain annoying structures that make life difficult for us.

\section{Comparable Elements in $k$-skeletons}

Knowing that an element is comparable with 
another is certainly helpful.
This section is devoted to proving the following theorem.

\begin{theorem}
\label{reduce-comparable}
Let $a$ and $b$ be elements of the $k$-skeleton $M$, 
where $a\less b$.
If either 
\begin{itemize}
\item[(i)] $a\prec b$, or
\item[(ii)] $a\cong b$, and $\fr^*(a)\geq 3$,
\end{itemize}
then $M\ba a$ is a $k$-skeleton.
\end{theorem}

\begin{proof}
Assume that either (i) or (ii) holds.
We first show that

\begin{sublemma}
\label{reduce-comparable1}
$M\ba a$ is $k$-coherent, and, if $a\prec b$, then 
both $M\ba a$ and $M/b$ are $k$-coherent.
\end{sublemma}

\begin{proof}
If condition (ii) holds this 
follows from Corollary~\ref{3-comparable}.
Assume that $a\prec b$. 
Assume that neither $M\ba a$ nor $M/b$ is $k$-coherent.
If $M\ba a$ is $3$-connected and $k$-fractured, then by 
Lemma~\ref{comparable},
$M/b$ is $k$-coherent and $b$ is 
cofixed in $M$, contradicting the fact
that $M$ is a $k$-skeleton. Thus 
$M\ba a$ is not $3$-connected, and  
similarly $M/b$ is not $3$-connected.

Assume that $b$ is in a triangle. 
Then, as $b\succ a$, there is a triangle
$T$ containing both $b$ and $a$. By Lemma~\ref{skeleton-triangle},
$M\ba a$ is $3$-connected, contradicting the fact 
that $M\ba a$ is not 3-connected.
Hence $b$ is not in a triangle. Thus $\si(M/b)$ is not $3$-connected.
As $a\prec b$ it follows from Corollary~\ref{freedom2}, that
$\si(M/a)$ is not $3$-connected.
But then $\co(M\ba a)$ is $3$-connected, and, as $a$
is not in a triad, we have contradicted the assumption that 
$M\ba a$ is not $3$-connected.

Thus either $M\ba a$ or $M/b$ is $k$-coherent.
Assume that $M/b$ is $k$-coherent. 
If $M\ba a$ is not $3$-connected, then,
by the dual of Corollary~\ref{freedom2},
$b$ is cofixed in $M$
contradicting the fact that $M$ is a skeleton.
Thus $M\ba a$ is $3$-connected. If $M\ba a$ is not $k$-coherent,
then by Lemma~\ref{comparable}, we again obtain the 
contradiction that $b$ is cofixed 
in $M$.
Thus $M\ba a$ is indeed $k$-coherent. 
A dual argument shows that, if $M\ba a$ is $k$-coherent,
then $M/b$ is $k$-coherent. Thus, both $M\ba a$ and $M/b$ are 
$k$-coherent.
\end{proof}

\begin{sublemma}
\label{b-not-cofixed}
$b$ is not cofixed in $M\ba a$.
\end{sublemma}

\begin{proof}
As $M/b$ is $k$-coherent, and $M$ is a
$k$-skeleton, $b$ is not cofixed in $M$.
It now follows from Lemma~\ref{cofree-at-last} that
$b$ is not cofixed in $M\ba a$.
\end{proof}

We now work towards showing that $M\ba a$ is a $k$--skeleton.

\begin{sublemma}
\label{reduce-comparable3}
If $z$ is fixed in $M\ba a$, then $M\ba a,z$ is not $k$-coherent.
\end{sublemma}

\begin{proof}
Assume that $z$ is fixed in $M\ba a$ and that $M\ba a,z$ is 
$k$-coherent.
Then $z$ is fixed in $M$ and, as $M$ is a $k$-skeleton, $M\ba z$
is not $k$-coherent. As $M\ba a,z$ is $k$-coherent, $M\ba z$ is 
$3$-connected. Hence $M\ba z$ is $k$-fractured.
Altogether we have $M$, $M\ba a$
and $M\ba a,z$ are $k$-coherent and $M\ba z$ is $k$-fractured.
As $z$ is fixed in $M$, $z\neq b$.
Let $(P_1,P_2,\ldots,P_k)$ be a $k$--fracture of $M\ba z$. By 
Lemma~\ref{comparable-in-petal}, we may assume that both
$a$ and $b$ are in $P_1$.

Consider $M\ba z,a$. By Lemma~\ref{gain-coherence}(i), 
$P_1-\{a\}$ is a set
of loose elements of the swirl-like quasi-flower 
$(P_1-\{a\},P_2,P_3,\ldots,P_k)$ of $M\ba z,a$.
By Lemma~\ref{loose2}, elements of $P_1-\{a\}$ are either fixed or 
cofixed in $M\ba z,a$. But $b$ is not fixed in $M$ and hence
not fixed in $M\ba z,a$. Therefore $b$ is cofixed in $M\ba z,a$.
Now $z$ is fixed in $M\ba a$, so it is not the case
that $z\more a$ in $M$. 
Therefore, by Corollary~\ref{born-free}(ii), $b$ is cofixed in 
$M\ba a$, contradicting \ref{b-not-cofixed}.
\end{proof}

\begin{sublemma}
\label{reduce-comparable4}
If $z$ is cofixed in $M\ba a$, then $M\ba a/z$ is not
$k$-coherent.
\end{sublemma}

\begin{proof}
Assume that $z$ is cofixed in $M\ba a$ and that $M\ba a/z$ is $k$-coherent.
By \ref{b-not-cofixed}, $z\neq b$.
Assume that  $z$ is not cofixed in $M$. Then, by 
Corollary~\ref{born-free}(ii), 
$z\less a$ in $M$, and hence
$z\less b$ in $M$ so that $z\less b$ in $M\ba a$. But then, 
$b$ is cofixed in $M\ba a$ contradicting \ref{b-not-cofixed}.
Hence $z$ is cofixed in $M$.

We now show that $M/z$ is $3$-connected. Assume not. As $M/z\ba a$
is $k$-coherent and hence $3$-connected, there is an element $p$ of $M$
such that $\{a,z,p\}$ is a triangle. As $z$ is cofixed in $M$,
$z\npreceq b$. Hence there is a cyclic flat of $M$
containing $b$ but not $z$.
As $a\less b$, this cyclic flat contains $a$. 
But now, by Lemma~\ref{freedom3}, $a$ is fixed in $M$, and
$M\ba a$ is $k$-coherent, contradicting the fact that 
$M$ is a skeleton. Thus $M/z$ is $3$-connected.

Since $z$ is cofixed in $M$ and $M/z$ is $3$-connected, it must
be the case that $M/z$ is $k$-fractured. 
By Lemma~\ref{comparable-in-petal},
there is a $k$--flower $(P_1,P_2,\ldots,P_k)$ of $M/z$ such that
$\{a,b\}\subseteq P_1$. But $M/z\ba a$ is $k$-coherent, so 
by Lemma~\ref{gain-coherence}(i), the elements of $P_1-\{a\}$ are all
fixed or cofixed in $M/z\ba a$, in particular, 
$b$ is either fixed or cofixed in
$M/z\ba a$. As $a\less b$ in $M/z$, and $M/z$ is 
$3$-connected, $b$ is not fixed in $M/z\ba a$. Thus 
$b$ is cofixed in $M/z\ba a$. But now, $b$ is cofixed in 
$M\ba a$, and we have contradicted \ref{b-not-cofixed}.
\end{proof}

By \ref{reduce-comparable1},\ref{reduce-comparable3} and
\ref{reduce-comparable4}, $M\ba a$ is a $k$-skeleton. 
In the case that $a\prec b$, it follows by duality that
$M/b$ is also a $k$-skeleton.
\end{proof}

As an easy consequence of Theorem~\ref{reduce-comparable} we 
have

\begin{corollary}
\label{stuck}
Let $a$ and $b$ be comparable elements in the $k$-skeleton
$M$. If neither $M\ba a$ nor $M/a$ is a $k$-skeleton, then
$\{a,b\}$ is a clonal class, and $\fr(a)=\fr^*(a)=2$.
\end{corollary}

\section{Bogan Couples}
\label{bogan-couples}

In this section we examine one of the problematic structures that
turns out to be of importance to us. 
Recall that a flower is {\em canonical\ } if it is the only member of
its equivalence class. If $f$ is an element of 
the $k$-skeleton $M$
such that $M\ba f$ is 3-connected but not $k$-coherent, then 
$M\ba f$ is certainly not a $k$-skeleton. But one might
hope that the fact that $M$ is a $k$-skeleton would impose
structure on a $k$-fracture in $M\ba f$. In particular it would
be useful to know that such a $k$-fracture is canonical.
In this section we show that this is almost always the case unless
$f$ is a very particular type of feral element.

Let $a,b$ be elements of the $k$-coherent matroid $M$. Then
a {\em bogan display\ }
\index{bogan display} 
for $\{a,b\}$ is a partition
$(R,S,T,\{a,b\})$ of $E(M)$ together
with  partitions $(R_1,R_2,\ldots,R_{k-2})$,
$(S_1,S_2,\ldots,S_r)$ and $(T_1,T_2,\ldots,T_{k-2})$ 
of $R$, $S$ and $T$
respectively such that the following hold.
\begin{itemize}
\item[(i)] $(R_1\cup\{b\},R_2,\ldots,R_{k-2},S_1,S_2,\ldots,S_r,T)$
is a maximal $k$-fracture of $M\ba a$, and 
$b$ is in the coguts of $R_1\cup\{b\}$ and $T$.
\item[(ii)] $(R,S\cup\{b\},T_1,T_2,\ldots,T_{k-2})$ is a maximal
$k$-fracture of $M/a$.
\item[(iii)] $(R,S_1,S_2,\ldots,S_r,T_1,T_2,
\ldots,T_{k-3},T_{k-2}\cup\{a\})$
is a maximal $k$-fracture of $M\ba b$, and $a$ is in the
coguts of $T_{k-2}\cup\{a\}$ and $R$.
\item[(iv)] $(R_1,R_2,\ldots,R_{k-2},S\cup\{a\},T)$ is a maximal
$k$-fracture of $M/b$.
\item[(v)] $(R_1,R_2,\ldots,R_{k-2},S\cup T\cup\{a,b\})$,
$(S_1,S_2,\ldots,S_r,R\cup T\cup\{a,b\})$ and 
$(T_1,T_2,\ldots,T_{k-2},R\cup S\cup\{a,b\})$ are 
maximal swirl-like flowers in $M$.
\item[(vi)] $a$ and $b$ are cofixed in $M$.
\end{itemize}

The pair $\{a,b\}$ is a {\em bogan couple\ }
\index{bogan couple} 
if it has a 
bogan display.
The property of being a bogan couple is not self dual. If $\{a,b\}$
is a bogan couple in $M^*$, then $\{a,b\}$ is a 
{\em cobogan couple}.
\index{cobogan couple}

\begin{figure}
\begin{tikzpicture}[thick,line join=round]
	\coordinate (o) at (0,0);
	\coordinate (p) at (0,2);
	\coordinate (q) at (4,2);
	\coordinate (r) at (4,0);
	\coordinate (s) at (6,0);
	\coordinate (t) at (6,2);
	\coordinate (n) at (0,1);
	\coordinate (a) at ($(p)+(45:1.414)$);
	\coordinate (b) at ($(q)+(135:1.414)$);
	\coordinate (c) at (-45:2.121);
	\coordinate (d) at ($(s)+(-135:2.121)$);
	\node[pattern color=lines,draw,circle through=(b),pattern=north east lines] at ($(b)!0.5!(q)$) {};
	\node[pattern color=lines,draw,circle through=(a),pattern=north east lines] at ($(a)!0.5!(b)$) {};
	\node[pattern color=lines,draw,circle through=(p),pattern=north east lines] at ($(p)!0.5!(a)$) {};
	\node[pattern color=olines,draw,circle through=(o),pattern=north east lines] at ($(o)!0.5!(p)$) {};
	\node[pattern color=rlines,draw,circle through=(c),pattern=north east lines] at ($(c)!0.5!(o)$) {};
	\node[pattern color=rlines,draw,circle through=(d),pattern=north east lines] at ($(d)!0.5!(c)$) {};
	\node[pattern color=rlines,draw,circle through=(s),pattern=north east lines] at ($(s)!0.5!(d)$) {};
	\filldraw[fill=white] (o) -- (p) -- (a) -- (b) -- (q) -- (t) -- (s) -- (d) -- (c) -- cycle;
	\coordinate[label=-135:$b$] (y) at (4,1);
	\coordinate[label=right:$a$] (z) at (6,1);
	\draw (p) -- (q) -- (r);
	\draw (n) -- (z);
	\draw (o) -- (s);
	\node at ($(a)!0.5!(b) + (0,0.5)$) [rectangle,fill=white,draw=white,inner sep=1pt] {\textcolor{glabels}{$R_2$}};
	\node at ($(a)!0.5!(p) + (-0.25,0.25)$) [rectangle,fill=white,draw=white,inner sep=0.5pt] {\textcolor{glabels}{$R_3$}};
	\node at ($(s)!0.5!(d) + (0.35,-0.35)$) [rectangle,fill=white,draw=white,inner sep=1pt] {\textcolor{rlabels}{$T_1$}};
	\node at (3.725,2.775) [rectangle,fill=white,draw=white,inner sep=0.5pt] {\textcolor{glabels}{$R_1$}};
	\node at (-0.5,1) [rectangle,fill=white,draw=white,inner sep=1pt] {\textcolor{olabels}{$S$}};
	\node at (0.3,-1.1) [rectangle,fill=white,draw=white,inner sep=1pt] {\textcolor{rlabels}{$T_3$}};
	\node at (3,-2.3) [rectangle,fill=white,draw=white,inner sep=1pt] {\textcolor{rlabels}{$T_2$}};
	\foreach \pt in {y,z} \fill[black] (\pt) circle (3pt);
\end{tikzpicture}
\caption{A Bogan Couple}\label{bogancouple}
\end{figure}
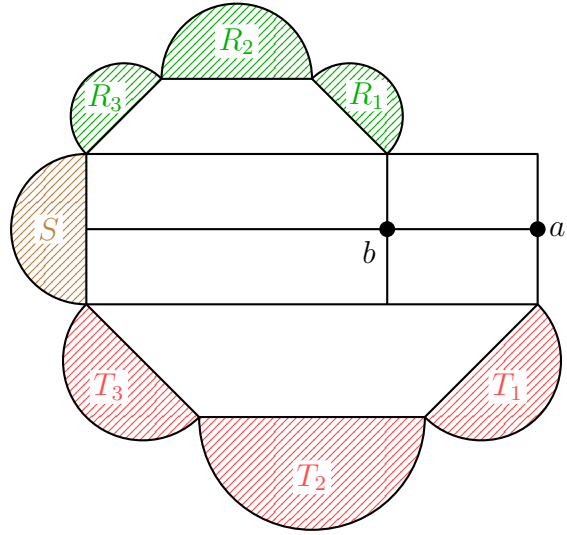

\begin{theorem}
\label{canon1}
Let $a$ be an element of the $k$-skeleton $M$ such that $M\ba a$ is 
$k$-fractured. If $b$ is a loose element of a $k$-fracture $\PP$
of $M\ba a$, then $b$ is a loose coguts element of $\PP$ and 
either 
\begin{itemize}
\item[(i)] $a\more b$, or
\item[(ii)] $\{a,b\}$ is a bogan couple of $M$.
\end{itemize}
\end{theorem}

\begin{proof}
Assume that ${\bf P}=(P_1\cup\{b\},P_2,\ldots,P_l)$ is a 
maximal $k$-fracture of $M\ba a$ and that $b$ is in the guts or the
coguts of $(P_1\cup\{b\},P_2)$.

\begin{sublemma}
\label{canon1.1}
$b$ is in the coguts of $(P_1\cup\{b\},P_2)$.
\end{sublemma}

\subproof
Assume that $b$ is in the guts of $(P_1\cup\{b\},P_2)$. Then, by
Lemma~\ref{loose2}, $b$ is fixed in $M\ba y$, 
so that $b$ is fixed in $M$.
By Lemma~\ref{red-hill}, $M\ba a,b$ is 
3-connected up to series pairs, so that
either $M\ba b$ is 3-connected, or $M\ba b$ has a series pair. 

We now eliminate the latter case.
Assume that $M\ba b$ has a series pair.
Then $b$ is in a triad $\{p,b,q\}$
of $M$ and $\co(M\ba b)$ is 3-connected. 
(Otherwise $M$ has a 4-element fan, 
contradicting Lemma~\ref{no-fan}.) 
As $b$ is fixed in $M$, this triad is not in a clonal triple,
so, by Lemma~\ref{skeleton-triangle}, it is $k$-wild.
But $\co(M\ba b)$ is 3-connected, so by 
Lemma~\ref{k-wild-different}(iii), it is not standard.
We need to show that it is not costandard. By 
Lemma~\ref{loosey-loose} $(p,b,q)$ is a consecutive set 
of loose elements in the flower
$(P_1\cup\{b\},P_2,\ldots,P_l)$ of $M\ba a$. But then,
$\si(M\ba a/p,q)$ is $3$-connected, as $p$ and $q$ are
loose coguts elements of this flower.
Thus $\si(M/p,q)$ is $3$-connected, contradicting 
Lemma~\ref{cowild-win}.

Therefore $M\ba b$ is 3-connected. As $b$ is fixed in $M$,
the matroid $M\ba b$ is $k$-fractured. 
By considering a $k$-fracture of $M\ba b$,
we see that there is a 3-separation $(R\cup\{a\},G)$ of 
$M\ba b$ such that $a\in\cl(R)$, and $(R\cup \{a\},G)$ is 
well blocked by $b$. Now $(R\cup\{b\},G)$ and
$(R,G\cup\{b\})$ are unsplit 4-separating sets of $M\ba a$.
If $(R,G)$ does not cross $P_1$, then $b$ does not block
$(R,G)$, as $b\in\cl(P_1)$. Thus $(R,G)$ crosses $P_1$,
and, similarly, $(R,G)$ crosses $P_2$.

Now Lemma~\ref{4-cross0} applies and we may assume that, for some
$T\in\{R,G\}$, we have $T\cup\{b\}\subseteq P_1\cup P_2\cup\{b\}$.
By considering the flower $(P_1\cup\{b\},P_2,\ldots,P_l)$,
Lemma~\ref{2-crossing} implies that 
$\lambda_{M\ba a}((P_1\cup\{b\})\cap T)=2$ 
and by considering the flower
$(P_1,P_2\cup\{b\},\ldots,P_l)$, Lemma~\ref{2-crossing} implies
that $\lambda_{M\ba a}(P_1\cap T)=2$. Hence
$b\in\cl^{(*)}(P_1\cap T)$. 
But $b\in\cl(P_2)$, so $b\notin\cl^*(P_1\cap T)$.
Therefore $b\in\cl(P_1\cap T)$, contradicting the assumption that 
$b$ blocks $(R,G)$.
\end{proof}

Assume that that $a\not\more b$. We need to prove that
$\{a,b\}$ is a bogan couple. We proceed by accumulating
properties of the pair $\{a,b\}$.

\begin{sublemma}
\label{canon1.15}
The element $b$ is cofixed in $M$ and $M/b$ is $3$-connected
and $k$-fractured.
\end{sublemma}

\subproof
As $b$ is 
a coguts element of $\PP$, the element
$b$ is cofixed in $M\ba a$. As
$a\not\more b$, it follows from Corollary~\ref{im-free} 
that $b$ is cofixed in
$M$. It follows from
\ref{canon1.1} that $b$ is the unique loose element between
$P_1$ and $P_2$. So, by Lemma~\ref{red-hill}, 
$M\ba a/b$ is $3$-connected.  Thus $M/b$ is $3$-connected
unless there is a triangle of $M$ containing $\{a,b\}$.

Assume that $T$ is such a triangle. As $b$ is cofixed in $M$,
the triangle $T$ is not a clonal triple and hence is 
$k$-wild. As $b$ is in the coguts of $(P_1\cup\{b\},P_2)$,
the matroid $\co(M\ba a,b)$ is not $3$-connected.
It follows from the definition of $k$-wild triangle that,
if $T$ is standard, then $\co(M\ba a,b)$ is 
3-connected. Moreover, by basic properties of the 
$\Delta-Y$ operation, the same conclusion holds if 
$T$ is costandard. 

It follows from the above contradiction that $M/b$ is $3$-connected.
As $b$ is cofixed in
$M$, the matroid $M/b$ is $k$-fractured.
\end{proof}

\begin{sublemma}
\label{canon1.2} We may assume that, 
for some $i_1,i_2\in\{1,2,\ldots,l\}$, with
$2 \leq i_1\leq i_2 \leq l$  the partition 
$(P_1,\ldots,P_{i_1-1},P_{i_1}\udots P_{i_2}\cup \{a\},
P_{i_2+1},\ldots,P_l)$
is a maximal $k$-fracture of $M/b$.
\end{sublemma}

\subproof
Let $(Q_1\cup\{a\},Q_2,\ldots,Q_m)$ be a maximal $k$-fracture of
$M/b$. By Lemma~\ref{loose-removable}, 
$(P_1,P_2,\ldots,P_l)$ is a maximal flower in
$M\ba a/b$.

Consider the partition $(Q_1,Q_2,\ldots,Q_m)$ of $E(M\ba a/b)$. This
will be a flower $\QQ$ in $M\ba a/b$ unless $|Q_1|=1$, in which 
case it is a quasi-flower. If 
$\QQ\less(P_1,P_2,\ldots,P_l)$, then the sublemma clearly holds. 
Assume otherwise. Then, by
Lemma~\ref{in-petal}, 
we may assume that there is a pair $i,j$, such
that $E(M\ba a/b)=P_i\cup Q_j$.

Assume that $j=1$. Then either $P_1\subseteq Q_1$ or 
$P_2\subseteq Q_1$. We lose no generality in assuming that 
$P_1\subseteq Q_1$, so that $P_1\cup\{a\}\subseteq Q_1\cup\{a\}$.
But $b\in\cl^*_{M\ba a}(P_1)$.
Hence $b\in\cl^*_M(Q_1\cup\{a\})$ and it follows that
$(Q_1\cup\{a,b\},Q_2,\ldots,Q_l)$ is a $k$-fracture of $M$,
contradicting the
fact that $M$ is $k$-coherent.

Hence we may assume that $j\neq 1$. 
Then either $Q_1\cup Q_2\subseteq P_i$ or
$Q_1\cup Q_m\subseteq P_i$ and we may assume that the former case holds.
As $M/b\ba a$ is 3-connected, $a$ is not a coguts element of
$Q_1$, so that $a\in\cl_{M/b}(Q_1\cup Q_2)$, that is,
$a\in\cl_{M/b}(P_i)$. Hence 
$(P_1,P_2,\ldots,P_i\cup\{a\},\ldots,P_l)$
fractures $M/b$ and the sublemma holds with $i_1=i_2=i$.
\end{proof}

Assume that $i_1$ and $i_2$ are chosen so that 
\ref{canon1.2} is satisfied and $i_2-i_1$ 
is maximum.

\begin{sublemma}
\label{canon1.3}
\begin{itemize}
\item[(i)] $i_1\geq 3$, $i_2\leq l$.
\item[(ii)] Either $i_1>3$ or $i_2<l$.
\item[(iii)] $a$ blocks both
$P_1\cup\{b\}$ and $P_2\cup\{b\}$.                          
\item[(iv)] $a$ blocks either $P_3\cup P_4\udots P_{l-1}$ or
$P_4\cup P_5\udots P_l$.
\item[(v)] $a$ does not block $P_1$ or $P_2$.
\end{itemize}
\end{sublemma}

\subproof
Say $i_1=2$. Then $a\in\cl_M(P_2\cup P_3\udots P_{i_2}\cup\{b\})$,
so that $(P_1\cup P_2\udots P_{i_2}\cup\{b\},P_{i_2+1},\ldots,P_l)$
is a $k$-fracture of $M$. Hence $i_1\geq 3$ and similarly
$i_2\leq l$ so that (i) holds. Part (ii) follows from the fact that 
$(P_1,P_2,\ldots,P_{i_1-1},P_{i_1}\udots P_{i_2}\cup \{a\},
P_{i_2+1},\ldots,P_l)$
is a $k$-fracture of $M/b$.

Assume that $a$ does not block $P_1\cup\{b\}$.
Certainly $a\notin\cl(P_1\cup\{b\})$, so 
$a\in\cl(E(M)-(P_1\cup\{a,b\}))$.
Then $b\in\cl^*(P_1)$ and 
$(P_1\cup\{b\},P_2,\ldots,P_{i_1}\udots P_{i_2}\cup\{a\},\ldots,P_l)$
is a $k$-fracture of $M$. Thus $a$ blocks $P_1\cup\{b\}$ and
similarly
$a$ blocks $P_2\cup\{b\}$ and (iii) holds.

Consider (iv). Assume that $a$ blocks neither 
$P_3\cup P_4\udots P_{l-1}$
nor $P_4\cup P_5\udots P_l$. As $a$ does not block $P_3\udots P_{l-1}$
and $a$ does block $P_1\cup \{b\}$, 
we see that $a\notin\cl(P_3\udots P_{l-1})$.
Hence $a\in\cl(P_l\cup P_1\cup \{b\}\cup P_2)$. Similarly
$a\in\cl(P_1\cup\{b\}\cup P_2\cup P_3)$ and it follows from 
Lemma~\ref{modular}
that $a\in\cl(P_1\cup\{b\}\cup P_2)$. By (ii) and symmetry we
may asume that $i_2<l$. As 
$a\in\cl(P_{i_1}\cup P_{{i_1}+1}\udots P_{i_2}\cup\{b\}$,
we have $a\in\cl(\{b\}\cup P_2\cup P_3\udots P_{i_2})$ so that by 
Lemma~\ref{modular}
$a\in\cl(\{b\}\cup P_2)$, contradicting the fact that $a$ blocks
this set. Thus (iv) holds.

Consider (v). As $i_1\geq 3$ and $i_2\leq l$, the 
set $\cl(P_{i_1}\cup P_{{i_1}+1}\udots P_{i_2}\cup\{b\})$ 
avoids $P_1$ and $P_2$. As 
$a\in \cl(P_{i_1}\cup P_{{i_1}+1}\udots P_{i_2}\cup\{b\})$,
neither $P_1$ nor $P_2$ is blocked by $a$.
\end{proof}

\begin{sublemma}
\label{canon1.4}
$a$ is cofixed in $M$.
\end{sublemma}

\subproof
By \ref{canon1.3}(v), we may assume that $a$ blocks 
$P_3\cup P_4\udots P_{l-1}$. By \ref{canon1.3}(iii), 
$a$ blocks $P_1\cup\{b\}$.
By \ref{canon1.3}(iv), 
$(P_1\cup\{b\},P_2,P_3\cup P_4\udots P_{l-1},P_l)$
is a swirl-like flower of $M\ba a$. It now follows from
Lemma~\ref{swirl-cofix} that $a$ is cofixed in $M$.
\end{proof}

\begin{sublemma} 
\label{canon1.5}
$M/a$ is $3$-connected.
\end{sublemma}

\subproof
Assume that $a$ is in a triangle. As $M/b\ba a$ is
$3$-connected, this triangle must contain $b$, contradicting
the fact that $M/b$ is 3-connected. Thus $a$ is not in a 
triangle.
Assume that $M/a$ is not 3-connected. 
As $a$ is not in a triangle, $\si(M/a)$ 
is not 3-connected, so there is 
a 3-separation $(R\cup\{a\},B)$ of $M$ such that $(R,B)$ is a
$3$-separation of $M\ba a$ and $a\in\cl(R),\cl(B)$. By
Lemma~\ref{cross1}, $(R,B)$ conforms with the maximal flower
$\PP$ of $M\ba a$. If either 
$R$ or $B$ is contained in a petal of a flower equivalent to $\PP$,
then, as $a\in\cl(R),\cl(B)$ we obtain the contradiction that
$M$ is $k$ fractured. Thus $(R,B)$ is displayed
by $\PP$. But then, either $a$ does not block $P_1\cup\{b\}$
or $P_2\cup\{b\}$, contradicting \ref{canon1.3}(iii).
\end{proof}

As $M/a$ is 3-connected, and $a$ is cofixed in $M$ we see that
$M/a$ is $k$-fractured. Thus $a$ is a feral element of 
$M$. Let $\QQ=(Q_1,\ldots,Q_m)$ be a 
maximal $k$-fracture of $M/a$. We may assume
that there is a labelling of $\PP$ and $\QQ$ that forms a 
feral display in either $M$ or $M^*$. In the proof 
of the next claim we assume that
the reader is familiar with terminology from 
Section~\ref{feral-elements} 
associated with feral elements.

\begin{sublemma}
\label{canon1.6}
\begin{itemize}
\item[(i)] $m=k$, $i_1=k$, and $i_2=l$.
\item[(ii)] Up to labels $Q_2=P_2\cup P_3\udots P_{k-1}$, 
$Q_1=P_k\cup P_{k+1}\udots P_{l}\cup\{b\}$,
and $\{Q_3,Q_4,\ldots,Q_m\}$ is a partition of $P_1$.
\item[(iii)] $(P_2,\ldots,P_{k-1},E(M)-Q_2)$,
$(P_k,\ldots,P_l,E(M)-(Q_1-\{b\}))$, and
$(Q_3,\ldots,Q_k,E(M)-P_1)$ are displayed flowers in $M$.
\end{itemize}
\end{sublemma}

\subproof
Consider the
$k$-fracture $\PP$ of $M\ba a$. By Lemma~\ref{feral1}(iv), 
$a$ is either
$1$-blocking for $\PP$ or is $2$-spanned by $\PP$. 

Our first, somewhat painful, task is to 
show that $a$ is 1-blocking for $\PP$.
Assume otherwise. Then $a$ is 2-spanned
by $\PP$. As $a$ is 2-spanned but not 1-blocking, 
$a$ blocks exactly two adjacent petals of a flower 
equivalent to $\PP$. We also have $l=k$
so that $P_l=P_k$. As $a$ is not 1-blocking for
$\PP$ we are in the case where we have a feral display for
$a$ in $M^*$. We know that $a$ blocks $P_1\cup\{b\}$ and,
by \ref{canon1.3}(v), $a$ does not block $P_2$. Thus it must
be the case that $a$ blocks $P_k$.
Moreover, as $a$ is $2$-spanned, $a\in\cl(P_k\cup P_1\cup\{b\})$.
We have committed to a labelling of $\PP$, but we are still free
to label $\QQ$. Given this, we can say that $(Q_1,Q_2,\ldots,Q_m)$
and $(P_k,P_1\cup\{b\},P_2,\ldots,P_{k-1})$ form a feral 
display in $M^*$. In other words,
$(Q_1,Q_2,\ldots, Q_m)$ and 
$(P_k,P_1\cup\{b\},P_2,\ldots,P_{k-1})$
play the roles played by $(P_1,P_2,\ldots,P_m)$
and $(Q_1,Q_2,\ldots,Q_k)$ in the definition of feral display.
With this in mind, we see from the properties of a 
feral display that $P_2\cup P_3\udots P_{k-1}\subseteq Q_1$,
and, for some $j\in\{2,3,\ldots,m-1\}$, we have
$Q_{j+1}\cup Q_{j+2}\udots Q_m\subseteq P_k$
and $Q_2\cup Q_3\udots Q_j\subseteq P_1\cup\{b\}$.

We next show that $i_1=i_2=k$. As 
$P_{i_1}\cup P_{{i_1}+1}\udots P_{i_2}\cup \{a\}$ is $3$-separating
in $M/b$ but not in $M$, we have 
$b\in\cl_M(P_{i_1}\cup P_{{i_1}+1}\udots P_{i_2}\cup \{a\})$.
But $a\notin\cl_M(P_{i_1}\cup P_{{i_1}+1}\udots P_{i_2})$,
as otherwise $M$ is $k$-fractured, so we also have
$a\in\cl_M(P_{i_1}\cup P_{{i_1}+1}\udots P_{i_2}\cup\{b\})$.
It follows from this and the fact that $a$ blocks $P_k$, that
$i_2=k$. Now $a\in\cl_{M/b}(P_1\cup P_k)$
and $a\in\cl_{M/b}(P_{i_1}\cup P_{{i_1}+1}\udots P_{i_2}\cup P_k)$.
So, by Lemma~\ref{modular}, $a\in\cl_{M/b}(P_k)$.
It now follows from \ref{canon1.2} that $i_1=k$.
Thus we do have $i_1=i_2=k$ and we can also conclude that
$b\in\cl_M(P_k\cup\{a\})$.

Assume that $b\in Q_1$. By \ref{canon1.3}(iii),
$a$ blocks $P_2\cup\{b\}$, so
$a\in\cl^*_M(P_2\cup\{b\})$. 
But $Q_1\supseteq P_2\cup\{b\}$,
so $a\in\cl^*_M(Q_1)$. This implies the contradiction
that $(Q_1\cup\{a\},Q_2,\ldots ,Q_m)$
is a $k$-fracture of $M$. Thus $b\notin Q_1$, so that 
$b\in Q_2\cup Q_3\udots Q_j$. Assume that $|Q_2\udots Q_j|\geq 3$.
By property (vii) of a feral display, $Q_2\cup Q_3\udots Q_j$
is 3-separating in 
$M$ and hence $M\ba a$. It is now easily seen that $b$ is in the
coguts of $Q_2\cup Q_3\udots Q_j$. Thus 
$b\notin\cl_M(E(M)-(Q_2\cup Q_3\udots Q_j))$, contradicting the fact 
that $b\in\cl_M(P_k\cup\{a\})$. 

We are left with the annoying
case where $|Q_2\cup Q_3\udots Q_j|=2$, so that
$j=2$, $m=k$ and $|Q_2|=2$. 
Say $Q_2=\{b,z\}$. Recall that there is  a set $Z_2$ such that
$P_k=Q_3\cup Q_4\udots Q_k\cup Z_2$.
We know that $b\in\cl_M(P_k\cup\{a\})$.
Thus $\lambda_M(P_k\cup\{a,b\})\leq \lambda_M(P_k\cup\{a\})=3$.
We now show that $\lambda_M(P_k\cup\{a,b,z\})\geq 3$.
Assume otherwise; then $\lambda_{M\ba a}(P_k\cup\{b,z\})=2$.
Consider the flower 
$(P_1\cup\{b\},P_2,\ldots,P_k)$ of $M\ba a$.
The petal $P_1\cup\{b\}$ is tight and $b$ is loose between
$P_1$ and $P_2$, so that, by Lemma~\ref{fine1},
$z\notin \fcl_{M\ba a}(P_k)$. But now it is clear that
$P_k\cup\{b,z\}$ does not have the form of Lemma~\ref{3-cross}.
Hence $\lambda_M(P_k\cup\{a,b,z\})\geq 3$. 
As $a$ coblocks $Q_2\cup Q_3\udots Q_k$, we also deduce
that $\lambda_M(\{a,b,z\}\cup Q_3\cup Q_4\udots Q_k)=3$.
We have shown that $\lambda_M(P_k\cup\{a,b\})\leq 3$
and $\lambda_M(\{a,b,z\}\cup Q_3\cup Q_4\udots Q_k)=3$. 
The union of these
sets is $P_k\cup\{a,b,z\}$, and 
$\lambda_M(P_k\cup\{a,b,z\})\geq 3$. But the intersection
of these sets is $\{a,b\}\cup Q_3\cup Q_4\udots Q_k$.
Hence $\lambda_M(\{a,b\}\cup Q_3\cup Q_4\udots Q_k)\leq 3$.
Now $a\notin\clstar(Q_3\cup Q_4\udots Q_k)$, as otherwise
one of $M\ba a$ of $M/a$ is not 3-connected. Hence
$\lambda_M(\{a\}\cup Q_3\cup Q_4\udots Q_k)=3$.
Thus $b\in\clstar_M(Q_3\cup Q_4\udots Q_k\cup\{a\})$.
But $b\in\cl_M(P_1\cup P_2)$ as $b$ is a loose 
coguts element between $P_1$ and $P_2$ in $\PP$. 
Therefore $b\notin\cl^*_M(Q_3\cup Q_4\udots Q_k\cup\{a\})$.
Hence $b\in\cl_M(Q_3\cup Q_4\udots Q_k\cup\{a\})$
so that $b\in\cl_{M/a}(Q_3\cup Q_4\udots Q_k)$ and it follows that
$\{z,b\}$ is a loose petal of $\QQ$ 
contradicting the fact that  it is not.

The above contradiction shows at last that $a$
is indeed 1-blocking for $\PP$ and we consider the structure
that arises in this case now.
By the symmetry between $P_1$ and $P_2$, we may assume that
$a$ blocks $P_1\cup\{b\}$ and no other petal of $\PP$.
Indeed, we may assume that $m=k$ and that
$(P_1\cup\{b\},P_2,\udots, P_l)$ and $(Q_1,Q_2,\ldots,Q_k)$
form a feral display. By the properties of a feral
display $P_1\cup\{b\}$ properly contains $Q_3\cup Q_4\udots Q_k$, 
and for some $i\in\{2,3,\ldots,l-1\}$, we have 
$P_2\cup P_3\udots P_i\subseteq Q_2$,
and $P_{i+1}\cup P_{i+2}\udots P_k\subseteq Q_1$. 
By \ref{canon1.3}(v), $P_1$ is not blocked
by $a$. Hence $P_1$ is 3-separating in $M/a$ and it follows easily
that $P_1$ is displayed by $\QQ$. But
$P_1\cup\{b\}$ properly contains $Q_3\cup Q_4\udots Q_k$. 
Hence $P_1=Q_3\cup Q_4\udots Q_k$.
Also by the properties of a feral display,
$\lambda_M(P_2\cup P_3\udots P_i)=2$ and 
$\lambda_M(P_{i+1}\cup P_{i+2}\udots P_l)=2$.
Note that, using the notation of the definition of feral
display we have $Z_1\cup Z_2=\{b\}$.

Either $b\in Q_1$ or $b\in Q_2$. 
Assume for a contradiction that $b\in Q_2$, so that 
$(Q_2,Q_1)=
(\{b\}\cup P_2\cup P_3\udots P_i,P_{i+1}\cup P_{i+2}\udots P_l)$.
By a property of feral elements $a$ coblocks either 
$Q_1$ or $Q_2$. But $Q_1=P_{i+1}\cup P_{i+2}\udots P_l$ and
it was observed above that this set is 3-separating in 
$M$. Hence $a$ coblocks $Q_2$.
Therefore $a\in\cl_M(Q_3\cup Q_4\udots Q_k\cup Q_1)$.
But $a\notin\cl_M(P_{i+1}\cup P_{i+2}\udots P_l)=\cl_M(Q_1)$.
As $a$ blocks $P_1\cup\{b\}$, we see that 
$a\notin\cl_M(P_1)$, and this set contains 
$Q_3\cup Q_4\udots Q_l$ so that $a\notin\cl_M(Q_3\cup Q_4\udots Q_k)$. 
Hence, by Lemma~\ref{pi-minor}
\begin{align*}
&\sqcap_{M/a}(Q_1,Q_3\cup Q_4\udots Q_k)\\
=&\sqcap_M(Q_1,Q_3\cup Q_4\udots Q_k)+1\\
\geq &\sqcap_M(P_l,P_1)+1\\
=&2,
\end{align*}
contradicting the fact that $(Q_1,Q_2,Q_3\cup Q_4\udots Q_k)$
is a swirl-like flower of $M/a$.

It follows that $b\in Q_1$. Thus 
$Q_2=P_2\cup P_3\udots P_i$, 
$Q_1=P_{i+1}\cup P_{i+2}\udots P_l\cup\{b\}$,
and $P_1=Q_3\cup Q_4\udots Q_k$. Recall that 
$(P_1,P_2,\ldots,P_l)$ is a maximal $k$-fracture of $M/b\ba a$.

We next show that $i_1=i+1$.
We first show that $i_1\geq i+1$. Assume otherwise,
so that $i_1<i+1$. We have
$a\in\cl_{M/b}(P_{i_1}\cup P_{{i_1}+1}\udots P_{i_2})$. Also, by 
a property of feral elements, $a$ coblocks either $Q_2$
or $Q_1$. But $Q_2=P_2\cup P_2\udots P_i$ and this 
set is 3-separating in $M$. Hence 
$a$ coblocks $Q_1$, so 
$a\in\cl_{M/b}(P_{i+1}\cup P_{i+2}\udots P_l)$.
Now, by Lemma~\ref{modular}, 
$a\in\cl_{M/b}(P_{i+1}\cup P_{i+2}\udots P_{i_2})$.
Therefore 
$(P_1,P_2,\ldots,P_{i_1},\ldots,P_{i+1}\cup P_{i+2}\udots 
P_{i_2}\cup \{a\},\ldots,P_l)$
is a flower of $M/b$ that 
properly refines 
$(P_1,P_2,\ldots,P_{i_1}\udots P_{i_2}\cup \{a\},\ldots, P_l)$
The former flower is clearly tight 
contradicting the fact established in \ref{canon1.2}
that the latter flower is maximal. 
Thus $i_1\geq i+1$.

Assume that $i_1>i+1$. As $a$ coblocks $Q_1$, we see that
$a\in\cl_M(P_1\cup P_2\udots P_i)$. 
Thus $a\in\cl_{M/b}(P_1\cup P_2\udots P_i)$.
Also $a\in\cl_{M/b}(P_{i_1}\cup P_{{i_1}+1}\udots P_l\cup P_1)$. 
Hence,
by Lemma~\ref{modular}, $a\in\cl_{M/b}(P_1)$ so that
$a\in\cl_M(P_1\cup \{b\})$ contradicting the fact 
that $a$ blocks this
set. Thus $i_1=i+1$.

We now prove that $i_2=l$. Assume otherwise, so that
$i_2<l$. We have 
$a\in\cl_{M/b}(P_{i_1}\cup P_{{i_1}+1}\udots P_{i_2})$
so that $a\in\cl_{M/b}(P_2\cup P_3\udots P_{i_2})$.
Hence $a\in\cl_M(\{b\}\cup P_2\cup P_3\udots P_{i_2})$.
Hence $\lambda_M(P_l\cup P_1)=2$, so that
$\lambda_{M/a}(P_l\cup P_1)=2$. But $P_1=Q_3\cup Q_4\udots Q_k$
and $P_l\subseteq Q_1$. There are several ways to see that
neither $P_l$ nor $Q_1-P_l$ is a set of loose elements of
$\QQ$. One way is to observe that both $|P_l|\geq 2$
and $|Q_1-P_l|\geq 2$, so that, if either of these sets is 
loose, one of the elements is a loose coguts element
and we have contradicited the dual of \ref{canon1.1}.
It follows that $P_l\cup P_1$ does not conform with 
$\QQ$ and we have contradicted the fact this this is a
maximal flower in $M/a$. Hence $i_2=l$.

We have now established that (i) and (ii) hold. Consider (iii).
The partition 
$(P_2,\ldots,P_l,P_1\cup\{b\})$ is a flower in $M\ba a$, and
$a\in\cl(Q_1\cup\{b\})$ so that 
$(P_2,\ldots,P_{k-2},E(M)-Q_2)$ is a flower in $M$ and the first
claim of (iii) holds. As $a$ coblocks $Q_1$, we see that 
$a\in\cl_M(E(M)-Q_1)$. Hence $a$ does not block 
$P_{k-1}\udots P_l$ so that the second claim of (iii) holds.
As $a$ blocks $P_1$, we have $a\in\cl^*(Q_1\cup Q_2)$.
Thus $a$ does not coblock $Q_1\cup Q_2$ and the third
claim of (iii) holds.
\end{proof}

To simplify notation a little, set $Q'_1=P_k\cup\udots P_l$,
so that $Q'_1=Q_1-\{b\}$. 
Summarising some of the information gained
so far, we see that the partition $\PP=(P_1\cup\{b\},P_2,\ldots,P_l)$ 
is a maximal $k$-fracture of
$M\ba a$, the partition $\QQ=(Q'_1\cup\{b\},Q_2,Q_3\udots Q_k)$ 
is a maximal $k$-fracture
of $M/a$, and the partition $(P_2,\ldots,P_{k-1},Q'_1\cup\{a\},P_1)$ 
is a maximal $k$-fracture of $M/b$. We next prove

\begin{sublemma}
\label{canon1.7}
The matroid $M\ba b$ is $3$-connected and 
$\OO=(Q_2,P_k,P_{k+1},\ldots,P_l,Q_k,Q_{k-1},\ldots,Q_3\cup\{a\})$
is a maximal $k$-fracture of $M\ba b$. Moreover 
$a$ is in the
coguts of $Q_3\cup\{a\}$ and $Q_2$ in $\OO$.
\end{sublemma}

\subproof
We omit the routine verification that $M\ba b$ is 3-connected.
We first show that $\OO$ is a flower in 
$M\ba b$. To do this, by 
\cite[Lemma~4.11(i)]{flower}, it
suffices to show the union of each pair of consecutive members 
of the partition $\OO$ in the linear order is 3-separating 
in $M\ba b$. In other words we need to check all consecutive
pairs, but we do not have to check that
$(Q_3\cup\{a\})\cup Q_2$ is 3-separating.

Consider $Q_2\cup P_k$. As $b$ is in the coguts of $P_1$ and $P_2$
in $M\ba a$, we have $\lambda_{M\ba a,b}(P_2\udots P_k)=1$, that is,
$\lambda_{M\ba a,b}(Q_2\cup Q_3\cup P_k)=1$, so that 
$\lambda_{M\ba b}(Q_2\cup P_k)\leq 2$ and, as $M\ba b$ is 
3-connected, equality holds here. If $i\in\{k,k+1,\ldots,l-1\}$,
then it follows from \ref{canon1.6}(iii) that 
$\lambda_{M\ba b}(P_i\cup P_{i+1})=2$. 

Consider $P_l\cup Q_k$. We have $\lambda_{M/a}(Q_1\cup Q_k)=2$,
so that $\lambda_{M/a\ba b}(Q_1'\cup Q_k)=2$. Also 
$b\in\cl^*_M(Q_2\cup\{a\})$, so $a\in\cl_M^*(Q_2\cup\{b\})$,
that is, $a\in\cl^*_{M\ba b}(Q_2)$. 
Hence $\lambda_{M\ba b}(Q_1'\cup Q_k)=2$.
We also know that $\lambda_{M\ba a,b}(P_l\cup P_1)=1$, so that 
$\lambda_{M\ba b}(P_l\cup P_1)=2$. As 
$(P_l\cup P_1)\cap (Q'_1\cup Q_k)=P_l\cup Q_k$, 
an uncrossing argument gives 
$\lambda_{M\ba b}(P_l\cup Q_k)=2$ as required. 
By \ref{canon1.6}(iii),
if $i\in\{k,k-1,\ldots,4\}$, 
then $\lambda_{M\ba b}(Q_i\cup Q_{i-1})=2$. 

It remains to prove that 
$\lambda_{M\ba b}(Q_4\cup Q_3\cup\{a\})=2$.
We will use an uncrossing argument in $M\ba a,b$ to
show that $\lambda_{M\ba a,b}(Q_4\cup Q_3)\leq 1$.
We first show that $\lambda_{M\ba a,b}(Q_2\cup Q_3\cup Q_4)\leq 2$.
We know that $\lambda_{M/a\ba b}(Q_2\cup Q_3\cup Q_4)=2$,
as $b\in\cl_{M/a}(Q'_1)$. By \ref{canon1.3}(iii), 
$a$ blocks $P_2\cup\{b\}$, so
$a\notin\cl_M(Q_5\cup Q_6\udots Q_k\cup Q'_1)$.
Hence $a$ does not coblock the 3-separation
$(Q_2\cup Q_3\cup Q_4,Q_5\cup Q_6\udots Q_k\cup Q'_1)$
of $M/a\ba b$. It follows from this fact that the above 3-separation
is induced in $M\ba b$ and hence 
$\lambda_{M\ba a,b}(Q_2\cup Q_3\cup Q_4)\leq 2$.
As $b$ is in the coguts of $P_1$ and $P_2$ in $M\ba a$,
we see that $\lambda_{M\ba a,b}(P_1)=1$,
in other words, $\lambda_{M\ba a,b}(Q_3\cup Q_4\udots Q_k)=1$.
Also, $b\in\cl_{M\ba a}(P_1\cup P_2)$, and $M\ba a$ is
3-connected, so $\lambda_{M\ba a,b}(Q'_1)=\lambda_{M\ba a}(Q'_1)=2$;
that is, $\lambda_{M\ba a,b}(Q_2\cup Q_3\udots Q_k)=2$.
Altogether we have;
$\lambda_{M\ba a,b}(Q_2\cup Q_3\cup Q_4)\leq 2$;
$\lambda_{M\ba a,b}(Q_3\cup Q_4\udots Q_k)=1$;
and $\lambda_{M\ba a,b}(Q_2\cup Q_3\udots Q_k)=2$. Uncrossing
these separations gives us the desired outcome that
$\lambda_{M\ba a,b}(Q_3\cup Q_4)\leq 1$. It follows from this
fact, and the fact that $M\ba b$ is 3-connected, that
$\lambda_{M\ba b}(Q_3\cup Q_4\cup\{a\})=2$, as required.

We have proved that $\OO$ is a flower in $M\ba b$. 
As $\PP$ is a swirl-like flower in
$M\ba a$, we have 
$\sqcap(Q_1',Q_2)=\sqcap(P_2\udots P_{k-1},P_k\udots P_l)=1$. 
Thus $\OO$ is a swirl-like or spike like. 
As $a\notin\cl(P_1)$, and $a\in \cl(Q_1\cup Q_2)$, we see that
$\sqcap_{M/a}(Q_3,Q_k)=\sqcap_M(Q_3,Q_k)$. 
But $\sqcap_{M/a}(Q_3,Q_k)=0$ and it follows that
$\OO$ is swirl-like.
\end{proof}

Relabel $(P_2,\ldots,P_{k-1})$ by $(R_1,\ldots,R_{k-2})$,
relabel $(P_k,\ldots,P_l)$ by $(S_1,\ldots S_r)$, and 
relabel $(Q_3,\ldots,Q_k)$
by $(T_1,\ldots,T_{k-2})$. It is now easily checked that 
these partitions form a bogan display for $\{a,b\}$, so that
$\{a,b\}$ is indeed a bogan couple. 
\end{proof}

\section{Gangs of Three}

We know that if $x$ is an element of the skeleton $M$ that
is comparable with some other element of $M$,
then either $M\ba x$ or $M/x$ is a $k$-skeleton unless $x$ belongs 
to a clonal pair. What if $x$ is not comparable to another
element and $M/x$ is $k$-coherent? In this situation it will
almost always be the case that $M/x$ is a $k$-skeleton. But there is
one exceptional structure for which this is not true and we
describe this structure now. Recall that a 
$3$-connected matroid $M$ is {\em uniquely} $k$-fractured it there
is a bloom $\FF$ of $M$ such that every $k$-fracture of
$M$ is displayed by this bloom.

Let $M$ be a $k$-coherent matroid and $\{r,s,t\}\subseteq E(M)$.
Then $\{r,s,t\}$ is a {\em gang of three}
\index{gang of three} 
in $M$ if there is a
partition $(R,S,T,Z,\{r,s,t\})$ of $E(M)$, and partitions
$(R_2,R_3,\ldots,R_{k-1})$, $(S_2,S_3,\ldots,S_{k-1})$ and
$(T_2,T_3,\ldots,T_{k-1})$, of $R$, $S$ and $T$ respectively
such that the following hold.
\begin{itemize}
\item[(i)] $M/r$, $M/s$ and $M/t$ are $k$-coherent.
\item[(ii)] $M\ba r$, $M\ba s$ and $M\ba t$ are 
3-connected.
\item[(iii)] $(R_2,R_3,\ldots,R_{k-1},E(M)-R)$,
$(S_2,S_3,\ldots,S_{k-1},E(M)-S)$ and $(T_2,T_3,\ldots,T_{k-1},E(M)-T)$
are tight swirl-like flowers in $M$.
\item[(iv)] $(\{s,t\},R_2,R_3,\ldots,R_{k-1},S\cup T\cup Z)$,
$(\{r,t\},S_2,S_3,\ldots,S_{k-1},R\cup T\cup Z)$ and
$(\{r,s\},T_2,T_3,\ldots,T_{k-1},R\cup S\cup Z)$ are 
canonical maximal $k$-fractures
of $M\ba r$, $M\ba s$ and $M\ba t$ respectively. Morever
these $k$-fractures are unique.
\item[(v)] $r$, $s$ and $t$ are fixed in $M$.
\end{itemize}

\begin{figure}
\begin{tikzpicture}[thick,line join=round]
	\coordinate (o) at (0,0);
	\coordinate (a) at (-3, 1);
	\coordinate (b) at (3,1);
	\coordinate (c) at ($(a) + (0,4.5)$);
	\coordinate (d) at (0,3);
	\coordinate (e) at ($(b) + (0,4)$);
	\coordinate (f) at ($(c)!0.6!(e)$);
	\coordinate (g) at ($(d)!0.35!(e)$);
	\coordinate (h) at ($(c)!0.35!(d)$);
	\coordinate (i) at ($(o) + (45:1.2)$);
	\coordinate (j) at ($(d) + (-45:1.2)$);
	\coordinate (k) at ($(a) + (120:1.5)$);
	\coordinate (l) at ($(c) + (-120:1.5)$);
	\coordinate (m) at ($(b) + (45:1.4)$);
	\coordinate (n) at ($(e) + (-45:1.4)$);
	\coordinate (en) at ($(e)!0.5!(n)$);
	\coordinate (nm) at ($(n)!0.5!(m)$);
	\coordinate (mb) at ($(m)!0.5!(b)$);
	\coordinate (b') at ($(b) + (45:-0.1)$);
	\coordinate (oi) at ($(o)!0.5!(i)$);
	\coordinate (ij) at ($(i)!0.5!(j)$);
	\coordinate (jd) at ($(j)!0.5!(d)$);
	\coordinate (o') at ($(o) + (45:-0.1)$);
	\coordinate (i') at ($(i) - (0,0.1)$);
	\coordinate (ak) at ($(a)!0.5!(k)$);
	\coordinate (kl) at ($(k)!0.5!(l)$);
	\coordinate (lc) at ($(l)!0.5!(c)$);
	\coordinate (k') at ($(k) + (120:0.1)$);
	\coordinate[label=-135:$r$] (r) at (intersection of e--h and d--f);
	\coordinate[label=-135:$s$] (s) at (intersection of d--f and c--g);
	\coordinate[label=above:$t$] (t) at (intersection of e--h and c--g);
	\draw (d) -- (f);
	\draw (c) -- (d) -- (e) -- cycle;
	\filldraw[pattern color=orlines,pattern=north west lines] (0,0.5) ellipse (4 and 1.1);
	\draw[white, line width=2mm,line join=round] (e) -- (b) -- (m);
	\fill[white] let \p1 = ($(m) - (mb)$),
		\n1 = {veclen(\x1,\y1)} in
			(b') arc (90-45-180:270-45-180:\n1+1mm);
	\filldraw[pattern color=glines,pattern=north west lines] let \p1 = ($(m) - (mb)$),
		\n1 = {veclen(\x1,\y1)} in
			(b) arc (90-45-180:270-45-180:\n1);
	\filldraw[pattern color=glines,pattern=north west lines] let \p1 = ($(e) - (en)$),
		\n1 = {veclen(\x1,\y1)} in
			(n) arc (-45:180-45:\n1);
	\filldraw[pattern color=glines,pattern=north west lines] let \p1 = ($(n) - (nm)$),
		\n1 = {veclen(\x1,\y1)} in
			(m) arc (-90:90:\n1);
	\draw[white, line width=2mm,line join=round] (d) -- (o) -- (i) -- (j);
	\fill[white] let \p1 = ($(o) - (oi)$),
		\n1 = {veclen(\x1,\y1)} in
			(o') arc (90-45-180:270-45-180:\n1+1mm);
	\fill[white] let \p1 = ($(i) - (ij)$),
		\n1 = {veclen(\x1,\y1)} in
			(i') arc (-90:90:\n1+1mm);
	\filldraw[pattern color=rlines,pattern=north west lines] let \p1 = ($(o) - (oi)$),
		\n1 = {veclen(\x1,\y1)} in
			(o) arc (90-45-180:270-45-180:\n1);
	\filldraw[pattern color=rlines,pattern=north west lines] let \p1 = ($(d) - (jd)$),
		\n1 = {veclen(\x1,\y1)} in
			(j) arc (-45:180-45:\n1);
	\filldraw[pattern color=rlines,pattern=north west lines] let \p1 = ($(i) - (ij)$),
		\n1 = {veclen(\x1,\y1)} in
			(i) arc (-90:90:\n1);
	\draw[white, line width=2mm,line join=round] (k) --(a) --(c);
	\fill[white] let \p1 = ($(a) - (ak)$),
		\n1 = {veclen(\x1,\y1)} in
			(k') arc (120:180+120:\n1+1mm);
	\filldraw[pattern color=olines,pattern=north west lines] let \p1 = ($(c) - (lc)$),
		\n1 = {veclen(\x1,\y1)} in
			(l) arc (180+60:60:\n1);
	\filldraw[pattern color=olines,pattern=north west lines] let \p1 = ($(l) - (kl)$),
		\n1 = {veclen(\x1,\y1)} in
			(k) arc (270:90:\n1);
	\filldraw[pattern color=olines,pattern=north west lines] let \p1 = ($(k) - (ak)$),
		\n1 = {veclen(\x1,\y1)} in
			(k) arc (120:180+120:\n1);
	\draw (g) -- (c) -- (l) -- (k) -- (a) -- (c) -- (d) -- (o) -- (i) -- (j) -- (d) -- (e) -- (b) -- (m) -- (n) -- (e) -- (h);
	\foreach \pt in {o,a,b,r,s,t} \fill[black] (\pt) circle (3pt);
	\node at ($(en) + (45:0.35)$) [rectangle,fill=white,draw=white,inner sep=0.5pt] {\textcolor{glabels}{$S_2$}};
	\node at ($(nm) + (0:0.5)$) [rectangle,fill=white,draw=white,inner sep=0.5pt] {\textcolor{glabels}{$S_3$}};
	\node at ($(mb) + (-45:0.35)$) [rectangle,fill=white,draw=white,inner sep=0.5pt] {\textcolor{glabels}{$S_4$}};
	\node at ($(jd) + (45:0.3)$) [rectangle,fill=white,draw=white,inner sep=0pt] {\textcolor{rlabels}{$T_2$}};
	\node at ($(ij) + (0:0.3)$) [rectangle,fill=white,draw=white,inner sep=0.5pt] {\textcolor{rlabels}{$T_3$}};
	\node at ($(oi) + (-45:0.3)$) [rectangle,fill=white,draw=white,inner sep=0pt] {\textcolor{rlabels}{$T_4$}};
	\node at ($(lc) + (160:0.35)$) [rectangle,fill=white,draw=white,inner sep=0.5pt] {\textcolor{olabels}{$R_2$}};
	\node at ($(kl) + (180:0.4)$) [rectangle,fill=white,draw=white,inner sep=0.5pt] {\textcolor{olabels}{$R_3$}};
	\node at ($(ak) + (-160:0.35)$) [rectangle,fill=white,draw=white,inner sep=0.5pt] {\textcolor{olabels}{$R_4$}};
	\node at (2,0.4) [rectangle,fill=white,draw=white] {\textcolor{orlabels}{$Z$}};
\end{tikzpicture}
\caption{A Gang of Three}\label{gangof3}
\end{figure}
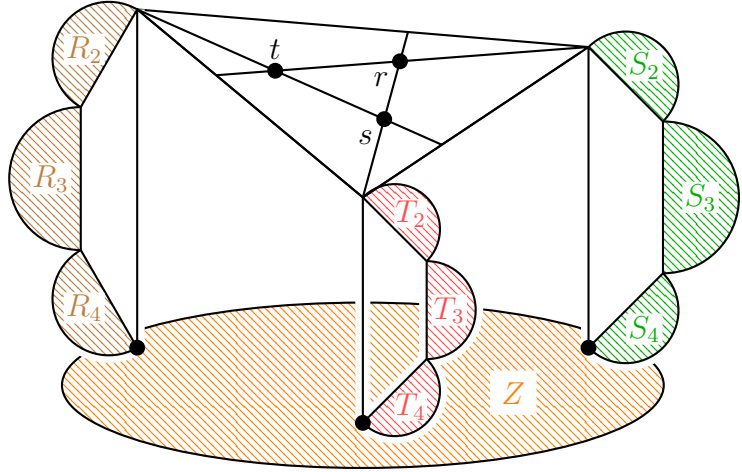

Figure~\ref{gangof3} illustrates a gang of three with notation as
in the definition.
Note that gangs of three are not self dual. A gang of three in $M^*$
is a {\em cogang of three}.
\index{cogang of three} 
The goal of this section is to prove

\begin{theorem}
\label{contract}
Let $x$ be an element of the $k$-skeleton $M$ such that 
$M/x$ is $k$-coherent
and  $x$ is not comparable with any other element of $M$. Then 
either $M/x$ is a $k$-skeleton, or $x$ is a member of a
gang of three in $M$.
\end{theorem}

\begin{figure}
\begin{tikzpicture}[thick,line join=round]
	\coordinate (Eo) at (0,-2);
	\coordinate (Ef) at ($(Eo) + (0,2.25)$);
	\coordinate (Ea) at ($(Ef) + (0:2.625)$);
	\coordinate (Eb) at ($(Ef) + (45:2.625)$);
	\coordinate (Ec) at ($(Ef) + (90:2.625)$);
	\coordinate (Ed) at ($(Ef) + (135:2.625)$);
	\coordinate (Ee) at ($(Ef) + (180:2.625)$);
	\coordinate (Eg) at ($(Ef) + (225:2.625)$);
	\coordinate (Eh) at ($(Ef) + (-45:2.625)$);
	\coordinate (Ei) at ($(Ee) + (0.75,-3)$);
	\coordinate (Ej) at ($(Ea) + (-0.75,-3)$);
	\coordinate (Eab) at ($(Ea)!0.5!(Eb)$);
	\coordinate (Ebc) at ($(Eb)!0.5!(Ec)$);
	\coordinate (Ecd) at ($(Ec)!0.5!(Ed)$);
	\coordinate (Ede) at ($(Ed)!0.5!(Ee)$);
	\coordinate (Eei) at ($(Ee)!0.5!(Ei)$);
	\coordinate (Eij) at ($(Ei)!0.5!(Ej)$);
	\coordinate (Eja) at ($(Ej)!0.5!(Ea)$);
	\coordinate (Ef') at ($(Eo)!0.5!(Ec)$);
	\node[pattern color=glines,pattern=north west lines,draw,circle through=(Eb)] at (Eab) {};
	\node[pattern color=glines,pattern=north west lines,draw,circle through=(Ec)] at (Ebc) {};
	\node[pattern color=glines,pattern=north west lines,draw,circle through=(Ed)] at (Ecd) {};
	\node[pattern color=glines,pattern=north west lines,draw,circle through=(Ee)] at (Ede) {};
	\node[pattern color=rlines,pattern=north west lines,draw,circle through=(Ei)] at (Eei) {};
	\node[pattern color=rlines,pattern=north west lines,draw,circle through=(Ej)] at (Eij) {};
	\node[pattern color=rlines,pattern=north west lines,draw,circle through=(Ea)] at (Eja) {};
	\filldraw[fill=white] (Ee) -- (Eg) -- (Eo) -- (Eh) -- (Ea) -- (Eb) -- (Ec) -- (Ed) -- (Ee) -- (Ei) -- (Ej) -- (Ea);
	\draw (Eo) -- (Ec);
	\filldraw[pattern=north west lines,pattern color=olines] (Ee) .. controls ($(Ee) + (1,0)$) and ($(Eg) + (1,1)$) .. (Eg);
	\filldraw[pattern=north west lines,pattern color=olines] (Eg) .. controls ($(Eg) + (1,1)$) and ($(Eo) + (0,1)$) .. (Eo);
	\filldraw[pattern=north west lines,pattern color=olines] (Ea) .. controls ($(Ea) + (-1,0)$) and ($(Eh) + (-1,1)$) .. (Eh);
	\filldraw[pattern=north west lines,pattern color=olines] (Eh) .. controls ($(Eh) + (-1,1)$) and ($(Eo) + (0,1)$) .. (Eo);
	\node at ($(Ef') + (0:0.4)$) {$f$};
	\fill[black] (Ef') circle (3pt);
	\node at ($(Ede) + (157.5:0.5)$) [rectangle,fill=white,draw=white,inner sep=1] {\textcolor{glabels}{$P_2$}};
	\node at ($(Ecd) + (112.5:0.5)$) [rectangle,fill=white,draw=white,inner sep=1] {\textcolor{glabels}{$P_3$}};
	\node at ($(Ebc) + (67.5:0.5)$) [rectangle,fill=white,draw=white,inner sep=1] {\textcolor{glabels}{$P_4$}};
	\node at ($(Eab) + (22.5:0.5)$) [rectangle,fill=white,draw=white,inner sep=1] {\textcolor{glabels}{$P_5$}};
	\node at ($(Eei) + (202.5:0.75)$) [rectangle,fill=white,draw=white,inner sep=1] {\textcolor{rlabels}{$Q_3$}};
	\node at ($(Eij) + (270:1)$) [rectangle,fill=white,draw=white,inner sep=1] {\textcolor{rlabels}{$Q_4$}};
	\node at ($(Eja) + (-22.5:0.75)$) [rectangle,fill=white,draw=white,inner sep=1] {\textcolor{rlabels}{$Q_5$}};
	\node at ($(Ee) + (0.75,-0.75)$) [rectangle,fill=white,draw=white,inner sep=1] {\textcolor{olabels}{$Z_1$}};
	\node at ($(Eo) + (-0.8,0.55)$) [rectangle,fill=white,draw=white,inner sep=1] {\textcolor{olabels}{$Z_2$}};
	\node at ($(Eo) + (0.8,0.55)$) [rectangle,fill=white,draw=white,inner sep=1] {\textcolor{olabels}{$Z_3$}};
	\node at ($(Ea) + (-0.75,-0.75)$) [rectangle,fill=white,draw=white,inner sep=1] {\textcolor{olabels}{$Z_4$}};
\end{tikzpicture}
\caption{$M\ba f$ is not uniquely fractured}
\label{not-unique}
\end{figure}
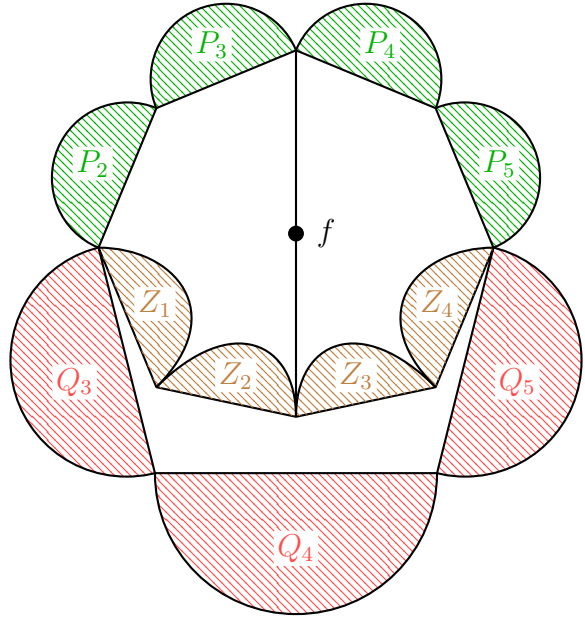

We first prove some preliminary lemmas.
One might expect that, if $f$ is feral, then $M\ba f$ would
be uniquely $k$-fractured. Perhaps surprisingly this is not always true;
a feral element for which $M\ba f$ is not uniquely $k$-fractured
is illustrated in Figure~\ref{not-unique}.
But it is true for members of bogan couples. 

\begin{lemma}
\label{uniquely-bogan}
If $M$ is $k$-coherent and $a$ is a member of a bogan couple
in $M$, then 
$M\ba a$ and $M/a$ are uniquely $k$-fractured.
\end{lemma}

\begin{proof}
Say that $\{a,b\}$ is a bogan couple of $M$. 
Then it has an associated bogan display. In what follows
we use the notation for a bogan display
given in the definition. We know that
$(R_1,R_2,\ldots,R_{k-2},S,T\cup\{b\})$ is a $k$-fracture of
$M\ba a$. Say $(P_1,P_2,\ldots,P_l)$ is another $k$-fracture.
By Lemma~\ref{in-petal}, we may assume that 
$P_1\cup P_2\udots P_{l-1}$ is contained in the full closure of a 
petal of $(R_1,R_2,\ldots,R_{k-2},S,T\cup\{b\})$. This 
petal must be blocked by $a$, otherwise $M$ is $k$-fractured.
The only possibilities are that 
$P_1\cup P_2\udots P_{l-1}\subseteq R_1\cup\{b\}$
or that $P_1\cup P_2\udots P_{l-1}\subseteq T\cup\{b\}$. 
But, in either case
we see that $b\in\cl^*(P_l)$ so that we may assume that
$P_1\cup P_2\udots P_{l-1}$ is contained in $T$ or $R_1$. As neither
$T$ nor $R_1$ is blocked by $a$ we again get the contradiction
that $M$ is $k$-fractured. Thus $M\ba a$ is uniquely $k$-fractured.
The argument for $M/a$ is similar and is omitted.
\end{proof}

The next lemma shows that, while an element can be loose in more
than one flower of a matroid, this can only happen in a very
restricted way. The first part of the lemma is true for any type
of flowers, but would take slightly longer to prove. Recall that,
if $Q_i$ is a petal of a swirl-like flower $\PP$, then 
$Q_i^+$ denotes the ordered sequences of loose 
elements between $Q_i$ and $Q_{i+1}$ and $Q_i^-$ denotes the 
ordered sequence of loose elements between $Q_{i-1}$
and $Q_i$.

\begin{lemma}
\label{one-only}
Let $\PP$ and 
$\QQ=(Q_1,Q_2,\ldots,Q_n)$ be inequivalent 
tight maximal swirl-like 
flowers of the $3$-connected matroid $M$ of order at least
three. Assume that the 
petal $P$ of $\PP$ is contained in $\hQ_i$ for some 
$i\in\{1,2,\ldots,n\}$. Then $P$ contains at most one element $z$ 
of $Q_i^+$. In the case that 
$P$ contains such an element $z$, then $z$ is the first
element of $Q_i^+$.
\end{lemma} 

\begin{proof}
By Lemma~\ref{in-petal}, $\hQ_i$ contains all but one petal of 
$\PP$. It follows that $P$ is not displayed by $\QQ$. 
The lemma now follows from Lemma~\ref{3-cross}.
\end{proof}

The following lemma is an easy consequence of the structure of
$k$-wild triangles.

\begin{lemma}
\label{not-awkward}
Let $T$ be a $k$-wild triangle of the $k$-coherent matroid
$M$, let $z_1,z_2$ be elements of $E(M)-T$. Assume that
$N\in \{M\ba z_1,z_2,M\ba z_1/z_2\}$ and that $N$
is $3$-connected. Then $T$ is not a set of
loose elements of a swirl-like flower of $N$.
\end{lemma}



The next two lemmas provide the bulk of the
proof of Theorem~\ref{contract}.
Lemma~\ref{del-del} below is the more straightforward. Note that
we use the {\em dual} of Lemma~\ref{del-del} for Theorem~\ref{contract}.

\begin{lemma}
\label{del-del}
Let $x$ be an element of the $k$-skeleton  $M$ 
such that $x$ is not comparable with any other element and
$M\ba x$ is $k$-coherent. Then there is no element 
$y\in E(M\ba x)$ such that $y$ is fixed in $M\ba x$ and 
$M\ba x,y$ is $k$-coherent.
\end{lemma}

\begin{proof}
Assume that $y$ is fixed in $M\ba x$ and that $M\ba x,y$
is $k$-coherent. Then $M\ba y$ is 3-connected.
Then
$y$ is fixed in $M$ so that $M\ba y$ is $k$-fractured. Let
$\FF=(P_1,P_2,\ldots,P_m)$ be a maximal 
$k$-fracture of $M\ba y$, where
$x\in P_1$. By Lemma~\ref{gain-coherence}, $m=k$, and 
$(P_1-\{x\},P_2,\ldots,P_k)$ is a swirl-like quasi-flower of 
$M\ba x,y$ and $P_1$ is a set of loose elements of this
quasi-flower.

\begin{sublemma}
\label{contract1}
$P_1$ contains no triangles. 
\end{sublemma}

\subproof
Assume that $T$ is a triangle in $P_1$. Say $x\notin T$. By 
Lemma~\ref{loose2}, 
$T$ contains an element that is fixed in $M\ba x,y$ and consequently
fixed in $M$. But, by Lemma~\ref{not-awkward}, $T$ is not
a $k$-wild triangle of $M$. 
But then, by Lemma~\ref{skeleton-triangle},
$T$ is a clonal triple in $M$, and does not contain an
element that is fixed in $M$.
On the other hand, if $x\in T$, then, as $x$ is not fixed in $M$,
we see, again by Lemma~\ref{skeleton-triangle}, that the triangle
$T$ is a clonal triple, contradicting the fact that $x$ is not 
comparable with any other element of $M$.
\end{proof}

By \ref{contract1}, and the fact that
$P_1-\{x\}$ is a fan in $M\ba x,y$, we have $|P_1|\leq 4$.

\begin{sublemma}
\label{contract2}
$|P_1|\neq 2$.
\end{sublemma}

\subproof
Assume that $|P_1|=2$, say $P_1=\{x,z\}$. 
As $x$ is not fixed in $M$, there is a matroid
$M'$ obtained by independently cloning
$x$ by $x'$. By Lemma~\ref{freedom4},
$\{x,x',z\}$ is a triangle in $M'\ba y$ and hence in
$M'$. This proves that 
$z\less x$ in $M$, contradicting the fact that $x$ 
is not comparable
with any other element of $M$.
\end{proof}

\begin{sublemma}
\label{contract3}
$|P_1|\neq 3$.
\end{sublemma}

\subproof
Say $|P_1|=3$. Then, as $P_1$ is $3$--separating in $M\ba y$ and
$M\ba x,y$ is $3$-connected, $P_1$ is a triangle of $M\ba y$ and
hence of $M$, contradicting \ref{contract1}.
\end{proof}

\begin{sublemma}
\label{contract4}
$|P_1|\neq 4$.
\end{sublemma}

\subproof
Assume that $|P_1|=4$.
By \ref{contract1} $P_1-\{x\}$ is a triad of $M\ba x,y$. Say
$P_1-\{x\}=\{a,b,c\}$, where $(a,b,c)$ is a fan between
$P_k$ and $P_2$ in $M\ba x,y$. By Lemma~\ref{loose2}, 
$a$ and $c$ are cofixed in $M\ba x,y$. As $y$ is fixed in
$M\ba x$, we see, by Lemma\ref{born-free}(ii), 
that $a$ and $c$ are cofixed in
$M\ba x$. But $x$ is not comparable with any element of
$M$, so, again by Corollary~\ref{born-free}(ii), we conclude
that $a$ and $c$ are cofixed in $M$.

Say $p\in\{a,c\}$. By Lemma~\ref{loose1}, 
$M\ba x,y/p$ is $k$-coherent. By 
\ref{contract1}, $x$ is not in a triangle of $M\ba y$, so
$M\ba y/p$ is $3$-connected. Now, by 
Lemma~\ref{gain-coherence}, either $M\ba y/p$
is uniquely fractured by $(P_1-\{p\},P_2,\ldots,P_k)$,
or $M\ba y/p$ is $k$-coherent and $P_1-\{p\}$ is a 
loose petal of $(P_1-\{p\},P_2,\ldots,P_k)$. It is conceivable that
$p$ is in a triangle $T$ of $M$. Assume that this is the case.
As $p$ is cofixed in $M$, the triangle
$T$ cannot be a clonal triple, so, by 
Lemma~\ref{skeleton-triangle}, $T$ must be $k$-wild.
Readers who find is obvious that $T$ cannot be $k$-wild
should skip the remainder of this paragraph. Note that
$T$ meets $P_1$ in $\{p\}$ as otherwise $P_1$ contains a triangle.
Observe that $p\in\cl_M(P_1-\{p\})$; that is, there is a 
3-element subset $Z$ of $E(M)-T$ such that $p\in\cl_M(Z)$.
It is easily checked from properties of a $k$-wild display
that this is only possible if $T$ is a standard $k$-wild
triangle. In this case $\si(M/p)$ is not $3$-connected
by Lemma~\ref{k-wild-different}(iii). But $M\ba y/p$
is $3$-connected, so $\si(M/p)$ is $3$-connected. Hence
$p$ is not in a triangle.

It follows that $M/p$ is $3$-connected. 
As $p$ is cofixed in $M$,
$M/p$ is not $k$-coherent. It now follows from 
either Lemma~\ref{gain-coherence} (in the case that 
$M\ba y/p$ is $k$-coherent) or Corollary~\ref{unique-fracture}
(in the case that $M\ba y/p$ is uniquely $k$-fractured) that
there is an $i$ and $j$ in $\{1,2,\ldots,k\}$ 
such that $y\in\cl_{M/a}(\hP_i)$
and $y\in\cl_{M/c}(\hP_j)$. Hence $y\in\cl_M(\{a\}\cup \hP_i)$
and $y\in\cl_M(\{c\}\cup\hP_j)$. If $i=1$, then
$y\in\cl_M(P_1)$ and $M$ is $k$-fractured. Hence
$i,j\in\{2,3,\ldots,k\}$.

We now show that $2\in\{i,j\}$. Assume otherwise.  Let 
$P'_3=P_3\cup P_4\udots P_k$. 
Then $y\in\cl(\{a\}\cup P'_3)$, and $y\in\cl(\{c\}\cup P'_3)$.
If $y\in\cl(P'_3)$, then, by Lemma~\ref{modular}, $y\in\cl(P_i)$,
contradicting the fact that $M$ is $k$-coherent. Hence
$c\in\cl(P'_3\cup\{y\})$ and $a\in\cl(P'_3\cup \{y\})$, so that
$c\in\cl(P'_3\cup\{a\})$.
Now consider
the swirl-like flower $(\{a,b,c\},P_2,\ldots,P_k)$ of $M\ba x,y$.
Here  $(a,b,c)$ is a fan between $P_k$
and $P_2$ in $M\ba x,y$, and $c$ is a rim element of this fan.
This implies that $b\in\cl_{M\ba x,y}(P_2\cup\{c\})$, so that 
$c\in\cl_M(P_2\cup \{b\})$.
Hence $c\notin\cl^*_M(\{a\}\cup P'_3)$.
From this contradiction it follows 
that $2\in\{i,j\}$, and, similarly, $k\in\{i,j\}$.
But now, $y\in\cl(P_1\cup P_2)$ and $y\in\cl(P_1\cup P_k)$,
so, by Lemma~\ref{modular}, $y\in\cl(P_1)$ implying that
$M$ is $k$-fractured.
\end{proof}

The lemma now follows from \ref{contract2}, \ref{contract3} 
and \ref{contract4}.
\end{proof}

We now come to the more substantial lemma.

\begin{lemma}
\label{con-del}
Let $x$ be an element of the $k$-skeleton $M$ such that 
$x$ is not comparable with any other element of $M$ and
$M/x$ is $k$-coherent. If there is an element $y\in E(M/x)$
such that $y$ is fixed in $M/x$ and $M/x\ba y$ is $k$-coherent,
then there is exactly one more such element $z$ and
$\{x,y,z\}$ is a gang of three in $M$.
\end{lemma}

\begin{proof}
Assume that $M/x$ is $k$-coherent, that
$y$ is fixed in $M/x$, and that $M/x\ba y$ is $k$-coherent.

\begin{sublemma}
\label{contract5}
$y$ is fixed in $M$, and $M\ba y$ is $3$-connected and $k$-fractured.
\end{sublemma}

\subproof
As $x$ is not comparable with any other element of $M$, 
we see,
by Corollary~\ref{born-free}(i), that $y$ is fixed in $M$. 
Thus $M\ba y$ is not 
$k$-coherent. If $M\ba y$ is not 3-connected, then, as 
$M\ba y/x$ is $3$-connected, $y$ is in a triad with $x$.
But $x$ is not cofixed in $M$, so, by 
Corollary~\ref{free-triangle} the triad is not $k$-wild, and 
$x$ is not comparable with any other element of $M$, so the triad
cannot be a clonal triple. It now follows by 
Lemma~\ref{skeleton-triangle} that $y$ is not in a triad
so that $M\ba y$ is $3$-connected
and $k$-fractured.
\end{proof}

Let $\PP=(P_1,P_2,\ldots,P_m)$ be a $k$-fracture of $M\ba y$ where
$x\in P_1$ and $P_1$ is fully closed. As $M\ba y/x$ is $k$-coherent,
it follows from Lemma~\ref{gain-coherence}(ii) 
that $m=k$ and that $(P_1,P_2,\ldots,P_k)$
uniquely fractures $M\ba y$. Also, as $M\ba y/x$ is $k$-coherent,
$P_1-\{x\}$ is a fan of loose elements in the quasi-flower
$(P_1-\{x\},P_2,\ldots,P_k)$ of $M\ba y/x$.

\begin{sublemma}
\label{contract6}
\begin{itemize}
\item[(i)] $\PP$ has no loose elements in $P_1$.
\item[(ii)] The fan $P_1-\{x\}$ in $M\ba y/x$ 
begins and ends with spoke elements.
\item[(iii)] $P_1$ has an even number of elements.
\end{itemize}
\end{sublemma}

\subproof
If $P_1$ has a loose element $z$, then,
by Theorem~\ref{canon1}, either $y\more z$
or $\{y,z\}$ is a bogan couple. The former
case contradicts \ref{contract5}. Hence
$\{y,z\}$ is a bogan couple. Again, by
Theorem~\ref{canon1}, $z$ is a loose coguts element
of $P_1$. It follows from the properties of a bogan display
that $y$ does not block $P_1-\{z\}$. As $M/x$ is 
$k$-coherent, $z\neq x$.

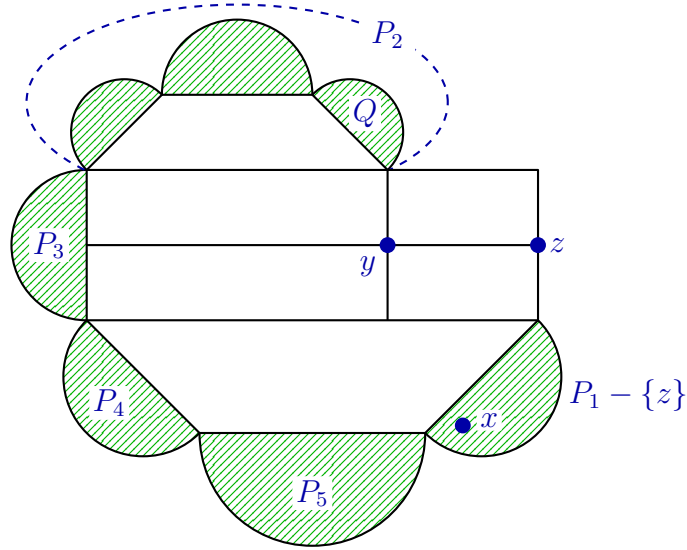
\begin{figure}
\begin{tikzpicture}[thick,line join=round]
	\coordinate (o) at (0,0);
	\coordinate (p) at (0,2);
	\coordinate (q) at (4,2);
	\coordinate (r) at (4,0);
	\coordinate (s) at (6,0);
	\coordinate (t) at (6,2);
	\coordinate (n) at (0,1);
	\coordinate (a) at ($(p)+(45:1.414)$);
	\coordinate (b) at ($(q)+(135:1.414)$);
	\coordinate (c) at (-45:2.121);
	\coordinate (d) at ($(s)+(-135:2.121)$);
	\node[pattern color=lines,draw,circle through=(b),pattern=north east lines] at ($(b)!0.5!(q)$) {};
	\node[pattern color=lines,draw,circle through=(a),pattern=north east lines] at ($(a)!0.5!(b)$) {};
	\node[pattern color=lines,draw,circle through=(p),pattern=north east lines] at ($(p)!0.5!(a)$) {};
	\node[pattern color=lines,draw,circle through=(o),pattern=north east lines] at ($(o)!0.5!(p)$) {};
	\node[pattern color=lines,draw,circle through=(c),pattern=north east lines] at ($(c)!0.5!(o)$) {};
	\node[pattern color=lines,draw,circle through=(d),pattern=north east lines] at ($(d)!0.5!(c)$) {};
	\node[pattern color=lines,draw,circle through=(s),pattern=north east lines] at ($(s)!0.5!(d)$) {};
	\draw[dashed,labels] ($(a)!0.5!(b)-(0,0.1)$) ellipse (2.8 and 1.3);
	\filldraw[fill=white] (o) -- (p) -- (a) -- (b) -- (q) -- (t) -- (s) -- (d) -- (c) -- cycle;
	\coordinate[label=-135:{\textcolor{labels}{$y$}}] (y) at (4,1);
	\coordinate[label=right:{\textcolor{labels}{$z$}}] (z) at (6,1);
	\draw (p) -- (q) -- (r);
	\draw (n) -- (z);
	\draw (o) -- (s);
	\node at (4,3.8) [rectangle,fill=white,draw=white] {\textcolor{labels}{$P_2$}};
	\node at (3.7,2.75) [rectangle,fill=white,draw=white,inner sep=1pt] {\textcolor{labels}{$Q$}};
	\node at (-0.5,1) [rectangle,fill=white,draw=white,inner sep=1pt] {\textcolor{labels}{$P_3$}};
	\node at (0.3,-1.1) [rectangle,fill=white,draw=white,inner sep=1pt] {\textcolor{labels}{$P_4$}};
	\node at (3,-2.3) [rectangle,fill=white,draw=white,inner sep=1pt] {\textcolor{labels}{$P_5$}};
	\coordinate (x) at (5,-1.4);
	\foreach \pt in {x,y,z} \fill[labels] (\pt) circle (3pt);
	\node at ($(x)+(10pt,2pt)$) [rectangle,fill=white,draw=white,inner sep=1pt] {\textcolor{labels}{$x$}};
	\node at (7.2,-1) [rectangle,fill=white,draw=white,inner sep=1pt] {\textcolor{labels}{$P_1-\{z\}$}};
\end{tikzpicture}
\caption{The Labelling for the Bogan Couple of \ref{contract6}}
\label{bogan-labelling}
\end{figure}

Consider a bogan display for $\{z,y\}$. 
By the properties of a bogan display, either 
$P_1-\{z\}$ or $P_2$ contains $k-2$ petals of a 
$k$-fracture for $M\ba z$. 
As $P_1-\{x\}$ is a fan in 
$M\ba y/x$, the latter case must hold.
An illustration of the bogan display that we now have is given
for the case $k=5$ in Figure~\ref{bogan-labelling}.
By considering the $k$-fracture of $M\ba z$ displayed by
the bogan display, and the properties of bogan displays as 
given in the definition, we observe that there is a
subset $Q$ of
$P_2$ such that $y\in\cl(P_4\cup P_5\udots P_k\cup P_1\cup Q)$
and $z\notin\cl(P_4\cup P_5\udots P_k\cup P_1\cup Q)$. Therefore
$y\in\cl_{M/x}(P_4\cup P_5\udots P_k\cup(P_1-\{x\})\cup Q)$
and $z\notin\cl_{M/x}(P_4\udots P_k\cup (P_1-\{x\})\cup Q)$,
so $y\not\more z$ in $M/x$. Moreover $z$ is clearly cofixed in
$M/x\ba y$, and now, by Corollary~\ref{born-free}, 
$z$ is cofixed in $M/x$.

By Lemma~\ref{uniquely-bogan}, 
$(P_1-\{z\},P_2,P_3\cup\{y\},\ldots,P_k)$
uniquely fractures $M/z$. Also $P_1-\{z,x\}$ 
is a set of loose elements in the quasi-flower 
$(P_1-\{x,z\},P_2,P_3\cup\{y\},\ldots,P_k)$
of $M/z,x$. Hence $M/z,x$ is $k$-coherent. 
As $z$ is cofixed in $M/x$,
we have contradicted the dual of Lemma~\ref{del-del}.
It follows that $P_1$ has no loose elements and (i) holds.

Consider (ii). We know that $P_1-\{x\}$ 
is a fan between $P_k$ and $P_1$
in $M\ba y/x$. If an end of this fan, say $z$, is a rim element,
then $z\in\cl^*_{M\ba y/x}(P_2)$ for some $i\in\{k,2\}$.
Hence $z\in\cl^*_{M\ba y}(P_2)$, contradicting (i). 
Thus (ii) holds. Part (iii)
follows immediately. 
\end{proof}

By \ref{contract6}(i), we may assume that $P_2$ and $P_k$ are fully closed.
The easy proof of the next sublemma is omitted.

\begin{sublemma}
\label{contract7}
\begin{itemize}
\item[(i)] $P_1$ contains no triangles in $M$.
\item[(ii)] If $T$ is a triad of $P_1-\{x\}$
in $M\ba y/x$, then $T$ is blocked by $y$.
\end{itemize}
\end{sublemma}

Consider the fan $P_1-\{x\}$ between $P_k$ and $P_2$ in
$M\ba y/x$. Denote this fan by $(z_k,f_1,f_2,\ldots,f_l,z_2)$. If
$|P_2|=2$, then $z_k=z_2$ and $\{f_1,f_2,\ldots,f_l\}=\emptyset$.
Otherwise, by \ref{contract6}(iii), $l$ is odd. By \ref{contract6}(ii),
$z_k\in\cl_{M\ba y/x}(P_k)$ and $z_2\in\cl_{M\ba y/x}(P_2)$.

\begin{sublemma}
\label{contract8}
$z_k\notin\cl_{M\ba y}(P_k)$ and $z_2\notin\cl_{M\ba y}(P_2)$.
\end{sublemma}

\subproof
If the sublemma fails, then either $z_k$ or $z_2$ is a loose element of
$P_1$ contradicting \ref{contract6}(i).
\end{proof}

\begin{sublemma}
\label{contract9}
Both $z_k$ and $z_2$ are fixed in $M$.
\end{sublemma}

\subproof
Say $i\in\{k,2\}$. As $z_i$ is a spoke element of a fan between 
$P_k$ and $P_2$ in $M\ba y/x$, we see, by Lemma~\ref{loose2}(iii),
that $z_i$ is fixed in $M\ba y/x$.
Hence $z_i$ is fixed in $M/x$. By hypothesis, 
$x$ is not comparable with any other element of $M$, so by 
Corollary~\ref{im-free},
$z_i$ is fixed in $M$.
\end{proof}

Our next task is to reduce to the case that $|P_1|=2$. Until further notice
assume that $|P_1|>2$.

\begin{sublemma}
\label{contract10}
$\{z_k,z_2,x\}$ is not a triad of $M\ba y$.
\end{sublemma}

\subproof
Assume that $\{z_k,z_2,x\}$ is a triad of $M\ba y$. 
We have $z_2\in\cl_{M\ba y/x}(P_2)$ and by \ref{contract8}, 
$z_2\notin\cl_{M\ba y}(P_2)$, so $x\in\cl_{M\ba y}(P_2\cup\{z_2\})$.
Say $z_k\in\cl_{M\ba y}(P_2\cup\{z_2\})$. Then, as
$\sqcap(P_2,\{z_k,z_2,x\})\leq \sqcap(P_2,P_1)=1$, 
an easy rank calculation shows that $\{z_k,z_2,x\}$ is a
triangle in $M$, contradicting \ref{contract7}(i). Thus
$z_k\notin\cl_{M\ba y}(P_2\cup\{z_2\})$ and also
$z_2\notin\cl_{M\ba y}(P_2\cup\{z_k\})$.

Next we show that $z_2$ and $z_k$ are cofixed in $M$.
As $x\in\cl_{M\ba y}(P_2\cup\{z_2\})$, but 
$z_k\notin\cl_{M\ba y}(P_2\cup\{z_2\})$, there is a 
cyclic flat of $M\ba y$ that contains $x$ but not $z_k$.
Hence there is a cyclic flat of $(M\ba y)^*$ that contains
$z_k$ but not $x$. But $\{z_k,z_2,x\}$ is a triangle of 
$(M\ba y)^*$, so by Lemma~\ref{freedom3}, $z_k$
is fixed in $(M\ba y)^*$. Thus $z_k$, and similarly $z_2$,
are cofixed in $M\ba y$. Since $y$ is fixed in $M$, 
Corollary~\ref{born-free} implies that 
$z_k$ and $z_2$ are cofixed in $M$.

We now show that $M/z_k$ is 3-connected. We know that $z_k$ and $z_2$ 
are spoke ends of a fan between $P_k$ and $P_2$ in $M\ba y/x$.
Moreover $P_1$ is fully closed in $M\ba y$ so that
$P_1-\{x\}$ is fully closed in $M\ba y/x$. Hence the fan
is maximal so 
that $M\ba y/x\ba z_k,z_2$ is 3-connected. As $\{x,z_k,z_2\}$
is a triad of $M\ba y$, we see that
$M\ba y/x\ba z_k,z_2=M\ba y/z_k\ba x,z_2$.
Hence $M\ba y/z_k\ba x,y_1,z_2$ is 3-connected and $M/z_k$
is 3-connected up to parallel classes. If $M/z_k$ is not 
$3$-connected, then $z_k$ is in a triangle $T$ of $M$.
By \ref{contract7}(i), and the fact that
$P_1$ is fully closed in $M\ba y$, the triangle $T$ contains $y$.
But $y$ is fixed in $M$, so $T$ is $k$-wild. As 
$\si(M/z_k)$ is 3-connected, $T$ is not standard. Assume that
$T$ is costandard. By Lemma~\ref{cowild-win},
$\co(M\ba y,z_k)$ is not $3$-connected. But 
$M\ba y,z_k/x$ is 3-connected, so $\co(M\ba y,z_k)$ is
$3$-connected and $T$ is not costandard. Hence
$M/z_k$, and 
similarly $M/z_2$, is 3-connected.
Since $z_k$ and $z_2$ are cofixed in $M$ and $M/z_k$ and $M/z_2$
are 3-connected, we deduce that 
both $M/z_k$ and $M/z_2$ are $k$-fractured.

Consider the flower $(P_1-\{z_k\},P_2,\ldots,P_k)$ of
$M\ba y/z_k$. We now show that $P_1-\{z_k\}$ is not a
loose petal of this flower. Assume that it is a loose petal.
Then there is an element $r\in P_1-\{z_k\}$ such that.
$r\in\cl^{(*)}_{M\ba y/z_k}(P_2)$.
Assume that $r\in\cl_{M\ba y/z_k}(P_2)$. Then
$r\in\cl_{M\ba y}(P_2\cup\{z_k\})$. As $\{z_k,z_2,x\}$ is a triad of
$M\ba y$, we see that $r\in\{z_2,x\}$. If $r=x$, then 
$z_k\in\cl_{M\ba y/x}(P_2)$, so $r\neq x$. Assume that $r=z_2$.
Then $z_2\in\cl_{M\ba y}(P_2\cup\{z_k\})$, and hence
$z_k\in\cl_{M\ba y}(P_2\cup\{z_2\})$. But also, 
$x\in\cl_{M\ba y}(P_2\cup\{z_2\})$ and we have already seen that 
this situation does not occur. Thus $r\notin\cl_{M\ba y/z_k}(P_2)$.
Then $r\in\cl^*_{M\ba y/z_k}(P_2)$ so that
$r\in\cl^*_{M\ba y}(P_2)$ contradicting \ref{contract6}(i).
Therefore $P_1-\{z_k\}$ is not a loose petal of 
$(P_1-\{z_k\},P_2,\ldots,P_k)$.

Assume that the flower $(P_1-\{z_k\},P_2,\ldots,P_k)$ 
of $M\ba y/z_k$ is not maximal.
In this case it is routinely seen that 
there is a  partition $(P',P'')$ of $P_1-\{z_k\}$,
such that $(P',P'',P_2,\ldots,P_k)$ is a tight swirl-like flower
in $M\ba y/z_k$. As $z_k\in\cl_{M\ba y/x}(P_k)$
and $z_k\notin\cl_{M\ba y/x}(P_k)$, we know that 
$x\in\cl_{M\ba y/z_k}(P_k)$. Say $z_2\in P'$.
Then, as $\sqcap_{M\ba y/z_k}(P_k\cup P')=0$,
we see that $z_2\notin\cl_{M\ba y/z_k,x}(P_2)$,
contradicting the fact that $z_2\in\cl_{M\ba y/x}(P_2)$.
Say $z_2\in P''$. Then, as $\sqcap(P_k\cup\{x\},P''\cup P_2)=0$,
and $z_2\notin\cl_{M\ba y/z_k}(P_2)$, we also have 
$z_2\notin\cl_{M\ba y/z_k,x}(P_2)$, and again we
contradict the fact that $z_2\in\cl_{M\ba y/x}(P_2)$.

Therefore $(P_1-\{z_k\},P_2,\ldots,P_k)$ is a maximal $k$-fracture
of $M\ba y/z_k$ and it is easily seen that this fracture is unique.
Similarly $(P_1-\{z_2\},P_2,\ldots,P_k)$ is  the unique 
maximal $k$-fracture of $M\ba y/z_2$. In what follows we discuss
related flowers in different matroids. By $\hP_i$ we will
always mean $\fcl_{M\ba y}(P_i)$. Observe that,
if $i\in\{2,3,\ldots,k\}$, then $\fcl_{M\ba y/z_k}(P_i)=\hP_i$
and, if $i\in\{1,2,\ldots,k-1\}$, then 
$\fcl_{M\ba y/z_2}(P_i)=\hP_i$.

As $M/z_k$ is uniquely $k$-fractured,
it follows by Corollary~\ref{unique-fracture} that 
$y\in\cl_{M/z_k}(\fcl_{M\ba y/z_k}(P_i))$
for some $i\in\{2,3,\ldots,k\}$, or 
$y\in\cl_{M/z_k}(\fcl_{M\ba y/z_k}(P_1-\{z_k\}))$.
If $y\in\cl_{M/z_k}(\fcl_{M\ba y/z_k}(P_1-\{z_k\}))$, then 
$y\in\cl_M(P_1)$ and it follows that $M$ is $k$-fractured.
Thus the former case occurs and we also have
$y\in\cl_{M/z_2}(\fcl_{M\ba y/z_2}(P_i))$
for some $i\in\{2,3,\ldots,k\}$.

Assume that $y\in\cl_{M/z_k}(\hP_i)$ for some 
$i\in\{2,3,\ldots,k-1\}$. 
Let $P'_2=\hP_2\cup\hP_3\udots\hP_{k-1}$.
Then $y\in\cl_M(P'_2\cup\{z_k\})$,
and, as $z_k\notin\cl_M(P'_2)$, we have
$z_k\in\cl_M(P'_2\cup\{y\})$. 
Observe that 
$\cl_M(P'_2\cup\{y\})\cap P_1=\cl_M(P'_2\cup\{z_k\})\cap P_1$.
As $z_2\in\cl_{M/x}(P'_2)$, and $z_2\notin\cl_M(P'_2)$,
we see that, if $z_2\in\cl_M(P'_2\cup\{y\})$,
then $x\in\cl_M(P'_1\cup\{y\})$. As $\sqcap(P'_2\cup\{y\},P_1)=2$,
we deduce that $\{x,z_2,z_k\}$ is a triangle in $M$. 
This contradicts a number of things, amongst which is the fact that
$M/x$ is 3-connected.

It follows tht $z_2\in\cl_M(\hP_k\cup\{y\})$. But we may apply
the previous argument using $z_2$ to deduce that 
$z_k\in\cl_M(\hP_2\cup\{y\})$.
Thus $y\in\cl_M(P_1\cup \hP_k)$ and $y\in\cl_M(P_1\cup \hP_2)$,
so by Lemma~\ref{modular}, $y\in\cl_M(P_1)$ contradicting the fact that
$M$ is $k$-coherent.
At last we can conclude that $\{z_k,z_2,x\}$ is not a triad
of $M\ba y$.
\end{proof} 

\begin{sublemma}
\label{contract11}
$M\ba y,z_k$ and $M\ba y,z_2$ are $3$-connected.
\end{sublemma}

\subproof 
Certainly $M\ba y/x\ba z_k$ is $3$-connected, so that
$M\ba y,z_k$ is 3-connected up to a series pair containing
$x$. Such a series pair must be $\{x,z_2\}$, otherwise
$z_2\in\cl_{M\ba y}(P_2)$ contradicting \ref{contract8}.
But, in this case, $\{x,z_k,z_2\}$ is a triad of $M\ba y$, 
contradicting \ref{contract10}. Hence $M\ba y,z_k$ has no
series pair containing $x$ and is therefore 3-connected.
\end{proof}

Consider the flower 
$(P_1-\{z_k\},P_2,\ldots,P_k)$ in $M\ba y,z_k$.

\begin{sublemma}
\label{contract12}
The ordered elements of $P_1-\{z_k\}$ that are loose elements between
$P_k$ and $P_1-\{z_k\}$ 
in $M\ba y,z_k$ form an initial segment of
$(f_1,\ldots,f_l,z_2,x)$.
\end{sublemma}

\subproof
Say $r\in P_1-\{z_k\}$ and $r\in\cl^{(*)}_{M\ba y,z_k}(P_k)$.
By \ref{contract6}(i), $r\notin\cl(P_k)$, so 
$r\in\cl^*_{M\ba y,z_k}(P_k)$. As $x\in\cl(P_2\cup\{z_2\})$,
and $P_2\cup\{z\}\subseteq E(M\ba y,z_k)-P_k$, we see that
$x\notin\cl^*_{M\ba y,z_k}(P_k)$ so that $r\neq x$. 
Hence $r\in\cl^*_{M\ba y,z_k/x}(P_k)$, that is,
$r\in\cl^*_{M\ba y/x}(P_k\cup\{z_k\})$. Thus $r=f_1$.

Assume that $l>1$ and that there is an element 
$r\in\cl_{M\ba y,z_k}(P_k\cup\{f_1\})$. If $r=x$, then 
$x\in\cl_M(P_k\cup\{f_1\})$, so that $f_1\in\cl_{M/x}(P_k)$.
But $z_2\in\cl_{M/x}(P_k)$ and we have contradicted the structure 
of loose elements in swirl-like flowers. Hence $r\neq x$.
Now $r\in\cl_{M\ba y,z_k/x}(P_k\cup\{f_1\})$, so $r=f_2$. 

An easy induction now proves that the loose elements of $P_k$
form an initial segment of $(f_1,f_2,\ldots,f_l,r,s)$, where
$(r,s)$ is a permutation of $\{z_2,x\}$. Assume that the 
elements $\{f_1,f_2,\ldots,f_l\}$ 
are all loose and that we can continue.
Note that there is a circuit $C$ such that
$\{x,z_2\}\subseteq C\subseteq P_2\cup\{x,z_2\}$.
Thus $r\notin \cl^*_{M\ba y,z_k}(P_k\cup\{f_1,f_2,\ldots,f_l\})$,
so that
$r\in\cl_{M\ba y,z_k}(P_k\cup\{f_1,f_2,\ldots,f_l\})$.
Say $r=x$, so that 
$x\in\cl_{M\ba y,z_k}(P_k\cup\{f_1,f_2,\ldots,f_l\})$.
But $z_k\in\cl_{M\ba y/x}(P_k\cup\{f_1,f_2,\ldots,f_l\})$
and $z_2\in\cl_{M\ba y/x}(P_k\cup\{f_1,f_2,\ldots,f_l,z_k\})$,
so $z_2\in\cl_{M\ba y/x}(P_k\cup\{f_1,f_2,\ldots,f_l\})$.
Hence $z_2\in\cl_{M\ba y,z_k}(P_k\cup\{f_1,f_2,\ldots,f_l,x\})$.
Thus 
$\{z_2,x\}\subseteq \cl_{M\ba y,z_k}(P_k\cup\{f_1,f_2,\ldots,f_l\})$,
contradicting the structure of loose elements in a swirl-like flower.
Hence $r=z_2$ and the claim holds.
\end{proof}

Now consider the flower $(P_1-\{z_2\},P_2,\ldots,P_k)$ of 
$M\ba y,z_2$.

\begin{sublemma}
\label{contract13}
The ordered subset of elements of $P_1-\{z_2\}$ that are loose
elements between $P_k$ and $P_1-\{z_2\}$
in $M\ba y,z_2$ is either empty or the first element
of the set is $x$, in which case $x$ is in the coguts of $P_k$
and $P_1-\{z_2\}$.
\end{sublemma}

\subproof
If $r\in P_k-\{z_2\}$ and $r\in\cl^{(*)}_{M\ba y,z_2}(P_k)$,
then $r\in\cl^*_{M\ba y,z_2}(P_k)$. Say $r\neq x$. Then
$r\in\cl^*_{M\ba y,z_2/x}(P_k)$ and 
$r\in\cl^*_{M\ba y/x}(P_k\cup\{z_2\})$. 
But $P_1-\{z\}$ is a fan in $M\ba y/x$ that begins and ends
at the spoke elements $z_k$ and $z_2$. 
Thus $(M\ba y/x)|(P_1-\{x\})$ is connected.
Also $z_2\in\cl_{M\ba y/x}(P_2)$. It follows that 
there is a circuit $C$ of 
$M\ba y/x$ such that $r\in C\subseteq P_2\cup(P_1-\{x\})$.
Hence $r\notin\cl^*_{M\ba y/x}(P_k)$. The claim follows from
this contradiction.
\end{proof}

\begin{sublemma}
\label{contract14}
If $P_1-\{z_2\}$ is a tight petal of $(P_1-\{z_2\},P_2,\ldots,P_k)$
in $M\ba y,z_2$, then the set of elements in $P_1-\{z_2\}$
that are loose between $P_k$ and $P_1-\{z_2\}$ is empty.
\end{sublemma}

\subproof
Assume that the set is nonempty. Then, by \ref{contract13}, 
$x$ is in the coguts of $P_1-\{z_2\}$ and $P_k$. Now 
$P_1-\{x,z_2\}$ is a tight petal of the
flower $(P_1-\{x,z_2\},P_2,\ldots,P_k)$ in the 3-connected matroid
$M\ba y,z_2/x$, so that $P_1-\{x\}$ is certainly not a fan
of loose elements
in the flower $(P_1-\{x\},P_2,\ldots,P_k)$ of $M\ba y/x$.
\end{proof}

As both $M\ba y,z_k$
and $M\ba y, z_2$ are 3-connected, both
$M\ba z_k$ and $M\ba z_2$ are 3-connected. As $z_k$ and $z_2$ 
are fixed in $M$, both $M\ba z_k$ and $M\ba z_2$ are
$k$-fractured.

\begin{sublemma}
\label{contract15}
We may assume that $y\notin\cl(P_1\cup P_2)$.
\end{sublemma}

\subproof
If $y\in\cl(P_k\cup P_1)$ and $y\in\cl(P_1\cup P_2)$, then
by Lemma~\ref{modular}, $y\in\cl(P_1)$, contradicting the fact that
$M$ is $k$-coherent. Thus, up to symmetry, we may assume that
$y\notin \cl(P_1\cup P_2)$.
\end{proof}

\begin{sublemma}
\label{contract16}
Say $i\in\{k,2\}$. Assume that $P_1-\{z_i\}$ is a fan of
loose elements between
$P_k$ and $P_2$ in the flower $(P_1-\{z_i\},P_2,\ldots,P_k)$ of
$M\ba y,z_i$, with initial element $\alpha$.
Then $y\in\cl(P_k\cup\{\alpha\})$.
\end{sublemma}

\subproof
Under the hypothesis of the sublemma it is clear that
$M\ba y,z_i$ is $k$-coherent. It is easily checked that if
$(P_1-\{z_i\},P_2,\ldots,P_k)$
induces a $k$-fracture in $M\ba z_i$, then it follows that 
$M$ is $k$-fractured. 
As $M\ba z_i$ is $k$-fractured, by Lemma~\ref{lose-coherence}
it must be the case that
for some petals $\hQ_1,\hQ_2$ of some 
other swirl-like flower of $M\ba y,z_i$ of order $k-1$,
we have $y\in\cl(\hQ_1),\cl(\hQ_2)$. By Lemma~\ref{in-petal}, 
we may assume that
$\hQ_1\subseteq \fcl_{M\ba y,z_i}(P_j)$ 
for some $j\in\{2,3,\ldots,k\}$. 
The only cases that 
do not quickly lead to a contradiction to the fact that $M$ is 
$k$-coherent are when either $\hQ_1\subseteq \fcl_{M\ba y,z_i}(P_k)$ 
or $\hQ_1\subseteq \fcl_{M\ba y,z_i}(P_2)$.
The latter case contradicts \ref{contract15}, so the former case 
holds. Now by Lemma~\ref{one-only} we have $\hQ_1\subseteq P_k\cup\{\alpha\}$,
so that $y\in\cl(P_k\cup\{\alpha\})$.
\end{proof}

\begin{sublemma}
\label{contract17}
Say $i\in\{k,2\}$. Assume that $P_1-\{z_i\}$ is a tight petal of
$(P_1-\{z_i\},P_2,\ldots,P_k)$ in $M\ba y,z_i$. Then
$y\in\cl(\fcl_{M\ba y,z_i}(P_k))$.
\end{sublemma}

\subproof 
We have $(P_1-\{z_i\},P_2,\ldots,P_k)$
uniquely fractures $M\ba y,z_i$ and $M\ba z_i$ is
$k$-fractured. By Lemma~\ref{in-closure}, either $y\in\cl(P_1-\{z_i\})$
or $y\in\cl(\fcl_{M\ba y,z_i}(P_j))$ for some $j\in\{2,3,\ldots,k\}$. 
The only 
case that does not either contradict \ref{contract15} or the 
$k$-coherence of $M$ is if $y\in\cl(\fcl_{M\ba y,z_i}(P_k))$.
\end{proof}

\begin{sublemma}
\label{contract18}
$x\in\cl(P_k\cup\{y\})$.
\end{sublemma}

\subproof
Say that $P_1-\{z\}$ is a tight petal of the flower
$(P_1-\{z_2\},P_2,\ldots,P_k)$ in $M\ba y,z_2$. 
By \ref{contract17} we have $y\in\cl(\fcl_{M\ba y,z_2}(P_k))$,
and by \ref{contract14}, we may assume that
$\fcl_{M\ba y,z_2}(P_k)=P_k$, so $y\in\cl(P_k)$,
giving the contradiction that $M$ is $k$-fractured.
Thus $P_1-\{z_2\}$ is a fan between 
$P_k$ and $P_1$ in $M\ba y,z_2$. 
By \ref{contract13} and \ref{contract16}, 
$y\in\cl(P_k\cup\{x\})$. As $y\notin\cl(P_k)$,
we have $x\in\cl(P_k\cup\{y\})$ as required.
\end{proof}

Assume that $P_1-\{z_k\}$ is a fan between $P_k$ and $P_1$
in $(P_1-\{z_k\},P_2,\ldots,P_k)$. Then by \ref{contract16},
$y\in\cl(P_k\cup\{f_1\})$ and now by \ref{contract18}
$x\in\cl(P_k\cup\{f_1\})$. But, by \ref{contract12}, $x$ is the last
element of the fan. Moreover, $f_1$ is an initial rim element of the
fan and it is easily seen that the fan structure implies the
contradiction that $x\notin\cl(P_k\cup\{f_1\})$.

Thus we may assume that $P_1-\{z_k\}$ is a 
tight petal. By \ref{contract12},
the loose elements of $P_1-\{z_k\}$ between $P_k$ and $P_1-\{z_k\}$ 
are of the form $(f_1,\ldots,f_i)$ for some $i\leq l$. Moreover,
$x\notin\cl(P_k\cup\{f_1,\ldots,f_i\})$, otherwise $x$ is also a loose
element. By \ref{contract16}, $y\in\cl(P_k\cup\{f_1,\ldots,f_i\})$,
so by \ref{contract18} $x\in\cl(P_k\cup\{f_1,\ldots,f_i\})$.

From this final contradiction we at last deduce that $|P_1|=2$ and we 
consider this case now. Assume that $P_1=\{x,z\}$. The task now
is to show that $\{x,y,z\}$ is a gang of three.

\begin{sublemma}
\label{contract19}
$x$, $y$ and $z$ are fixed in $M$.
\end{sublemma}

\subproof
We already know that $y$ is fixed in $M$. Say $\{a,b\}=\{x,z\}$.
Then $b$ is fixed in $M\ba y/a$ as $a$ is in the guts of a pair of 
petals of a swirl-like flower in $M\ba y/a$. 
Thus $b$ is fixed in $M/a$.
As $a$ and $b$ are not comparable in $M$, it follows from
Corollary~\ref{im-free} that $b$ is fixed in $M$. Therefore both
$x$ and $z$ are fixed in $M$.
\end{proof}

The next claim is clear.

\begin{sublemma}
\label{contract20}
$M\ba x$, $M\ba y$ and $M\ba z$ are $3$-connected and 
$k$-fractured.
\end{sublemma}

Let $\OO=(O_1,O_2,\ldots,O_l)$ and 
$\QQ=(Q_1,Q_2,\ldots,Q_m)$ be $k$-fractures
of $M\ba x$ and $M\ba z$ respectively. 
As $\PP$ uniquely fractures $M\ba y$,
it follows by Lemma~\ref{2-element-fracture} 
that $M\ba y,x$ and $M\ba y,z$ are 
$k$-coherent. Moreover $(P_2,P_3,\ldots,P_k\cup\{z\})$ and
$(P_2,P_3,\ldots,P_k\cup\{x\})$ are maximal 
swirl-like flowers in $M\ba y,x$
and $M\ba y,z$ respectively. At this stage we almost have symmetry
between $x$ and $z$ except that we do not yet know that
$M/z$ is $k$-coherent. 

By Lemma~\ref{in-closure}, $y\in\cl_M(\fcl_{M\ba y,x}(P_i))$
for some $i\in\{2,3,\ldots,k\}$. 
It is easily seen that we contradict the
$k$-coherence of $M$ unless, up to labels we have
$y\in\cl(P_k\cup\{z\})$, $y\notin\cl(P_k)$, and 
$y\notin\cl(P_2\cup\{z\})$ A similar
conclusion holds by considering the flower 
$(P_2,P_3,\ldots,P_k\cup\{x\})$, establishing the next claim.

\begin{sublemma}
\label{contract21}
We may assume that
$y\in\cl(P_k\cup\{z\})$, $y\in\cl(P_k\cup\{x\})$,$y\notin\cl(P_k)$,
$y\notin\cl(P_2\cup\{z\})$ and $y\notin\cl(P_2\cup\{x\})$.
\end{sublemma}

Consider the $k$-fracture $\OO=(O_1,O_2,\ldots,O_l)$ of $M\ba x$.
We may assume that $O_1$ is fully closed and that $y\in O_1$.
Consider the quasi-flower $(O_1-\{y\},O_2,\ldots,O_l)$
in $M\ba x,y$. As $M\ba x,y$ is $k$-coherent, it follows
from Lemma~\ref{gain-coherence} that $l=k$ and 
that $O_1-\{y\}$ is a loose petal of this quasi-flower.
As $O_1$ is fully closed, $O_1-\{y\}$ is a maximal fan
in $M\ba x,y$. By Lemma~\ref{lose-coherence}
$y\in\cl((O_1-\{y\})\cup O_2)$ and 
$y\in\cl((O_k\cup(O_1-\{y\})$. By Lemma~\ref{in-petal}
we may assume, up to labels, that $(O_1-\{y\})\cup O_2$
is contained in the full closure of some petal of the flower
$(P_2,P_3,\ldots,P_k\cup\{z\})$ of $M\ba x,y$.
As $y\in\cl((O_1-\{y\})\cup O_2)$, this petal must
be $P_k\cup\{z\}$.

\begin{sublemma}
\label{contract21.5}
$P_2\subseteq\fcl_{M\ba x,y}(O_k)$, 
$z\in\cl^*_{M\ba x,y}(O_k)$, and $z\in O_1$.
\end{sublemma}

\subproof
By Lemma~\ref{in-petal}, $P_2$ is contained in the full closure
of a petal $O_i$ of $(O_1-\{y\},O_2,\ldots,O_k)$ in
$M\ba x,y$, where $i\in \{3,4,\ldots,k\}$.
But $z\in\cl^*_{M\ba x,y}(P_2)$, so, 
$z\in\cl^*_{M\ba x,y}(O_i)$.
By the fact that
$(O_1-\{y\},O_2,\ldots,O_k)$ is a swirl-like quasi-flower
and the fact that $z\in\cl^*_{M\ba x,y}(O_i)$, we see
that $i=k$. Thus $P_2\subseteq\fcl_{M\ba x,y}(O_k)$
and $z\in\cl^*_{M\ba x,y}(O_k)$.

If $z\notin O_1-\{y\}$, then $z\in O_{k-1}$. But, in this case,
$z$ is a loose coguts element of the flower $(O_1,O_2,\ldots,O_k)$
in $M\ba x$, contradicting the fact that $M\ba x,z$ is $3$-connected.
Hence $z\in O_1-\{y\}$.
\end{proof}

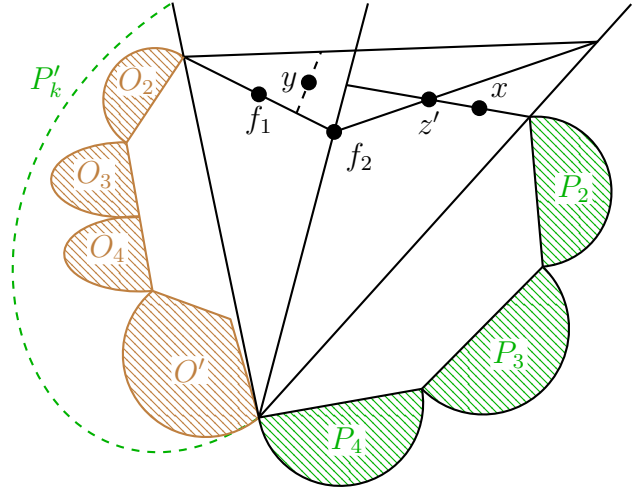
\begin{figure}
\begin{tikzpicture}[thick,line join=round]
	\coordinate (o) at (0,0);
	\coordinate[label=-45:$f_2$] (f2) at (1,3.8);
	\coordinate (a) at (4.5,5);
	\coordinate (c) at (-1,4.8);
	\coordinate (b) at ($(o)!0.8!(a)$);
	\coordinate[label=below:$z'$] (z') at (intersection of f2--a and b--c);
	\coordinate (d) at (intersection of f2--o and b--c);
	\coordinate[label=below:$f_1$] (f1) at ($(f2)!0.5!(c)$);
	\coordinate (e) at ($(f1)!0.5!(f2)$);
	\coordinate (f) at ($(c)!0.8!(intersection of f2--d and a--c)$);
	\coordinate[label=180:$y$] (y) at ($(e)!0.5!(f)$);
	\coordinate[label=45:$x$] (x) at ($(b)!0.5!(z')$);
	\coordinate (h) at (10:2.2);
	\coordinate (i) at ($(b) + (-85:2)$);
	\coordinate (j) at ($(o)!1.15!(c)$);
	\coordinate (oh) at ($(o)!0.5!(h)$);
	\coordinate (hi) at ($(h)!0.5!(i)$);
	\coordinate (bi) at ($(b)!0.5!(i)$);
	\coordinate (k) at (130:2.2);
	\coordinate (l) at ($(k) + (100:1)$);
	\coordinate (m) at ($(k) + (100:2)$);
	\coordinate (cm) at ($(c)!0.5!(m)$);
	\coordinate (ok) at ($(o)!0.5!(k)$);
	\coordinate (n) at ($(k)!0.9!(intersection of k--h and o--c)$);
	\draw[dashed,color=glines] (o) .. controls (-3,-2) and (-5,3) .. (j);
	\filldraw[pattern color=olines,pattern=north west lines,draw=olines] let \p1 = ($(o) - (ok)$),
		\n1 = {veclen(\x1,\y1)} in
			(o) arc (130+180:-50+180:\n1);
	\filldraw[pattern color=olines,pattern=north west lines,draw=olines] let \p1 = ($(c) - (cm)$),
		\n1 = {veclen(\x1,\y1)},
		\n2 = {atan(\x1/\y1)} in
			(c) arc (-\n2-90+180:-\n2+90+180:\n1);
	\filldraw[pattern color=olines,pattern=north west lines,draw=olines] (k) arc (100+180:-90+180:1 and 0.5);
	\filldraw[pattern color=olines,pattern=north west lines,draw=olines] (l) arc (100+180:-90+180:1 and 0.5);
	\filldraw[pattern color=glines,pattern=north west lines] let \p1 = ($(o) - (oh)$),
		\n1 = {veclen(\x1,\y1)} in
			(o) arc (-170:10:\n1);
	\filldraw[pattern color=glines,pattern=north west lines] let \p1 = ($(h) - (hi)$),
		\n1 = {veclen(\x1,\y1)},
		\n2 = {atan(\x1 / \y1)} in
			(h) arc (\n2-180:45.32:\n1);
	\filldraw[pattern color=glines,pattern=north west lines] let \p1 = ($(b) - (bi)$),
		\n1 = {veclen(\x1,\y1)} in
			(i) arc (-85:180-85:\n1);
	\filldraw[pattern color=olines,pattern=north west lines,draw=olines] (o) -- (n) -- (k);
	\draw (j) -- (o) -- (h) -- (i) -- (b) -- (d);
	\draw ($(o)!1.1!(a)$) -- (o) -- ($(o)!1.45!(f2)$);
	\draw (a) -- (c) -- (f2) -- cycle;
	\draw[dashed] (e) -- (f);
	\draw[olines] (k) -- (m) -- (c);
	\foreach \pt in {f1,f2,y,z',x} \fill[black] (\pt) circle (3pt);
	\node at ($(bi) + (5:0.5)$) [rectangle,fill=white,draw=white,inner sep=0.5pt] {\textcolor{glabels}{$P_2$}};
	\node at ($(hi) + (-45:0.5)$) [rectangle,fill=white,draw=white,inner sep=0.5pt] {\textcolor{glabels}{$P_3$}};
	\node at ($(oh) + (-80:0.5)$) [rectangle,fill=white,draw=white,inner sep=0.5pt] {\textcolor{glabels}{$P_4$}};
	\node at ($(ok) + (-135:0.3)$) [rectangle,fill=white,draw=white,inner sep=0.5pt] {\textcolor{olabels}{$O'$}};
	\node at ($(l) + (-135:0.6)$) [rectangle,fill=white,draw=white,inner sep=0.5pt] {\textcolor{olabels}{$O_4$}};
	\node at ($(m) + (-135:0.6)$) [rectangle,fill=white,draw=white,inner sep=0.5pt] {\textcolor{olabels}{$O_3$}};
	\node at ($(cm) + (135:0.35)$) [rectangle,fill=white,draw=white,inner sep=0.1pt] {\textcolor{olabels}{$O_2$}};
	\node at ($(cm) + (135:0.35) - (1.2,0)$) [rectangle,fill=white,draw=white,inner sep=0.1pt] {\textcolor{glabels}{$P'_k$}};
\end{tikzpicture}
\caption{Getting Technical}\label{technical}
\end{figure}

The disturbing possibility that $O_1\neq \{y,z\}$ 
needs to be eliminated. Figure~\ref{technical} is an illustration
for the proof of \ref{contract22}. In the diagram, 
$P'_k=P_k-\{f_1,f_2\}$ and $O'=O_k\cap P_k$.

\begin{sublemma}
\label{contract22} 
$O_1=\{y,z\}$.
\end{sublemma}

\subproof
If the sublemma fails, then $|O_1-\{y\}|>1$.
In this case, by Lemma~\ref{lose-coherence},
$y\in\cl(O_1-\{y\})$. If $|O_1-\{y\}|\neq 3$,
then $O_1$ contains a triangle of $M$.
We omit the routine verification that this cannot happen.
Hence $|O_1-\{y\}|=3$. By \ref{contract21.5},
$z\in\cl^*_{M\ba x,y}(O_k)$, so that there is an
ordering $(z,f_2,f_1)$ of $O_1-\{y\}$ that gives a maximal
fan of loose elements between $O_k$ and $O_2$ in the 
flower $(O_1-\{y\},O_2,\ldots,O_k)$ of $M\ba x,y$.
Note that $f_1$ is a rim element of this fan, so that,
by Lemma~\ref{loose2}(iv), $f_1$ is cofixed in $M\ba x,y$.
We next show that $f_1$ is cofixed in $M$.
Now $x\in\cl(P_2\cup\{z\})$ and $P_2\subseteq O_k$.
By the fact that $(z,f_2,f_1)$ is a fan of loose elements
between $O_k$ and $O_2$, we see that $f_1\notin \cl(O_k\cup \{z\})$.
Thus $f_1\notin \cl(P_2\cup\{z\})$. 
Hence $x\not\more f_1$ in $M\ba y$. Therefore, by 
Corollary~\ref{born-free}, $f_1$ is cofixed
in $M\ba y$. As $y$ is fixed in $M$ it follows, again by
Corollary~\ref{born-free}, that  
$f_1$ is cofixed in $M$.

The element $f_1$ is a terminal rim element of a 
maximal fan between $O_k$ and $O_2$ in the 
$k$-coherent matroid $M\ba x,y$.
Thus, by 
Lemma~\ref{loose1},
$M\ba y,x/f_1$ is $k$-coherent. But 
$f_1\notin\cl_{M}(P_2\cup\{x\})$,
so $x\notin\cl_{M/f_1}(P_2)$. Therefore,
$$x\in\cl_{M/f_1}(P_2\cup\{z\})-\cl_{M/f_1}(P_2).$$ 
As $f_1\in P_k$, by \ref{contract21}, we have 
$$x\in\cl_{M/f_1}((P_k-\{f_1\})\cup\{z\})-\cl_{M/f_1}(P_k-\{f_1\}).$$
It is now easily checked that $(\{x,z\},P_2,\ldots,P_k-\{f_1\})$
uniquely fractures $M\ba y/f_1$. We know that 
$y\in\cl_{M/f_1}((P_k-\{f_1\})\cup\{x,z\})$. But
$y\notin\cl_{M/f_1}(P_1-\{f_1\})$. Assume for a contradiction
that $y\in\cl_{M/f_1}(\{x,z\})$. Then $x\in\cl_{M/f_1}(\{y,z\})$
so that $x\in\cl_M(\{y,z,f_1\})$. 
Recall that $\OO$ is a $k$-fracture of $M\ba x$ and that
$O_1=\{y,z,f_1,f_2\}$. It follows that 
$(O_1\cup\{x\},O_2,\ldots,O_k)$ is a $k$-fracture of $M$.
From this contradiction we deduce that
$y\notin \cl_{M/f_1}(\{x,z\})$.
Now, by Corollary~\ref{unique-fracture}, 
$M/f_1$ is $k$-coherent. 
But $f_1$ is cofixed in $M$ and we have contradicted the fact that 
$M$ is a $k$-skeleton.
\end{proof}

Evidently the same argument works for $\QQ$ so that we have
$\OO=(\{y,z\},O_2,\ldots,O_k)$ and $\QQ=(\{x,y\},Q_2,\ldots,Q_k)$.

\begin{sublemma} 
\label{contract23}
\begin{itemize}
\item[(i)] $\OO$, $\PP$ and $\QQ$ uniquely fracture $M\ba x$,
$M\ba y$ and $M\ba z$ respectively.
\item[(ii)] Up to labels $(O_2,O_3,\ldots,O_k\cup\{x,y,z\})$,
$(P_2,P_3,\ldots,P_k\cup\{x,y,z\})$ and 
$(Q_2,Q_3,\ldots,Q_k\cup\{x,y,z\})$
are swirl-like flowers in $M$.
\item[(iii)] $O_2\cup O_3\udots O_{k-1}$, 
$P_2\cup P_3\udots P_{k-1}$ and
$Q_2\cup Q_3\udots Q_{k-1}$ are mutually disjoint.
\end{itemize}
\end{sublemma}

\subproof
We already know that $\PP$ uniquely fractures $M\ba y$, 
so that $M\ba y,x$
is $k$-coherent. Consider $\OO$. If $\OO$ did not uniquely fracture
$M\ba x$, then $M\ba x,y$ would not be $k$-coherent by 
Lemma~\ref{gain-coherence}.
Thus $\OO$, and similarly $\QQ$ uniquely fracture $M\ba x$ and 
$M\ba z$ respectively.

Consider (ii). As $y\in\cl(P_k\cup\{x,z\})$, we see that
$(P_2,P_3,\ldots,P_k\cup\{x,y,z\})$ is a swirl-like flower in $M$.
Similar arguments to those that establish that
$y\in\cl(P_k\cup\{x,z\})$ prove that, up to labels,
$x\in\cl(O_k\cup\{y,z\})$ and $z\in\cl(Q_k\cup\{x,y\})$.
Thus (ii) holds.

Consider (iii). Consider the inequivalent flowers 
$(O_2,O_3,\ldots,O_k\cup\{z\})$ and $(P_2,P_3,\ldots,P_k\cup\{z\})$
in $M\ba x,y$. We have already observed that, up to labels,
$O_2\subseteq P_k$. Note that $P_k\cup\{z\}$ is fully 
closed in $M\ba x,y$. Thus, by Lemma~\ref{in-petal}, $P_k$
contains all but one member of $(O_2,O_3,\ldots,O_k)$.
As $z$ is in the coguts of $O_2$ and
$O_k\cup\{z\}$ in $M\ba x,y$ we see that $z\in\cl(O_2\cup O_k)$.
so that $y\in\cl(O_2\cup O_k)$. But $y\notin\cl(P_k)$,
so $O_2\cup O_k\not\subseteq P_k$. Hence $O_k\not\subseteq P_k$.
Thus $O_2\udots O_{k-1}\subseteq P_k$, so that
$O_2\cup O_3\udots O_{k-1}$ and 
$P_2\cup P_3\udots P_{k-1}$ are disjoint.
The rest of (iii) follows from similar arguments.
\end{proof}

Note that we now have symmetry between $y$ and $z$.

\begin{sublemma}
\label{contract24}
$M/y$ and $M/z$ are $k$-coherent.
\end{sublemma}

\subproof
By the symmetry between $y$ and $z$ is suffices to prove that
$M/z$ is $k$-coherent. Consider $M\ba x/z$. Clearly
$(P_2,\ldots,P_k\cup\{y\})$ is a swirl-like flower in this matroid.
But $x\in\cl(P_2\cup\{z\})$, so $x\in\cl_{M/x\ba z}(P_2)$.
By \ref{contract21}, $y\in\cl(P_k\cup\{x\})$,
so $x\in\cl(P_k\cup\{y,z\})$. Thus
$(P_2\cup\{x\},P_3,\ldots,P_{k-1},P_k\cup\{y\})$ is a 
swirl-like flower
in $M/z$ and $x$ is in the guts $P_2\cup\{x\}$ and
$P_k\cup\{y\}$. It now follows by Lemma~\ref{loose1}
that $M/z$ is $k$-coherent if and only if $M\ba x/z$ is.
As $M\ba x/z$ is $k$-coherent, the sublemma follows.
\end{proof}

Relabel $(x,y,z)$ by $(r,s,t)$; relabel
$(O_2,O_3,\ldots,O_{k-1})$ by $(R_2,R_3,\ldots,R_{k-1})$; 
relabel $(P_2,P_3,\ldots,P_{k-1})$ by 
$(S_2,S_3,\ldots,S_{k-1})$; relabel $(Q_2,Q_3,\ldots Q_{k-1})$
by $(T_2,T_3,\ldots,T_{k-1})$; and 
let $Z=O_k\cap P_k\cap Q_k$. With this relabelling it is
a matter of routine checking of the defnition to confirm that
$\{x,y,z\}$ is indeed a gang of three.
\end{proof}

At last we can perform the ritual incantation that completes the
proof of Theorem~\ref{contract}.

\begin{proof}[Proof of Theorem~\ref{contract}]
Assume that $M/x$ is not a $k$-skeleton. 
Then there is an element $y\in E(M/x)$
such that either
\begin{itemize}
\item[(i)] $y$ is cofixed in $M/x$ and $M/x,y$ is coherent, or
\item[(ii)] $y$ is fixed in $M/x$ and $M/x\ba y$ is coherent.
\end{itemize}
By the dual of Lemma~\ref{del-del} (i) does not occur so that (ii)
holds. But then, by Lemma~\ref{con-del}
$x$ is a member of a gang of
three and the theorem follows.
\end{proof}

\section{Removing a Gang of Three}

Proving that $k$-skeletons in a class can be found by an
inductive search is the topic of the next chapter.
Elements in gangs of three are problematic in that we cannot 
remove them individually to keep the property of being a 
$k$-skeleton. In this section we show that we
can remove the whole gang of three. 
Note that, while this is a 3-element move, viewed from a ``bottom
up'' perspective, it is not particularly complicated. It amounts
to extending a $k$-skeleton by adding two elements into the
guts of petals of certain swirl-like flowers of order $k-1$
and performing a
single coextension on the resulting matroid.

\begin{theorem}
\label{gogang}
Let $\{r,s,t\}$ be a gang of three in the $k$-skeleton $M$.
Then $M/r\ba s,t$ is a $k$-skeleton.
\end{theorem}

Given the gang of three $\{r,s,t\}$ there is an associated
canonical partition of the ground set of $M$. 
In what follows we use the 
same labelling for this associated  partition as 
the one given in the definition.

\begin{lemma}
\label{gogang1}
Let $\{r,s,t\}$ be a gang of three in the $k$-skeleton $M$.
\begin{itemize}
\item[(i)] $(S_2\cup\{t\},S_3,\ldots, S_{k-1},R\cup T\cup Z)$
is a swirl-like flower in $M/r$ with $t$ in the guts of 
the petals $S_2\cup \{t\}$
and $R\cup T\cup Z$. Also 
$(T_2\cup \{s\},T_3,\ldots,T_{k-1},R\cup S\cup Z)$
is a swirl-like flower in $M/r$ with $s$ in the guts of 
$T_2\cup\{s\}$ and $R\cup S\cup Z$.
\item[(ii)] $M/r\ba s,t$ is $k$-coherent.
\item[(iii)] If $\{r,p,q\}$ is another gang of three in $M$, then
$\{r,p,q\}=\{r,s,t\}$.
\end{itemize}
\end{lemma}

\begin{proof}
Part (i) is a routine consequence of the definition. 
Part (ii) follows from (i) and Lemma~\ref{clean}. Consider (iii). 
Say that $\{r,s',t'\}$ is a gang of three 
where $\{s',t'\}\neq \{s,t\}$.
Then $\{s',t'\}$ is a 2-element petal in a tight
swirl-like flower
$(\{s',t'\},Q_2,\ldots,Q_{k-1},E(M)-(Q_2\udots Q_{k-1}\cup \{r,s',t'\})$
of $M\ba r$. As $\{s,t\}$ also has this property, we see that
$\{s,t\}\cap \{s',t'\}=\emptyset$. As $\{r,s',t'\}$ is a gang of
three, $\{s',t'\}$ is not a petal of a displayed swirl-like flower
of $M$ and it follows easily that $\{s',t'\}\subseteq Z$.
It is also easily seen that $Q_2\cup Q_3\udots Q_{k-1}\subseteq Z$.
But  $r\notin\cl(E(M)-(Q_2\cup Q_3\udots Q_{k-1}\cup\{r,s',t'\})$
and we obtain the contradiction that $M$ is $k$-fractured.
\end{proof}

We omit the easy proof of the next lemma.

\begin{lemma}
\label{gogang2}
Let $a$ and $b$ be elements of the $3$-connected matroid
$M$. Assume that $a$ is in the guts of a pair of petals of 
a swirl-like flower of $M$ if order at least $4$.
\begin{itemize}
\item[(i)] If $b$ is fixed or cofixed in $M\ba a$, then
$b$ is respectively fixed or cofixed in $M$.
\item[(ii)] If $M\ba a,b$ is $k$-coherent, then 
$M\ba b$ is $k$-coherent.
\item[(iii)] If  $M\ba a/b$ is $k$-coherent and  
$M$ has no triangle containing $a$ and $b$, then $M/b$
is $k$-coherent. 
\end{itemize}
\end{lemma}

\begin{proof}[Proof of Theorem~\ref{gogang}]
By Lemma~\ref{gogang1}(ii) $M/r\ba s,t$ is $k$-coherent. 
Assume that there is an element
$y$ that is fixed in $M/r\ba s,t$ such that 
$M/r\ba s,t,y$ is $k$-coherent.
Then, by Lemma~\ref{gogang2}, $y$ is fixed in 
$M/r$ and $M/r\ba y$ is $k$-coherent.
Now by Lemma~\ref{con-del}, there is a gang of three in $M$ containing 
$r$ and $y$, contradicting Lemma~\ref{gogang1}(iii). 
Thus, if $y$ is fixed in $M/r\ba s,t$, then 
$M/ r\ba s,t,y$ is not $k$-coherent.

Assume that
$y$ is cofixed in $M/r$ and $M/r\ba s,t/y$ is $k$-coherent. 
By Lemma~\ref{gogang2}(i), $y$ is cofixed in $M/r$ and hence
cofixed in $M$. Assume that $M/r$ has no triangle containing
$y$ and either $s$ or $t$.
By Lemma~\ref{gogang2}(iii), $M/r,y$ is $k$-coherent. All up,
$M/r$ and $M/r,y$ are $k$-coherent and $y$ is cofixed in $M/r$,
and we have contradicted the dual of Lemma~\ref{del-del}. 
Otherwise, it is easily seen that $M/r,y$ is $k$-coherent up to
a single parallel pair containing either $s$ or $t$. Moreover,
this parallel pair is coblocked by $r$. We omit the routine argument
that proves that, in this case, $M/y$ is $k$-coherent, 
contradicting the
fact that $y$ is cofixed in $M$.
Thus, if $y$ is cofixed in $M/r\ba s,t$, then $M/r\ba s,t/y$
is not $k$-coherent. Hence $M/r\ba s,t$ is indeed a $k$-skeleton.
\end{proof}

\chapter{A Chain Theorem For Skeletons}
\label{chain-gang}

The goal of this chapter is to prove that skeletons do not
occur sporadically. In particular we prove the following theorem.

\begin{theorem}
\label{not-sporadic}
If $M$ is a $k$-skeleton and $|E(M)|>4$, then $M$ has a $k$-skeleton
minor $M'$ with $|E(M)-E(M')|\leq 4$.
\end{theorem}

In other words all $k$-skeletons in a minor-closed
class can be found by an
inductive search beginning with $U_{2,4}$ 
using 1-, 2-, 3-, or 4-element moves. It
turns out that, in the case that 3- or 4-element moves are required,
very specific structure arises, the full details of which are
made explicit in Theorem~\ref{inductive-search}. 

We also prove in Corollary~\ref{skeleton-free-swirl}
that a swirl-like flower of order $l$ in a $k$-skeleton
has a rank-$l$ free swirl minor. This fact is used later in this paper.
Apart from this, and from some of the easy lemmas,
the results of this chapter are not used later in this paper.
But they are of some independent interest; for example it is possible
that they could be used to get explicit bounds for the 
number of inequivalent representations of $4$-connected
matroids over small prime fields such as $GF(7)$. Moreover, the 
material is so tightly interwoven with other material in this
paper that it would be unwieldy to write it up separately.  
In any case the reader should feel perfectly
relaxed about skipping this chapter.

We will say that a
$k$-skeleton $M$ is 1-{\em reduced} 
\index{$1$-reduced $k$-skeleton}
if $|E(M)|>4$ and  
$M$ has no element $e$ such that
either $M\ba e$ or $M/e$ is a $k$-skeleton. 
We know from Corollary~\ref{stuck} 
that, if $p$ and $q$ are comparable elements
in a 1-reduced $k$-skeleton, then $\{p,q\}$ form a clonal
class of size 2. In this case 
it  often true that  $M\ba p/q$
is a $k$-skeleton, but unfortunately this is not always so.
We say that a $k$-skeleton $M$ is
2-{\em reduced} 
\index{$2$-reduced $k$-skeleton}
if it is 1-reduced and has no clonal pair
$\{p,q\}$ such that $M\ba p/q$ is a $k$-skeleton.
We begin by getting a more explicit structural description of 
1- and $2$-reduced skeletons. To do this it will help to have a 
slightly revised version of $3$-trees of matroids.

\section{Augmented $3$-trees}

Let $M$ be a 3-connected matroid and $(D,C)$ be a
$3$-separation of $M$, where $D$ is a quad.
Then for any partition $(A,B)$ of $D$ into 2-element subsets
the partition
$(A,B,C)$ is a tight swirl-like flower of $M$.
Such a flower has order two as the 3-separation 
$(D,C)$ displays all of the non-sequential 3-separations 
displayed by the flower. 
This is a somewhat degenerate situation, and
in the definition of $3$-trees given in Section~\ref{3tree},
these flowers were specifically excluded from being displayed
unless $(A,B,C)$ refines to a 
larger swirl-like or spike-like 
flower $(A,B,P_1,P_2,\ldots,P_l)$
in $M$, in which case, the refined flower certainly is displayed. 
An example
of this situation arises when $M$ is a swirl, 
in which case
each consecutive pair of petals of the swirl is a quad of $M$.

The partition of $D$ is somewhat less arbitrary 
if we are told that $D$ partitions
into $2$-element clonal classes $A$ and $B$. In this case it
is easily seen that any $3$-separation that crosses $D$ crosses
neither $A$ nor $B$. Choosing to display 
the flower $(A,B,C)$ now seems somewhat less arbitrary. 
We next give a modified version of 3-trees in which such flowers 
are displayed. 
Recall appropriate definitions from Section~\ref{3tree}. The
following definition is sufficiently 
lengthy to challenge even a robust digestive
system, but it is just a modification of the definition of
$3$-tree that guarantees that dags that partition into
2-element clonal classes get displayed by flower vertices.
Let $\pi$ be a partition of $E(M)$. We say that 
the $\pi$-labelled tree $T$ is an {\em augmented $3$-tree}
\index{augmented $3$-tree}
for $M$ if the following hold.
\begin{itemize}
\item[(A1)] For each edge $e$ of $T$, the partition 
$(X,Y)$ of $E$ displayed by
$e$ is $3$-separating, and, if $e$ is incident 
with two bag vertices,
then $(X,Y)$ is a non-sequential $3$-separation.
\item[(A2)] Every non-bag vertex $v$ is 
labelled either $D$ or $A$;
if $v$ is labelled $D$, then there is a cyclic
ordering on the edges incident with $v$.
\item[(A3)] If a vertex $v$ is labelled $A$, 
then either the partition of
$E$ displayed by $v$ is a tight maximal 
anemone of order at least 3 or it has the form $(A,B,C)$
where $A\cup B$ is a quad and $A$ and $B$ are $2$-element
clonal classes.
\item[(A4)] If a vertex $v$ is labelled $D$, 
then either the partition of $E$
displayed by $v$, with the cyclic order 
induced by the cyclic ordering on
the edges incident with $v$, is a tight maximal 
daisy of order at least 3, or it has the form $(A,B,C)$
where $A\cup B$ is a quad and $A$ and $B$ are $2$-element
clonal classes.
\item[(A5)] For every tight maximal flower of $M$ of order 
three, there is an equivalent flower that is displayed by
a vertex of $T$.
\item[(A6)] For every partition $(A,B,C)$ of $E(M)$,
where $A\cup B$ is a quad and $A$ and $B$ are $2$-element
clonal classes, there is a vertex of $T$ that displays a,
possibly trivial,
refinement of 
$(A,B,C)$.
\item[(A7)] If a vertex is incident with two edges $e$ and 
$f$ that display equivalent 3-separations, then the other
ends of $e$ and $f$ are flower vertices, $v$ has degree two,
and $v$ labels a non-empty bag.
\end{itemize}

\begin{lemma}
\label{augmented}
Every $3$-connected matroid has an augmented $3$-tree.
\end{lemma}

\begin{proof}
To facilitate the proof we say that a $\pi$-labelled tree
is a {\em semi-augmented} 3-tree for $M$ if it satisfies all
properties of an augmented 3-tree except that
$(A6)$ need not hold for all partitions $(A,B,C)$ of $M$
where $A\cup B$ is a quad and $A$ and $B$ are 2-element
clonal classes. Note that 3-trees are semi-augmented
3-trees. Thus every 3-connected matroid has a
semi-augmented 3-tree.

Let $T'$ be a semi-augmented 3-tree for $M$. Let $D$ be
a quad of $M$ consisting of 2-element clonal classes
$A$ and $B$ and let $C=E(M)-D$. Assume that no refinement
of $(A,B,C)$ is displayed by $T'$. We claim
that $T'$ can be enhanced to a semi-augmented 3-tree
that displays a refinement of $(A,B,C)$.

Assume that the 3-separation $(D,C)$
is sequential. In this case, $T'$ consists of a single vertex
and it is routinely seen that 
one obtains an augmented 3-tree for $M$
by replacing $T'$ with a tree $T$ containing
a single flower vertex displaying the flower $(A,B,C)$.
Thus we may assume that $(D,C)$ is non-sequential.

Assume that $M$ has a $3$-separation $(X,Y)$ crossing $D$
that is not equivalent to $(D,C)$. Then
we may assume that $A\subseteq X$ and $B\subseteq Y$. 
In this case it is readily checked that
$(X-A,A,B,Y-B)$ is a swirl-like, spike-like 
or V\'amos-like flower of 
$M$ of order $4$. A flower equivalent
to this is displayed by $T'$, 
and it is readily checked that we may modify $T'$
to obtain a semi-augmented $3$-tree in which a refinement of
$(X-A,A,B,Y-B)$ is displayed. Thus the claim holds in this case.

Assume that $M$ has no $3$-separation crossing $D$ that is not
equivalent to $(D,C)$.
Then a $3$-separation equivalent to $(D,C)$ is displayed
by $T$ and we see that $D$ is contained in a bag associated
with a leaf $v$ of $T$. 
Add new vertices $v'$, $v_A$ and $v_B$
to $T'$ to obtain a tree $T$ where $v'$ is incident with
$v$, $v_A$ and $v_B$. In the new tree, $v'$ is a flower vertex,
the new bag at $v$ is the bag in $T'$ at $v$ with $D$ removed,
and the bags at $v_A$ and $v_B$ are $A$ and $B$ respectively.
Label this flower vertex $A$ or $D$ arbitrarily.
The tree $T$ is our required enhanced semi-augmented 3-tree
and the claim holds in this case too.

The lemma now follows by induction.
\end{proof}

Note that an augmented $3$-tree need not be be a $3$-tree,
as it may violate property (ii) of 3-trees.
Nonetheless we extend terminology for $3$-trees to 
augmented $3$-trees in an obvious way.

The new leaf bags that we create in an augmented tree are certainly
petals of a spike-like or swirl-like flower. In fact this is
true for any 2-element leaf bag.

\begin{lemma}
\label{2-element-bag}
Let $T$ be an augmented $3$-tree for a 
$3$-connected matroid $M$ and let $P_1$ be a $2$-element 
leaf bag
of $M$. Then $M$ has a tight 
spike-like, swirl-like or V\'amos-like
flower $(P_1,P_2,\ldots,P_m)$
with at least three petals.
\end{lemma}

\begin{proof}
Certainly $P_1$ is a petal of a tight flower $(P_1,P_2,\ldots,P_m)$
with at least three petals as this is true for any sequential 
3-separating set that is a leaf bag of $T$. Say that 
$(P_1,P_2,\ldots,P_m)$ is not a swirl-like, spike-like 
or V\'amos-like flower. 
Then, up to duality, we may
assume that the flower is a paddle. But then $\sqcap(P_1,P_2)=2$
so that $r(P_1\cup P_2)=r(P_1)+r(P_2)-2=r(P_2)$ and it follows
that $P_1\subseteq \cl(P_2)$, contradicting the fact that the flower
$(P_1,P_2,\ldots,P_m)$ is tight. 
\end{proof}

\section{$1$-reduced $k$-skeletons}

In this section we develop some properties of $k$-skeletons
where we have no single-element deletion or contraction
available to maintain the
property of being a $k$-skeleton. We first need to define
another type of element.

Let $M$ be a $k$-coherent matroid. Then the element $e$ of $M$
is {\em semi-feral}
\index{semi-feral element} 
if either
\begin{itemize}
\item[(i)] $M\ba f$ is 3-connected and $k$-fractured, and
$M/f$ is not 3-connected, or
\item[(ii)] $M/f$ is 3-connected and $k$-fractured, and
$M\ba f$ is not 3-connected.
\end{itemize}

Note that elements of $k$-wild triangles and triads are semi-feral.
The goal of this section is to prove 

\begin{theorem}
\label{1-reduced}
Let $M$ be a $1$-reduced $k$-skeleton were $|E(M)|>4$.
Then the ground set of $M$ consists of feral elements,
semi-feral elements, members of gangs or cogangs of three,
and $2$-element clonal classes. Moreover,
if $T$ is an augmented $3$-tree for $M$, then 
each leaf bag of $T$ is a union of  $2$-element clonal classes 
of $M$.
\end{theorem}

One part of Theorem~\ref{1-reduced} is clear.

\begin{lemma}
\label{ground-set}
Let $M$ be a $1$-reduced $k$-skeleton.
Then the ground set of $M$ consists of 
feral elements, semi-feral elements, members of gangs
or cogangs of three and $2$-element clonal classes.
\end{lemma}

\begin{proof}
Say $e\in E(M)$. If $M$ has an element $f$ that is comparable with
$e$, then, by Corollary~\ref{stuck}, $e$ is in a 
2-element clonal class.

Assume that $e$ is not comparable to any other element. If $e$
is in a triangle or a triad $T$, then $T$ is $k$-wild
by Lemma~\ref{skeleton-triangle} and therefore $e$ is semi-feral.
Assume that $e$ is not in a triangle or a triad. Then either
$M\ba e$ or $M/e$ is 3-connected. If neither $M\ba e$ nor $M/e$
is $k$-coherent, then, $e$ is either feral or semi-feral.
If $M/e$ is $k$-coherent, then, by Theorem~\ref{contract}, $e$ is a 
member of a gang of three and if $M\ba e$ is $k$-coherent,
then, by the dual of Theorem~\ref{contract}, $e$ is a member
of a cogang of three.
\end{proof}

Let $x$ be an element of the $3$-connected matroid $M$.
Recall that $x$ is {\em peripheral} if it belongs to a leaf 
bag of some 3-tree for $M$. We say that $x$ is {\em strongly
peripheral}
\index{strongly peripheral element} 
if it belongs to a leaf bag of some augmented
$3$-tree for $M$.

It is shown in Lemmas~\ref{wheel-wild} and \ref{feral-internal}
that feral elements and members of $k$-wild triangles are
not peripheral. Note that if $x$ is strongly peripheral,
then $x$ is peripheral. We omit the routine proof of 
the next lemma.

\begin{lemma}
\label{no-gang}
If $M$ is a $k$-coherent matroid and $x$ is a member of a
gang or cogang of three in $M$. Then $x$ is not peripheral.
\end{lemma}

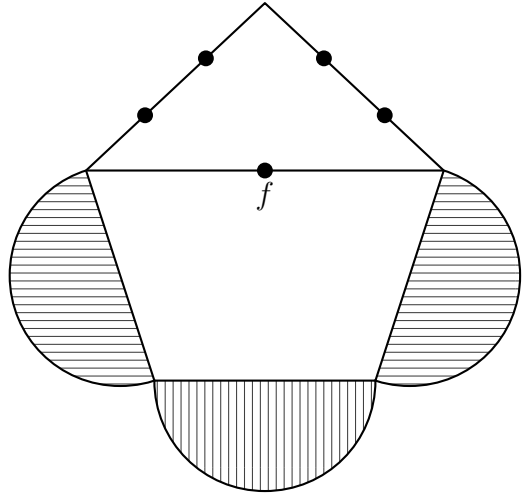
\begin{figure}[h]
\begin{center}
\begin{tikzpicture}[thick,line join=round]
	\coordinate (a) at (378:2.5);
	\coordinate (b) at (450:3);
	\coordinate (c) at (522:2.5);
	\coordinate (d) at (594:2.5);
	\coordinate (e) at (666:2.5);
	\coordinate (cd) at ($(c)!0.5!(d)$);
	\coordinate (de) at ($(d)!0.5!(e)$);
	\coordinate (ea) at ($(e)!0.5!(a)$);
	\node[pattern color=darkgray,draw,circle through=(c),pattern=horizontal lines] (CD) at (cd) {};
	\node[pattern color=darkgray,draw,circle through=(d),pattern=vertical lines] (DE) at (de) {};
	\node[pattern color=darkgray,draw,circle through=(e),pattern=horizontal lines] (EA) at (ea) {};
	\filldraw[fill=white] (a) -- (b) -- (c) -- (d) -- (e) -- cycle;
	\draw (a) -- (c);
	\coordinate[label=below:$f$] (ac) at ($(a)!0.5!(c)$);
	\coordinate (x) at ($(b)!0.33!(a)$);
	\coordinate (y) at ($(b)!0.33!(c)$);
	\coordinate (z) at ($(a)!0.33!(b)$);
	\coordinate (t) at ($(c)!0.33!(b)$);
	\foreach \pt in {x,y,z,t,ac} \fill[black] (\pt) circle (3pt);
\end{tikzpicture}
\label{semi-feral-fig}
\caption{A semi-feral element that is peripheral but not strongly
peripheral}
\end{center}
\end{figure}

On the other hand, it is possible for a 
semi-feral element of a $k$-coherent matroid
to be peripheral.
Figure~\ref{semi-feral-fig} illustrates an example. Indeed,
for an 
arbitrary $k$-coherent matroid it is possible for a semi-feral
element to be strongly peripheral.
Our  task is to  show that this cannot 
happen for 1-reduced $k$-skeletons.

\begin{lemma}
\label{semi-feral-incomparable}
If $f$ is a semi-feral element of the $k$-skeleton 
$M$ then $f$ is not comparable with any other element of $M$.
\end{lemma}

\begin{proof}
Say $f$ is semi-feral and is comparable with the element
$g$ of $M$. If $f$ and $g$ are not clones, then, 
by Corollary~\ref{stuck},
either $M\ba f$ or $M/f$ is $k$-coherent, contradicting the definition
of semi-feral element.
Assume that $f$ and $g$ are clones.
Then either $M\ba f$ or $M/f$ is $k$-coherent by
Corollary~\ref{clone-coherent}, contradicting the defintion
of semi-feral elements.
\end{proof}

Except in special circumstances, semi-feral elements belong to the
guts of some vertical 3-separation.

\begin{lemma}
\label{semi-feral-in-guts}
If $f$ is a semi-feral element of the $k$-skeleton 
$M$ and $M\ba f$ is $3$-connected, then
$\si(M/f)$ is not $3$-connected unless $f$ is in a
costandard $k$-wild triangle.
\end{lemma}

\begin{proof}
By the definition of semi-feral element $M/f$ is not
3-connected. If $\si(M/f)$ is $3$-connected, then
$f$ is in a triangle $T$ that is $k$-wild 
by Lemma~\ref{skeleton-triangle}. If 
$T$ is costandard, then the claim holds. 
Otherwise $T$ is standard.
In this case 
it follows from Lemma~\ref{k-wild-different}(iii)
that $f$ is in the guts of a 
vertical 3-separation so that $\si(M/f)$ is not 
3-connected.
\end{proof}

\begin{lemma}
\label{strong-peripheral}
If $f$ is a semi-feral element of the $1$-reduced $k$-skeleton
$M$, then $f$ is not strongly peripheral.
\end{lemma}

\begin{proof}
Assume that $f$ is semi-feral, where $M\ba f$
is 3-connected. If $f$ is in a $k$-wild triangle, then the
lemma holds by Lemma~\ref{wheel-wild}. Thus we
may assume that $f$ is not in a $k$-wild triangle.
By Lemma~\ref{semi-feral-in-guts}, $f$  is in the guts of a 
vertical $3$-separation $(X\cup\{f\},Y)$ in $M$. Let
$\PP=(P_1,P_2,\ldots,P_n)$ be a $k$-fracture of $M\ba f$.
By Lemma~\ref{semi-feral-incomparable}, Theorem~\ref{canon1}
and the fact that $f$ is not feral, and hence not in a bogan couple,
we deduce that $\PP$ is canonical. If either $X$ or $Y$ is
contained in a petal of $\PP$, then $M$ is not
$k$-coherent, so, up to labels, there is an $i\in \{2,3,\ldots,n-2\}$
such that $X=P_1\cup P_2\udots P_i$. If $i\notin \{2,n-2\}$, 
then it is clear that $f$ is not peripheral. Therefore we
may assume that $i\in\{2,n-2\}$ and, indeed, that $i=2$.

The flower $(P_1\cup P_2\cup\{f\},P_3,\ldots,P_n)$ of $M$
certainly needs to be displayed in $T$. Thus, if 
$f$ is strongly peripheral, then it must be the case that
$P_1\cup P_2\cup\{f\}$ is a leaf bag of $T$. Assume that 
$P_1$ is non-sequential. Then the swirl-like flower
$(P_1,P_2,\{f\}\cup P_3\udots P_n)$ is displayed by $T$,
as $P_1$ is not equivalent to $P_1\cup P_2$. Hence
$P_1$ and $P_2$ are sequential. 
If $|P_1|>2$, then $P_1$ contains either a triangle or triad.
By Lemma~\ref{wheel-wild} this triad or triangle is a clonal
triple and we contradict the assumption that $M$ is 1-reduced.
Hence $|P_1|=|P_2|=2$. Moreover, $(P_1,P_2,\{f\}\cup P_3\udots P_n)$
is a swirl-like or spike-like flower, so that $P_1\cup P_2$ is a 
quad of $M$.

Say $p\in P_1$. As $P_1\cup P_2$ is a quad, 
$M\ba p$ and $M/p$ are $k$-coherent by Lemma~\ref{dag-skeleton}. 
Neither of these matroids is a $k$-skeleton, and $p$ is 
peripheral so is not in 
a gang or cogang of three. Therefore, by Theorem~\ref{contract},
$p$ is comparable with some other element of $M$.  It
follows from Theorem~\ref{reduce-comparable} that
$p$ is in a clonal pair. Evidently this clonal pair must
be $P_1$. Similarly $P_2$ is a clonal pair. Indeed
$P_1$ and $P_2$ are 2-element clonal classes.
By the definition of augmented 3-tree the flower 
$(P_1,P_2,\{f\}\cup P_3\udots P_n)$ is displayed in $M$
leading to the conclusion that $f$ is indeed
not strongly peripheral.
\end{proof}

We can now complete the proof of Theorem~\ref{1-reduced}
which is nothing more than a summary of the facts that
we have garnered so far.

\begin{proof}[Proof of Theorem~\ref{1-reduced}]
By Lemma~\ref{ground-set} the ground set of a 
1-reduced $k$-skeleton consists of feral elements,
semi-feral elements, members of gangs and cogangs of three
and 2-element clonal classes. By Lemmas~\ref{wheel-wild},
\ref{feral-internal}, \ref{no-gang} and \ref{strong-peripheral},
no feral element, semi-feral element, or member of a 
gang or cogang of three is strongly peripheral. Hence 
every leaf bag of an augmented 3-tree for $M$ consists of 
2-element clonal classes of $M$.
\end{proof}

We also note another easy consequence of the results of this section.

\begin{lemma}
\label{cool-peripheral}
If $x$ is a strongly peripheral element of the $k$-skeleton
$M$ and $x$ does not belong to a clonal pair, then
either $M\ba x$ or $M/x$ is a $k$-skeleton.
\end{lemma} 

One consequence of Theorem~\ref{1-reduced} is that
swirl-like flowers in $k$-skeletons have free swirl minors.
We turn attention to this now.
We first prove a lemma that is a consequence of
Tutte's Linking Theorem. In this chapter we only apply the 
swirl-like case. The spike-like case and the copaddle
cases see applications
in Chapter~\ref{unavoidable-minors-of-large}.
Let $\PP=(P_1,P_2,\ldots,P_l)$ be a flower in the connected
matroid $M$. A clonal pair $\{p_i,p'_i\}$ in 
$P_i$ is $\PP$-strong if 
$\kappa(\{p_i,p'_i\},P_1\udots P_{i-1}\cup P_{i+1}\udots P_l)=2$.

\begin{lemma}
\label{free-link}
Let $M$ be a connected matroid with 
a flower $\PP=(P_1,P_2,\ldots,P_l)$ where $l\geq 3$
that is swirl-like (respectively spike-like, a paddle, or
a copaddle).
Assume that, for all $i\in\{1,2,\ldots,l\}$,
the petal $P_i$ contains a $\PP$-strong
clonal pair $\{p_i,q_i\}$. Then
$M$ contains a $\Delta_l$-minor (respectively $\Lambda_l$-,
$U_{2,2l}$- or $U_{2l-2,2l}$-minor). Moreover, in
the the spike-like or swirl-like case, 
the legs of the $\Delta_l$- or $\Lambda_l$-minor are 
$\{p_1,q_1\},\{p_2,q_2\},\ldots,\{p_l,q_l\}$.
\end{lemma}

\begin{proof}
Let $s$ be the 
first element of $\{1,2,\ldots,l\}$ such that $|P_s|>2$.
Let $Z=E(M)-P_s$.
By Tutte's Linking Theorem, there is a minor $M'$ on
$Z\cup\{p_s,q_s\}$
such that 
$M'|Z=M|Z$ and such that
$\lambda_{M'}(\{p_s,q_s\})=2$. One routinely checks
that $M'$ is $3$-connected and that
$\PP'=(P_1,P_2,\ldots,P_{s-1},\{p_s,q_s\},P_{s+1},\ldots,P_l)$
is a flower in $M'$ having the same type that it does
in $M$. It is easily checked that $\{p_i,q_i\}$ is 
a $\PP'$-strong clonal pair for all $i\in\{1,2,\ldots,l\}$.
The lemma now follows routinely.
\end{proof}

\begin{corollary}
\label{skeleton-free-swirl}
Let $M$ be a $k$-skeleton and let $l\geq 4$ be an integer.
If $M$ contains a swirl-like flower of 
order $l$, then $M$ has a $\Delta_l$-minor.
\end{corollary}

\begin{proof}
Let $\PP=(P_1,P_2,\ldots,P_l)$ be a swirl-like flower of order
$l$ in $M$. Assume that $M$ is not 1-reduced. Then there is an
element $x$ in a petal $P_j$ of $\PP$ such that, up to
duality, $M\ba x$ is a $k$-skeleton. By Lemma~\ref{clean},
$P_j-\{x\}$ is not a loose petal of the swirl-like quasi-flower
$(P_1,\ldots,P_{j-1},P_j-\{x\},P_{j+1},\ldots,P_l)$ of 
$M\ba x$. In other words, the above partition is a tight
swirl-like flower of order $l$ in $M\ba x$. It follows
that we lose no generality in assuming that $M$ is 
1-reduced. In this case, by Theorem~\ref{1-reduced}, each petal of
$\PP$ contains a clonal pair and it follows from
Lemma~\ref{free-link} that $M$ has a $\Delta_l$-minor.
\end{proof}

\section{A Miscellany}

We have seen that 1-reduced $k$-skeletons are quite structured.
The next task is to impose further structure on 2-reduced 
$k$-skeletons. Before doing that we develop
some further terminology and prove a few 
lemmas that will be used later in this chapter.
Let $\PP=(P_1,P_2,\ldots,P_m)$ be a flower in a matroid $M$ and $e$
be an element of the petal $P_i$ of $\PP$. Let $N$ be a
$3$-connected matroid in $\{M\ba e, M/e\}$.
Then $e$ {\em opens the flower $\PP$
in $N$}
\index{opens a flower}
or {\em opens the petal $P_i$ of $\PP$}
\index{opens a petal} 
if, for some partition 
$(P_{i_1},P_{I_2},\ldots,P_{i_t})$ of  $P_i-\{e\}$,
the partition 
$(P_1,\ldots,P_{i-1},P_{i_1},\ldots,P_{i_t},P_{i+1},\ldots,P_m)$
is a flower in $N$ whose order is greater than that of 
$\PP$.  

The first lemma is an easy consequence of Lemma~\ref{in-petal}.

\begin{lemma}
\label{either-or}
Let $\PP=(P_1,P_2,\ldots,P_m)$ be a swirl-like flower of the 
$k$-coherent matroid $M$, where $m\geq 3$. 
Assume that $x\in P_1$ and that
$M\ba x$ is $3$-connected and $k$-fractured. Let 
$\QQ=(Q_1,Q_2,\ldots,Q_n)$ be a $k$-fracture of $M\ba x$.
Then either
\begin{itemize}
\item[(i)] The flower $\QQ$ is obtained by opening the petal $P_i$
of $\PP$, or
\item[(ii)] $P_2\cup P_3\udots P_m\subseteq \hQ_j$
for some $j\in\{1,2,\ldots,n\}$.
\end{itemize}
\end{lemma}

\begin{lemma}
\label{either-or-more}
Let $(X,Y)$ be a $3$-separation in a $k$-coherent matroid
$M$ where $|X|\geq 4$.
Assume that $(X,Y)$ does not refine to a swirl-like flower
with at least three petals. If $x\in X$ and $M\ba x$ is $3$-connected
and $k$-fractured by $\PP$, then all but one petal of $\PP$
is contained in $X$.
\end{lemma}

\begin{proof}
Let $\PP=(P_1,P_2,\ldots,P_n)$. We may assume that 
$\PP$ is maximal.
Consider the 3-separation $(X-\{x\},Y)$ in $M\ba x$.
If $X-\{x\}\subseteq \hP_i$ for some $i$, we obtain a 
contradiction to the fact that $M$ is $k$-coherent.
Otherwise, on the assumption that the lemma fails, we may
assume, by possibly moving to a flower equivalent to $\PP$,
that there is an $i\in \{2,3,\ldots,n-2\}$ such that
$X-\{x\}=P_1\cup P_2\udots P_i$. But then
$(X,P_{i+1},P_{i+2},\ldots,P_n)$ is a flower in $M$ that refines
$(X,Y)$.
\end{proof}

Recall that a $3$-connected matroid $M$ is 
{\em uniquely $k$-fractured}
if there is a flower $\QQ$ such that, for every $k$-fracture
$\PP$ of $M$, we have $\PP\less\QQ$. In such a case $\QQ$ may
have order greater than $k$.

\begin{lemma}
\label{open-flower}
Let $M$ be a matroid with a pair of elements $p$ and $q$ such that
$\PP=(\{p,q\},P_2,\ldots,P_m)$ is a maximal tight swirl-like
flower where $m\geq 3$. Say $f\in P_i$ for some 
$i\in\{2,3,\ldots,m\}$, and $f$ has the property that $M\ba f$ 
is not $k$-coherent, but $M\ba p/q\ba f$ is $k$-coherent.
Then the following hold.
\begin{itemize}
\item[(i)] $M\ba f$ is $3$-connected and has a unique $k$-fracture
$\QQ$.
\item[(ii)] $\QQ$ has order $k$.
\item[(iii)] $\QQ$ is obtained by opening the petal $P_i$ of $\PP$.
\end{itemize}
\end{lemma}

\begin{proof}
Assume that $M\ba f$ is not 3-connected. Then, as $\{p,q\}$ is a 
clonal pair and $M\ba p/q\ba f$ is $k$-coherent and therefore 
$3$-connected, we see that $\{p,q\}$ is a series pair of $M\ba f$.
Thus $\{p,q,f\}$ is a triad of $M$ and hence $f$ is a loose element
of $\PP$
in the the coguts of $\{p,q\}$, which is, up to labels, contained
in $P_2$. But then $(P_2-\{f\},P_3\cup P_4\udots P_m)$ is a
2-separation of $M\ba p/q\ba f$, contradicting the fact that
$M\ba p/q\ba f$ is 3-connected. Hence $M\ba f$ is 
$3$-connected.

Say $f\in P_i$. Observe that 
$(\{p,q\},P_2,\ldots,P_{i-1},P_i-\{f\},P_{i+1},\ldots,P_m)$
is a flower in $M\ba f$ which refines to a maximal flower
$\QQ=(\{p,q\},Q_2,\ldots,Q_s)$ in $M\ba f$. Assume that
$\RR=(R_1,R_2,\ldots,R_t)$ is a maximal $k$-fracture of
$M \ba f$ that is not equivalent to $\QQ$. By Lemma~\ref{in-petal},
up to equivalence and labels in $\RR$, 
there is an $i\in\{2,3,\ldots,s\}$
such that
$Q_{i+1}\cup Q_{i+2}\udots Q_s\cup\{p,q\}
\cup Q_2\udots Q_{i-1}\subseteq R_1$.
But, in this case $(R_1-\{p,q\},R_2,\ldots,R_t)$ is clearly a 
$k$-fracture of $M\ba p/q\ba f$. Hence the only possible $k$-fractures
of $M\ba f$ are flowers equivalent to $\QQ$.

Pert (i) of the lemma follows immediately. If $\QQ$ has order greater
than $k$, then $(Q_2,Q_3,\ldots,Q_s)$ is a $k$-fracture of 
$M\ba p/q\ba f$, so (ii) holds. 
Certainly $\QQ$ does not have the property given by 
Lemma~\ref{either-or}(ii), so, by that lemma, 
$\QQ$ is obtained by
opening the petal $P_i$ of $\PP$. Hence part (iii) also holds.
\end{proof}

While we may not be able to remove elements from a 1-reduced
$k$-skeleton to keep a $k$-skeleton, 
we can always remove peripheral elements to 
keep $k$-coherence.

\begin{lemma}
\label{keep-coherent-pq}
Let $M$ be a $1$-reduced $k$-skeleton and let
$\{p,q\}$ be a clonal pair in a leaf bag $B$ 
of an augmented $3$-tree for $M$. Then $M\ba p$, $M/p$,
and $M\ba p/q$ are $k$-coherent.
\end{lemma}

\begin{proof}
By Theorem~\ref{1-reduced} $B$ consists of 
$2$-element clonal classes of $M$. By 
Corollary~\ref{clone-coherent} we may assume, up to
duality, that $M\ba p$ is $k$-coherent. Assume that
$M\ba p/q$ is not $k$-coherent. As $M$ is $1$-reduced,
there is an element $f$ such that either
\begin{itemize}
\item[(a)] $f$ is fixed in $M\ba p$ and 
the matroid $M\ba p,f$ is $k$-coherent, or
\item[(b)] $f$ is cofixed in $M\ba p$ and 
the matroid $M\ba p/f$ is $k$-coherent.
\end{itemize}
As $q$ is not fixed in 
$M\ba p$ and $M\ba p/q$ is not $k$-coherent, 
we see that in either case $f\neq q$. By 
Theorem~\ref{1-reduced},
$f$ is not a clone of $q$. Now, if $f$ is comparable with $q$,
then by Corollary~\ref{stuck} 
$M$ is not $1$-reduced, so $f$ is not comparable with $q$. 
In case (a), $f$ is clearly
fixed in $M$. In case (b) by Corollary~\ref{im-free}, 
$f$ is cofixed in $M$.
By Theorem~\ref{1-reduced}, $B$ consists of 2-element
clonal classes of $M$. Thus,
in either case (a) or (b), we see that $f\notin B$. 

Let $N=M\ba p$ if case (a) holds and let $N=M/p$ if case (b) holds.
Let $\PP=(P_1,P_2,\ldots,P_n)$ be a $k$-fracture in $N$. 
Assume that $B\neq \{p,q\}$ so that $|B|\geq 4$. Assume that, for
some $i\in\{2,3,\ldots,n-2\}$, we have $B=P_1\cup P_2\udots P_i$.
Then $(P_1,P_2,\ldots,P_i,P_{i+1}\udots P_n\cup \{f\})$
is a swirl-like flower in $M$, and this is easily seen to 
contradict the assumption that $B$ is a peripheral bag of 
an augmented 3-tree for $M$. Thus we may assume that
$B\subseteq P_1$. But now $P_1-\{p\}$ contains a clonal pair
and cannot be a set of loose elements of 
the flower $(P_1-\{p\},P_2,\ldots,P_n)$ of $N\ba p$,
contradicting the assumption that $N\ba p$ is $k$-coherent.

Hence $B=\{p,q\}$. Now $B$ is a $2$-element bag of an augmented
$3$-tree for $M$. In this case, by Lemma~\ref{2-element-bag},
there exists a swirl-like, 
spike-like or V\'amos-like 
flower $\QQ=(\{p,q\},Q_2,\ldots,Q_m)$ in 
$M$. For some $i\in\{2,3,\ldots,m\}$, we have $f\in Q_i$.
By Lemma~\ref{open-flower}, there is a partition
$(Q_{i_1},\ldots,Q_{i_t})$ of $Q_i$ such that 
$$(\{p,q\},Q_2,\ldots,Q_{i-1},Q_{i_1},
\ldots,Q_{i_t},Q_{i+1},\ldots,Q_m)$$ 
is a $k$-fracture of $N$. But $N\ba p/q$ is $k$-coherent so 
$$(Q_2,\ldots,Q_{i-1},Q_{i_1},
\ldots,Q_{i_t},Q_{i+1},\ldots,Q_m)$$ 
is a swirl-like flower of order $k-1$ in $N\ba p/q$. However 
$M\ba p/q$ is not $k$-coherent, so by Lemma~\ref{in-closure},
$p\in\clstar(\hQ)$ for some petal $Q$ of the above flower.
This is easily seen to contradict the maximality of $\QQ$ in $M$.

Thus $N\ba p/q$ is $k$-coherent and it is routinely verified that so
too is $N/q$.
\end{proof}

The last two lemmas of this section 
are perhaps oddly placed, but they are close to their
first application.

\begin{lemma}
\label{2-elt-feral}
Let $M$ be a $k$-coherent matroid with an element $z$
such that $M\ba z$ is $3$-connected with a 
$k$-fracture $(\{p_1,p'_1\},P_2,\ldots,P_l)$.
If $\{p_1,p'_1\}$ is fully closed in $M\ba z$,
then $M\ba p_1$ is $k$-coherent.
\end{lemma}

\begin{proof}
By Lemma~\ref{2-3-4-remove}, $M\ba z,p_1$ is 3-connected.
Therefore $M\ba p_1$ is 3-connected.

Assume that $M\ba z$ is not uniquely $k$-fractured.
Let $(T_1,T_2,\ldots,T_m)$ be another $k$-fracture of 
$M\ba z$. By Lemma~\ref{in-petal}, there is an
$i\in\{1,2,\ldots,l\}$ and a $j\in\{1,2,\ldots,m\}$
such that 
$P_1\cup P_2\udots P_{i-1}\cup P_{i+1}\udots P_l\subseteq \hat{T}_j$,
and 
$T_1\cup T_2\udots T_{j-1}\cup T_{j+1}\udots T_m\subseteq \hP_i$.
Evidently $i\neq 1$, so that $p_1\in T_j$.
Certainly $M\ba z,p_1$ is 3-connected and
$(T_1,\ldots,T_{j-1},T_j-\{p_1\},T_{j+1},\ldots,T_m)$
is a $k$-fracture of $M\ba z,p_1$. If this $k$-fracture 
induces a $k$-fracture $(T'_1,T'_2,\ldots,T'_k)$
in $M\ba p_1$, then, as $T_j-\{p_1\}\subseteq T'_\mu$,
for some $\mu$ in $\{1,2,\ldots,k\}$, we obtain
the contradiction that
$(T'_1,T'_2,\ldots,T'_{\mu-1},T'_\mu
\cup\{p_1\},T'_{\mu+1},\ldots,T'_k)$
is a $k$-fracture of $M$.

From the above we deduce that any $k$-fracture of 
$M\ba p_1$ is induced by the quasi-flower
$(\{p'_1\},P_2,\ldots,P_l)$ of $M\ba z,p_1$ which 
we may assume is maximal. Thus, on the assumption
that the lemma fails, there is a quasi-flower
$\OO=(O_1,O_2,\ldots,O_k)$, displayed by the above 
quasi-flower such that 
$(O_1,O_2,\ldots,O_{i-1},O_i\cup\{z\},O_{i+1},\ldots,O_k)$
is a $k$-fracture of $M\ba p_1$ for some $i\in\{1,2,\dots,k\}$.
As $(\{p'_1\},P_2,\ldots,P_l)$ is maximal,
and $p'_1\in\cl^*_{M\ba p_1,z}(P_2)$,
we see that $P_2\cup\{p_1'\}\subseteq \hO_j$ for some
$j\in\{1,2,\ldots,k\}$, and also that
$P_l\cup\{p'_1\}\subseteq \hO_s$ for some $s\in\{1,2,\ldots,k\}$.
Note that $p_1\in \cl(\{p'_1\}\cup P_2)$
and $p_1\in\cl(\{p'_1\}\cup P_l)$.
The only cases that are not almost immediately seen to lead to the
contradiction that $M$ is $k$-fractured are when, up to labels,
we have, either (a) $j=s=i-1$ or (b) $j=i-1$ and $s=i+1$.
By Lemma~\ref{modular}, case (b) leads to the contradiction
that $(O_1,O_2,\ldots,O_{i-1},O_i\cup\{z,p\},O_{i+1},\ldots,O_k)$
is a $k$-fracture of $M$. On the other hand, by considering the
structure of swirl-like flowers, we may assume in case (a),
that $(P_2\cup\{p'_1\})\cap \hO_i=\emptyset$,
so that, for some flower $(O'_1,O'_2,\ldots,O'_k)$,
equivalent to $\OO$, the partition
$(O'_1,O'_2,\ldots,O'_{i-1}\cup\{p_1\},
O'_1\cup\{z\},O'_{i+1},\ldots,O'_k)$ 
is a $k$-fracture of $M$. The lemma follows from this final
contradiction.
\end{proof} 

We omit the proof of the next lemma which amounts to little
more than observing the properties of a feral display.

\begin{lemma}
\label{display-way}
Let $f$ be a feral element of the $k$-coherent matroid
$M$.
\begin{itemize}
\item[(i)] If $f$ blocks two petals of a $k$-fracture
of $M\ba f$, then there is a feral display for $f$ in $M^*$.
\item[(ii)] If $f$ coblocks two petals of a $k$-fracture
of $M/f$, then there is a feral display for $f$ in $M$.
\end{itemize}
\end{lemma}

\section{$2$-reduced Skeletons}

For $2$-reduced skeletons we can strengthen the outcome
of Theorem~\ref{1-reduced} somewhat. 
We say that a clonal pair $\{p,q\}$ of a 3-connected matroid $M$
is {\em strongly peripheral}
\index{strongly peripheral clonal pair}
if it is contained in a leaf bag of an augmented 3-tree for $M$.
Let $\{p,q\}$ be a strongly peripheral
clonal pair of the $2$-reduced $k$-skeleton $M$.
For the remainder of this chapter, if $\{p,q\}$
is a clonal pair of $M$, then the matroid $M\ba p/q$
will be denoted by $M_{pq}$. By Lemma~\ref{keep-coherent-pq},
$M_{pq}$ is $k$-coherent. As $M_{pq}$ is not a $k$-skeleton,
there is an element $f$ such that either 
$f$ is fixed in $M_{pq}$ and $M_{pq}\ba f$ is $k$-coherent, in which
case we say that $f$ is {\em $pq$-annoying for deletion},
\index{$pq$-annoying for deletion} 
or $f$ is cofixed in $M_{pq}$ and 
$M_{pq}/f$ is $k$-coherent, in which case we say that
$f$ is {\em $pq$-annoying for contraction}.
\index{$pq$-annoying for contraction} 
If either one of the
cases holds we say that $f$ is {\em $pq$-annoying}.
\index{$pq$-annoying}

\begin{lemma}
\label{annoying-fracture}
Let $M$ be a $2$-reduced $k$-skeleton and let $\{p,q\}$
be a strongly peripheral clonal pair of $M$. 
\begin{itemize}
\item[(i)] If $f$ is $pq$-annoying for deletion, then $f$
is fixed in $M$ and $M\ba f$ is $3$-connected and $k$-fractured.
\item[(ii)] If $f$ is $pq$-annoying for contraction, then $f$
is cofixed in $M$ and $M/f$ is $3$-connected and $k$-fractured.
\end{itemize}
\end{lemma}

\begin{proof}
Assume that $M\ba f$ is not 3-connected. Then, as $M_{pq}\ba f$
is $3$-connected, $\{p,q,f\}$ is a triad of $M$ and hence a
clonal triple so that $M$ is not $2$-reduced. Hence $M\ba f$ is
$3$-connected. Certainly $f$ is not comparable with 
either $p$ or $q$ so, by Corollary~\ref{born-free}, $f$
is fixed in $M$. It now follows from the definition of 
$k$-skeleton that $M\ba f$ is $k$-fractured. Thus (i) holds. Part
(ii) is the dual of (i).
\end{proof}

\begin{lemma}
\label{2-reduce}
Let $M$ be a $2$-reduced $k$-skeleton and let $B$ be a 
leaf bag of an augmented $3$-tree for $M$. 
Then the following hold.
\begin{itemize}
\item[(i)] $B$ consists of a single clonal pair $\{p,q\}$.
\item[(ii)] $M$ has a tight, maximal,
swirl-like flower $\PP=(B,P_2,\ldots,P_m)$ for some $m\geq 3$..
\item[(iii)] If $f$ is $pq$-annoying for deletion, then $f$
opens the flower $\PP$ in $M\ba f$.
\end{itemize}
\end{lemma}

\begin{proof}
By Theorem~\ref{1-reduced}, $B$ consists of clonal
classes of size $2$. Let $\{p,q\}$ be a clonal pair contained in 
$B$. As $M$ is 2-reduced, there is an element $f$ 
of $M_{pq}$ that is $pq$-annoying.
We lose no generality in assuming that $f$ is $pq$-annoying
for deletion. By Lemma~\ref{annoying-fracture},
$M\ba f$ is 3-connected and $k$-fractured.
Let  $\QQ=(Q_1,Q_2,\ldots,Q_n)$
be a maximal $k$-fracture of $M\ba f$. 

Assume that $B\neq \{p,q\}$ so that $|B|\geq 4$. Assume that, for
some $i\in\{2,3,\ldots,n-2\}$ we have $B=Q_1\cup Q_2\udots Q_i$.
Then $(Q_1,Q_2,\ldots,Q_i,Q_{i+1}\udots Q_n\cup \{f\})$
is a swirl-like flower in $M$, and this is easily seen to 
contradict the assumption that $B$ is a peripheral bag of 
an augmented 3-tree for $M$. Thus we may assume that
$B\subseteq Q_1$. But now $Q_1-\{p,q\}$ contains a clonal pair
of the 3-connected matroid $M_{pq}\ba f$
and cannot be a set of loose elements of 
the flower $(Q_1-\{p,q\},Q_2,\ldots,Q_n)$ of $M_{pq}\ba f$
contradicting the fact that $M_{pq}\ba f$ is $k$-coherent.
Therefore $B=\{p,q\}$ so that (i) holds.

By Lemma~\ref{2-element-bag} $M$ has a tight  
flower $\PP=(\{p,q\},P_2,\ldots,P_m)$ where $m\geq 3$. Say
$f\in P_i$. By Lemma~\ref{open-flower}, any $k$-fracture of
$M\ba f$ is obtained by opening $\PP$. This implies that
$\PP$ is swirl-like. (If $m=3$, then $\PP$ will be ambiguous.)
Parts (ii) and (iii) of the lemma follow from these observations.
\end{proof}



\section{Removing a Bogan Couple}

Our goal is to show that in a $2$-reduced $k$-skeleton
we can always find a $3$- or $4$-element move
that preserves the property of being a $k$-skeleton.
In this section we show that bogan couples lead
to a good outcome in that, if a $pq$-annoying element 
belongs to a bogan couple, then we have a 4-element move.

\begin{lemma}
\label{bogan-reduce}
Let $\{p,q\}$ be a strongly-peripheral clonal pair of the
$2$-reduced $k$-skeleton $M$ and let $a$ be an element of $M_{pq}$
that is $pq$-annoying. If $a$ belongs to a bogan or cobogan
couple $\{a,b\}$, then $M_{pq}\ba a/b$ is a $k$-skeleton.
\end{lemma}

\begin{proof}
Up to duality we may assume that $\{a,b\}$ is a bogan couple.
Associated to the bogan couple $\{a,b\}$ is a partition
$(R,S,T,\{a,b\})$ together with partitions 
$(R_1,R_2,\ldots,R_{k-2})$, $(S_1,S_2,\ldots,S_r)$ and
$(T_1,T_2,\ldots,T_{k-2})$ of $R$, $S$ and $T$ respectively that
form a bogan display for $\{a,b\}$. Note that we have labelled the
display just as in the definition of bogan couple. The case where
$a$ is $pq$-annoying for deletion is not quite symmetric to the
case where $a$ is $pq$-annoying for contraction---or at least
it requires argument to show that it is---but essentially the same
proof works in either case. We prove the lemma in the case that 
$a$ is $pq$-annoying for deletion. In this case
$(R_1\cup\{b\},R_2,\ldots,R_{k-2},S_1,\ldots,S_r,T)$
is a maximal swirl-like flower in $M\ba a$ and $b$ is in the
coguts of the petals $R_1\cup\{b\}$ and $T$. Note that, if $r>1$,
then the above partition with $\{p,q\}$ removed induces a 
$k$-fracture of $M_{pq}\ba a$. Hence $r=1$, so that $S=S_1$
and the flower is $\FF=(R_1\cup\{b\},R_2,\ldots, R_{k-2},S,T)$.

\begin{sublemma}
\label{bogan-reduce1}
$\{p,q\}= R_i$ for some $i\in\{1,2,\ldots,k-2\}$.
\end{sublemma}

\subproof
As $\{p,q\}$ is a clonal pair, this set is contained in a petal of
$\FF$. Moreover this petal becomes a, possibly empty, set of loose
elements in the induced flower in $M_{pq}\ba a$. Given this,
it is clear that 
$\{p,q\}$ is not contained in $T$. 
As $\{p,q\}$ is strongly peripheral and $M$ is 2-reduced,
by Lemma~\ref{2-reduce}, there is a flower $\OO$ displayed
in an augmented 3-tree for $M$ in which $\{p,q\}$ is a petal.
Say $\{p,q\}\subseteq S$.
The flower $\OO$ conforms with the maximal flowers displayed by the
bogan display for $\{a,b\}$ and it follows that $E(M)-S$
is contained in a petal of $\OO$. But then $E(M)-S$
contains a strongly peripheral set other than $\{p,q\}$ and,
again by Lemma~\ref{2-reduce}, we see that $E(M)-S$ contains a 
clonal pair of $M$ other than $\{p,q\}$ and is 
hence not a set of loose elements of a swirl-like flower
in $M_{pq}\ba a$. Thus $\{p,q\}\subseteq R_i$ for some
$i\in\{1,2,\ldots,k-2\}$. If $R_i$ is not a petal of $\FF$,
then the previous argument applies with $S$ replaced by $R_i$.
Thus $R_i$ is a petal of $\OO$, that is $\{p,q\}=R_i$.
\end{proof}

\begin{sublemma}
\label{bogan-reduce2}
$M_{pq}\ba a/b$ is $3$-connected.
\end{sublemma}

\subproof
Certainly $M_{pq}\ba a$ is $3$-connected. Let
$(R'_1,R'_2,\ldots,R'_{k-3})=
(R_1,\ldots,R_{i-1},R_{i+1},\ldots,R_{k-2})$.
Note that $b$ is in the coguts of the petals $R_1'$ and $T$ of
the swirl-like flower $(R'_1\cup\{b\},R'_2,\ldots,R'_{k-3},S,T)$
of $M_{pq}\ba a/b$. Thus $M_{pq}\ba a/b$ is 3-connected unless
$b$ is in a triangle $\{b,r,t\}$ where $r\in R'_1$ and $t\in T$.
If $R_1\neq \{p,q\}$, then $R'_1=R_1$ and we deduce that $b$
is in a triangle of $M$ contradicting the fact that $M/b$
is $3$-connected. If $R_1=\{p,q\}$, then $r\in\cl_M(\{p,q\})$
and $\{p,q,r\}$ is a triangle of $M$ contradicting the fact that
$M$ is 2-reduced.
\end{proof}

Assume that $M_{pq}\ba a/b$ is not a $k$-skeleton. Then there is
an element $h$ such that either (i) $h$ is fixed in $M_{pq}\ba a/b$
and $M_{pq}\ba a/b\ba h$ is $k$-coherent, of (ii) $h$ is cofixed
in $M_{pq}\ba a/b$ and $M_{pq}\ba a/b,h$ is $k$-coherent.
In what follows we assume that (i) holds noting that the argument
in the case that (ii) holds is essentially identical.

\begin{sublemma}
\label{bogan-reduce3}
$h$ is fixed in $M$.
\end{sublemma}

\subproof
Assume that $h$ is not fixed in $M$. Let $M'$ be a matroid
obtained by independently cloning $h$ by $h'$. For a set
$Z\subseteq E(M)$, let $Z'$ denote $Z\cup\{h'\}$ if $h\in Z$ and
otherwise $Z'=Z$. It is easily checked that the bogan display for
$\{a,b\}$ in $M$ extends to a bogan display in $M'$ where each
member $Z$ of the partitions of the display in $M$ is replaced by
$Z'$. But then one readily checks that $\{h,h'\}$ is independent
in $M'\ba p/q\ba a/b$ so that $h$ is not fixed in $M_{pq}\ba a/b$.
Thus $h$ is indeed fixed in $M$.
\end{proof}

We omit the routine verification of the next claim.

\begin{sublemma}
\label{bogan-reduce4}
$M\ba h$ is $3$-connected.
\end{sublemma}

As $h$ is fixed in $M$ and $M\ba h$ is 3-connected, $M\ba h$ is
$k$-fractured. We will obtain a contradiction by showing that there
is no sensible place for $h$. Assume first 
that $h\in R_j$ for some $j\in\{1,\ldots,k-2\}$.
Assume that a $k$-fracture for $M\ba h$ is obtained by opening the
petal $R_j$ of the flower $(R_1,R_2,\ldots,R_{k-2},S\cup T\cup\{a,b\})$
of $M$. Say 
$(R_1,\ldots,R_{j-1},R_{j1},\ldots,R_{jl},R_{j+1},\ldots R_{k-2},
S\cup T\cup\{a,b\})$ is such a $k$-fracture. In this case
$(R_1,\ldots,R_{j-1},R_{j1},\ldots,R_{jl},R_{j+1},\ldots R_{k-2},
S,T\cup\{b\})$ is a $k$-fracture of $M\ba a$ whose order is greater 
than $k$ and it follows that $M_{pq}\ba a/b\ba h$ is 
$k$-fractured. Otherwise $E(M)-R_j$ is contained in a petal 
of a $k$-fracture for $M\ba h$ and again it follows that 
$M_{pq}\ba a/b\ba h$ is $k$-fractured.

Thus $h\notin R$ and, similarly, $h\notin T$. Hence $h\in S$.
We next show that $|S|\geq 3$. Assume for a contradiction that
$|S|=2$, say $S=\{h,h'\}$. Then
$(R_1\cup\{b\},R_2,\ldots,R_{k-2},\{h,h'\},T)$ is a $k$-fracture
of $M\ba a$. One readily checks that $\{h,h'\}$ is fully
closed in $M\ba a$. Lemma~\ref{2-elt-feral} now applies and
it follows from that lemma that $M\ba h$ is $k$-coherent.
Thus $|S|\geq 3$ as desired. As $M\ba h$ is 
3-connected and $S$ is 3-separating in $M$ we have
$h\in\cl(S-\{h\})$.

If $E(M)-S$ is contained in a petal of a $k$-fracture for
$M\ba h$, then again it follows easily that $M_{pq}\ba a/b\ba h$ is
$k$-fractured. Otherwise there is a partition $(Z_1,Z_2)$ of
$E(M)-S$ such that $\lambda_{M\ba h}(Z_1)=\lambda_{M\ba h}(Z_2)=2$.
As $h\in\cl(S-\{h\})$, neither $Z_1$ nor $Z_2$ is blocked 
by $h$ so that 
$\lambda_M(Z_1)=\lambda_M(Z_2)=2$, that is $(S,Z_1,Z_2)$ is a flower
in $M$. Such a partition contradicts the maximality of the flowers
displayed in a bogan display. The lemma follows from this
final contradiction.
\end{proof}

\section{Life in  $2$-reduced Skeletons}
\label{life}

In this section we obtain more information about the structure
associated with a strongly peripheral clonal pair in
a 2-reduced $k$-skeleton.
Throughout this section $\{p,q\}$ will denote a peripheral
clonal pair of the 2-reduced $k$-skeleton $M$.
As $M$ is $2$-reduced, there is an element $f$
of $M$ that is $pq$-annoying. Up to duality we may assume that
$f$ is $pq$-annoying for deletion.
By Lemma~\ref{2-reduce}, 
there is a swirl-like flower
$\QQ=(Q_1,Q_2,\ldots,Q_k)$ in $M\ba f$, and an $s\in\{2,3,\ldots,k-2\}$
such that 
\begin{itemize}
\item[(i)] $Q_t=\{p,q\}$ for some $t\in\{s+1,\ldots,k\}$, and
\item[(ii)] $(Q_1\udots Q_s\cup\{f\},Q_{s+1},\ldots,Q_k)$
is a maximal swirl-like flower of $M$.
\end{itemize}
Let $Q'_1=Q_1\cup Q_2\udots Q_s\cup\{f\}$, and
let $\QQ'=(Q'_1,Q_{s+1},\ldots,Q_k)$.
As with $M$ and $\{p,q\}$, the element $f$ and the flowers $\QQ$
and $\QQ'$
will be fixed throughout this section.

We say that the element
$g$ is $f$-{\em bad for deletion}
\index{$f$-bad for deletion} 
if $f$ is fixed in $M_{pq}\ba f$
and $M_{pq}\ba f,g$ is $k$-coherent, and $g$ is 
$f$-{\em bad for contraction}
\index{$f$-bad for contraction} 
if $f$ is cofixed in $M_{pq}\ba f$
and $M_{pq}\ba f/g$ is $k$-coherent. The element
$g$ is $f$-{\em bad}
\index{$f$-bad} 
if it is $f$-bad for either
deletion or contraction. The next lemma follows from the 
definition of $k$-skeleton.

\begin{lemma}
\label{bad-boy1}
If $M_{pq}\ba f$ is not a $k$-skeleton, then $M_{pq}\ba f$
has an element that is $f$-bad.
\end{lemma} 

\begin{lemma}
\label{bad-boy2}
If $g$ is $f$-bad for deletion, then $g$ is fixed in $M$,
and if $g$ is $f$-bad for contraction, then 
$g$ is cofixed in $M$.
\end{lemma}

\begin{proof}
Assume that $g$ is $f$-bad for deletion.
Then $g$ is fixed in $M_{pq}\ba f$ so that $g$ is certainly
fixed in $M/p$. But $p$ is not comparable with $g$ so 
$g$ is fixed in $M$. Assume that $g$ is $f$-bad for contraction.
Then $g$ is cofixed in $M_{pq}\ba f$ so that 
$g$ is cofixed in $M\ba p,f$. If $g$ is not
cofixed in $M\ba f$, then, by Corollary~\ref{born-free}, 
$p\more g$ in $M\ba f$ implying that $\{p,q,g\}$ is a
triangle of $M$. Hence $g$ is  cofixed in $M\ba f$.
But as $g$ is not comparable with any element of $M$, we see,
by Corollary~\ref{im-free}, that $g$ is cofixed in $M$. 
\end{proof}

Assume that $g$ is $f$-bad. If $g$ is $f$-bad for deletion, then
the symbol $*$ will denote deletion and 
if $g$ is $f$-bad for contraction, then $*$  will denote 
contraction. 

\begin{lemma}
\label{bad-boy3}
If $g$ is $f$-bad, then $M*g$ and $M\ba f*g$ are $3$-connected.
\end{lemma}

\begin{proof}
Assume that $M\ba f*g$ is not 3-connected.
Consider the 3-connected matroid $M\ba f$. There is a
$3$-separation $(X\cup\{g\},Y)$ of $M\ba f$
for which $g$ is either in the guts or the coguts according
as to whether $*$ is contraction or deletion. 
For some $i\in\{1,2,\ldots,k\}$, we have $g\in Q_i$.
We may assume that $Q_i$ is fully closed. 
As $\QQ$ is maximal, $(X\cup\{g\},Y)$ conforms with
$\QQ$. Thus we may assume that $X\cup\{g\}\subseteq Q_i$.
But then $(X,Y-\{p,q\})$ is clearly a 2-separation of
$M_{pq}\ba f*g$ contradicting the fact that this matroid
is 3-connected. It follows from this contradiction that 
$M\ba f*g$ is 3-connected.

Consider $M*g$. Assume that $M*g$ is not 3-connected. As
$M\ba f*g$ is 3-connected, $M*g$ has a single parallel pair
containing $f$. Thus $*$ is contraction and $M$ has a triangle
$T$ containing $f$ and $g$. As $g$ is cofixed in $M$, the 
triangle $T$ is $k$-wild. Now $\si(M/g)$ is a 3-connected matroid,
so, by Lemma~\ref{k-wild-different}, $T$ is costandard.
But it follows easily from the definition of $k$-wild display that,
in this case, $M\ba f/g$ has a swirl-like flower whose order is 
greater than $k$. But this implies that $M_{pq}\ba f/g$
is not $k$-coherent. Therefore $M*g$ is also 3-connected.
\end{proof}

The next lemma follows immediately from Lemmas~\ref{bad-boy2}
and \ref{bad-boy3}.

\begin{lemma}
\label{bad-boy4}
If $g$ is $f$-bad, then $M*g$ is $3$-connected and $k$-fractured.
\end{lemma}

Next we gain a little more information about the location of
$f$-bad elements.

\begin{lemma}
\label{bad-boy5}
If $g$ is $f$-bad, then $g\in Q_1\cup Q_2\udots Q_s$.
\end{lemma}

\begin{proof}
Recall that $\QQ'$ denotes the flower
$(Q'_1,Q_{s+1},\ldots,Q_k)$ of $M$, where 
$Q'_1=Q_1\udots Q_s\cup\{f\}$.
Assume that $g\in Q_j$, where $j\in\{s+1,s+2,\ldots,k\}$.  Assume that 
a $k$-fracture of $M*g$ is obtained 
by opening the petal $Q_j$ of $\QQ'$. Then the flower
$(Q_1,Q_2,\ldots,Q_{j-1},Q_j-\{g\},Q_{j+1},\ldots,Q_k)$ of $M\ba f*g$
expands to a swirl-like flower whose order is strictly greater than
$k$. It follows from this that $M_{pq}\ba f*g$ is $k$-fractured.
On the other hand, by Lemma~\ref{either-or}, $E(M)-Q_j$
is contained in a petal $P_1$ of a $k$-fracture 
$(P_1,P_2,\ldots,P_k)$ of  $M*g$. In this case
it is clear that $P_1-\{p,q,f\}$ is not a set of loose elements
of the swirl-like flower $(P_1-\{p,q,f\},P_2,\ldots,P_k)$
of $M_{pq}\ba f*g$ and again we see that this matroid 
is $k$-fractured. Thus $g\in Q'_1$, that is, 
$g\in Q_1\cup Q_2\udots Q_s$.
\end{proof}

One possibility that leads to a good outcome is when 
$f$ is semi-feral.

\begin{lemma}
\label{bad-boy6}
If $f$ is semi-feral, then $M_{pq}\ba f$ is a $k$-skeleton.
\end{lemma}

\begin{proof}
Assume that $f$ is semi-feral and assume that $M_{pq}\ba f$
is not a $k$-skeleton. By Lemma~\ref{bad-boy1},
there is an element $g$ that is $f$-bad.

We first prove that 
$f$ is not in a costandard $k$-wild triangle. For ease of
notation we relabel $f$ to $a$ and assume that $a$ belongs
to the costandard $k$-wild triangle $\{a,b,c\}$. 
With notation as in the definition of costandard
$k$-wild triangle, we see that 
$M\ba a$ has a unique $k$-fracture 
$(A_1,A_2,\ldots,A_{k-2},B\cup\{b\},C\cup\{c\})$.
As $\{a,b,c\}$ is costandard, the 
elements $b$ and $c$ are in the coguts of
$B\cup\{b\}$ and $C\cup\{c\}$ respectively. As $\QQ'$
and $(A_1,A_2,\ldots,A_{k-2},B\cup C\cup\{a,b,c\})$ are 
flowers in $M$ and both have petals that open to 
$k$-fractures in the uniquely
$k$-fractured matroid $M\ba a$, we deduce that $s=2$, and that 
$(A_1,A_2,\ldots,A_{k-2},B\cup\{b\},C\cup\{c\})=
(Q_3,Q_4,\ldots,Q_k,Q_1,Q_2)$. Hence $\{p,q\}=A_i$
for some $i\in\{1,2,\ldots,k-2\}$.  By Lemma~\ref{bad-boy5},
$g\in B\cup C\cup\{b,c\}$. Up to symmetry, the $f$-bad 
(that is $a$-bad) element
$g$ belongs to $B\cup\{b\}$.
Say $g=b$. In this case $M/g$ is not 3-connected, so,
by Lemma~\ref{bad-boy4}, $g$ must be $f$-bad for deletion.
But $g$ is in the coguts of $B$ in $M\ba a$, so $M\ba a,g$ is not
$3$-connected, contradicting Lemma~\ref{bad-boy3}. Thus $g\neq b$.
Recall that there is a partition $(B_1,B_2,\ldots,B_k)$ 
of $B$ such that
$(B_1,B_2,\ldots,B_k,A\cup C\cup\{a,b,c\})$ is a swirl-like flower
of $M$. We have $g\in B_i$ for some 
$i\in\{1,2,\ldots,k-2\}$. By Lemma~\ref{either-or}, a $k$-fracture
for $M*g$ is either obtained by opening the petal $B_i$
of the flower $(B_1,B_2\ldots,B_k,A\cup C\cup\{a,b,c\})$ of $M$,
or has the property that $E(M)-B_i$ is contained in a petal of
this $k$-fracture. In either case it is easily deduced that
$M_{pq}\ba a*f=M\ba g\ba p/q\ba a$ is either not $3$-connected or
is $k$-fractured. It follows that $f$ is not in a costandard
$k$-wild triangle.

As $f$ is not in a costandard $k$-wild triangle, 
by Lemma~\ref{semi-feral-in-guts}, $f$ is in the
guts of a vertical $3$-separation $(X,Y)$ of $M$, where $f\in X$.
If $(X,Y)$ crosses $(Q'_1,Q_{s+1}\udots Q_k)$, then it is easily
seen that we either contradict the maximality of the flower $\QQ'$
in $M$, or we find that $\QQ'$ is not canonical, contradicting
Lemma~\ref{clean}. If $X$ properly contains $Q'_1$, then, again
as $\QQ'$ is maximal, we deduce that $Y=Q_i\cup Q_{i+1}\udots Q_j$,
where, up to labels, $s<i\leq j\leq k$. But then, as $x\in \cl(Y)$,
we see again that $f$ is loose in $\QQ'$, contradicting 
Lemma~\ref{clean}. Thus $X\subseteq Q'_1$. Now 
$\lambda_{M\ba f}(X-\{f\})=2$. If $X-\{f\}\subseteq Q_i$ for some
$i\in\{1,2,\ldots,s\}$, then we obtain the contradiction that 
$M$ is $k$-fractured. Otherwise, as $\QQ'$ is maximal, we see that,
up to labels, $X-\{f\}=Q_i\cup Q_{i+1}\udots Q_j$, 
where $1\leq i\leq j<s$,
so that $(Q_1\cup Q_2\udots Q_{s-1},Q_s,Q_{s+1},\ldots Q_k)$ 
is a flower
in $M$, contradicting the maximality of $\QQ'$. Therefore $X=Q'_1$.
It follows that $f\in\cl(Q_{s+1}\cup Q_{s+2}\udots Q_k)$ and hence
$(Q_1,Q_2,\ldots,Q_s,Q_{s+1}\udots Q_k\cup\{f\})$ is a swirl-like 
flower in $M$.

Consider the $f$-bad element $g$. By Lemma~\ref{bad-boy5}
$g\in Q_i$ for some $i\in\{1,2,\ldots,s\}$. By 
Lemma~\ref{either-or}, a $k$-fracture for 
$M*g$ either opens the petal $Q_i$ of 
$(Q_1,Q_2,\ldots,Q_s,Q_{s+1}\udots Q_k\cup\{f\})$
or has the property that
$E(M)-Q_i$ is contained in a petal of a $k$-fracture of $M*g$.
Routine arguments show that both cases lead to the contradiction
that $M_{pq}\ba f*g$ is not $k$-coherent.
\end{proof}

The next lemma follows from Lemma~\ref{bad-boy6} 
and Theorem~\ref{1-reduced}.

\begin{lemma}
\label{bad-boy7}
If $M_{pq}\ba f$ is not a $k$-skeleton, then
$f$ is either a feral element or is in a gang-of-three.
\end{lemma}

It is good news
if $M$ has a gang of three
since, by Theorem~\ref{gogang}, we can remove the whole gang-of-three
and keep the property of being a $k$-skeleton. 
We also know that if $f$ belongs to a bogan couple
we can obtain a 4-element reduction. The assumption 
in the next lemma that
$M$ has no gangs or cogangs of three and that $f$ does not belong
to a bogan couple is probably not necessary, but it does
simplify the argument a little.

\begin{lemma}
\label{bad-boy8}
Assume that $M$ has no gangs of cogangs of three and that 
$f$ is not in a bogan couple. If $g$ is $f$-bad, then 
any $k$-fracture of $M*g$ is obtained by opening the petal
$Q'_1$ of $\QQ'$. Moreover $M*g$ is uniquely $k$-fractured.
\end{lemma}

\begin{proof}
By Lemma~\ref{bad-boy4}, $M*g$ is 3-connected and $k$-fractured.
By Lemma~\ref{bad-boy5}, $g\in Q'_1$.
Let $\RR=(R_1,R_2,\ldots,R_u)$ be a $k$-fracture of 
$M*g$. Assume that $\RR$ is not obtained by opening the 
petal $Q_1'$ of $\QQ$. By Lemma~\ref{either-or}, 
we may assume, up to labels,
that $R_1$ is fully closed and that 
$R_1\supseteq Q_{s+1}\cup Q_{s+2}\udots Q_k$. 
Moreover it is routinely seen that $u=k$.

For some $i\in\{1,2,\ldots,s\}$, we have $g\in Q_i$.

\begin{sublemma}
\label{bad-boy8.1}
$|Q_i|>2$.
\end{sublemma}

\subproof
Assume that $|Q_i|=2$.
As $f$ is not in a bogan couple, the flower $\QQ$ is 
canonical, so that $Q_i$ is fully closed.
Now, by Lemma~\ref{2-elt-feral}, $M\ba g$ is $k$-coherent.
The claim follows from this contradiction.
\end{proof}

Assume that $f\in R_j$.

\begin{sublemma}
\label{bad-boy8.2}
The operation $*$ is contraction and $|R_j|=2$.
\end{sublemma}

\subproof
By Lemma~\ref{bad-boy3}, $M\ba f*g$ is 
$3$ connected and by \ref{bad-boy8.1},
$|Q_i|>2$.
Hence 
$(Q_1,Q_2,\ldots,Q_{i-1},Q_i-\{g\},Q_{i+1},\ldots,Q_k)$
is a flower in $M\ba f*g$. It follows easily
that, if $*$ is deletion, then $g\in\cl_{M\ba f}(Q_i)$
and, if $*$ is contraction, then $g\in\cl^*_{M\ba f}(Q_i)$.

Note that 
$(R_1,R_2,\ldots,R_{j-1},R_j-\{f\},R_{j+1},\ldots,R_k)$
is a quasi-flower in $M\ba f*g$.
As $R_1\supseteq Q_{s+1}\udots Q_k$, and
$M_{pq}\ba f*g$ is $k$-coherent, we see that
$R_j-\{f\}$ is a set of loose elements of
$(R_1,R_2,\ldots,R_{j-1},R_j-\{f\},R_{j+1},\ldots,R_k)$.
Assume that $Q_i-\{g\}$ is not contained in $R_1$.
Then $(R_2\cup R_3\udots R_k)-\{f\}\subseteq Q_i$.
But then one readily concludes that
$f\in\cl_M(Q_i)$ giving the contradiction that 
$M$ is $k$-fractured. Therefore 
$Q_i-\{g\}\subseteq R_1$.

By the observation that 
$g\in\cl_{M\ba f}(Q_i)$, or 
$g\in\cl^*_{M\ba f}(Q_i)$, according as to whether 
$*$ is deletion or contraction, we deduce that
$(R_1\cup\{g\},R_2,\ldots,R_{j-1},R_j-\{f\},R_{j+1},\ldots,R_k)$
is a flower in $M\ba f$.

Now $f\in\cl_{M*g}(R_j-\{f\})$. If 
$f\in\cl_M(R_j-\{f\})$, then 
$(R_1\cup\{g\},R_2,\ldots,R_k)$ is a $k$-fracture of $M$.
Hence $*$ is contraction.

It is also the case that $R_j-\{f\}$ is a set of loose
elements in the flower
$(R_1\cup\{g\},R_2,\ldots,R_{j-1},R_j-\{f\},R_{j+1},\ldots,R_k)$
of $M\ba f$. Let $l$ be an initial element of $R_j-\{f\}$
in this flower. For
some $h$, we have $l\in Q_h$. 
Assume that $l$ is a guts element. Then, by 
Lemma~\ref{loose-removable}, $l$ does not expose any
3-separations in $M\ba f,l$. This means that
$(Q_1,Q_2,\ldots,Q_{h-1},Q_h-\{l\},Q_{h+1},\ldots,Q_k)$
is a maximal flower in $M\ba f,l$. Moreover, either
$M\ba f,l$ is $k$-coherent or this is a unique 
$k$-fracture of $M\ba f,l$. It follows from either
Lemma~\ref{in-closure} or Corollary~\ref{unique-fracture}
that $M\ba l$ is $k$-coherent. But it is readily seen that
$l$ is fixed in $M$ contradicting the definition of 
a $k$-skeleton. Thus $l$ is a coguts element.

Assume that $|R_j|>2$. Note that $f\in\cl_{M/g}(R_j-\{f\})$.
Again, by Lemma~\ref{loose-removable}, $l$ does not
expose any 3-separations in $M\ba f/l$, so that
$(Q_1,Q_2,\ldots,Q_{h-1},Q_h-\{l\},Q_{h+1},\ldots,Q_k)$
is a maximal flower in $M\ba f/l$. But 
$l$ is cofixed in $M\ba f$ and $f$ is not comparable with
any element of $M$, so $l$ is cofixed in $M$. Therefore
$M/l$ is $k$-fractured. Again, by 
Lemma~\ref{in-closure} or Corollary~\ref{unique-fracture},
we see that $f$ is in the closure in $M/l$ of a petal of the
above flower. Recall that $g\in Q_i$
and that $f\in\cl_{M/g}(R_j-\{f\})$.
From these facts it follows that we must have 
$f\in \cl_{M/l}(Q_i)$, so that $l\in\cl_M(Q_i\cup\{f\})$.
Let $l'$ be the other end of $R_j-\{f\}$,
then it is also the case that $l'\in \cl_M(Q_i\cup\{f\})$
and, indeed, that $l'\in\cl_M(Q_i\cup\{l\})$. But
$Q_i-\{g\}\subseteq R_1$, and
$(R_1,R_2,\ldots,R_{j-1},R_j-\{f\},R_{j+1},\ldots,R_k)$
is a swirl-like flower. By Lemma~\ref{fine1},
$l'\notin\cl_{M/g}(R_1\cup\{l\})$, so
$l'\notin\cl_M(Q_1\cup\{l\})$. From this contradiction
we can finally deduce that $|R_j|=2$.
\end{proof}

By \ref{bad-boy8.2}, $f$ is in a 
fully closed 2-element petal of 
a $k$-fracture of $M/g$. But 
$f$ is feral and  we have a contradiction 
to Lemma~\ref{2-elt-feral}.
\end{proof}

\section{The Chain Theorem}

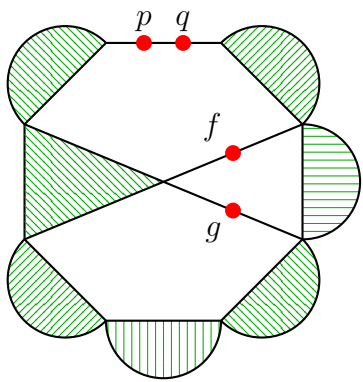
\begin{figure}[h]
\begin{center}
\begin{tikzpicture}[thick,line join=round]
	\coordinate (a) at (0,0);
	\coordinate (s) at (22.5:2);
	\coordinate (t) at (90-22.5:2);
	\coordinate (u) at (135-22.5:2);
	\coordinate (v) at (180-22.5:2);
	\coordinate (w) at (225-22.5:2);
	\coordinate (x) at (270-22.5:2);
	\coordinate (y) at (315-22.5:2);
	\coordinate (z) at (360-22.5:2);
	\coordinate (st) at ($(s)!0.5!(t)$);
	\coordinate[label=above:$p$] (p) at ($(u)!0.33!(t)$);
	\coordinate[label=above:$q$] (q) at ($(t)!0.33!(u)$);
	\coordinate (uv) at ($(u)!0.5!(v)$);
	\coordinate (wx) at ($(w)!0.5!(x)$);
	\coordinate (xy) at ($(x)!0.5!(y)$);
	\coordinate (yz) at ($(y)!0.5!(z)$);
	\coordinate (zs) at ($(z)!0.5!(s)$);
	\node[pattern color=lines,draw,circle through=(s),pattern=north east lines] at (st) {};
	\node[pattern color=lines,draw,circle through=(u),pattern=north west lines] at (uv) {};
	\node[pattern color=lines,draw,circle through=(w),pattern=north east lines] at (wx) {};
	\node[pattern color=lines,draw,circle through=(x),pattern=vertical lines] at (xy) {};
	\node[pattern color=lines,draw,circle through=(y),pattern=north west lines] at (yz) {};
	\node[pattern color=lines,draw,circle through=(z),pattern=horizontal lines] at (zs) {};
	\filldraw[fill=white] (s) -- (t) -- (u) -- (v) -- (w) -- (x) -- (y) -- (z) -- cycle;
	\fill[pattern color=lines,pattern=north west lines] (a) -- (v) -- (w) -- cycle;
	\draw (s) -- (w);
	\draw (v) -- (z);
	\coordinate[label=135:$f$] (f) at ($(a)!0.5!(s)$);
	\coordinate[label=225:$g$] (g) at ($(a)!0.5!(z)$);
	\foreach \pt in {p,q,f,g} \fill[rfill] (\pt) circle (3pt);
\end{tikzpicture}
\label{3-reduced-fig}
\caption{Illustration of a $4$-element move}
\end{center}
\end{figure}

At last we are in a position to prove the more detailed
version of Theorem~\ref{not-sporadic}.
Figure~\ref{3-reduced-fig} illustrates one situation
where a 4-element move is needed. This is the
type of situation that
the proof of Theorem~\ref{inductive-search}
converges to.

\begin{theorem}
\label{inductive-search}
Let $M$ be a nonempty $k$-skeleton. Then at 
least one of the following holds.
\begin{itemize}
\item[(i)] There is an element $e$ such that either 
$M\ba e$ or $M/e$ is a 
$k$-skeleton.
\item[(ii)] There is a clonal pair $\{p,q\}$ 
such that $M\ba p/q$ is a $k$-skeleton.
\item[(iii)] There is a gang of three $\{r,s,t\}$ 
such that $M/r\ba s,t$ is a $k$-skeleton.
\item[(iv)] There is a cogang of three $\{r,s,t\}$ 
such that $M\ba r/s,t$ is a $k$-skeleton.
\item[(v)] The ground set of $M$ consists of clonal 
classes of size two, feral elements and
semi-feral elements.
\end{itemize}
Moreover, in the case that (v) holds, but (i)--(iv) do not,
then every leaf bag of an augmented $3$-tree for
$M$ contains exactly one
clonal pair, and for any such clonal pair $\{p,q\}$  
at least one of the following holds.
\begin{itemize}
\item[(vi)] There is a feral or semi-feral element $f$ such that either
$M\ba p/q\ba f$ or $M\ba p/q/f$ is a $k$-skeleton.
\item[(vii)] There is a pair $f,g$ of feral elements of $M$ such that
$M\ba p/q\ba f/g$ is a $k$-skeleton.
\end{itemize}
\end{theorem}

For the remainder of this section we assume that we are
under the hypotheses of Theorem~\ref{inductive-search}.
We also assume that none of (i), (ii), (iii) or (iv) holds. 
By Theorems~\ref{1-reduced} and \ref{gogang}, (v) holds. 
Moreover $M$ is 2-reduced, so by Lemma~\ref{2-reduce},
every leaf bag of an augmented 3-tree for $M$ consists of a 
single clonal pair. Let $\{p,q\}$ be a peripheral
clonal pair of $M$. It remains to prove that either (vi) or
(vii) hold. As (ii) does not hold
there is an element that is $pq$-annoying. By
Lemmas~\ref{bad-boy6} and \ref{bogan-reduce}
we may also assume the following.
\begin{itemize}
\item[(a)] Any element that is $pq$-annoying is feral.
\item[(b)] No element that is $pq$-annoying belongs
to a bogan couple.
\end{itemize}

Let $f$ be an element of $M$ that is $pq$-annoying. Up to
duality we may assume that $f$ is 
$pq$-annoying for deletion. We are now
effectively under the hypotheses of Section~\ref{life} and in
what follows we use the notational conventions of that section,
so that there is a swirl-like flower
$\QQ=(Q_1,Q_2,\ldots,Q_k)$ in $M\ba f$, and an $s\in\{2,3,\ldots,k-2\}$
such that 
\begin{itemize}
\item[(i)] $Q_t=\{p,q\}$ for some $t\in\{s+1,\ldots,k\}$, and
\item[(ii)] $(Q_1\udots Q_s\cup\{f\},Q_{s+1},\ldots,Q_k)$
is a maximal swirl-like flower of $M$.
\end{itemize}
We also let $Q'_1=Q_1\cup Q_2\udots Q_s\cup\{f\}$
and $\QQ'=(Q'_1,Q_{s+1},\ldots,Q_k)$.

Assume that
(vi) does not hold. Then there is an element $g$ that is
$f$-bad. By, Lemma~\ref{bad-boy5}, 
$g\in Q_i$ for some $i\in\{1,2,\ldots,s\}$. Let
$\QQ''=(Q_1,Q_2,\ldots,Q_{i-1},Q_i-\{g\},Q_{i+1},\ldots,Q_k)$

\begin{lemma}
\label{inductive-search1}
$\QQ''$ uniquely $k$-fractures $M\ba f*g$.
\end{lemma}

\begin{proof}
Consider the quasi-flower $\QQ''$
of $M\ba f*g$. If this refines to a flower of
order greater than $k$, then is is clear that 
$M_{pq}\ba f*g$ is $k$-fractured. Hence
either $\QQ''$
is a maximal flower of order $k$ in $M\ba f*g$,
or $Q_i-\{g\}$ is a set of loose elements of this flower.
Consider the former case. Then, either the lemma holds, of
there is another $k$-fracture of $M\ba f*g$. But then
$M_{pq}\ba f*g$ is certainly not $k$-coherent.
Thus, if the lemma fails, the latter case holds. Assume that 
we are in this case.

The flower 
$\QQ''=(Q_1,Q_2,\ldots,Q_{i-1},Q_i-\{g\},Q_{i+1},\ldots,Q_k)$
has order $k-1$. By Lemma~\ref{bad-boy8},
a $k$-fracture of $M*g$ is obtained by opening the petal
$Q'_1$ of $\QQ'$. 
Thus $\QQ''$ induces a $k$-fracture in $M*g$.
By Lemma~\ref{lose-coherence} there is a set $Q'$ of loose
elements of $\QQ''$ that has the property that 
$Q'\cup\{f\}$ is a tight petal of the induced $k$-fracture
in $M*g$. As $f$ is not in a bogan couple, the flower
$\QQ$ is canonical, so the only loose elements of 
$\QQ''$ are in $Q_i-\{g\}$. Therefore $Q'=Q_i-\{g\}$,
so that 
$(Q_1,Q_2,\ldots,Q_{i-1},(Q_i-\{g\})\cup\{f\},Q_{i+1},\ldots,Q_k)$
is a $k$-fracture of $M*g$, giving the contradiction that 
$(Q_1,Q_2,\ldots,Q_{i-1},Q_i\cup\{f\},Q_{i+1},\ldots,Q_k)$
is a $k$-fracture of $M$. The lemma follows from this 
contradiction.
\end{proof}

As $\QQ''$ uniquely fractures $M\ba f*g$, and 
$M*g$ is $k$-fractured, there is a petal $Q$
of a flower equivalent to $\QQ'$ such that 
$f\in\cl_{M*g}(Q)$.  Certainly $Q\neq \{p,q\}$ and it follows
that $M_{pq}*g$ is $k$-coherent. 
This establishes the next lemma.

\begin{lemma}
\label{inductive-search2}
The element $g$ is $pq$-annoying and  hence
$g$ is feral and not in a bogan couple. Consequently a 
$k$-fracture of $M*g$ is canonical.
\end{lemma}

\begin{lemma}
\label{inductive-search3}
The element $g$ is $f$-bad for contraction.
\end{lemma}

\begin{proof}
Assume that $g$ is $f$-bad for deletion.
By Lemma~\ref{inductive-search1}, 
$\QQ''$ uniquely $k$-fractures $M\ba f,g$. 
But $M\ba g$ is $k$-fractured, so
$f$ is in the closure of a petal of some flower equivalent to
$\QQ''$. One consequence of this is that $s=2$ (otherwise either 
$f\in\cl(Q_1\cup Q_2\udots Q_{s-1})$ or 
$f\in\cl(Q_2\cup Q_3\udots Q_{s})$ so that
$\QQ'=(Q_1\cup Q_2\udots Q_s\cup\{f\},Q_{s+1},\ldots,Q_k)$ 
is not a maximal
flower in $M$.) We may now also assume that $i=2$, so that
$\QQ''=(Q_1,Q_2-\{g\},Q_3,\ldots,Q_k)$ and $Q'_1=Q_1\cup Q_2$.
Another consequence is that for some maximal fan 
$H$ between $Q_1$ and $Q_2-\{g\}$, we have $f\in\cl(Q_1\cup H)$.

We now consider possibilities for the fan $H$. Let $Q''_1=Q_1-H$,
$Q''_2=Q_2-(H\cup\{g\})$, and assume that $H=(h_1,h_2,\ldots,h_t)$
is ordered from $Q_1$ to $Q_2-\{g\}$ in $\QQ''$. 
Thus the quasi-flower 
$(Q''_1,H,Q''_2,Q_3,\ldots,Q_k)$ in $M\ba f,g$ is equivalent
to $\QQ''$. Moreover $f\in\cl(Q''_1\cup H)$ 
and $g\in\cl(Q''_2\cup H)$. Note also that we have symmetry between
$f$ and $g$.

\begin{sublemma}
\label{subsearch1}
The elements $h_1$ and $h_t$ are coguts elements of $H$.
\end{sublemma}

\subproof
Assume that $h_1$ is a guts element. Then
$h_1\in\cl(Q''_1)$ so that $h_1$ is a loose element of the
$k$-fracture $(Q''_1,H\cup Q''_2\cup\{g\},Q_3,\ldots,Q_k)$
of $M\ba f$. This contradicts the fact that $f$ is not in a
bogan or cobogan couple. Thus $h_1$ is a coguts element of $H$,
and so too is $h_t$.
\end{proof}

\begin{sublemma}
\label{subsearch2}
$f\notin\cl(Q''_1\cup\{h_1,h_2,\ldots,h_{t-1}\})$ and
$g\notin\cl(Q''_2\cup\{h_2,h_3,\ldots,h_t\})$.
\end{sublemma}

\subproof
If $f\in\cl(Q''_1\cup\{h_1,h_2,\ldots,h_{t-1}\})$,
then $h_t$ is a loose element of a $k$-fracture of $M\ba g$
contradicting the fact that $g$ is not in a bogan 
or cobogan couple. The
claim follows from this observation and the symmetry between $f$
and $g$.
\end{proof}

\begin{sublemma}
\label{subsearch3}
$f\in\cl(Q''_1\cup\{h_t\})$ and $g\in\cl(Q''_2\cup\{h_1\})$.
\end{sublemma}

\subproof
Assume that $f\notin\cl(Q''_1\cup\{h_t\})$. Note that this means
that $t>1$. Consider $M\ba f,g/h_t$. 
It follows from Lemma~\ref{loose-removable} and the fact that
$M\ba f,g$ is uniquely $k$-fractured by $\QQ''$, that
$M\ba f,g/h_t$ is uniquely $k$-fractured by
$(Q''_1,\{h_1,h_2,\ldots,h_{t-1}\},Q''_2,Q_3,\ldots,Q_k)$.
Now $g\in\cl_{M/h_t}(Q''_2\cup\{h_1,h_2,\ldots,h_{t-1}\})$, but
$g\notin \cl_{M/h_t}(Q''_2\cup\{h_2,\ldots,h_{t-1}\})$. 
Therefore there is a circuit $C$ in $M/h_t$ such that
$h_1\in C$ and $C\subseteq \{h_1,h_2,\ldots,h_{t-1}\}\cup Q''_2$.
It follows that the flower 
$(Q''_1,\{g,h_1,h_2,\ldots,h_{t-1}\}\cup Q''_2,Q_3,\ldots,Q_k)$
of $M\ba f/h_t$ is canonical. Evidently it is also the unique
$k$-fracture of $M\ba f/h_t$. Now $f\notin\cl_M(Q''_1\cup\{h_t\})$,
so $f\notin\cl_{M/h_t}(Q''_1)$, and $f$ is certainly not in the closure
of any other petal of the flower. We
conclude, by Corollary~\ref{unique-fracture}, that 
$M/h_t$ is $k$-coherent. But $M$ is a 2-reduced $k$-skeleton
with no gangs or cogangs of three. So this means that
$h_t$ is in a clonal pair. But $h_t$ evidently has no clone and
the claim follows from this contradiction.
\end{proof}

\begin{sublemma}
\label{subsearch4}
$|H|\in\{1,3\}$.
\end{sublemma}

\subproof
Otherwise $F$ contains a triangle $T$ in $M\ba f,g$ and hence in 
$M$. By Lemma~\ref{not-awkward} $T$ cannot be $k$-wild. But $T$
certainly cannot be a clonal triple of $M$ and the claim follows.
\end{proof}

\begin{sublemma}
\label{subsearch5}
The elements of $H$ are feral.
\end{sublemma}

\subproof
It is easily seen that the elements of $H$ are either feral or
semi-feral. We omit the routine verification that they 
are not semi-feral.
\end{proof}

We first consider that case that $|H|=3$.
Let $\QQ'''=(Q''_1\cup\{f,h_3\},Q''_2\cup\{g,h_1\},Q_3,\ldots,Q_k)$.

\begin{sublemma}
\label{subsearch6}
$\QQ'''$ is a $k$-fracture of $M\ba h_2$.
\end{sublemma}

\begin{proof}
Consider $M\ba f,g,h_2/h_3$. Observe that 
$(Q''_1,Q''_2\cup\{h_3\},Q_3,\ldots,Q_k)$ is a 
swirl-like flower in this matroid and that $h_3$ is in the 
coclosure of $Q''_1$. But $\{h_1,h_3\}$ is a series pair
in $M\ba f,g,h_2$ so that 
$(Q''_1,Q''_2\cup\{h_1,h_3\},Q_3,\ldots,Q_k)$
is a swirl like flower in $M\ba f,g,h_2$.
Moreover, $h_3\in\cl^*_{M\ba f,g,h_2}(Q''_1)$,
so that the partition
$(Q''_1\cup\{h_3\},Q''_2\cup\{h_1\},Q_3,\ldots,Q_k)$
has the property that the union of every consecutive pair
of sets in the cyclic order is exactly 3-separating. (This
partition is not
technically a flower according to our definition as the series
pair is split between two petals.) By \ref{subsearch3},
$f\in\cl(Q''_1\cup\{h_3\})$ and $g\in\cl(Q''_1\cup\{h_1\})$.
By \ref{subsearch5}, $M\ba h_2$ is 3-connected.
Therefore 
$\QQ'''$ is a flower in $M\ba h_2$. It now follows easily
that it is swirl-like and tight, establishing the claim.
\end{proof}

Note that $h_2$ blocks both $Q''_1\cup\{f,h_3\}$
and $Q''_2\cup\{g,h_1\}$. Thus, by Lemma~\ref{display-way},
there is a $k$-fracture $\PP=(P_1,P_2,\ldots,P_m)$
of $M/h_2$ such that $\QQ'''$ and $\PP$ form a feral display
for $f_2$ in $M^*$. Mimicking the notation in the definition of 
feral display as closely as possible, we may assume, for some
$i\in\{3,4,\ldots,m-1\}$, that 
$P_2\cup P_3\udots P_i\subseteq Q''_2\cup\{g,h_1\}$, and
$P_{i+1}\cup P_{i+2}\udots P_m\subseteq Q''_1\cup\{f,h_3\}$.
Let $Z_1=P_1\cap(Q''_1\cup\{f,h_3\})$ and
$Z_2=P_1\cap(Q''_2\cup\{g,h_1\})$.

\begin{sublemma}
\label{subsearch7}
$\{f,h_3\}\subseteq Z_1$ and $\{g,h_1\}\subseteq Z_2$.
\end{sublemma}

\begin{proof}
Say $f\in P_j$ for some $j\in\{i+1,i+2,\ldots,m\}$.
As $f$ is feral it follows from Lemma~\ref{2-elt-feral}
that $|P_j|\geq 3$. Using this fact and the fact that $f$
is feral, we deduce that $f$ is not on the guts of  $P_j$.
It now follows from Lemma~\ref{either-or}, or 
Lemma~\ref{either-or-more}, that any $k$-fracture
of $M\ba f$ is either obtained by opening the petal $P_i$
in the flower 
$(P_1\cup P_2\udots P_i\cup \{h_2\},P_{i+1},\ldots,P_m)$
or has the property that $E(M)-P_j$ is contained in a petal.
Such a $k$-fracture is certainly not $\QQ$, the flower
that we know uniquely $k$-fractures $M\ba f$. Hence $f\in Z_1$
and, symmetrically, $g\in Z_2$.

Assume that $h_3\in P_j$. Arguing as before we have $|P_j|>2$.
Also $h_3$ is not in the coguts of $P_j$ as otherwise it is not feral.
Thus $h_3\notin\cl^*(E(M)-P_j)$. But then, as $f\notin P_j$
we see that $h_3\notin\cl^*_M(\{h_1,h_2,f,g\})$ so that
$h_3\notin\cl^*_{M\ba f,g}(\{h_1,h_2\})$, contradicting the
fact that $\{h_1,h_2,h_3\}$ is a triad of $M\ba f,g$. Therefore
$h_3\in Z_1$ and symmetrically $h_1\in Z_2$ as required.
\end{proof}

By \ref{subsearch7}, $\{h_1,h_3,f,g\}\subseteq P_1$.
But, by the properties of the feral display in $M^*$, the element
$h_2$ coblocks $P_1$ in $M$, so that $h_2\notin\cl^*(P_1)$.
Therefore $h_2\notin\cl^*_M(\{h_1,h_3,f,g\})$ and again we 
contradict the fact that $\{h_1,h_2,h_3\}$
is a triad in $M$.

Assume that $|H|=1$ so that $H=\{h_1\}$. Clearly
$(Q''_1\cup\{f\},Q''_2\cup\{g\},Q_3,\ldots,Q_k)$
is a $k$-fracture of $M/h_1$. Moreover, 
one readily checks that it is unique. It follows that $M_{pq}/h_1$
is $k$-coherent and hence $h_1$ is $pq$-annoying. Therefore
$h_1$ is not in a bogan or cobogan couple.
By \ref{subsearch5}, $h_1$
is feral. Moreover $h_1$ coblocks both
$Q''_1\cup\{f\}$ and $Q''_2\cup \{g\}$.
Thus, by Lemma~\ref{display-way}, there is a feral display for 
$h_1$ in $M$. By examining the definition of feral display
we observe that there are petals $P_i$ and $P_{i+1}$
of a $k$-fracture of $M\ba h_1$ such that 
$P_i\subseteq Q''_1\cup\{f\}$, $P_{i+1}\subseteq Q''_2\cup\{g\}$,
$\lambda_M(P_i)=\lambda_M(P_{i+1})=2$, and 
$\sqcap_M(P_i,P_{i+1})=1$. But $\sqcap_M(Q''_1,Q''_2)=0$,
so we may assume without loss of generality that
$f\in P_i$. But $f\notin\cl(Q''_1)$, so that 
$f$ is a coloop of $M|P_i$ and hence 
$\lambda_{M\ba f}(P_i-\{f\})=1$. As $M\ba f$
is 3-connected, it must be the case that $|P_i|=2$.
As $h_1$ is not in a bogan or cobogan couple, the flower
$\PP$ is canonical. Therefore $P_i$ is fully closed in 
$M\ba h_1$. But now Lemma~\ref{2-elt-feral} implies 
the contradiction that
$h_1$ is not feral. All cases lead to a contradiction and the
lemma follows.
\end{proof}

From now on we assume that $g$ is $f$-bad for contraction. Recall
that $Q'_1=Q_1\cup Q_2\udots Q_s$ and, for some 
$i\in\{1,2,\ldots,s\}$, we have $g\in Q_i$. We may assume, up to 
labels, that $g\notin Q_1$.

\begin{lemma}
\label{get2}
$i=s=2$.
\end{lemma}

\begin{proof}
Assume that $s\neq 2$. By Lemma~\ref{inductive-search1},
$\QQ''=(Q_1,Q_2,\ldots,Q_{i-1},Q_i-\{g\},Q_{i+1},\ldots,Q_k)$
uniquely $k$-fractures $M\ba f/g$.
But $M/g$ is $k$-fractured. Thus, by Lemma~\ref{keep-fracture},
$f$ is in the span of a petal of a flower equivalent to $\QQ''$.
Certainly the petal is not equivalent to $Q_i-\{g\}$ as 
otherwise $M$ is $k$-fractured. An elementary argument,
based on the fact that $\QQ'$ is a maximal flower in $M$
shows that we may now assume that $i=s$ and that
$f$ is in the span of a petal equivalent to $Q_1$.
As $\QQ$ is a canonical flower in $M\ba f$, and $s\neq 2$,
the petal $Q_1$ is fully closed in $\QQ''$. Hence
$(Q_1\cup\{f\},Q_2,\ldots,Q_{s-1},Q_s-\{g\},Q_{s+1},\ldots,Q_k)$
is a $k$-fracture of $M/g$.

Consider the $k$-fracture $\QQ=(Q_1,Q_2,\ldots,Q_k)$
of $M\ba f$.  
Certainly $f$ blocks $Q_s$ and, as 
$\QQ'=(Q_1\cup Q_2\udots Q_s\cup\{f\},Q_{s+1},\ldots,Q_k)$ 
is a maximal swirl-like flower in 
$M$, it also must be the case that $f$ blocks 
$Q_2\cup Q_3\udots Q_s$. Thus
$f$ is not 2-spanned by $\QQ$. 
Now, by the definition of a feral display, one sees that
$f$ has a feral display in $M$ rather than $M^*$.
As $f$ blocks $Q_s$, we see that, with the labelling given in 
the definition of feral display, we have
$(P_1,P_2,\ldots,P_i,P_{i+1},\ldots,P_k)
=(Q_s,Q_{s+1},\ldots,Q_k,Q_1,\ldots, Q_k)$, where $P_i=Q_k$.
It now follows from the properties of a feral display that
$f\in\cl(Q_s\cup Q_{s+1}\udots Q_k)$.
Therefore $f\in\cl_{M/g}((Q_s-\{g\})\cup Q_{s+1}\udots Q_k)$.
But this implies that $f$ is a loose element of the $k$-fracture
$(Q_1\cup\{f\},Q_2,\ldots,Q_{s-1},Q_s-\{g\},Q_{s+1},\ldots,Q_k)$,
of $M/g$ contradicting the fact that $g$ is not in a bogan couple.
The lemma follows from this contradiction.
\end{proof}

We now know that $\QQ''=(Q_1,Q_2-\{g\},Q_3,\ldots,Q_k)$
is a $k$-fracture in $M\ba f/g$. As in the proof of the previous
lemma we have a, possibly empty, fan $H$ of loose elements
between $Q_1$ and $Q_2-\{g\}$ such that
$f\in\cl_{M/g}(Q_1\cup H)$. Let 
$Q''_1=Q_1-H$ and $Q''_2=Q_2-(H\cup \{g\})$. Assume that 
$H=(h_1,h_2,\ldots,h_t)$ is ordered from $Q_1$ to $Q_2-\{g\}$ in 
$\QQ''$. In other words, we have a  quasi-flower
$(Q''_1,H,Q''_2,Q_3,\ldots,Q_k)$ in $M\ba f/g$ is equivalent to 
$\QQ''$. Moreover 
$(Q''_1\cup H\cup\{f\}, Q_2'',Q_3,\ldots,Q_k)$
and $(Q''_1,Q''_2\cup H\cup\{g\},Q_3,\ldots,Q_k)$
are canonical $k$-fractures of $M/g$ and $M\ba f$
respectively. As $f$ and $g$ are both $pq$-annoying,
these $k$-fractures are canonical.

Our next goal is to get rid of the irritating fan $H$.

\begin{lemma}
\label{fan-free}
The flower $\QQ''$ is canonical. In particular, the
fan $H$ is empty.
\end{lemma}

\begin{proof}
It is easily seen that the only possible loose elements
of $\QQ''$ are in $H$. 
Assume that $H\neq\emptyset$.

\begin{sublemma}
\label{fan-free1}
The element $h_1$ is a guts element of $H$ and $h_2$ is a coguts
element. Consequently $|H|\geq 2$.
\end{sublemma}

\subproof
Consider $h_t$, the last element of $H$. If 
$h_t$ is a guts element of $H$, then $h_t\in\cl_{M\ba f/g}(Q''_2)$
so that $h_t\in\cl_{M/g}(Q''_2)$, contradicting the fact that
the $k$-fracture 
$(Q''_1\cup H\cup\{f\},Q''_2,Q_3,\ldots,Q_k)$ of $M/g$ is canonical.
Hence $h_t$ is a guts element of $H$ and, dually, 
$h_1$ is a guts element of $H$.
\end{proof}

We omit the easy proof of the next claim.

\begin{sublemma}
\label{fan-free2}
The elements of $H$ are feral.
\end{sublemma}

Consider possible $k$-fractures of $M\ba h_1$.

\begin{sublemma}
\label{fan-free3}
There is an $i\in\{1,2,\ldots,t\}$ such that
$$(Q''_1\cup\{f,h_2,h_3,\ldots,h_i\},Q''_2\cup 
\{g,h_{i+1},h_{i+2},\ldots,h_t\},Q_3,\ldots,Q_k)$$ 
is a $k$-fracture of $M\ba h_1$.
\end{sublemma}

\subproof
Consider the matroid $M\ba h_1/g$.
It is easily seen that $M\ba h_1/g$ is 3-connected.
As $h_1\in\cl_{M/g}(Q''_1)$, and 
$f\in\cl_{M/g}(Q''_1\cup H)$, we have 
$f\in\cl_{M/g}(Q''_1\cup\{h_2,h_3,\ldots,h_t\})$.
Thus
$(Q''_1\cup\{f,h_2,h_3,\ldots,h_t\},Q''_2,Q_3,\ldots,Q_k)$
is a $k$-fracture of $M\ba h_1/g$. Let 
$i$ be the least integer such that
$(Q''_1\cup\{f,h_2,h_3,\ldots,h_i\},Q''_2\cup 
\{h_{i+1},h_{i+2},\ldots,h_t\},Q_3,\ldots,Q_k)$
is also a $k$-fracture of $M\ba h_1/g$.

As $h_1$ is feral, $M\ba h_1$ is $k$-fractured.
Let $(O_1\cup \{g\},O_2,\ldots, O_k)$ be a 
$k$-fracture of $M\ba h_1$. By the structure of swirl-like flowers
either $g\in\cl(O_k\cup O_1)$ or $g\in\cl(O_1\cup O_2)$.
Up to labels we may assume that $g\in\cl(O_1\cup O_2)$.
The partition $(O_1\cup O_2,O_3,\ldots,O_k)$ is a swirl-like
flower in the 3-connected matroid $M\ba h_1/g$.
Assume that this flower is not comparable with 
$(Q''_1\cup\{f,h_2,h_3,\ldots,h_i\},Q''_2\cup
\{h_{i+1},h_{i+2},\ldots,h_t\},Q_3,\ldots,Q_k)$.
By Lemma~\ref{in-petal}, either $O_1\cup O_2$ is contained
in the full closure of a petal of 
$(Q''_1\cup\{f,h_2,h_3,\ldots,h_i\},Q''_2\cup
\{h_{i+1},h_{i+2},\ldots,h_t\},Q_3,\ldots,Q_k)$
or $O_3\cup O_4\udots O_k$ is contained in the full closure
of a petal of this flower. The latter case routinely leads to a 
contradiction of the fact that $M_{pq}/g$ is $k$-coherent.
The only non-contradictory possibility in the former case is
if $O_1\subseteq Q''_2\cup\{h_{i+1},h_{i+2},\ldots,h_t\}$.
But, in this case, the claim holds.
\end{proof}

Consider the $k$-fracture
$(Q''_1\cup\{f,h_2,h_3,\ldots,h_i\},Q''_2\cup 
\{g,h_{i+1},h_{i+2},\ldots,h_t\},Q_3,\ldots,Q_k)$
of $M\ba h_1$ given by \ref{fan-free3}.
As $h_1$ is a feral element of $M$, there is a $k$-fracture
$\PP=(P_1,P_2,\ldots,P_m)$ of $M/h_1$ such that
$\PP$ and $(Q''_1\cup\{f,h_2,h_3,\ldots,h_i\},Q''_2\cup 
\{g,h_{i+1},h_{i+2},\ldots,h_t\},Q_3,\ldots,Q_k)$
form a feral display in $M$ or $M^*$. Observe that
$h_1$ is spanned in $M$ by
$(Q''_1\cup\{f,h_2,h_3,\ldots h_i\})\cup 
(Q''_2\cup\{g,h_{i+1},h_{i+2},\ldots,h_t\})$ but not by
either $Q''_1\cup\{f,h_2,h_3,\ldots h_i\}$
or $Q''_2\cup\{g,h_{i+1},h_{i+2},\ldots,h_t\}$, as otherwise
$M$ is $k$-fractured. This shows that we have a feral display in 
$M$. By the definition of feral display, we may assume
that there is an $l\in\{2,3,\ldots,m-1\}$ such that
$P_2\cup P_3\udots P_l\subseteq Q''_2\cup
\{g,h_{i+1},h_{i+2},\ldots,h_t\}$,
$P_{l+1}\cup P_{l+2}\udots P_m\subseteq Q''_1\cup
\{f,h_1,h_2,\ldots,h_i\}$, and
$Q_3\cup Q_4\udots Q_k\subseteq P_1$.
As usual we let
$Z_1=(Q''_1\cup\{f,h_1,h_2,\ldots,h_i\})\cap P_1$,
and 
$Z_2=(Q''_2\cup\{g,h_{i+1},h_{i+2},\ldots,h_t\})\cap P_1$.

\begin{sublemma}
\label{fan-free4}
$Q''_1\cup\{f,h_1,h_2,\ldots,h_i\}\subseteq Z_1$ and
$Q''_2\cup\{g,h_{i+1},h_{i+2},\ldots,h_t\}\subseteq Z_2$.
\end{sublemma}

\subproof
Say $f\notin Z_1$. Then $f\in P_u$ for some
$u\in\{l+1,l+2,\ldots,m\}$. But then,
by Lemma~\ref{either-or}, or Lemma~\ref{either-or-more},
$(Q_1,Q_2,\ldots,Q_k)$ is not a $k$-fracture of $M\ba f$.
Thus $f\in Z_1$ and similarly $g\in Z_2$.

If $\{h_2,h_3,\ldots,h_t\}\subseteq Z_1\cup Z_2$,
then the claim holds. Assume otherwise.
Let $\omega$ be the first element of $\{1,2,\ldots,t\}$
such that $h_\omega\notin Z_1\cup Z_2$.
Say $\omega\leq i$. Then
$Q''_1\cup\{f,g,h_1,h_2,\ldots,h_{\omega-1}\}
\subseteq E(M)-(P_{l+1}\cup P_{l+2}\udots P_m)$,
and $h_\omega\in\clstar_{M\ba f/g}(Q''_1\cup
\{f,g,h_1,h_2,\ldots, h_{\omega-1}\})$, so that
$h_\omega\in\clstar_M(E(M)-(P_{l+1}\cup P_{l+2}\udots P_m))$.
It follows from this that, if $m\neq l+1$, 
then  $h_\omega$ is a loose element of a flower
of $M$ of order at least three, contradicting
Lemma~\ref{clean}.
Say $l+1=m$. By Lemma~\ref{2-elt-feral},
$|P_m|>2$. Thus $h_w$ is in the guts or coguts of a 
3-separation of $M$, contradicting the fact that 
$h_w$ is feral. The same argument works in the case
that $\omega >i$ and the sublemma follows.
\end{proof}

By the properties of a feral display,
$\sqcap_M(P_l,P_{l+1})=1$. But,
by \ref{fan-free4}, $P_l\subseteq Q''_1$
and $P_{l+1}\subseteq Q''_2$. Hence 
$\sqcap_M(Q''_1,Q''_2)\geq 1$, that is,
$\sqcap_{M\ba f}(Q''_1,Q''_2)\geq 1$.
But $g\notin\cl(Q''_1)\cap \cl(Q''_2)$,
so $\sqcap_{M\ba f/g}(Q''_1,Q''_2)\geq 1$.
However the existence of the coguts element $h_t$ of the fan
$H$ in $M\ba f/g$ ensures that 
$\sqcap_{M\ba f/g}(Q''_1,Q''_2)=0$.

It follows from the above contradiction that a feral display
does not exist for $h_1$ so that $h_1$ is not feral,
and we can at last conclude that $H\neq \emptyset$.
\end{proof}

The next lemma shows, thankfully, that we do not have to dig any
deeper.

\begin{lemma}
\label{search-over}
There is no element $h$ such that $h$ is fixed in
$M_{pg}\ba f/g$ and $M_{pq}\ba f/g\ba h$ is
$k$-coherent.
\end{lemma}

\begin{proof}
Assume that the lemma fails so that we have an element 
$h$ that is fixed in
$M_{pg}\ba f/g$ and $M_{pq}\ba f/g\ba h$ is $k$-coherent.
It is clear that $h$ is either feral or semi-feral, and,
in the semi-feral case, $M\ba f$ is 3-connected. In other words,
$M\ba h$ is 3-connected and $k$-fractured.
By Lemma~\ref{fan-free}, 
$\QQ''=(Q_1,Q_2-\{g\},Q_3,\ldots,Q_k)$ canonically
$k$-fractures $M\ba f/g$, and this $k$-fracture is unique.
Moreover $(Q_1\cup\{f\},Q_2-\{g\},Q_3,\ldots,Q_k)$
and $(Q_1,Q_2,\ldots,Q_k)$ are unique canonical $k$-fractures of
$M/g$ and $M\ba g$ respectively. 
Recall also that $Q_t=\{p,q\}$ for some $t\in\{3,4,\ldots,k\}$.

We omit the argument that shows that either $h\in Q_1$
or $h\in Q_2-\{g\}$ which is similar to, but easier than, the arguments
below.

Consider the case that  that $h\in Q_2-\{g\}$. We know that
$(Q_1\cup\{f\},Q_2-\{g\},Q_3,\ldots,Q_k)$ is a canonical 
unique $k$-fracture
of $M/g$ and $h\in Q_2-\{g\}$. Consider $M/g\ba h$. 
This matroid is clearly $3$-connected. 
If, for some partition $(Q'_2,Q''_2)$ of $Q_2-\{g\}$, the 
partition  $(Q_1\cup\{f\},Q'_2,Q''_2,Q_3,\ldots,Q_k)$ is a 
tight swirl-like flower of $M/g\ba h$, 
then it is clear that $M_{pq}\ba f/g\ba h$
is not $k$-coherent. Hence 
$(Q_1\cup\{f\},Q_2-\{g,h\},Q_3,\ldots,Q_k)$
is a maximal swirl-like quasi-flower in $M/g\ba h$.
This flower is either a $k$-fracture, or else $Q_2-\{g,h\}$
is a set of loose elements of the flower. 
In the latter case we have $M/g\ba h$
is $k$-coherent, as any other fracture would again contradict the fact
that $M_{pq}\ba f/g\ba h$ is $k$-coherent. As $M\ba h$ is $k$-fractured,
it follows that 
$(Q_1\cup\{f\},Q_2-\{h\},Q_3,\ldots,Q_k)$ is a maximal $k$-fracture
of $M\ba h$. Note that $g\in\cl(Q_1\cup\{f\})$,
as otherwise $(Q_1\cup\{f\},Q_2,\ldots,Q_k)$ is a $k$-fracture
of $M$. But $g\in Q_2-\{h\}$, so tht $g$ is loose in 
$(Q_1\cup\{f\},Q_2-\{h\},Q_3,\ldots,Q_k)$.
This gives the contradiction that $\{g,h\}$ is a bogan couple.

Consider the case that $h\in Q_1$. 
Assume that $M\ba f,h$ is not 3-connected.
Then, as $M\ba f,h/g$ is 3-connected, there is a series pair
$\{g,g'\}$ in $M\ba f,h$ so that $\{g,g',h\}$ is a triad of
$M\ba f$. But $\QQ$ is a canonical flower of 
$M\ba f$ and it follows that any triad must be contained in a single
petal of $\QQ$, contradicting the fact that $h\in Q_1$ and 
$g\in Q_2$.
Thus $M\ba f,h$ is $3$-connected. Arguing as in the previous
case we see that  
$(Q_1-\{h\},Q_2,\ldots,Q_k)$ either uniquely 
$k$-fractures $M\ba f,h$
or  $Q_1-\{h\}$ is a set of loose elements of this flower. But
$M\ba h$ is $k$-fractured. In the case that $M\ba f,h$
is $k$-fractured, it must be the case that this 
$k$-fracture is $(Q_1-\{h\},Q_2,Q_3,\ldots,Q_k)$,
as otherwise it is easily checked that
$M_{pq}\ba f/g\ba h$ is $k$-fractured. In either case
we deduce that $M_{pq}\ba f,h$ is $k$-coherent and we have contradicted
Lemma~\ref{inductive-search3}.
\end{proof}

Assume that $M_{pq}\ba p/q$ is not a $k$-skeleton. 
By Lemma~\ref{search-over}, there is no element $h$
such that either $h$ is fixed in $M_{pq}\ba f/g$ and
$M_{pq}\ba f/g\ba h$ is $k$-coherent. By the symmetry
between $f$ and $g$ under duality, it also follows
from Lemma~\ref{search-over}, that there is no element
$h$ such that $h$ is cofixed in $M_{pq}\ba f/g$ and
$M_{pq}\ba f/g,h$ is $k$-coherent. Thus
$M_{pq}\ba f/g$ is indeed a $k$-skeleton. This completes the
proof of Theorem~\ref{inductive-search}.

\chapter{Paths of $3$-separations}
\label{paths-of-3-separations}

\section{Introduction}

Recall that if $A$ and $B$ are disjoint sets 
of elements of the matroid
$M$, then $\kappa(A,B)$ denotes the minimum of 
$\lambda(X',Y')$ over 
all partitions $(X',Y')$ of 
$E(M)$ with $X\subseteq X'$ and $Y\subseteq Y'$.
A {\em path} of $3$-separations
\index{path of $3$-separations} 
in a matroid
$M$ is a partition $\PP=(P_0,P_1,\ldots,P_l)$ of $E(M)$ into  
subsets such that $\kappa(P_0,P_l)=2$ and
$\lambda(P_0\udots P_i)=2$ , for all 
$i\in\{0,1,\ldots,l-1\}$.  
If $i\in\{0,1,\ldots,l\}$, then 
$P_i$ is a {\em step} 
\index{step of a path of $3$-separations}
of the path; $P_0$ and 
$P_l$ are {\em end steps}
\index{end step of a path of $3$-separations}
and otherwise $P_i$ is an {\em internal step}.
\index{internal step of a path of $3$-separations} 

We allow the possibility that internal steps 
can be empty. A {\em solid} 
path has no empty steps. If $\PP$ has 
$l+1$ nonempty steps, then the 
{\em length} of $\PP$ is $l$.

Recall that $\Lambda_m$ denotes the rank-$m$
free spike and recall that $\eq$ denotes the class
of matroids that has no $U_{2,q+2}$-, $U_{q,q+2}$-
or $\Lambda_q$-minor. In this chapter we begin
the task of controlling the structure of a $k$-skeleton
in $\eq$. Our primary goal is to prove 
that $k$-skeletons in $\eq$ cannot 
have arbitrarily long paths of $3$-separations.
In other words, we prove

\begin{theorem}
\label{main1}
Let $k\geq 5$ and $q\geq 2$ be integers. 
Then there is a function $f_{\ref{main1}}(k,q)$
such that, if $M$ is a $k$-skeleton with a path of $3$-separations
of length $f_{\ref{main1}}(k,q)$, 
then $M\notin \eq$.
\end{theorem}

We also prove a number of other lemmas that will be of use
later in the paper.
Before ploughing on into the technicalities of the proof we note a 
corollary of Theorem~\ref{main1}.

\begin{corollary}
\label{main2}
Let $k\geq 5$ and $q\geq 2$  be integers. 
Let $M$ be a $k$-coherent matroid with a path of $3$-separations
$\PP$ of length $f_{\ref{main1}}(k,q)+2$ such that each step of
$\PP$ contains an element that is neither fixed nor cofixed. 
Then $M\notin \eq$. 
\end{corollary}

\begin{proof}
Assume that $M$ satisfies the hypotheses of the corollary. Let 
$n=f_{\ref{main1}}(k,q)$.
Then $M$ has a path 
$\PP=(P_0,P_1,\ldots,P_n)$ such that each step of $\PP$
contains an element that is neither fixed nor cofixed and 
both $P_0$
and $P_n$ contain at least two such elements.
If $M$ is a $k$-skeleton, the result follows immediately, so assume
that $M$ is not a $k$-skeleton. 
Then, up to duality, we may assume that
there is a fixed element $z$ such that
$M\ba z$ is $k$-coherent. For some 
$i\in\{0,1,\ldots,l\}$, we have $z\in P_i$.
If $x$ is not fixed in $M$,
then $x$ is not fixed in $M\ba z$. Assume that $x$ is not cofixed
in $M$. Then, as $z$ is fixed in $M$, it is not the case that
$z$ is freer than $x$ in $M$. Thus, by Corollary~\ref{born-free},
$x$ is not cofixed in $M\ba x$. It is now clear that
$(P_0,\ldots,P_{i-1},P_i-\{z\},P_{i+1},\ldots,P_n)$ is a
path of $3$-separations in $M\ba z$ each step of which contains
an element that is not fixed or cofixed. The result now
follows from an obvious induction.
\end{proof}

\section{Some Preliminaries}

We begin by developing 
some more terminology for paths of $3$-separations.
A nonempty internal step $P_i$ of $\PP=(P_0,P_1,\ldots,P_l)$
is {\em prime}
\index{prime step} 
if there is no partition $(P_{i1},P_{i2})$
of $P_i$ into nonempty subsets such that 
$(P_0,P_1,\ldots,P_{i-1},P_{i1},P_{i2},P_{i+1},\ldots,P_l)$
is a path of $3$-separations of $M$. The path $\PP$
is {\em maximal}
\index{maximal path of $3$-separations} 
if all of its nonempty internal steps are prime.
It may be that an internal step $P_i$ contains a single element
$p_i$, in which case it is a {\em singleton} step.
\index{singleton step}
A singleton step $P_i=\{p_i\}$ is a {\em guts} 
step
\index{guts step}
or a {\em coguts} step
\index{coguts step} 
if it is respectively in the guts or coguts of the 3-separation
$(P_0\cup P_1\cup P_1\udots P_i,
P_{i+1}\cup P_{i+2}\cup P_{i+2}\udots P_l)$ 
of $M$. To simplify notation
we often denote the singleton step $\{p_i\}$ by $p_i$.
We also use $P_i^-$ to denote the set 
$P_0\cup P_1\udots P_{i-1}$ and
$P_i^+$ to denote the set $P_{i+1}\cup P_{i+2}\udots P_l$.

The next lemma follows from the fact that the connectivity function
of a matroid is monotone under minors.

\begin{lemma}
\label{induced-path}
Let $\PP=(P_0,P-1,\ldots,P_l)$ be a path of 
$3$-separations in the matroid
$M$ and let $M'$ be a minor of $M$ on the set $E'$. 
If $\kappa_{M'}(P_0\cap E',P_l\cap E')=2$, then
$\PP'=(P_0\cap E',P_1\cap E',\ldots,P_l\cap E')$ is a path of $3$-separations in $M'$.
\end{lemma}

There are no surprises in the next elementary lemma. It 
essentially says that if we keep the same local connectivity,
then we keep the same guts.


\begin{lemma} 
\label{basic2}
Let $Y$ and $Z$ be disjoint sets of a matroid
$M$ and let $y$ be an element of $Y$. Assume that $Z'\subseteq Z$,
and $\sqcap(Y,Z')=\sqcap(Y,Z)$. Then $y\in\cl(Z')$ if and only
if $y\in\cl(Z)$.
\end{lemma} 

\begin{proof}
One direction is clear. For the other direction, say
that $y\in\cl(Z)$. Assume, for a contradiction, that
$y\notin \cl(Z')$.
We consider two cases. For the
first assume that $y\notin\cl(Y-\{y\})$. Then
$\sqcap(Z\cup\{y\},Y-\{y\})=\sqcap(Z,Y)-1$,
and $\sqcap(Z'\cup\{y\},Y-\{y\})
=(r(Z')+1)+(r(Y)-1)-r(Z'\cup Y)
=\sqcap(Z',Y)$. This contradicts the fact that
$\sqcap(Z'\cup\{y\},Y-\{y\})\leq \sqcap(Z\cup\{y\},Y-\{y\})$.

For the other case, assume that $y\in\cl(Y-\{y\})$.
Then $\sqcap(Z\cup\{y\},Y-\{y\})=\sqcap(Z,Y)$,
and $\sqcap(Z'\cup\{y\},Y-\{y\})
=(r(Z')+1)+r(Y)-r(Z'\cup Y)=\sqcap(Z',Y)+1$.
Again we contradict the fact that
$\sqcap(Z'\cup\{y\},Y-\{y\})\leq \sqcap(Z\cup\{y\},Y-\{y\})$.
\end{proof}

We will also need the following slight 
strengthening of Tutte's Linking
Theorem. 

\begin{corollary}
\label{linking}
Let $(A,Z,B)$ be a partition of the ground set of a matroid $M$
where $\kappa(A,B)=k$. Then there is a partition $\{I,J\}$ of
$Z$ such that $\sqcap(A,I)=\sqcap(B,I)=0$, and
$\lambda_{M/I\ba J}(A,B)=k$.
\end{corollary}

\begin{proof}
By Tutte's Linking Theorem there is a partition $(I,J)$
of $E(M)-(A\cup B)$ such that 
$\lambda_{M/I\ba J}(A,B)=\kappa_M(A,B)$. Assume that
amongst all partitions with this property we have chosen
$I$ to have least cardinality. Let $N=M\ba J$.
Then, for all $z\in I$
there is a $k$-separation $(A',B')$ of $N$ such that
$A\subseteq A'$, $B\subseteq B'$ and $z$ is in the coguts of
$(A',B')$. A routine uncrossing argument shows that
there is an ordering $(z_1,z_2,\ldots,z_l)$ of $I$,
with the property that 
$(A\cup\{z_1,\ldots,z_i\},\{z_{i+1},\ldots,z_l\}\cup B)$
is a $k$-separation of $N$ with $z_i$ in the 
coguts for all $i\in\{1,2,\ldots,l\}$. But then 
$r(A\cup I)=r(A)+|I|$ and $r(B\cup I)=r(B)+|I|$ giving the
corollary.
\end{proof}


Let ${\mathcal U}(q)$ denote the class of matroids with 
no $U_{2,q+2}$-minor and ${\mathcal U}^*(q)$ denote the class
of matroids with no $U_{q,q+2}$-minor.
We use the following result of Kung \cite{ku93}.

\begin{lemma}
\label{kung}
If $M$ is a  simple rank-$r$ matroid in ${\mathcal U}(q)$,
then $M$ has at most $(q^r-1)/(q-1)$
elements.
\end{lemma}

An easy corollary of Lemma~\ref{kung} is

\begin{corollary}
\label{kung2}
If $M$ is a simple matroid in ${\mathcal U}(q)$, 
and $A\subseteq E(M)$
has $\lambda(A)=r$, then there are at most $(q^r-1) /(q-1)$
elements in $\cl(A)-A$.
\end{corollary}

\section{Strands and Paths of Clonal Pairs}

Inn this section we develop  some straightforward properties of
paths of $3$-separations. In doing this we typically 
do not need strong
connectivity assumptions about the underlying matroid.
Recall that in
Chapter~\ref{flower-chapter} flowers were defined for
connected, but not necessarily 3-connected matroids. 

\subsection*{Strands}

Let $(A,B)$ be a partition of the ground set of the matroid
$M$. A $B$-{\em strand} 
\index{$B$-strand}
\index{strand}
is a  minimal subset $X$ of $B$
such that $\sqcap(A,X)=1$. Say $a\in A$.
Note that, if $X$ is a $B$-strand, 
and $a\in\cl(X)$, then $X\cup\{a\}$ is a circuit.
Note also, that if $Y$ is a subset of $B$ 
with $\sqcap(Y,A)=1$ and $a\in\cl(Y)$, then there
is a strand $X\subseteq Y$ such that $a\in\cl(X)$.

\begin{lemma}
\label{strand1}
Let $(A,x,B)$ be a path of $3$-separations in the matroid
$M$. If there exists a $B$ strand $X$ such that $x\in\cl(X)$,
then $x$ is fixed in $M$.
\end{lemma}

\begin{proof}
Assume that $X$ is a $B$ strand with $x\in\cl(X)$.
Then $x\in\cl(B)$. If $x\notin\cl(A)$, then
$(A,B\cup\{x\})$ is a 2-separation in $M$ contradicting the 
definition of path of $3$-separations. Therefore
$x\in\cl(A)$. 

Let $M'$ be a matroid obtained by cloning $x$ by $x'$.
Then $x'\in\cl_{M'}(A)$, so
$\sqcap(A\cup\{x,x'\},X)=\sqcap(A\cup\{x\},X)=1$.
But $\{x,x'\}\subseteq\cl_{M'}(X)$, so that $x$ and $x'$ are
parallel in $M'$. Hence $x$ is fixed in $M$.
\end{proof}

Let $(A,x,B)$ be a path of 3-separations of 
$M$, where $x$ is a guts singleton.
Then $x$ is {\em fixed from the right}
\index{fixed from the right} 
if there is a $B$-strand
$X$ such that $x\in\cl(X)$ and $x$ is {\em fixed from the left}
\index{fixed from the left}
if there is an $A$-strand $X$ such that $x\in\cl(X)$.

The next lemma follows from Theorem~6.2 and Lemma~6.3 of 
\cite{gegewh06}. See also \cite{be06}.

\begin{lemma}
\label{free-in-guts}
Let $(A,B)$ be an exact $3$-separation of the matroid $M$.
Then there exists an extension $M'$ of $M$
by the element $x'$ such that $x'\in\cl_{M'}(A)\cap \cl_{M'}(B)$
and $x'$ is not fixed in $M'$.
\end{lemma}

\begin{lemma}
\label{strand2}
Let $(A,x,B)$ be a path of $3$-separation, 
where $x$ is a guts singleton. Then
$x$ is fixed in $M$ if and only if $x$ is fixed from the left or
fixed from the right.
\end{lemma}

\begin{proof}
If $x$ is fixed from the left or right, then $x$ is fixed in $M$
by Lemma~\ref{strand1}. Consider the converse. Assume that
$x$ is fixed in $M$.
By Lemma~\ref{free-in-guts} it is possible to extend $M$ by an element $x'$
to obtain a matroid $M'$ with a path of 3-separations
$(A,x,x',B)$ such that
$x'$ is a guts singleton and $x'$ is not fixed in $M'$. 
As $x$ is fixed in
$M'$, there is a circuit $C$ of $M$ containing $x$ such that
$x'\notin\cl_{M'}(C)$. Consider $M'/x'$. As 
$x'\notin\cl_{M'}(C)$, we see that $C$ is a circuit of
$M'/x'$. Moreover, $x$ is in the guts of the 2-separation
$(A\cup\{x\},B)$ of
$M'/x'$. But now $(A,B)$ is a separation
of $M'/x',x$ and $C-\{x\}$ is a circuit of this matroid.
Thus $C-\{x\}\subseteq A$ or $C-\{x\}\subseteq B$. Assume the 
latter holds. As $x\in \cl(C-\{x\})$, we see that
$\sqcap(A\cup\{x\},C-\{x\})\geq 1$.
But $x'\notin \cl(C-\{x\})$, so, by Lemma~\ref{basic2},
$\sqcap(A\cup\{x\},C-\{x\})<2$. Hence $C-\{x\}$ is a $B$-strand,
and $x$ is fixed from the right.
\end{proof}

\begin{lemma}
\label{strand3}
Let $(A,x,B)$ be a path of $3$-separations of the matroid $M$, 
where $x$ is a guts singleton. Assume that
$B$ has a partition $\{B_1,B_2,\ldots,B_l\}$ such that
\begin{itemize}
\item[(i)] $(A\cup\{x\},B_1,B_2,\ldots,B_l)$ is a swirl-like or 
spike-like flower with at least
three petals, and
\item[(ii)] $B_1$ is a clonal pair.
\end{itemize}
Then $x$ is fixed from the right in $(A,x,B)$ if and only 
if either $x\in\cl(B_1)$ or $x\in\cl(B_l)$.
\end{lemma}

\begin{proof}
Say $x\in\cl(B_1)$ or $x\in\cl(B_l)$. Then, as
$\sqcap(B_1,A)=\sqcap(B_l,A)=1$, we see that $x$ is fixed from the right.

Assume that $x$ is fixed from the right. Then there is a $B$-strand 
$X$ whose closure contains $x$ so that
$X\cup\{x\}$ is a circuit. For the first case, assume that
$X\cap B_1\neq \emptyset$. Then, as $X\cup\{x\}$ is a circuit,
and $B_1$ is a clonal pair, we see that $B_1\subseteq \cl(X)$.
Now
$$1=\sqcap(B_1,A)\leq \sqcap(\cl(X),A)=1,$$
so that $\sqcap(B_1,A)= \sqcap(\cl(X),A)$. By Lemma~\ref{basic2},
$x\in\cl(B_1)$, as required. For the other case assume that
$X\cap B_1=\emptyset$. Then
$$1=\sqcap(B_2\cup B_3\udots B_l,A)=\sqcap(X,A)=\sqcap(B_l,A)=1,$$
and again by Lemma~\ref{basic2}, it follows that $x\in\cl(B_l)$.
\end{proof}

\begin{lemma}
\label{strand4}
Let $(A,x,B)$ be a path of $3$-separations of the matroid $M$, 
where $x$ is a guts
singleton. Assume that $x$ is not fixed from the right. Say
$I\subseteq B$, and $\sqcap(A,I)=0$. Then, in $M/I$, the element
$x$ is not parallel to any element of $A$.
\end{lemma}

\begin{proof}
Assume that $y\in B$ and $\{x,y\}$ is a parallel pair in $M/I$. Then
there is a circuit $C$ of $M$ such that 
$\{x,y\}\subseteq C\subseteq \{x,y\}\cup I$.
Now $1\leq \sqcap(C-\{x\},A)\leq \sqcap(I,A)+1=1$, so
that $\sqcap(C-\{x\},A)=1$. Thus $C-\{x\}$ is a $B$-strand
and $x\in\cl(C-\{x\})$. 
This contradicts the fact that $x$ is not fixed from
the right.
\end{proof}

\begin{lemma}
\label{strand5}
Let $\PP$ be a path of $3$-separations with a set
$\{x_1,\ldots,x_l\}$ of singleton guts elements, none of which
are fixed from the right. Then $M$ has a $U_{2,l}$-minor
on the set $\{x_1,\ldots,x_l\}$.
\end{lemma}

\begin{proof}
By an appropriate concatenation we may assume that
$$\PP=(P_0,x_1,P_1,x_2,P_2,\ldots,P_{l-1},x_l,P_l).$$
Note that $\PP$ need not be a solid path in that it may
have some empty steps.
By Corollary~\ref{linking},
there is a partition $\{I,J\}$ of $P_{l-1}$ such that
$\sqcap(P_0\cup\{x_1\}\cup P_1\udots P_{l-2}\cup\{x_{l-1}\},I)=
\sqcap(I,\{x_l\}\cup P_l)=0$, and 
$\lambda_{M/I\ba J}(P_0\cup\{x_1\}\cup\udots P_{l-2}
\cup\{x_{l-1}\},\{x_l\}\cup P_l)=2.$ Evidently
$(P_0,x_1,P_1,\ldots,P_{l-2},x_{l-1},x_l,P_l)$ is a path of 
$3$-separations in $M/I\ba J$. By Lemma~\ref{strand4},
$x_{l-1}$ is not parallel to $x_l$. Moreover,
each element of $\{x_1,x_2,\ldots,x_l\}$
is a guts singletons in this path and is not fixed from the right.
Repeating this process we obtain a minor $M'$ with a 
path of $3$-separations $P'=(P_0,x_1,x_2,\ldots,x_l,P_l)$, where
the members of $\{x_1,x_2,\ldots,x_l\}$ are guts singletons and 
$x_i$ is not parallel to $x_j$ for all distinct 
$i,j\in\{1,2,\ldots,l\}$.
Evidently $M'|\{x_1,x_2,\ldots,x_l\}\cong U_{2,l}$.
\end{proof}

\subsection*{Displayed Flowers}

Let $\PP=(P_0,P_1,\ldots,P_l)$ be a path of 3-separations 
of the matroid $M$ and
let $\QQ$ be a flower of $M$. Then $\QQ$ is {\em displayed} 
by $\PP$ if 
each petal of $\QQ$ is a union of steps
of $\PP$. 
\index{flower displayed by a path}
The flowers displayed by $\PP$ are 
partially ordered by refinement.
The flower $\QQ$ is a {\em maximal} displayed flower of $\PP$
\index{maximal displayed flower} 
if it is maximal in this
partial order. Note that a maximal displayed flower 
need not be a maximal
flower of $M$ and flowers that are incomparable in this partial order
may be comparable in the usual partial order of flowers.

\begin{lemma}
\label{display1}
Let $\PP=(P_0,P_1,\ldots,P_l)$ be a path of $3$-separations of the 
connected matroid
$M$ with the property that $\lambda(P_i)\geq 2$ for all
$i\in\{0,1,\ldots,l\}$.
Let $\QQ=(Q_1,Q_2,\ldots,Q_r)$ be a maximal
displayed
flower of $\PP$. Then there exist $i,j$ such that
$0\leq i\leq j<l$ , and 
$\{Q_1,\ldots,Q_r\}=\{P_i^-,P_i,P_{i+1},\ldots,P_j,P_j^+\}$.
\end{lemma}

\begin{proof}
Certainly $\QQ$ has petals $Q_\alpha$
and $Q_\beta$ containing $P_0$ and $P_l$ respectively.
We first show that any other petal of $\QQ$ contains exactly
one step of $\PP$. Assume otherwise.
Say that $Q$ is a petal of $\QQ$ other than $Q_\alpha$
or $Q_\beta$ that contains more than
one step of $\PP$. Assume that $P_i,P_j\subseteq Q$, where $i<j$.
Then $(P_i^-\cup P_i,P_i^+)$ is not displayed by $\QQ$
and a routine uncrossing argument shows that there is a flower
displayed by $\PP$ that refines $\QQ$ and displays this 3-separation,
contradicting the fact that $\QQ$ is a maximal displayed 
flower of $\PP$.

Consider $Q_\alpha$, the petal containing $P_0$. Let $i$ be the
greatest integer such that $P_i\subseteq Q_\alpha$ We now 
show that $Q_\alpha=P_0\cup P_1\udots P_i$. Assume that this is
not the case. Then there is a $j$ with
$0<j<i$ such that 
$P_j\not\subseteq Q_\alpha$. In this case an uncrossing argument using
the 3-separation $(P_j^-\cup P_j,P_j^+)$ again shows that
$\QQ$ is not a maximal displayed flower of $\PP$. 
Hence $Q_\alpha$ and $Q_\beta$ respectively contain
an initial and a terminal sequence of steps of $\PP$ and the lemma
follows.
\end{proof}

A straightforward uncrossing argument also proves

\begin{lemma}
\label{display2}
Let $\PP=(P_0,P_1,\ldots,P_l)$ be a path of $3$-separations of the 
connected matroid
$M$ with the property that $\lambda(P_i)\geq 2$ for all
$i\in\{0,1,\ldots,l\}$.
If $P$ is a step of $\PP$, then $P$ is 
a petal of at most one maximal displayed flower.
\end{lemma}

\subsection*{Special Paths of $3$-separations}

We say that a flower in a matroid is a {\em spiral} 
\index{spiral}
if it is
either swirl-like or spike-like.
A path of $3$-separations 
$(P_0,P_1,\ldots,P_l)$ is {\em special}
if 
\begin{itemize}
\item[(i)] $(P_i^-,P_i,P_i^+)$ is a spiral
for all
$i\in\{1,2,\ldots,l-1\}$, and
\item[(ii)] $(P_0,P_1,\ldots,P_l)$ displays no $4$-petal flowers.
\end{itemize}
\index{special path of $3$-separations}

We begin by characterising special paths of length 4. We first
note an elementary operation that can be performed on special paths
of 3-separations.

\begin{lemma}
\label{special1}
Let $(A,B,C,D)$ is a special path of $3$-separations. Then 
so too is $(B,A,C,D)$. 
\end{lemma}

\begin{proof}
As $(A,B,C\cup D)$ is a spiral, $\lambda(B)=2$, so that
$(B,A,C,D)$ is a path of $3$-separations. The flowers displayed
by $(B,A,C,D)$ are the same as the flowers displayed by
$(A,B,C,D)$, so that the path displays no 4-petal flowers. 
To show that 
the path $(B,A,C,D)$ is special, 
we need to show that (i) $(B\cup A,C,D)$
and (ii) $(B,A,C\cup D)$ are spirals. As $(A,B,C,D)$
is special $(A\cup B,C,D)$ is a spiral, so that (i) holds.
We have already observed that $(A,B,C\cup D)$ is a spiral and
this spiral is equivalent or $(B,A,C\cup D)$ so that (ii) holds.
\end{proof}

We say that the special path $(B,A,C,D)$ 
is obtained from $(A,B,C,D)$ by
{\em switching}. 
\index{switching}
Evidently $(D,C,B,A)$ is also a special path
of 3-separations
and we say that it is obtained from $(A,B,C,D)$ by {\em reversal}.
\index{reversal}
The special path $(A,B,C,D)$ has
\begin{description}
\item[Type I] if $\sqcap(A,C)=\sqcap(B,C)=1$ and $\sqcap(B,D)=0$,
\item[Type II] if $\sqcap(B,C)=1$ and $\sqcap(A,C)=\sqcap(B,D)=0$, and
\item[Type III]  if $\sqcap(A,C)=\sqcap(B,C)=\sqcap(B,D)=0$.
\end{description}
It is easily seen that no sequence of 
switches and reversals can convert
a special path of one type into another type. 

\begin{lemma}
\label{special2}
If $(A,B,C,D)$ is a special path of $3$-separations, 
then some sequence of
switches and reversals converts $(A,B,C,D)$ into a path of 
Type I, II, or III.
\end{lemma}

\begin{proof}
As $(A,B,C\cup D)$ is a spiral in $M$, we have
$\sqcap(A,B)=1$ and similarly $\sqcap(C,D)=1$. Moreover
$\sqcap(A,C)\leq \sqcap(A,C\cup D)=1$. Thus
$\sqcap(A,C)$, $\sqcap(A,D)$, $\sqcap(B,C)$ and 
$\sqcap(C,D)$ are all at most 1.

Say that $\sqcap(A,C)=\sqcap(B,C)=1$. 
Assume that $\sqcap(B,D)=1$. 
Then
\begin{align*}
\lambda(A\cap C)&=r(A\cup C)+r(B\cup D)-r(M)\\
&=(r(A)+r(C)-1)+(r(B)+r(D)-1)-r(M)\\
&=r(A\cup B)+r(C\cup D)-r(M)\\
&=2.
\end{align*}
It follows from this that $(A,B,C,D)$ is a flower and it 
is easily checked that this flower is spike-like,
contradicting the definition of special path.
Thus $\sqcap(B,D)=0$, and $(A,B,C,D)$ has 
Type~I.

Assume that $(A,B,C,D)$ cannot be converted into a path 
of  Type~I or Type~III by a sequence of switchings and reversals. 
Then we may assume,
up to switchings and reversals, that $\sqcap(A,C)=0$ and
$\sqcap(B,C)=1$. If $\sqcap(B,D)=1$, 
then $(A,B,C,D)$ converts into a path of  
Type~I. Assume that $\sqcap(B,D)=0$. If $\sqcap(A,D)=1$, then it is
easily checked that $(A,B,C,D)$ is a swirl-like flower.
Hence $\sqcap(A,D)=0$ and $(A,B,C,D)$ has Type~II.
\end{proof}

We omit the straightforward rank calculation that proves the next
lemma.

\begin{lemma}
\label{special3}
If $(A,B,C,D)$ is a special path of $3$-separations 
of Type I, II or III
in $M$, then $(B,A,D,C)$ is a special path of 
$3$-separations of Type I, II or
III respectively in $M^*$.
\end{lemma}

\subsection*{Paths of Clonal Pairs}

Let $\{p,q\}$ be a clonal pair in a matroid $M$. Then
$\{p,q\}$ is $M$-{\em strong} if $\lambda_M(\{p,q\})=2$.
\index{$M$-strong clonal pair}
Let $\PP=(P_0,\ldots,P_l)$ be  a path of 
$3$-separations of the connected matroid $M$. Then $\PP$
is a {\em path of clonal pairs} if 
\begin{itemize}
\item[(i)] $\PP$ is maximal;
\item[(ii)] $P_i$ is an $M$-strong clonal pair for
all $i\in\{1,2,\ldots,l-1\}$.
\item[(iii)] $(P_i^-,P_i,P_i^+)$ is a flower for
all $i\in\{1,2,\ldots,l-1\}$.
\end{itemize}
\index{path of clonal pairs}

Note that (iii) follows from (ii) except for the fact that there is
nothing in the definition of a path of 3-separations to prevent the
possibility of $M$ having a 2-separating set that 
crosses $P_0\cup P_l$. Such a 2-separating set causes no
difficulties structurally, but, given our definition of
flower it does prevent $\PP$ from having displayed flowers.
This triviality having been dealt with we move on.

Let $\PP=(P_0,P_1,\ldots,P_l)$ be a path of clonal pairs.
Say $i\in\{1,2,\ldots,l-1\}$ and $P_i=\{p_i,q_i\}$.
If the flower
$(P_i^-,P_i,P_i^+)$ is
a paddle or copaddle, then
$(P_0,P_1,\ldots,P_{i-1},p_i,q_i,P_{i+1},\ldots,P_l)$ 
is a path of 3-separations contradicting
the maximality of $\PP$. Thus  $(P_i^-,P_i,P_i^+)$ is a spiral.
The next corollary follows from  this observation and 
Lemmas~\ref{special2} and \ref{special3}.

\begin{lemma}
\label{clonal1}
Let $\PP=(P_0,P_1,\ldots,P_l)$ be a path of clonal pairs. Then every flower displayed by $\PP$ is a spiral. Moreover, there are integers
$0\leq i_1\leq\cdots\leq i_m<l$ such that $(Q_1,Q_2,\ldots,Q_r)$ 
is a maximal displayed
flower of $\PP$ if and only if 
$$\{Q_1,Q_2,\ldots,Q_r\}=\{P_{i_j+1}^-,P_{i_{j}+1},
 \ldots,P_{i_{j+1}},P_{i_{j+1}}^+\}$$
for some $j\in\{1,2,\ldots,m\}$.
\end{lemma}

Via Lemma~\ref{clonal1}, the maximal displayed flowers of a path of clonal
pairs can be canonically associated with a partition of 
$(P_0,P_1,\ldots,P_{l-1})$ into consecutive sets of steps. 
We call this partition the
{\em flower partition} of $\PP$ in $M$. We omit the routine proof
of the next lemma.

\begin{lemma}
\label{clonal2}
Let $\PP=(P_0,P_1,\ldots,P_l)$ be a 
path of clonal pairs of the matroid
$M$. If $t\in\{1,2,\ldots,l-1\}$ and 
$P_t=\{p_t,q_t\}$ is a petal of a displayed 
flower of $\PP$ with at least four petals, then the following hold.
\begin{itemize}
\item[(i)] $(P_0,P_1,\ldots,P_{t-1},P_{t+1},\ldots,P_l)$ 
is a path of clonal pairs in $M\ba p_t/q_t$.
\item[(ii)] If 
$({\mathcal P}_1,{\mathcal P}_2,\ldots,{\mathcal P}_s)$ is the flower
partition of $\PP$ in $M$, where $P_t\in{\mathcal P}_j$, then 
$({\mathcal P}_1,{\mathcal P}_2,\ldots,
{\mathcal P}_{j-1},{\mathcal P}_j-\{P_t\},
{\mathcal P}_{j+1},\ldots,{\mathcal P}_l)$ is 
the flower partition of the path of clonal pairs 
$(P_0,P_1,\ldots,P_{t-1},P_{t+1},\ldots,P_l)$ in $M\ba p_t/q_t$.
\end{itemize}
\end{lemma}

The next corollary follows from Lemma~\ref{clonal2} and an
elementary calculation.

\begin{corollary}
\label{clonal3} 
Let $s\geq 1$  and $k\geq 5$ be integers and
let  $\PP=(P_0,P_1,\ldots,P_l)$ be a 
path of clonal pairs of the connected matroid $M$. Assume that 
$\PP$ does not display
any spiral with $k$ petals. If 
$l> (k-3)(s-1)$, then $M$ has a  minor with a 
special path of clonal pairs of length $s$.
\end{corollary}

\subsection*{Special Paths of Clonal Pairs}

At the risk of overkill we point out that a {\em
special path of clonal pairs} is a special path of 
3-separations that is also a path of clonal pairs.
We omit the easy proof of the next lemma.

\begin{lemma}
\label{special10}
Let $\PP=(P_0,P_1,\ldots,P_l)$ be a special path 
of clonal pairs of the connected matroid
$M$. If $P_i=\{p_i,q_i\}$ is an internal step of $\PP$, then
$(P_0,\ldots,P_{i-1},\{q_i\}\cup P_i^+)$ and
$(P_i^-\cup\{q_i\},P_{i+1},\ldots,P_l)$ are both special 
paths of clonal pairs in both $M\ba p_i$ and $M/p_i$.
\end{lemma}

\begin{lemma}
\label{special11}
Let $(P_0,\{p_1,q_1\},\{p_2,q_2\},P_3)$ be a special 
path of clonal pairs
in the connected matroid $M$. 
Then for some $N\in\{M,M^*\}$, either 
\begin{itemize}
\item[(i)] $q_1$ is not fixed from the right in the 
path $(P_0,q_1,\{p_2,q_2\},P_3)$
of $N/p_1$, or
\item[(ii)] $q_2$ is not fixed from the left in the 
path $(P_0,\{p_1,q_1\},q_2,P_3)$
of $N/p_2$.
\end{itemize}
\end{lemma}

\begin{proof}
Consider the path $(P_0,\{p_1,q_1\},\{p_2,q_2\},P_3)$. 
By Lemma~\ref{special2} a sequence
of switches and reversals converts it into a path or 
Type I, II or III. If the lemma
holds for a given path, then it certainly holds for 
the reversal of that path.
It is now easily seen that we lose no generality in assuming that
there is a path $(A,B,C,D)$ of Type I, II or III 
such that $\{A,B\}=\{P_0,\{p_1,q_1\}\}$
and $\{C,D\}=\{\{p_2,q_2\},P_3\}$. By Lemma~\ref{special3}, 
we may further assume that
$D=\{p_2,q_2\}$. Then 
$\sqcap_M(P_0,\{p_2,q_2\})=\sqcap_M(\{p_1,q_1\},\{p_2,q_2\})=0$.
Hence $q_2\notin \cl_{M/p_2}(P_0)$ and 
$q_2\notin\cl_{M/p_2}(\{p_1,q_1\})$.
Now, by Lemma~\ref{strand3}, $q_2$ is not 
fixed from the left in 
$(P_0,\{p_1,q_1\},q_2,P_3)$.
\end{proof}

\begin{lemma}
\label{special12}
Let $s$ be a positive integer and let 
$\PP=(P_0,P_1,\ldots,P_l)$ be a special 
path of clonal pairs of the  connected matroid $M$. If $l>8s$,
then $M$ has a $U_{2,s}$- or a $U_{s-2,s}$-minor.
\end{lemma}

\begin{proof}
Assume that $l>8s$. 
By Corollary~\ref{clonal3}, we may assume that, 
up to reversal and taking duals,
that $\PP$ has a subsequence 
$(P_{i_1}=\{p_{i_1},q_{i_1}\},P_{i_2}=\{p_{i_2},q_{i_2}\},
\ldots,P_{i_s}=\{p_{i_s},q_{i_s}\})$ 
of internal
steps with the following property: 
for $j\in\{1,2,\ldots,s\}$, the element $q_{i_j}$
is not fixed from the right in the path
$(P_0,P_1,\ldots,P_{i_j-1},q_{i_j},P_{i_j+1},\ldots,P_l)$ 
of $3$-separations
of $M/p_{i_j}$. By Lemma~\ref{special10}, 
$(P_0\cup P_1\udots P_{i_j-1}\cup\{q_{i_j}\},P_{i_j+1},\ldots,P_l)$
is a special path of clonal pairs in this matroid. Define 
$\QQ=(Q_0,Q_1,\ldots,Q_l)$ as follows. 
If $r\in\{i_1,i_2,\ldots,i_s\}$, 
then $Q_r=\{q_r\}$
and otherwise $Q_r=P_r$. This is clearly a path of 3-separations in 
$M/p_{i_1},p_{i_2},\ldots,p_{i_s}$ and the elements 
$\{q_{i_1},q_{i_2}\ldots,q_{i_s}\}$ are
guts singletons that are not fixed from the right. 
By Lemma~\ref{strand5}, $M$ has a $U_{2,s}$-minor.
\end{proof}

\begin{corollary}
\label{special13}
There is a function $f_{\ref{special13}}(k,q)$ such that, if 
$n\geq f_{\ref{special13}}(k,q)$ and 
$\PP=(P_0,P_1,\ldots,P_n)$ is a path
of clonal pairs in a connected matroid $M$ that does not display
any swirl-like flower of order $k$,
then $M\notin \eq$. 
\end{corollary}

\begin{proof}
Let $s=\max\{k,q\}$ and let 
$f_{\ref{special13}}(k,q)=((s+2)-3)8s$.
Certainly $\PP$ does not display any
swirl-like flower with $s+2$ petals. If $\PP$
displays a spike-like flower with $s+2$ petals,
then all but at most two of the petals of this flower are
clonal pairs so that $M$ has a $\Lambda_s$- and hence
a $\Lambda_q$-minor. Therefore we may assume that
$\PP$ does not display any spiral with $s+2$ petals.

By Corollary~\ref{clonal3}, $M$ has a minor with a special path of 
clonal pairs of length $8s+1$. By Lemma~\ref{special12},
$M$ has a $U_{2,q+2}$- or a $U_{q,q+2}$-minor and we conclude that
$M\notin \eq$.
\end{proof}

\section{Feral Elements in Paths}

To prove Theorem~\ref{main1} we would like to reduce a long path of
$3$-separations in a $k$-skeleton to one 
in which the internal steps are
either singletons or clonal pairs. 
But there are obstacles to doing this
caused by the presence of feral elements.

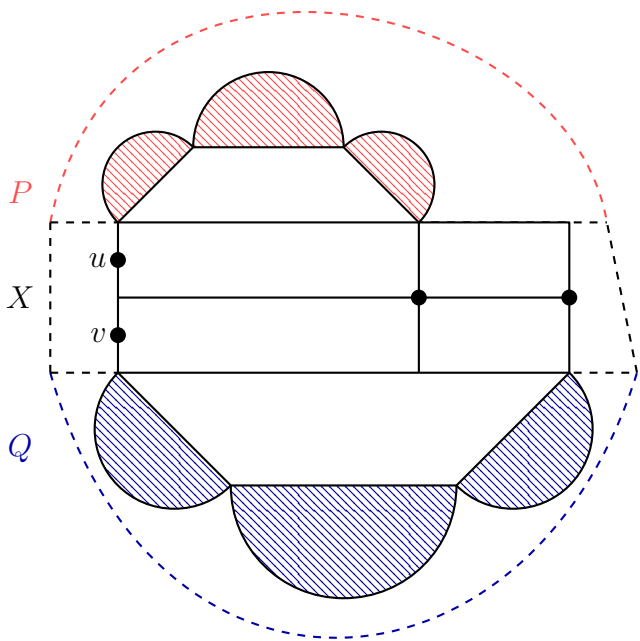
\begin{figure}
\begin{tikzpicture}[thick,line join=round]
	\coordinate (o) at (0,0);
	\coordinate (p) at (0,1);
	\coordinate (q) at (0,2);
	\coordinate (r) at (4,2);
	\coordinate (s) at (6,2);
	\coordinate (t) at (4,1);
	\coordinate[label=left:$u$] (u) at ($(p)!0.5!(q)$);
	\coordinate[label=left:$v$] (v) at ($(o)!0.5!(p)$);
	\coordinate (w) at (6,1);
	\coordinate (x) at (4,0);
	\coordinate (y) at (6,0);
	\coordinate (a) at ($(q) + (45:1.414)$);
	\coordinate (b) at ($(r) + (135:1.414)$);
	\coordinate (c) at (-45:2.121);
	\coordinate (d) at ($(y) + (-135:2.121)$);
	\coordinate (e1) at ($(q) + (-0.9,0)$);
	\coordinate (e2) at ($(s) + (0.5,0)$);
	\draw[dashed,rlabels] (e1) .. controls ($(a) + (-1,3.4)$) and ($(b) + (3,2)$) .. (e2);
	\coordinate (f1) at (-0.9,0);
	\coordinate (f2) at ($(y) + (0.9,0)$);
	\draw[dashed,labels] (f1) .. controls ($(c) + (-1,-3.4)$) and ($(d) + (1,-3)$) .. (f2);
	\draw[dashed] (f1) -- (e1);
	\draw[dashed] (f2) -- (e2);
	\node (X) at ($(p) + (-1.3,0)$) {$X$};
	\node at ($(X) + (0,1.4)$) {\textcolor{rlabels}{$P$}};
	\node at ($(X) + (0,-2)$) {\textcolor{labels}{$Q$}};
	\draw[dashed] (f1) -- (f2);
	\draw[dashed] (e1) -- (e2);
	\node[pattern color=rlines,draw,circle through=(q),pattern=north west lines] at ($(q)!0.5!(a)$) {};
	\node[pattern color=rlines,draw,circle through=(a),pattern=north west lines] at ($(b)!0.5!(a)$) {};
	\node[pattern color=rlines,draw,circle through=(b),pattern=north west lines] at ($(b)!0.5!(r)$) {};
	\node[pattern color=dblines,draw,circle through=(o),pattern=north west lines] at ($(o)!0.5!(c)$) {};
	\node[pattern color=dblines,draw,circle through=(c),pattern=north west lines] at ($(c)!0.5!(d)$) {};
	\node[pattern color=dblines,draw,circle through=(d),pattern=north west lines] at ($(d)!0.5!(y)$) {};
	\filldraw[fill=white] (o) -- (c) -- (d) -- (y) -- (s) -- (r) -- (b) -- (a) -- (q) -- cycle;
	\draw (p) -- (w);
	\draw (o) -- (y);
	\draw (x) -- (r) -- (q);
	\foreach \pt in {u,v,t,w} \fill[black] (\pt) circle (3pt);
\end{tikzpicture}
\caption{A feral pack in the path $(P,X,Q)$}\label{feral-pack}
\end{figure}

Let $X$ be a prime step in a path of $3$-separations  in the
$k$-coherent matroid $M$. Then
$X$ is a {\em feral pack}
\index{feral pack} 
if no element of
$X$ is in a triangle or triad, $X$ contains a single clonal pair
$\{u,v\}$, and $X-\{u,v\}$ consists of feral elements of $M$.
Figure~\ref{feral-pack} illustrates a feral pack. Note that $X-\{u,v\}$
is a bogan couple. We know of no examples of feral packs where
$X-\{u,v\}$ is neither a bogan nor a cobogan couple. 

Recall that if $Z\subseteq E(M)$, then $\coh(Z)$ denotes the 
set $Z-\fcl(E(M)-Z)$.

\begin{lemma}
\label{pack0}
Let $Z$ be a $3$-separating set of the $k$-coherent matroid
$M$ such that no member of $Z$ is in a triangle or a triad
and $|Z|\geq 4$. Then the following hold.
\begin{itemize}
\item[(i)] There is an element $z\in Z$ such that either $M\ba z$ or $M/z$
is $k$-coherent.
\item[(ii)] If $Z$ contains a clonal pair $\{u,v\}$, then there is
an element $z\in Z-\{u,v\}$ such that either $M\ba z$
or $M/z$ is $k$-coherent.
\end{itemize}
\end{lemma}

\begin{proof}
As no member of $Z$ is in a triangle or a triad, 
$Z$ is non-sequential. Thus a $3$-separation 
equivalent to $Z$ is displayed in 
a 3-tree for $M$. It follows that $Z$ contains a peripheral set $P$.
Assume that  $P$ is non-sequential, then by Lemma~\ref{wheel3.5}(ii), if 
$z\in\coh(P)$, either $M\ba z$ or $M/z$ is $k$-coherent.
As $|\coh(P)|\geq 4$ there is an element 
$z\in\coh(P)-\{u,v\}$ and the lemma
holds. Assume that $P$ is sequential. 
Then, as no element of $Z$ is in a triangle
or a triad, $|P|=2$. By Corollary~\ref{2-elt-win}, 
if $p\in P$, then either $M\ba p$ or $M/p$ 
is $k$-coherent and the lemma holds unless $P=\{u,v\}$. If $P=\{u,v\}$,
then it is clear that $Z$ contains another peripheral 
set $P'$ and the lemma
holds by considering this peripheral set.
\end{proof}

It follows immdiately from Lemma~\ref{pack0} that if
$X$ is a feral pack of the $k$-coherent matroid $M$ and
$Z\subseteq X$ is $3$-separating, then $|Z|\leq 2$.
Because of this the next lemma applies to feral elements of
feral packs.

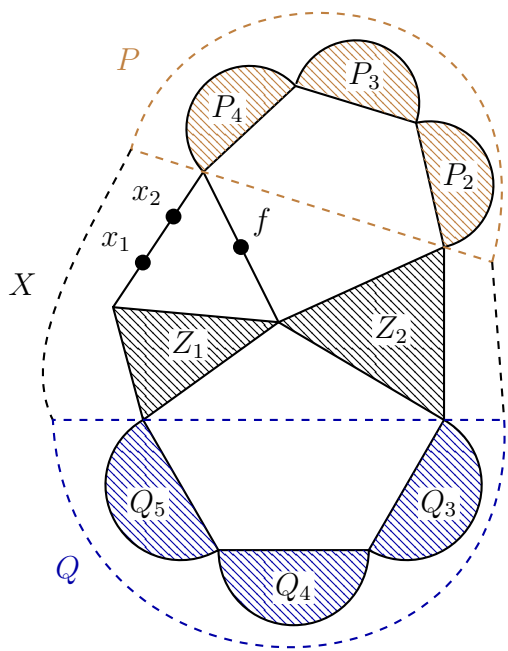
\begin{figure}
\begin{tikzpicture}[thick,line join=round]
	\coordinate (o) at (0,0);
	\coordinate (a) at (-1.8,-1.3);
	\coordinate (b) at (-2.2,0.2);
	\coordinate (c) at (-1,2);
	\coordinate (d) at (2.2,1);
	\coordinate (e) at (2.2,-1.3);
	\coordinate[label=45:$f$] (f) at ($(o)!0.5!(c)$);
	\coordinate (g) at ($(a)+(-60:2)$);
	\coordinate (h) at ($(g)+(2,0)$);
	\coordinate (i) at ($(c)+(43:1.676)$);
	\coordinate (j) at ($(i)+(-17:1.676)$);
	\coordinate[label=135:$x_1$] (x1) at ($(b)!0.33!(c)$);
	\coordinate[label=135:$x_2$] (x2) at ($(b)!0.67!(c)$);
	\draw (o) -- (a) -- (b) -- (o) -- (c) -- (b);
	\draw (o) -- (d) -- (e) -- cycle;
	\draw (a) -- (g) -- (h) -- (e);
	\draw (c) -- (i) -- (j) -- (d);
	\fill[pattern=north west lines] (o) -- (a) -- (b) -- (o) -- (d) -- (e) -- (o);
	\filldraw[pattern color=dblines,pattern=north west lines] let \p1 = ($(g) - (a)$), 
		\n1 = {veclen(\x1,\y1)} in
			(a) arc (90+30:180+120:\n1/2)
			(h) arc (0:-180:\n1/2)
			(e) arc (60:-180+60:\n1/2);
	\filldraw[pattern color=olines,pattern=north west lines] let \p2 = ($(i) - (c)$),
		\n2 = {veclen(\x2,\y2)} in
			(i) arc (43:180+43:\n2/2);
	\filldraw[pattern color=olines,pattern=north west lines] let \p3 = ($(i) - (j)$),
		\n3 = {veclen(\x3,\y3)} in
			(j) arc (43+43-90-14:43+90+43-14:\n3/2);
	\filldraw[pattern color=olines,pattern=north west lines] let \p7 = ($(j) - (d)$),
		\n7 = {veclen(\x7,\y7)} in
			(d) arc (17+43+43-180:17+43+43:\n7/2);
	\foreach \pt in {f,x1,x2} \fill[black] (\pt) circle (3pt);
	\node at (1.5,-0.1) [rectangle,fill=white,draw=white,inner sep=0.5pt] {$Z_2$};
	\node at (-1.2,-0.3) [rectangle,fill=white,draw=white,inner sep=0.5pt] {$Z_1$};
	\node at ($(a)!0.5!(g) + (-150:0.5)$) [rectangle,fill=white,draw=white,inner sep=0.5pt] {$Q_5$};
	\node at ($(h)!0.5!(g) + (-90:0.5)$) [rectangle,fill=white,draw=white,inner sep=0.5pt] {$Q_4$};
	\node at ($(h)!0.5!(e) + (-30:0.5)$) [rectangle,fill=white,draw=white,inner sep=0.5pt] {$Q_3$};
	\node at ($(c)!0.5!(i) + (180-43:0.4)$) [rectangle,fill=white,draw=white,inner sep=0.5pt] {$P_4$};
	\node at ($(j)!0.5!(i) + (90-17:0.4)$) [rectangle,fill=white,draw=white,inner sep=0.5pt] {$P_3$};
	\node at ($(j)!0.5!(d) + (17:0.4)$) [rectangle,fill=white,draw=white,inner sep=0.5pt] {$P_2$};
	\coordinate (t1) at ($(c)!1.2!(d)$);
	\coordinate (t2) at ($(d)!1.3!(c)$);
	\coordinate (b1) at ($(a)!1.2!(e)$);
	\coordinate (b2) at ($(e)!1.3!(a)$);
	\draw[dashed,olines] (t1) .. controls ($(j) + (2,2)$) and ($(i) + (-1.4,2)$) .. (t2) -- cycle;
	\draw[dashed] (t1) -- (b1);
	\draw[dashed] (t2) .. controls ($(t2) + (-0.5,-1)$) and ($(b2) + (-0.5,1)$) .. (b2);
	\draw[dashed,dblines] (b1) .. controls ($(h) + (2,-2.3)$) and ($(g) + (-2,-2.3)$) .. (b2) -- cycle;
	\node at (-2,3.5) {\textcolor{olabels}{$P$}};
	\node at (-3.4,0.5) {$X$};
	\node at (-2.8,-3.3) {\textcolor{dblabels}{$Q$}};
\end{tikzpicture}
\caption{Illustration for Lemma~\ref{pack1}}\label{picture_pack1}
\end{figure}

\begin{lemma}
\label{pack1}
Let $M$ be a $k$-coherent matroid and let 
$(P,X,Q)$ be a path of $3$-separations in $M$ where $X$ is a
prime step.
Assume that no element of $X$ is in a triangle or a triad, and that
there is no $3$-separating set $Z\subseteq X$ with $|Z|\geq 3$.
Let $f$ be a feral element of $X$. Let $(P_1,P_2,\ldots,P_m)$ and
$(Q_1,Q_2,\ldots,Q_n)$ be maximal $k$-fractures of $M\ba f$ and $M/f$
respectively. Then there is a partition $\{\{x_1,x_2\},Z_1,Z_2\}$
of $X-\{f\}$  such that, up to labels, duality and equivalence of flowers,
the following hold.
\begin{itemize}
\item[(i)] $m=k$, $P_1=Q\cup Z_1\cup Z_2$,
$P_2\cup P_3\udots P_{k-1}=P$, $P_k=\{x_1,x_2\}$, and 
$(Q\cup X,P_2,P_3,\ldots,P_{k-1})$ is a swirl-like flower of 
order $k-1$ in $M$.
\item[(ii)] $n=k$, $Q_3\cup Q_4\udots Q_k=Q$, 
$Q_1=Z_1\cup\{x_1,x_2\}$,
$Q_2=Z_2\cup P$, and $(P\cup X,Q_3,Q_4,\ldots,Q_k)$
is a swirl-like flower of order $k-1$ in $M$.
\item[(iii)] $M\ba x_1$ is $k$-coherent.
\end{itemize}
\end{lemma}

\begin{proof}
Up to duality we may assume that $(P_1,P_2,\ldots,P_m)$ and
$(Q_1,Q_2,\ldots,Q_n)$ form a feral display for $f$. It 
follows that $n=k$ and that
$(Q_1\cup Q_2\cup\{f\},Q_3,\ldots,Q_k)$ is a maximal 
swirl-like flower of $M$.  Now consider
the way that the $3$-separations $(P,X\cup Q)$ and $(P\cup X,Q)$
interact with this flower. If neither $P$ nor $Q$ is contained in a
petal of the flower, then one deduces that $X$ contains a 
$3$-separating set equivalent to the non-sequential 3-separating
set $Q_1\cup Q_2\cup\{f\}$. Thus, up to switching the labels of $P$
and $Q$, we may assume that $P$ is contained in a petal $Q'$.
Assume that $Q'\neq Q_1\cup Q_2\cup\{f\}$.Then we deduce that
$Q_1\cup Q_2\cup\{f\}\subseteq X$, contradicting the assumption that
$X$ contains no $3$-separating set $Z$ with $|Z|\geq 3$.
Hence $P\subseteq Q_1\cup Q_2\cup\{f\}$. If 
$P\cup X\not\subseteq Q_1\cup Q_2\cup\{f\}$, we contradict the 
assumption that $X$ is prime. 
Thus $P\cup X\subseteq Q_1\cup Q_2\cup\{f\}$.
Now $(Q_1,Q_2,\ldots,Q_k)$ is a maximal flower in $M/f$, and
$((P\cup X)-\{f\},Q)$ is a 3-separation. We know that
$(P\cup X)-\{f\}\subseteq Q_1\cup Q_2$. If 
$(P\cup X)-\{f\}$ is contained in either $\fcl(Q_1)$ or $\fcl(Q_2)$,
we obtain the contradiction that $M$ is $k$-fractured. 
Therefore $(P\cup X)-\{f\}$ is equivalent to $Q_1\cup Q_2$,
and we may assume that $(P\cup X)-\{f\}=Q_1\cup Q_2$. 

Consider the flower $(P_1,P_2,\ldots,P_n)$ of $M\ba f$. 
As we have a feral display for $f$, we have
$P_1\supseteq Q_3\cup Q_4\udots Q_k$ and, for some 
$i\geq 3$, the partition $(P_2,\ldots,P_i,P_{i+1}\udots P_n\cup\{f\})$
is a swirl-like flower of order $i$ in $M$. 
Arguing as before we deduce that
$P_2\cup P_3\udots P_i=P$. Now 
$P_{i+1}\cup P_{i+2}\udots P_n\subseteq X$
and $\lambda_M(P_{i+1}\cup P_{i+2}\udots P_n)=2$. 
By the assumption that
$X$ contains no $3$-separating set $Z$ with $|Z|\geq 3$, we
must have $n=k$, $i=k-1$, and $|P_n|=2$. Let $P_k=P_n=\{x_1,x_2\}$.
From the properties of a feral display we 
see that there is a partition
$(Z_1,Z_2)$ of $X-\{x_1,x_2\}$ such that $Q_1=Z_1\cup\{x_1,x_2\}$
and $Q_2=Z_2\cup P$. This shows that parts (i) and (ii) of the lemma
hold.

Consider part (iii). The partition
$(P_1,P_2,\ldots,P_k)$ is a $k$-fracture of $M\ba f$. 
Moreover, $\{x_1,x_2\}$ is a fully-closed petal
of this $k$-fracture. It now follows from
Lemma~\ref{2-elt-feral} that $M\ba x_1$ is $k$-coherent. 
\end{proof}

The next lemma is elementary.

\begin{lemma}
\label{display-sequential}
Let $S$ be a sequential $3$-separating set in a 
$3$-connected matroid $M$ where
$|S|,|E(M)-S|\geq 3$ and let $(U,V)$ be a non-sequential 
$3$-separation of $M$.
Then either $S\subseteq \fcl(U)$ or $S\subseteq \fcl(V)$.
\end{lemma}

\begin{proof} 
As $S$ is sequential it contains a triangle or triad $T$ such that
$S\subseteq \fcl(T)$. If $|U\cap T|\geq 2$, then 
$S\subseteq \fcl (U\cap T)$,
so that $S\subseteq \fcl(U)$. The lemma follows from this observation.
\end{proof}

The next lemma shows that feral packs are essentially the only obstruction
to obtaining a satisfactory simplification of a path of $3$-separations. 

\begin{lemma}
\label{path1}
Let $M$ be a $k$-coherent matroid and let $(P,X,Q)$ be a
path of $3$-separations where $X$ is a prime step
and $|P|,|Q|\geq 3$. Let 
$X'=(\fcl(P\cup X)\cap Q)\cup(\fcl(Q\cup X)\cap P)$. 
Assume that there is no fixed element in $X'$ such that $M\ba x$ is
$k$-coherent and no cofixed element in $X'$ such that $M/x$
is $k$-coherent. 
\begin{itemize}
\item[(i)] If $|X|\geq 2$, then there exists $x\in X$ such that
either $M\ba x$ or $M/x$ is $k$-coherent.
\item[(ii)] If $|X|\geq 3$, $X$ contains a clonal pair
$\{u,v\}$, and $X$ is not a feral pack,
then there exists $x\in X-\{u,v\}$ such that
either $M\ba x$ or $M/x$ is $k$-coherent.
\end{itemize}
\end{lemma}

\begin{proof}
If we are in case (i) of the lemma,  let 
$Y=X$ and if we are in case (ii), let $Y=X-\{u,v\}$.
Assume that there is no element $y\in Y$ such that either
$M\ba y$ or $M/y$ is $k$-coherent. Note that, as $X$ is prime,
both $P$ and $Q$ are fully closed.

\begin{sublemma}
\label{subpath1}
No element of $Y$ is in both a triangle and a triad.
\end{sublemma}

\subproof
If $y\in Y$ is in both a triangle and a triad, then $y$ is in a fan $F$ with
at least four elements. As $|X|>1$ and $X$ is prime, 
it is easily seen
that $M$ is not a wheel or a whirl. 
Hence $F$ has ends. If $e$ is an end of $F$,
then either $M\ba e$ or $M/e$ is $k$-coherent by 
Corollary~\ref{remove-sequential} and 
$e$ is respectively fixed or cofixed.
Evidently either $e\in X'$ or $e\in X$. 
The former case contradicts our assumption about the elements of $X'$.
Consider the latter case. The element $e$ is either fixed or cofixed,
so $e\notin\{u,v\}$. Thus $e\in Y$, contradicting our assumption about
the elements of $Y$.
\end{proof}

\begin{sublemma}
\label{subpath2}
No  element of $X$ belongs to a triangle or triad.
\end{sublemma}

\subproof
If the sublemma fails, then we may assume that there is a triangle 
$T=\{a,b,c\}$ that meets $X$.

Assume that $T$ is not $k$-wild. 
Note that the elements of $T$ are in $X\cup X'$.
Say $a\in X'$. We may assume that $a\in P$. Then $a\in\cl(P-\{a\})$,
but $T\not\subseteq \cl(P-\{a\})$, so by Lemma~\ref{freedom3},
$a$ is fixed in $M$. By \ref{subpath1}, and 
the fact that $T$ is not $k$-wild, 
we see that there is an element $z\in T$ such that
$M\ba z$ is $k$-coherent. The only way that this does not
contradict assumptions about the members of $X'\cup Y$ is to 
have $z\in\{u,v\}$. But in this  case 
$a\notin \{u,v\}$. By Corollary~\ref{free-triangle},
$M\ba a$ is $k$-coherent again contradicting our
assumptions about elements of $Y\cup X'$. 

Assume that $T$ is $k$-wild. By \ref{subpath1}, 
$T$ is either a standard or
costandard $k$-wild triangle.
Let
$(A_1,A_2,\ldots,A_{k-2},B_1,B_2,\ldots,B_{k-2},C_1,C_2,\ldots,C_{k-2})$
be a $k$-wild display for $T$, where
$A=A_1\cup A_2\udots A_{k-2}$, $B=B_1\cup B_2\udots B_{k-2}$
and $C=C_1\cup C_2\udots C_{k-2}$. Let $(U,V)$ be a non-sequential
$3$-separation of $M$. If neither $U$ nor $V$ is contained
in any of $\fcl(A)$, $\fcl(B)$ or $\fcl(C)$, then one deduces by
uncrossing that one of the flowers of $M$ displayed by the 
$k$-wild display is not maximal. This contradicts the fact that
$M$ is $k$-coherent.
This means that
for any non-sequential 3-separation $(U,V)$, up to labels,
$U$ is contained in either $\fcl(A)$, $\fcl(B)$ or 
$\fcl(C)$. If $|P\cup Q|=4$, the sublemma is clear. 
Thus we may assume that $|P|>2$. Note that
$T\not\subseteq \fcl(P)$. If $P$ is non-sequential, then we see that we
can assume that $P\subseteq \fcl(A)$ by the above observation. 
If $P$ is sequential, the same conclusion follows from 
Lemma~\ref{display-sequential}.
If $|Q|>2$, then we may assume that $Q\subseteq \fcl(B)$ so that
$X$ contains a set equivalent to $C$, 
and we see that $X$ contains at least
three peripheral sets and it is clear 
that $Y$ contains a peripheral set $Z\subseteq C$ 
whose full closure avoids $\{u,v\}$. 
If $p\in\fcl(Z)\cap P$, then $p$ is a loose element of the swirl-like flower 
$(A_1,A_2,\ldots,A_{k-2},B\cup C\cup T)$ of $M$. 
But then, by Lemma~\ref{loose-removable},
up to duality, $X'$ has  a fixed element whose deletion 
preserves $k$-coherence.
Hence $\fcl(Z)$ avoids $P$ and similarly avoids $Q$.
Thus $\fcl(Z)\subseteq Y$. But now, by Corollary~\ref{upcor}, $Y$
contains an element that can be either 
deleted or contracted to preserve
$k$-coherence. We omit the easy analysis for the case that $|Q|=2$.
\end{proof}

\begin{sublemma}
\label{subpath3}
There is no $3$-separating set $Z\subseteq X$ with $|Z|\geq 3$.
\end{sublemma}

\subproof
Assume that $Z$ is such a set. If $Z$ is sequential then it contains a
triangle or triad contradicting \ref{subpath2}. 
Thus $Z$ is non-sequential and $|Z|\geq 4$.
In this case it follows from
Lemma~\ref{pack0} $Z$ contains an element $z\notin\{u,v\}$ 
such that $M\ba z$ or $M/z$ is $k$-coherent. This contradiction establishes the sublemma.
\end{proof}

\begin{sublemma}
\label{subpath4}
For all $y\in X$, both $M\ba y$ and $M/y$ are $3$-connected.
\end{sublemma}

\subproof
If the lemma fails, then, as no element of $X$ is  
in a triangle or triad, we 
may assume up to duality that, 
for some $y\in X$, the element $y$ is in the 
guts of a vertical $3$-separation 
$(Y_1\cup\{y\},Y_2)$. If $Y_1\subseteq P$, then $y\in\cl(P)$,
contradicting the fact that $P$ is fully closed. Hence 
neither $Y_1$ nor $Y_2$ is contained in either $P$ or $Q$.
If $P\subseteq Y_i$ and $Q\subseteq Y_j$ for some permutation
$(i,j)$ of $(1,2)$, then $(P,Y_i-P,\{y\},Y_j-Q,Q)$ is a path of
$3$-separations contradicting the fact that $P$ is prime. Thus
$(Y_1,Y_2)$ crosses at least one of $(P,X\cup Q)$ and $(P\cup X,Q)$.

Up to symmetry, we may now assume that $(Y_1,Y_2)$ crosses $(P,X\cup Q)$
and further that, if $(Y_1,Y_2)$ crosses $(P\cup X,Q)$, then all
vertical $3$-separations of $M$ with $y$ in the guts cross both
$(P,X\cup Q)$ and $(P\cup X,Q)$. Without loss of generality,
$|P\cap Y_1|\geq 2$. So, by uncrossing, we see that
$\lambda((X\cup Q)\cap Y_2)\leq 2$ and 
$\lambda((X\cup Q)\cap(Y_2\cup \{y\}))\leq 2$. If 
$r((X\cup Q)\cap Y_2)\geq 3$, then 
$(P\cup Y_1\cup\{y\},(X\cup Q)\cap Y_2)$ is a vertical
3-separation with $y$ in the guts that does not cross 
$(P,X\cup Q)$. Moreover this 3-separation crosses
$(P\cup X,Q)$ if and only if $(Y_1,Y_2)$ does and we
have contradicted the assumption about the choice
of $(Y_1,Y_2)$.
Thus $r((X\cup Q)\cap Y_2)\leq 2$. If
$r((X\cup Q)\cap Y_2)=2$, then $y\in\cl((X\cup Q)\cap Y_2)$ and
$y$ is in a triangle. Hence $|(X\cup Q)\cap Y_2|=1$, so that 
$|P\cap Y_2|>1$. We may now repeat the argument to conclude that 
$|Y_1\cap (X\cup Q)|=1$, so that $|X\cup Q|=3$ giving numerous
contradictions, one of which is to the fact that $y$ is not in a triangle or triad.
\end{proof}

Assume that $Y=X$. Then each element of $X$ is feral. But this contradicts
Lemma~\ref{pack0}. This contradiction shows that part (i) holds. Assume
that $Y=X-\{u,v\}$. Then each element of $X-\{u,v\}$ is a feral
element and $X$ is a feral pack. This contradiction
shows that (ii) holds.
\end{proof}

The next lemma shows that feral packs cannot hunt in packs in 
a $k$-skeleton. In fact they cannot hunt in packs in arbitrary paths
of $3$-separations, but it is easier to prove for $k$-skeletons.

\begin{lemma}
\label{pack6}
Let $M$ be a $k$-skeleton with a maximal path 
$\PP=(P_0,P_1,\ldots,P_l)$
of $3$-separations. Assume that $P_i$ is a feral pack for some 
$i\in\{1,2,\ldots,l-2\}$. Then $\lambda(P_{i+1})=2$ and 
$(P_{i+1}^-,P_{i+1},P_{i+1}^+)$ is a spiral. In particular,
$P_{i+1}$ is not a feral pack.
\end{lemma}

\begin{proof}
By Lemma~\ref{pack1}, $M$ has a maximal swirl-like flower
$\RR=(R_1,R_2,\ldots,R_{k-1})$, where $R_1=P_0\cup P_1\udots P_i$.
By Lemma~\ref{clean}, swirl-like flowers of order 
$(k-1)$ in $k$-skeletons
are canonical. Consider the $3$-separation
$(P_{i+1}^-\cup P_{i+1},P_{i+1}^+)$ of $M$. If $P_{i+1}^+$
is contained in a petal of $\RR$, then we contradict the
fact that $\PP$ is maximal. It follows that
an equivalent $3$-separation equivalent
to this is displayed by $\RR$, and, as $\RR$ is canonical, this
$3$-separation must be precisely 
$(P_{i+1}^-\cup P_{i+1},P_{i+1}^+)$. Thus $P_{i+1}^-\cup P_{i+1}$
is a union
of consecutive petals of $\RR$ and it follows that
$(P_{i+1}^-,P_{i+1},P_{i+1}^+)$ is a concatenation of $\RR$ 
so that $(P_{i+1}^-,P_{i+1},P_{i+1}^+)$ is a spiral and 
$\lambda(P_{i+1})=2$. It follows routinely that 
$P_i$ is not a feral pack.
\end{proof}

\section{Proof of the Main Theorem} 

We are almost in a position to achieve the primary 
goal of this chapter and prove Theorem~\ref{main1}. 
We need three more lemmas.

\begin{lemma}
\label{main99}
Let $M$ be a $k$-skeleton and $(P_0,p_1,P_2,P_3)$ 
be a path of $3$-separations
in $M$ where $p_1$ is a fixed guts 
singleton and $P_2$ is prime. Then
$(P_0\cup\{p_1\},P_2,P_3)$ is a tight spiral in $M$.
\end{lemma}

\begin{proof} 
As $M$ is a skeleton, no element of $M$ is in 
both a triangle and a triad.
We first prove that $M\ba p_1$ is 3-connected. Assume that $r(P_0)=2$.
Then $P_0\cup\{p_1\}$ is a triangle. As $p_1$ is fixed, 
this triangle is
$k$-wild by Lemma~\ref{skeleton-triangle}. 
But $p_1$ is not in a $4$-element
fan, so this triangle is either standard or costandard in which case
$M\ba p_1$ is $3$-connected.

Assume that $r(P_0)>2$. Then $\co(M\ba p_1)$ 
is 3-connected by Bixby's Lemma.
Assume that $p_1$ is in a triad $T$. As $M/ p_1$ is not 3-connected,
$T$ is not $k$-wild and
therefore $T$ is a clonal triple contradicting the fact that 
$p_1$ is fixed in $M$.

Thus $M\ba p_1$ is $3$-connected and, as $p_1$ is fixed in $M$, 
we see that $M\ba p_1$ is $k$-fractured.
Let $\RR=(R_0,R_1,\ldots,R_m)$ be a maximal 
$k$-fracture of $M\ba p_1$.
It is easily seen that, for some $j\in\{0,1,\ldots,m-2\}$, we have, 
up to labels in 
$\RR$, that $P_0=R_0\udots R_j$ and
that $(P_0\cup\{p_1\},R_{j+1},\ldots,R_m)$ is a 
maximal swirl-like flower of $M$.
By Lemma~\ref{clean}, this flower has no loose elements so that
$P_0\cup\{p_1\}$ is fully closed. This means that $P_2$
is not a singleton and it follows routinely 
that  $(P_0\cup\{p_1\},P_2,P_3)$
is a tight spiral in $M$. 
\end{proof}

Recall that ${\mathcal U}_q$ and ${\mathcal U}^*_q$ denote the
classes of matroids with no $U_{2,q+2}$= and $U_{q,q+2}$-minor 
respectively.

\begin{lemma}
\label{main4}
Let $M$ be a $k$-coherent matroid in $\youq$ with a path
$\PP=(P_0,P_1,\ldots,P_l)$ of $3$-separations. 
Then there 
is a function $f_{\ref{main4}}(m,q)$ such that, if at least
$f_{\ref{main4}}(m,q)$ internal steps of $\PP$ contain clonal 
pairs, then $M$ contains a $k$-coherent minor with a path
of clonal pairs of length $m$.
\end{lemma}

\begin{proof}
Let $m\geq 2$ be an integer and let $f_{\ref{main4}}(m,q)=2(q+2)(m+1)$.
Assume that $n\geq f_{\ref{main4}}(m,q)$ and that, for this value
of $m$, the matroid $M$ is a minor-minimal counterexample 
to the lemma.
Given the hypotheses of the lemma, we have a distinguished set $S$
of $n$ clonal pairs such that each step of the path contains at most
one clonal pair in $S$.
By taking an appropriate concatenation we may assume that
$P_0$ and $P_l$ each contain a clonal pair in $S$.
By Lemma~\ref{path1} and the minimality assumption, 
we see that each internal step of 
$\PP$ is either a singleton, a clonal pair or a feral pack.
A clonal pair $P_i=\{p_i,q_i\}$ may be non-prime so that 
$(P_0,P_1,\ldots,P_{i-1},p_i,q_i,P_{i+1},\ldots,P_l)$ is also
a path of 3-separations, or else it may be prime. For the moment we
do not refine the path by splitting non-prime clonal pairs.

If $M$ has a
fixed element $x$ such that $M\ba x$ is $k$-coherent, then it is
clear that the induced path in $M\ba x$ has the properties of $\PP$
described above. By this observation, its dual and the minimality
assumption we deduce that $M$ is a $k$-skeleton.
It now follows by Lemma~\ref{pack6}, that if $i\in\{2,3,\ldots,l-2\}$
and $P_i$ is a feral pack,
then both $P_{i-1}$ and $P_{i+1}$ are prime clonal pairs. 
Thus at most $n/2$ of the steps that contain 
clonal pairs are feral packs.
We now sacrifice the clonal pairs in feral packs. 
Let $S'$ denote the set of
clonal pairs in $S$ that are not in feral packs. Then $S'$ contains
at least $n/2$ clonal pairs. Let  $M'$ be a 
minimal-minimal $k$-coherent minor of 
$M$ whose ground set contains $S'$ and let
$\QQ=(Q_0,Q_1,\ldots,Q_{l'})$ be the path of $3$-separations 
induced by $\PP$
in $M'$. Note that $\QQ$ is well-defined since both 
$P_0$ and $P_l$ contain
clonal pairs in $S'$. By Lemma~\ref{path1}, the internal steps of
$\QQ$ are either clonal pairs or singletons. 
It is also clear that $M'$
is a $k$-skeleton and that 
if $Q_i=\{q_i\}$ is a singleton internal step,
then neither $M'\ba q_i$ nor $M'/q_i$ is $k$-coherent.

\begin{sublemma}
\label{main4.2}
If $2\leq i\leq l'-2$ and $Q_i=\{q_i\}$ is a guts singleton, 
then $q_i$ is not fixed in $M'$.
\end{sublemma} 

\subproof
Assume that $q_i$ is a fixed guts singleton. 
By Lemma~\ref{main99}
$(Q_{i}^+\cup\{q_i\},Q_{i+1},Q_{i+1}^+)$
is a tight spiral flower in $M$, so that 
$Q_{i+1}$ is a prime clonal pair. Evidently
 $q_i\notin\cl(Q_{i+1})$
and $q_i\notin\cl(Q_{i+1}^+)$. Thus, by Lemma~\ref{strand3}, $q_i$
is not fixed from the right in $M'$. A symmetric argument shows that 
$q_i$ is not fixed from the left in $M'$. By Lemma~\ref{strand2},
$q_i$ is not fixed in $M$.
\end{proof}

Refine $\QQ$ to a maximal path by splitting the non-prime clonal pairs 
in internal steps to
obtain a path $\QQ'=(Q_0',Q'_1,\ldots,Q_l')$. By 
\ref{main4.2} and its dual, all the singletons
of this path are either unfixed guts singletons 
or uncofixed coguts singletons.
If there are at least $2(q+2)$ such singletons, 
then, by Lemma~\ref{strand5},
$M$ has either a $U_{2,q+2}$- or a $U_{q,q+2}$-minor 
contradicting the
fact that $M\in\youq$. 
Thus there are at most $2(q+2)-1$ such singletons, so that $\QQ'$ has
a section of at least $f_{\ref{main4}}(q,m)/2(q+2)=m+1$ consecutive prime clonal pairs. We conclude that
$M$ has a $k$-coherent minor with a path of clonal pairs of length $m$ 
contradicting the assumption that $M$ was a counterexample.
\end{proof}

Finally, a lemma about the way the gangs of three can appear in paths.
The lemma follows easily from the structure of the 3-separations of
$M$ around a gang of three and we omit the proof.

\begin{lemma}
\label{gangway}
Let $(P,X,Q)$ be a path of $3$-separations in the $k$-coherent matroid
$M$. If $\{x,y,z\}$ is a gang of three in $M$ and
$\{x,y,z\}\cap X\neq\emptyset$, then $\{x,y,z\}\subseteq X$ 
and $|X|\geq 4$.
\end{lemma}

We are at last in a position to prove Theorem~\ref{main1} which
we restate here for convenience.

\begin{theorem}
\label{main3}
Let $k\geq 5$ and $q\geq 2$ be integers. Then there is a function
$f_{\ref{main3}}(k,q)$
such that, if $M$ is a $k$-skeleton with a path of $3$-separations
of length $f_{\ref{main3}}(k,q)$, then $M\notin \eq$.
\end{theorem}

\begin{proof}
Let
$f_{\ref{main3}}(k,q)=2f_{\ref{main4}}
(f_{\ref{special13}}(k,q),q)+2(q+2)+1.$
Assume that $M$ is a minor-minimal counterexample to the theorem in that
$M$ is a $k$-skeleton in $\eq$ with a path of
$3$-separations of length $n=f_{\ref{main3}}(k,q)$, but no
proper $k$-skeleton minor of $M$ has such a path. Let 
$\PP=(P_0,P_1,\ldots,P_n)$ be a path of $3$-separations of 
length $n$ in $M$. We may assume that each internal step of $\PP$
is prime.

\begin{sublemma}
\label{main3.1}
Each internal step of $\PP$ is either a 
singleton or contains a clonal pair.
\end{sublemma}

\subproof
Assume that $P_i$ is an internal step that has at least two elements but
no clonal pair. As $M$ is a $k$-skeleton, $M$ has no fixed element $x$ 
such that $M\ba x$ is $k$-coherent and dually. Thus, by 
Lemma~\ref{path1}(i), there is an element $x\in P_i$ such that either 
$M\ba x$ or $M/x$ is $k$-coherent. If either $M\ba x$ or 
$M/x$ is a $k$-skeleton, then 
$(P_0,P_1,\ldots,P_{i-1},P_i-\{x\},P_{i+1},\ldots,P_n)$ is a path of 
$3$-separations in a proper $k$-skeleton minor of $M$, contradicting the choice
of $M$. So neither $M\ba x$ nor $M/x$ is a $k$-skeleton. If 
$x$ is comparable with another element, then, by 
Theorem~\ref{reduce-comparable},
either $M\ba x$ or $M/x$ is a $k$-skeleton. Thus $x$ is not comparable
with any other element of $M$. In this case, by Theorem~\ref{contract},
it must be the case that, up to duality, $x$ is in a gang of three. Say
$\{x,y,z\}$ is such a gang of three. By Lemma~\ref{gangway}, 
$\{x,y,z\}\subseteq P_i$ and $|P_i|\geq 4$. 
But now, by
Theorem~\ref{gogang}, $M/x\ba y,z$ is a $k$-skeleton, and
$(P_0,P_1,\ldots,P_{i-1},P_i-\{x,y,z\},P_{i+1},\ldots,P_n)$
is a path of $3$-separations of length $n$
in a proper $k$-skeleton minor
of $M$, again contradicting the choice of $M$.
\end{proof}

If there are at least $q+2$ unfixed guts singletons or 
$q+2$ uncofixed coguts singletons in $\PP$, then, 
by Lemma~\ref{strand5},
$M$ has a $U_{2,q+2}$- or a $U_{q,q+2}$-minor. 
Thus there are at 
least $n-2(q+2)$ steps that are either fixed guts singletons,
cofixed coguts singletons or contain a clonal pair.
By Lemma~\ref{main99} and its dual,
if $2<i<n-1$, and $P_i=\{p_i\}$ is a fixed guts singleton or cofixed
coguts singleton, then neither
$P_{i-1}$ nor $P_{i+1}$ are singletons and therefore both of these
sets contain clonal pairs. Hence there are at least 
$(n-2(q+2)-1)/2$ internal steps that contain clonal pairs, that is,
at least $f_{\ref{main4}}(f_{\ref{special13}}(k,q),q)$ internal steps
contain clonal pairs. By Lemma~\ref{main4}, $M$ has a $k$-coherent
minor with a path of clonal pairs of length $f_{\ref{special13}}(k,q)$.
But now, by Corollary~\ref{special13}, $M$ is not in $\youq$,
contradicting the assumption that $M$ was a counterexample to the theorem.
\end{proof}

\section{A Last Lemma}

We conclude this chapter with a lemma on paths of 
3-separations that focusses on displayed swirl-like flowers 
without making assumptions about the $k$-coherence of the
matroid. We will use the following straightforward result.

\begin{lemma}
\label{shrink-flower}
Let $(P_1,P_2,\ldots,P_l)$ be a maximal swirl-like flower
of order $l\geq 5$ in the $3$-connected matroid $M$. Assume
that $P_1$ is fully closed. Then the following hold.
\begin{itemize}
\item[(i)] $M$ has a $3$-connected minor on
$P_2\cup P_3\udots P_l$.
\item[(ii)] If $M'$ is a $3$-connected minor of $M$
on $P_2\cup P_3\udots P_l$, and $(X,Y)$ is a $3$-separation
of $M'$, then either $(X\cup P_1,Y)$ or $(X,Y\cup P_1)$
is a $3$-separation of $M$.
\end{itemize}
\end{lemma}

\begin{proof}
Part (i) follows from Theorem~\ref{upgrade} and an easy
induction. Consider (ii). It is easily seen that
$(P_2,P_3,\ldots,P_l)$ is a maximal flower in
$M'$. Say $(X,Y)$ is a 3-separation of $M'$. 
If $(X,Y)$ is displayed by $(P_2,P_3,\ldots,P_l)$,
then (ii) follows easily. Otherwise we may assume that
$X\subseteq \fcl_{M'}(P_i)$ for some $i\in\{2,\ldots,l\}$.
But then $(X,Y\cup P_1)$ is clearly a 3-separation of $M$.
\end{proof}

Let $\PP=(P_0,P_1,\ldots,P_n)$ be a path of $3$-separations in
the connected matroid $M$ and let $\{p_i,q_i\}$ be a clonal
pair contained in $P_i$ for some $i\in\{0,1,\ldots,n\}$.
Then $\{p_i,q_i\}$ is $\PP$-{\em strong}
\index{$\PP$-strong clonal pair} 
if 
$\kappa(\{p_i,q_i\},P_0\cup P_1\udots P_{i-1}\cup P_{i+1}\udots P_n=2$.
In other words, $\{p_i,q_i\}$ is $\PP$-strong if there is no
2-separating set $A$ with the property that 
$\{p_i,q_i\}\subseteq A\subseteq P_i$.

\begin{lemma}
\label{at-last}
Let $\PP=(P_0,P_1,\ldots,P_n)$ be a path of $3$-separations in
the connected matroid $M$ each step of which contains a 
$\PP$-strong
clonal pair and let $k$ be an integer. Then there is a function
$f_{\ref{at-last}}(k,q)$ such that, if $\PP$ displays no
swirl-like flower of order $k$, and $n\geq f_{\ref{at-last}}(k,q)$,
then $M\notin \eq$.
\end{lemma}

\begin{proof}
Assume that $k\geq 5$.
Let $f_{\ref{at-last}}(k,q)=f_{\ref{main1}}(k,q)+2$.
We claim that the lemma holds with this definition of 
$f_{\ref{at-last}}(k,q)$. 

Assume not.
Let $M$ be a counterexample to the lemma, so that
$M$ contains a path $\PP=(P_0,P_1,\ldots,P_l)$ of 3-separations
satisfying the hypotheses of the lemma where 
$l\geq f_{\ref{at-last}}(k,q)$, and $M\in \eq$. 

Assume that $M$ is chosen to have a ground set of minimum
cardinality. 
It follows routinely from the definitions of
path of $3$-separations and $\PP$-strong clonal pair
that, under this assumption, $M$ is 3-connected. We now show that
$M$ is $k$-coherent. Assume not. Then there is a swirl-like flower
$\QQ=(Q_1,Q_2,\ldots,Q_t)$ of order $t\geq k$ in $M$. A subflower
of $\QQ$, displayed by $\PP$ has order at most $k-1$. It
follows from this fact and 
Lemma~\ref{display1}, that, 
up to flowers equivalent to $\QQ$, there is a step $P_i$ of 
$\PP$ such that $Q_j\cup Q_{j+1}\subseteq P_i$
for some $j\in\{1,2,\ldots,t\}$. Let $\{p,p'\}$ be a $\PP$-strong 
clonal pair in $P_i$. Then we may assume that $Q_j$
avoids $\{p,p'\}$. It is also easily seen that we may
assume that $P_i$ is fully closed. By Lemma~\ref{shrink-flower}(i),
we may remove elements from $P_i$ to obtain
a $3$-connected minor $M'$ of $M$. But it follows from
Lemma~\ref{shrink-flower}(ii) that
$(P_1,P_2,\ldots,P_{i-1},P_i-Q_j,P_{i+1},\ldots,P_l)$ is a path
of $3$-separations in $M'$ satisfying the hypotheses of the 
lemma and we have contradicted the minimality of the choice of
$|E(M)|$. 

Thus $M$ is $k$-coherent. This contradicts
Corollary~\ref{main2}.
In the case that $k\leq 5$, the lemma holds by letting
$f_{\ref{at-last}}(k,q)=f_{\ref{main1}}(5,q)+2$.
\end{proof}

Finally we not a result that is a more-or-less immediate
corollary of Lemma~\ref{at-last} and Lemma~\ref{free-link}.

\begin{corollary}
\label{get-free-clonal}
Let $\PP$ be a path of $3$-separations of length $n$
in the connected matroid $M$ 
in $\eq$, each step of which contains a $\PP$-strong
clonal pair.
Then there is a function $f_{\ref{get-free-clonal}}(m,q)$
such that, if $n\geq f_{\ref{get-free-clonal}}(m,q)$,
then $M$ has a $\Delta_m$ minor, each leg of which
is equal to $\{p_i,q_i\}$ for some $i$.
\end{corollary}

\chapter{Taming a Skeleton}
\label{taming}

\section{Introduction}

In the previous chapter we proved that a $k$-skeleton in
$\eq$ cannot have an arbitrarily long path of 3-separations.
In this chapter we continue the process of controlling structure
in skeletons. The goal is to prove the following theorem.

\begin{theorem}
\label{tame}
There is a function $f_{\ref{tame}}(m,k,q)$ such that, if $M$ is 
a $k$-skeleton in $\eq$ with at least $f_{\ref{tame}}(m,k,q)$
elements, then $M$ has a $4$-connected minor whose ground set
contains a set of $m$ pairwise disjoint clonal pairs.
\end{theorem}

With this theorem in hand we can forget about $k$-skeletons and focus on
$4$-connected matroids with many clonal pairs. We begin by 
learning more about $3$-trees associated with $k$-skeletons in $\eq$.
We already know that they cannot contains arbitrarily long
induced paths. If $v$ is a flower vertex of such a tree,
then, as $M$ is a $k$-skeleton the associated flower cannot
have high degree as it induces a path of $3$-separations in 
$M$. But so, far there is nothing to control the degree of a bag vertex.

\section{Potatoes}

In this section we prove the following theorem.

\begin{theorem}
\label{small-tree}
Let $M$ be a $k$-skeleton in $\eq$, and let $T$ be a $3$-tree for
$M$. Then
there is a function $f_{\ref{small-tree}}(m,k,q)$ such that, if
$T$ has at least $f_{\ref{small-tree}}(m,k,q)$ vertices,
then $M$ has a $4$-connected minor with a set of $m$ pairwise-disjoint
clonal pairs.
\end{theorem}

As noted above we need to control the degree of bag vertices
in a $3$-tree for a $k$-skeleton. Such a vertex $v$ typically
exemplifies a more highly connected part of the matroid,
attached to which are the 3-separating sets displayed
by $v$. This motivates the
definition of ``potato'' that we now give.

Let $M$ be a $3$-connected matroid and $n\geq 0$ be an integer. 
A {\em potato}
\index{potato} 
of $M$ is a
partition $\{P_1,P_2,\ldots,P_n\}$ of a subset of  
$E(M)$ such that
the following hold.
\begin{itemize}
\item[(i)] $\lambda(P_i)=2$, $|P_i|\geq 3$, and $|E(M)-P_i|\geq 5$,
for all $i\in\{1,\ldots,n\}$.
\item[(ii)] If $(X,Y)$ is a 3-separation of $M$, then, for some 
$i\in\{1,\ldots,n\}$, either $X$ or $Y$ is contained in $P_i$.
\item[(iii)] $\lambda(P_i\cup P_j)\geq 3$ for all distinct $i$ and $j$
in $\{1,2,\ldots,n\}$.
\end{itemize}
The set $E(M)-(P_1\cup P_2\udots P_n)$ is the {\em core} of the potato
\index{core of a potato}
and each $P_i$ is an {\em eye} of the potato.
\index{eye of a potato} 
A potato may have an empty core.
If $M$ is 4-connected, then the empty set of eyes defines a potato.
If $n\geq 4$, then condition (iii) is redundant, but if $n\in\{2,3\}$,
then (iii) says that a potato cannot be a 3-petal flower.
Note that the eyes of a potato are fully-closed sets.

Readers familiar with matroid 
tangles (see for example \cite{gegewh09}) 
will recognise that, except for trivially small matroids,
a tangle $\tau$
of order 4 can be associated with a potato 
$\PP$ of the 3-connected matroid
$M$. If $\lambda(X,Y)\leq 2$, then 
$(X,Y)\in\tau$ if either $|X|\leq 2$,
or $X$ is contained in an eye of $\PP$. 
Of course, not all tangles of order
4 can be derived from a potato in this way. 

Our main interest is in potatoes in $k$-coherent matroids
and $k$-skeletons.
If 
$\PP$ is a potato of the $k$-coherent matroid $M$, then one annoying 
possibility is that one of its eyes is a $k$-wild triangle or triad. 
Let $\{a,b,c\}$ be 
a $k$-wild triangle of $M$. Assume that $\{a,b,c\}$ is an eye of 
$\PP$. Then $\{a,b,c\}$ is certainly not in a 
$4$-element fan, so that it is either standard or costandard,
in which case there is a partition $(A,B,C,\{a,b,c\})$
of $E(M)$ such that $A$, $B$, and $C$ are each non-trivial $3$-separating
sets. If $\{a,b,c\}$ is standard, 
then $(A\cup\{a\},B\cup C\cup\{b,c\})$ 
is a $3$-separation that crosses $\{a,b,c\}$, so it must be the case that
$\{a,b,c\}$ is costandard. In this case there must be eyes that contain
$A$, $B$ and $C$, so the potato must be $\{A,B,C,\{a,b,c\}\}$.  
Using the above notation we have
proved the following lemma.

\begin{lemma}
\label{potato3}
If $\{a,b,c\}$ is a $k$-wild triangle of the $k$-coherent matroid $M$
with associated partition $(A,B,C,\{a,b,c\})$,
and $\{a,b,c\}$ is an eye of the potato $\PP$, then $\{a,b,c\}$
is costandard  and $\PP=\{A,B,C,\{a,b,c\}\}$.
In particular $\PP$ has four eyes, an empty core, and  at
most one eye of $\PP$ is a $k$-wild triangle or triad.
\end{lemma}

For
an arbitrary $k$-coherent matroid the converse
of Lemma~\ref{potato3} need not hold as, for example,
there may be a single element in $\cl(A)\cap B$. But the converse
does hold for $k$-skeletons. We omit the easy proof.

\begin{lemma}
\label{potato3-5}
If $\{a,b,c\}$ is a costandard $k$-wild triangle of the $k$-skeleton
$M$ with associated partition $(A,B,C,\{a,b,c\})$, 
then $\{A,B,C,\{a,b,c\}\}$ is a potato in $M$.
\end{lemma}

The next lemma helps us to find potatoes in $k$-skeletons.

\begin{lemma}
\label{potato2}
Let $M$ be a $k$-skeleton and let $A_1$ 
and $A_2$ be $3$-separating sets each having 
at least three elements. 
If $\lambda(A_1\cup A_2)>2$ and $A_1\cap A_2\neq\emptyset$,
then either  $A_1$ or $A_2$ is a $k$-wild triangle or triad.
\end{lemma}

\begin{proof}
Assume that $A_1\cap A_2\neq\emptyset$. As 
$\lambda(A_1)=\lambda(A_2)=2$ and $\lambda(A_1\cup A_2)>2$, 
we have $\lambda(A_1\cap A_2)<2$.
Thus $|A_1\cap A_2|=1$. Say $A_1\cap A_2=\{a\}$. As $A_1$ and 
$E(M)-A_2$ are 3-separating and the union of these sets avoids at least 
two elements of $M$, we have $\lambda(A_1\cap(E(M)-A_2))\leq 2$,
that is, $\lambda(A_1-\{a\})=2$ and similarly $\lambda(A_2-\{a\})=2$.
Now $a\in\clstar(A_1-\{a\})$ and, up to duality, we may assume that
$a\in\cl(A_1-\{a\})$. Evidently $a\notin\cl^*(A_2-\{a\})$, so that we
also have $a\in\cl(A_2-\{a\})$.

We now show that $a$ is fixed in $M$. Say that $\sqcap(A_1,A_2)=2$.
Then
$$\lambda(A_1\cup A_2)=\lambda(A_1)+\lambda(A_2)-
\sqcap(A_1,A_2)-\sqcap^*(A_1,A_2)\leq 2.$$
Thus $\sqcap(A_1,A_2)=1$. But 
$\sqcap(A_1-\{a\},A_2-\{a\})=\sqcap(A_1,A_2)$, and it follows by 
Lemma~\ref{freedom1} that $a$ is fixed in $M$.

If either $A_1$ or $A_2$ is a triangle, then, by 
Lemma~\ref{skeleton-triangle}, that
triangle is $k$-wild and the lemma holds. Therefore we may assume
that neither $A_1$ nor $A_2$ is a triangle and indeed that 
$r(A_1),r(A_2)\geq 3$. Thus $\si(M/a)$ is not 3-connected.
Now by Bixby's Lemma and the fact that $M$ has no
$4$-element fans, the matroid 
$M\ba a$ is 3-connected and therefore $k$-fractured;
that is, $a$ is semi-feral. Let $\PP=(P_1,P_2,\ldots,P_m)$
be a maximal fracture of $M\ba a$. By Corollary~\ref{clonal3}, there
exists $j\in\{2,3,\ldots,m-1\}$ such that $a$ is in the guts of
$(P_1\cup P_2\udots P_j,P_{j+1}\cup P_{j+2}\udots P_m)$. 
As $M$ is a $k$-skeleton
and $a$ is fixed in $M$, it follows from Theorem~\ref{canon1} that
$\PP$ is canonical. Note that, if $A_1-\{a\}$ is contained in a
petal of $\PP$, then $a\notin\cl(A_1-\{a\})$. Thus $A_1-\{a\}$ and
similarly $E(M\ba a)-A_1$ is not contained in a petal of $\PP$.
Thus, by Lemma~\ref{cross1}, $A_1-\{a\}$ is equal to a union of
petals of $\PP$. It is easily seen that this can only hold if, up
to labels, $A_1-\{a\}=P_1\cup P_2\udots P_j$. Similarly we must have 
$A_2-\{a\}=P_{j+1}\cup P_{j+2}\udots P_m$. 
But then we obtain the contradiction that $\sqcap(A_1,A_2)=2$.
\end{proof}

The next two lemmas show that we may shrink the eyes of
a potato $\PP$ to find
a 4-connected minor of a $k$-skeleton with a clonal pair for each eye
of $\PP$.

\begin{lemma}
\label{potato4}
Let $\PP=\{P_1,P_2,\ldots,P_n\}$ be a potato of the $k$-skeleton
$M$. Then there is a $k$-skeleton minor $M'$ of $M$ with a potato
$\PP'=\{P_1',P'_2,\ldots,P_n'\}$ such that the following hold.
\begin{itemize}
\item[(i)] The core of $\PP'$ is equal to the core of $\PP$.
\item[(ii)] $P'_i\subseteq P_i$ for $i\in\{1,\ldots,n\}$.
\item[(iii)] At most one eye of $\PP'$ does not contain a clonal pair.
\end{itemize}
\end{lemma}

\begin{proof}
Let $P$ be an eye of $\PP$. 
Assume that $|P|\geq 4$.
If $P$ is sequential with at least four elements, then 
$P$ contains a triangle or triad that is certainly not $k$-wild
and again $P$ contains a clonal pair. 
Assume that $P$ is non-sequential.
Then $P$ contains a strongly peripheral set $Z$. 
If $Z$ contains no clonal pairs,
then, by Lemma~\ref{cool-peripheral}, 
there is an element $z\in Z$ such that
$M\ba z$ or $M/z$ is a $k$-skeleton. Up to duality we may assume that
$M\ba z$ is a $k$-skeleton. Assume that $P=P_1$. We now
show that $\{P_1-\{z\},P_2,\ldots,P_n\}$ is  potato of $M\ba z$.

Condition (iii)  of a potato clearly holds.
Consider condition (i). If this fails, then 
$|E(M\ba z)-P_i|<5$ for some $i>1$, 
that is $|E(M)-P_i|=5$. But, in this case,
$n=2$ and the core of $\PP$ has at most one element, meaning that
$\lambda_M(P_1\cup P_2)<3$ contradicting the definition of a potato.
 
For the more substantial case consider property (ii) of
a potato. 
Note that, as $M\ba z$ is 
3-connected, $z\in\cl(P_1-\{z\})$.
Let $(X,Y)$ be a 3-separation of 
$M\ba z$. If $X$ contains $P_1-\{z\}$, then $(X\cup\{z\},Y)$
is a 3-separation of $M$ and, by the definition of potato,
$Y\subseteq P_i$ for some $i\in\{1,2,\ldots,n\}$, so that (ii) holds
in this case. Thus we may assume that $(X,Y)$ 
crosses $(P_1-\{z\},E(M)-P_1)$. We now
consider subcases of this case. Assume that
$\lambda_{M\ba z}(X-P_1)=\lambda_{M\ba z}(Y-P_1)=2$.
These 3-separating sets are not blocked by $z$ and it follows that 
$(X-P_1,Y-P_1,P_1)$ is a partition of $E(M)$ into 3-separating sets.
From the definition of potato, $|E(M)-P_1|\geq 5$, so 
we may assume that $|X-P_1|\geq 3$.
Consider the 3-separation $(X-P_1,Y\cup P_1)$. 
By property (ii) of a potaot, we may assume that
$X-P_1$ is contained in $P_2$. But then we have the contradiction that
$Y-P_1$ is contained in both $P_3$ and $P_4$.
Therefore it is not the case that
$\lambda_{M\ba z}(X-P_1)=\lambda_{M\ba z}(Y-P_1)=2$.

If $\lambda_{M\ba z}(X-P_1)<2$, then $|X-P_1|=1$ so that 
$P_1-\{z\}$ is not fully closed in $M\ba z$ contradicting the fact that
$P_1$ is fully closed in $M$. 

For the last, somewhat irritating, subcase  we may assume that
$\lambda_{M\ba z}(X-P_1)>2$. 
This means that $|Y\cap (P_1-\{z\})|=1$. But now, $P_1-\{z\}$ and $Y$
are 3-separating sets in the $k$-skeleton $M\ba z$ that meet, 
and $\lambda((P_1-\{z\}) \cup Y)>3$. By Lemma~\ref{potato2}
we deduce that $Y$ is a $k$-wild triangle or triad of $M\ba z$.
Say $Y=\{a,b,c\}$, where $\{a\}=Y\cap(P_1-\{z\})$. Note that,
as $\PP$ is a potato in $M$, the set $\{a,b,c\}$ is blocked in
$M$, so $\{a,b,c\}$ is a triad, and $a\in\cl^*(P_1-\{z\})$. Let
$(A_1,A_2,\ldots,A_{k-2},B_1,B_2,\ldots,B_{k-2},
C_1,C_2,\ldots,C_{k-2})$
be a $k$-wild display for $\{a,b,c\}$ and let
$A=A_1\cup A_2\udots A_{k-2}$, $B=B_1\cup B_2\udots B_{k-2}$ and
$C=C_1\cup C_2\udots C_{k-2}$.
We have $a\in\cl^*_{M\ba z}(P_1-\{z,a\})$, 
but $P_1-\{z\}$ is fully closed in $M\ba z$.. 
It is readily checked that, since $M$ is a $k$-skeleton,
this implies that
$P_1-\{z\}=A\cup\{a\}$. In this case $\{a,b,c\}$ is a standard 
$k$-wild triad. Therefore $B\cup\{b\}$ and $C\cup\{c\}$ are 
$3$-separating.
Certainly they are not triangles or triads. 
Thus we may assume that for some distinct
$i$ and $j$, we have $B\cup\{b\}\subseteq P_i$ and 
$C\cup\{c\}\subseteq P_j$.
Hence $\PP=(P_1,P_i,P_j)=(A\cup\{a,z\},B\cup\{b\},C\cup\{c\})$.
But now $\lambda_M(P_1\cup P_i)=2$, contradicting the definition of 
a potato. Therefore 
$\{P_1-\{z\},P_2,\ldots,P_n\}$ is indeed a potato in $M\ba z$. 

The process may be repeated to obtain a $k$-skeleton minor 
$M'$ of $M$ with a potato $\PP'=\{P'_1,P'_2\ldots,P'_n\}$ such that
each eye either contains a clonal pair of $M'$
or is a triangle or triad.
If one of these is $k$-wild 
then the lemma holds by Lemma~\ref{potato3}.
Otherwise each triangle
or triad is a clonal triple  
and the lemma holds in this case too.
\end{proof}

\begin{lemma}
\label{potato5}
Let $\PP=\{P_1,P_2,\ldots,P_n\}$ be a potato of the $3$-connected matroid
$M$ with the property that at most one eye of $\PP$ does not contain
a clonal pair. Then $M$ contains a $4$-connected
minor containing the core of $\PP$ and a set of 
$n-1$ pairwise-disjoint clonal pairs.
\end{lemma}

\begin{proof}

Assume that the potato $\PP=\{P_1,P_2,\ldots,P_n\}$ 
has core $C$. Let $\{p,q\}$ be a clonal pair in 
$P_1$. By Tutte's Linking Theorem there is a 3-connected minor $N$
of $M$ on $\{p,q\}\cup P_2\udots P_n\cup C$ such that
$N|(P_2\cup P_3\udots P_n\cup C)=M|(P_2\cup P_3\udots P_n\cup C)$.
We show that $\{P_2,\ldots,P_n\}$ is a potato in $N$.

Properties (i) and (iii) of a potato are clear. Consider (ii).
Let $(X,Y)$ be a 3-separation in $N$. Assume that $\{p,q\}\subseteq X$.
Then $(X\cup P_1,Y)$ is a 3-separation in $M$. But $X\neq \{p,q\}$,
so $X\cup P_1$ is not contained in an eye of $\PP$. Therefore
$Y\subseteq P_i$ for some $i\in\{2,3,\ldots,n\}$ as required.
Assume that $p\in X$ and $q\in Y$. As $p$ and $q$ are clones, we may
assume, up to duality, that both $p$ and $q$ are in the guts of $(X,Y)$.
Thus $\{p,q\}\subseteq\cl(Y-\{q\})$. If $|Y-\{q\}|=2$, then
$N|(Y\cup\{p\})\cong U_{2,4}$ and, if $y\in Y-\{q\}$, then 
$(P_1\cup\{y\},E(M)-(P_1\cup\{y\}))$ is a 3-separation in $M$
that violates condition (ii). Therefore $(X\cup\{q\},Y-\{q\})$
is a 3-separation of $N$. Arguing as above we have $Y\subseteq P_i$
for some $i\in\{2,3,\ldots,n\}$, so that $\{p,q\}\subseteq\cl_N(P_i)$.
It is now easily checked that $\lambda_M(P_1\cup P_i)=2$,
contradicting the fact that $\PP$ satisfies (iv).

Thus $\{P_2,P_3,\ldots,P_n\}$ is indeed a potato in $N$.
The core of this potato is $C\cup\{p,q\}$.
Repeating this process $n$ times gives a 3-connected matroid $M'$
with a potato having an empty set of eyes and a core consisting of
$C$ together with a disjoint set of $n$ clonal pairs. If 
$(X,Y)$ is a $3$-separation then either $X$ or $Y$ is contained in an
eye of this potato. This impossibility shows 
that $M'$ is 4-connected.
\end{proof}

\subsection*{Sequential $3$-separating Sets}

Here we simply recall for convenience some elementary facts
and prove an easy lemma.
Let $Z$ be a sequential 3-separating set
of the 3-connected matroid $M$. Recall that
a {\em sequential ordering} for $Z$ is an
ordering $(z_1,z_2,\ldots,z_t)$ of $Z$ such that, 
$\{z_1,z_2,\ldots,z_i\}$
is $3$-separating for all $i\in\{1,2,\ldots,t\}$. 
Recall also that  subset $Z'$ of $Z$ {\em generates} 
$Z$ if some ordering of $Z'$ is an initial segment of
a $3$-sequence for $Z$. If $|Z'|\geq 2$, then it is easily seen
that $Z'$ generates $Z$ if and only if $Z\subseteq \fcl(Z')$.
The first three elements of a 3-sequence for $Z$ either form a triangle
or triad, so every sequential 3-separating set with at least three elements
is generated by some triangle or triad $T$. 
It is then immediate that
any 2-element subset of $T$ generates $Z$. 

\begin{lemma}
\label{seq-gen}
Let $Z$ be a sequential $3$-separating set of the $3$-connected
matroid $M$ and let $X$ and $Y$ be $3$-separating subsets of $Z$.
where $Z\subseteq X\cup Y$.  Then either $X$ or $Y$ generates $Z$.
\end{lemma}

\begin{proof}
Note that either
$X$ or $Y$ contains at least two elements of a generating
triangle or triad for $Z$. Assume the former. Then
$Z\subseteq \fcl(T\cap X)$ so that $Z\subseteq \fcl(X)$. By 
Lemma~\ref{sub-seq-3-sep} $X$ is sequential. Thus $X$ has
a sequential ordering. As $Z\subseteq \fcl(X)$
this extends to a sequential ordering of $Z$.
\end{proof}

\subsection*{Strong Bag Vertices}
Let $T$ be a 3-tree of the 3-connected matroid $M$, and let $v$ be a bag
vertex of $T$. Certain bag vertices identify more highly connected parts
of the matroids---in other words tangles of order 4---others do not.
If the vertex has degree at least 4, then it is certainly the case
that $v$ identifies a more connected part of the matroid.
For the main result of this section we are only interested in vertices
of high degree, but later we will care about the low degree case.

In the low degree case, the issue is 
to identify whether there is substance in the bag.
In general, having a lot of elements in $B_v$ does not suffice
as we could, for example,
just be finding a large sequential 3-separating set. 
Say that $(X,Y)$ is a non-sequential 3-separation, then there is a
path $(\coh(X),z_1,z_2,\ldots,z_n,\coh(Y))$ 
of $3$-separations of $M$ and the internal elements
of this path
can all be in a bag of low degree. If $M$ is a $k$-skeleton in $\eq$,
then, by Theorem~\ref{main1}, this number is bounded
by $f_{\ref{main1}}(k,q)$. In fact a much more modest bound is 
straightforwardly seen to hold, but we take the lazy way out.

\begin{lemma}
\label{lazy-bound}
Let $M$ be a $k$-skeleton in $\eq$ and 
let $(X,Y)$ be a $3$-separation of $M$.
\begin{itemize}
\item[(i)] If $X$ is sequential, the $X$ has at most 
$f_{\ref{main1}}(k,q)$ elements.
\item[(ii)] If $(X,Y)$ is non-sequential, 
then $E(M)-(\coh(X)\cup\coh(Y))$
has at most $f_{\ref{main1}}(k,q)$ elements.
\end{itemize}
\end{lemma}

Let $M$ be a $k$-skeleton in $\eq$ and let $T$ be a 3-tree for $M$.
If $v$ is a bag vertex of $T$, then we denote the subset
of $E(M)$ that labels $v$ by $B_v$. If we say that $B$ is a {\em bag}
of $M$, then we mean that $B=B_v$ for some bag vertex of a 
$3$-tree for $M$.
We now define what it means for a bag vertex to be {\em strong}.
\index{strong bag vertex}
If $v$ has degree at least $3$, then $v$ is strong. If 
$v$ has degree at most $2$, then $v$ is strong if 
$|B_v|\geq 4f_{\ref{main1}}(k,q)$.

Let $v$ be a strong bag vertex of $T$. Define the 
collection ${\mathcal E}_v$ of
3-separating sets of $M$ as follows. 
A $3$-separating set $X$ is in ${\mathcal E}_v$
if $|X|\geq 3$ and $X$ is equivalent to 
a $3$-separating set displayed by
$v$ or is equivalent to one contained in $B_v$. Note that as 
3-trees do not attempt to display sequential 3-separating sets, there
may well be many of the latter type. 
The maximal members of ${\mathcal E}_v$
are the {\em eyes} of $v$. 
\index{eyes of a bag vertex}

\begin{lemma}
\label{potato6}
Let $T$ be a $3$-tree for the $k$-skeleton $M$ and let $v$ be a strong
bag vertex of $T$. Then the eyes of $v$ form a potato of $M$.
\end{lemma}

\begin{proof}
Let $\PP=\{P_1,\ldots,P_n\}$ be the collection of eyes of $v$.

\begin{sublemma}
\label{potato6.1}
The lemma holds if an eye of $\PP$ is a $k$-wild triangle.
\end{sublemma}

\subproof
Say $P_1=\{a,b,c\}$ is $k$-wild with associated partition
$(A,B,C,\{a,b,c\})$ obtained from a $k$-wild display. 
Evidently $\{a,b,c\}$ is costandard.
The sets $A$, $B$, $C$ and $\{a,b,c\}$ are clearly the eyes of 
$v$ and by Lemma~\ref{potato3-5} they form a potato.
\end{proof}

In what follows we may assume that no eye of $\PP$ is a
$k$-wild triangle or triad.

\begin{sublemma}
\label{potato6-1}
$|E(M)-P_i|\geq 5$ for all $i\in\{1,2,\ldots,n\}$
\end{sublemma}

\subproof The claim is clear if $n\geq 3$. If $n\leq 2$, then
the sublemma is a consequence 
of the definition of strong bag vertex and
Lemma~\ref{lazy-bound}.
\end{proof}

\begin{sublemma}
\label{potato6-2}
If $P_i$ and $P_j$ are distinct eyes of $\PP$, then 
$\lambda(P_i\cup P_j)\geq 3$.
\end{sublemma}

\subproof
Consider the eyes
$P_1$ and $P_2$ of $\PP$. 
Assume that $\lambda(P_1\cup P_2)\leq 2$. Consider the case that
$n=2$. Then $d(v)\leq 2$, so that $|B_v|\geq 4f_{\ref{main1}}(k,q)$.
Let $B'_v=B_v-(\fcl(P_1)\cup \fcl(P_2))$. Certainly
$\lambda(B'_v)\leq 2$, and, by Lemma~\ref{lazy-bound},
$|B'_v|\geq 2f_{\ref{main1}}(k,q)$. But then $B'_v$ is a non-sequential
3-separating set that is not displayed by $T$, contradicting the 
definition of a 3-tree.

Assume that $n\geq 3$. Then $\lambda(P_1\cup P_2)\geq 2$ so that
$\lambda(E(M)-(P_1\cup P_2))=2$. If this set is sequential, then
$d(v)\leq2$ and we again obtain the contradiction that 
$B_v-(\fcl(P_1)\cup\fcl(P_2))$ is a sequential 
$3$-separating set with more than $f_{\ref{main1}}(k,q)$ elements.
Hence $E(M)-(P_1\cup P_2)$ is non-sequential.
If $P_1\cup P_2$ is sequential, then by 
Lemma~\ref{seq-gen}, $P_1\cup P_2$ is generated by either
$P_1$ or $P_2$ contradicting the definition of eyes.
So $P_1\cup P_2$ is non-sequential. But, again, it is easily seen that
$(P_1,P_2-P_1,E(M)-(P_1\cup P_2))$ is a flower in $M$. As this
flower is not displayed in $T$, at least one petal, say $P_1$ must be
sequential and $v$ has degree at most 2, so that $B_v$ has at least
$4f_{\ref{main1}}(k,q)$ elements. But, again by Lemma~\ref{lazy-bound}
we see that $B_v$ is covered by three sets whose union has at most
$3f_{\ref{main1}}(k,q)$ elements. 
\end{proof}

\begin{sublemma}
\label{potato6-3}
The eyes of $\PP$ are pairwise disjoint.
\end{sublemma}

\subproof
This follows from \ref{potato6-2},
Lemma~\ref{potato2}
and the assumption that no eye of $\PP$ is a 
$k$-wild triangle or triad.
\end{proof}

\begin{sublemma}
\label{potato6-4}
If $(X,Y)$ is a $3$-separation of $M$, then either
$X$ or $Y$ is contained in an eye of $\PP$.
\end{sublemma}

\subproof
If $(X,Y)$
is non-sequential then $(X,Y)$ 
is equivalent to a 3-separation displayed by
$T$ and it is clear that either $X$ or $Y$ is contained in an eye
of $v$. Say $X$ is sequential. Then $X=\fcl(T)$ for some triangle or
triad $T$. If $T$ is $k$-wild, then $X=T$ contradicting the 
assumption that no eye of $\PP$ is a $k$-wild triangle or triad.
Thus  $T$ is a clonal triple in which case the argument
that $T$ (and hence $\fcl(T)$) is 
contained in $P_i$ for some $i$ is clear. 
\end{proof}

The lemma follows 
from \ref{potato6-1},
\ref{potato6-2}, \ref{potato6-3} and \ref{potato6-4}. 
\end{proof}

We are now in a position to prove the main result of this section.

\begin{proof}[Proof of Theorem~\ref{small-tree}]
Let $\mu(m,k,q)={\rm max}\{m,f_{\ref{main1}}(k,q)\}$. Now
let $f_{\ref{small-tree}}(m,k,q)=\mu(m,k,q)^{f_{\ref{main1}}(k,q)}$.
Assume that $T$ has at least $f_{\ref{small-tree}}(m,k,q)$ vertices.

We may assume that $m\geq 5$. Say  $T$ has a bag vertex $v$ of degree
$m$. Then $v$ is strong and has at least $m$ eyes.
So, by Lemma~\ref{potato6}, $M$ has a potato with $m$ eyes,
and by Lemma~\ref{potato4}, $M$ has a 3-connected minor
$N$ with a potato with
$m$ eyes such that each eye contains a clonal pair. Then, by
Lemma~\ref{potato5}, $M$ has a $4$-connected minor with $m$ clonal
pairs. Thus the lemma holds if $T$ has a bag vertex of degree $m$,
so we may assume that each bag vertex has degree less than $m$.

Each flower vertex of degree $l$ induces a path of $3$-separations
in $M$ of length $l$, so by Theorem~\ref{main1}, $T$ has no flower
vertex of degree at least $f_{\ref{main1}}(k,q)$. Hence the degree of
each vertex is less than $\mu(m,k,q)$. But as $T$ has at least
$f_{\ref{small-tree}}(m,k,q)$ vertices, it has a path of length 
$f_{\ref{main1}}(k,q)$ contradicting the fact that $M\in \eq$.
\end{proof}

\section{Bounding Fettered Elements}

An element of a matroid is {\em fettered}
\index{fettered element}
if it is either fixed or cofixed.
Otherwise it is {\em unfettered}.
\index{unfettered element} 
A matroid is {\em unfettered}
\index{unfettered matroid} 
if it has no fettered elements.

In the previous section it was shown that if a 
3-tree $T$ for a $k$-skeleton
$M$ in $\eq$ has  sufficiently many vertices, 
then $M$ has a $4$-connected minor with many 
clonal pairs. This turns the focus on to 3-trees of bounded size.
The task of this section is to show that in this case we gain control 
over the
number of fettered elements. 

\begin{theorem}
\label{control-fetter}
Let $M$ be a $k$-skeleton in $\eq$ and let $T$ be a $3$-tree
for $M$ with $\nu$ vertices.
Then there is a function $f_{\ref{control-fetter}}(\nu,k,q)$ 
such that  $M$ has at most $f_{\ref{control-fetter}}(\nu,k,q)$
fettered elements.
\end{theorem}

There are no surprises in the next lemma.

\begin{lemma}
\label{local-swirl}
Let $(L,P,R,Q)$ be a swirl-like flower in the $3$-connected matroid
$M$. Let $(L',R')$ be a partition of $L\cup P\cup R$ such
that $L\subseteq L'\subseteq L\cup P$. Then
$\lambda_{M\ba Q}(L',R')=\lambda_M(L'\cup Q,R')-1$.
\end{lemma}

\begin{proof}
As $L\subseteq L'\subseteq L\cup P$, we see that
$\sqcap(L,Q)\leq \sqcap(L',Q)\leq \sqcap(L\cup P,Q)$, so that
$\sqcap(L',Q)=1$. Also $\lambda_M(Q)=2$
so that $r(Q)-r(M)=2-r(L\cup P\cup R)$. Therefore

\begin{align*}
\lambda_M(L'\cup Q,R')-1&=r(L'\cup Q)+r(R')-r(M)-1\\
&=r(L')+r(Q)+r(R')-r(M)-2\\
&= r(L')+r(R')-r(L\cup P\cup R)\\
&=\lambda_{M\ba Q}(L',R')
\end{align*}
as required.
\end{proof}

The next lemma is essentially the dual of Lemma~\ref{swirl-cofix}.

\begin{lemma}
\label{fetter4}
Let $e$ be an unfixed element of the $3$-connected matroid 
$M$ and $l\geq 4$ be an integer.
Let $\PP=(P_1,P_2,\ldots,P_l)$ be a swirl-like flower of $M/e$. 
Then, up to labels, 
$(P_1\cup P_2\cup\{e\},P_3,\ldots,P_l)$ is a flower in $M$.
\end{lemma}

The next corollary follows immediately from Lemma~\ref{fetter4}.

\begin{corollary}
\label{fetter5}
Let $M$ be a $k$-coherent matroid and $e$ be an element of $M$
that is not fixed and has the property that $M/e$ is $3$-connected.
Then $M/e$ is $(k+1)$-coherent, and, if $\PP$ is a flower in
$M$ that opens to a $k$-fracture in $M/e$, then $\PP$
has order $k-1$.
\end{corollary}

\begin{lemma}
\label{fetter6}
Let $M$ be a $k$-skeleton in $\eq$ and let $P$ 
be a petal of the swirl-like 
flower $\PP$ of $M$. Let $\{a_1,a_2,\ldots,a_n\}$ 
be a set of fettered elements
in $\PP$ such that the following hold for all $i\in\{1,2,\ldots,n\}$.
\begin{itemize}
\item[(i)] $M\ba a_i$ is $k$-coherent.
\item[(ii)] $M/a_i$ is $3$-connected.
\item[(iii)] There is a $k$-fracture of $M/a_i$ that opens $\PP$.
\end{itemize}
Then $n\leq q+1$.
\end{lemma}

\begin{proof}
Say $i\in\{1,2,\ldots,n\}$. As $M\ba a_i$ is $k$-coherent, $a_i$ is not
fixed in $M$. By Corollary~\ref{fetter5}, $\PP$ is a swirl-like
flower  of order
$k-1$ in $M$ and the flower in $M/a_i$ that opens $\PP$ has order $k$.
Say $\PP=(P,P_2,\ldots,P_{k-1})$. 
Then for each $i\in\{1,\ldots,n\}$,
there is a partition $(L_i,R_i)$ of $P$ such that 
$\PP_i=(L_i,R_i,P_2,\ldots,P_k)$ is a maximal swirl-like
flower of $M/a_i$. 
We now show that $\PP_i$ is canonical.

As $M\ba a_i$ is $k$-coherent, 
$a_i$ is not in a bogan couple, so, if $p$ is a loose element
of $\PP$, then, by Theorem~\ref{canon1}, $a\less p$. It is 
easily checked that 
$p$ can only be a guts element in this flower, so $p$ is fixed in
$M/a_i$. This means that either $a_i$ is fixed in $M$, or 
$a_i$ and $p$ are a clonal pair. As $M\ba a_i$ is $k$-coherent,
$a_i$ is not fixed and the latter case contradicts the fact that $a_i$
is cofixed. Hence $\PP_i$ is indeed canonical. 
Let $R=P_3$, $Q=P_4$ and $L=P_5\cup P_6\udots P_{k-1}$, so that
$(L,P,R,Q)$ is a (cyclically shifted) concatenation of $\PP$.

\begin{sublemma}
\label{fetter6.1}
If $i\neq j$ and $a_j\in R_i$, then $L_i\subseteq L_j$.
\end{sublemma}

\subproof
We know that $(L,L_i,R_i,R,Q)$ is a canonical flower in $M/a_i$. 
If $((L\cup L_j\cup\{a_j\})-\{a_i\},(R_j\cup R\cup Q)-\{a_i\})$
is a 3-separation in $M/a_i$, we obtain the contradiction that $a_j$ is 
a loose element of the above flower. Hence
$((L\cup L_j\cup\{a_j\})-\{a_i\},(R_j\cup R\cup Q)-\{a_i\})$
is an exact $4$-separation in $M/a_i$.
It now follows
by Lemma~\ref{4-cross0} that there is a partition $(C,R_j)$ of $R_i$
such that $L_j-\{a_i\}=L_i\cup C$. But $a_i\in\cl_M(L\cup L_i)$,
and it follows readily that $L_j=L_i\cup\{a_i\}\cup C$, so that 
$L_i\subseteq L_j$.
\end{proof}

A consequence of \ref{fetter6.1} is that we may assume that indices
are chosen for elements of $\{a_1,a_2,\ldots,a_n\}$ so that 
$L_1\subseteq L_2\subseteq \cdots\subseteq L_n$. 
Consider $M\ba Q$. By Lemma~\ref{local-swirl},
$\lambda_{M\ba Q}(L\cup L_i\cup\{a_i\},R_i\cup R)
=\lambda_{M\ba Q}(L\cup L_i,\{a_i\}\cup R_i\cup R)=2$.
Again, by Lemma~\ref{local-swirl}, and the fact that $\PP$ is maximal,
$\kappa_{M\ba Q}(L\cup L_1,R\cup R_n)=2$. For $i\in\{2,3,\ldots,n\}$,
let $L'_i=L_i-(L_{i-1}\cup\{a_i\})$. Then
$$(L\cup L_1,a_1,L_2',a_2,L_3',\ldots,L'_n,a_n,R_n\cup R)$$
is a path of 3-separations in $M\ba Q$.

Say $i\in\{1,2,\ldots,n\}$. It is easily seen that 
$a_i$ is not a coguts singleton in the path described above.
Hence $a_i$ is a guts
singleton. Now $M\ba a_i$ is $k$-coherent,
so that $a_i$ is not fixed in $M$ and hence is not fixed
$M\ba Q$.
Thus Lemma~\ref{strand5} applies and we deduce that  
$M$ has a $U_{2,n}$-minor. As $M$ is in $\eq$, we conclude that
$n\leq q+1$.
\end{proof}

\subsection*{Bounding Feral Elements}

The task of this subsection is to bound the number of feral elements
in a bag $B_v$ of a 3-tree of a 
$k$-skeleton as a function of the degree
of $v$. We first develop some terminology specific to this task.

Let $M$ be a 3-connected matroid and let $\PP=(P_1,P_2,\ldots,P_m)$
be a tight swirl-like flower in $M$, where $m\geq 3$. An element
$f\in P_1$ {\em expands $\PP$ at $P_2$} if $M\ba f$ is 
3-connected and there is a partition 
$(P_1',P'_{m+1})$ of $P_1-\{f\}$ such that 
\begin{itemize}
\item[(i)] $(P_1',P_2,P_3,\ldots,P_m,P'_{m+1})$
is a tight swirl-like flower in $M\ba f$, and
\item[(ii)] $P_1'$ is 3-separating in $M$.
\end{itemize}
The element $f$
{\em expands\ } $\PP$ if it expands $\PP$ at either $P_2$ or $P_m$.
Assume that $f$ expands $\PP$ at $P_2$ and that
$(P_1',P_2,\ldots,P_m,P'_{m+1})$ has  properties 
(i) and (ii) above. Assume that, 
in addition,  $(P_1',P_2,\ldots,P_m,P'_{m+1})$
has the property that 
whenever $(Q_1',P_2,\ldots,P_m,Q'_{m+1})$
is another flower in $M\ba f$
obtained by expanding $\PP$ at $P_2$ for which
$P_1'\subseteq Q_1'$, then $P_1'=Q_1'$. 
In this case we will say that 
$(P'_1,P_2,\ldots,P_m,P'_{m+1})$ is an 
{\em expansion of $\PP$ by $f$ at $P_2$}.
\index{expansion of $\PP$}
It is clear that if $f$ expands $\PP$, then an expansion of
$\PP$ by $f$ at either $P_2$ or $P_m$ exists.

\begin{lemma}
\label{tame-feral}
Let $M$ be a $k$-skeleton in $\eq$ and let 
$v$ be a strong bag vertex of 
a $3$-tree $T$ for $M$. Let $\PP=(P_1,P_2,\ldots,P_m)$ be a 
maximal swirl-like
flower of $M$, where $m\geq 3$ and $B_v\subseteq P_1$.
Assume that $\PP$ is displayed by a vertex adjacent to $v$.
Then there are at most $8(q^2+q+1)$ feral elements in $B_v$ that
expand $\PP$.
\end{lemma}

\begin{proof}
Assume that the flowers $\PP'=(P_1',P_2,\ldots,P_m,P'_{m+1})$
and $\QQ'=(Q'_1,P_2,\ldots,P_m,Q'_{m+1})$ are expansions of
$\PP$ at $P_2$ by the distinct feral elements $f$ and $g$ in $B_v$.
By the maximality of the choice of $P'_1$ and $Q'_1$ and the fact
that $M$ is a $k$-skeleton we deduce that $P'_1$ and $Q'_1$ are 
fully closed $3$-separating sets of $M$.

Let $\RR'=(R'_1,P_2,\ldots,P_m,R'_{m+1})$ be an expansion of $\PP$
at $P_2$ by an element in $P_1$. We will say that $\RR'$ is 
{\em equivalent}
\index{equivalent expansions}
to $\PP'$ if $P'_1=R'_1$.

\begin{sublemma}
\label{tame-feral1}
There are at most $q^2+q+1$ distinct
expansions of $\PP$ that are equivalent
to $\PP$.
\end{sublemma}

\subproof
Note that $\lambda_M(P'_{m+1}\cup\{f\})=3$, as otherwise $\PP$ is not a 
maximal flower in $M$. Thus $f$, and any other element, that expands
$\PP$ to a flower equivalent to $\PP'$ is in the coclosure in $M$
of $P'_1\cup P_2\udots P_m$. As $M$ is cosimple,
the dual of Lemma~\ref{kung2} implies that this set has at most 
$q^2+q+1$ elements.
\end{proof}

\begin{sublemma}
\label{tame-feral2}
If $|P'_1|\geq 3$, then $P'_1$ is an eye of $v$ and
$g\notin P'_1$.
\end{sublemma}

\subproof
Assume that $P'_1$ is not sequential. Then a 3-separating set
equivalent to $P'_1$ is displayed in $T$. The maximality of the choice of 
$P'_1$ ensures that $P'_1$ is displayed by an edge of $T$. If this 
edge is not incident with $v$, then we again contradict 
the choice of $P'_1$. Thus $P'_1$ is an eye of $v$. Elements in 
$P'_1\cap B_v$ are either in the guts or coguts of a 3-separation
and are therefore not feral. Thus $g\notin P'_1$. Assume that
$P'_1$ is sequential. It is easily checked that $P'_1$ is not a 
$k$-wild
triangle or triad. Therefore $P'_1=\fcl(T)$ for some clonal triangle or
triad of $M$ and again it is  evident that
if $P'_1$ is not an eye of $v$ we contradict the maximality of the 
choice of $P'_1$. Again it is clear that $g\notin P'_1$.
\end{proof}

\begin{sublemma}
\label{tame-feral3}
If $\PP'$ is not equivalent to $\QQ'$ and
$|P'_1|\geq 3$, then $|Q'_1|=2$.
\end{sublemma}

\subproof
Assume that $P'_1$ and $Q'_1$ both have at least three elements.
By \ref{tame-feral2} $P'_1$ and $Q'_1$ are eyes of $v$ and by  
Lemma~\ref{potato6} the eyes of $v$ form a potato. Thus $P'_1$
and $Q'_1$ are disjoint. But $Q'_1$ does not contain a feral element
in $B_v$. Hence $Q'_1\subseteq P'_{m+1}$. But then
$\sqcap(Q'_1,P_2)=0$ contradicting the fact that these are adjacent
petals in the swirl-like flower $\QQ'$.
\end{proof}

\begin{sublemma}
\label{tame-feral4}
If $|P'_1|=|Q'_1|=2$, and $P'_1\neq Q'_1$, then $f\in Q'_1$.
\end{sublemma}

\subproof
Say that $P'_1=\{f_1,f_2\}$ and $Q'_1=\{g_1,g_2\}$.
As $\sqcap(Q'_1,P_2)=1$ we deduce that $Q'_1\not\subseteq P'_{m+1}$.
Therefore we may assume that $g_1\in\{f_1,f_2,f\}$. Assume for a 
contradiction that $g_1=f_1$. Then $\{f_2,g_2\}\subseteq\cl(P_2\cup\{f_1\})$,
so that $r(P_2\cup\{f_1,f_2,g_2\})=r(P_2)+1$. But
$\sqcap(P_1,P_2)=1$ as these are adjacent petals of swirl-like flower in
$M$. Thus $\sqcap(P_2,\{f_1,f_2,g_2\})=1$ and it follows that
$r(\{f_1,f_2,g_2\})=2$ meaning that this set is a triangle. This 
contradicts the maximality of the choice of $P'_1$. We conclude
that $g_1=f$.
\end{proof}

A consequence of \ref{tame-feral4} is that there
are at most three inequivalent expansions 
of $\PP$ at $P_2$ where the petal
of the expansion adjacent to $P_2$ has two elements. Combining this 
with \ref{tame-feral3} we deduce that there 
are at most four inequivalent
expansions of $\PP$ at $P_2$ by a feral element in $B_v$. If 
we consider expansions at $P_m$ we obtain at most eight inequivalent
expansions altogether. The lemma follows from this fact and
\ref{tame-feral1}.
\end{proof}

Let $f$ be a feral element of the $k$-coherent matroid $M$.
By Theorem~\ref{feral}, $f$ has a feral display in either 
$M$ or $M^*$. If the former case holds we say that
$f$ is a {\em standard}
\index{standard feral element}
feral element and if the later case holds we say that
$f$ is a {\em costandard} feral element.
\index{costandard feral element}

\begin{corollary}
\label{bound-feral}
Let $v$ be bag vertex of a $3$-tree for the $k$-skeleton $M$
in $\eq$. If $v$ has degree $l$, then there are at most
$16l(q^2+q+1)$ feral elements in $B_v$.
\end{corollary}

\begin{proof}
Say that $v$ is not strong. If $d(v)\leq 1$, then the elements of
$B_v$ are peripheral and by Lemma~\ref{feral-internal},
$B_v$ contains no feral elements.
Say $d(v)=2$. Let $X$ and $Y$ be the $3$-separating
sets displayed by $v$. If $(X,B_v\cup Y)$ is equivalent to 
$(X\cup B_v,Y)$, then the elements of $B_v$ are not feral
as they are either in the guts or coguts of a $3$-separation.
On the other hand, if $(X,B_v,Y)$ is a flower and $B_v$ is not
displayed in a 3-tree, then $B_v$ is sequential
and the elements of $B_v$ are not peripheral  unless possibly $|B_v|=2$,
in which case the bound of the lemma certainly holds.

Assume that $v$ is strong. Let $f$ be a feral element in $v$.
Assume that $f$ is standard. Let $\PP=(P_1,P_2\ldots,P_l)$
be a maximal $k$-fracture of $M\ba f$. Then, using the definition
of feral element, there is an $i\in\{2,3,\ldots,l\}$ such that
$(P_2,P_3,\ldots,P_i,P_{i+1}\udots P_l\cup P_1\cup\{f\})$
and $(P_{i+1},P_{i+2},\ldots,P_i,P_1\cup P_2\udots P_i\cup\{f\})$ 
are maximal swirl-like flowers of $M$, 
one of which may have only two petals.
We lose no generality in assuming that $i>2$. In this case
the flower $\PP$ is an expansion of 
$(P_2,P_3,\ldots,P_i,P_{i+1}\udots P_l\cup P_1\cup\{f\})$ by $\{f\}$.
It is readily seen that this flower is  displayed by a 
vertex adjacent to $v$. There are at most $l$ such flowers and, by 
Lemma~\ref{tame-feral}, each such flower is expanded by at most
$8(q^2+q+1)$ feral elements. Hence there are at most 
$8l(q^2+q+1)$ standard feral elements By duality there are at most
$8l(q^2+q+1)$ costandard feral elements and the corollary follows.
\end{proof}

\begin{lemma}
\label{seq-unfetter}
Let $A$ be a sequential $3$-separating set of the $k$-skeleton
$M$. If $|A|\geq 4$, then no element of $A$ is fettered.
\end{lemma}

\begin{proof}
Let $B=E(M)-A$. Clearly we may assume that $A$ is fully closed.
Let $(a_1,a_2,\ldots,a_n)$ be a sequential ordering of $A$.
Up to duality we may assume that $\{a_1,a_2,a_3\}$
is a triangle. 

\begin{sublemma}
\label{seq-unfetter1}
$\{a_1,a_2,a_3\}$ is a clonal triple and $M\ba a_1$ is a $k$-skeleton.
\end{sublemma}

\subproof
As $M$ is a $k$-skeleton $\{a_1,a_2,a_3\}$ is not in a 
4-element fan, so by Tutte's Triangle Lemma, $M\ba a_i$
is $3$-connected for some $i\in \{1,2,3\}$. By
Corollary~\ref{remove-sequential}, $M\ba a_i$ is $k$-coherent.
Hence $\{a_1,a_2,a_3\}$ is not $k$-wild and is therefore a 
clonal triple. By Corollary~\ref{stuck} 
$M\ba a_1$ is a $k$-skeleton.
\end{proof}

\begin{sublemma}
\label{seq-unfetter2}
The lemma holds if $|A|=4$.
\end{sublemma}

\subproof
Say $A=\{a_1,a_2,a_3,a_4\}$. By \ref{seq-unfetter1},
$a_1$, $a_2$ and $a_3$ are unfettered. Consider $a_4$.
If $A$ is a line, then the claim
holds as $A$ is a clonal set. In the other case
$a_4\in\cl^*(\{a_1,a_2,a_3\})$. Thus $\co(M\ba a_4)$
is not $3$-connected. By Bixby's Lemma and the fact that
$M$ has no 4-element fans, $\si(M/a_4)$ is $3$-connected.
By Corollary~\ref{remove-sequential}, $M/a_4$ is $k$-coherent.
Thus $a_4$ is not cofixed in $M$. 

It remains to show that $a_4$ is not fixed in $M$. 
Note that $\{a_2,a_3,a_4\}$ is a triad of $M\ba a_1$.
By
Corollary~\ref{born-free}(i), $a_4$ is not cofixed in $M\ba a_1$.
Hence $\{a_2,a_3,a_4\}$ is not a $k$-wild triad in $M\ba a_1$.
As $M\ba a_1$ is a $k$-skeleton $\{a_2,a_3,a_4\}$ is a 
clonal triple of $M\ba a_1$. Thus $a_4$ is not
fixed in $M\ba a_1$. By Lemma~\ref{freedom21},
$a_4$ is not fixed in $M$.
\end{proof}

Assume that $n>4$ and, for induction, that the lemma holds 
if $|A|=n-1$. Then all the elements of $A-\{a_1\}$ are unfettered
in $M\ba a_1$. By Lemma~\ref{freedom21}, they are unfettered in 
$M$.
\end{proof}

\begin{proof}[Proof of Theorem~\ref{control-fetter}]
Let $T$ be a $3$-tree for $M$. Assume that $T$ has at most $\mu$
vertices. It is possible that a fettered element $z$ of
$M$ is in the guts or
coguts of a $3$-separation. Consider this case now.
The element $z$ may be in a $k$-wild triangle
or triad. As a generous bound for the number of such elements we
observe that

\begin{sublemma}
\label{control-wild}
There are at most $3\nu$ elements of $M$ in $k$-wild triads or
triangles.
\end{sublemma}

Assume that $z$ is not in a $k$-wild triangle or triad.
By Lemma~\ref{seq-unfetter}, $z$ is in the guts or
coguts of a non-sequential 3-separation $(X,Y)$. Note that there
is a path $(\coh(X),z_1,z_2,\ldots,z_m,\coh(Y)$ of
$3$-separations in $M$. By Theorem~\ref{main1}, there are at
most $f_{\ref{main1}}(k,q)$ elements of $M$ in such a path.
Moreover, a 3-separation equivalent to $(X,Y)$ is displayed by $T$.
Non sequential 3-separations are displayed by either flower vertices
or edges of $T$. There are at most $2^\nu$ 3-separations displayed in 
such a way. Putting this information together we get another generous
bound.

\begin{sublemma}
\label{control-guts}
There are at most $2^\nu f_{\ref{main1}}(k,q)$ fettered elements in the
guts or coguts of a non-sequential $3$-separation in $M$.
\end{sublemma}

If $z$ is a fettered element not covered by the previous cases,
then $M\ba z$ and $M/z$ are both 3-connected. 
By Corollary~\ref{bound-feral}
there are at most $16\nu(q^2+q+1)$ feral elements in any vertex bag of
$T$. It follows that

\begin{sublemma}
\label{control-feral}
There are at most $16\nu^2(q^2+q+1)$ feral elements in $M$.
\end{sublemma}

Any fettered element not covered by the previous cases
has the property that either
(a) $M\ba z$ is $k$-coherent and $M/z$ is $3$-connected
and $k$-fractured, or (b) $M/z$ is $k$-coherent and 
$M\ba z$ is $3$-connected and $k$-fractured. 
Assume that $z$ has property (a). By 
Corollary~\ref{fetter5}, there is a swirl-like flower
$(P_1,P_2,\ldots,P_{k-1})$ of $M$, where $x\in P_1$ such that,
for some partition $(P'_1,P''_1)$ of $P_1-\{z\}$,
the partition $(P'_1,P''_1,P_2,\ldots,P_{k-1})$ is a maximal 
$k$-fracture of $M/z$. By Lemma~\ref{fetter6}, there are at
most $\nu(q+1)$ elements of this type. Dually, there are
at most $\nu(q+1)$ fettered elements of type (b).

From the above fact, \ref{control-wild}, \ref{control-guts} and
\ref{control-feral}, we deduce that the theorem holds by letting
$$f_{\ref{control-fetter}}(\nu,k,q)=3\nu+2^\nu f_{\ref{main1}}(k,q)
+16\nu^2(q^2+q+1)+2\nu(q+1).$$
\end{proof}

\section{Finding an Unfettered Minor}

We are now in a position to simplify structure by obtaining
a large 4-connected unfettered matroid from a sufficiently
large $k$-skeleton.

\begin{theorem}
\label{unfetter1}
There is a function $f_{\ref{unfetter1}}(m,\nu,k,q)$ such that,
if $M$ is a $k$-skeleton in $\eq$
having a $3$-tree with at most $\nu$
vertices, then $M$ has a $4$-connected unfettered minor with at
least $m$ elements.
\end{theorem}

Our strategy will be to move to find a bounded size $k'$
such that a sufficiently large $k$-skeleton has a large
unfettered $k$-coherent minor.

\begin{lemma}
\label{unfetter4}
Let $M$ be a $k$-coherent matroid that is not a wheel or a whirl.
If $M$ has $l$ fettered elements, then there is an element $z$
of $M$ such that either $M\ba z$ or $M/z$ is $2k$-coherent with 
at most $(l-1)$ fettered elements.
\end{lemma}

\begin{proof}
Let $x$ be a fixed element in $M$. 
Assume that  the element $z$ is fettered 
in $M\ba x$. If $z$ is fixed in $M\ba x$, 
then $z$ is clearly fixed in $M$. Say
that $z$ is cofixed in $M\ba x$. 
Then, by Corollary~\ref{born-free} and the
fact that is $x$ is fixed in $M$, we see $z$
is cofixed in $M$. From this we deduce

\begin{sublemma}
\label{unfetter4.1}
If $x$ is fixed in $M$, and $y$ is fettered in $M\ba x$, then $y$
is fettered in $M$.
\end{sublemma}

Let $x$ be a fettered element of $M$. If $x$ is in a 4-element fan,
then, as $M$ is not a wheel or a whirl, this fan has an end $z$. 
Moreover, $z$ is either fixed and $M\ba z$ is
$k$-coherent, or cofixed and $M/z$ is $k$-coherent.
By \ref{unfetter4.1}, the lemma is satisfied with this choice
of $z$. We may thus assume that $M$ has no 4-element fans.

Assume that $x$ is in a triangle $T$. Certainly $T$ is not a clonal
triple so it has a fixed element $z$. If $T$ has an unfixed element,
then by Corollary~\ref{free-triangle}, $M\ba z$ is $k$-coherent and
the lemma holds by \ref{unfetter4.1}. Assume that all elements of
$T$ are fixed. Then, as $M$ has no $4$-element fans, 
$T$ has an element $z$ such that $M\ba z$ is $3$-connected. 
It remains to show that $M\ba z$ is $2k$-coherent. Let
$(P_1,P_2,\ldots,P_n)$ be a maximal swirl-like flower in $M\ba z$.
Assume that the other elements of $T$ are in $P_1$ and $P_i$.
Then 
$(P_2,P_3,\ldots,P_{i-1},P_i\cup P_{i+1}\udots P_n\cup P_1\cup\{z\})$
and 
$(P_1\cup P_2\udots P_i\cup\{z\},P_{i+1},P_{i+2}\ldots,P_n)$ 
are swirl-like
flowers in $M$ of order $i-1$ and $n-i+1$ respectively. Thus
$i-1<k$ and $n-i+1<k$ so that $n<2k$ and $M\ba z$ is $2k$-coherent.

We may now assume that $x$ is not in a triangle or a triad. Assume
that $x$ is in the guts of a vertical 3-separation. Then
$M\ba x$ is 3-connected and there is a 3-separation $(A,B)$
of $M\ba x$ such that $x\in\cl(A)$ and $x\in\cl(B)$.
Let $\PP=(P_1,P_2,\ldots,P_n)$ be a maximal swirl like flower of 
$M\ba x$. By Lemma~\ref{cross1}, up to labels and flowers equivalent
to $\PP$ there is an $i\in\{1,2,\ldots,n\}$ such
that either $A\subseteq P_i$ or 
$(A,B)=(P_1\cup P_2\udots P_i,P_{i+1}\cup P_{i+2}\udots P_n)$.
As $x\in\cl(A),\cl(B)$ we again deduce that $M\ba x$ is $2k$-coherent.

If $x$ is fixed in $M$, then, by the \ref{unfetter4.1}, $M\ba x$
has at most $l-1$ fettered elements, but it may be the case that
$x$ is not fixed but cofixed. We consider this case now. Assume that the
element $z$ is fettered in $M\ba x$, but not in $M$. Then, by 
Corollary~\ref{born-free}, $x\more z$.
As $x$ is cofixed in $M$ we see that $x\not\cong z$. Since $x\in\cl(A)$
and $x\in\cl(B)$, it is also the case that
$z\in\cl(A)$ and $z\in\cl(B)$. Therefore 
$z$ is in the guts of the $3$-separation $(A,B)$. Thus $z$ is fixed in
$M$. In this case $M\ba z$ is $2k$-coherent with at most $l-1$
fettered elements.

We may now assume that both $M\ba x$ and $M/x$ are 3-connected. 
Assume that
$x$ is fixed in $M$. Then the lemma holds for $M\ba x$ unless 
$M\ba x$ is not $2k$-coherent. In this case, by 
Corollary~\ref{fetter5},
$x$ is cofixed. If $M/x$ is $k$-coherent, then the lemma
holds for $M/x$. Otherwise $x$ is a feral element and has a 
feral display in $M$ or $M^*$. 
It follows immediately from the definition of a feral display 
that $M/x$ is $2k$-coherent,
and the lemma holds in this case too.
\end{proof}

As an immediate consequence of Lemma~\ref{unfetter4} we obtain

\begin{corollary}
\label{unfetter5}
Let $M$ be a $k$-coherent matroid with $l$ fettered elements and
a nonempty set of unfettered elements, then $M$ has an
unfettered $2^lk$-coherent minor with at least $|E(M)-l|$ elements.
\end{corollary}

\begin{lemma}
\label{unfetter21}
Let $(A,B)$ be a $3$-separation of the matroid $M$, where $A$ is fully
closed and let $a_1$ and $a_2$ be elements of $A$. Let $N$ be a 
$3$-connected minor of $M$ on $B\cup\{a_1,a_2\}$. If the element
$b$ of $B$ is unfettered in $M$, then $b$ is unfettered in $N$.
\end{lemma}

\begin{proof}
Say $a\in A$ and $N$ is a minor of $M/a$. Then $a\notin\cl(B)$, as
otherwise $N$ is not $3$-connected. As $A$ is fully closed,
$b\notin\cl^*(A)$, so there is a circuit $C\subseteq B$ containing $b$.
As $a\notin\cl(B)$ we have $a\notin\cl(C)$. Therefore it is not the
case that $a\less b$. Thus, by Corollary~\ref{born-free}, 
$b$ is not fixed in
$M/a$. Also $b$ is not cofixed in $M/a$ as $b$ is not cofixed in $M$.
The lemma now follows by duality and induction.
\end{proof}

An unfettered $k$-coherent matroid is evidently a $k$-skeleton. We
use this fact without comment from now on.

\begin{lemma}
\label{unfetter22}
Let $M$ be a $k$-coherent unfettered matroid in $\eq$,
and let $v$ be a strong
bag vertex of a $3$-tree for $M$. If 
$|B_v|\geq mf_{\ref{main1}}(k,q)$, then $M$ has a 
$4$-connected unfettered minor with at least $m$ elements.
\end{lemma}

\begin{proof}
By Lemma~\ref{potato6}, the eyes of $v$ form a potato $\PP$ of $M$.
Say $\PP$ has $l$ eyes. Note that $M$ has no $k$-wild triangles
or triads so that an obvious minor perturbation of 
Lemmas~\ref{potato4} and \ref{potato5} proves that 
$M$ has a $4$-connected minor $N$ containing the core of $\PP$ and
a set of $2l$ other elements. Moreover, by Lemma~\ref{unfetter21}
this matroid is unfettered. If $l\geq m/2$, then the lemma follows.
On the other hand say $l<m/2$. By Lemma~\ref{lazy-bound},
each eye contains at most $f_{\ref{main1}}(k,q)$
elements of $B_v$. As $f_{\ref{main1}}(k,q)\geq 2$, we see that
in this case the core of $\PP$ has at least 
$m$ elements. Thus the lemma holds in either case.
\end{proof}

\begin{proof}[Proof of Theorem~\ref{unfetter1}]
Let $k'=2^{f_{\ref{control-fetter}}(\nu,k,q)}k$.
Let $\nu'=f_{\ref{small-tree}}(m,\nu,k',q)$.
Let $f_{\ref{unfetter1}}(m,k,q)=(mf_{\ref{main1}}(k',q))^{\nu'}
+f_{\ref{control-fetter}}(\nu,k,q)$.

Assume that $M$ as at least $f_{\ref{unfetter1}}(m,k,q)$ elements.
By Theorem~\ref{control-fetter}, $M$ has at most
$f_{\ref{control-fetter}}(\nu,k,q)$ fettered elements.
By Corollary~\ref{unfetter5}, $M$ has an unfettered
$k'$-coherent minor $N$ with at least 
$(mf_{\ref{main1}}(k',q))^{\nu'}$ elements.
Let $T'$ be a 3-tree for $N$. If $T'$ has at least
$\nu'$ vertices, then, by Theorem~\ref{small-tree}, $N$
has a $4$-connected minor with a set of $m$ pairwise disjoint clonal
pairs. It is easily deduced from the construction of the minor
and Lemma~\ref{unfetter21}, that the minor is unfettered,
so the theorem holds in this case. Thus we may assume that
$T'$ has at most $\nu'$ vertices. Then $T$ has a vertex $v$ for
which $|B_v|\geq mf_{\ref{main1}}(k',q)$. We may assume that 
$m\geq 4$, so that $v$ is strong. The theorem now follows from
Lemma~\ref{unfetter22}.
\end{proof}

\section{Unfettered Matroids}

We now know that a sufficiently large $k$-skeleton in 
$\eq$ has a large $4$-connected unfettered minor. 
The next task is to  extract
a $4$-connected minor having many clonal pairs from such a matroid.
While we only need the result for matroids in $\eq$, free
spikes play no role in the arguments, so we focus on the
class of matroids with no $U_{2,q+2}$- or $U_{q,q+2}$-minor,
that is, the class $\youq$.

\begin{theorem}
\label{4clonal1}
There exists a function $f_{\ref{4clonal1}}(q,t)$
such that, if $M$ is a $4$-connected unfettered matroid in 
$\youq$ with at least $f_{\ref{4clonal1}}(q,t)$ elements, 
then $M$ has a $4$-connected minor with at least $n$ 
pairwise disjoint clonal pairs.
\end{theorem}

As usual we prepare for the proof with a series of lemmas.
We begin with ones that focus on connectivity.
A matroid $M$ is {\em $4$-connected up to rank-$k$ $3$-separators},
\index{$4$-connected up to rank-$3$ separators}
if it is 3-connected and, whenever $(A,B)$ is a 3-separation of $M$,
either $r(A)\leq k$ or $r(B)\leq k$. A vertically
$4$-connected matroid is one that is $4$-connected up to 
rank-2 3-separators. If $M$ is $4$-connected up to
rank-3 3-separators we will say that $M$ is 
{\em $4$-connected up to planes}.
\index{$4$-connected up to planes}

If $M$ is $3$-connected and $(A,B)$ is a 3-separation where 
$r(A)=3$, then $\lambda(A-\cl(B))=2$, $r(A-\cl(B))=3$ and
$A-\cl(B)$ is a cocircuit. It follows that, if $M$ is 4-connected
up to planes, then $M$ is vertically $4$-connected if and
only if $M$ has no 3-separating cocircuits. 
The next lemma is clear.

\begin{lemma}
\label{4clonal2}
Let $x$ be an element of the $3$-connected matroid $M$. If
$M\ba x$ is $4$-connected up to rank-$k$ $3$-separators, then $M$
is also $4$-connected up to rank-$k$ $3$-separators.
\end{lemma}

\begin{lemma}
\label{4clonal3}
Let  $M$ be an unfettered matroid that is $4$-connected up to planes. 
If $r(M)\geq 6$, then the 
$3$-separating cocircuits of $M$ 
partition a subset of the ground set of $M$.
\end{lemma}

\begin{proof}
Assume that $C_1$ and $C_2$ are 3-separating cocircuits of
$M$. Say $x\in C_1\cap C_2$. As $x$ is not fixed, there
is a matroid $M'$ obtained
by independently cloning $x$ by $x'$.
Then $x'\notin \cl_{M'}(E(M)-C_1)$, as otherwise $x$ and $x'$
are not clones. Hence $C_1\cup\{x'\}$ and $C_2\cup\{x'\}$
are 3-separating cocircuits of $M'$. 
Now $r(C_1\cup\{x'\})\cup (C_2\cup\{x'\})|\geq 2$, so that
$\lambda_{M'}(C_1\cup C_2\cup\{x\})=2$ and, by an easy rank
argument, we also have $r_{M'}(C_1\cup C_2\cup\{x'\})=4$.
Thus
$r_M(C_1\cup C_2)=4$ and $\lambda_M(C_1\cup C_2)=2$.
As $M$ has rank at least six, 
$r_M(E(M)-(C_1\cup C_2))\geq 4$ and we have contradicted
the fact that $M$ is 4-connected up to planes.
\end{proof}



The next lemma is a consequence of Lemma~\ref{freedom3}.

\begin{lemma}
\label{4clonal-triangle}
If $T$ is a triangle of an unfettered matroid, then the members
of $T$ are clones.
\end{lemma}

\begin{lemma}
\label{4clonal4}
Let $M$ be a vertically $4$-connected unfettered matroid. If 
the element $x$ of $M$
does not have a clone, then both $M\ba x$ and $M/x$ are 
$3$-connected and unfettered.
\end{lemma}

\begin{proof}
Assume that $x$ has no clone.
Clearly $M\ba x$ is $3$-connected and, by Lemma~\ref{4clonal-triangle}
so too is $M/x$.
Certainly no element of $M\ba x$ is fixed. Assume that the
element $z$ is cofixed in $M\ba x$. Then, by Corollary~\ref{born-free},
either $x$ is cofixed in $M$, contradicting the fact that
$M$ is unfettered, or $x$ and $z$ are clones, contradicting the fact
that $x$ does not have a clone. Thus $M\ba x$ is 
unfettered and dually, so too is $M/x$.
\end{proof}

\begin{lemma}
\label{4clonal8}
Let $t$ be an element of a triangle of a vertically $4$-connected
unfettered matroid. Then $M\ba t$ is vertically $4$-connected
and unfettered.
\end{lemma}

\begin{proof}
Assume that $t$ belongs to the triangle $T$ of $M$.
By Lemma~\ref{4clonal-triangle}, 
$T$ is a clonal triple of $M$ and it is easily seen
that $M\ba t$ is vertically $4$-connected. No element of
$M\ba t$ is fixed in this matroid. It also follows routinely
from Corollary~\ref{born-free} that $M\ba t$ has no cofixed elements.
\end{proof}

\begin{lemma}
\label{4clonal5}
Let $M$ be an unfettered matroid that is 
$4$-connected up to planes, and 
let $A$ be a $3$-separating cocircuit of $M$. Then the following hold.
\begin{itemize}
\item[(i)] If $l$ is a non-trivial line that meets $A$, then 
$l\subseteq A$.
\item[(ii)] If $a\in A$ and $a$ has no clone, then $M/a$ is 
$4$-connected up to planes.
\item[(iii)] If the element $a$ of $A$ 
is in a non-trivial line $l$ of $M$ and
$|A-l|\geq 2$, then $\si(M/a)$ is $4$-connected up to planes.
\end{itemize}
\end{lemma}

\begin{proof}
Let $l$ be a non-trivial line of $M$
that meets $A$. Assume that
$l$ is not contained in $A$. Then $|(E(M)-A)\cap l|=1$; say 
$(E(M)-A)\cap l=\{x\}$. As $x$ is in the guts of the
$3$-separation $(E(M)-A,A)$,
we have $x\in\cl(E(M)-(A\cup\{x\}))$. But then 
$x$ is fixed in $M$ by Lemma~\ref{freedom3}, contradicting the
fact that $M$ is unfettered. Thus (i) holds.

Consider (ii). The result is routine if $A$ is not a triad.
Assume that $A$ is a triad $\{a,b,c\}$. By the dual of 
Lemma~\ref{4clonal-triangle}, $\{a,b,c\}$ is a clonal triple
so that $\{b,c\}$ is a clonal pair in $M/a$. Assume that
$(B,C)$ is a $3$-separation of $M/a$. If 
$\{b,c\}\subseteq B$, then $(B\cup\{a\},C)$ is a $3$-separation
of $M$, and it follows that either $B$ or $C$ has rank at most
three in $M/a$. Assume that $b\in B$ and $c\in C$. If
there is a circuit $Z$ contained in $B$ that contains $c$,
then $c\in\cl_{M/a}(Z)$, and we may apply the previous
argument to conclude that either $B$ or $C$ has rank at most
$3$ in $M/a$. Otherwise, both $b$ and $c$ are in the coguts of
the $3$-separation $(B,C)$. Consider this case.
Assume for a contradiction that
$r_{M/a}(B)\geq 4$ and $r_{M/a}(C)\geq 4$. Then
$(B-\{b\},C\cup\{b\})$ is a $3$-separation in $M/a$
and $b\in\cl^*_{M/a}(B-\{b\})$. As $\{a,b,c\}$ is a triad
of $M$ we have $a\in\cl^*_M(C\cup\{b\})$ so that
$(B-\{b\},C\cup\{a,b\})$ is a 3-separation in $M$.
But $b\in\cl^*_{M/a}(B-\{b\})$, so
$b\in\cl^*_M(B-\{b\})$. Hence $(B,C\cup\{a\})$ is a 3-separation of
$M$. But $r_M(B),r_M(C\cup\{a\})\geq 4$ and we have contradicted
the fact that $M$ is 4-connected up to planes.

We omit the easy proof of (iii).
\end{proof}

A 3-separation $(A,B)$ of a matroid $M$
is {\em cyclic}
\index{cyclic $3$-separation} 
if both $A$ and $B$
contain circuits of $M$. The matroid $M$ is {\em cyclically
$4$-connected}
\index{cyclically $4$-connected}
if it 3-connected and
has no cyclic 3-separations. Note that
$M$ is cyclically 4-connected if and only if $M^*$ is vertically
4-connected. Note also that cyclically 4-connected matroids
do not contain triangles unless the matroid
is degenerately small. 

\begin{lemma}
\label{4clonal9}
Let $p$ be an element of the plane $P$ of the
cyclically $4$-connected unfettered matroid $M$ where
$|P|\geq 5$. Assume
that $M\ba p$ is not cyclically $4$-connected. Then $|P|=5$ and
$M\ba z$ is cyclically $4$-connected for all $z\in P-\{p\}$.
\end{lemma}

\begin{proof}
Let $(A,B)$ be a cyclic 3-separation of $M\ba p$, where
$B$ is coclosed. If $|A\cap(P-\{p\})|\leq 1$, then
$(A,B\cup\{p\})$ is a cyclic 3-separation of $M$.
Thus $|A\cap(P-\{p\})|> 1$ and symmetrically
$|B\cap (P-\{p\})|> 1$. But $r(A\cap(P-\{p\})\leq 2$
and $r(B\cap (P-\{p\})\leq 2$ as otherwise $p$ is in the
closure of either $A$ or $B$ implying that $M$ is not cyclically
$4$-connected. As $M$ has no triangles we deduce that
$|P|=5$ and that  $|A\cap(P-\{p\})|=|B\cap (P-\{p\})|=2$.
Let $\{a_1,a_2\}=A\cap(P-\{p\})$ and let
$\{b_1,b_2\}=B\cap(P-\{p\})$.
Observe that $A$ and $\cl(P)$ are a modular pair of cyclic flats of
$M$ whose intersection does not span $P$. 
Thus any clone of 
$a_1$ or $a_2$ lies on the line $\{a_1,a_2\}$, and by 
symmetry any clone of $b_1$ or $b_2$ lies on the line
$\{b_1,b_2\}$. It follows that 
$\{a_1,a_2\}$ and $\{b_1,b_2\}$ are clonal lines of $M$.

Assume that $M\ba a_1$ is not cyclically 4-connected.
Then, arguing as above, we deduce that, 
for some $x\in\{b_1,b_2,p\}$, the pair
$\{a_2,x\}$ is a clonal line of $M$  contradicting the fact that 
$\{a_1,a_2\}$ is a clonal line of $M$.
\end{proof}

The proof of the next lemma is entirely analogous to the proof
of Bixby's Lemma.

\begin{lemma}
\label{4clonal10}
Let $x$ be an element of the vertically $4$-connected matroid
$M$ that is not in a triangle. Assume that $M/x$ is not
$4$-connected up to planes. Then $M\ba x$ is $4$-connected
up to rank-$4$ $3$-separators. Moreover, if $(D_1,D_2)$ is a vertical
$4$-separation of $M\ba x$ such that $r(D_2)\leq 4$,
then $D_2$ is covered by a pair of lines of $M$.
\end{lemma}

\begin{proof}
Assume that $M/x$ is not 4-connected up to planes. Then there
is a 4-separation $(C_1\cup\{x\},C_2)$ of $M$
with $x$ in the guts such that $r(C_1),r(C_2)\geq 5$.

Let $(D_1,D_2)$ be a $3$-separation of $M\ba x$. We
need only consider the case that  
$r(D_1),r(D_2)\geq 3$. Then, as neither $(D_1\cup\{x\},D_2)$
nor $(D_1,D_2\cup\{x\})$ is a $3$-separation of $M$, we see
that $(D_1\cup\{x\},D_2)$ is a $4$-separation of $M$ with
$x$ in the coguts. 

Without loss of generality $r(C_1\cap D_1)\geq 3$.
Thus $\lambda(C_1\cap D_1)\geq 3$. By uncrossing
$\lambda((C_2\cap D_2)\cup\{x\})\leq 3$. But
$x\in\cl(C_1\cup D_1)$ and $x\in\cl^*(C_1\cup D_1)$.
Hence $r(C_1\cup D_1\cup\{x\})=r(C_1\cup D_1)$
and $r(C_2\cap D_2)=r((C_2\cap D_2)\cup\{x\})-1$.
Therefore $\lambda(C_2\cap D_2)\leq 2$.
As $M$ is vertically 4-connected, we have $r(C_2\cap D_2)\leq 2$.
But $r(C_2)\geq 5$ so that $r(C_2\cap D_1)\geq 3$.
Repeating the above argument shows that $r(C_1\cap D_2)\leq 2$.
Therefore $r(D_2)\leq 4$ and $D_2$ is covered by a pair of lines
of $M$.
\end{proof}

Let $S$ be a set of elements of the matroid $M$. Define $P_M(S)$
to be the set $\cl(\cl^*(S))$.

\begin{lemma}
\label{4clonal7}
If $z\notin P_M(S)$, then $P_{M\ba z}(S)\supseteq P_M(S)$.
\end{lemma}

\begin{proof}
Observe that $\cl^*_{M\ba z}(S)=\cl^*_M(S\cup\{z\})-\{z\}$.
It follows that $\cl^*_{M\ba z}(S)\supseteq\cl^*_M(S)$.
As $z\notin\cl_M(\cl^*_M(S))$, it follows that
$\cl_M(\cl^*_M(S))\subseteq\cl_{M\ba z}(\cl^*_{M\ba z}(S))$
as required.
\end{proof}

Recall that it follows from Lemma~\ref{kung} that a simple 
rank-$t$ matroid
in $\youq$ has at most $(q^t-1)/(t-1)$ elements.
Let $h(q,t)=(q^t-1)/(t-1)$. By duality, if $M$ is a cosimple matroid
in $\youq$ of corank $t$, then $M$ has at most $h(q,t)$ elements.

\begin{lemma}
\label{kung3}
Let $M$ be a simple, cosimple matroid in  $\youq$.
\begin{itemize}
\item[(i)] If $S\subseteq E(M)$ and $|S|\leq t$, then
$|P_M(S)|\leq h(q,h(q,t))$.
\item[(ii)] If $r^*(M)\geq h(q,t)$, then $r(M)\geq t$.
\end{itemize}
\end{lemma}

\begin{proof}
Consider (i). Assume that 
$|S|\leq t$. As $r^*(S)\leq |S|$, and $r^*(\cl^*(S))=r^*(S)$,
we see that $|\cl^*(S)|\leq h(q,t)$. But 
$r(\cl(\cl^*(S))=r(\cl^*(S))$ so that $|P_M(S)|\leq h(q,h(q,t))$.

Consider (ii). Assume that $r^*(M)\geq h(q,t)$. Then
$|E(M)|\geq h(q,t)$ so that, by Lemma~\ref{kung}, $r(M)\geq t$.
\end{proof}

For an unfettered matroid $M$, let $L_M$ denote the set of 
elements in non-trivial clonal classes of $M$. Not surprisingly,
a {\em clonal line}
\index{clonal line} 
of $M$ is a line whose elements form a 
clonal set. Of course, if the line is non-trivial its elements will
form a clonal class.

\begin{lemma}
\label{4clonal12}
Let $M$ be a vertically $4$-connected unfettered matroid
in $\youq$. Then there is a function
$f_{\ref{4clonal12}}(q,t)$ such that such that if 
$r^*(M)\geq f_{\ref{4clonal12}}(q,t)$, 
then at least one of the following holds.
\begin{itemize}
\item[(i)] $M$ has a vertically $4$-connected restriction
$N$ such that $|L_N|\geq t$.
\item[(ii)] There is an element $z\in E(M)-P_M(L_M)$ such that
$M/z$ is $4$-connected up to planes.
\end{itemize}
\end{lemma}

\begin{proof}
To simplify notation we inductively define $h^i(q,t)$ for 
$i\geq 2$ by $h^i(q,t)=h(h(q,h^{i-1}(q,t))$.
Assume that $r^*(M)>h^6(q,t)$.
For a set $S$ of elements of a matroid $N$, set
$Q_N(S)=P_N(P_N(S))$. Let $Z=(z_1,\ldots,z_{l-1})$
be a maximal sequence of elements of $M$ such that, 
for $i\in\{1,\ldots,l-1\}$, the following hold.
\begin{itemize}
\item[(i)]  $M\ba z_1,z_2,\ldots,z_i$ is vertically 4-connected.
\item[(ii)] $z_i\notin L_{N\ba z_1,z_2,\ldots,z_i}$.
\end{itemize}

Clearly $r(M\ba Z)=r(M)$.
By Lemma~\ref{4clonal4}, $M\ba Z$ is unfettered.
If $|L_{M\ba Z}|\geq t$, then the lemma holds. 
Assume that $|L_{M\ba Z}|<t$. By Lemma~\ref{kung3}(ii),
$r(M)> h^5(q,t)$. As $r(M\ba Z)=r(M)$, it also
follows from Lemma~\ref{kung3}(ii) that
$r^*(M\ba Z)>h^4(q,t)$.
By Lemma~\ref{kung3}(i),
there is an element $z_l\in E(M\ba Z)-Q_{M\ba Z}(L_M)$.

If $M\ba Z/z_l$ is
$4$-connected up to planes, then, by Lemma~\ref{4clonal2},
$M/z_l$ is $4$-connected up to planes and the lemma holds.
Assume otherwise. Then the conclusions of Lemma~\ref{4clonal10}
hold for $M\ba Z\ba z_l$. Let $M'=M\ba Z\ba z_l$.
It follows from 
the choice of $Z$ that $M'$ is not vertically $4$-connected.
Moreover, if $A$ is a 3-separator of $M'$ with $r(A)\in\{3,4\}$,
then $z_l\in\cl^*_{M\ba Z}(A)$.

Assume that $M'$ is not 4-connected up to planes. 
Let $A$ be a rank-4 3-separator of $M'$. Then there are lines 
$l_1$ and $l_2$ such that $l_1\cup l_2=A$. As $r(A)=4$,
these lines are disjoint. If $A\subseteq P_{M\ba Z}(L_{M\ba Z})$,
then $z_l\in Q_{M\ba Z}(L_{M\ba Z})$, 
contradicting the choice of $z_l$.
Thus we may assume that there is an element $a\in l_1$
such that $a\notin P_{M\ba Z}(L_{M\ba Z})$. Certainly $a$ is not
in a triangle, so $|l_1|=2$. Let $L$ be a subset of $l_2$
with $|L|=|l_2|-2$. By Lemma~\ref{4clonal8}, $M\ba Z\ba L$
is vertically $4$-connected, so that $(M\ba Z\ba L)^*$ is 
cyclically 4-connected. Moreover, $\{z\}\cup l_1\cup (l_2-L)$
is a 5-point plane of this matroid. It now follows by the dual
of Lemma~\ref{4clonal9}, that, if $z\in l_1$, then
$M\ba Z\ba L/z$ is vertically $4$-connected. By Lemma~\ref{4clonal2},
$M/z$ is vertically 4-connected and by Lemma~\ref{4clonal7}, 
$z\notin P_M(L_M)$.

Assume that $M'$ is $4$-connected up to planes and let $A$ be a 
rank-3, 3-separating cocircuit of $M'$. Arguing as above we deduce
that there is an element $z\in A$ that is not in 
$P_{M\ba Z}(L_{M\ba Z})$.
Now $z$ is not in a triangle of $M'$, as otherwise, 
by Lemma~\ref{4clonal-triangle},
$z\in L_M$. So, by Lemma~\ref{4clonal5}, $M'/z$ is 4-connected
up to planes. Again it follows from Lemma~\ref{4clonal7} that
$z\notin P_M(L_M)$.

We conclude that the lemma holds by letting
$f_{\ref{4clonal12}}(q,t)=h^6(q,t)$.
\end{proof}

\begin{lemma}
\label{4clonal13}
Let $M$ be a vertically $4$-connected unfettered matroid
in $\youq$. Assume that $r^*(M)\geq f_{\ref{4clonal12}}(q,t)$.
Then at least one of the following holds.
\begin{itemize}
\item[(i)] $M$ has a vertically $4$-connected restriction
$N$ such that $|L_N|\geq t$.
\item[(ii)] There is a nonempty 
set  $J\subseteq(E(M)-L_M)$ such that
$M/J$ is unfettered and vertically $4$-connected.
\end{itemize}
\end{lemma}

\begin{proof}
Assume that the lemma fails. 
Then (i) does not hold, so,
by Lemma~\ref{4clonal12}, there is an element
$z\in E(M)-P(L_M)$ such that $M/z$ is 4-connected up to planes.
If $M/z$ is vertically 4-connected, then, by Lemma~\ref{4clonal4},
$M/z$ is unfettered
so that the lemma is satisfied with $J=\{z\}$. It follows that $M/z$
is not  vertically 4-connected. Let $\{Z_1,Z_2,\ldots,Z_m\}$ be the
3-separating cocircuits of $M/z$. Say $i\in\{1,2,\ldots,m\}$.
If $Z_i\subseteq P(L_M)$, then $z\in P(L_M)$. Thus there is an element
$z_i\in(Z_i-P(L_M))$. 

\begin{sublemma}
\label{4clonal13.1}
There is an $i\in\{1,2,\ldots,m\}$ such that every 
element of $Z_i-L_M$
is contained in a triangle in $M/z$.
\end{sublemma}

\subproof
Assume otherwise. For $i\in\{1,2,\ldots,m\}$, choose $z_i\in Z_i-L_M$
so that $z_i$ is not in a triangle in $M/z$. 
It is straightforwardly seen
that now the lemma holds with $J=\{z,z_1,z_2,\ldots,z_m\}$,
contradicting the assumption that the lemma fails. 
\end{proof}

Let $Z$ be a 3-separating cocircuit of $M/z$ satisfying \ref{4clonal13.1}
and let $H=E(M)-Z$.

\begin{sublemma}
\label{4clonal13.2}
$Z$ is coclosed in $M$ and $M/z$.
\end{sublemma}

\subproof
If $z\in\cl^*(Z)$, then $r(H-\{z\})=r(M)-2$, 
so that $(H-\{z\},Z\cup\{z\})$
is a vertical 3-separation of $M$, contradicting the fact that $M$ is 
vertically 4-connected. 
Thus $z\notin\cl^*(Z)$. Say that $t\in H-\{z\}$,
and $t\in\cl^*(Z)$. Then $t\in\cl^*_{M/z}(Z)$ and $Z\cup\{t\}$
is a rank-4 3-separation of $M/z$ contradicting the fact that $M/z$
is $4$-connected up to planes. Therefore $Z$ is coclosed in
$M$ and hence in $M/z$.
\end{proof}

\begin{sublemma}
\label{4clonal13.3}
$Z$ partitions into two clonal lines of $M/z$. One of these,
say $l'_z$, is a clonal line of $M$; the other, say $l_z$, spans a 
plane of $M$ and contains at least two points $u$ and $v$ that 
are not in $L_M$.
\end{sublemma}

\subproof
If $Z\subseteq L_M$, then $z\in\cl(L_M)$ contradicting the choice of
$z$. Thus there is a point $u\in Z-L_M$. By \ref{4clonal13.1},
$u$ is in a line $l_z$ of $M/z$, where $|l_z|\geq 3$. As $M/z$
is unfettered, it follows from Lemma~\ref{4clonal5}(i),
that $l_z\subseteq Z$. As $u\notin L_M$, we see that
$l_z$ is not a line of $M$. Hence $r_M(l_z)=3$, and $z\in\cl_M(l_z)$.

Say $y\in Z$. Assume that $y$ is not in a triangle in $M/z$.
By \ref{4clonal13.1}, $y\in L_m$ and hence has a clone in $M$
and therefore in $M/z$. Such a clone must be in $Z$.
Thus $Z$ partitions into clonal lines of $M/z$.

Let $l'_z$ be a clonal line of $M/z$ in $Z-l_z$.
Assume that $r_M(l'_z)=3$. Then $z\in\cl_M(l'_z)$, and
$z\in\cl_M(H-\{z\})$ (by \ref{4clonal13.2}). It follows routinely that
$z$ is fixed in $M$. Thus $l'_z$ is a line of $M$. If this line is not a
clonal line, then it has two points and is not contained in $L_M$, and 
we obtain a contradiction to \ref{4clonal13.1}. Thus $Z-l_z$
partitions into clonal lines of $M$.

If $u$ is the only element of $l_z$ not in $L_M$, the 
$u$ is the only element of $Z$ not in $L_M$. Thus 
$u\in\cl^*_{M/z}(L_M)$ so that $u\in\cl^*(L_M)$.
This gives the contradiction that $z\in P(L_M)$. 
Hence there is another point $v\in l_z$ that is not
in $L_M$.

Assume that there is more than one clonal line in $Z-l_z$. If these
lines are skew in $M$, then $z\in\cl(L_M)$, 
contradicting the choice of $z$.
Thus $Z-l_z$ spans a plane in $M$. If both $u$ and $v$ are on this 
plane, then any clone $u'$ of $u$ must lie on both this plane and the 
plane $l_z\cup\{z\}$ (certainly $u$ is not a coloop of this plane).
Thus $\{u',u,v\}$ is a triangle and we deduce that $\{u,v\}$ is a 
clonal line of $M$, contradicting the fact that these elements are 
not in $L_M$. 

It follows that we may assume that $u\notin\cl_M(Z-l_z)$ and it is
now easily checked that, in $M/u$, the set $Z\cup\{z\}$ contains a
$U_{3,5}$ restriction (indeed a $U_{3,6}$ restriction, but we do not 
need this). Thus, if $(P,Q)$ is a vertical 3-separation of $M/u$,
we may assume that $Z\cup\{z\}\subseteq P$. But then 
$u\notin \cl_M(Q)$, so this 3-separation is not coblocked by $u$,
so that $(P\cup\{u\},Q)$ is a vertical 3-separation of $M$,
contradicting the fact that $M$ is vertically $4$-connected.
Therefore $l'_z$ is the unique clonal line of $M\ba l_z$, and the
claim holds.
\end{proof}

\begin{sublemma}
\label{4clonal13.4}
$M/u$ and $M/v$ are $4$-connected up to planes.
\end{sublemma}

\subproof
Let $P=(l_z-\{u\})\cup\{z\}$. By Lemma~\ref{4clonal5}, 
$\si(M/u,z)$ is 4-connected up to planes 
apart from the single parallel class $P$. 
If $M/u$ is not 4-connected up to
planes, then there is a 3-separation $(X,Y)$ of $M/u,z$
with $r_{M/u,z}(X)=3$, such that $z\in\cl_{M/u}(X)$ and
$z\notin\cl_{M/u}(Y)$. We may assume that $X$ is a cocircuit of
$M/u,z$. If $P\subseteq \cl_{M/u,v}(Y)$, then $z\in\cl_{M/u}(Y)$,
so $P\subseteq X$. But, by Lemma~\ref{4clonal3}, the 3-separating
cocircuits of $M/z$ are all disjoint from $P$.
\end{proof}

\begin{sublemma}
\label{4clonal13.5}
We may assume that neither $u$ nor $v$ is in $P(M_L)$.
\end{sublemma}

\subproof
If all but one element $t$ of $l_z$ is in $P(M_L)$, then the
remaining element is clearly in the closure of $Z-\{t\}$, so that
$z\in\cl(Z-\{t\})$, and $z\in\cl(P(M_L))$, that is, $z\in P(M_L)$.
\end{proof}

If $M/u$ is vertically 4-connected, then the lemma holds with $J=\{u\}$.
Thus we may assume that neither $M/u$ nor $M/v$ is vertically
4-connected. We now have symmetry between $u,v$ and $z$ in
that there exist sets $U$ and $V$ with partitions $\{l_u,l'_u\}$
and $\{l_v,l'_v\}$ respectively such that the conclusions established
above hold with $(z,Z,l_z,l'_z)$ replaced by $(u,U,l_u,l'_u)$
or $(v,V,l_v,l'_v)$.

\begin{sublemma}
\label{4clonal13.6}
$l_z\cup\{z\}=l_u\cup\{u\}=l_v\cup\{v\}$.
\end{sublemma}

\subproof
We first show that $l_u\cap l_z\neq\emptyset$. Assume otherwise.
As $l_u$ contains a circuit containing $u$, 
we have $|l_u\cap Z|\geq 2$.
Thus there is an element of $l'_z$ contained in $l_u$. As $l'_z$
is a clonal line it follows that $l'_z\subseteq l_u$. 
Note that any other
element of $l_u$ is in $H$. Now $r(l'_z\cup\{u\})=3$, so 
$l'_z\cup\{u\}$ spans $l_u$. Let $p$ and $q$ be elements of
$l_u-(l'_z\cup\{u\})$ and let $M'$ be the matroid obtained by 
independently cloning $p$ by $p'$. Then 
$\{p,p',q\}\subseteq \cl(l'_z\cup\{u\})$. Also, by \ref{4clonal13.2},
$p\in\cl(H-\{p\})$ so that $p'\in\cl(H)$.
But $\sqcap(l'_z\cup\{u\},H)=2$, so $\{p,p',q\}$ is a triangle.
This shows that $l_u-(l'_z\cup\{u\})$ is a clonal line of $M$
contradicting the fact, established by \ref{4clonal13.3} and symmetry,
that this set has at least two elements not in $L_M$. Therefore 
$l_u\cap l_z\neq\emptyset$.

We may now assume that there is an element $p\in l_z\cap(l_u-\{u\})$.
Then $p$ is in $U$. But $l_u$ contains a circuit containing $p$
and $U$ is a cocircuit. Hence there is another element of $U$
contained in $l_z$. If such an element is contained in $l_u$, then,
as $u$ is not in a triangle, we see that 
$\cl(l_u\cup\{u\})=\cl(l_z\cup\{z\})$ and it follows easily that
$l_u\cap\{u\}=l_z\cap\{z\}$. Thus we may assume that the element
is in $l'_u$, and as $l'_u$ is a clonal line we deduce that
$l'_u\subseteq l_z$. By symmetry we also have $l'_z\subseteq l_u$
and it follows easily that $U\subseteq \cl(Z)$. But then $\cl(Z)$
contains two distinct cocircuits, $U$ and $Z$, so that 
$r(E(M)-\cl(Z))\leq r-2$, and hence $\lambda(Z)=2$, contradicting
the fact that the matroid $M$ is vertically 4-connected.
\end{proof}

Let $A$ denote the common set given by \ref{4clonal13.6},
that is, $A=l_z\cup\{z\}$.

\begin{sublemma}
\label{4clonal13.7}
There is a pair $\{s,t\}\subseteq\{u,v,z\}$ such that
$\sqcap(l'_s\cup l'_t,A)=2$.
\end{sublemma}

\subproof
Assume that the sublemma fails. 
Then $\sqcap(l'_u\cup l'_v\cup l'_z,A)=1$.
Moreover, elementary rank calculations establish that 
$\sqcap(l'_u\cup l'_v,A)=1$ and that 
$\sqcap_{M/l_z}(l'_u\cup l'_v,A)=0$. From this latter fact we deduce
that $\sqcap^*(l'_u\cup l'_v)\geq 3$. Using this and 
Lemma~\ref{lambda-meet}
we obtain
\begin{align*}
\lambda(l'_u\cup l'_v\cup A)&=\lambda(l'_u\cup l'_v)
+\lambda(A)-\sqcap(l'_u\cup l'_v,A)-\sqcap^*(l'_u\cup l'_v,A)\\
&\geq 3+3-1-0=5.
\end{align*}
But by uncrossing the 4-separations $A\cup l'_u$ and $A\cup l'_v$,
we see that $\lambda(l'_u\cup l'_v\cup A)\leq 3$. This contradiction
establishes the sublemma.
\end{proof}

By \ref{4clonal13.7}, we may assume that $\sqcap(l'_u\cup l'_v,A)=2$.
Thus $\sqcap(l'_u\cup l'_v,l_z)=2$. 
But $l'_u\cup l'_v\subseteq H-\{z\}$,
and $\sqcap_{M/z}(H-\{z\},l_z)\leq 1$, so that 
$\sqcap_{M/z}(l'_u\cup l'_v,l_z)<\sqcap_M(l'_u\cup l'_v,l_z)$.
By Lemma~\ref{pi-minor}, $z\in \cl_M(l'_u\cup l'_v)$
contradicting the fact that $z\notin L_M$. This contradiction at last 
completes the proof of the lemma.
\end{proof}

\begin{corollary}
\label{4clonal14}
Let $M$ be an unfettered vertically $4$-connected matroid in 
$\youq$. If $r^*(M)\geq f_{\ref{4clonal12}}(q,t)$, then 
$M$ has a
vertically $4$-connected minor $N$ with 
$|L_N|\geq t$.
\end{corollary}

\begin{proof}
Assume that $r^*(M)\geq f_{\ref{4clonal12}}(q,t)$.
Let $Z=\{z_1,z_2,\ldots,z_l\}$ be a maximal set of elements of
$M$ such that, for all $i\in\{1,2,\ldots,l\}$, the matroid
$M/z_1,z_2,\ldots,z_i$ is unfettered and vertically 4-connected and 
such that
$z_i\notin L_{M/z_1,z_2,\ldots,z_{l-1}}$. By
Lemma~\ref{4clonal4} $M/Z$ is unfettered. Certainly
$Z$ is independent. Hence 
$r^*(M/Z)=r^*(M)\geq f_{\ref{4clonal13}}(q,t)$.
By the definition of $Z$, part (ii) of Lemma~\ref{4clonal12}
does not hold for $M/Z$. Thus part (i) of that lemma holds 
for $M/Z$ and gives the required minor.
\end{proof}

Finally we can achieve the purpose of this section.

\begin{proof}[Proof of Theorem~\ref{4clonal1}]
Let $f_{\ref{4clonal1}}(q,t)=h(q,f_{\ref{4clonal12}}(q,(q+2)t))$. 
Let $M$ be an unfettered 
$4$-connected matroid in $\youq$ with 
at least $f_{\ref{4clonal1}}(q,t)$ elements..
Then $r^*(M)\geq f_{\ref{4clonal12}}(q,(q+2)t)$. 
By Corollary~\ref{4clonal14},
$M$ has a vertically 4-connected minor $N$ with the property that
$|L_N|\geq (q+2)t$. 
Let $N'$ be the matroid obtained by deleting 
all but two elements from each non-trivial clonal line of $N$.
Such lines have at most $q+1$ points.
By Lemma~\ref{4clonal8}, $N'$ is an unfettered vertically 
4-connected matroid. 
But $N'$ has no triangles, so that
$N'$ is 4-connected. Moreover $|L_{N'}|\geq 3t$. Hence
$N'$ has at least $t$ pairwise disjoint clonal pairs.
\end{proof}

\chapter[Unavoidable Minors of Clonal Matroids]{Unavoidable 
Minors of Large $3$-connected Clonal Matroids}
\label{unavoidable-minors-of-large}

\section{Introduction}
Let $M$ be a matroid and $\mathcal C$ be a partition of $E(M)$, then 
$\mathcal C$ is a {\em clonal partition}
\index{clonal partition} 
of $M$ if $C$ is a partition
into clonal pairs.
The matroid $M$ is a 
{\em clonal matroid}
\index{clonal matroid} 
if $E(M)$ has a clonal partition. 

Let $N$ be a minor of $M$. Then a
clonal subset $A$ of $N$ is an $M$-{\em clonal} 
\index{$M$-clonal subset}
subset if $A$ is a 
clonal set in $M$. In particular, the clonal pair $\{a,a'\}$ is 
an {\em $M$-clonal pair}
\index{$M$-clonal pair} 
if $\{a,a'\}$ is a clonal pair in $M$.
We say that $N$ is a {\em clonal minor}
\index{clonal minor}
of $M$ if 
$E(N)$ has a partition
into $M$-clonal pairs. 

Let $M$ be a clonal matroid with associated clonal 
partition $\mathcal P$.
If $\{p,p'\}$ is a member of $\mathcal P$, then we say that $p'$ is
the {\em clonal mate}
\index{clonal mate}
of $p$. We may do this at times without mentioning
the underlying partition, but only when no danger of ambiguity arises.
In any unexplained context,
we indicate the clonal mate of an element by adding a prime symbol.
Thus $a'$ will denote the clonal mate of $a$. If the subset $A$
of elements of $M$ contains no member of $\mathcal P$, then
$A'$ will denoted the set $A'=\{a':a\in A\}$.

The goal of this chapter is to prove.

\begin{theorem}
\label{klonal1}
There exists a function $f_{\ref{klonal1}}(m,q)$
such that, if $M$ is a $4$-connected matroid in
$\eq$ with at least $f_{\ref{klonal1}}(m,q)$ pairwise disjoint
clonal pairs,
then $M$ has a clonal $\Delta_m$-minor.
\end{theorem}

The next lemma gives the first step towards proving 
Theorem~\ref{klonal1}.

\begin{lemma}
\label{3-klonal}
Let $M$ be a $3$-connected matroid with a nonempty set $A$ of elements that
has a partition into pairwise disjoint clonal pairs. Then 
$M$ has a $3$-connected clonal minor on $A$.
\end{lemma}

\begin{proof}
Certainly $M$ is not a wheel or a whirl. Thus, if $M$ has a $4$-element
fan, then this fan has an end. The end of the fan gives an element
$f$ that is not in a clonal pair. The element $f$ is not in $A$ and
can be either deleted or contracted to preserve 3-connectivity. 
Thus we may
assume that $M$ has no $4$-element fans. Say $x\in E(M)-A$. 
Assume that $x$ is
in a triangle $T$. If $T$ contains a member of $A$, then by 
Corollary~\ref{free-triangle},
$M\ba x$ is $3$-connected. Otherwise, by Tutte's Triangle Lemma,
there is an element of $T$ that can be deleted to preserve 
3-connectivity. Assume that $M$ has no triangles or triads, then,
by Bixby's Lemma, either $M\ba x$ or $M/x$ is 3-connected. The lemma
follows by induction on $|E(M)-A|$.
\end{proof}

Given that 4-connected matroids are 3-connected, Theorem~\ref{klonal1}
is an immediate corollary of Lemma~\ref{3-klonal} and the next
theorem.

\begin{theorem}   
\label{unavoidable-swirl}
There is a function $f_{\ref{unavoidable-swirl}}(m,q)$ 
such that if $M$ is a
$3$-connected clonal matroid in $\eq$ whose ground set has at least
$f_{\ref{unavoidable-swirl}}(m,q)$ elements, then $M$ has a clonal
$\Delta_m$-minor.
\end{theorem}

Proving Theorem~\ref{unavoidable-swirl} is the task of this chapter. 
This theorem has a similar flavour to the next important theorem
of Ding, Oporowski, Oxley and Vertigan \cite{doov}. 
A matroid is a {\em whorl}
\index{whorl} 
if it is either a whirl or the cycle matroid of a wheel.

\begin{theorem}
\label{doov}
There is a function $f_{\ref{doov}}(m)$ such that, 
if $M$ is a $3$-connected
matroid with at least $f_{\ref{doov}}(m)$ elements, 
then $M$ has one of the
following as a minor:
$U_{m,m+2}$, $U_{2,m+2}$, $M(K_{3,m})$, $M^*(K_{3,m})$,
a rank-$m$ whorl, or a rank-$m$ spike.
\end{theorem}

We use the following immediate corollary of Theorem~\ref{doov}.

\begin{corollary}
\label{doov1}
There is a function $f_{\ref{doov1}}(m)$ such that, 
if $M$ is a $3$-connected
matroid in $\eq$
with at least $f_{\ref{doov1}}(m)$ elements, 
then $M$ has one of the
following as a minor:
$M(K_{3,m})$, $M^*(K_{3,m})$,
a rank-$m$ whorl, or a rank-$m$ spike.
\end{corollary}

A brief outline of our path to 
Theorem~\ref{unavoidable-swirl} follows.
We begin by using Corollary~\ref{doov1} to find a
$M(K_{3,m})$, $M^*(K_{3,m})$,
rank-$m$ whorl, or rank-$m$ spike as a minor of our
large 3-connected clonal matroid $M$ in $\eq$. 
In doing so we have lost our clones and these need to be
recovered. Up to duality we
obtain a series extension of our minor which has many 
$M$-clonal series pairs. In Section~\ref{bridging-m-clonal} 
we show that we can
find such a series extension where the series pairs are
bridged in a particularly
simple way. After that, it is a matter of inspecting each
type of minor in turn and demonstrating that in each case
the bridging process either produces a violation to the
assumption that we are in $\eq$ or 
produces a large clonal free-swirl minor.

Theorem~\ref{klonal1} serves the needs of this paper,
but it probably not the strongest possible result. Given a large
3-connected clonal matroid, we are happy to obtain a large
line, coline or free-spike minor. The only minor we care about being
clonal is the free swirl. We make the following conjecture.

\begin{conjecture}
\label{clonal-doov}
There is a function $f_{\ref{clonal-doov}}(m)$ such that if a 
$3$-connected clonal matroid has at least $f_{\ref{clonal-doov}}(m)$
elements, then it has one of the following as a clonal 
minor: $U_{2,2m}$, $U_{2m-2,2m}$, $\Lambda_m$, or $\Delta_m$.
\end{conjecture}

\section{Bridging $M$-clonal Series Pairs}
\label{bridging-m-clonal}

As noted above, our first objective is to produce a highly
structured minor with many $M$-clonal series pairs that
are bridged in a particularly simple way. Recall that a 
matroid $M$ is {\em $3$-connected up to series pairs} if,
whenever $(X,Y)$ is a $2$-separation of $M$, then either 
$X$ or $Y$ is a series pair.

\begin{lemma}
\label{clonal-series}
Let $M$ be a $3$-connected clonal matroid in $\eq$. Then
there is a function $f_{\ref{clonal-series}}(m,t,q)$ such that,
if $M$ has at least
$f_{\ref{clonal-series}}(m,t,q)$ elements, then either
$M$ has a clonal $\Delta_m$-minor or, up to duality,
$M$ has a $3$-connected minor $N$ with a coindependent set $J$
such that the following hold.
\begin{itemize}
\item[(i)] $N\ba J$ is $3$-connected up to a set of at least
$t$ series pairs.
\item[(ii)] Each series pair of $N\ba J$ is $M$-clonal.
\item[(iii)] $\co(N\ba J)$ is either a spike, a whorl, or
for some integer $l$, is isomorphic to 
$M(K_{3,l})$ or $M^*(K_{3,l})$.
\end{itemize}
\end{lemma}

The reader should now recall material on bridging sequences from 
Chapter~\ref{background} 
Section~\ref{bridging-sequences}.
It is shown in \cite{gehlwh,leox03} that a bridging sequence
for a 2-separation has at most five elements.
For an $M$-clonal
series pair it is not hard to do better.

\begin{lemma}
\label{bridge-klonal2}
Let $N$ be a connected minor of the matroid $M$, 
let $\{p,p'\}$ be an
$M$-clonal parallel pair in $N$, and let 
$V$ be a  minimal bridging sequence for $\{p,p'\}$.
Then $|V|\leq 2$.
\end{lemma}

\begin{proof}
Say that $V=(v_0,v_1,\ldots,v_t)$. Assume that $|V|>2$,
that is, assume that $\{p,p'\}$ is not bridged
in $N[v_0,v_1]$. Let $Z=E(N)-\{p,q\}$. Consider $N[v_0]$. If $v_0$
is an extension element of $V$, then $v_0\in\cl_{N[v_0]}(\{p,p'\})$,
so that $\{p,p',v_0\}$ is a parallel set in $N[v_0]$.
But then $v_0\in\cl_{N[v_0]}(Z)$, contradicting 
Lemma~\ref{bridge-klonal1}. Thus $v_0$ is a coextension element
of $V$.

\begin{sublemma}
\label{sub-bridge-clonal2}
$\{p,p',v_0\}$ is both a triangle and a triad in $N[(v_0)]$.
\end{sublemma}

\subproof 
This seems clearer in the dual. Here $\{p,p'\}$ is a series pair
of $N^*$. By Lemma~\ref{bridge-basic}, 
$v_0\in\cl_{N^*[v_0]}(\{p,p'\})$.
As $\{p,p'\}$ is a clonal pair, $\{p,p',v_0\}$
is a triangle in $N^*[v_0]$.  If $v_0\in\cl_{N^*[v_0]}(Z)$,
then we violate Lemma~\ref{bridge-klonal1}.
Therefore $\{p,p',v_0\}$ is also a triad.
\end{proof}

Consider $N[v_0,v_1]$. Then $v_1$ is an extension element of 
$V$. By this fact and Lemma~\ref{bridge-basic}, 
in $N[v_0,v_1]$ we have  $v_1\in\cl(\{p,p',v_0\})$, and $v_1$
is not parallel to $p$ or $p'$. Thus $\{p,p',v_1\}$ is a triangle
in $N[v_0,v_1]$ But $v_1\notin\cl_{N[v_0,v_1]}(Z)$, so that
$N[v_0,v_1]\ba v_0\cong N[v_0,v_1]\ba v_1$. Therefore $N[V]/v_0$
has a $N$ as a minor, contradicting Lemma~\ref{bridging-0}. 
\end{proof}

It is perhaps surprising that we have not used the next easy
fact about  spikes and swirls earlier in this paper.

\begin{lemma}
\label{clones-in-spike}
Let $M$ be a spike or a swirl of rank at least four. If $M$ contains
a clonal pair, then $M$ is a free spike or  free swirl.
\end{lemma}

\begin{proof}
Let $\{a,a'\}$ be a clonal pair of $M$. Then $\{a,a'\}$
is certainly a leg of the spike or swirl.
Assume that $M$ is not a free spike or free swirl. Then $M$ contains
a circuit-hyperplane $H$ that is a transversal of the legs.
But then $H$ contains exactly one of $\{a,a'\}$, contradicting the
assumption that $\{a,a'\}$ is a clonal pair.
\end{proof}

\begin{proof}[Proof of Lemma~\ref{clonal-series}]
We lose no generality in assuming that $m>q$ as otherwise,
we can define $f_{\ref{clonal-series}}(m,t,q)$ to be equal to 
$f_{\ref{clonal-series}}(q+1,t,q)$. 
Set $f_{\ref{clonal-series}}(m,t,q)=f_{\ref{doov1}}(m+2t,q)$.
Assume that $|E(M)|\geq f_{\ref{clonal-series}}(m,t,q)$. 
By Corollary~\ref{doov1}, we see that $M$ has minor
$N$ that is either a rank-$(m+2t)$ whorl, a rank-$(m+2t)$ 
spike, or is
isomorphic to $M(K_{3,m+2t})$ or $M^*(K_{3,m+2t})$.

\begin{sublemma}
\label{clonal-series-1}
If $\{a,a'\}$ is a clonal pair of $M$ and $a\in E(N)$, then 
$a'\notin E(N)$. 
\end{sublemma}

\subproof
Elements of whorls,
$M(K_{3,t+2m})$ and $M^*(K_{3,t+2m})$ are either fixed or cofixed
so these matroids contain no clonal pairs. If a spike
contains a clonal pair, then by Lemma~\ref{clones-in-spike},
that spike is a free spike contradicting the assumptions that
$M\in\eq$.
\end{proof}

Choose a partition $\mathcal C$ of $E(M)$ into clonal pairs.
Assume that $N=M/I\ba J$ where $I$ is independent and $J$ is coindependent in $M$. 
Each element $e$ of $N$ has a clonal mate $e'\in E(M)-E(N)$.
As $|E(N)|\geq 2(m+2t)$, we may assume,
up to duality,  that at least $m+2t$ of these clonal mates
are in $J$. Let $K$ be the set of elements of $J$ 
that are not clonal mates of elements of $N$. 
Consider $M/I\ba K$. Say $z\in E(N)$ has a clonal mate
$z'\in J$. Then $\{z,z'\}$ must be  a parallel pair
in $M\ba I/K$ unless $z$ is not fixed in $N$. This can only happen
if either (a) $N$ is a whirl or (b) $N$ is a spike. 
Assume that we are in one of these cases and that 
there are at least $m$ members of $E(N)$ that have clonal mates
in $M\ba I/K$ that are not in parallel pairs.
In case (a) we routinely see that $M\ba I/K$ has a $\Delta_m$-minor,
the legs of which are $M$-clonal, so that
the lemma holds. In case (b) we routinely see that 
$M\ba I/K$ has a $\Lambda_m$-minor and we have contradicted the
assumption that $M\in\eq$.

It follows from the argument of the previous paragraph
that we may assume from now on that
$M/I\ba K$ is $3$-connected up
to a set of at least $2t$ parallel pairs and 
each of these parallel pairs is a  
member of $\mathcal C$. If $\{x,x'\}$ is such a 
parallel pair, then, by Lemma~\ref{bridge-klonal2}, its
corresponding 2-separation has a 
bridging sequence of length at most 2.

Let $X$  be the union of the clonal pairs in 
$M/I\ba K$ that are bridged with a 1-element bridging sequence
and let $Y$ be the union of the clonal pairs that are bridged 
by a 2-element minimal bridging sequence. Let $A=X\cap E(N)$,
$A'=X-A$, $B=Y\cap E(N)$ and $B'=Y-B$.

\begin{sublemma}
\label{clonal-series-2}
The lemma holds if $|A|\geq t$.
\end{sublemma}

\subproof
For each $a\in A$, let $a''$ denote an element of $I$ that bridges
$\{a,a'\}$ and let $A''=\{a'':a\in A\}$. Note that $A''$ could be
a small set, indeed it may only have one element. 
Consider $M/(I-A'')\ba K$. Note that, if
$b\in B$, then $\{b,b'\}$ is a parallel pair in this matroid
as otherwise $\{b,b'\}$ has a 1-element bridging sequence.
Hence each parallel pair of $M/I\ba (J-A')$ is bridged in
$M/(I-A'')\ba (J-A')$ and the claim follows by taking the dual.
\end{proof}

We now consider the case that $|B|\geq t$.
Say $B=\{b_1,b_2,\ldots,b_n\}$. Then 
$B'=\{b_1',b'_2,\ldots,b_n'\}$. 
For $i\in\{1,2,\ldots,n\}$, let $b_i''$ and 
$b_i'''$ be the first and second
elements of a 2-element bridging sequence for $\{b_i,b'_i\}$.
Set $B''=\{b_1'',b''_2,\ldots,b_n''\}$ and 
$B'''=\{b_1''',b'''_2,\ldots,b_n'''\}$.
As $\{b_i,b'_i\}$ is a parallel pair in $M/I\ba K$, we see that
$b''_i$ is a coextension element of the bridging sequence $(b''_i,b'''_i)$
and hence that $b'''_i$ is an extension element of this bridging sequence.
Therefore $B''\subseteq I$ and $B'''\subseteq J$.

\begin{sublemma}
\label{clonal-series-3}
For each $i\in\{1,\ldots,n\}$, the pair $(b''_i,b'''_i)$ is a minimal
bridging sequence for $\{b_i,b'_i\}$ in $M/I\ba (K\cup A)$.
\end{sublemma}

\subproof
Say that $a\in A$. Then $b''_i\notin\cl_{(M/I\ba K)[b''_i]}(\{a,a'\})$,
as otherwise $b''_i$ coblocks $\{b_i,b'_i\}$ contradicting the 
definition of bridging sequences. Thus $\{a,a'\}$ is a parallel
pair in $(M/I\ba K)[b''_i]$ and also in $(M/I\ba K)[b''_i,b'''_i]$.
The claim follows easily from this observation.
\end{proof}

The effect of \ref{clonal-series-3} 
is that we can ignore $A'$ and we have
the following setup.
Let $I''=I-B''$ and $J''=J-(B'\cup B''')$. Then 
$N=\si(M/(I''\cup B'')\ba (J''\cup B'''))$. 
Moreover, for $i\in\{1,\ldots,n\}$, the pair
$\{b_i,b'_i\}$ is parallel in $M/(I''\cup B'')\ba (J''\cup B''')$
and is bridged in $M/I''\ba J''$ by the bridging sequence 
$(b''_i,b'''_i)$.

By the definition of bridging sequence, $\{b_i,b'_i,b''_i\}$ is 
2-separating in $M/(I''\cup (B''-\{b''_i\}))\ba (J''\cup B''')$,
and, in this matroid, $b''_i\in\cl(\{b_i,b'_i\})$. This shows that
$\{b_i,b'_i,b''_i\}$ is a triangle in this matroid and, as
$\{b_i,b'_i\}$ is a clonal pair, $\{b_i,b'_i\}$ is a series pair
in $M/(I''\cup(B''-\{b_i''\}))\ba (J''\cup B'''\cup\{b_i''\})$.

Consider $M/I''\ba (J''\cup B''')$. If the 2-separating set 
$\{b_i,b'_i,b''_i\}$
of $M/(I''\cup (B''-\{b''_i\}))\ba (J''\cup B''')$ 
is bridged in this matroid,
then the parallel pair $\{b_i,b'_i\}$ of 
$M/(I''\cup B'')\ba (J''\cup B''')$
is also bridged. But, as $B''$ is independent in 
$M/I''\ba (J''\cup B''')$,
a minimal bridging sequence must have size one, contradicting the 
assumptions that $\{b_i'',b_i'''\}$ is a minimal bridging sequence
for $\{b_i,b'_i\}$. From this we deduce that
$\{b_i,b'_i,b''_i\}$ is a $2$-separating triangle of 
$M/I''\ba (J''\cup B''')$, and that $\{b_i,b'_i\}$ is a series
pair of $M/I''\ba (J''\cup B'''\cup B'')$.

Moreover, it is easily seen that $\{b_i,b_i',b''_i\}$ 
is a clonal triple
of $M/I''\ba (J''\cup B''')$, so that
$M/(I''\cup B')\ba (J''\cup B'''\cup B'')\cong M/(I''\cup B'')
\ba (J''\cup B'''\cup B')$. But the latter matroid is equal to $N$
and the former is $\co(M/I''\ba (J\cup B'''\cup B''))$. 
Each triangle $\{b_i,b'_i,b''_i\}$ of
$M/I''\ba(J''\cup B''')$ is blocked by $b_i'''$ in $M/I''\ba J''$,
and it follows that each series pair $\{b_i,b'_i\}$ of
$M/I''\ba (J''\cup B'''\cup B'')$ is blocked by $b'''_i$
in $M/I''\ba (J''\cup B'')$. Thus the lemma holds in this case too.
\end{proof}

\section{The Whorl Case.}

We begin by examining the case when we have a large whorl minor.
The next lemma is the goal of this section. Its
proof is surprisingly lengthy.

\begin{lemma}
\label{whorl-case}
Let $M$ be a $3$-connected matroid in $\eq$ with a coindependent
set $J$ such that the following hold.
\begin{itemize}
\item[(i)] $M\ba J$ is $3$-connected up to a set of 
$n$ series pairs that are $M$-clonal.
\item[(ii)] $\co(M\ba J)$ is a whorl.
\end{itemize}
Then there is a function $f_{\ref{whorl-case}}(m,q)$ such
that, if $n\geq f_{\ref{whorl-case}}(m,q)$, then $M$
has a $\Delta_m$-minor, each leg of which is a 
series pair of $M\ba J$.
\end{lemma}

\subsection*{Preliminary Results}

The next three lemmas are easily proved and are certainly well
known. Note that Lemmas~\ref{easy-hyp} and \ref{easy-hyp2}
follow from Ramsey-theoretic results on matrices given in 
\cite{doov1}.
A vertex of a hypergraph is
{\em isolated}
\index{isolated vertex} 
if it is not incident with any edges. Edges of a hypergraph are
{\em parallel} if they are incident with the same set of vertices.

\begin{lemma}
\label{easy-hyp}
Let $l$ and $n$ be integers and 
let $H=(V,E)$ be a hypergraph with no isolated vertices
and no parallel pairs of edges.
Assume that no edge of
$H$ is incident with more than $l$ vertices and let 
$U$ be an $n$-element set of vertices of $H$. Then there is a function
$f_{\ref{easy-hyp}}(m,l)$ such that, if 
$n\geq f_{\ref{easy-hyp}}(m,l)$,
then there is a subset $\{u_1,u_2,\ldots,u_m\}$ of $V$ and a set
$\{e_1,e_2,\ldots,e_m\}$ of edges of $H$ such that, 
for $i\in\{1,2,\ldots,m\}$,
$u_i$ is incident with $e_j$ if and only if $i=j$.
\end{lemma}

The next lemma is a strengthening of Lemma~\ref{easy-hyp} that
gives a somewhat more specific outcome.

\begin{lemma}
\label{easy-hyp2}
Let $l$ be an integer and let $H=(V,E)$ be a hypergraph with no
isolated vertices and no parallel pairs of edges. 
Assume that no edge of $H$ is incident with more
than $l$ vertices. Let $\phi:V\rightarrow E$ be a function
such that, for all $v\in V$, the vertex $v$ is incident with
$\phi(v)$. Then there is a function
$f_{\ref{easy-hyp2}}(m,l)$ such that, if 
$n\geq f_{\ref{easy-hyp2}}(m,l)$,
then there is a subset $\{u_1,u_2,\ldots,u_m\}$ of $V$ 
such that, 
for $i\in\{1,2,\ldots,m\}$, the vertex
$u_i$ is incident with $\phi(u_j)$ if and only if $i=j$.
\end{lemma}

\begin{lemma}
\label{ham}
Let $G=(V,E)$ be a graph where $V$  is cyclically ordered
and $E$ is a matching.
Then there is a function $f_{\ref{ham}}(m)$ such that, if $|E|=n$ and 
$n\geq f_{\ref{ham}}(m)$, then, 
for some labelling $(v_1,v_2,\ldots,v_l)$
of $V$ that respects the cyclic order, there is a set
of $m$ edges in $E$ that can be directed and ordered
$((v_{i_1},v_{j_1}),(v_{i_2},v_{j_2}),\ldots,(v_{i_m},v_{j_m}))$ 
such that one of the following holds.
\begin{itemize}
\item[(i)] $i_1<j_1<i_2<j_2<\cdots<i_m<j_m$.
\item[(ii)] $i_1<i_2<\cdots <i_m<j_m<j_{m-1}<\cdots<j_1$.
\item[(iii)] $i_1<i_2<\cdots<i_m<j_1<j_2<\cdots<j_m$.
\end{itemize}
\end{lemma} 

Recall that a flower $\PP$ in a connected matroid $M$ has the
property that if $(X,Y)$ is a 2-separation, then either 
$X$ or $Y$ is contained in a petal of $\PP$. Let 
$\PP=(P_1,P_2,\ldots,P_m)$
be a flower in the connected matroid $M$. 
Recall that a clonal pair
$\{p_i,q_i\}$ contained in the petal $P_i$ is {\em $\PP$-strong} if 
$\kappa(\{p_i,q_i\},P_1\cup P_2\udots P_{i-1}\cup P_{i+1}\udots P_m)=2$.
Equivalently $\{p_i,q_i\}$ is $\PP$-strong if 
there is no $2$-separating
set $X$ of $M$ with $\{p_i,q_i\}\subseteq X\subseteq P_i$. 
For convenience we restate here a special case of 
Lemma~\ref{free-link}.

\begin{lemma}
\label{get-free-swirl}
Let $M$ be a connected matroid with a swirl-like flower
$\PP=(P_1,P_2,\ldots,P_n)$ such that for all $i\in\{1,2,\ldots,n\}$
the petal $P_i$ contains a $\PP$-strong clonal pair $\{p_i,q_i\}$ 
Then $M$ contains a $\Delta_n$-minor with associated flower
$(\{p_1,q_1\},\{p_2,q_2\},\ldots,\{p_n,q_n\})$.
\end{lemma}

We will also use the following technical but elementary lemma.

\begin{lemma}
\label{blocking-lemma}
Let $M$ be a matroid with an element $z''$ such that $M\ba z''$
is connected with
an $M$-clonal series pair $\{z,z'\}$. Assume that $M\ba z''$ has 
an exact $3$-separation $(X,Y)$ where $\{z,z'\} \subseteq Y$.
Assume further that $(X,Y)$ is blocked in $M$ and that
$\lambda_{M\ba z''}(X\cup\{z,z'\})>2$. Then the $3$-separation 
$(X,Y-\{z'\})$ of $N\ba z''/z'$ is blocked by 
$z''$ in $M/z'$.
\end{lemma}

\begin{proof}
Evidently $(X,Y-\{z'\})$ is an exact 3-separation of $M\ba z''/ z'$.
If this is not blocked by $z''$
in $M/z'$, then either $z''\in\cl_{M/z'}(Y-\{z'\})$
or $z''\in\cl_{M/z'}(X)$. In other words,
$z''\in\cl_M(Y)$ or $z''\in\cl_M(X\cup\{z'\})$.
The former case does not occur.
Consider the latter. In this case $z'\in\cl_M(X\cup\{z''\})$ and,
as $\{z,z'\}$ are clones, $z\in\cl_M(X\cup\{z''\})$. Thus 
$r_M(X\cup\{z,z'\})\leq r_M(X)+1$. But 
$r_M(Y-\{z,z'\})\leq r_M(Y)-1$. Therefore 
$\lambda_{M\ba z''}(X\cup\{z,z'\})\leq 2$, 
contradicting a hypothesis of the lemma.
\end{proof}

It is perhaps surprising the we have not needed the following
lemmas on freedom until now.

\begin{lemma}
\label{more-free}
Let $\{a,b,c\}$ be a triad of the matroid $M$, 
where $b$ and $c$ are clones.
Then $a$ is freer than $b$ in $M/c$.
\end{lemma}

\begin{proof}
Let $F$ be a cyclic flat of $M/c$ that contains $a$. Then either $F$ or
$F\cup\{c\}$ is a cyclic flat of $M$. In the latter case $F$ contains
$a$ as $\{a,c\}$ is a clonal pair. In the former case $F$ contains
either $a$ or $c$ as $\{a,b,c\}$ is a triad of $M$. But again, as
$a$ and $c$ are clones we deduce that $F$ contains $a$.
\end{proof}

\begin{lemma}
\label{keep-free}
Let $a$ and $b$ be elements of a matroid $M$, 
where $a$ is freer than $b$.
If $N$ is a minor of $M$ whose ground set contains 
both $a$ and $b$, then 
$a$ is freer than $b$ in $N$.
\end{lemma}

\begin{proof}
Say $c\in E(M)-\{a,b\}$. Consider $M\ba c$.
If $F$ is a cyclic flat of $M\ba c$, then either 
$F$ or $F\cup\{c\}$ is a cyclic flat of $M$. In either case
we deduce that $b\in F$. Thus $a$ is freer than $b$ in $M\ba c$.

As $a$ is freer than $b$ we have $b$ is freer than $a$ in
$M^*$, so that $b$ is freer than $a$ in $M^*\ba c$.
Thus $a$ is freer than $b$ in $M/c$. The lemma follows from
these observations.
\end{proof}

\subsection*{Cleanly-blocked Coextended Whorls.} 
The matroid $M$ is a
{\em coextended whorl}
\index{coextended whorl}
if it is $3$-connected up to series pairs
and $\co(M)$ is a whorl.
We now consider the case where the series 
pairs of a coextended whorl are blocked 
one at a time by the blocking elements. 
We begin by developing some terminology for this case.

Let $N$ be a rank-$n$ whorl
with rim elements labelled $R=\{r_1,r_2,\ldots,r_n\}$ 
and spoke elements
labelled $S=\{s_1,s_2,\ldots,s_n\}$. 
We say this labelling is {\em standard}
\index{standard labelling of a whorl}
if, for all $i\in\{1,2,\ldots,n\}$,
the sets $\{s_i,r_i,s_{i+1}\}$ and $\{r_i,s_{i+1},r_{i+1}\}$
are respectively triangles and triads of $N$, where indices
are taken modulo $n$. 
Let $M$ be a matroid. Then $M$ is a {\em cleanly blocked}
\index{cleanly blocked coextended whorl} 
coextended whorl of {\em order}
$n$ with {\em distinguished $5$-tuple}
\index{distinguished $5$-tuple} 
$(R,S,T,T',T'')$
if we have labellings $R=\{r_1,r_2,\ldots,r_n\}$, 
$S=\{s_1,s_2,\ldots,s_n\}$,
$T=\{t_1,t_2,\ldots,t_n\}$, $T'=\{t'_1,t'_2,\ldots,t'_n\}$
and $T''=\{t''_1,t''_2,\ldots,t''_n\}$ such that the following hold.
\begin{itemize}
\item[(i)] $E(M)$ consists of the union of the 
disjoint sets $R$, $S$, $T'$ and 
$T''$. The set $T$ is contained in $R\cup S$.
\item[(ii)] $M\ba T''/T'$ is a whorl for 
which $R$ and $S$ give a standard labelling.
\item[(iii)] $M\ba T''$ is a coextended whorl with series pairs 
$\{\{t_i,t'_i\}:i\in\{1,2,\ldots,n\}\}$.
\item[(iv)] The series pairs of $M\ba T''$ are $M$-clonal.
\item[(v)] For all $i\in\{1,2,\ldots,n\}$, 
the element $t''_i$ blocks the
series pair $\{t_i,t'_i\}$ of $M\ba T''$, but 
blocks no other series pair.
\end{itemize}

Our goal for this case is to prove.

\begin{lemma}
\label{clean-block}
Let $M$ be a cleanly-blocked coextended whorl in $\eq$ with 
distinguished $5$-tuple $(R,S,T,T',T'')$.
Then there is a function $f_{\ref{clean-block}}(m,q)$ such that,
if $|T|\geq f_{\ref{clean-block}}(m,q)$, then $M$ has
a $\Delta_m$-minor, each leg of which is a series pair of $M\ba T''$.
\end{lemma}

Note that if $M$ is a cleanly blocked coextended whorl, then $M$ is 
$3$-connected, the set $\{t_i,t'_i,t''_i\}$ is a triad
for all $i\in\{1,2,\ldots,n\}$ and $\co(M\ba T'')=M\ba T''/T'$.
We call a triad of the form $\{t_i,t'_i,t''_i\}$ a {\em flap}
\index{flap} 
of $M$. The next lemma is immediate.

\begin{lemma}
\label{flap-gone}
If $\{t_i,t'_i,t''_i\}$ is a flap of $M$, then $M\ba t''_i/t'_i$
is a cleanly blocked coextended whorl of order $n-1$ with distinguished
$5$-tuple $(R,S,T-\{t_i\},T'-\{t'_i\},T''-\{t''_i\})$.
\end{lemma}

The minor $N$ of $M$ is obtained by {\em removing flaps}
\index{flap removal} 
if it is obtained
by a sequence of operations of the form 
described in Lemma~\ref{flap-gone}.
The matroid $M$ is {\em rim based}
\index{rim-based coextended whorl} if all the series pairs of
$M\ba T''$ are based at rim elements, that is, if $T\subseteq R$. It
is {\em spoke based}
\index{spoke-based coextended whorl} 
if all the series pairs 
are based at spoke elements,
that is, if $T\subseteq S$. We consider the two cases in turn.

\subsection*{The Rim Case.}

Our goal in the rim case is to prove

\begin{lemma}
\label{rim-case}
Let $M$ be a rim-based cleanly-blocked coextended whorl in $\eq$ with 
distinguished $5$-tuple $(R,S,T,T',T'')$.
Then there is a function $f_{\ref{rim-case}}(m,q)$ such that,
if $|T|\geq f_{\ref{rim-case}}(m,q)$, then $M$ has
a $\Delta_m$-minor, each leg of which is a series pair of $M\ba T''$.
\end{lemma}

Throughout this subsection we assume that $M$ is a rim-based cleanly blocked coextended whorl with distinguished 5-tuple $(R,S,T,T',T'')$. 
Note that a flap of $M$ will have the form $\{r_i,r'_i,r''_i\}$
for some $r_i\in R$.

Say $|R|=n$. 
If $1\leq k<n$, then the set $\{s_n,r_n,s_1,r_1,\ldots,r_{k-1},s_k\}$
is 3-separating in $M\ba T''/T'$ and the corresponding $3$-separation is 
clearly induced in $M\ba T''$. We denote this induced
3-separation by $(L_k,K_k)$. We say that the element $r''_h$
of $T''$ {\em blocks} $(L_k,K_k)$ {\em from the right}
\index{blocks from the right}
if $h\in\{k,k+1,\ldots,n-1\}$ and $r''_h$ blocks 
$(L_k,K_k)$. On the other hand, $r''_h$ {\em blocks} $(L_k,K_k)$ 
{\em from the left}
\index{blocks from the left} 
if $h\in\{1,2,\ldots,k-1\}$ and $r''_h$ blocks $(L_k,K_k)$.

\begin{lemma}
\label{bound-block}
If $M\in \eq$ and $2\leq k< n$, then at most 
$\frac{q^k-1}{q-1}$ elements of $T''$ block $(L_k,K_k)$
from the right.
\end{lemma}

\begin{proof}
Let $U''$ be the set of elements of $T''$ that block $(L_k,K_k)$
from the right. Let $U=\{r_i\in T:r''_i\in U''\}$.
Let $(L,K)=(\{s_n,r_n,s_1,r_1,\ldots,r_{k-1},s_k\},\{r_k,s_{k+1},\ldots,
s_{n-1},r_{n-1}\})$. By Lemma~\ref{blocking-lemma},
we have

\begin{sublemma}
\label{bound-block-1}
If $u''\in U''$, then $u''$ blocks the $3$-separation $(L,K)$
of $M\ba T''/T'$.
\end{sublemma}

Define the equivalent relation $\sim$ on $\{s_k,s_{k+1},\ldots,s_n\}$
as follows. If $k\leq i\leq j\leq n$, then $s_i\sim s_j$ if the set 
$\{r_i,r_{i+1},\ldots,r_{j-1}\}$ contains no member of $U$. Otherwise
$s_i\sim s_j$ if $s_j\sim s_i$. Evidently $\sim$ has $|U|+1$
equivalence classes. Let
$$V=U''\cup\{r_i:i\in\{r_k,r_{k+1}\ldots,r_{n-1}\};r_i\notin U\}.$$

\begin{sublemma}
\label{bound-block-2}
$V$ is independent in $M/T'$.
\end{sublemma}

\subproof
As $M\ba T''/T'$ is a whorl, 
the set $\{r_k,r_{k+1},\ldots,r_{n-1}\}$ of 
rim elements is independent in $M/T'$. It is elementary that if $x$ is 
an element of an independent set $X$ of a 
matroid and $x''$ is freer than
$x$, then $(X-\{x\})\cup\{x''\}$ is also independent. The claim now
follows from the above facts, Lemma~\ref{more-free} and
Lemma~\ref{keep-free}.
\end{proof}

We now focus on the rank-$k$ matroid $M/T'/V$, 
the goal being to show that it has
at least $|U''|$ parallel classes.

\begin{sublemma}
\label{bound-block-3}
If $i\in\{1,2,\ldots,k\}$, then $\{s_i\}$ is independent in $M/T'/V$.
\end{sublemma}

\subproof
As $M\ba T''/T'$ is a whorl, $\{s_i,r_k,r_{k+1},\ldots,r_{n-1}\}$
is independent in $M\ba T''/T'$. By Lemma~\ref{keep-free} we see that
$\{s_i\}\cup V$ is independent in $M/T'$, proving the claim.
\end{proof}

\begin{sublemma}
\label{bound-block-4}
If $i,j\in\{k,k+1,\ldots,n\}$, and $s_i\not\sim s_j$, 
then $\{s_i,s_j\}$ is 
independent in $M/T'/V$.
\end{sublemma}

\subproof
Assume otherwise. Then, by \ref{bound-block-3}, 
$\{s_i,s_j\}$ is a circuit
of $M/T'/V$ and there is a subset $W$ of $V$ 
such that $\{s_i,s_j\}\cup W$
is a circuit of $M/T'$. Let $X$ be the set consisting of those elements
$r_t$ of $\{r_k,r_{k+1},\ldots,r_{n-1}\}$ 
such that either $r_t\in W$ or 
$r''_t\in W$. Note that, as $s_i\not\sim s_j$, there is at least one
element $r''_\alpha$ in $W\cap U''$.

Say $r''_t\in W$. Then 
$r''_t\in\cl_{M/T'}((W-\{r''_t\})\cup \{s_i,s_j\})$
and, as $r''_t$ is freer than $r_t$ in this matroid, we see that 
$r_t\in\cl_{M/T}(W\cup\{s_i,s_j\})$. 
We conclude that $\cl_{M/T}(W\cup\{s_i,s_j\})$ contains
$X\cup\{s_i,s_j\}$. But, in $M\ba T''/T'$, 
this is a set of at most $r-2$ rim elements
together with two spoke elements and contains at most one circuit. Thus,
the rank of $X\cup\{s_i,s_j\}$ in $M/T'$ is equal to the rank of 
$W\cup\{s_i,s_j\}$ in this matroid. Therefore $X\cup\{s_i,s_j\}$
spans $\cl_{M/T'}(W\cup\{s_i,s_j\})$, 
so that $r''_\alpha\in\cl(X\cup\{s_i,s_j\})$
contradicting the fact that $r''_\alpha$ blocks this set.
\end{proof}

As $r(M/T'/V)=k$, and the relation
$\sim$ has $|U|+1$ equivalence classes,  it follows from
\ref{bound-block-4} that this matroid has at least $|U|=|U''|$ parallel
classes. The elements of $U''$ block $(L_u,K_u)$ from the right.
As $M\in \eq$, by Lemma~\ref{kung}, there are at most 
$\frac{q^k-1}{q-1}$ of them.
\end{proof}

It may be that an element $r_i\in R-T$ has the property that
$M/r_i\ba s_i$ is a cleanly-blocked coextended whorl with
distinguished sets $(R-\{r_i\},S-\{s_i\},T,T',T'')$. If this is
the case we say that $M/r_i\ba s_i$ is a {\em reduction} of 
$M$, and that $M$ is {\em reduced} if it has no reductions.
The next lemma is clear.

\begin{lemma}
\label{easy-reduce}
Assume that $M$ is reduced. If $r_i\in R-T$, then there is a 
flap $\{r_j,r'_j,r''_j\}$ such that $\{r_i,r_j,r'_j,r''_j\}$
is a circuit.
\end{lemma}

In the reduced case we can bound the number of elements that are 
in $R-T$.

\begin{lemma}
\label{reduce-r-t}
Assume that $M\in \eq$ is reduced. Then there is a function
$f_{\ref{reduce-r-t}}(m,q)$ such that, if
$|R-T|\geq f_{\ref{reduce-r-t}}(m,q)$, then $M$ has a 
$\Delta_m$-minor all of whose legs are series pairs of $M\ba T''$.
\end{lemma}

\begin{proof}
Let $f_{\ref{reduce-r-t}}(m,q)=
(q+2)f_{\ref{ham}}(2f_{\ref{get-free-clonal}}(m,q)+2).$

Define a graph $G$ with vertex set $R$ such that, for 
$r_i\in R-T$ and $r_j\in T$, the pair  $\{r_i,r_j\}$ is an edge
if $\{r_i,r_j,r_j',r''_j\}$ is a circuit. Note that,
if $r_j\in T$, then the degree of $r_j$ is at most one in
this graph so that it decomposes into stars.

\begin{sublemma}
\label{reduce-r-t-1}
If $r_i\in R$, then $d(r_i)\leq q+2$.
\end{sublemma}

\subproof
Up to labels we may assume that $i=1$, that is $r_i=r_1$. 
Assume that $d(r_1)>1$. Then $r_i\in R-T$. By 
Lemma~\ref{bound-block} at most $q+1$ elements block 
$(L_2,K_2)$ from the right. If $\alpha\neq s$, and 
$\{r_1,r_\alpha\}$ is an edge of $G$, then $r_\alpha\in T$,
and $r''_\alpha$ blocks $(L_2,K_2)$ from the right. 
There are at most 
$q+1$ such elements. It is possible that $\{r_1,r_n\}$
is an edge. Altogether we have $d(r_i)\leq q+2$.
\end{proof}

Thus $G$ has a matching of size 
$f_{\ref{ham}}(2f_{\ref{get-free-clonal}}(m,q)+2).$
By Lemma~\ref{ham}, $G$ has a collection $A$ of 
$2f_{\ref{get-free-clonal}}(m,q)+2$ edges whose indices
can be labelled to satisfy one of the cases of Lemma~\ref{ham}.
Let $N$ be the matroid obtained from $M$ by removing all
flaps $\{r_i,r'_i,r''_i\}$ for which $r_i\notin A$. If either Case~(i)
or (ii) of Lemma~\ref{ham} holds, then it is easily seen
that $N$ has a path of $3$-separations of length at least
$f_{\ref{get-free-clonal}}(m,q)$ each step of which contains a 
series pair of $M\ba T''$. Thus each step of
$\PP$ contains a $\PP$-strong clonal pair, 
so that in these cases the lemma
follows from Corollary~\ref{get-free-clonal}. 
Note that, if Case~(i) holds
this minor is found more directly via Lemma~\ref{get-free-swirl}. 

Consider Case~(iii). In this case we have a sequence of indices
$i_1<i_2<\cdots<j_1<j_2<\cdots$ such that, 
for $1\leq l\leq 2f_{\ref{get-free-clonal}}(m,q)+2$, 
either $r_{i_l}\in T$ and 
$\{r_{i_l},r'_{i_l},r''_{i_l},r_{j_l}\}$ is a circuit,
or $r_{j_l}\in T$ and 
$\{r_{j_l},r'_{j_l},r''_{j_l},r_{i_l}\}$ is a circuit. 
Up to labels we may assume that the former case occurs at
least $t=f_{\ref{get-free-clonal}}(m,q)+1$ times. (We do this simply 
for notational convenience; not for any structural reason.)
After another round of flap removals, reductions and label resetting,
we obtain a cleanly bridged coextended whorl $N'$ with 
distinguished 5-tuple
$(R',S',U,U',U'')$ such that $|R'|=2t$,
$U=\{r_1,r_2,\ldots,r_t\}$ and, for all $i\in\{1,2,\ldots,t\}$,
the set $\{r_i,r'_i,r''_i,r_{t+i}\}$ is a circuit.
Let $N''=N'\ba \{r_1,r'_1,r''_1,r_{t+1}\}$.
For $k\in\{2,3,\ldots,t\}$, let 
$P_k=\{s_k,s_{t+k},r_k,r'_k,r'_{t+k}\}$ and set 
$\PP=(P_2,P_3,\ldots,P_{t-1},P_t\cup\{s_1,s_{2t}\})$.

\begin{sublemma}
\label{reduce-r-t-2}
$N''$ is $3$-connected,
and $\PP$ is a path of $3$-separations in this matroid.
\end{sublemma}

\subproof
Consider the matroid $N''\ba U''/U'$. 
As this matroid is obtained by deleting
two rim elements of a whorl, it is a matter of 
elementary graph theory to 
verify the following facts. Let
$S_1=(s_2,r_2,s_3,\ldots,r_t,s_{t+1})$ and 
$S_2=(s_{t+2},r_{t+2},s_{t+3},\ldots,s_{2t},r_{2t},s_1)$.
\begin{itemize}
\item[(i)] The unique separation of $N''\ba U''/U'$
is $(S_1,S_2)$.
\item[(ii)] If $(X,Y)$ is a 2-separation of $N''\ba U''/U'$
then either $X$ or $Y$ is an initial or terminal 
segment of $S_1$ or $S_2$.
\item[(iii)] For $i\in\{2,3,\ldots,t-1\}$, the set 
$Z_i=\{s_2,r_2,\ldots,s_i,r_i,s_{t+2},r_{t+2},\ldots,s_{t+i},r_{t+i}\}$
is $3$-separating in $N''\ba U''/U'$.
\end{itemize}

We omit the routine verification of the following claims. 
The separation $(S_1,S_2)$ of $N''\ba U''/U'$ is not induced in 
$N''/U'$. Indeed $\kappa_{N''}(S_1,S_2)\geq 3$. Also all 2-separations
of $N''\ba U''/U'$ are bridged in $N''$. 
From these claims we deduce that
$N''$ is $3$-connected.

Consider a $3$-separating set $Z_i$ as described in (iii). Say
$j\in\{1,2,\ldots,i\}$. Then $r_j\in\cl^*_{N''\ba U''}(\{r_j\})$
as $\{r_j,r'_j\}$ is a series pair of this matroid. 
Thus $\lambda_{N''\ba U''}(Z_i\cup\{r'_2,r'_3,\ldots,r'_i\})=2$.
Also $\{r_j,r'_j,r''_j,r_{t+j}\}$ is a circuit of $N''$.
Hence 
$Z_i\cup\{r_2',r'_3,\ldots,r_i'\}\cup\{r''_2,r''+3\ldots,r''_i\}$
is $3$-separating in $N''$. But this set is equal to 
$P_2\cup P_3\udots P_i$.
\end{proof}

The lemma now follows from \ref{reduce-r-t-2},
the fact that $\PP$ has length 
$t-1= f_{\ref{get-free-clonal}}(m,q)$,
and Corollary~\ref{get-free-clonal}.
\end{proof}

We develop some notation for the next lemma. If $M_t$
is a minor of $M$ obtained by flap removal and reductions, we denote
its distinguished 5-tuple by $(R_t,S_t,T_t,T'_t,T''_t)$,
where for some integer $l$, 
we have $R_t=\{r^t_1,r^t_2,\ldots,r^t_l\}$
and $S_t=\{s^t_1,s^t_2,\ldots,s^t_l\}$. 
For an integer $i$, we define the
$3$-separation $(L^t_i,K^t_i)$ in a way that is precisely analogous 
to the way that we defined $(L_i,K_i)$ in $M$.

\begin{lemma}
\label{right-clean}
Assume that $M\in \eq$ and $M$ has no $\Delta_m$-minor
each leg of which is a clonal pair of $M\ba T''$. Then there is a
function $f_{\ref{right-clean}}(h,m,q)$ with the property that,
if $|T|\geq f_{\ref{right-clean}}(h,m,q)$, then 
$M$ has a minor $M_h$ obtained by flap removals 
and reductions such that the following hold.
\begin{itemize}
\item[(i)] There is a sequence of indices $(h_1,h_2,\ldots,h_h)$
such that, for $i\in\{1,2,\ldots,h\}$, the $3$-separation
$(L^h_{h_i},K^h_{h_i})$ is not blocked from the right.
\item[(ii)] $L^h_{h_1}$ contains exactly one member of $T_h$ and,
if $1<i\leq h$, then $L^h_{h_i}-L^h_{h_{i-1}}$ contains
exactly one member of $T_h$.
\end{itemize}
\end{lemma}

\begin{proof}
For $0\leq u\leq h$, 
let $\mu(u,m,q)=(q^{(u+1+f_{\ref{reduce-r-t}}(m,q))}-1)/(q-1)$,
and let $f_{\ref{right-clean}}(h,m,q)=\sum_{u=0}^h\mu(u,m,q)$.

Let $M_0=M$ and assume that, for some $u\in\{0,1,\ldots,h-1\}$,
the matroid $M_u$, obtained from $M$ by flap removals and reductions,
satisfies the conclusion of the lemma with $h$ replaced by $u$.

\begin{sublemma}
\label{right-clean-1}
If $|K^u_{u_u}\cap T_u|\geq 1$, then $M$ has a minor $M_{u+1}$,
obtained from $M$ by flap removals and reductions that 
satisfies the conclusions for the lemma with $m$ replaced by 
$u+1$, and such that 
$$|K^{u+1}_{(u+1)_{u+1}}\cap T_{u+1}|\geq |K^u_{u_u}|-\mu(u,m,q).$$
\end{sublemma}

\subproof
It is clear that reductions preserve the desired properties of
$M_u$, so that we may assume that $M_u$ is reduced. Let 
$i$ be the least integer such that $r^u_i\notin L^u_i$, 
and $r^u_i\in T_u$.
Consider $L^u_i$. This set contains $u+1$ elements of 
$R_u\cap T_u$ and, by Lemma~\ref{reduce-r-t},
at most $f_{\ref{reduce-r-t}}(m,q)$ elements of $R_u-T_u$.
Hence $i\leq u+1+f_{\ref{reduce-r-t}}(m,q)$.
It now follows from Lemma~\ref{bound-block} that at most
$\mu(u,m,q)$ elements of $T''_u$ block $(L^u_i,K^u_i)$ from the 
right. If we remove the flaps associated with the blocking 
elements, then we obtain the desired minor.
\end{proof}

The lemma follows from \ref{right-clean-1} and induction.
\end{proof}

\begin{proof}[Proof of Lemma~\ref{rim-case}]
Let
$$f_{\ref{rim-case}}(m,q)=
f_{\ref{right-clean}}
(f_{\ref{right-clean}}((f_{\ref{get-free-clonal}}(m,q),m,q),m,q)).$$
Assume that $|T|\geq f_{\ref{right-clean}}(m,q)$. 
Let $h=f_{\ref{right-clean}}(f_{\ref{get-free-clonal}}(m,q),m,q)$.
Assume that the lemma fails. Then
by Lemma~\ref{right-clean} and a reversal of the indices, we obtain
a minor $M_h$ of $M$ with a sequence $\{h_1,h_2,\ldots,h_h\}$ of 
indices such that, for $i\in\{1,2,\ldots,h\}$, the 3-separation
$(L^h_{h_i},R^h_{h_i})$ is not blocked from the 
left and such that, after
possibly removing some extra flaps, 
has the property that $|L_1^h\cap T_h|=1$,
and, for $i\in \{2,3,\ldots,h\}$, the set $L^h_i-L^h_{i-1}$ 
contains one element
of $T_h$. We would like to apply Lemma~\ref{right-clean} again, but
we are not quite in a position to do this 
as we have distinguished indices to
worry about. We omit the details  of the obvious upgrade
of Lemma~\ref{right-clean} that covers this
and conclude that, for some
$l\geq f_{\ref{get-free-clonal}}(m,q)$, we have a minor $M_l$
of $M$, obtained by flap removal 
and reductions that has the property that
there is a sequence of indices, $\{l_1,\ldots,l_l\}$ such that,
for $i\in\{1,2,\ldots,l\}$, the separation $(L^l_i,K^l_i)$ is neither
blocked from the left nor the right, that is, is 3-separating
in $M_l$. 

Let $P_1=L^l_1$,  let $P_i=L^l_i-L^l_{i-1}$ 
for $i\in\{2,3,\ldots,h-1\}$, and
let $P_h=K^l_{h-1}$. Then the path $(P_1,P_2,\ldots,P_h)$ of 
$3$-separations of $M_l$ has the property 
that each step contains a clonal
pair that is a series pair of $M\ba T''$. 
By Corollary~\ref{get-free-clonal}
$M_l$ has a $\Delta_m$-minor each leg of which is a clonal series pair
of $M\ba T''$ and the lemma follows.
\end{proof}

\subsection*{The Spoke Case} The goal here is to prove

\begin{lemma}
\label{big-spoke}
Let $M$ be a spoke-based cleanly-blocked coextended whorl  in 
$\eq$ with 
distinguished $5$-tuple $(R,S,T,T',T'')$. Then there is a function
$f_{\ref{big-spoke}}(m,q)$ such that, if 
$|T|\geq f_{\ref{big-spoke}}(m,q)$, then $M$ has a $\Delta_m$-minor,
each leg of which is a series pair of $M\ba T''$.
\end{lemma}

Now for the series of lemmas that lead to Lemma~\ref{big-spoke}.

\begin{lemma}
\label{spoke0}
Let $M$ be a matroid whose ground set contains
disjoint sets $T=\{t_1,t_2,\ldots,t_n\}$, 
$T'=\{t'_1,t'_2,\ldots,t'_n\}$,
and $T''=\{t''_1,t''_2,\ldots,t''_n\}$ such that the following hold.
\begin{itemize}
\item[(i)] For all $i\in\{1,2,\ldots,n\}$, the
pair $\{t_i,t'_i\}$ is an $M$-clonal series 
class of $M\ba T''$, and these are the only non trivial series
classes of $M\ba T''$. 
\item[(ii)] $T$ is a parallel class of $\co(M\ba T'')=M\ba T''/T'$.
\item[(iii)] Up to the parallel class $T$, the matroid 
$M\ba T''/T'$ is a whorl with rim elements $R$, and $T$ is at 
a spoke of this whorl.
\item[(iv)] There is an injective function 
$\phi:\{1,2,\ldots,n\}\rightarrow R$ such that, for all
$i\in\{1,2,\ldots,n\}$, the set $\{t_i,t'_i,t''_i,\phi(i)\}$
is a circuit of $M$.
\end{itemize}
Then there is a function $f_{\ref{spoke0}}(q)$ such that, if 
$n\geq f_{\ref{spoke0}}(q)$, then $M\notin \eq$. 
\end{lemma}

\begin{proof}
Let $f_{\ref{spoke0}}(q)=f_{\ref{main1}}(5,q)+3$.

Let $(s_1,r_1,s_2,r_2,\ldots,s_t,r_t)$ be a labelling of 
$\si(M\ba T''/T')=M\ba T''/T'\ba (T-\{t_1\})$, 
where $R=\{r_1,r_2,\ldots,r_t\}$, $S=\{s_1,s_2,\ldots,s_t\}$,
$t_1=s_1$ and triples of the
form $(s_i,r_i,s_{i+1})$ and $(r_i,s_{i+1},r_{i+1})$ are 
triangles and triads respectively. We may assume that, if 
$i>j$, then $\phi(i)>\phi(j)$. Define a partition
$(P_1,P_2,\ldots,P_n)$ of $E(M\ba T''/T')$ as follows: let
$P_1=\{t_1,r_1,s_2,r_2,\ldots,r_{\phi(1)}\}$; let 
$P_k=\{t_k,s_{\phi(k-1)+1},r_{\phi(k-1)+1},\ldots,s_{\phi(k)},
r_{\phi(k)}\}$ for $k\in\{2,\ldots,n-1\}$; and let
$P_n=\{t_n,s_{\phi(n-1)+1},r_{\phi(n-1)+1},\ldots,s_t,r_t\}$.
This is a path of 3-separations in $M\ba T''/T$. Moreover,
$(P_1\cup\{t_1',t''_1\},P_2\cup\{t_2',t_2''\},\ldots,
P_n\cup\{t_n',t_n''\})$ is a path of $3$-separations in $M$ of length
$n-1$. 

Furthermore, it is readily checked that if $(A,B)$ is a 
$3$-separation of $M$, where $|A|,|B|>4$, then either 
$A$ or $B$ is a fan of $M$. Thus $M$ is 5-coherent. 
It now follows from Corollary~\ref{main2} that, if 
$n \geq f_{\ref{spoke0}}(q)$, then $M$ is not in $\eq$.
\end{proof}

\begin{lemma}
\label{spoke1}
Let $M$ be a spoke-based cleanly-blocked coextended whorl in 
$\eq$ with distinguished $5$-tuple $(R,S,T,T',T'')$.
Assume that there is an injection $\phi$ from the indices
of members of $T$ to the indices of members of $R$ such that,
for all $s_i\in T$, the set
$\{s_i,s'_i,s''_i,r_{\phi(i)}\}$ is a circuit of $M$. Then there
is a function $f_{\ref{spoke1}}(m,q)$ such that, if 
$|T|\geq f_{\ref{spoke1}}(m,q)$, then $M$ has a $\Delta_m$-minor,
each leg of which is a series pair of $M\ba T''$.
\end{lemma}

\begin{proof}
Let $f_{\ref{spoke1}}(m,q)=\max\{f_{\ref{ham}}(2m),
f_{\ref{ham}}(2f_{\ref{spoke0}}(q))\}$.
Assume that $n\geq f_{\ref{spoke1}}(m,q)$. Then
$n\geq f_{\ref{ham}}(2t)$, where $t\geq \max\{m,f_{\ref{spoke0}}(q)\}$.

By Lemma~\ref{ham}, a majority argument, and an appropriate cyclic
ordering of the indices, we see that there is a minor $N$ of $M$,
obtained by flap removal, with $t$ flaps,
$$\{\{s_{i_1},s'_{i_1},s''_{i_1},r_{\phi(i_1)}\},
\{s_{i_2},s'_{i_2},s''_{i_2},r_{\phi(i_2)}\},\ldots,
\{s_{i_t},s'_{i_t},s''_{i_t},r_{\phi(i_t)}\}\}$$
where either
\begin{itemize}
\item[(i)] $i_1<\phi(i_1)<i_2<\phi(i_2)<\cdots<i_t<\phi(i_t)$,
\item[(ii)] $i_1<i_2<\cdots <i_t<\phi(i_t)<\cdots<\phi(i_1)$, or
\item[(iii)] $i_1<i_2<\cdots <i_t<\phi(i_1)<\phi(i_2)<\phi(i_t)$.
\end{itemize}
Assume that case (i) holds. For $k\in\{1,2,\ldots, t-1\}$, let 
$P_k=\{s_{i_k},r_{i_k},
s_{i_k+1},r_{i_k+1},\ldots,s_{i_{k+1}-1},r_{i_{k+1}-1},s'_{i_k},
s''_{i_k}\}$ and let $P_t=E(N)-(P_1\cup P_2\udots P_{t-1})$. 
Note that, for $k\in\{1,2,\ldots,t\}$, 
the set $P_k$ contains the clonal
pair $\{s_{i_k},s'_{i_k}\}$ and that 
$(P_1,P_2,\ldots,P_t)$ is a swirl-like flower in $N$. It now 
follows from Lemma~\ref{get-free-swirl} that in this case, as 
$t\geq m$,  the lemma is satisfied by producing
a $\Delta_m$-minor each of whose legs are 
clonal pairs of the form $\{s_i,s'_i\}$.

Assume that either (ii) or (iii) hold. Let 
$R'=\{r_j\in R:i_1\leq j\leq i_k\}$. 
Consider $N/R'$.
Apart from a single
parallel class at the common basepoint of the series pairs
$\{\{s_{i_1},s'_{i_1}\},\{s_{i_2},s'_{i_2}\},
\ldots,\{s_{i_t},s'_{i_t}\}\}$
in $N/R'\ba\{s_{i_1}'',s''_{i_2},\ldots,s_{i_t}''\}$
that does not affect the argument, the hypotheses of Lemma~\ref{spoke0}
hold for $N/R'$ so that, as $t\geq f_{\ref{spoke0}}(m)$, 
we obtain the contradiction that $M\notin \eq$.
\end{proof}

\begin{proof}[Proof of Lemma~\ref{big-spoke}]
Define $f_{\ref{big-spoke}}(m,q)$ by 
$$f_{\ref{big-spoke}}(m,q)=
\frac{m^{f_{\ref{spoke0}}(q)f_{\ref{spoke1}}(m,q)}-1}
{2(m-1)}.$$
Let $Z$ be a maximal subset of $R$ with the property that
each series pair of $M\ba T''/Z$ is blocked
in $M/Z$. Note that, up to parallel classes at spokes,
$M\ba T''/(T'\cup Z)$ is a whorl.

Let $\lambda$ denote the maximum number of elements contained in a single
parallel class of $M\ba T''/(T'\cup Z)$. The next claim follows
from an elementary bookkeeping argument and Lemma~\ref{spoke0}.

\begin{sublemma}
\label{sub-big-spoke1}
If $\lambda\geq f_{\ref{spoke0}}(q)$, then $M\notin \eq$.
\end{sublemma}

Thus $\lambda<f_{\ref{spoke0}}(m)$. Let $\mu$
denote the rank of $M\ba T''/(T'\cup Z)$. 

\begin{sublemma}
\label{sub-big-spoke2}
The lemma holds if 
$\mu\geq f_{\ref{spoke0}}(q)f_{\ref{spoke1}}(m,q)$. 
\end{sublemma}

\begin{proof}
Let $R'$ denote the rim elements of $M\ba T''/(T'\cup Z)$. 
If $r\in R'$, then, by the definition of
$Z$, there is an element $t_i$ of $T$ such that $\{t_i,t'_i,t''_i,r\}$ is 
a circuit. There may be more than one. Arbitrarily choose one for each 
rim element to define a function $\rho:R'\rightarrow T$. Note that this
function is injective.  
As $\lambda<f_{\ref{spoke0}}(q)$ and 
$\mu\geq f_{\ref{spoke0}}(q)f_{\ref{spoke1}}(m,q)$,
we may choose a subset $R''$ of $R'$ with
$|R''|\geq f_{\ref{spoke1}}(m,q)$ such that, if $r_i$ and $r_j$
are elements of $R''$, then $\rho(r_i)$ is not parallel to $\rho(r_j)$
in $M\ba T''/(T'\cup Z)$. It is now straightforward to take an appropriate
minor of $M$ and apply Lemma~\ref{spoke1} to prove the sublemma.
\end{proof}

We may now assume that $\mu <f_{\ref{spoke0}}(q)f_{\ref{spoke1}}(m,q)$.
Assume that $n\geq f_{\ref{big-spoke}}(m,q)$. Then 
$n\geq \frac{m^\mu-1}{2(m-1)}$. Consider $M/(Z\cup T'')$. 
This matroid has
rank $\mu$. Moreover, each clonal pair $\{t_i,t'_i\}$ is independent
in this matroid. Hence it has at 
least $2n$ parallel classes, that is, it
has at least $\frac{m^\mu-1}{m-1}$ parallel classes. 
By Lemma~\ref{kung} 
we obtain the contradiction that $M\notin\eq$.
\end{proof}

\subsection*{The General Case} 
First observe that Lemma~\ref{clean-block} follows routinely 
from Lemmas~\ref{rim-case} and \ref{big-spoke}. 
We omit the ritual incantation that establishes it.

Lemmas~\ref{rim-case}
and \ref{big-spoke} deal with
one case that arises in Lemma~\ref{whorl-case} we now consider
a complementary case.

\begin{lemma}
\label{elder}
Let $M$ be a $3$-connected matroid with an element $x$ such that 
$M\ba x$ is $3$-connected up to an $n$-element set of 
$M$-clonal series pairs and such that $\co(M\ba x)$ is a whorl.
Then there is a function $f_{\ref{elder}}(q)$ such that,
if $n\geq f_{\ref{elder}}(q)$, then $M\notin\eq$.
\end{lemma}

\begin{proof}
Let $f_{\ref{elder}}(q)=f_{\ref{main1}}(5,q)+4$ and let
$t=f_{\ref{main1}}(5,q)+2$.

Let $r_0$ be a rim element of $\co(M\ba x)$. It may be that
$r_0$ is in a series pair $\{r_0,r'_0\}$ of $M\ba x$. In this case
let $N=M\ba r_0/r'_0$. Otherwise let $N=M\ba r_0$. Observe that
$N\ba x$ has a path $(P_0,P_1,\ldots,P_t)$ of $2$-separations
each step of which contains a series pair of $M\ba x$.
We omit the routine verification of the following claims.
The path $(P_0\cup\{x\},P_1,\ldots,P_t)$ is a well-defined path of
$3$-separations in $N$. If $(A,B)$ is a 2-separation, then
either $A$ or $B$ is contained in either $P_0$ or $P_t$. The
underlying 3-connected matroid $M'$ 
obtained by appropriately removing all but one element of
each maximal 2-separating subset of $M$ contained in 
$P_0$ or $P_t$ is 5-coherent. Each step of the 
path in $M$ induced by $\PP$ 
contains a clonal pair of $M'$.
By Corollary~\ref{main2}, $M\notin \eq$.
\end{proof}

At last we can achieve the purpose of thes section.

\begin{proof}[Proof of Lemma~\ref{whorl-case}]
Recall the hypotheses of the lemma.
Let $l=f_{\ref{clean-block}}(m,q)$ and let
$f_{\ref{whorl-case}}(m,q)=f_{\ref{easy-hyp}}(l,f_{\ref{elder}}(q))$.
Assume that $n\geq f_{\ref{whorl-case}}(m,q)$.

Let $H$ be the hypergraph whose vertex set is the collection of series
pairs of $M\ba J$ and whose edge set is $J$ where 
an element $x\in J$ is
incident with a series pair if $x$ blocks that series pair. 
If $H$ has an edge containing at least $f_{\ref{elder}}(q)$
vertices, then it is easily seen that $M$ has a minor
satisfying the hypotheses of Lemma~\ref{elder} giving
the contradiction that $M\notin \eq$. Thus, by Lemma~\ref{easy-hyp}
and the definition of $f_{\ref{whorl-case}}(m,q)$, there is a
set $S=\{\{s_1,s_1'\},\{s_2,s'_2\},\ldots,\{s_l,s_l'\}\}$ 
of series pairs and a collection
$T=\{t_1,t_2,\ldots,t_l\}$ of elements of $J$ such that,
if $i,j\in\{1,2,\ldots,l\}$, then $t_i$ blocks
$\{s_i,s_j\}$ if and only if $i=j$.

We ignore the elements of $J$ not in $T$ and focus on $M\ba (J-T)$.
Let $\{u,u'\}$ be a series pair of $M\ba J$ that is not in $S$.
Consider $M\ba (J-T)/u'$. Assume, for a contradiction, that,
for some $i\in\{1,2,\ldots,l\}$, the series pair $\{s_i,s'_i\}$ is not
blocked in this matroid. Then we have
$t_i\in\cl_{M\ba (J-T)\ba u'}(\{s_i,s'_i\})$, that is,
$t_i\in\cl_M(\{s_i,s'_i,u'\})$, so that 
$u'\in\cl_M(\{s_i,s'_i,t_i\})$. But $\{u,u'\}$ is a clonal pair of
$M$, so that $\{u,u'\}\subseteq \cl_M(\{s_i,s'_i,t_i\})$. 
However, this implies that
$\lambda_{M\ba J}(\{u,u',s_i,s'_i\})=1$, contradicting the fact that
$M\ba J$ is 3-connected up to series pairs.

Therefore each series pair of $S$ is blocked in $M\ba (J-T)/u'$.
It follows from this fact and an obvious induction
that $M$ has a minor $N$ with a set $T$ of elements that is 
a cleanly-blocked coextended whorl, where the series pairs of
$N\ba T$ are series pairs of $M\ba J$. Moreover,
$N\ba T$ has at least $l=f_{\ref{clean-block}}(m,q)$ series pairs.
By Lemma~\ref{clean-block} $N$, and hence $M$, has a 
$\Delta_m$-minor each leg of which is a series pair of
$N\ba T$, and hence of $M\ba T$.
\end{proof}

\section{The $M(K_{3,n})$ Case}

We now consider the cases where the underlying matroid is 
$M(K_{3,n})$. The goal of this section is to prove

\begin{lemma}
\label{k3n-case}
Let $M$ be a $3$-connected matroid with a coindependent 
set $J$ such that the following hold.
\begin{itemize}
\item[(i)] $M\ba J$ is $3$-connected up to a set of $n$ 
series pairs that are $M$-clonal.
\item[(ii)] $\co(M\ba J)\cong M(K_{3,t})$ for some integer $t$..
\end{itemize}
Then there is a function $f_{\ref{k3n-case}}(q)$ such that, if
$n\geq f_{\ref{k3n-case}}(q)$, then $M\notin \eq$.
\end{lemma}

As usual we develop a series of lemmas. The next lemma
is just a special case of Lemma~\ref{free-link} that we recall
here for convenience.

\begin{lemma}
\label{get-clonal-flower}
Let $M$ be a connected matroid with a paddle $\PP=(P_1,P_2,\ldots,P_n)$
such that, for all $i\in\{1,2,\ldots,n\}$ the petal $P_i$ contains
a $\PP$-strong clonal pair. Then $M$ contains a 
$U_{2,2n}$-minor.
\end{lemma}

To make life easier in this section we develop some local terminology.
Let $N$ be a connected
matroid such that $\co(N)\cong M(K_{3,t})$ for some $t\geq 3$.
Just as there is a unique flower of order $n$ in $\co(N)$, so too
is there a unique flower of order $n$ in $N$ and, extending existing terminology,
we say that this is the {\em canonical flower associated with $N$}.

A matroid $M$ is a {\em blocked coextended $M(K_{3,t})$}
\index{blocked coextended $M(K_{3,t})$} 
of {\em order}
$n$ with {\em blocking set $J$} if the following hold.
\begin{itemize}
\item[(i)] $M\ba J$ is 3-connected up to a 
set of $n$ series pairs that are $M$-clonal.
\item[(ii)] $\co(M\ba J)\cong M(K_{3,t})$.
\item[(iii)]  If $j\in J$, then there is a series pair of $M\ba J$
that is blocked by $j$.
\end{itemize}
Let $M$ be a blocked coextended $M(K_{3,t})$
with blocking set $J$ and let $\PP=(P_1,P_2,\ldots,P_t)$
be the canonical flower associated with $M\ba J$. For $x\in J$,
and $i\in\{1,2,\ldots,t\}$, we say that 
$x$ is {\em incident} with $P_i$ if
either $x$ blocks a series pair in $P_i$, or $x$ blocks $P_i$.
Note that the former case is redundant 
unless $x\in \cl(P_i)$. The above
incidence relation defines a hypergraph whose vertex set is 
the set of petals of $\PP$ and whose edge set is $J$. We denote
this hypergraph by $H(\PP,J)$.

\begin{lemma}
\label{k3n-case-1}
Let $M$ be a blocked coextended $M(K_{3,t})$ with blocking set
$J$ and canonical associated flower $(P_1,P_2,\ldots,P_t)$. Let
$\{e_i,e'_i\}$ be a series pair of $M\ba J$ that is contained in
$P_i$. Say that the element $x$ of $J$ blocks $P_j$, where 
$P_j\neq P_i$. Then $x$ blocks $P_j$ in $M/e$.
\end{lemma}

\begin{proof}
Assume otherwise. Then $x\in\cl_M(P_j\cup\{e\})$ so that
(i) $e\in\cl_M(P_j\cup\{x\})$ and (ii) $x\in\cl_{M/e}(P_j)$.
From (i) and the fact that $\{e,e'\}$ are clones we deduce that
$\{e,e'\}\subseteq \cl_M(P_j\cup\{x\})$. By this fact and (ii),
we see that $e'\in\cl_{M/e}(P_j)$, contradicting the fact that
$\co(M\ba J)\cong M(K_{3,t})$.
\end{proof}

\begin{lemma}
\label{k3n-case-2}
Let $M$ be a blocked coextended $M(K_{3,t})$ with blocking set
$J$ and canonical flower $\PP=(P_1,P_2,\ldots,P_t)$. If $H(\PP,V)$
has an edge that is incident with at least $q$ vertices, then
$M$ has a $\Lambda_q$-minor.
\end{lemma}

\begin{proof}
Assume that $x$ is incident with at least $q$ vertices.
We may assume that $J=\{x\}$ and, by Lemma~\ref{k3n-case-1},
we may assume that $M\ba x$ has a single series pair 
$\{a_1,a'_1\}$ that is contained in $P_1$ and is blocked by $x$. 
We may further assume
that $x$ blocks every petal of $\PP$.

Let $A=\{a_1,a_2,\ldots,a_n\}$, $B=\{b_1,b_2,\ldots,b_n\}$ and
$C=\{c_1,c_2,\ldots,c_n\}$ be a partition of 
$(P_1-\{a_1'\})\cup P_2\udots P_t$ such
that the following hold: $P_1=\{a_1,a'_1,b_1,c_1\}$; for each 
$i\in\{2,3,\ldots,t\}$, 
we have $P_i=\{a_i,b_i,c_i\}$; 
and $A$, $B$, and $C$ are cocircuits of $M\ba x/a'_1$. In other words,
$A$, $B$ and $C$ are stars of the underlying $K_{3,t}$. Consider 
$N=\si(M/b_1,c_1)$. The parallel pairs of 
$M/b_1,c_1$ are $\{b_2,c_2\},\{b_3,c_3\},\ldots,\{a_t,b_t\}$,
so we may assume that $E(N)=A\cup (B-\{b\})\cup\{a_1',x\}$. Note that
$N\ba x\cong M(K_{2,t})$, that for all $i\in\{2,3,\ldots,t\}$, the
set $\{a_1,a'_1,a_i,b_i\}$ is a circuit of $N$, 
and that each series pair
of $N\ba x$ is blocked by $x$. It follows that $N/x$ is a spike.
But the clonal pair $\{a_1,a'_1\}$ is a leg of this spike.
By Lemma~\ref{clones-in-spike},
the only spikes with clonal pairs are free spikes. Hence 
$N/x\cong\Lambda_t$. As $t$ is the degree of $x$ and 
$t\geq q$, the lemma follows.
\end{proof}

\begin{lemma}
\label{k3n-case-3}
Let $M$ be a blocked coextended $M(K_{3,t})$ of order $n$ with
blocking set $J$ and canonical flower $\PP=(P_1,P_2,\ldots,P_t)$
and let $l$ be an integer. Then there is a function
$f_{\ref{k3n-case-3}}(l,q)$ such that, if each edge of 
$H(\PP,J)$ has at most $l$ vertices and 
$n\geq f_{\ref{k3n-case-3}}(l,q)$, then $M$ has a 
$U_{2,q+2}$-minor.
\end{lemma}

\begin{proof}
Let $w=\lceil{\frac{q+2}{2}}\rceil$ and  
let $f_{\ref{k3n-case-3}}(1,q)=w$. For $l>1$ let 
$f_{\ref{k3n-case-3}}(l,q)=f_{\ref{easy-hyp2}}(s,l)$
where $s=w(q+2)(f_{\ref{k3n-case-3}}(l-1,q))+w(q+3)$.

Assume that $n\geq f_{\ref{k3n-case-3}}(l,q)$. Say that
$l=1$. Up to labels we may assume that $P_1,P_2,\ldots,P_w$ all contain
series pairs of $M\ba J$. But then, 

for each $i\in\{1,2,\ldots,w\}$,
there is an element $x_i\in J$ such that $x_i$ blocks a series
pair in $P_i$. By the fact that $l=1$ we have
$x_i\in\cl(P_i)$. Thus 
$(P_1\cup\{x_1\},P_2\cup\{x_2\},\ldots,P_w\cup\{x_w\})$ 
is a paddle in 
$M|(P_1\udots P_n\cup\{x_1,\ldots,x_w\}$
each petal of which contains a clonal pair
that is strong relative to this flower. By 
Lemma~\ref{get-clonal-flower} $M$ has a $U_{2,q+2}$-minor.

Assume that $l>1$ and, for induction, assume that the lemma
holds for smaller values of $l$. Up to labels of the petals of 
$\PP$, there is, by Lemma~\ref{easy-hyp2}, a subset
$J'=\{x_1,x_2,\ldots,x_s\}$ of $J$ and a set
$\{\{a_1,a'_1\},\{a_2,a'_2\},\ldots,\{a_s,a'_s\}\}$ of series pairs of 
$M\ba J$ such that $\{a_i,a'_i\}\subseteq P_i$ for 
$i\in\{1,2,\ldots,s\}$,  with the properties that
\begin{itemize}
\item[(i)] $x_i$ blocks $\{a_i,a'_i\}$ 
for all $i\in\{1,2,\ldots,s\}$; and
\item[(ii)]  if $i\in\{1,2,\ldots,s\}$, $j\in\{1,2,\ldots,t\}$,
and $x_i$ blocks $P_j$, then either $j=i$ or $j>s$.
\end{itemize}
It may be that there exist $i\in\{1,2,\ldots,s\}$ and
$j\in\{s+1,s+2,\ldots,t\}$ such that, for some 
$z\in P_j$, the set $\{a_i,a'_i,x_i,z\}$
is a circuit. For terminology restricted to this proof, we say in
this case that $x_i$ is a {\em fragile} element of $J'$
and $z$ is a {\em fragile} element of $P_j$.

Let $y$ be an element of $P_j$ for some 
$j\in\{s+1,s+2,\ldots,t\}$
and assume that $y$ is not fragile. Then, for all 
$i\in\{1,2,\ldots,s\}$, the series pair $\{a_i,a'_i\}$
of $[M/z\ba(P_j-\{z\})]\ba J$ is blocked by $x_i$. It follows
from this that if the lemma holds in the case that every element of 
$P_{s+1}\cup P_{s+2}\udots P_t$ is fragile, then it holds in general.
Thus we may assume that every element of this set is fragile.

Let $K$ be the set of fragile elements of $J'$.

\begin{sublemma}
\label{k3n-case-3.1}
The lemma holds if $|K|>w(q+2)$.
\end{sublemma}

\begin{proof}
Note that the elements of $K$ are all 2-element edges
of $H(\PP,J)$. Indeed the subhypergraph induced by $K$
is a union of stars. Thus, it either has a vertex of
degree $q+2$ or a matching of size $w$. The latter case implies
that $M$ has a $3$-connected minor with a paddle containing $w$
petals, each petal of which contains a clonal pair. In
this case, by Lemma~\ref{get-clonal-flower}, $M$ has a
$U_{2,q+2}$-minor. Consider the former case. Assume that $P_j$
has degree $q+2$ in $H(\PP,J)$. 
Up to labels we may assume that the members of
$\{x_1,x_2,\ldots,x_{q+2}\}$ are fragile elements 
incident with $P_j$ in $H(\PP,J)$. Say $i\in\{1,2,\ldots,q+2\}$.
Then, as $\sqcap(\{a_i,a'_i\},P_j)=0$, there is at most one
element of $P_j$ in the closure of $\{a_i,a'_i,x_i\}$.
Thus, again up to labels, we may assume that there is an element
$z$ in $P_j$ such that 
$z\notin\cl(\{a_i,a'_i,x_i\})$ for $i\in\{1,2,\ldots,w\}$.
This means that, if $i\in\{1,\ldots,w\}$, the series pair
$\{a_i,a'_i\}$ of $M/z\ba J$ is blocked by $x_i$. But
$x_i\in\cl(P_i\cup P_j)$ and it follows from properties of $M(K_{3,n})$
that $x_i\in\cl_{M/z}(P_i)$. Indeed
$(M/z)|(P_1\cup P_2\udots P_w\cup\{x_1,\ldots,x_w\})$ is a 
3-connected matroid with a paddle
$(P_1\cup\{x_1\},P_2\cup\{x_2\},\ldots,P_w\cup\{x_w\})$,
each petal of which contains a clonal pair. Again, by
Lemma~\ref{get-clonal-flower} $M$ has a $U_{2,q+2}$-minor.
\end{proof}

We may now assume that $J'$ has at most $w(q+2)$ fragile elements.
If $J'$ contains at least $w$ elements of degree 1 in
$H(\PP,J)$, then we again find a $U_{2,q+2}$-minor, so we may
assume that $J'$ has at least $s-(w(q+2)+w)=s-w(q+3)$
elements that are not fragile with degree at least 2.
But, each element of $\{P_{s+1},P_{s+2},\ldots,P_t\}$ is incident
with a fragile element and, somewhat crudely, we see that
$t-s\leq w(q+2)$. From these fact and the definition of
$s$, we deduce that, for some $j\in\{s+1,s+2,\ldots,t\}$,
at least $f_{\ref{k3n-case-3}}(l-1,q)$ 
members of $J'$ are incident with 
$P_j$. Up to labels we may assume that the members
of $J''=\{x_1,x_2,\ldots,x_p\}$ are incident with $P_j$,
where $p\geq f_{\ref{k3n-case-3}}(l-1,q)$. Say 
$z\in P_j$. Consider 
$M/z|(P_1\cup P_2\udots P_p\cup J'')$. This is a blocked coextended
$M(K_{3,t})$ of order $f_{\ref{k3n-case-3}}(l-1,q)$ with 
blocking set $J''$ such that each edge
of $H((P_1,P_2,\ldots,P_p),J'')$ is incident with at most
$l-1$ vertices. By the induction assumption this matroid,
and therefore $M$, has a $U_{2,q+2}$-minor.
\end{proof}

\begin{proof}[Proof of Lemma~\ref{k3n-case}]
Let $f_{\ref{k3n-case}}(q)=f_{\ref{k3n-case-3}}(q,q)$. 
Let $\PP$ be the canonical flower associated with $M\ba V$.
If $H(\PP,V)$ has an edge incident with $q$ vertices,
then, by Lemma~\ref{k3n-case-2}, $M\notin \eq$. 
Otherwise $M\notin \eq$ by Lemma~\ref{k3n-case-3}.
\end{proof}

\section{The Spike Case}

We now turn to spikes. In this case we prove

\begin{lemma}
\label{big-spike-case}
Let $M$ be a $3$-connected matroid with a 
coindependent set $J$ such that the
following hold.
\begin{itemize}
\item[(i)] $M\ba J$ is $3$-connected up 
to a set of $n$ series pairs that are $M$-clonal.
\item[(ii)] $\co(M\ba J)$ is a spike.
\end{itemize}
Then there is a function $f_{\ref{big-spike-case}}(q)$  
such that, if $n\geq f_{\ref{big-spike-case}}(q)$,
then $M$ has a $\Lambda_q$-minor.
\end{lemma}

Lemma~\ref{big-spike-case} is 
a routine consequence of the next
lemma, which is indeed somewhat stronger.
Let $x$ be an element of the matroid $M$ and $\PP$
be a flower in $M\ba x$. Then $\PP$ is 
{\em well-blocked}
\index{well-blocked flower} 
by $x$ if 
$x$ blocks $(P,Q)$ whenever 
$(P,Q)$ is a $3$-separation of $M\ba x$ such that
both $P$ and $Q$ are unions of at least two petals 
of $\PP$.

\begin{lemma}
\label{small-spike-case}
Let $M$ be a $3$-connected matroid with a triad $\{x,a_1,a'_1\}$ such that $M\ba x/a'_1$ is a rank-$n$ spike and $\{a_1,a'_1\}$
is an $M$-clonal series pair in $M\ba x$.
Then, if $n\geq q^2-1$, the matroid $M$ has a $\Lambda_q$-minor.
\end{lemma}

\begin{proof}
Let $(P_1,P_2,\ldots,P_n)$ be the spike-like flower
associated with $M\ba x$, where $\{a_1,a'_1\}\subseteq P_1$
and $(P_1-\{a_1'\},P_2,\ldots,P_n)$ is a spike in
$M\ba x/a_1'$. Say $P_1=\{a_1,a'_1,b_1\}$ and,
for $i\in\{2,3,\ldots,n\}$, let $P_i=\{a_i,b_i\}$.

\begin{sublemma}
\label{small-spike1}
If $\PP$ is well-blocked by $x$ and $n\geq q+1$, then
$M$ has a $\Lambda_q$-minor.
\end{sublemma}

\subproof 
First observe that $M\ba x/b_1\ba P_n\cong M(K_{2,n-1})$.
We omit the easy rank calculations that establish this fact.
Let $M'=M/b_1\ba P_n$. 
We now show that every member
of $(P_1-\{a'_1\},P_2,\ldots,P_{n-1})$ in $M'\ba x$ is 
blocked by $x$ in $M'$.

Consider $P_1-\{x\}=\{a_1,a'_1\}$. If $x$ does not block
$P_1-\{x\}$, then $x\in\cl_{M'}(P_2\cup P_3\udots P_{n-1})$,
so that $x\in\cl_M(E(M)-\{a_1,a'_1\})$, contradicting the fact
that $x$ blocks $\{a_1,a'_1\}$ in $M$.
Say $i\in\{2,3,\ldots,n-1\}$ and assume that
$x$ does not block $P_i$. Then 
$x\in\cl_{M'}((P_1-\{a_1'\})\cup P_2\udots P_{i-1}\cup P_{i+1}
\udots P_{n-1})$ so that
$x\in\cl_{M\ba P_n}(P_1\cup P_2\udots P_{i-1}\cup P_{i+1}
\udots P_{n-1})$ and therefore $x$ does not block $P_i\cup P_n$
in $M$, contradicting the assumption that $\PP$ is 
well-blocked by $x$ in $M$. Therefore every member
of $(P_1-\{a'_1\},P_2,\ldots,P_{n-1})$ in $M'\ba x$ is 
indeed blocked by $x$ in $M'$.

We now have an anemone in $M'$. 
But, if $P$ and $P'$ are petals of this flower, then
$\sqcap(P,P')=1$. It follows that the flower is spike-like,
that is,
$M'$ is a spike with cotip $x$. Thus $M'/x$ is a spike.
But $\{a_1,a'_1\}$ is a clonal pair, so 
$M'/x\cong \Lambda_{n-1}$ by Lemma~\ref{clones-in-spike}. As
$n-1\geq q$, the claim follows.
\end{proof}

Let $\calP$ and $\calQ$ be maximal sets of petals such that
$x\in\cl(E(M)-\cup_{P\in\calP}(P))$ and
$x\in\cl(E(M)-\cup_{Q\in\calQ}(Q))$. 
If $\calP\cap\calQ\neq\emptyset$, then, by Lemma~\ref{modular},
we deduce that 
$x\in\cl(E(M)-(\cup_{P\in\calP}(P)\cup \cup_{Q\in\calQ}(Q)))$
and we have contradicted the assumption that $\calP$
and $\calQ$ were maximal with the given property.
Observe that $x\notin\cl(E(M)-(P_1\cup\{x\}))$,
but, if $i\in\{2,3,\ldots,n\}$,
then $x\in\cl(E(M)-(P_i\cup\{x\}))$. This establishes the
next claim.

\begin{sublemma}
\label{small-spike2}
There is a partition $(\{P_1\},\calP_2,\ldots,\calP_s)$
of the petals of $\PP$ such that, for all $i\in\{2,3,\ldots,s\}$,
the set $\calP_i$ is a maximal collection of petals with the 
property that $x\in\cl_M(E(M)-\cup_{P\in\calP_i}(P))$.
\end{sublemma}

Consider the partition given by \ref{small-spike2}.

\begin{sublemma}
\label{small-spike3}
If $s\geq q+1$, then $M$ has a $\Lambda_q$-minor.
\end{sublemma}

\subproof
Up to labels we may assume that $(P_2,P_3,\ldots,P_s)$
is a transversal of $(\calP_2,\calP_3,\ldots,\calP_s)$.
Let $M'=M\ba\{a_{s+1},a_{s+2},\ldots,a_n\}/\{b_{s+1},b_{s+2}
\ldots,b_n\}$. Then $(P_1,P_2,\ldots,P_s)$ is a spike-like
flower in $M'\ba x$. Elementary rank calculations show
that $(P_1,P_2,\ldots,P_s)$ is well-blocked by $x$ in $M'$
and that $\{a_1,a'_1\}$ is a series pair in $M'\ba x$ that is 
blocked by $x$. By \ref{small-spike1}, $M'$ has a 
$\Lambda_{s-1}$-minor and the sublemma follows.
\end{proof}

\begin{sublemma}
\label{small-spike4}
If $i\in\{2,3,\ldots,s\}$ and
$|\calP_i|\geq q-1$, then $M$ has a $\Lambda_q$-minor.
\end{sublemma}

\subproof
Say $|\calP_i|=l$. Up to labels
we may assume that $\calP_i=P_2\cup P_3\udots P_{l+1}$.
Let $P_1'=P_1\cup P_{l+2}\cup P_{l+3}\udots P_n$.
Observe that $\PP'=(P'_1,P_2,\ldots,P_{l+1})$ is a spike-like
flower in $M$ and that $\{a_1,a'_1\}$ is a $\PP'$-strong
clonal pair in $P'_1$. It follows easily from Tutte's Linking
Lemma and the fact that $\{a_1,a'_1\}$ is a clonal pair, that
$M$ has a 3-connected minor on
$\{a_1,a'_1\}\cup P_2\cup P_3\udots P_{l+1}$ such that 
$(\{a_1,a'_1\},P_2,\ldots,P_{l+1})$ is a spike in $M''$.
Now $\{a_1,a'_1\}$ is a clonal pair and $M''\cong\Lambda_q$
by  Lemma~\ref{clones-in-spike}.
\end{proof}

As $n\geq (q+1)(q-1)$, either \ref{small-spike3}
or \ref{small-spike4} applies and the lemma follows.
\end{proof}

\begin{corollary}
\label{small-spike-cor}
Let $M$ be a $3$-connected matroid with an element
$x$ such that $M\ba x$ is $3$-connected up to a nonempty
set of $M$-clonal series pairs and that
$\co(M\ba x)$ is a rank-$n$ spike.
If $n\geq q^2-1$, then $M$ has a $\Lambda_q$-minor.
\end{corollary}

\begin{proof}
The goal is to contract elements from $M$ that are 
in series pairs of $M\ba x$ keeping a matroid satisfying the 
hypotheses of the lemma until we have only one series pair remaining,
in which case we can apply Lemma~\ref{small-spike-case}.
No difficulties arise with this strategy unless we have two
clonal pairs $\{s_i,s'_i\}$ and $\{s_j,s'_j\}$
such that $\{s_j,s'_j\}$ is not blocked in $M/s_i$,
that is, if $x\in\cl_{M/s_i}(\{s_j,s'_j\})$. 
In this case $s_i\in\cl_M(\{s_j,s'_j,x\})$ and,
as $\{s_i,s'_i\}$ is a clonal pair, 
$\{s_i,s'_i\}\subseteq \cl_M\{s_j,s'_j,x\}$.
Therefore $r_M(\{s_i,s'_i,s_j,s'_j\})=3$. But this means
that $\{s'_i,s'_j\}$ are parallel in $M/s_i,s_j$ contradicting the
fact that $\co(M)$ is 3-connected.
It follows that our strategy is fine and that the corollary 
indeed follows from Lemma~\ref{small-spike-case}.
\end{proof}

Finally we observe that Lemma~\ref{big-spike-case} is an
almost immediate consequence of Corollary~\ref{small-spike-cor}.

\section{The $M^*(K_{3,n})$ Case}

We now come to the final case. This seems to be more
difficult that the previous cases --- although we may be missing
an easy argument. Rather than solve it independently of the other
cases we combine it with the fact that the previous cases have
been resolved. We begin by developing some specialist terminology.

Let $M$ be a matroid and $l$ be an integer. A set $X$ of elements of
$M$ is a {\em star}
\index{star} 
of {\em order} $l$ if:
\begin{itemize}
\item[(i)] $X$ is closed;
\item[(ii)] $X$ is the union of $l$ 2-element series classes of $M$; 
and
\item[(iii)] $X\cap E(M)$ is a parallel set in $\co(M)$.
\end{itemize}
Evidently a star is {\em maximal} if it is not contained in
a larger star.
Insisting that $X$ be closed 
is just a non-triviality condition that we do not really
need, but it helps to make certain arguments cleaner.
The matroid $M$ is $3$-connected {\em up to stars},
\index{$3$-connected up to stars} 
if whenever
$(X,Y)$ is a $2$-separation of $M$, either $X$ or $Y$ is 
a star. 

A matroid $M$ is a {\em star extended}
\index{star extended $\kstar$} 
$\kstar$ if
$M$ is 3-connected up to stars, and $\si(\co(M))\cong \kstar$. 
Let $M$ be a star extended $\kstar$. Then there is a flower
$\PP=(P_1,P_2,\ldots,P_t)$ such that
$(P_1\cap E(\si(\co(M)),P_2\cap E(\si(\co(M)),\ldots,
P_t\cap E(\si(\co(M)))$ is the canonical flower associated with
$\kstar$ in $\si(\co(M))$. We call this flower the
copaddle {\em associated} with $M$.

Let $M$ be a star extended $\kstar$, then 
$M$ is an $(n,l)$-{\em copaddle}
\index{$(n,l)$-copaddle}  
if $n\leq t$ and the flower
$\PP=(P_1,P_2,\ldots,P_t)$ associated with $M$ has the property
that $P_i$ contains a unique star of order $l$ for each
$i\in\{1,2,\ldots,n\}$.

Let $M$ be a connected matroid.
The triple  $(M,\PP,J)$ is a {\em blocked} $(n,l)$-copaddle,
\index{blocked $(n,l)$-copaddle}
with {\em blocking set} $J$ and {\em associated flower} 
$\PP=(P_1,P_2,\ldots,P_t)$
if $J$ is a coindependent set of $M$ such that $M\ba J$
is an $(n,l)$-copaddle  with associated flower $\PP$
and the following properties hold.
\begin{itemize}
\item[(i)] Every series pair of $M\ba J$ in
$P_1\cup P_2\udots P_n$ is blocked in $M$.
\item[(ii)] All series pairs of
$M\ba J$ are $M$-clonal. 
\end{itemize}
If we say that a matroid $M$ is a 
blocked $(n,l)$-copaddle, we mean that there exist $\PP$ and $J$
such that $(M,\PP,J)$ is a blocked $(n,l)$-copaddle.

Let $M$ be a blocked $(n,l)$-copaddle with blocking set
$J$ and associated copaddle $\PP=(P_1,P_2,\ldots,P_t)$.
If we say that a set $S\subseteq P_1\cup P_2\udots P_t$
is a series pair or a star of $\PP$ we will always mean that
it has this property in $M\ba J$ and if we say that an element
$x$ of $J$ blocks $S$ we will always mean that 
this set is blocked by $x$ in $M\ba (J-\{x\})$.

For the time being we focus on blocked $(n,1)$-copaddles
so that our stars are series pairs.
The more general structures will eventually play a role.
We begin by showing that, given a sufficiently large 
3-connected clonal matroid, we can always either achieve
a desirable outcome or produce a blocked $(n,1)$-copaddle.

\begin{lemma}
\label{star1}
Let $M$ be a $3$-connected clonal matroid in $\eq$. Then there is a 
function $f_{\ref{star1}}(s,m,q)$ such that, if 
$|E(M)|\geq f_{\ref{star1}}(s,m,q)$, then at
least one of the following holds.
\begin{itemize}
\item[(i)] $M$ has a clonal $\Delta_m$-minor.
\item[(ii)] Either $M$ or $M^*$
has a minor $N$ with an associated triple
$(N,\PP,J)$ such that $(N,\PP,J)$ is a blocked
$(s,1)$-copaddle. Moreover, all of the series pairs of
$N\ba J$ are $M$-clonal.
\end{itemize}
\end{lemma}

\begin{proof}
Let $\mu=\max\{3s,f_{\ref{whorl-case}}(m,q),f_{\ref{k3n-case}}(m,q),
f_{\ref{big-spike-case}}(m,q)\}$. Let
$f_{\ref{star1}}(m,q)=f_{\ref{clonal-series}}(\mu,q)$.

By Lemma~\ref{clonal-series}, $M$ or $M^*$
has a 3-connected minor
$N$ with a coindependent set $J$ such that one  the 
following holds.
\begin{itemize}
\item[(i)] $N\ba J$ is $3$-connected up to a set of at least
$\mu$ series pairs.
\item[(ii)] Each series pair of $N\ba J$ is $M$-clonal.
\item[(iii)] $\co(N\ba J)$ is either a spike, a whorl, or
for some integer $l$, is isomorphic to 
$M(K_{3,t})$ or $M^*(K_{3,t})$.
\end{itemize}
If $\co(N\ba J)$ is a whorl, then, by Lemma~\ref{whorl-case},
$M$ has an clonal $\Delta_m$-minor and the lemma holds.
If $\co(N\ba J)$ is a spike or is isomorphic to 
$M(K_{3,t})$, then, by Lemma~\ref{big-spike-case}
or Lemma~\ref{k3n-case}, we contradict the assumption that
$M\in \eq$. Hence $\co(N\ba J)\cong \kstar$.
In this case it is easily seen that $N$ has a minor
$N'$ that has a representation as a blocked 
$(\mu/3,1)$-copaddle where all of the series pairs
in the copaddle associated with $N'$ are $M$-clonal.
By the definition of $\mu$, we have $\mu/3\geq s$, 
and the lemma follows.
\end{proof}

We now define a more highly structured type of 
blocked $(n,1)$-copaddle. Let $(M,\PP,J)$ be a blocked 
$(n,1)$-copaddle where $\PP=(P_1,P_2,\ldots,P_t)$.
For all $i\in\{1,2,\ldots,n\}$ there is a unique series
pair in $P_i$ which we denote by $\{a_i,a'_i\}$. The triple
$(M,\PP,J)$ is an $n$-{\em blockage}
\index{$n$-blockage} 
if there is a function 
$\phi:\{1,2,\ldots,n\}\rightarrow J$
such that, for all $i\in\{1,2,\ldots,n\}$, the
element $\phi(i)$ blocks $\{a_i,a'_i\}$ and
$\phi(i)$ blocks $P_{n+1}\cup P_{n+2}\udots P_t$ but 
$\phi(i)$ does not block $P_j$
for any $j\in\{1,2,\ldots,i-1,i+1,\ldots,n\}$.
We call the function $\phi$ the {\em blocking function}
\index{blocking function}
of $(M,\PP,J)$.
The remainder of this section is devoted to producing a
large $m$-blockage from a very large blocked $(n,1)$-copaddle.

\begin{theorem}
\label{star3}
Let $(M,\PP,J)$ be a blocked $(n,1)$-copaddle in $\eq$. Then there
is a function $f_{\ref{star3}}(m,q)$ such that, if 
$n\geq f_{\ref{star3}}(m,q)$, then there is a permutation
$\QQ$ of the petals of $\PP$ and a partition
$(J',J'')$ of $J$ such that $(M\ba J',\QQ,J'')$
is an $m$-blockage.
\end{theorem}

This is a technical result and probably does not deserve to
be called a theorem, but doing so helps clarify the 
organisation of this section.
Before proving Theorem~\ref{star3}, we note a corollary of it.

\begin{corollary}
\label{star-col}
Let $M$ be a $3$-connected clonal matroid in $\eq$.
Then there is a function $f_{\ref{star-col}}(s,m,q)$
such that, if $|E(M)|\geq f_{\ref{star-col}}(s,m,q)$,
then either
\begin{itemize}
\item[(i)] $M$ has a clonal $\Delta_m$-minor, or
\item[(ii)] $M$ or $M^*$ has a minor $M'$ with an associated
triple $(M',\PP,J)$ such that $(M',\PP,J)$
is an $s$-blockage. Moreover, all of the series pairs
in $M'\ba J$ are $M$-clonal.
\end{itemize}
\end{corollary}

\begin{proof}
Let 
$f_{\ref{star-col}}(s,m,q)=
f_{\ref{star1}}(f_{\ref{star3}}(s,q),m,q)$.
With this function, the corollary follows immediately from
Lemma~\ref{star1} and Theorem~\ref{star3}.
\end{proof}

Let $\calP$ be a collection of subsets of the ground set of a matroid.
To simplify life we make some convenient abbreviations.
We will often say $\cl(\calP)$ to refer to $\cl(\cup_{P\in\calP}P)$.
This convention is extended to $r(\calP)$, $\lambda(\calP)$,
and so on. 

\begin{lemma}
\label{k*2a}
Let $\PP$ be a copaddle in the matroid $M$ and let $\calP$ and
$\calQ$ be disjoint subsets of petals of $\PP$ whose union does
not contain all the petals of $\PP$. Then
$r(\calP)+r(\calQ)=r(\calP\cup\calQ)$.
\end{lemma}

Part (i) of the next lemma is true for any petal of a copaddle
that has no elements in the coguts of a petal. This is clearly true
in our case.

\begin{lemma}
\label{k*2b}
Let $M$ be an  
$(n,l)$-copaddle with associated flower $(P_1,P_2,\ldots,P_t)$.
Then the following hold. 
\begin{itemize}
\item[(i)] $M|P_i$ is connected for all $i\in\{1,2,\ldots,t\}$.
\item[(ii)] If $i\in\{1,2,\ldots,n\}$
then $M/P_i$ is an $(n-1,l)$-copaddle with
associated  flower $(P_1,\ldots,P_{i-1},P_{i+1},\ldots,P_t)$.
\item[(iii)] If $i\in\{n+1,n+2,\ldots,t\}$, then
$M/P_i$ is an $(n,l)$-copaddle with
associated  flower $(P_1,\ldots,P_{i-1},P_{i+1},\ldots,P_t)$.
\end{itemize}
\end{lemma}

The next lemma is also elementary.

\begin{lemma}
\label{k*2x}
Let $M$ be a matroid with an element $x$ such that
$M\ba x$ is connected with a copaddle $(P_1,P_2,\ldots,P_n)$.
Say that $x$ blocks $P_2$ in $M$, but not in $M/P_1$,
then
\begin{itemize}
\item[(i)] $x\in \cl(P_1\cup P_2)$.
\item[(ii)] If $\mathcal P$ is any set of petals that contains
$P_1$ but not $P_2$, then $x$ blocks $\cup_{P\in\calP}(P)$.
\end{itemize}
\end{lemma}

The straightforward properties described above enable
us to prove the next useful lemma.

\begin{lemma}
\label{k*2c}
Let $M$ be a matroid with an element $x$ such that $M\ba x$
is an $(n,l)$-copaddle with associated flower
$\PP=(P_1,P_2,\ldots,P_t)$. 
Let $\{p,p'\}\subseteq P_1$ be an $M$-clonal series pair of 
$M\ba x$ and let $\PP'=(P_1-\{p'\},P_2,\ldots,P_t)$.
Let $(\calP,\calQ)$ be a partition of the petals
of $\PP'$ where $|\calP|\geq 2$ and $(P_1-\{p'\})\in \calP$.
Then $x$ blocks $(\calP\cup\{p'\},\calQ)$ in $M$ if and only
if $x$ blocks $(\calP,\calQ)$ in $M/p'$.
\end{lemma}

\begin{proof}
For the not completely trivial direction assume that $x$ blocks
$(\calP\cup\{p'\},\calQ)$ in $M$ and, for a contradiction,
that $x$ does not block $(\calP,\calQ)$ in
$M/p'$. In this case it is clear that 
$x\in\cl_{M/p'}(\calQ)$, so that $x\in\cl_M(\calQ\cup\{p'\})$.
But $x\notin\cl_M(\calQ)$. Hence $p'\in\cl_M(\calQ\cup\{x\})$.
By Lemma~\ref{k*2b}, $p'\in\cl_M(P_1-\{p'\})$. By
Lemma~\ref{k*2a}, $\sqcap_M(\calQ\cup\{x\},P_1-\{p'\})\leq 1$.
Now, by Lemma~\ref{freedom1}, $p'$ is fixed in $M$,
contradicting the fact that the pair
$\{p,p'\}$ is $M$-clonal.
\end{proof}

The previous four lemmas will be used freely without reference
throughout this section.

\subsection*{The Well-blocked Case}
Let $M$ be a matroid with an element $x$ such that $M\ba x$
is connected with a copaddle $\PP$. Then $\PP$ is 
{\em well-blocked} by $x$ if, $P_i\cup P_j$ is blocked in $M$ for all 
distinct petals $P_i$ and $P_j$ of $\PP$. This situation
corresponds to a straightforward case in our analysis. 

\begin{lemma}
\label{k*1}
Let $M$ be a connected matroid with an element $x$ such that
$M\ba x$ is $3$-connected up to $M$-clonal
series pairs and 
$\co(M\ba x)\cong M^*(K_{3,n})$. Assume that the
canonical copaddle associated with $M\ba x$ is well-blocked by
$x$.
Then there is
a function $f_{\ref{k*1}}(q)$ such that, if 
$n\geq f_{\ref{k*1}}(q)$, the matroid $M$ has a $U_{q,q+2}$-minor.
\end{lemma} 

Lemma~\ref{k*1} is a consequence of the next lemma.
We prove the dual version as it seems intuitively clearer.

\begin{lemma}
\label{real-doovy}
Let $M$ be a matroid with an element $x$ 
such that $M/x\cong M(K_{3,n})$
for some $n\geq 4$. Assume that $x$ coblocks 
every pair of triads of 
$M(K_{3,n})$. Then there is a function 
$f_{\ref{real-doovy}}(q)$
such that, if $n\geq f_{\ref{real-doovy}}(q)$, 
then $M$ has a $U_{2,q+2}$-minor.
\end{lemma}

\begin{proof}
Let $\PP=(P_1,P_2,\ldots,P_n)$ be the canonical maximal flower 
in  $M(K_{3,n})$ and let $A=\{a_1,a_2,\ldots,a_n\}$ 
be a transversal of the petals of $\PP$. Let $P_i'=P_i-\{a_i\}$.
By assumption
$x$ coblocks $P_i\cup P_j$. 
Therefore $r_M(P_i\cup P_j)=5$ for all distinct 
$i,j\in\{1,2,\ldots,n\}$. Consider
$M/x/A$. Note that $r(M/x/A)=2$ and, for all 
$i\in\{1,2,\ldots,n\}$, we have
$r_{M/x/A}(P'_i)=r_{M/A}(P'_i)=2$.
As $x$ coblocks 
$P_i\cup P_j$, we see that
$x\in\cl_M(P_i\cup P_j)$, so $x\in\cl_{M/A}(P'_i\cup P'_j)$.
Hence, for all distinct $i,j\in\{1,2,\ldots,n\}$, 
we have $r_{M/A}(P'_i\cup P'_j)=3$.
Thus $P'_i$ and $P'_j$ span distinct lines of $M/A$. Hence $M/A$ has
at least $n$ distinct lines. By Lemma~\ref{kung},
a rank-3 matroid with no 
$U_{2,q+2}$-minor
has at most $q^2+q+1$ parallel classes and hence at most
$\binom{q^2+q+1}{2}$ lines. The lemma holds by letting 
$f_{\ref{real-doovy}}(q)=\binom{q^2+q+1}{2}$.
\end{proof}

\begin{proof}[Proof of Lemma~\ref{k*1}]
Let $f_{\ref{k*1}}(q)=f_{\ref{real-doovy}}(q)$
and assume that $n\geq f_{\ref{k*1}}(q)$.
It follows from Lemmas~\ref{k*2x} and \ref{k*2c} that
the flower associated with $\kstar$ in $\co(M\ba x)$
is well blocked by $x$.
Now $M$ has
a $U_{q,q+2}$-minor by the dual of
Lemma~\ref{real-doovy}.
\end{proof}

For a blocked $(n,1)$-copaddle $(M,\PP,J)$ we would not expect to find
elements giving us the easy win of Lemma~\ref{k*1}.
Nonetheless, for $x\in J$,
we can find a partition which effectively gives us
a concatenation of $\PP$ which is well-blocked by $x$.

\begin{lemma}
\label{partition}
Let $M$ be a connected matroid with an element $x$ such that
$M\ba x$ is connected with 
a copaddle $\PP=(P_1,P_2,\ldots,P_t)$. Then there is a partition
$(\calQ,\calP_1,\calP_2,\ldots,\calP_s)$ of the petals of $\PP$,
where possibly $\calQ=\emptyset$
such that the following property holds: if $\calP$ is a set of
petals of $\PP$ such that $x\in\cl(\calP)$, then
$(\{P_1,P_2,\ldots,P_t\}-\calP)\subseteq \calP_i$ for some 
$i\in\{1,2,\ldots,s\}$.
\end{lemma}

\begin{proof}
Assume that $\{\calQ_1,\calQ_2,\ldots,\calQ_s\}$ 
are the distinct minimal subsets of petals whose closure
spans $x$. For $i\in\{1,2,\ldots,s\}$,
let $\calP_i=\{P_1,P_2,\ldots,P_t\}-\calQ_i$,
and let
$\calQ=\calQ_1\cap\calQ_2\udots\calQ_s$.
By Lemma~\ref{modular}, if $i$ and $j$ are distinct elements
of $\{1,2,\ldots,s\}$, then $\calQ_i\cup\calQ_j$
contains all of the petals of $\PP$, that is,
$\calP_i$ and $\calP_j$ are disjoint. It follows from this
fact that
$(\calQ,\calP_1,\calP_2,\ldots,\calP_s)$ satisfies the lemma.
\end{proof}

Note that, in the partition given by Lemma~\ref{partition},
$\calQ$ consists of the petals of $\PP$ that are
blocked by $x$.
Let $(M,\PP,J)$ be a blocked $(n,1)$-copaddle. 
For $x\in J$, let $\pi(x)$
denote the partition of the petals of $\PP$
given by Lemma~\ref{partition}.

\begin{lemma}
\label{partition1}
Let $(M,\PP,J)$ be a blocked $(n,1)$-copaddle where $M\in\eq$, and let 
$x$ be an element of $J$ with 
$\pi(x)=(\calQ,\calP_1,\calP_2,\ldots,\calP_k)$.
Then $|\calQ|+k\leq f_{\ref{k*1}}(q)$.
\end{lemma}

\begin{proof}
Let $\{Q_1,Q_2,\ldots,Q_k\}$ be a transversal of 
$(\calP_1,\calP_2,\ldots,\calP_k)$, and let
$Q$ be the union of the set of petals of $\PP$
not in $\calQ\cup \{Q_1,Q_2,\ldots,Q_k\}$.
Then $\calQ\cup \{Q_1,Q_2,\ldots,Q_k\}$ is the set
of petals of a copaddle in $M\ba J/Q$.
It is easily checked that his flower is well-blocked
by $x$.
The lemma now follows from Lemma~\ref{k*1}.
\end{proof}

As immediate consequences of Lemma~\ref{partition1}
we get the next two results.

\begin{corollary}
\label{partition-cor}
Let $(M,\PP,J)$ be a blocked $(n,1)$-copaddle where 
$M\in \eq$.  If $x\in J$, then $x$ blocks at
most $f_{\ref{k*1}}(q)$ petals of $\PP$.
\end{corollary}

\begin{corollary}
\label{partition2}
Let $(M,\PP,J)$ be a blocked $(n,1)$-copaddle where $M\in\eq$,
and let $\calP$ be a subset of $\PP$ with
$|\calP|\geq 2f_{\ref{k*1}}(q)$. If $x\in J$,
then there is a subset $\calQ$ of $\calP$
such that:
\begin{itemize}
\item[(i)] $|\calQ|\geq|\calP|/f_{\ref{k*1}}(q)$, and
\item[(ii)] $x$ is spanned by the union of the petals
of $\calP$ that are not in $\calQ$.
\end{itemize}
\end{corollary}

\subsection*{Good $3$-paths} Let $(M,\PP,J)$ be a blocked
$(n,1)$-copaddle where $\PP=(P_1,P_2,\ldots,P_t)$.
A {\em good $3$-path of length $s$}
\index{good $3$-path} 
in $(M,\PP,J)$ is a 
partition $(\calP_1,\calP_2,\ldots,\calP_s)$
of $\{P_{n+1},P_{n+2},\ldots,P_t\}$, and a subset
$\{x_1,x_2,\ldots,x_t\}$ of $J$ such that 
\begin{itemize}
\item[(i)] $x_i\in\cl(\calP_i)$ 
for all $i\in\{1,2,\ldots,s\}$, and 
\item[(ii)] $x_i$ blocks a series pair of a member of
$\calP_i$ for all
$i\in\{1,2,\ldots,s\}$.
\end{itemize}

It is not hard to see that having a sufficiently long
good 3-path is a clear win.

\begin{lemma}
\label{copaddle-win}
Let $(M,\PP,J)$ be a blocked $(n,1)$-copaddle.
Then there is a function $f_{\ref{copaddle-win}}(q)$
such that, if $M$ has a good $3$-path of length $s$, where
$s>f_{\ref{copaddle-win}}(q)$, then $M\notin\eq$.
\end{lemma}

\begin{proof}
Let $f_{\ref{copaddle-win}}(q)=f_{\ref{at-last}}(3,q)$.
Assume that $s>f_{\ref{copaddle-win}}(q)$.
Assume that the partition $(\calP_1,\calP_2,\ldots,\calP_s)$
of $\{P_{n+1},P_{n+2},\ldots,P_t\}$, and the subset
$\{x_1,x_2,\ldots,x_t\}$ of $J$ define a good 3-path
in $(M,\PP,J)$.
Let $Q_1$ be the union of the sets in 
$\{P_1,P_2,\ldots,P_n,\{x_1\}\}\cup \calP_1$, and for 
$i\in\{2,3,\ldots,s\}$, let $Q_i=\cup_{P\in\calP_i}(P)\cup\{x_i\}$.
Then 
$\QQ=(Q_1,Q_2,\ldots,Q_s)$
is a path of 3-separations in $M\ba (J-\{x_1,x_2,\ldots,x_s\})$,
each step of which contains a $\QQ$-strong clonal pair.
Assume that $\QQ$ displays a swirl-like flower
$(R_1,R_2,R_3)$. Let $E'=E(\co(M\ba J))$.
By Lemma~\ref{keep-flower-type}
$(R_1\cap E',R_2\cap E',R_3\cap E')$ is a swirl-like flower
in $\co(M\ba J)$. But, also by Lemma~\ref{keep-flower-type}
$(R_1\cap E',R_2\cap E',R_3\cap E')$ is a copaddle in
$\co(M\ba J)$. Thus $\QQ$ does not display a 3-petal
swirl-like flower. By Lemma~\ref{at-last}, 
$M\notin\eq$. 
\end{proof}

\subsection*{Near $n$-blockages} Let $(M,\PP,J)$
be a blocked $(n,1)$-copaddle. For $i\in\{1,2,\ldots,n\}$
there is an element of $J$ that 
blocks the series pair $\{a_i,a'_i\}$ in
$P_i$. Of course there may be more than one such element. Let 
$\phi:\{1,2,\ldots,n\}\rightarrow J$ be a function with 
the property that $\phi(i)$ blocks $\{a_i,a'_i\}$.
If $(M,\PP,J)$ is endowed with such a function we will
say that $(M,\PP,J)$ is a blocked $(n,1)$-copaddle with
{\em blocking function} $\phi$.
Evidently deleting the members of $J$ that are not in the range
of $\phi$ preserves the property of being a blocked
$(n,1)$-copaddle.

Let $(M,\PP,J)$ be a blocked $(n,1)$-copaddle with 
blocking function $\phi$ where $\PP=(P_1,P_2,\ldots,P_t)$.
Then $(M,\PP,J)$ is a {\em near} $n$-blockage if
\index{near $n$-blockage}
$P_{n+1}\cup P_{n+2}\udots P_t$ blocks $\phi(i)$ for all
$i\in\{1,2,\ldots,n\}$.

The next task is to find a near $n$-blockage. Note that,
if $(M,\PP,J)$ is a blocked $(n,1)$-copaddle of order $n$, then
it is also a blocked $(n,1)$-copaddle of any order
$m\leq n$.

\begin{lemma}
\label{almost-blockage}
Let $(M,\PP,J)$ be a blocked $(n,1)$-copaddle with
blocking function $\phi$, where $\PP=(P_1,P_2,\ldots,P_t)$.
Assume that $(M,\PP,J)$ has a good $3$-path of length $s$ and
that $(M,\PP,J)$ is not a near $n$-blockage. If 
$n\geq mf_{\ref{k*1}}(q)$, then there is an ordering of the
elements of $\{P_1,P_2,\ldots,P_n\}$ such that, relative to
this ordering, $(M,\PP,J)$ is a blocked $(m,1)$-copaddle
with a good $3$-path of length $s+1$.
\end{lemma}

\begin{proof}
Assume that $n\geq mf_{\ref{k*1}}(q)$. Assume that
the partition $(\calP_1,\calP_2,\ldots,\calP_s)$ of 
$\{P_{n+1},\ldots,P_t\}$ together
with the subset $\{x_1,x_2,\ldots,x_s\}$ of $J$ give a 
good path of length $s$.
As $(M,\PP,J)$ is not a near $n$-blockage, there is
an $i\in\{1,2,\ldots,t\}$ such that $\phi(i)$ is not
blocked by $P_{n+1}\cup P_{n+2}\udots P_t$. By 
Corollary~\ref{partition2}, there is a subset 
$\calP$ of $\{P_1,P_2,\ldots,P_n\}$ of size
$m$ such that $x$ is in 
the closure of the union of the petals of $\PP$
that are not in $\calP$. Up to
labels we may assume that this subset is $\{P_1,P_2,\ldots,P_m\}$.
As $\phi(i)$ blocks $\{a_i,a'_i\}$, we see that
$i\notin \{1,2,\ldots,m\}$. Let 
$\calP_{s+1} =\{P_{m+1},P_{m+2},\ldots,P_n\}$ and let 
$x_{s+1}=\phi(i)$. It is now clear that
$(\calP_1,\calP_2,\ldots,\calP_{s+1})$ and 
$(x_1,x_2,\ldots,x_{s+1})$ define a good path of length
$s+1$ in the blocked $(m,1)$-copaddle $(M,\PP,J)$.
\end{proof}

\begin{corollary}
\label{get-almost}
Let $(M,\PP,J)$ be a blocked $(n,1)$-copaddle
with blocking function $\phi$ where  
$M\in\eq$. Then there is a function $f_{\ref{get-almost}}(m,q)$
such that, if $n\geq f_{\ref{get-almost}}(m,q)$, then
there is an ordering of
the petals of $\PP$ such that, with respect to this ordering,
$(M,\PP,J)$ is a near $m$-blockage whose blocking function
is $\phi$ restricted to
$\{1,2,\ldots,m\}$.
\end{corollary}

\begin{proof}
Let $f_{\ref{get-almost}}(m,q)=
mf_{\ref{k*1}}(q)^{f_{\ref{copaddle-win}}(q)}$.
With this function, the corollary follows immediately from
Lemmas~\ref{copaddle-win} and \ref{almost-blockage}.
\end{proof}

\subsection*{Cleaning a Near $n$-blockage}
At last we are able to achieve the goal of this section.

\begin{lemma}
\label{get-blockage}
Let $(M,\PP,J)$ be a near $n$-blockage with blocking function
$\phi$, where $\PP=(P_1,P_2,\ldots,P_t)$. Assume that
$M\in\eq$.
Then there is a function $f_{\ref{get-blockage}}(m,q)$
such that if $n\geq f_{\ref{get-blockage}}(m,q)$,
then there is an ordering of the petals
of $\PP$ 
such that, with respect to this ordering
$(M,\PP,J)$ is an $m$-blockage.
\end{lemma}

\begin{proof}
Let $f_{\ref{get-blockage}}(m,q)=
f_{\ref{easy-hyp2}}(f_{\ref{k*1}}(q),m)$.

Define a hypergraph whose vertices are $\{P_1,P_2,\ldots,P_n\}$
and whose edges are $\{\phi(1),\phi(2),\ldots,\phi(n)\}$
as follows: for $i\in\{1,2,\ldots,n\}$, the vertices
incident with $\phi(i)$ are the members of
$\{P_1,P_2,\ldots,P_n\}$ that $\phi$ blocks. It follows from
Corollary~\ref{partition-cor} that each edge is incident with at most 
$f_{\ref{k*1}}(q)$ vertices. It follows from the definition
of $\phi$ that
$P_i$ is incident with $\phi(i)$
for all $i\in\{1,2,\ldots,n\}$.
It now follows from Lemma~\ref{easy-hyp2} that,
up to an appropriate permutation of $\{P_1,P_2,\ldots,P_n\}$,
we may assume, for all $i,j\in\{1,2,\ldots,m\}$,
that $\phi(i)$ blocks $P_j$ if and only if $i=j$. The
lemma now follows by observing that
the members of $\{\phi(i),\phi(2),\ldots,\phi(m)\}$ have the
properties required for an $m$-blockage.
\end{proof}

\begin{proof}[Proof of Theorem~\ref{star3}]
Let $f_{\ref{star3}}(m,q)=
f_{\ref{get-almost}}(f_{\ref{get-blockage}}(m,q),q)$.
Let $\phi$ be a blocking function for $(M,\PP,J)$.
By Corollary~\ref{get-almost}, there is an ordering of 
$\PP$ such that, with respect to this ordering,
$(M,\PP,J)$ is a near 
$f_{\ref{get-almost}}(m,q)$-blockage whose blocking
function is $\phi$ restricted to 
$\{1,2,\ldots,m\}$. Life is now good as we can observe that,
by Lemma~\ref{get-blockage}, there is an ordering of the
petals of $\PP$, relative to which, $(M,\PP,J)$ is
an $m$-blockage.
\end{proof}

\section{Building a Blocked $(n,l)$-copaddle}

In this section we prove that given a very large $3$-connected
clonal matroid in $\eq$ then we can either find a large free swirl
minor or we can build blocked $(n,l)$-copaddles, where $l$
is large.

\begin{theorem}
\label{star-block}
Let $M$ be a $3$-connected clonal matroid in 
$\eq$. Then there is a function $f_{\ref{star-block}}(s,l,m,q)$
such that, if $|E(N)|\geq f_{\ref{star-block}}(s,l,m,q)$,
then at least one of the following holds.
\begin{itemize}
\item[(i)] $M$ has a clonal $\Delta_m$-minor.
\item[(ii)] Either $M$ or $M^*$
has a minor $N$ with an associated triple
$(N,\PP,J)$ such that $(N,\PP,J)$ is a blocked $(s,l)$-copaddle
where all of the series pairs of $N\ba J$ are $M$-clonal.
\end{itemize}
\end{theorem}

Of course we will need some preliminary lemmas. 
The next technical lemma is straightforward but crucial.

\begin{lemma}
\label{crucial}
Let $M$ be a connected matroid with an element $x$ such that
$M\ba x$ is $3$-connected up to stars and let $S$ be a maximal
star of $M\ba x$. Assume that $x\in\cl(S)\cap\cl(E(M)-S)$ and
that $x$ is cofixed in $M\ba S$. Let $M'$ be a matroid
obtained by cocloning $x$ by $x'$.
If $\{x,x'\}$ is not a series pair of $M'$ then there is a
partition $(S_1,S_2)$ of $S$ such that the following hold.
\begin{itemize}
\item[(i)] $\{x,x'\}$ is a parallel pair in $M'\ba S_1/S_2$.
\item[(ii)] $M'\ba S_1/S_2\ba x'=M\ba S$.
\end{itemize}
\end{lemma}

\begin{proof}
Let $T=E(M)-(S\cup\{x\})$. We first observe,

\begin{sublemma}
\label{crucial1}
$T$ is not coblocked by $x'$.
\end{sublemma}

\begin{proof}
Otherwise $x'\in \cl_{M'}(T)$ so that $\{x,x'\}\subseteq \cl_{M'}(T)$.
This means that $\{x,x'\}$ is a 
coindependent clonal pair in $M'|(T\cup\{x,x'\})$. 
But $M'|(T\cup\{x,x'\})/x'\cong M\ba S$, and we have
contradicted
the fact that $x$ is cofixed in $M\ba S$.
\end{proof}

Note that $\{x,x'\}$ is independent in $M'$, as otherwise
$x$ is a loop of $M=M'/x$.
Let $S'\subseteq S$ be a transversal of the series pairs of $S$ and
let $S''=S-S'$.
Each series pair in $S$ is also a series pair in $M'$ and 
it follows that $\{x,x'\}$ is independent in $M'/S'$.
As $S$ is a star of $M$, we have
$r_{M'/S'/x'}(S''\cup \{x\})=1$. Thus $r_{M'/S'}(S''\cup \{x,x'\})=2$.
If $S''$ consists of a single parallel set in the closure of
$T$, then $\{x,x'\}$ is a series pair in $M'/S'$ and hence,
also in $M'$. Thus there is an element $s\in S''$ such that 
$s\notin\cl_{M'/S'}(T)$. As $\{x,x'\}$ is a clonal pair,
$\{s,x,x'\}$ is a triangle. Let 
$S_1=S''-\{s\}$ and $S_2=S'\cup \{s\}$. 
It is now clear that the lemma is satisfied 
by this choice of $S_1$ and $S_2$.
\end{proof}

Let $(M,\PP,J)$ be an $n$-blockage with blocking function
$\phi$, where $\PP=(P_1,P_2,\ldots,P_t)$.
Then $(M,\PP,J)$ is a {\em minimal} $n$-blockage
\index{minimal $n$-blockage}
if $J=\{\phi(i):i\in\{1,2,\ldots,n\}\}$, and $M$ is $3$-connected,
that is, if there are no unblocked series pairs in 
$P_{n+1}\cup P_{n+2}\udots P_t$. Clearly there is no loss
of generality in focussing on minimal $n$-blockages.
Note that, for any $n$-blockage, $\phi$ is injective, so that
for a minimal $n$-blockage $\phi$ is a bijection.

In what follows, in the $n$-blockage $(M,\PP,J)$,
we denote the series pair in the petal $P_i$ of $\PP$ 
by $\{a_i,a'_i\}$ for $i\in\{1,2,\ldots,n\}$. 
We denote the other 
two elements of $P_i$ by $\{b_i,c_i\}$. Let
$B=\{b_1,b_2,\ldots,b_n\}$ and $C=\{c_1,c_2,\ldots,c_n\}$.

\begin{lemma}
\label{find-4fans}
Let $(M,\PP,J)$ be a minimal $n$-blockage with blocking function
$\phi$. Then the following hold.
\begin{itemize}
\item[(i)] $M/C$ is $3$-connected.
\item[(ii)] If $i\in\{1,2,\ldots,n\}$, then
$(b_i,a_i,a'_i,\phi(i))$ is a maximal fan in $M/C$ where
$\{b_i,a_i,a'_i\}$ is a triangle and $\{a_i,a'_i,\phi(i)\}$
is a triad.
\end{itemize}
\end{lemma}

\begin{proof}
Consider part (i). Evidently $M/C$ is connected.
Assume that $M/C$ is not 3-connected.
Let $(X',Y')$ be a 2-separation of $M/C$. 
Let $(X,Y)=(X'-J,Y'-J)$. Assume that $|X|\leq 1$.
As $J$ is coindependent in $M$, and hence in $M/C$, we 
see that $r_{M/C}(X)\in\{1,2\}$. In either case 
we deduce that there is an element $z\in E(M)$
and an $\{i,j\}\subseteq \{1,2,\ldots,n\}$ such that
$C\cup\{\phi(i),\phi(j),z\}$ contains a circuit. A routine 
check shows that all possible cases 
contradiction the definition of minimal $n$-blockage.
Therefore $(X,Y)$ is a 2-separation in $M/C\ba J$.
It is straightforwardly verified that, up to labels,
we may assume that either 
\begin{itemize}
\item[(a)] $X=\{a_1,a'_1\}$ or
\item[(b)] $X=\{\{a_1,a'_1,b_1\},\{a_2,a'_2,b_2\},
\ldots,\{a_i,a'_i,b_i\}\}$ for some $i\in\{1,2,\ldots,n\}$.
\end{itemize}
The verification of this 
is particularly routine if one considers the dual.
Recall that $\phi(1)$ blocks both $\{a_1,a'_1\}$
and $P_{n+1}\cup P_{n+2}\udots P_t$ in $M$.
Thus, in either of the above cases, $\phi(1)$ blocks
$X$ in $M/C$. Therefore $M/C$ is $3$-connected
so that (i) holds.

Consider (ii). Say $i\in\{1,2,\ldots,n\}$. 
Evidently $\{a_i,a'_i,\phi(i)\}$ is a triad in 
$M/C$. As $\{a_i,a'_i,b_i,c_i\}$ is a circuit in 
$M$ and $M/C$ is $3$-connected we see that $\{a_i,a'_i,b_i\}$
is a triangle in $M/C$. Thus (ii) holds.
\end{proof}

Let $C$ and $D$ be disjoint subsets of a matroid and
let $N=M\ba D/C$. Say $d\in D$ and $c\in C$. 
In the remainder of this section we will at times refer
to the matroid $M\ba(D-\{d\})/C$ as being obtained from
$N$ by {\em undeleting}
\index{undeletion} 
$d$ and the matroid
$M\ba D/(C-\{c\})$ as being obtained from $M$ by {\em uncontracting}
\index{uncontraction}
$c$.

\begin{lemma}
\label{starry-night}
Let $M$ be a $3$-connected clonal matroid in $\eq$
and let $N$ be a $3$-connected minor of $M$ with the following
properties.
\begin{itemize}
\item[(i)] There is a set $B=\{b_1,b_2,\ldots,b_n\}$
of elements of $N$ such that $N\ba B$ is an 
$(n,l)$-copaddle, with associated flower 
$(P_1,P_2,\ldots,P_t)$.
\item[(ii)] The series pairs of $N\ba B$ are $M$-clonal.
\item[(iii)] For $i\in\{1,2,\ldots,n\}$, the maximal
star $S_i\in P_i$ of $N\ba B$ has the property that
$b_i\in\cl_N(S_i)\cap\cl_N(E(N)-S_i)$.
\end{itemize}
Then there is a function $f_{\ref{starry-night}}(m,q)$
such that, if $n\geq f_{\ref{starry-night}}(m,q)$,
then $M$ has an $(m,l+1)$-copaddle minor all of whose series
pairs are $M$-clonal.
\end{lemma}

\begin{proof}
Let $f_{\ref{starry-night}}(m,q)=f_{\ref{k3n-case}}(q)+m$.
Assume that $n\geq f_{\ref{starry-night}}(m,q)$.
As $M$ is clonal, each member of $\{b_1,b_2,\ldots,b_n\}$
has a clone in $M$. For $i\in\{1,2,\ldots,n\}$, denote the
clone of $b_i$ by $b'_i$. Let $C$ and $D$ be disjoint subsets
of $E(M)$ such that $N=M\ba D/C$.

Consider $P_1$. There are two cases to consider. For the first
assume that $b_1'\in D$. Let $N'$ denote the matroid obtained
by undeleting $b'_1$. As $b_1$ is fixed in $N$, we see that
$\{b_1,b'_1\}$ is a parallel pair in $N'$. In this case
let $P'_1=(P_1-S_1)\cup\{b'_1\}$ and let $N_1=N'\ba S_1$.
Note that we have replaced the star $S_1$ of $N\ba B$
by an $M$-clonal parallel pair at the base point of $S_1$.

For the second case, assume that $b'_1\in C$. In this case
let $N'$ be the matroid obtained from $N$ by uncontracting 
$b'_1$. Assume that $\{b_1,b'_1\}$ is not a series pair 
in $N'$. Note that
$b_1$ is cofixed in $N\ba S_1$ so that Lemma~\ref{crucial}
applies. By this lemma there is a partition
$(X,Y)$ of $S_1$ such that $\{b_1,b'_1\}$ is a parallel
pair in $N'\ba X/Y$ and $N'\ba X/Y\ba b'_1=N$.
In this case, let $N_1=N'\ba X/Y$ and let
$P'_1=(P_1-S_1)\cup\{b'_1\}$. Observe that we have
again replaced the star $S_1$ by an $M$-clonal parallel pair 
at the basepoint of the star. On the other hand, if
$\{b_1,b'_1\}$ is a series pair, then let $N_1=N'$ and
let $P'_1=P_1\cup\{b'_1\}$.
Observe that $S_1\cup\{s_1,s'_1\}$ is a star of order
$l+1$ in $P'_1$.

Repeat the above process for the remaining elements of 
$B$ to obtain a matroid $N_n$ with associated flower
$\PP'=(P'_1,P'_2,\ldots,P'_t)$, where, for $i\in\{n+1,n+2,\ldots,t\}$
we have $P'_i=P_i$.
Assume that for some $s\geq f_{\ref{k3n-case}}(q)$
the are $s$ members of $\{P'_1,P'_2,\ldots,P'_n\}$
that contain $M$-clonal parallel pairs.
Up to labels we may assume that these are 
$\{P'_1,P'_2,\ldots,P'_s\}$. Observe that
$(N_n/(P'_{s+1}\cup P'_{s+2}\udots P'_t))^*$ satisfies the
hypotheses of Lemma~\ref{k3n-case}, giving the contradiction that
$N_n\notin \eq$. Otherwise, by the definition of
$f_{\ref{starry-night}}(m,q)$, at least $m$ petals in 
$\{P'_1,P'_2,\ldots,P'_n\}$ have stars of order $l+1$.
The series pairs in these stars are $M$-clonal.
The lemma follows from these observations.
\end{proof}

\begin{lemma}
\label{star71}
Let $M$ be a $3$-connected matroid in $\eq$ with a minor
$N$ that is an $(n,l)$-copaddle. Assume that all series pairs
in $N$ are $M$-clonal. Then there is a function 
$f_{\ref{star71}}(m,q)$ such that,
if $n\geq f_{\ref{star71}}(m,q)$, then $M$ has a minor
$M'$ that is a blocked $(n,l)$-copaddle $(M',\PP',J)$. Moreover,
all of the series pairs of $\PP'$ are $M$-clonal.
\end{lemma}

\begin{proof}
Let $f_{\ref{star71}}(m,q)=f_{\ref{k3n-case}}(q)+m$.
Assume that $n\geq f_{\ref{star71}}(m,q)$.
Let $\PP=(P_1,P_2,\ldots,P_t)$ be the flower associated with
$N$. For $i\in\{1,2,\ldots,n\}$, let $S_i$ denote the maximal
star in $P_i$. By Lemma~\ref{bridge-klonal2}, if $\{s_i,s'_i\}$
is a series pair in $S_i$, then $\{s_i,s'_i\}$
has a 1- or 2-element bridging
sequence. We will say that $P_i$ has {\em type-$1$} if 
every series pair in $S_i$ has
a 1-element bridging sequence and has {\em type}-2 otherwise. 

Note that the first element of a bridging sequence for a series
pair consists of a delete element. It follows 
routinely that
if there are $m$ members of $\{P_1,P_2,\ldots,P_n\}$ that
are of type-1, then the lemma holds. Otherwise, we may assume
that for some $\mu\geq f_{\ref{k3n-case}}(q)$, there are
$\mu$ members of $P_i$ that are of type-2. Up to 
labels we may assume
that $\{P_1,P_2,\ldots,P_\mu\}$ are all of type-2.
For $i\in\{1,2,\ldots,\mu\}$, let 
$\{s_i,s'_i\}$ be a series pair in $S_i$ with a 2-element
bridging sequence $(t_i,u_i)$. Let $T=\{t_1,t_2,\ldots,t_\mu\}$
and let $N'$ be the matroid obtained from
$N$ by undeleting the elements of $T$. Say $t_i\in T$.
As $(t_i,u_i)$ is a bridging sequence and 
$t_i$ is a delete element of this bridging sequence,
$t_i\in\cl_{N'}(\{s_i,s'_i\})$.  Note that
this means that $\{s_i,s'_i\}$ is a series pair
in $N'\ba t_i$ and that
$(P_1\cup\{t_1\},P_2\cup\{t_2\},\ldots,P_n\cup\{t_n\},P_{n+1},
\ldots,P_t)$ is a copaddle in $N'$. As $\{s_i,s'_i\}$
is $M$-clonal, $\{t_i,s_i,s'_i\}$ is a triangle in 
$N'$. If $\{s_i,s'_i\}$ is a series pair in $N'$,
then $u_i$ cannot coblock $\{s_i,s'_i\}$.
Hence $\{s_i,s'_i\}$ is not a series pair in $N'$
so that $\{t_i,s_i,s'_i\}$ is also a triad in
$N'$.

For $i\in\{1,2,\ldots,n\}$ let $P'_i=P_i-(S_i-\{s_i,s'_i\})$.
Let $N''=N'/T\ba((S_1-\{s_1,s'_1\})\udots(S_\mu-\{s_\mu,s'_\mu\}))
/(P_{n+1}\cup P_{n+2}\udots P_t)$.
Then $\PP'=(P'_1,P'_2,\ldots,P'_n)$ is a copaddle
in $N''$. 
Moreover, each petal of $\PP'$ contains a single 
$M$-clonal parallel pair and $\si(N'')\cong M^*(K_{3,n})$.
Thus $\PP'$ satisfies the hypotheses of Lemma~\ref{k3n-case}
in $M^*$ and, by that lemma we obtain the contradiction that
$M\notin\eq$.
\end{proof}

\begin{proof}[Proof of Theorem~\ref{star-block}]
The proof is by induction on $l$.
If $l>1$, assume that $f_{\ref{star-block}}(s,l-1,m,q)$ has been
defined and, for induction, assume that the theorem holds 
with this definition of $f_{\ref{star-block}}(s,l-1,m,q)$.
Let $n_4=f_{\ref{star71}}(s,q)$. Let
$n_3=f_{\ref{starry-night}}(n_4,q)$. Let
$n_2=n_3+f_{\ref{k3n-case}}(q)$. Let
$n_1=f_{\ref{star1}}(n_2,m,q)$ if $l=1$
and let $n_1=f_{\ref{star-block}}(n_2,l-1,m,q)$
if $n>1$. 
Finally 
let $f_{\ref{star-block}}(s,l,m,q)=f_{\ref{star-col}}(n_1,m,q)$.

Assume that $n\geq f_{\ref{star-block}}(s,l,m,q)$.
Assume that $M$ does not have a clonal $\Delta_m$-minor.
By Corollary~\ref{star-col}, 
$M$ or $M^*$ has a minor $M_1$ with an associated triple
$(M_1,\PP,J)$ such that 
$(M_1,\PP,J)$ is an $n_1$-blockage. 
We now apply Lemma~\ref{find-4fans}.
Using the notation of that lemma we have a set $C$ such that
$M_1/C$ is $3$-connected, and every set in
$\{\{a_1,a'_1\},\{a_2,a'_2\},\ldots,\{a_{n_1},a'_{n_1}\}\}$
is contained in both a triangle and a triad.
By Lemma~\ref{3-klonal} $M_1/C$ has a 
3-connected minor $N$ whose ground set is
$\{a_1,a'_1\}\cup \{a_2,a'_2\}\udots \{a_{n_1},a'_{n_1}\}$. 

If $l=1$ we apply Lemma~\ref{star1} and if $l>1$ we apply
the induction assumption to deduce that either $M_1/C$
or $(M_1/C)^*$ has a minor $M_2$ with an associated triple
$(M_2,\PP_2,J_2)$ such that $(M_2,\PP_2,J_2)$
is a blocked $(n_2,l)$-copaddle. Moreover,
all of the series pairs of $M_2\ba J_2$ are $M_1$-clonal
and hence $M$-clonal. Let $M_1'=M_1/C$ if $M_2$ is a minor of 
$M_1$ and otherwise let $M_1'=(M_1/C)^*$. Every series pair
of $M_2\ba J_2$ is in a 4-element fan in $M'_1$. Thus every
series pair in $M_2\ba J_2$ is in a triangle of $M_1'$.
By a slight abuse of notation assume that
$\PP_2=(P_1,P_2,\ldots,P_t)$. For $i\in\{1,2,\ldots,n_2\}$,
let $S_i$ denote the maximal star
in $P_i$ and let $\{s_i,s'_i\}$ be a series pair in $S_i$.
Let $B'=\{b_1,b_2,\ldots,b_{n_2}\}$ be elements
of $E(M'_1)-E(M_2)$ such that
$\{s_i,s'_i,b_i\}$ is a triangle in $M'_1$
for $i\in\{1,2,\ldots,n_2\}$.

Observe that the members of $B'$ were all deleted from
$M_1'$ as otherwise, for some $i\in\{1,2,\ldots,n_2\}$, 
we contradict the fact that $\{s_i,s'_i\}$
is independent in $M_2$.
Let $M_3$ be the matroid obtained by undeleting the elements
of $B'$ from $M_2$. Observe that $\{s_i,s'_i,b_i\}$
is a triangle in $M_3$ for all $i\in\{1,2,\ldots,n_2\}$.
Applying Lemma~\ref{k3n-case} and arguing just as in
the proof of Lemma~\ref{star71}, we deduce, up to the labels of
petals in $\{P_1,P_2,\ldots,P_{n_2}\}$, that we may assume that
$\{s_i,s'_i\}$ is a series pair in $M_3$ for all
$i\in\{1,2,\ldots,n_3\}$. It follows from this that,
for $i\in\{1,2,\ldots,n_3\}$, the element
$b_i$ is in the guts of a 2-separation of $M_3$
one side of which is $S_i$. Let $B=\{b_1,b_2,\ldots,b_{n_3}\}$
and let $M_4$ be the matroid obtained by undeleting the elements
of $B$ from $M_2$.

We may now apply Lemma~\ref{starry-night} and deduce that
$M$ or $M^*$ has an $(n_4,l+1)$-copaddle minor all of whose
series pairs are $M$-clonal. Finally, by Lemma~\ref{star71},
$M$ or $M^*$ has a minor $M'$ with an associated triple
$(M',\PP',J')$ such that $(M',\PP',J')$
is a blocked $(s,l+1)$-copaddle
such that all series pair of $M'\ba J'$ are $M$-clonal,
completing the proof of the theorem.
\end{proof}

\section{Excluding a Blocked $(n,l)$-copaddle}

We now show that a matroid in $\eq$ cannot contain a blocked 
$(n,l)$-copaddle for large values of $l$. In striking contrast
to Theorem~\ref{star-block}, the bound here is quite modest.
We will use the following theorem of  Lemos and Oxley
\cite{oxlem01}.

\begin{theorem}
\label{oxley-lemos}
Let $M$ be a connected matroid in which a largest circuit and 
cocircuit have $c$ and $c^*$ elements respectively. Then
$M$ has at most $\frac{1}{2}cc^*$ elements.
\end{theorem}

\begin{lemma}
\label{kill-star}
Let $M$ be a $3$-connected matroid with a set $S$ such that 
$M|S\cong M(K_{2,l})$. Assume that all of the series pairs
in $M|S$ are $M$-clonal. Then there is a function 
$f_{\ref{kill-star}}(q)$ such that, if 
$l\geq f_{\ref{kill-star}}(q)$, then $M$ has a 
$\Lambda_q$- or $U_{2,q+2}$-minor.
\end{lemma}

\begin{proof}
Let $f_{\ref{kill-star}}(q)=\left(\frac{1}{2}q^2+1\right)\lfloor
\frac{q+2}{2}\rfloor$. Assume that
$l\geq f_{\ref{kill-star}}(q)$.

Say that the series pairs in $M|S$ are
$\mathcal S=\{\{s_1,s'_1\},\{s_2,s'_2\},\ldots,\{s_l,s'_l\}\}$.
Then, by Lemma~\ref{3-klonal}, $M$ has a $3$-connected minor 
$M_1$ on $S$. It follows from Lemma~\ref{pi-minor} that
$\sqcap_{M_1}(\{s_i,s'_i\},\{s_j,s'_j\})\geq 1$ for all distinct
$i$ and $j$ in $\{1,2,\ldots,l\}$. Consider the lines of 
$M_1$ that are spanned by members of $\mathcal S$. If any
of these lines has at least $q+2$ points, then the lemma holds,
so we may assume that each line contains at most
$\lfloor\frac{q+2}{2}\rfloor$ members of
$\mathcal S$. Thus there is an integer $n_1$, where
$n_1+1\geq  \frac{1}{2}q^2+1$ such that, up to labels,
$\sqcap_{M_1}(\{s_i,s'_i\},\{s_j,s'_j\})=1$ for all
$i\in\{1,2,\ldots,n_1,n_1+1\}$. 
Let 
$M_2=M_1|(\{s_1,s'_1\},\{s_2,s'_2\},\ldots,\{s_{n_1+1},s'_{n_1+1}\})$.
Evidently $M_2$ is 3-connected. 

\begin{sublemma}
\label{kill-star1}
There is an $i\in\{1,2,\ldots,n_1+1\}$ such that
$M_2/s_i$ is $3$-connected.
\end{sublemma}

\subproof
Assume that the sublemma fails. 
Then, for all  $i\in\{1,2,\ldots,n_1+1\}$ there is a
path $(X_i,\{s_i,s'_i\},Y_i)$ of $3$-separations in $M_2$,
where $\{s_i,s'_i\}\subseteq \cl(X_i)$ and, of course,
$\{s_i,s'_i\}\subseteq \cl(Y_i)$.
Amongst all such paths assume that we have chosen $s_i$
and the associated path so that $|X_i|$ is minimal.
Note that we may assume that $X_i$ is a union of clonal pairs.
If $|X_i|=2$, then $\sqcap_{M_1}(X_i,\{s_i,s'_i\})=2$ and
$X_i$ is a clonal pair so that we have contradicted the fact
that $\sqcap_{M_1}(\{s_i,s'_i\},\{s_j,s'_j\})=1$ for all
$i\in\{1,2,\ldots,n_1,n_1+1\}$. Thus $|X_i|\geq 4$.

Choose $\{s_j,s'_j\}\subseteq X_i$ and let 
$(X_j,\{s_j,s'_j\},Y_j)$ be an associated path of 
3-separations where $\{s_i,s'_i\}\subseteq Y_j$.
Assume that $X_j\cap Y_i\neq \emptyset$.
An easy uncrossing argument shows that
$\lambda(X_j\cap Y_i)=\lambda((X_j\cap Y_i)\cup \{s_j\})=2$.
Thus $s_j$ is in the guts or coguts of
$(X_j\cap Y_i,E(M)-(X_j\cap Y_i))$. But $s_j\in\cl(Y_j)$
and $Y_j\subseteq E(M)-(X_j\cap Y_i)$, so that
$s_j$, and hence also $s'_j$ is in $\cl(X_j\cap Y_i)$.
Therefore $\{s_i,s'_i,s_j,s'_j\}\subseteq \cl(Y_i)$
and $M|\{s_i,s'_i,s_j,s'_j\}\cong U_{2,4}$, contradicting the
fact that $\sqcap_{M_1}(\{s_i,s'_i\},\{s_j,s'_j\})=1$.
We conclude that $X_j\cap Y_i$ is empty, so that
$X_j$ is a subset of $X_i$. As $\{s_j,s'_j\}\subseteq X_i-X_j$,
we have contradicted the minimality assumption and the sublemma
follows.
\end{proof}

By \ref{kill-star1} we may assume that $M/s_{n_1+1}$
is $3$-connected. Relabel $s'_{n_1+1}$ by $t$
and let $M_3=M_2/s_{n_1+1}$. Observe that
$\{s_i,s'_i,t\}$ is a triangle in $M_3$ for all 
$i\in\{1,2,\ldots,n_1\}$. Consider $M_3/t$.
This matroid is connected, so $\si(M_3/t)$ is connected.
By Theorem~\ref{oxley-lemos} $\si(M_3/t)$ either
has a circuit of size at least $q$ or a cocircuit 
of size at least $q$. Assume that $C$ is such a set.
Up to labels we may assume that
$C=\{s_1,s_2,\ldots,s_\mu\}$, where $\mu\geq q$.

Assume that $C$ is a circuit of $\si(M_3/t)$ of size at 
least $q$. 
As $\{s_1,s'_1\}$ is a clonal
pair in $M_3$, we have $s'_1\in\cl_{M_3}(C)$ so that
$t\in\cl_{M_3}(C)$. Thus $C\cup \{t\}$ is a circuit in $M_3$.
It now follows that 
$M_3|(\{s_1,s'_1\}\cup\{s_2,s'_2\}\udots \{s_\mu,s'_\mu\})
\cong \Lambda_\mu$. Hence $M$ has a $\Lambda_q$-minor.

Assume that $C$ is a cocircuit
of $\si(M_3/t)$ of size at least $q$. Then 
$C'=\{s_1,s'_1\}\cup\{s_2,s'_2\}\udots \{s_\mu,s'_\mu\}$
is a cocircuit of $M_3$. 
Let $H=E(M_3)-C'$ and let $B$ be a basis of $M_3|H$ 
containing $t$. Consider $\si(M_3/(B-\{t\}))$. 
We may assume that $t\in \si(M_3/(B-\{t\}))$.
Evidently $r(\si(M_3/(B-\{t\}))=2$. Moreover, 
if $i\in\{1,2,\ldots,\mu\}$, the pair $\{s_i,s'_i\}$
is independent in $\si(M_3/(B-\{t\})$, as otherwise
we have, at some stage, contracted a point on a line spanned
by $\{s_i,s'_i\}$. As $C'$ is a cocircuit, such a point
must have been parallel with $t$ and we contradict the
fact that $t$ is not a loop of $\si(M_3/(B-\{t\}))$.
Hence $\si(M_3/(B-\{t\}))\cong U_{2,2\mu+1}$ and in this, the final
case, we also conclude that $M$ has a $U_{2,q+2}$-minor.
\end{proof}

\section{Summing Up}

At last we are able to achieve the goals of this chapter.
We first consider Theorem~\ref{unavoidable-swirl}.

\begin{proof}[Proof of Theorem~\ref{unavoidable-swirl}]
Let $f_{\ref{unavoidable-swirl}}(m,q)
=f_{\ref{star-block}}(1,f_{\ref{kill-star}}(q),m,q)$.
Assume that $n\geq f_{\ref{unavoidable-swirl}}(m,q)$ and that
the 3-connected matroid
$M$ has a partition into $n$ pairwise-disjoint clonal pairs.
Assume that $M$ does not have a clonal $\Delta_m$-minor.
By Theorem~\ref{star-block},
either $M$ or $M^*$
has a minor $N$ with an associated triple
$(N,\PP,J)$ such that $(N,\PP,J)$ is a blocked 
$(1,f_{\ref{kill-star}}(m,q))$-copaddle
where all of the series pairs of $N\ba J$ are $M$-clonal.
Say $\PP=(P_1,P_2,\ldots,P_t)$. Let $S$ be the maximal
star of order $f_{\ref{kill-star}}(m,q)$ in $P_1$.
Evidently $N|S\cong M(K_{2,f_{\ref{kill-star}}}(m,q))$.
The hypotheses of Lemma~\ref{kill-star} are satisfied and
by that lemma $N\notin\eq$. This contradicts the assumption that
$M\in \eq$. It follows from this contradiction that $M$ has a 
clonal $\Delta_m$-minor.
\end{proof}

We conclude by observing again that Theorem~\ref{klonal1} is
an immediate corollary of Theorem~\ref{unavoidable-swirl}
and Lemma~\ref{3-klonal}.

\chapter{Strict Paths of $2$-separations}
\label{strict-2-paths}

\section{Introduction}

Let $M$ be a connected matroid. A {\em path of $2$-separations}
\index{path of $2$-separations}
in $M$ is a partition $\PP=(P_0,P_1,\ldots,P_l)$ of $E(M)$ into
subsets such that 
$\lambda(P_0\cup P_1\udots P_i)=1$ for all $i\in\{0,1,\ldots,l-1\}$.
Adapting
terminology from Chapter~\ref{paths-of-3-separations}, we say, 
for $i\in\{0,1,\ldots,l\}$,
that $P_i$ is a {\em step}
\index{step of a path of $2$-separations} 
of $\PP$, that $P_0$ and $P_l$ are 
{\em end}
\index{end step of a path of $2$-separations} 
steps and otherwise $P_i$ is an {\em internal} step.
\index{internal step of a path of $2$-separations}
A 2-separation $(X,Y)$ of $M$ is {\em displayed}
\index{displayed $2$-separation} 
by  $\PP$ if 
$(X,Y)=(P_0\cup P_1\udots P_i,P_{i+1}\cup P_{i+2}\udots P_l)$
for some $i\in\{0,1,\ldots,l-1\}$. Let $\QQ$ be a path of 
$2$-separations in $M$. Then $\QQ$ is a {\em concatenation} 
of $\PP$ if every 2-separation displayed by $\QQ$ is also 
displayed by $\PP$. If $\PP$ has $l$ steps,
then $\PP$ has {\em length} $l-1$.

So far we have simply extended terminology from paths of 
3-separations. The next two definitions are of specific 
importance to this chapter. Let $\PP$ be a path of 
2-separations in the connected matroid $M$.
A $2$-separation $(X,Y)$ of $M$ is {\em $\PP$-relevant}
\index{$\PP$-relevant $2$-separation} 
if there is a $j\in\{1,2,\ldots,n-1\}$ such that either
$X$ or $Y$ has the form $P_0\cup P_1\udots P_{j-1}\cup P'_j$
for some subset $P'_j$ of $P_j$.
The path $\PP$ is {\em strict}
\index{strict path of $2$-separations} 
if
\begin{itemize} 
\item[(i)]  $\lambda(P_i)=2$ for all $i\in\{1,2,\ldots,l-1\}$, and
\item[(ii)] if $(X,Y)$ is a $2$-separation of $M$,
that is not $\PP$-relevant, then either $X$ or $Y$ is contained
in a step of $\PP$.
\end{itemize}
We are particularly interested in paths whose steps
contain clonal pairs. A clonal pair $\{p_i,q_i\}$
contained in the internal step $P_i$ is {\em $\PP$-strong}
\index{$\PP$-strong clonal pair}
if $\kappa_{M}(\{p_i,q_i\},E(M)-P_i)=2$.

The results of this chapter focus on bridging 
strict paths of 2-separations 
whose steps contain strong clonal pairs. 
It is clearly possible to have
a matroid in $\eq$ in which all the steps of such a path
are bridged. An example is a free swirl.
But if constraints are placed on the way that the 2-separations
are bridged, then we do obtain certificates that prove that
a matroid is not in $\eq$. In this chapter we obtain a sequence
of such results leading to Lemma~\ref{get-simple}. This lemma
will, in turn, be used in Chapter~\ref{last-rites} to prove that
large free swirls that are bridged in certain ways cannot be in
$\eq$.

\section{Basic Facts}

We begin by obtaining some basic facts on strict paths of 
$2$-separations. 

\begin{lemma}
\label{2path0}
Let $\PP=(P_0,P_1,\ldots,P_l)$ be a strict path of 
$2$-separations in the connected matroid $M$. If
$i\in\{1,2,\ldots,l\}$, then the following hold.
\begin{itemize}
\item[(i)] $\sqcap(P_0\cup P_2\udots P_{i-1},P_{i+1}\cup P_{i+2}
\udots P_l)=0$.
\item[(ii)] $\sqcap^*(P_0\cup P_2\udots P_{i-1},P_{i+1}\cup P_{i+2}
\udots P_l)=0$.
\end{itemize}
\end{lemma}

\begin{proof}
Let $Q_0=P_0\cup P_1\udots P_{i-1}$
and $Q_l=P_{i+1}\cup P_{i+2}\udots P_l$.
By Lemma~\ref{lambda-meet}
$\lambda(Q_0\cup Q_l)=\lambda(Q_0)+\lambda(Q_l)
-\sqcap(Q_0,Q_l)-\sqcap^*(Q_0,Q_l)$. The lemma follows
from this and the fact that $\lambda(Q_0)=\lambda(Q_l)=1$.
\end{proof}

We omit the easy arguments establish the properties
of the next lemma.

\begin{lemma}
\label{2path1}
Let $\PP=(P_0,P_1,\ldots,P_l)$ be a strict path of $2$-separations
in the connected
matroid $M$. Then the following hold.
\begin{itemize}
\item[(i)] $\PP$ is a strict path of $2$-separations in $M^*$.
\item[(ii)] If $0<i\leq j<n$, then 
$\lambda(P_i\cup P_{i+1}\udots P_j)=2$.
\item[(iii)] If $\QQ$ is a concatenation of $\PP$,
then $\QQ$ is a strict path of $2$-separations in $M$.
\end{itemize}
\end{lemma}

Presumably the next lemma is well known. We omit the easy
proof. It is a generalisation of the fact that if 
$x$ is an element of the connected matroid $M$, then
either $M\ba x$ or $M/x$ is connected.

\begin{lemma}
\label{closed-coclosed}
Let $(A,B)$ be a $2$-separation of the connected matroid
$M$. Then either $M\ba A$ or $M/A$ is connected. 
\end{lemma}

The next lemma follows from Lemma~\ref{closed-coclosed}.

\begin{lemma}
\label{reduce-2path}
If $\PP=(P_0,P_1,\ldots,P_n)$ is a strict path
of $2$-separations in the connected matroid
$M$, and $1\leq i<n$, then, $(P_{i+1},P_{i+2},\ldots,P_n)$
is a strict path of $2$-separations
in either $M\ba (P_0\cup P_1\udots P_i)$
or $M^*\ba (P_0\cup P_1\udots P_i)$.
\end{lemma}

\begin{lemma}
\label{2path-lambda}
Let $\PP=(P_0,P_1,\ldots,P_l)$ be a strict path of $2$-separations
of the connected matroid $M$, let $J$ be a 
proper subset of $\{1,\ldots,l\}$ and let 
$Z=\cup_{j\in J}P_j$. Then the following hold.
\begin{itemize}
\item[(i)] If there exists $i\in\{1,2,\ldots,l-1\}$ such that
$J=\{1,2,\ldots,i\}$ or $J=\{i+1,i+2,\ldots,l\}$, then $\lambda(Z)=1$.
\item[(ii)] If there exist $i,j\in\{2,3,\ldots,l-1\}$ with $i\leq j$
and $J=\{i,i+1,\ldots,j\}$ or $J=\{1,2,\ldots,i-1,j+1,j+2,\ldots,l\}$,
then $\lambda(Z)=2$.
\item[(iii)] If neither (i) nor (ii) hold, then $\lambda(Z)>2$.
\end{itemize}
\end{lemma}

\begin{proof}
Part (i) follows from the definition. Part (ii) is just 
Lemma~\ref{2path1}(ii). Consider part (iii). 
Assume that $Z$ satisfies neither (i) nor (ii).
The result is vacuous if $\PP$ has length one. Assume for induction
that (iii) holds if $\PP$ has length at most $l-1$.
Up to complementation
we may assume that $1\notin J$. Let $i+1$ be the least integer in 
$J$.  Let $Y=E(M)-Z$ and let
$Y'=P_1\cup P_2\udots P_i$. By Lemma~\ref{2path1}(i) and 
Lemma~\ref{reduce-2path}, we may assume that 
$M\ba Y'$ is connected and that 
$(P_{i+1}\cup P_{i+2}\udots P_n)$ is a strict path of $2$-separations
in this matroid. By the choice of $J$ and either (ii) or the
induction assumption we have 
$\lambda_{M\ba Y'}(Z,Y-Y')\geq 2$. Note that
$\sqcap(Y',Y-Y')=0$, so that $r(Y)=r(Y-Y')+r(Y')$.
But $r(M)=r(M\ba Y')+r(Y')-1$. The lemma now follows from an
easy calculation.
\end{proof}

Let $\PP=(P_0,P_1,\ldots,P_l)$ be a path of 
$2$-separations in the matroid $M$. A subset $X$ of $E(M)$
is 
\index{initial set of steps}
\index{terminal set of steps}
\index{consecutive set of steps}
\index{coconsecutive set of steps}
\begin{itemize}
\item[(i)] an {\em initial} set of $\PP$ of 
$X=P_0\cup P_1\udots P_i$ for some $i\in\{0,1,\ldots,l\}$;
\item[(ii)] a {\em terminal} set of $\PP$ if
$E(M)-X$ is an initial set;
\item[(iii)] a {\em consecutive} set of $\PP$ if
$X=P_i\cup P_{i+1}\udots P_j$ where $0\leq i\leq j\leq l$;
and a 
\item[(iv)] a {\em coconsecutive} set of $\PP$
if $E(M)-X$ is a consecutive
set.
\end{itemize}

The next lemma is a version of our old friend
Lemma~\ref{modular} for paths rather than flowers.

\begin{lemma}
\label{2-path-modular}
Let $e$ be an element of the matroid $M$ such that 
$M\ba e$ is connected with a strict path 
$\PP=(P_0,P_1,\ldots,P_n)$ of $2$-separations.
Let $X$ and $Y$ be subsets of $E(M\ba e)$. Assume that
both $X$ and $Y$ are unions of steps and  
that $e\in\cl(X)$ and $e\in\cl(Y)$. 
\begin{itemize}
\item[(i)] If $X$ is the union of an initial set of steps,
$Y$ is  the union of a terminal set of steps,
and $X\cap Y\neq\emptyset$, then $e\in\cl(X\cap Y)$.
\item[(ii)] If $X$ and $Y$ are the union of either a 
consecutive or coconsecutive set of petals,
$X\cup Y\neq E(M)-\{e\}$, and $X\cap Y\neq \emptyset$,
then $e\in\cl(X\cap Y)$.
\end{itemize}
\end{lemma}

\begin{proof} In either case it suffices to show that 
$X$ and $Y$ form a modular pair. In case (i) we have
$\lambda(X)=\lambda(Y)=1$, $\lambda(X\cap Y)=2$ and 
$\lambda(X\cup Y)=0$. Thus (i) holds.

Consider (ii). Since
$X\cup Y\neq E(M)-\{e\}$ and $X\cap Y\neq \emptyset$,
it is readily deduced that 
both $X\cup Y$ and $X\cap Y$ are unions of consecutive or
coconsecutive sets of steps of $\PP$. It now follows from
Lemmas~\ref{2path1}(iii) and \ref{2path-lambda}(ii), that
$X$ and $Y$ form a modular pair.
\end{proof}

The next lemma establishes a connection between 
strict paths of 2-separations and swirl-like flowers. 
We omit the easy proof.

\begin{lemma}
\label{get-swirl}
Let $e$ be an element of the matroid $M$ such that $M\ba e$
is connected with a strict path 
$\PP=(P_0,P_1,\ldots,P_l)$ of $2$-separations. If 
$e\in\cl(P_0\cup P_l)$, but $e\notin \cl(P_0)$ and
$e\notin \cl(P_l)$, then $(P_l\cup P_0\cup\{e\},P_1,\ldots,P_{l-1})$
is a swirl-like flower in $M$.
\end{lemma}

If $\PP=(P_1,P_2,\ldots,P_n)$ is a swirl-like flower in the
connected matroid $M$, then it is easily seen that, if
$P_1$ is coclosed, then $M\ba P_1$ is connected.
We also omit the easy proof of the next lemma.

\begin{lemma}
\label{get-2-path}
Let $\PP=(P_1,P_2,P_3,\ldots,P_n)$ be a swirl-like flower
in the connected matroid $M$. If $P_1$ is coclosed,
then $(P_2,P_3,\ldots,P_n)$ is a strict path of $2$-separations
in the
connected matroid $M\ba P_1$.
\end{lemma}

\begin{lemma}
\label{kappa-kappa}
Let $\PP=(P_0,P_1,\ldots,P_l)$ be a strict path of 
$2$-separations in the connected matroid $M$, and
let $\{p_i,q_i\}$ be a $\PP$-strong clonal pair. Then,
$\kappa(\{p_i,q_i\},P_0\cup P_l)=2$.
\end{lemma}

\begin{proof}
If the lemma does not hold, then there is a 2-separation 
$(X,Y)$ of $M$ with $P_0\cup P_l\subseteq X$ and
$\{p_i,q_i\}\subseteq Y$. Such a separation is not 
$\PP$-relevant so it must be the case that
there is a petal $P_i$ such that $Y\subseteq P_i$.
But this contradicts the definition of $\PP$-strong
clonal pair.
\end{proof}

\begin{lemma}
\label{connect-2paths-3paths}
Let $N$ be a connected minor of the connected matroid $M$
and let $\PP=(P_0,P_1,\ldots,P_l)$ be a strict path of
$2$-separations in $N$. Assume that every $\PP$-relevant
$2$-separation of $N$ is bridged in $M$, and that
$\QQ=(Q_0,Q_1,\ldots,Q_m)$ is a partition of $E(M)$
such that, $P_0\subseteq Q_0$, 
$P_m\subseteq Q_m$ and $\lambda_M(Q_0\cup Q_1\udots Q_i)=2$ 
for $i\in\{0,1,\ldots,m-1\}$. Then the following hold.
\begin{itemize}
\item[(i)] $\QQ$ is a path of $3$-separations in $M$.
\item[(ii)] If $\{p,q\}$ is an $M$-clonal pair of 
$N$ that is $\PP$-strong, then it is also $\QQ$-strong.
\end{itemize}
\end{lemma}

\begin{proof}
Assume that $(Q'_0,Q'_m)$ is a $2$-separation of $M$
with $Q_0\subseteq Q'_0$ and $Q_m\subseteq Q'_m$. Then
$(Q'_0,Q'_m)$ is induced by a $2$-separation
$(P'_0,P'_l)$ with $P_0\subseteq P'_0$ and $P_l\subseteq P'_l$.
But by the definition of strict path, such $2$-separations
are $\PP$-relevant and hence bridged in $M$. Thus 
(i) holds. Part (ii) follows from Lemma~\ref{kappa-kappa}.
\end{proof}

\section{A First Certificate}

We now develop some certificates for showing that a matroid
is not in $\eq$. Interpreted in this way the next lemma
says that if we block all the 2-separations displayed by
a strict path of $2$-separations, and we do not keep a large
displayed swirl-like flower, then we are not in $\eq$. 

\begin{lemma}
\label{2path2}
Let $M$ be a matroid in $\eq$ with an element $b$ such that
$M\ba b$ is connected with a strict path 
$\PP=(P_0,P_1,\ldots,P_n)$
of $2$-separations each internal step of which contains a
$\PP$-strong $M$-clonal pair. 
Assume that $b$ blocks both $P_0$ and $P_n$.
Then there is a function $f_{\ref{2path2}}(m,q)$ such that,
if $n\geq f_{\ref{2path2}}(m,q)$, then the following holds.
There exist an  $s\in\{1,2,\ldots,n-2\}$ such that $s+m<n$
with the property that 
$$\RR=(P_0\cup P_1\udots P_s\cup P_{s+m}\cup P_{s+m+1}
\udots P_n\cup\{b\},P_{s+1},\ldots,P_{s+m-1})$$ 
is a swirl-like flower in $M$ of order $m$. Moreover, each petal
of $\RR$ contains an $\RR$-strong $M$-clonal pair.
\end{lemma}

\begin{proof}
Let $f_{\ref{2path2}}(m,q)=f_{\ref{at-last}}(m+1,q)+2$.
Assume that $n\geq f_{\ref{2path2}}(m,q)$.
Note that,
if $(X,Y)$ is a 2-separation of $M\ba b$ where 
$P_0\subseteq X$ and $P_n\subseteq Y$, then $b$
blocks $(X,Y)$. Hence
$\PP'=(P_0\cup\{b\},P_1,\ldots,P_n)$ 
is a path of 3-separations in $M$ as otherwise $b$
fails to block both $P_0$ and $P_n$. Say $i\in\{1,2,\ldots,n-1\}$.
Then, by the definition of $\PP$-strong for strict paths of
2-separations,  the step
$P_i$ of $\PP$ contains an $M$-clonal pair $\{p_i,p'_i\}$ such that
$\kappa_{M\ba b}(\{p_i,p'_i\},P_0\udots P_{i-1}
\cup P_{i+1}\udots P_n)=2$. Hence
$\kappa_{M}(\{p_i,p'_i\},\{b\}\cup 
P_0\udots P_{i-1}\cup P_{i+1}\udots P_n)=2$. It 
follows that the $\{p_i,p'_i\}$ is $\PP'$-strong in the sense
of $\PP'$-strong defined for paths of 3-separations.

Let 
$\PP''=(\{b\}\cup P_0\cup P_1,P_2,\ldots,P_{n-2},P_{n-1}\cup P_n)$.
Then each step of $\PP''$ contains a $\PP''$ strong
clonal pair and $\PP''$ has length $n-2$.
By Lemma~\ref{at-last}, $\PP''$ displays a 
swirl-like flower of order $m+1$. Let $\QQ$ be a maximal such
flower, say $\QQ$ has order $\mu+1\geq m+1$. 
Let $\QQ=(Q_1,Q_2,\ldots,Q_{\mu+1})$. Then, by 
Lemma~\ref{display1}, there is an integer $s\geq 2$, with
$s+\mu\leq n-2$ such that 
$\{Q_1,Q_2,\ldots,Q_{\mu+1}\}=\{\{b\}\cup P_0\cup P_1\udots P_s,
P_{s+1},\ldots, P_{s+\mu-1},P_{s+\mu}\cup P_{s+\mu+1}
\udots P_{n-1}\cup P_n\}$.
Say $\{i,i+1\}\subseteq \{s+1,s+2,\ldots,s+\mu-1\}$,
then, by the definition of strict 2-path, 
$\sqcap_M(P_i,P_{i+1})=1$, so that $P_i$ and $P_{i+1}$ are adjacent
in $\QQ$. Moreover, 
$\sqcap_M(\{b\}\cup P_0\cup P_1\udots P_s,P_{s+1})\geq 1$, so these
two petals of $\PP$ are adjacent in $\QQ$. Also
$\sqcap(P_{s+\mu-1},P_{s+\mu}\cup P_{s+\mu+1}
\udots P_{n-1}\cup P_n)\geq 1$,
and these two petals are adjacent in $\QQ$. By elimination
$\{b\}\cup P_0\cup P_1\udots P_s$ and 
$P_{s+\mu-1}\cup P_{s+\mu}\udots P_{n-1}\cup P_n$ are adjacent.
Hence
$\QQ=(\{b\}\cup P_0\cup P_1\udots P_s,
P_{s+1},\ldots, P_{s+\mu-1},P_{s+\mu}\udots P_{n-1}\cup P_n)$.
Let
$\QQ'=(\{b\}\cup P_0\cup P_1\udots P_s\cup 
P_{s+\mu}\udots P_{n-1}\cup P_n,
P_{s+1},P_{s+2},\ldots,P_{s+\mu-1})$
Every petal of this flower
contains a $\QQ'$-strong clonal pair. As $\mu\geq m$,
the lemma is satisfied by taking $\RR$ to be
an appropriate concatenation of $\QQ'$.
\end{proof}

\section{A Second Certificate}

Our next certificate requires somewhat more work.

\begin{lemma}
\label{back-block}
Let $M$ be a matroid with a set $B=\{b_2,b_3,\ldots,b_n\}$ of elements
such that $M\ba B$ is connected with a strict path
$\PP=(P_0,P_1,\ldots,P_n)$ of $2$-separations each internal step
of which contains a $\PP$-strong $M$-clonal pair. Assume that, for 
$i\in\{2,3,\ldots n\}$, the element $b_i$ is in 
$\cl(P_i\cup P_0)$, but not in $\cl(P_0)$ or $\cl(P_i)$.
Then there is a function $f_{\ref{back-block}}(q)$ such that,
if $n\geq f_{\ref{back-block}}(q)$, then $M\notin\eq$.
\end{lemma}

We begin by proving Lemma~\ref{back-block} in two special cases.
For the first case was assume that $P_0$ consists of a single
element.

\begin{lemma}
\label{2path3}
Let $M$ be a matroid with a set $B=\{b_2,b_3,\ldots,b_n\}$ of elements
such that $M\ba B$ is connected with a strict path
$\PP=(\{p_0\},P_1,\ldots,P_n)$ of $2$-separations each internal step
of which contains a $\PP$-strong $M$-clonal pair. Assume that, for 
$i\in\{2,3,\ldots n\}$, the element $b_i$ is in 
$\cl(P_i\cup\{p_0\})$, but not in $\cl(\{p_0\})$ or $\cl(P_i)$.
Then there is a function $f_{\ref{2path3}}(q)$ such that,
if $n\geq f_{\ref{2path3}}(q)$, then $M\notin\eq$.
\end{lemma}

\begin{proof}
Consider the partition
$$\PP'=(\{p_0\}\cup P_1,P_2\cup\{b_2\},\ldots,P_{n-1}
\cup\{b_{n-1}\},P_n\cup\{b_n\})$$
of $E(M)$. 
We first prove that $\PP'$ is a path of 
$3$-separations in $M$.

Say that $j\in\{1,2,\ldots,n-1\}$. 
Note that $p_0\notin \cl(P_{j+1}\cup P_{j+2}\udots P_n)$, as otherwise
$\lambda_{M\ba B}(\{p_0\}\cup P_{j+1}\udots P_n)=1$,
contradicting Lemma~\ref{2path-lambda}. 
Also $b_i\in\cl(P_i\cup\{p_0\})$ for all $i\in\{1,2,\ldots,n-1\}$.
Thus
$$r(\{p_0\}\cup P_1\cup P_2\cup\{b_2\}\udots P_{j-1}\cup\{b_{j-1}\})
=r(\{p_0\}\cup P_1\udots P_{j-1}).$$
and 
$$r(P_j\cup\{b_j\}\udots P_n\cup\{b_n\})
=r(P_j\udots P_n\cup\{p_0\})=r(P_j\udots P_n)+1.$$
Hence 
$$\lambda_M(\{p_0\}\cup P_1\cup P_2\cup\{b_2\}\udots 
P_{j-1}\cup \{b_{j-1}\})=2.$$
To prove that $\PP'$ is a path of 3-separations in $M$.
It remains to show that 
$\kappa_M(\{p_0\}\cup P_1,P_n\cup\{b_n\})=2$. 
Assume not. Then, for some $i\in\{1,\ldots,n\}$
and subset $P_i'$ of $P_i\cup\{b_i\}$, we have
$\lambda_M(\{p_0\}\cup P_1\cup P_2\cup\{b_2\}\udots 
P_{i-1}\cup \{b_{i-1}\}\cup P'_i)=1$.
But 
$\lambda_{M\ba B}(\{p_0\}\cup P_1\cup P_2\udots P_{i-1}
\cup P'_i)=1$,
so that $b_n$ does not block this 2-separation. Hence 
either 
$b_n\in\cl(\{p_0\}\cup P_1\cup P_2\udots P_{i-1}
\cup P'_i)$ or 
$b_n\in\cl((P_i-P'_i)\cup P_{i+1}\udots P_n)$.
But $b_n\in\cl(\{p_0\}\cup P_n)$. By 
Lemma~\ref{2-path-modular}, we obtain the contradiction
that either $b_n\in\cl(\{p_0\})$ or $b_n\in\cl(P_n)$.

If $i\in\{2,3,\ldots,n\}$, then $b\in\cl(P_i\cup\{b_i\})$.
Therefore, if $P'_i$ and $P'_j$ are distinct steps of the path
$\PP'$, then $\sqcap(P'_i,P'_j)>0$. This shows that 
$\PP'$ does not display any 4-petal swirl-like flower.
Moreover, it is evident that if $i\in\{2,3,\ldots,n-1\}$,
then any $\PP$-strong clonal pair in $P_i$ is 
$\PP'$-strong.  The 
lemma now follows by letting 
$f_{\ref{2path3}}(q)= f_{\ref{at-last}}(4,q)+2$.
\end{proof}

For the second special case we assume that $P_0$ consists
of a series class that is blocked in a particular way.

\begin{lemma}
\label{2path4}
Let $M$ be a matroid with a set $B=\{b_1,b_2,\ldots,b_n\}$
such that $M\ba B$ is connected with a strict path 
$\PP=(P_0,P_1,\ldots,P_n)$ of $2$-separations,
where $P_0=\{s_1,s_2,\ldots,s_n\}$ is a series class and
each internal step of $\PP$ contains a $\PP$-strong
$M$-clonal pair.
Assume that, for $i\in\{2,3,\ldots,n\}$, the element 
$b_i$ is in $\cl(P_i\cup\{s_i\})$, but $b_i\notin\cl(P_i)$
and $b_i\notin\cl(\{s_i\})$. Then there is a function
$f_{\ref{2path4}}(q)$ such that, if $n\geq f_{\ref{2path4}}(q)$,
then $M\notin \eq$.
\end{lemma}

\begin{proof}
We first note.

\begin{sublemma}
\label{2path4.0}
\begin{itemize}
\item[(i)] If $2\leq i\leq n$, then 
$b_i\notin\cl(P_0\cup P_1\udots P_{i-2}).$ In particular
$\cl(P_0)\cap B=\emptyset$.
\item[(ii)] Say $P'_0\subseteq P_0$ and $|P'_0|\geq 2$.
If $s_i\in P'_0$, then 
$b_i\notin\cl((P_0-P'_0)\cup P_1\cup P_2\udots P_n)$.
\end{itemize}
\end{sublemma}

\subproof
Assume that $b_i\in\cl(P_0\cup P_1\udots P_{i-2})$. 
Consider the pair of sets $P_0\cup P_1\udots P_{i-2}$
and $P_i\cup\{s_i\}$. As $s_i$ is in a non-trivial
series class of $M\ba B$ that does not meet $P_i$, we have
$r(P_i\cup\{s_i\})=r(P_i)+1$. The union of the two sets
is $P_0\cup P_1\udots P_{i-2}\cup P_i$. By Lemma~\ref{2path0}
the rank of this set is $r(P_0\cup P_1\udots P_{i-2})+r(P_i)$.
Moreover, the intersection of the sets is $\{s_i\}$.
We deduce that the two sets form a modular pair. As
$b_i$ is in the span of each we obtain the contradiction that
$b_i\in\cl(\{s_i\})$. Hence (i) holds.

Consider (ii). Assume that $|P'_0|\geq 2$ and that $s_i\in P'_0$.
As $b_i\in\cl(P_i\cup\{s_i\})$ but not in $\cl(P_i)$,
we have $s_i\in\cl(P_i\cup\{b_i\})$. Thus,
if $b_i\in\cl((P_0-P'_0)\cup P_1\udots P_n)$, then
we have $s_i\in\cl((P_0-P'_0)\cup P_1\udots P_n)$,
contradicting the fact that $P'_0$ is a series set in 
$M\ba B$.
\end{proof}

\begin{sublemma}
\label{2path4.1}
If $(X,Y)$ is a $2$-separation of $M$, then,
for some $Z\in\{X,Y\}$, either $Z\subseteq P_i$ for an
$i\in\{2,3,\ldots,n-2\}$ or $Z\subseteq P_{n-1}\cup P_n$.
\end{sublemma}

\subproof
Consider the 2-separation $(X,Y)$. Let
$(X',Y')=(X-B,Y-B)$. Assume that $|X'|\leq 1$.
Then, as $r(M\ba B)=r(M)$, and $M\ba B$ is connected,
we see that $X\subseteq \cl(Y)$. Hence $r(X)=1$.
But it is evident that no member of $B$ is in a non-trivial
parallel class. Hence $|X-B|>1$ and also $|Y-B|>1$,
so that
$(X',Y')$ is a 2-separation of $M\ba B$.

Assume that $(X',Y')$ is $\PP$-relevant, so that,
up to labels, 
$(X',Y')=(P_0\cup P_1\udots P_{i-1}\cup P'_i,
P''_i\cup P_{i+1}\udots P_n)$. If $i\neq n-1$, then,
by \ref{2path4.0}, $(X',Y')$ is blocked by a member of 
$B$. Thus $i=n-1$, and $Y'\subseteq P_{n-1}\cup P_n$.
In this case, it also follows from \ref{2path4.0} that
$\cl_M(Y')\cap B=\emptyset$, so that $Y'=Y$ and the claim
holds in this case. 

On the other hand, if $(X',Y')$ is not $\PP$-relevant then,
by the definition of strict path of 2-separation, 
we may assume that
$X'\subseteq P_i$ for some $i\in\{0,1,\ldots,n\}$. Again by
\ref{2path4.0}, $\cl_M(X')\cap B=\emptyset$, so that $X'=X$.
If $X'\subseteq P_0$, then it follows from \ref{2path4.0}
that $X'$ is blocked by at least one member of $B$. Thus
$X'\not\subseteq P_0$ and the sublemma follows.
\end{proof}

\begin{sublemma}
\label{2path4.2}
For all $i\in\{1,2,\ldots,n-1\}$, we have
$\lambda_M(P_1\cup P_2\udots P_i\cup\{s_1,s_2,\ldots,s_i\}
\cup\{b_1,b_2\ldots,b_i\})=2$.
\end{sublemma}

\subproof
By the definition of paths of 2-separations,
$\lambda_{M\ba B}(P_1\cup P_2\udots P_i\cup\{s_1,s_2,\ldots,s_n\})=1$.
Since $\{s_1,s_2,\ldots,s_n\}$ is a series class,
and $\sqcap(\{s_1,s_2,\ldots,s_n\},P_1)=1$, we have
$r(P_1\cup P_2\udots P_i\cup\{s_1,s_2\ldots,s_n\})
=r(P_1\cup P_2\udots P_i)+n-1$.
But $r(P_1\cup P_2\udots P_i\cup\{s_1,s_2,\ldots,s_i\})
=r(P_1\cup P_2\udots P_i)+i$,
and $r(P_{i+1}\cup P_{i+2}\udots P_n\cup\{s_{i+1},s_{i+2}\ldots,s_n\})
=r(P_{i+1}\cup P_{i+2}\udots P_n)+n-i$.
Hence 
$\lambda_{M\ba B}(P_1\cup P_2\udots P_i\cup\{s_1,s_2,\ldots,s_i\},
P_{i+1}\cup P_{i+2}\udots P_n\cup
\{s_{i+1},s_{i+2}\ldots,s_n\})=2$. The sublemma is
a straightforward consequence of this observation and the
fact that $b_j\in\cl(P_j\cup\{s_j\})$ for all $j\in\{1,2,\ldots,n\}$.
\end{proof}

By \ref{2path4.2} and \ref{2path4.1}, 
$\PP'=(P_1\cup\{s_1,b_1\},P_2\cup\{s_2,b_2\},\ldots,
P_{n-2}\cup\{s_{n-2},b_{n-2}\},P_{n-1}\cup 
P_n\cup\{s_{n-1},s_n,b_{n-1},b_n\})$
is a path of 3-separations in $M$. If this path displays a 
swirl-like flower with at least four petals, then there 
exists an $i\in\{2,3,\ldots,n-2\}$ such that 
$\lambda_M(P_i\cup\{s_i,b_i\})=2$, so that
$\lambda_{M\ba B}(P_i\cup\{s_i\})=2$.
But $s_i\notin\cl(P_i)$ and 
$s_i\in\cl(P_1\cup\{s_1,s_2,\ldots,s_{i-1},s_{i+1},\ldots,s_n\})$
so that 
$\lambda_{M\ba B}(P_i\cup\{s_i\})=\lambda_{M\ba B}(P_i)+1$.
By Lemma~\ref{2path1}, $\lambda_{M\ba B}(P_i)=2$.
Hence $\lambda_{M\ba B}(P_i\cup\{s_i\})=3$
and $\PP'$ does not display any 4-petal swirl-like
flowers. 

To complete the proof we need to show that our $M$-clonal
pairs are $\PP'$-strong. Say $i\in\{2,3,\ldots,n-2\}$,
and $\{p_i,p'_i\}\subseteq P'_i$ is a $\PP$-strong clonal
pair. If
$\kappa_M(\{p_i,p'_i\},E(M)-(P_i\cup\{b_i,s_i\}))<2$,
then there is a $2$-separating set $Z$ of $M$ such that
$\{p_i,p'_i\}\subseteq Z\subseteq P_i\cup\{b_i,s_i\}$.
By \ref{2path4.1}, $Z\subseteq P_i$ and is 2-separating in 
$M\ba B$ contradicting the definition of $\PP$-strong.

The lemma now 
follows by letting $f_{\ref{2path4}}(q)=f_{\ref{at-last}}(4,q)+3$.
\end{proof}

Before proving Lemma~\ref{back-block} we note an elementary
fact.

\begin{lemma}
\label{keep-circuit}
Let $C$ be a circuit of the connected matroid $M$.
If $x\in E(M)-C$ then there is an $N\in \{M\ba x,M/x\}$ 
such that $N$ is connected and
$C$ is a circuit of $N$.
\end{lemma}

\begin{proof}
If $M\ba x$ is connected, let $N=M\ba x$.
Otherwise say $(X,Y)$ is a separation of $M\ba x$, 
where $C\subseteq X$.
If $x\in \cl(C)$, then $(X\cup\{x\},Y)$ is a separation of
$M$. Thus $x\notin \cl(C)$ and the lemma holds by setting $N=M/x$.
\end{proof}

We also recall the theorem of Lemos and 
Oxley \cite{oxlem01} that we have
stated in this paper as Theorem~\ref{oxley-lemos}.
The next result is a straightforward consequence of 
this theorem. We omit the routine proof.

\begin{corollary}
\label{ox-lem-cor}
Let $(A,B)$ be a $2$-separation of the connected
matroid $M\in \eq$. Then there is a function 
$f_{\ref{ox-lem-cor}}(m,q)$ such that, if 
$|A|\geq f_{\ref{ox-lem-cor}}(m,q)$, then
$A$ contains either a circuit, a cocircuit, a parallel set,
or a series set  with at least $m$ elements.
\end{corollary}

\begin{proof}[Proof of Lemma~\ref{back-block}]
Let $\rho=\max\{f_{\ref{2path3}}(q),f_{\ref{2path4}}(q)\}$ 
and let 
$f_{\ref{back-block}}(q)=f_{\ref{2path3}}(q)
f_{\ref{ox-lem-cor}}((q+1)\rho).$
Assume that $n\geq f_{\ref{back-block}}(q)$.

To facilitate the proof, we make some local definitions. 
Let $M'$ be a minor of $M$ obtained by removing a proper subset
of elements of $P_0$, let $N'=M'\ba B$, 
and let $P'_0=P_0\cap E(N')$. Then $M'$ is an {\em allowable}
\index{allowable minor}
minor of $M$ if $N'$ is connected and $B$ is coindependent in $M'$.

We will always use the convention that,
if $M'$ is an allowable minor of $M$, then 
$P'_0=P_0\cap E(M')$ and $N'=M'\ba B$. The next claim is evident.

\begin{sublemma}
\label{back-block1}
Let $M'$ be an allowable minor of $M$.
Say $e\in P'_0$ and $|P'_0|>1$. If $N'\ba e$ is 
connected, then $M'\ba e$ is allowable and if $N'/e$ is connected,
then $M'/e$ is allowable. Moreover, either $M'\ba e$ or $M'/e$
is allowable.
\end{sublemma}

Let $M''$ be an arbitrary minor of $M$ obtained by removing 
elements of $P_0$ and let $P''_0=E(M'')\cap P_0$. For an
element $x$ of $P''_0$, define
$$d_{M''}(x)=|\{P_i:i\in\{1,2,\ldots,m\},
x\in\cl_{M''}(P_i\cup\{b_i\})\}|.$$
and define
$$d(M'')=|\{P_i:i\in\{1,2,\ldots,n\},b_i\in
\cl_{M''}(P_0\cup P_i),b_i\notin\cl_{M''}(P_i)\}|.$$
Let $x$ and $z$ be distinct elements of $P''_0$.
The next sublemma is elementary.

\begin{sublemma}
\label{back-block2}
\begin{itemize}
\item[(i)] $d_{M''\ba z}(x)=d_{M''}(x)$ and, if 
$z$ is not a coloop of $M''$, then $d(M''\ba z)=d(M'')$.
\item[(ii)]  $d_{M''/z}(x)\geq d_{M''}(x)$ and  $d(M''/z)=d(M'')-d(z)$.
\end{itemize}
\end{sublemma}

\begin{sublemma}
\label{back-block3}
Let $M'$ be an allowable minor of 
$M$, and let $p$  be an element of $P'_0$. 
If $d_{M'}(p)>f_{\ref{2path3}}(q)$, then
$M\notin {\mathcal E}(q)$.
\end{sublemma}

\subproof
Consider the partition
$$\PP'=(\{p\},(P'_0-\{p\})\cup P_1,P_2,\ldots,P_n).$$ 
It is 
easily checked that this is a strict path of 2-separations 
each internal step  of which contains a 
$\PP'$-strong, $M'$-clonal pair. But now the hypotheses of
Lemma~\ref{2path3} hold for $M'$, $N'$ and $\PP'$, so 
by that lemma $M'\notin\eq$. Hence $M\notin \eq$.
\end{proof}

\begin{sublemma}
\label{back-block4}
There is an allowable minor $M'$ of $M$ that
has the following properties: 
$P'_0$ has no non-trivial parallel classes; $d(M')=d(M)$;
and  $d_{M'}(p)>0$ for all $p\in P'_0$.
\end{sublemma}

\subproof
Let $M'$ be an allowable minor of $M$ that has the properties
that $d(M')=d(M)$ and that $|E(M')|$ is minimal. It follows from
\ref{back-block1} and \ref{back-block2} that
$M'$ satisfies the claim.
\end{proof}

From now on we assume that $M'$ is an allowable minor of $M$
satisfying \ref{back-block4}.

\begin{sublemma}
\label{back-block5}
If $r_{M'}(P'_0)\leq f_{\ref{ox-lem-cor}}((q+1)\rho)$, 
then $M\notin \eq$.
\end{sublemma}

\subproof
Recall that 
$n\geq f_{\ref{2path3}}(q)f_{\ref{ox-lem-cor}}((q+1)\rho)$. 
If $M'$ contains an allowable
minor $M''$ with an element $z\in E(M'')\cap P_0$
such that $d_{M''(x)}\geq f_{\ref{2path3}}(q)$, then the claim
follows from \ref{back-block3}. Thus we may assume that this 
never occurs. In this case we apply \ref{back-block1} and
\ref{back-block2} to obtain an allowable minor 
$M''$, where $E(M'')\cap P_0=\{p_0\}$. But, in this case
$d_{M''}(p_0)=d(M'')$ and $d(M'')\geq f_{\ref{2path3}(q)}$.
In this case $M\notin\eq$ by Lemma~\ref{2path3}.
\end{proof}

\begin{sublemma}
\label{back-block5.5}
If $P_0'$ contains a series class of $N'$ of size at least
$f_{\ref{2path4}}(q)$, then $M\notin \eq$.
\end{sublemma}

\begin{proof}
Let $P_0'=S$, $P'_1=(P_0-P'_0)\cup P_1$,
and $\PP'=(P'_0,P'_1,P_2,\ldots,P_n)$. Then $\PP'$, $N'$ and $M'$
satisfy the hypotheses of Lemma~\ref{2path4}, so that, by that
lemma $M\notin \eq$.
\end{proof}

\begin{sublemma}
\label{back-block6}
If $P_0'$ contains a circuit $C$ of $N'$ with at least 
$(q+1)f_{\ref{2path4}}(q)$ elements, then $M\notin \eq$. 
\end{sublemma}

\subproof
Say $P'_0\neq C$. Then by Lemma~\ref{keep-circuit} and
\ref{back-block1}, there is an element
$z\in P'_0$ such that either $M'\ba z$ or $M'/z$ is allowable
with $C$ as a circuit. By this fact and \ref{back-block2},
we lose no generality in assuming that $P'_0=C$. 
As $\lambda_{N'}(C)=1$  we have $r^*_{N'}(C)=2$.  If 
$M\in\eq$, then $C$ contains at most $q+1$ series classes
of $N'$. As $|C|\geq (q+1)f_{\ref{2path4}}$, 
there is a series class $S$
in $C$ of size at least $f_{\ref{2path4}}(q)$. 
By \ref{back-block5.5} $M\notin\eq$.
\end{proof}

An easy argument that we omit shows that

\begin{sublemma}
\label{back-block7}
If $M''$ is an allowable minor of $M'$ and the element 
$p\in P'_0\cap E(M'')$ belongs to a parallel class of size $l$,
then $d_{M''}(p)\geq l$.
\end{sublemma}

\begin{sublemma}
If $P'_0$ has a cocircuit $C$ of $N'$ with at least 
$(q+1)f_{\ref{2path3}}(q)$ elements,
then $M\notin \eq$.
\end{sublemma}

\subproof
Arguing as in \ref{back-block6}, but using the dual of 
Lemma~\ref{keep-circuit}, we obtain an allowable minor
$M''$ with $P'_0\cap E(M'')=C$ such that $C$ is a cocircuit in 
$M''$. But then, $r_{N'}(C)\leq 2$. So that $C$ has a parallel
set of size at least $f_{\ref{2path3}}(q)$. 
Say that $p$ is an element of such
a parallel set. Then, by \ref{back-block7}, 
$d_{M''}(p)\geq f_{\ref{2path3}}(q)$.
Let $\PP''=(\{p\},P_1\cup P_0-\{p\},P_2,\ldots,P_n)$.
With this path of 2-separations in $M''\ba B$ the hypotheses of
Lemma~\ref{2path3} are satisfied. Thus $M\notin \eq$.
\end{proof}

If $r(P'_0)\leq f_{\ref{ox-lem-cor}}((q+1)\rho)$,
then $M\notin \eq$ by \ref{back-block5}. Thus we may assume that
$r(P'_0)> f_{\ref{ox-lem-cor}}((q+1)\rho)$. By definition 
$P'_0$ contains no non-trivial parallel sets of $N'$,
so, by  the definition of $\rho$ and Corollary~\ref{ox-lem-cor}, 
$P'_0$ contains one of the following:
a series set of $N'$ of size at least
$f_{\ref{2path4}}(q)$, a circuit of size at least
$(q+1)f_{\ref{2path4}}(q)$, or a cocircuit of $N'$ of 
size at least $(q+1)f_{\ref{2path3}}(q)$. By
\ref{back-block5.5}, \ref{back-block6} and \ref{back-block7}
respectively, we deduce in each case that $M\notin \eq$.
\end{proof}

\section{A Third Certificate}

Let $\PP$ be a strict path of 2-separations. 
Recall that a  2-separation $(X,Y)$ 
is {\em $\PP$-relevant} if either $X$ or $Y$ is of the form
$P_0\cup P_1\udots P_{i-1}\cup P_i'$ for some subset $P'_i$ of $P_i$.

\begin{lemma}
\label{2path47}
Let $M$ be a matroid with a set coindependent set
$B$ such that $M\ba B$ has a 
strict path $\PP=(P_0,P_1,\ldots,P_n)$ of $2$-separations, each 
internal step of which contains a $\PP$-strong
$M$-clonal pair. Assume that
the following hold.
\begin{itemize}
\item[(i)] For all $b\in B$, there exists $i\in\{1,2,\ldots,n\}$ 
such that $b\in\cl(P_{i-1}\cup P_i)$.
\item[(ii)] Every $\PP$-relevant $2$-separation 
of $M\ba B$ is bridged in $M$.
\end{itemize}
Then there is a function $f_{\ref{2path47}}(q)$ such that,
if $n\geq f_{\ref{2path47}}(q)$, then $M\notin \eq$.
\end{lemma}

\begin{proof}
We may assume that $B$ is minimal in that, for all $b\in B$, there
is a $\PP$-relevant 2-separation that is induced in $M\ba b$.
For $i\in\{1,2,\ldots,n\}$, let 
$Z_i=\{b\in B: b\in\cl(P_0\cup P_1\udots P_i),
b\notin\cl(P_0\cup P_1\udots P_{i-1})\}$.

\begin{sublemma}
\label{2path47-partition}
$(Z_1,Z_2,\ldots,Z_n)$ partitions $B$ into
nonempty subsets.
\end{sublemma}

\subproof
It is immediate from the definition that the members of
$(Z_1,Z_2,\ldots,Z_n)$ are pairwise disjoint. Say 
$i\in\{1,2,\ldots,n\}$. We now prove that $Z_i\neq\emptyset$.
Consider $(P_0\cup P_1\udots P_{i-1},P_i\cup P_{i+1}\udots P_n)$. 
This 2-separation of $M\ba B$ is not induced in $M$, so some member
$b$ of $B$ blocks it. For such an element
$b$ we have
$b\notin\cl(P_0\cup P_1\udots P_{i-1})$ and 
$b\notin\cl(P_i\cup P_{i+1}\udots P_n)$. But
$b\in\cl(P_{j-1}\cup P_j)$ for some $j$, so we have
$j=i$ and it follows that $b\in\cl(P_0\cup P_1\udots P_i)$
so that $b\in Z_i$.
\end{proof}

\begin{sublemma}
\label{2path47-trivial}
If $b\in Z_k$, then $b\in\cl(P_{k-1}\cup P_k)$.
\end{sublemma}

\subproof
There is a $j\in\{1,2,\ldots,n\}$ such that
$b\in\cl(P_{j-1}\cup P_j)$. By the definition of 
$Z_k$ we have $j>k-1$. If $j=k$ the claim holds. Say $j>k$.
Then $b\in\cl(P_k\cup P_{k+1}\udots P_n)$, 
so by Lemma~\ref{2-path-modular},
$b\in\cl(P_k)$. Hence $b\in\cl(P_{k-1}\cup P_k)$ in this case
too.
\end{proof}

\begin{sublemma}
\label{2path47-1}
Say $j\in\{1,2,\ldots,n\}$ and 
$b\in Z_j$. Let $(X,Y)$ be a $\PP$-relevant 
$2$-separation of $M\ba B$ that induces a $2$-separation 
$(X',Y')$ of $M\ba b$. Then the following hold.
\begin{itemize}
\item[(i)] $P_0\cup P_1\udots P_{j-2}\subseteq X$ and 
$P_{j+1}\cup P_{j+2}\udots P_n\subseteq Y$.
\item[(ii)] $Z_1\cup Z_2\udots Z_{j-3}\subseteq X'$
and $Z_{j+2}\cup Z_{j+3}\udots Z_n\subseteq Y'$.
\end{itemize}
\end{sublemma}

\subproof
By the definition of $\PP$-relevant,
$$(X,Y)=(P_0\cup P_1\udots P_{i-1}
\cup P'_i,P''_i\cup P_{i+1}\udots P_n)$$
for some $i\in\{1,2,\ldots,n-1\}$. If $i<j-1$, then 
by \ref{2path47-trivial}, $z\in\cl(Y)$. If 
$i>j$, then $x\in\cl(X)$ by the definition of $Z_j$.
In either case we contradict the fact that $(X,Y)$ is not induced
in $M$. Thus $i\in\{j-1,j\}$ and (i) holds.

Consider (ii). Say $z\in Z_1\cup Z_2\udots Z_{j-3}$. 
Assume that $z\in Y'$. Then $z\in \cl(P_{j-1}\cup P_j\udots P_n)$
as otherwise $(X,Y)$ is not induced in $M\ba b$. But,
by definition, $z\in\cl(P_0\cup P_1\udots P_{j-3})$.
Now, by Lemma~\ref{2path0}, 
$\sqcap(P_0\cup P_1\udots P_{j-3},P_{j-1}\cup P_j\udots P_n)=0$,
contradicting the fact that $b$ is not a loop of $M$.
Thus $z\in X'$. A similar argument proves that if
$z\in Z_{j+2}\cup Z_{j+3}\udots Z_n$, then $z\in Y'$.
\end{proof}

Let $t=\lfloor4(n-1)\rfloor$. For $i\in\{1,2,\ldots,t\}$,
let $b_{4i}$ be an element of $Z_{4i}$. Let 
$(X'_{4i},Y'_{4i})$ be a 2-separation of $M\ba b_{4i}$
that is induced by a $\PP$-relevant 2-separation 
$(X_{4i},Y_{4i})$ of $M\ba B$. Note that 
$(X'_{4i},\{b_{4i}\}\cup Y'_{4i})$ is a 3-separation of $M$.

\begin{sublemma}
\label{2path47-2}
If $1\leq i<j\leq t$, then $X'_{4i}\subseteq X'_{4j}$.
\end{sublemma}

\subproof
By \ref{2path47}(i),
$P_{4i+1}\cup P_{4i+2}\udots P_n\subseteq Y_{4i}\subseteq Y'_{4i}$.
By \ref{2path47}(ii),
$Z_{4i+2}\cup Z_{4i+3}\udots Z_n\subseteq Y'_{4i}$. Similarly
$P_1\cup P_2\udots P_{4j-2}\subseteq X'_{4j}$
and $Z_1\cup Z_2\udots Z_{4j-3}\subseteq X'_{4j}$.
But $j>i$, so $4j-3\geq 4i+1$. Thus
$Y'_{4i}\cup X'_{4j}=E(M)$ and hence
$X'_{4i}\subseteq X'_{4j}$.
\end{proof}

Let $R_1=X'_4$, $R_t=\{b_t\}\cup Y'_t$, and for 
$i\in\{2,3,\ldots,t-1\}$, let $R_i=X'_{4i}-X'_{4i-4}$.
Note that, $\lambda_M(R_1\cup R_2\udots R_i)=2$
for $i\in\{1,2,\ldots,t-1\}$. By this fact and \ref{2path47-2},
we deduce that 
$\RR=(R_1,R_2,\ldots,R_t)$ is a path of 
3-separations in $M$. Each internal step of $\RR$ contains
an internal step of $\PP$ and hence contains a
$\PP$-strong,
$M$-clonal pair. Such a clonal pair is clearly $\RR$-strong.

\begin{sublemma}
\label{2path47-3}
If $i,j\in\{1,2,\ldots,t\}$,
and $i<j-2$, then $R_1\cup R_2\udots R_i$
and $R_j\cup R_{j+1}\udots R_t$ are skew.
\end{sublemma}

\begin{proof}
Consider $X_{4i}$ and $Y_{4j-4}$. Then $4i<4j-4$,
so, by Lemma~\ref{2path0}, $X_{4i}$ and 
$Y_{4j-4}$ are skew in $M\ba B$.
Now $X'_{4i}\subseteq\cl(X_{4i})$ and
$Y'_{4j}\cup\{b_{4j}\}\subseteq \cl(Y_{4j-4})$.
Moreover $X'_{4i}=R_1\cup R_2\udots R_i$
and $Y'_{4j}\cup\{b_{4j}\}\subseteq \cl(Y_{4j-4})$.
This establishes the claim.
\end{proof}

It is a straightforward consequence of 
\ref{2path47-3} that $\RR$ does not display any
swirl-like flowers of order 4. The lemma now 
follows by letting 
$f_{\ref{2path47}}(q)=4(f_{\ref{at-last}}(4,q)+2)$.
\end{proof}

\section{Simply-Bridged Paths of $2$-separations}

Let $\PP=(P_0,P_1,\ldots, P_n)$ be a strict path of
2-separations of the connected matroid $N$ and let 
$M$ be a matroid with an $N$-minor. Then we say that
$M$ {\em simply bridges}
\index{simply bridged path of $2$-separations} 
$\PP$ if there is a labelling
$V=\{v_0,\ldots,v_n\}$ of $E(M)-E(N)$
and a sequence  
$N=M_0,M_1,\ldots,M_n=M$ of minors of 
$M$ such that the following hold.
\begin{itemize}
\item[(i)] For all $i\in \{0,1,\ldots,n-1\}$, the matroid
$M_{i+1}$ is either a single-element extension or coextension
of $M_i$ by $v_i$.
\item[(ii)] For all $i\in \{0,1,\ldots,n-1\}$,
the partition 
$(P_0\cup P_1\udots P_i\cup
\{v_0,v_1,\ldots,v_{i-1}\},P_{i+1}\cup P_{i+2}\udots P_n)$
is a 2-separation of $M_i$ that is bridged by $v_i$ in $M_{i+1}$.
\item[(iii)] Every $\PP$-relevant $2$-separation of
$N$ is bridged in $M$.
\end{itemize}

The goal of this section is to prove the following lemma.

\begin{lemma}
\label{simple-bridge}
Let $\PP=(P_0,P_1,\ldots, P_n)$ be a strict path of 
$2$-separations in the connected matroid $N$ that is simply
bridged by the matroid 
$M$. Assume that each internal step of $\PP$ contains a
$\PP$-strong $M$-clonal pair. Then there is a function
$f_{\ref{simple-bridge}}(q)$ such that, if 
$n\geq f_{\ref{simple-bridge}}(q)$, then $M\notin \eq$.
\end{lemma}

For the remainder of this section we assume that we are 
under the hypotheses
of Lemma~\ref{simple-bridge} 
with labelling for the sequence of minors and set $V$ of
bridging elements as given in the definition of ``simply bridges''.

For $i\in\{0,1,\ldots,n-1\}$, we say that $v_i$ is a {\em delete}
(respectively {\em contract}) element of $V$ if $M_{i+1}$
is an extension (respectively coextension) of $M_i$. We denote
the set of delete and contract elements by $C$ and $D$ respectively.

Say that $S=\{s_1,s_2,\ldots,s_\omega\}$ is a subset of $V$.
We use the notation $N[S]$
or $N[s_1,s_2,\ldots,s_\omega]$ for the matroid
$M/(C-S)\ba (D-S)$. Thus we have
$M_i=N[v_1,v_2,\ldots,v_{i-1}]$. For clarity we will
use the latter notation from now on. The next lemma follows 
easily from the definition of simply bridges.

\begin{lemma}
\label{ind-coind}
$C$ is independent in $M$ and $D$ is coindependent in $M$.
Hence, if $S$ is a subset of $V$, then $C\cap S$ is independent
in $N[S]$ and $D\cap S$ is coindependent in $N[S]$.
\end{lemma}

Define the function $\alpha$ on 
$\{0,1,\ldots,n-1\}$ as follows: $\alpha(i)=0$ if $v$ bridges
$P_0$ in $N[v_i]$. Otherwise 
$\alpha(i)$ is the least integer $j$ such that the
2-separating set
$$P_0\cup P_1\udots P_{j+1}\cup \{v_0,v_1,\ldots,v_j\}$$
of $N[v_0,v_1,\ldots,v_j]=M_{j+1}$ is bridged in 
$N[v_0,v_1,\ldots,v_j,v_i]$. It follows from the 
definition of simply bridged that $\alpha$
is well defined. Moreover, if $\alpha(i)>0$,
then $v_i$ does not bridge the 2-separating set
$$P_0\cup P_1\udots P_{\alpha(i)}\cup
\{v_0,v_1,\ldots,v_{\alpha(i)-1}\}$$
of $N[v_0,\ldots,v_{\alpha(i)-1}]$ in
$N[v_0,\ldots,v_{\alpha(i)-1},v_i]$.

For an integer $i$ we will denote the collection
of delete elements whose indices are less than $i$ by
$D_i^-$ and the set whose indices are greater than or equal to $i$
by $D_i^+$. The sets $C_i^-$ and $C_i^+$ are defined 
analogously. 

\begin{lemma}
\label{clean-bridge-1}
Say that $v_i$ is a delete element of $V$. Then
\begin{itemize}
\item[(i)] $v_i\in\cl_M(P_{\alpha(i)+1}\cup P_{\alpha(i)+2}
\udots P_n\cup C^+_{\alpha(i)})$.
\item[(ii)] $v_i\in\cl_M(P_0\cup P_1\udots P_{i+1}\cup C^-_{i+1})$.
\end{itemize}
\end{lemma}

\begin{proof}
As the 2-separation 
$(P_0\cup P_1\udots P_{\alpha(i)}\cup\{v_0,v_1,\ldots,v_{\alpha(i)-1}\},
P_{\alpha(i)+1}\cup P_{\alpha(i)+2}\udots P_n)$
of $N[v_0,v_1,\ldots,v_{\alpha(i)-1}]$ is induced in
$N[v_0,v_1,\ldots,v_{\alpha(i)-1},v_i]$, we see that
$v_i\in\cl_{N[v_0,v_1,\ldots,v_{\alpha(i)-1},v_i]}
(P_{\alpha(i)+1}\cup P_{\alpha(i)+2}\udots P_n)$. 
Part (i) follows from this fact. We omit the easy proof of (ii).
\end{proof}

Let $\QQ=(Q_0,Q_1,\ldots,Q_t)$ be a concatenation of $\PP$.
Then $\QQ$ is {\em tidy} if for all $i\in\{0,1,\ldots,n-1\}$ the
following holds: whenever $P_i\subseteq Q_j$, then 
$P_{\alpha(i)}\subseteq Q_{j-1}\cup Q_j$. We now consider the
situation when we can find a tidy concatenation of sufficient
length.

\begin{lemma}
\label{tidy}
Let $\QQ=(Q_0,Q_1,\ldots,Q_t)$ be a tidy concatenation of $\PP$.
Then there is a function $f_{\ref{tidy}}(q)$ such that 
if $t\geq f_{\ref{tidy}}(q)$, then $M\notin\eq$.
\end{lemma}

\begin{proof}
We first show that
if $t$ is sufficiently large, then
we can get $s$  large such that we 
have one of  two cases in $M/C$.

\begin{sublemma}
\label{subtidy1}
There is a function $f_{\ref{subtidy1}}(s)$ such that if 
$t\geq f_{\ref{subtidy1}}(s)$, then one the following holds.
\begin{itemize}
\item[(i)] There are indices $\beta$, $\gamma$ with 
$\gamma\geq s$ and a concatenation
$$\QQ'=(Q_0\udots Q_{\beta},
Q_{\beta+1},Q_{\beta+2},\ldots,Q_{\beta+\gamma-1},
Q_{\beta+\gamma}\udots Q_t)$$
of $\QQ$ such that every $\QQ'$-relevant $2$-separation of
$N$ is bridged in $M/C$.
\item[(ii)] There is a path
$\RR'=(R'_0,R'_1\ldots,R'_s)$ of $\QQ$-relevant
$2$-separations of $N$
such that,
\begin{itemize}
\item[(a)]  for all $i\in\{0,1,\ldots,s\}$,
there is a $j\in \{0,1,\ldots,t\}$ such that 
$Q_j\cup Q_{j+1}\subseteq R'_i$, and
\item[(b)] $\RR'$ induces a strict path
$\RR=(R_0,R_1,\ldots,R_s)$ of $2$-separations in $M/C$.
\end{itemize}
\end{itemize}
\end{sublemma}

\subproof
Let $f_{\ref{subtidy1}}(s)=s(4s+3)$ and assume that
$t\geq f_{\ref{subtidy1}}(s)$.
Assume that (i) does not hold. Then, in any concatenation
$$\QQ_k=(Q_0\udots Q_k,Q_{k+1},\ldots,Q_{s+k-1},Q_{s+k}\udots Q_t)$$
of $\QQ$, there is a $\QQ_k$-relevant 2-separation that is
induced in $M/C$. Thus there is a path
$\RR''=(R_0'',R_1'',\ldots,R''_{4s+3})$ of 
$\QQ$-relevant 2-separations such that
\begin{itemize}
\item[(1)] for all $i\in\{0,1,\ldots,4s+2\}$,
if 
$$R_0''\cup R_1''\udots R''_i=Q_0\cup Q_1\udots Q_{j-1}\cup Q'_j,$$
and 
$$R_0''\cup R_1''\udots R''_{i+1}=
Q_0\cup Q_1\udots Q_{k-1}\cup Q'_k,$$
then $j<k$; and
\item[(2)] every displayed $\RR''$-relevant 2-separation is
induced in $M/C$.
\end{itemize}
For $i\in\{0,1,\ldots,s\}$, let 
$R'_i=R''_{4i}\cup R''_{4i+1}\cup R''_{4i+2}\cup R''_{4i+3}$,
and let $\RR'=(R'_0,R'_1,\ldots,R'_s)$. 
Evidently each step of $\RR'$ contains 
two consecutive steps of $\QQ$ so that (a) holds.
Let
$R_0=\cl_{M/C}(R'_0)$ and, for $i\in\{0,1,\ldots,s\}$,
let $R_i=\cl_{M/C}(R'_0\cup R'_1\udots R'_i)-
\cl_{M/C}(R'_0\cup R'_1\udots R'_{i-1})$.
It follows from the definition of $\RR'$, 
Lemma~\ref{clean-bridge-1}, and Lemma~\ref{2-path-modular}, 
that $\RR$ is a path of 
2-separations in $M/C$.

There remains the irritating possibility that $\RR$ is not
a strict path of 2-separations. 
We now show that we may choose $\RR'$ and $\RR$ so that
$\RR$ is indeed strict. If $\RR$ is not strict,
then there exists a 
2-separation $(X,Y)$  
of $M/C$ that is not $\RR$ relevant such that neither 
$X$ not $Y$ is contained in a step of $\RR$. 
By Lemma~\ref{ind-coind}, 
$D$ is a coindependent set of $M/C$. 
If $X\subseteq D$, then $X$ is a parallel class and is contained
in a step of $\RR$ by the definition of $\RR$. Thus we may assume
that neither $X$ nor $Y$ is contained in $D$. 
Let $X'=X-D$ and $Y'=Y-D$. 
As $D$ is coindependent in $M/C$ and $N$ is connected,
we have  
$\cl_{M/C}(X')\supseteq X$, $cl_{M/C}(Y')\supseteq Y$,
and $(X',Y')$ is a 2-separation in $N$.
As $\QQ$ is a strict path in $N$, we may assume,
up to labels, that $X\subseteq Q_{\omega}$ for some
$\omega\in\{1,2,\ldots,t-1\}$. Then there is an 
$i\in\{1,2,\ldots,s\}$ such that 
$Q_{\omega}\subseteq R'_i\cup R'_{i+1}$. Moreover,
by uncrossing,
$\lambda_{M/C}(R'_0\cup R'_1\udots R'_i\cup X)=1$.
Consider the path
$$(R'_0,R'_1,\ldots,R'_{i-1},R'_i\cup X',R'_{i+1}-X',R'_{i+2},
\ldots, R'_s).$$
Observe that the 2-separation
$(R'_0\cup R'_1\udots R'_i\cup X',
(R'_{i+1}-X')\cup R'_{i+2}\udots R'_s)$
is induced in $M/C$. Repeat the process. As we are always
moving sets from a set of higher index to one of lower index,
the process must terminate. When it does, and we perform
an appropriate relabelling we have found the required 
paths $\RR'$ and $\RR$ satisfying the sublemma.
\end{proof}

Assume that we are in Case~(ii) of \ref{subtidy1}
and consider the path
$\RR=(R_0,R_1,\ldots,R_s)$. 
For $i\in\{1,\ldots,t\}$, let $\beta(i)$ denote the least positive
integer such that $Q_i\subseteq R_{\beta(i)}\cup R_{\beta(i)-1}$.
Evidently $\beta$ is well defined.

\begin{sublemma}
\label{subtidy2}
If $v_i\in D$, and $P_{i+1}\subseteq Q_j$, then 
$v_i\in R_{\beta(j)}\cup R_{\beta(j)-1}$.
\end{sublemma}

\subproof
From 
the definition of $\beta$ and the fact that each member
of $\RR$ contains a member of $\QQ$ we see that 
$Q_j\cup Q_{j-1}\subseteq R_{\beta(j)}\cup R_{\beta(j)-1}$.
By Lemma~\ref{clean-bridge-1} and the definition of $\QQ$,
we see that $v_i\in\cl_M(Q_0\cup Q_1\udots Q_j\cup C)$ and
$v_i\in \cl_M(Q_{j-1}\cup Q_j\udots Q_s\cup C)$.
Therefore $v_i\in\cl_{M/C}(Q_0\cup Q_1\udots Q_j)$
and $v_i\in\cl_{M/C}(Q_{j-1}\cup Q_j\udots Q_s)$.
By Lemma~\ref{2-path-modular} 
$v_i\in\cl_{M/C}(Q_{j-1}\cup Q_{j})$.
Hence $v_i\in\cl_{M/C}(R'_{\beta(j)}\cup R'_{\beta(j)-1})$.

We now prove that $v_i\in R_{\beta(j)}\cup R_{\beta(j)-1}$.
If this is not the case, then 
$v_i\in\cl_{M/C}(R'_0\cup R'_1\udots R'_{\beta(j)-2})$.
Observe that, by construction, 
$Q_{j-2}\cup Q_{j-1} \subseteq R_{\beta(j)}\cup R_{\beta(j)-1}$.
But, by the fact that $\QQ$ is a tidy concatenation of
$\PP$, we have 
$v_i\in\cl_{M/C}(Q_{j-1}\udots Q_{j})$.
This yields the contradiction that $v_i$ is a loop of $M/C$
and the sublemma follows.
\end{proof}

\begin{sublemma}
\label{subtidy3}
If $v_i\in C$ and $P_{i+1}\subseteq Q_j$, then
$v_i\in\cl^*_M(R_{\beta(j)}\cup R_{\beta(j)-1})$. 
\end{sublemma}

\subproof
If $P_{i+1}\subseteq Q_j$, then 
$P_{\alpha(i)}\subseteq Q_j\cup Q_{j-1}$.
By the dual of Lemma~\ref{clean-bridge-1}
$v_i\in\cl^*_M(Q_j\cup Q_{j+1}\udots  Q_t\cup D^+_{\alpha(i)})$.
By \ref{subtidy2},
$D^+_{\alpha(i)}\subseteq R_{\beta(j)-1}\cup R_{\beta(j)}\udots R_s$.
Hence $v_i\in\cl^*_M(R_{\beta(j)-1}\cup R_{\beta(j)}\udots R_s)$.

Again by the dual of Lemma~\ref{clean-bridge-1} we have
$v_i\in\cl^*_M(Q_1\cup Q_2\udots Q_j\cup D_i^-)$,
so that 
$v_i\in\cl^*_M(R_1\cup R_2\udots R_{\beta(j)})$.

This proves that 
$v_i$ does not coblock either 
$R_{\beta(j)-1}\cup R_{\beta(j)}\udots R_s$
or 
$R_1\cup R_2\udots R_{\beta(j)}$. By the dual of 
Lemma~\ref{2-path-modular} we conclude
that $v_i\in\cl^*_M(R_{\beta(j)}\cup R_{\beta(j)-1})$
as required.
\end{proof}

Assume that (i) holds in \ref{subtidy1}.
By Lemma~\ref{2path47}, if 
$s\geq f_{\ref{2path47}}(q)$,
then $M\notin \eq$. On the other hand if 
(ii) holds in \ref{subtidy1} then we 
may apply the dual of Lemma~\ref{2path47}
and conclude that if $s\geq f_{\ref{2path47}}(q)$,
then $M\notin \eq$. Thus the lemma holds by letting
$f_{\ref{tidy}}(q)=f_{\ref{subtidy1}}(f_{\ref{2path47}}(q))$.
\end{proof}

\begin{proof}[Proof of Lemma~\ref{simple-bridge}]
Consider the path $\PP=(P_1,P_2,\ldots,P_n)$
satisfying the hypotheses of the lemma.

\begin{sublemma}
\label{subsimple}
Assume that $M\in \eq$.
Then there is a function $f_{\ref{subsimple}}(m,l,q)$ such that,
if $n\geq f_{\ref{subsimple}}(m,l,q)$, then there exists
an integer $t(m)\in\{1,2,\ldots,n\}$, and a 
concatenation 
$\QQ=(Q_0,Q_1,\ldots,Q_{m-1},Q_m)$ of $\PP$, where
$Q_m=P_{t(m)}\cup P_{t(m)+1}\udots P_n$, such that the 
following hold.
\begin{itemize}
\item[(i)] $(P_1,P_2,\ldots,P_{t(m)-1},Q_m)$
is simply bridged by $\{v_1,\ldots,v_{t(m)}\}$
in $N[v_1,\ldots,v_{t(m)}]$.
\item[(ii)] For all $i$, if $P_i\subseteq Q_j$,
then $P_{\alpha(i)}\subseteq Q_{j-1}\cup Q_j$.
\item[(iii)] $Q_{m-1}$ contains at least $l$ steps of $\PP$.
\end{itemize}
\end{sublemma}

\subproof
Let $N=2f_{\ref{back-block}}(q)$. 

We define $f_{\ref{subsimple}}(m,l,q)$ inductively. To begin
let $f_{\ref{subsimple}}(3,l,q)=(l+1)N$. 
Assume that $n\geq (l+1)N$.
Let $S$ be the subset of $V$ with indices $i$ for which 
$\alpha(i)=0$. If $|S|\geq N$, then there are either at
least $f_{\ref{back-block}}(q)$ members of $D$ in $S$ or at
least $f_{\ref{back-block}}(q)$ members of $C$ in $S$.
In either case by Lemma~\ref{back-block}, or its dual, 
we contradict the assumption that $M\in \eq$.
Thus there are at most $N$ indices $i$ such that 
$\alpha_i=0$. Thus, if $n\geq N(l+1)$, there are indices $s$
and $t$ with $t-s\geq m$ such that, for all 
$j\in\{s,s+1,\ldots,t\}$ we have $\alpha_j\geq 1$.
Let $Q_1=P_1$, $Q_2=P_2\cup P_3\udots P_{s-1}$,
$Q_3=P_s\cup P_{s+1}\udots P_{t-1}$ and 
$Q_4=P_t\cup P_{t+1}\udots P_n$.
Clearly $(Q_0,Q_1,Q_2,Q_3)$ satisfies the claim with $m=3$.

For $m\geq 4$, let 
$f_{\ref{subsimple}}(m,l,q)=f_{\ref{subsimple}}(m-1,N(l+1),q)$.
A repeat of essentially the same argument as in the base case
establishes the claim.
\end{proof}

Let $m=f_{\ref{tidy}}(q)$ and let $f_{\ref{simple-bridge}}(q)=
f_{\ref{subsimple}}(m,1,q)$. 
If $M\in\eq$, then \ref{subsimple} implies that $\PP$ has a tidy
concatenation $\QQ$ of length $f_{\ref{tidy}}(q)$.
By Lemma~\ref{tidy}, $M\notin \eq$.
\end{proof}

\section{Sequentially Bridged Paths of $2$-separations}

Let $\PP=(P_0,P_1,\ldots,P_n)$ be a strict path of 
2-separations in  the connected matroid
$N$ and let $M$ be a matroid having $N$ as a minor. Then
$\PP$ is {\em sequentially bridged}
\index{sequentially bridged path of $2$-separations} 
in $M$ if 
$M$ has a sequence $M_0=N,M_1,\ldots,M_n$ of minors such that 
the following hold.
\begin{itemize}
\item[(i)] For all $i\in\{0,1,\ldots,n-1\}$, the matroid
$M_i$ is a minor of $M_{i+1}$ and $M_n$ is a minor of $M$. 
Denote $E(M_{i+1})-E(M_i)$ by $V_i$.
\item[(ii)] For all $i\in \{0,1,\ldots,n-1\}$, the 2-separation
$(P_0\cup P_1\udots P_{i-1},P_i\udots P_n)$ of $N$
is bridged in
$M_i$, and $\lambda_{M_i}(P_1\udots P_i
\cup V_0\udots V_{i-1})= 1$.
\item[(iii)] Every $\PP$-relevant $2$-separation of $N$ is 
bridged in $M$.
\end{itemize}

\begin{lemma}
\label{get-simple}
Let $N$ be a matroid with a strict path $\PP=(P_0,P_1,\ldots,P_n)$
of $2$-separations that is sequentially bridged in the matroid
$M$. Assume that each step of $N$ is an
$M$-clonal pair. Then there is a function $f_{\ref{get-simple}}(q)$
such that, if $n\geq f_{\ref{get-simple}}(q)$,
then $M\notin \eq$.
\end{lemma}

\begin{proof}
We use the notation of the definition of sequentially bridges.
Thus we have matroids 
$N=M_0,M_1,\ldots,M_n$ such that, for all $i\in\{0,1,\ldots,n-1\}$
the matroid $M_i$ is a minor of $M_{i+1}$ and $M_n$ is a minor of $M$. 
We also hve $V_i=E(M_{i+1})-E(M_i)$.
It is straightforwardly seen that we lose no generality in assuming
that for all $i\in\{0,1,\ldots,n-1\}$, the set
$V_i$ is a bridging sequence
for the 2-separation 
$(P_0\udots P_i\cup V_0\udots V_{i-1},P_{i+1}\udots P_n)$
of $M_i$.
For $i\in\{0,\ldots,n\}$, let $v_i$ denote the first element 
of $V_i$. Let $V=\{v_i:i\in\{0,\ldots,n-1\}\}$ and let
$C$ and $D$ denote the elements of $V$ that are 
respectively contract and delete
elements of their bridging sequences. Let 
$N'=M_n/C\ba D$. Observe that  
$$\PP'=(P_0,P_1\cup (V_0-\{v_0\}),\ldots,P_n\cup 
(V_{n-1}-\{v_{n-1}\}))$$
is a strict path of 2-separations in $N'$. Moreover, this path is
simply bridged in $M_n$. 
Evidently, if $i\in\{1,2,\ldots,n-1\}$, then $P_i$ is a 
$\PP'$-strong $M$-clonal pair.
The lemma now follows
by letting  $f_{\ref{get-simple}}(q)=f_{\ref{simple-bridge}}(q)$.
\end{proof}

\chapter{Last Rites}
\label{last-rites}

It will take only a few more lemmas and we will finally be
able to prove the main theorems of the paper.

\section{Sequentially Bridged Paths of $3$-separations}

We have defined ``sequentially bridges'' for paths of 2-separations.
The extension to paths of 3-separations has no surprises.

Let $\PP=(P_0,P_1,\ldots,P_n)$ be a path of 
3-separations in the matroid $M$.
A $3$-separation $(X,Y)$ of $M$ is {\em $\PP$-relevant}
\index{$\PP$-relevant $3$-separation}
if there is a $j\in\{1,2,\ldots,n-1\}$ such that either
$X$ or $Y$ has the form $P_0\cup P_1\udots P_{j-1}\cup P'_j$
for some subset $P'_j$ of $P_j$.

Let $\PP=(P_0,P_1,\ldots,P_n)$ be a  path of 
3-separations in  the connected matroid
$N$ and let $M$ be a matroid having $N$ as a minor. Then
$\PP$ is {\em sequentially bridged}
\index{sequentially bridged path of $3$-separations} 
in $M$ if 
$M$ has a sequence $M_0=N,M_1,\ldots,M_n$ of minors such that 
the following hold.
\begin{itemize}
\item[(i)] For all $i\in\{0,1,\ldots,n-1\}$, the matroid
$M_i$ is a minor of $M_{i+1}$ and $M_n$ is a minor of $M$. 
Denote $E(M_{i+1})-E(M_i)$ by $V_i$.
\item[(ii)] For all $i\in \{0,1,\ldots,n-1\}$, the 3-separation
$(P_0\cup P_1\udots P_{i-1},P_i\udots P_n)$ of $N$
is bridged in
$M_i$, and $\lambda_{M_i}(P_1\udots P_i
\cup V_0\udots V_{i-1})\leq 1$.
\item[(iii)] Every $\PP$-relevant $3$-separation of $N$ is 
bridged in $M$.
\end{itemize}

\begin{lemma}
\label{getting-closer}
Let $M$ be a matroid  with a connected minor $N$ that
has a swirl-like flower $\PP=(P_0,P_1,\ldots,P_n)$
such that each petal of $\PP$, apart from $P_0$ and $P_n$
are $M$-clonal pairs, and such that the path $\PP$
is sequentially bridged in $M$. Then there is a function
$f_{\ref{getting-closer}}(q)$ such that, if
$n\geq f_{\ref{getting-closer}}(q)$, then $M\notin\eq$.
\end{lemma}

\begin{proof}
As $\PP$ is sequentially bridged, there is a sequence of 
minors $N=M_0,M_1,\ldots,M_n$ of $M$ satisfying the definition 
of sequentially bridged. Let $V_i=E(M_i)-E(M_{i-1})$.
For $i\in\{0,1,\ldots,n-1\}$ let
$A_i=V_0\cup V_2\udots V_i$. Observe that, as both $P_{n-1}$ and 
$P_{n-3}$ are clonal pairs, then $P_{n-2}$ is coclosed
in $N$, so that $N\ba P_{n-2}$ is connected.
Let $P'_0=P_0\cup P_n\cup P_{n-1}$,
and let $\PP'=(P'_0,P_1,\ldots,P_{n-3})$.
By Lemma~\ref{get-2-path} $\PP'$ is a strict path of 
$2$-separations in $N\ba P_{n-2}$.

\begin{sublemma}
$\PP'$ is sequentially bridged in $M$.
\end{sublemma}

\subproof
Say $i\in\{1,2,\ldots,n-4\}$. Then 
$\lambda_{M_i}(P_0\cup P_1\udots P_i\cup A_{i-1})=2$
and it follows that 
$(P_0\cup P_1\cup P_i\cup A_{i-1},P_{i+1},\ldots,P_n)$
is a swirl-like flower in $N_i$. Hence 
$(P'_0\cup P_1\cup P_i\cup A_{i-1},P_{i+1},\ldots,P_{n-3},P_{n-2})$
is a swirl-like flower in $M_i$ so that
$(P'_0\cup P_1\cup P_i\cup A_{i-1},P_{i+1},\ldots,P_{n-3})$
is a path of 2-separations in $M_i\ba P_{i-2}$.
Thus $\lambda_{M_i\ba P_{i-2}}(P'_0\cup P_1\udots P_i\cup A_{i-1})=1$.

We now prove that 
the 2-separation $(P'_0\cup P_1\udots P_{i-1},P_i\udots P_{n-3})$ 
of $\PP'$ is bridged in $M_i\ba P_{i-2}$. Assume otherwise.
Then there is a partition $(A',A'')$ of $A_{i-1}$ such that
$$\lambda_{M_i\ba P_{n-2}}(P'_0\cup P_1\udots P_{i-1}\cup A',
P_i\udots P_{n-3}\cup A'')=1.$$
But it is easily seen that 
$$r_{M_i}(P_i\udots P_{n-3}\cup A''\cup P_{n-2})
\leq r_{M_i}(P_i\udots P_{n-3}\cup A'')+1.$$
Thus $\lambda_{M_i}(P'_0\cup P_1\udots P_{i-1}\cup A')=2$.
By definition $\lambda_{N_i}(P_1\udots P_i\cup A_{i-1})=2$.
Uncrossing these two 3-separating sets proves that
$\lambda_{M_i}(P_1\udots P_{i-1}\cup A')=2$. Hence
$P_1\udots P_{i-1}$ is not bridged in $M_i$, contradicting 
the definition of sequentially bridges.

The claim follows from the above facts and the fact that
every $\PP'$ relevant $2$-separation is displayed by $\PP'$.
\end{proof}

Let $f_{\ref{getting-closer}}(q)=f_{\ref{get-simple}}(q)+3$.
The lemma now follows from Lemma~\ref{get-simple}.
\end{proof}

\section{Bridged Swirls}

The goal of this section is to prove the following theorem.

\begin{theorem}
\label{trapped}
Let $M$ be a $4$-connected matroid in $E(q)$
and let $N$ be a free-swirl minor of $M$ all of whose
clonal pairs are clonal in $M$. Then there is a function
$f_{\ref{trapped}}(q)$ such that 
$|E(N)|\leq f_{\ref{trapped}}(q)$.
\end{theorem}

Given Lemma~\ref{getting-closer}, it is clear that 
the task is to prove the following lemma.

\begin{lemma}
\label{gotcha}
Let $M$ be a $4$-connected matroid in $\eq$ with a 
$\Delta_n$-minor $N$ whose ground set consists of $M$-clonal
pairs. Then there is a function $f_{\ref{gotcha}}(m,q)$
such that, if $n \geq f_{\ref{gotcha}}(m,q)$,
then $M$ has a minor with a 
swirl-like flower $\QQ=(Q_0,Q_1,\ldots,Q_m)$,
each petal of which apart from $Q_0$ is an
$M$-clonal pair, having the property that 
the path $\QQ$ is sequentially bridged in $M$. 
\end{lemma}

In fact this task is quite straightforward, although notationally
unwieldy. In the light of Lemma~\ref{2path2}, there are no
surprises in the next lemma.

\begin{lemma}
\label{1-block-swirl}
Let $M$ be a matroid in $\eq$ with an element $b$
such that $M\ba e$ is connected with
a swirl-like flower
$\PP=(P_0,P_2,\ldots,P_n)$. Assume that
all petals of $\PP$, apart from $P_0$, are $M$-clonal pairs
and that $P_0$ is blocked by $b$.
Then there is a function $f_{\ref{1-block-swirl}}(m,q)$
such that, if $n\geq f_{\ref{1-block-swirl}}(m,q)$,
then the following holds. There exist 
$i$ and $j$ in $\{1,2,\ldots,l\}$ with $j-i\geq m$
such that
$$(P_0\cup P_1\udots P_{i-1}\cup P_{j+1}\udots P_n\cup \{b\},
P_i,P_{i+1},\ldots,P_j)$$
is a swirl-like flower in $M$.
\end{lemma}

\begin{proof}
Let $f_{\ref{1-block-swirl}}(m,q)=f_{\ref{2path2}}(m,q)+m+2$.
Observe that $P_{n-1}$ is coclosed in $M\ba b$,
so by Lemma~\ref{get-2-path}, 
$(P_0\cup P_n,P_1,\ldots,P_{n-2})$ is a strict path of
2-separations in $M\ba b\ba P_{n-1}$. 
Assume that $P_{n-m-2}\cup P_{n-m-1}\udots P_{n-2}$
is not blocked by $b$ in $M\ba P_{n-1}$.
Then the partition 
$(\{b\}\cup P_0\cup P_1\udots P_{n-m-3}\cup P_{n-1}\cup P_n,
P_{n-m-2},P_{n-m-1},\ldots,P_{n-2})$ of $E(M)$
clearly satisfies the lemma.
Assume then, that  
$P_{n-m-2}\cup P_{n-m-1}\udots P_{n-2}$
is blocked by $b$ in $M\ba P_{n-1}$.
Consider the path
$(P_n\cup P_0,P_1,\ldots,P_{n-m-2}\cup P_{n-m-1}\udots P_{n-2})$
of $2$-separations in $M\ba b\ba P_2$. This path satisfies the
hypotheses of Lemma~\ref{2path2}. Applying that lemma gives
a swirl-like flower in $M\ba P_2$
of the form
$$(P_0\cup P_1\udots P_{i-1}\cup P_{j+1}\udots P_{n-2}\cup 
P_n\cup \{b\},
P_i,P_{i+1},\ldots,P_j),$$
where $j-i\geq m$.
To complete the proof of the lemma it suffices to observe that
$P_{n-1}\subseteq \cl(P_{n-2}\cup P_n)$.
\end{proof}

Extending from a single blocking element to a bridging sequence
we obtain:

\begin{lemma}
\label{gotcha1}
Let $M$ be a matroid in $\eq$ with a minor $N$ that is connected
and has a swirl-like flower
$\PP=(P_0,P_1,\ldots, P_n)$. Assume that
all petals of $\PP$ except $P_0$ are $M$-clonal pairs. Let
$V_0$ be an $M$-bridging sequence for $P_0$.
Then there is a function $f_{\ref{gotcha1}}(m,q)$
such that if $n\geq f_{\ref{gotcha1}}(m,q)$, then the following
holds. There exist $i$ and $j$ in $\{1,2,\ldots,l\}$ with $j-i\geq m$
such that 
$$(P_0\cup P_1\udots P_{i-1}\cup P_{j+1}\udots P_n\cup V_0,
P_i,P_{i+1},\ldots,P_j)$$
is a swirl-like flower in $N[V_0]$.
\end{lemma}

\begin{proof}
Say that bridging sequence is $(v_0,v_1,\ldots,v_t)$.
Up to duality we may assume that $v_t$ is a deletion element.
Observe that 
$(P_0\cup \{v_0,v_1,\ldots,v_{t-1}\},P_1,P_2,\ldots,P_n)$
is a swirl-like flower in $N[V_0]\ba v_t$ all of whose
petals are $N[V_0]$-clonal pairs apart from the initial
petal. Moreover, the initial petal is blocked by $v_t$.
The lemma follows from Lemma~\ref{1-block-swirl} by letting 
$f_{\ref{gotcha1}}(m,q)=f_{\ref{1-block-swirl}}(m,q)$.
\end{proof}

To facilitate the proof of Lemma~\ref{gotcha} we introduce some
terminology.
Let $\PP=(P_0,P_1,\ldots,P_n)$ be a swirl-like flower in a matroid 
$N$ all of whose petals except $P_0$ consist of clonal pairs. 
If $i\in\{1,2,\ldots,n\}$, and $P_i=\{p_i,q_i\}$,
then $\PP=(P_0,\ldots,P_{i-1},P_{i+1},\ldots,P_n)$ is a
swirl-like flower in $N'=N\ba p_i/q_i=N/p_i\ba q_i$. In this
case we say that
$N'$ and $\PP'$ are obtained by {\em removal of $P_i$}. More
generally, if $N'$ and $\PP'$ are
 obtained by a sequence of such operations,
then we say that they  are obtained from $M$ and $\PP$ by 
{\em petal removal}.
\index{petal removal}

The next lemma, while lengthy to state, is an 
immediate corollary of Lemma~\ref{last-bridge}.

\begin{lemma}
\label{gotcha2}
Let $\PP=(P_0,P_1,\ldots, P_n)$ be a swirl-like flower in the matroid
$N$. Assume that
all petals of $\PP$ except $P_0$ are clonal pairs. Let
$V_0$ be an $M$-bridging sequence for $P_0$. Assume that
$1\leq i<j\leq n$ and that 
$$(P_0\cup P_1\udots P_{i-1}\cup P_{j+1}\udots P_n\cup V_0,
P_i,P_{i+1},\ldots,P_j)$$
is a maximal 
swirl-like flower in $N[V_0]$. Then either $i>1$ or $j<n$.

If $i>1$ and $N/C\ba D$ is the matroid obtained
by removing the petals $P_1,P_2,\ldots,P_{i-2},P_{j+1},P_{j+2}
\ldots,P_n$
from $\PP$  
then the following hold.
\begin{itemize}
\item[(i)] The  petal $P_0$ of  the swirl-like flower
$(P_0,P_{i-1},P_i,\ldots,P_{j-1},P_j)$ of the matroid $N/C\ba D$
is bridged in $N[V_0]/C\ba D$.
\item[(ii)] $\lambda_{N[V_0]/C\ba D}(P_i\cup P_{i+1}
\udots  P_j)=2$.
\end{itemize}
On the other hand, if $j<n$ and $N'$ is the matroid obtained 
by removing the petals $P_1,\ldots,P_{i-1},P_{j+2},\ldots, P_n$
from $\PP$ (say $N'=N/C'\ba D'$), then the following hold.
\begin{itemize}
\item[(iii)] The petal $P_0$ of  the swirl-like flower
$(P_0,P_i,\ldots,P_j,P_{j+1})$ of the matroid $N/C'\ba D'$
is bridged in $N[V_0]/C'\ba D'$.
\item[(iv)] $\lambda_{N[V_0]/C'\ba D'}(P_i\cup P_{i+1}
\udots  P_j)=2$.
\end{itemize}
\end{lemma}

In the next lemma we use the notation $N\less M$ 
to indicate that  $N$ is a minor of $M$.
The lemma is somewhat stronger than we need, but it is 
set up to facilitate an inductive proof. Note that $n+1=0\mod n+1$.

\begin{lemma}
\label{gotcha3}
Let $M$ be a $4$-connected matroid in $\eq$ and let
$N$ be a minor of $M$ with a swirl-like flower
$\PP=(P_0,P_1,\ldots,P_n)$ all of whose petals apart from
$P_0$ are $M$-clonal pairs.
Then there is a function $f_{\ref{gotcha3}}(m,t,q)$
such that, if $n\geq f_{\ref{gotcha3}}(m,t,q)$, then
$N$ has a minor $N'$ with a swirl-like flower
$(P_0,Q_1,Q_2,\ldots,Q_m)$ obtained by removing petals from
$\PP$ such that there are indices
$$0\leq l_1\leq l_2\leq \cdots \leq l_m <r_m\leq r_{m-1}
\leq \cdots\leq r_1\leq n+1$$
and minors $N'=N_0\less N_1\less \cdots \less N_m\less M$,
where $E_i=E(N_i)-E(N_0)$, 
such that the following hold.
\begin{itemize}
\item[(i)] $r_m-l_m\geq t$.
\item[(ii)] If $i\in\{2,3,\ldots,m\}$, then 
$(l_i-l_{i-1})+(r_{i-1}-r_i)=1$, and
$(l_1,r_1)\in\{(0,n),(1,0)\}$.
\item[(iii)] $\lambda_{N_i}(P_0\cup P_{l_1}\udots P_{l_{i}}
\cup P_{r_1}\udots P_{r_{i}}\cup E_{i})=2$
for $i\in\{1,2,\ldots,m\}$.
\item[(iv)] $P_0\cup P_{l_1}\udots P_{l_{i-1}}\cup
P_{r_1}\udots P_{r_{i-1}}\cup E_{i-1}$ is bridged in $N_{i}$
for $i\in\{1,2,\ldots,m\}$.
\end{itemize}
\end{lemma}

\begin{proof}
Define $f_{\ref{gotcha3}}(m,t,q)$, by
$f_{\ref{gotcha3}}(1,t,q)=f_{\ref{gotcha1}}(t,q)$ and,
for $m>1$ by
$f_{\ref{gotcha3}}(m,t,q)=
f_{\ref{gotcha3}}(m-1,f_{\ref{gotcha1}}(t,q),q)$.

It follows from Lemmas~\ref{gotcha1} and \ref{gotcha2} that
the lemma holds for $m=1$ 
with $N_1=N[V_0]$ as defined in Lemma~\ref{gotcha2}.
Say $k>1$ and assume that the lemma holds for all $t$ whenever
$m\leq k$. Let $t'=f_{\ref{gotcha1}}(t,q)$ and say that
$n\geq f_{\ref{gotcha3}}(k,t',q)$. Let $N_0$ be a minor of 
$M$ with a swirl-like flower $(P_1,Q_1,\ldots,Q_n)$
for which there exists a sequence of minors 
$N_0\less N_1\less\cdots \less N_k\less M$ with indices
$l_1,l_2,\ldots,l_k,r_1,r_2,\ldots,r_k$ 
that satisfy the conclusion of the
lemma where $r_k-l_k\geq t'$.
Let
$$Q_0=P_0\cup Q_1\udots Q_{l_k}\cup P_{r_k}\cup R_{r_k+1}
\udots P_n \cup (E(N_k)-E(N_0)),$$
and let
$$\QQ=(Q_0,Q_{l_k+1},Q_{l_k+2},\ldots,Q_{r_k-1}).$$
Then $\QQ$ is a swirl-like flower in $N_k$, each petal of which,
apart from $Q_0$ is an $M$-clonal pair. Therefore, by
Lemmas~\ref{gotcha1}, \ref{gotcha2}, and the choice of $t'$,
there is a minor $N_k/C\ba D$ of $N_k$, 
obtained by removing petals from $\QQ$, with a swirl-like flower
$\QQ'=(Q_0,Q'_1,\ldots,Q'_u)$, where $u\geq t+1$, for which the 
following hold.
\begin{itemize}
\item[(i)] $N_k/C\ba D\less N_{k+1}\less M$.
\item[(ii)] $Q_0$ is bridged in $N_{k+1}$.
\item[(iii)] Either 
$\lambda_{N_{k+1}}(Q_1'\cup Q'_2\udots Q'_{u-1})=2$
or $\lambda_{N_{k+1}}(Q'_2\cup Q'_3\udots Q'_u)=2$.
\end{itemize}
Consider the sequence
$$N_0/C\ba D,N_1/C\ba D,\ldots,N_k/C\ba D, N_{k+1},$$
and the swirl-like flower in $N_0/C\ba D$ obtained by removing the 
petals of $(P_0,Q_1,\ldots,Q_n)$ that are contained in $C\cup D$.
It is now easily verified that the conclusions of the lemma hold
for this choice of flower and sequence of minors.
\end{proof}

\begin{proof}[Proof of Lemma~\ref{gotcha}]
Let $f_{\ref{gotcha}}(m,q)=f_{\ref{gotcha3}}(2m,1,q)$.
By Lemma~\ref{gotcha3} we have the following: a sequence
of minors
$$N_0\less N_1\less\cdots\less N_{2m}\less M,$$
where $E_i=E(N_i)-E(N_0)$;
a swirl-like flower $\PP=(P_0,P_1,\ldots,P_n)$ in $N_0$,
all of whose petals apart from $P_0$ are $M$-clonal pairs;
and a sequence of indices
$$0\leq l_1\leq l_2\leq \cdots\leq l_{2m}
<r_{2m}\leq \cdots \leq r_1\leq 2m+1$$
such that the following hold.
\begin{itemize}
\item[(i)] If $i\in\{2,3,\ldots,2m\}$, then 
$(l_i-l_{i-1})+(r_{i-1}-r_i)=1$, and
$(l_1,r_1)\in\{(0,n),(1,0)\}$.
\item[(ii)] $\lambda_{N_i}(P_0\cup P_{l_1}\udots P_{l_{i}}
\cup P_{r_1}\udots P_{r_{i}}\cup E_{i})=2$
for $i\in\{1,2,\ldots,2m\}$.
\item[(ii)] $P_0\cup P_{l_1}\udots P_{l_{i-1}}\cup
P_{r_1}\udots P_{r_{i-1}}\cup E_{i-1}$ is bridged in $N_{i}$
for $i\in\{1,2,\ldots,2m\}$.
\end{itemize}

Let ${\mathcal N}=\{N_0,N_1,\ldots,N_{2m}\}$
For $i\in\{0,1\ldots,2m-1\}$ say that $N_i\in{\mathcal N}$ 
is {\em bridged on the left} if $l_{i+1}=l_i+1$; otherwise it
is {\em bridged on the right}. Let $l$ denoted the number of
matroids in $\mathcal N$ that are bridged on the left.
We lose no generality in assuming that $l\geq m$. Let
$N'_0,N'_1,\ldots,N'_l$ denoted the matroids that are bridged
on the left, where 
$N'_0\less N'_1\less\cdots\less N'_l$. Let 
$E'_i=E(N'_i)-E(N'_{i-1})$. The next
claim is an immediate 
consequence of the definitions.

\begin{sublemma}
\label{subgotcha1}
Say $i\in\{1,2,\ldots,l\}$. Then there is a sequence
$n>j_1\geq j_2\geq \cdots \geq j_l\geq l+2$
such that
\begin{itemize}
\item[(i)] $\lambda_{N'_i}(P_0\cup P_1\udots P_i\cup P_n\cup
P_{n-1}\udots P_{j_i}\cup E'_i)=2$,
but 
\item[(ii)] $P_0\cup P_1\udots P_i\cup P_n\cup
P_{n-1}\udots P_{j_{i-1}}\cup E'_{i-1}$ is bridged
in $N'_i$.
\end{itemize}
\end{sublemma}

Let 
$P'_0=P_{l+2}\cup P_{l+3}\udots P_n\cup P_0\cup (E(N'_0)-E(N_0))$
An easy uncrossing argument proves

\begin{sublemma}
\label{subgotcha2}
Say $i\in\{0,1,\ldots,l\}$. Then
\begin{itemize}
\item[(i)] $\lambda_{N'_i}(P'_0\cup P_1\udots P_i\cup E'_i)=2$,
but 
\item[(ii)] $P'_0\cup P_1\udots P_{i-1}\cup E'_{i-1}$ 
is bridged in $N'_{i}$.
\end{itemize}
\end{sublemma}

Let
$\PP'=(P'_0,P_1,P_2,\ldots,P_{l+1})$. Then $\PP'$ is a swirl-like
flower in $N'_0$ all of whose petals, apart from $P'_0$, are
$M$-clonal pairs. It follows from \ref{subgotcha2} that
the path $\PP'$ is sequentially bridged by the matroids
$N'_0,N'_1,\ldots,N'_l$. The lemma now follows from the fact that
$l\geq m$.
\end{proof}

\begin{proof}[Proof of Theorem~\ref{trapped}]
Let $f_{\ref{trapped}}(q)
=f_{\ref{gotcha}}(f_{\ref{getting-closer}}(q)+1,q))$.
Assume that $M$ has a $\Delta_n$-minor all
minor of $M$ all of whose petals are $M$-clonal pairs,
where $n\geq f_{\ref{trapped}}(q)$.
Then, by Lemma~\ref{gotcha}, for some 
$m\geq f_{\ref{getting-closer}}(q)+1$, the matroid
$M$ has a minor 
with a swirl-like flower $\QQ=(Q_0,Q_1,\ldots,Q_m)$
all of whose petals, apart from $Q_0$, are $M$-clonal
pairs with the property that the path $\QQ$ is 
sequentially bridged in $M$. By Lemma~\ref{getting-closer}
$M\notin \eq$.
\end{proof}

\section{The Interment}

At long last we are able to bury our skeletons and prove that there are
only a finite number of $k$-skeletons in $\eq$.

\begin{theorem}
\label{control-skeleton}
Let $q$ and $k\geq 5$ be  integers. Then there is a function 
$f_{\ref{control-skeleton}}(k,q)$ such that, if 
$M$ is a $k$-skeleton in $\eq$, then $M$ has at most 
$f_{\ref{control-skeleton}}(k,q)$ elements.
\end{theorem}

\begin{proof}
Let $n_1=f_{\ref{trapped}}(q)+1$, let $n_2=f_{\ref{klonal1}}(n_1,q)$,
and let $f_{\ref{control-skeleton}}(k,q)=f_{\ref{tame}}(n_2,k,q)$.

Let $M$ be a $k$-skeleton and 
assume that $|E(M)|> f_{\ref{control-skeleton}}(k,q)$.
By Theorem~\ref{tame}, $M$ has a 4-connected minor $N$ with 
a set of $n_2$ pairwise disjoint clonal pairs. By
Theorem~\ref{klonal1}, $N$ has an
$N$-clonal $\Delta_{n_1}$-minor. By Theorem~\ref{trapped},
$M\notin\eq$.
\end{proof}

\chapter{Applications to Matroid Representability}
\label{fields-at-last}

We are finally in a position to obtain consequences for 
inequivalent representations of matroids.
We begin with results that bound the number of inequivalent
representations. Two representations
of a matroid over a field $\mathbb F$ are 
{\em equivalent}
\index{equivalent representations}
if one can be obtained from another by elementary row
operations and column scalings. 
This differs from the definition given in \cite{ox92}
where field automorphisms are also allowed.
For the results presented here the difference is not
significant---for a finite field, a bound on the number
of inequivalent representations with respect to one notion
of equivalence implies a bound with respect to the other.

\section{Bounding Inequivalent Representations}

Our main theorems bounding inequivalent representations 
are corollaries of the next theorem.

\begin{theorem}
\label{master}
Let $k\geq 5$ and $q\geq 3$ be integers and let 
$\mathbb F$ be a finite field. Then there is a function
$f_{\ref{master}}(k,q,{\mathbb F})$ such that
a $k$-coherent member of
$\eq$ has at most $f_{\ref{master}}(k,q,{\mathbb F})$
inequivalent ${\mathbb F}$-representations.
\end{theorem}

To prove Theorem~\ref{master}, we need a few more 
easy facts.
Let $\mathbb F$ be a field. 
Suppose that $z$ is fixed in $M$, and consider two 
${\mathbb F}$-representations of $M$ of the form $[A,y]$ and 
$[A,y']$, where $A$ represents $M\ba z$. Now
$[A,y,y']$ represents a single-element extension of 
$M$ and it is easily checked that $y$ and $y'$ are clones.
Since $z$ is fixed, $\{y,y'\}$ is a parallel pair.
Thus $[A,y]$ and $[A,y']$ are equivalent.
This shows that, up to equivalence, any representation
of $M\ba z$ extends to at most one representation of $M$.
Part (i) of the next lemma follows from this argument.
Part (ii) is the dual of part (i).

\begin{lemma}
\label{extend-fixed}
Let $z$ be an element of the matroid $M$ and 
${\mathbb F}$ be a finite field. 
\begin{itemize}
\item[(i)] If $z$ is fixed in $M$, then the number of 
inequivalent ${\mathbb F}$-representations of 
$M$ is at most the number of inequivalent 
${\mathbb F}$-representations of $M\ba z$.
\item[(ii)] If $z$ is cofixed in $M$, then the number of 
inequivalent ${\mathbb F}$-representations of 
$M$ is at most the number of inequivalent 
${\mathbb F}$-representations of $M/z$.
\end{itemize}
\end{lemma}

Recall
that wheels have one $k$-skeleton minor, 
namely $U_{2,3}$ and whirls have 
an additional one, namely $U_{2,4}$.

\begin{lemma}
\label{field1}
Let $M$ be a $k$-coherent matroid 
and let $\mathbb F$
be a finite field. Then the number of inequivalent 
$\mathbb F$-representations
of $M$ is bounded above by the maximum of the 
number of inequivalent $\mathbb F$-representations
of members of the set of $k$-skeleton minors of $M$.
\end{lemma}

\begin{proof}
Assume that $M$ is not a $k$-skeleton. If $M$ is a wheel, then
$M$ is uniquely $\mathbb F$-representable, as is $U_{2,3}$.
If $M$ is a whirl,
then it is well known and easily seen that the number of inequivalent 
$\mathbb F$-representations of $M$ is equal to that of $U_{2,4}$.
So the lemma holds in these trivial cases. Assume that $M$ is not a 
wheel or a whirl. Then, up to duality, there is an element
$x\in E(M)$ such that $x$ is fixed in $M$ and $M\ba x$ is 
$k$-coherent and the lemma holds by Lemma~\ref{extend-fixed}
and an obvious induction.
\end{proof}

\begin{proof}[Proof of Theorem~\ref{master}] 
By Theorem~\ref{control-skeleton}
there are a finite number of $k$-skeletons in 
$\eq$. Let $f_{\ref{master}}(k,q,{\mathbb F})$
denote the maximum of the number of inequivalent
${\mathbb F}$-representations of a $k$-skeleton
in $\eq$. It follows from Lemma~\ref{field1}
that a $k$-coherent matroid in $\eq$
has at most $f_{\ref{master}}(k,q,{\mathbb F})$
as required.
\end{proof}

It follows from \cite[Lemma~11.6]{totally-free} that, if 
$p$ is a prime that exceeds $3$, then $\Lambda_p$
is not $GF(p)$ representable. Certainly neither
$U_{2,p+2}$ nor $U_{p,p+2}$ is $GF(p)$-representable.
We therefore have

\begin{lemma}
\label{exclude-spikes}
Let $p\geq 3$ be a prime. Then the class of $GF(p)$-representable
matroids is contained in ${\mathcal E}(p)$.
\end{lemma}

The next corollary follows immediately from 
Theorem~\ref{master} and Lemma~\ref{exclude-spikes}.

\begin{corollary}
\label{biggie2}
Let $k\geq 5$ be an integer, $p$ be a prime number
and $\mathbb F$ be a finite field. Then a 
$k$-coherent $GF(p)$-representable matroid
has at most 
$f_{\ref{master}}(k,p,{\mathbb F})$ inequivalent representations over 
$\mathbb F$.
\end{corollary}

A special case of Corollary~\ref{biggie2} is Theorem~\ref{biggie}.
We restate it here for convenience.

\begin{theorem}
\label{biggy}
Let $k\geq 5$ be an integer and $p$ be a prime number. Then there
is a function $f_{\ref{biggy}}(k,p)$ such that a $k$-coherent matroid
has at most $f_{\ref{biggy}}(k,p)$ 
inequivalent representations over $GF(p)$.
\end{theorem}

Finally Theorem~\ref{4-conn}, stated in the introduction,
follows from Theorem~\ref{biggie} as 4-connected matroids 
are $k$-coherent.

\section{Excluding a Free Swirl}

If, as well as excluding a free spike we exclude a free swirl,
we bound the number of inequivalent representations of 
3-connected matroids over a finite field. 

A 3-connected matroid $M$ is {\em totally free}
\index{totally free matroid} 
if it is not 
a wheel or a whirl of rank at least three, 
and has the properties that,
for all $x\in E(M)$, if $M\ba x$
is 3-connected, then $x$ is not fixed in $M$ and if $M/x$
is 3-connected, then $x$ is not cofixed in $M$. 
If, in addition to  excluding $U_{2,q+2}$, $U_{q,q+2}$ and
$\Lambda_q$, we also exclude $\Delta_q$, 
then, as well as having only
a finite number of $k$-skeletons, we have only a finite
number of totally-free matroids. It is proved in
\cite{totally-free}, and easily seen, that the 
maximum number of inequivalent representations of a
3-connected
matroid over a finite field is bounded above by that
of its totally free minors. Let 
$E(U_{2,q+2},U_{q,q+2},\Lambda_q,
\Delta_q)$, denote the class of matroids with no
$U_{2,q+2}$-, $U_{q,q+2}$-, $\Lambda_q$- or $\Delta_q$-minor.

\begin{theorem}
\label{totally-free}
Let $q$ be a positive integer. Then there are a finite number of 
totally free matroids in $E(U_{2,q+2},U_{q,q+2},\Lambda_q,
\Delta_q)$. 
\end{theorem}

\begin{proof}
Let $M$ be a totally free matroid. Observe that, if $M$
is $k$-coherent, then $M$ is also a $k$-skeleton.
By Theorem~\ref{control-skeleton}, there are a finite number
of $k$-skeletons in $\eq$ for any fixed $k$. Thus, if the theorem
fails, there is a totally free matroid $M\in
E(U_{2,q+2},U_{q,q+2},\Lambda_q,\Delta_q)$ that is not
$q$-coherent. Let $l=|E(M)|$. Then any 3-connected
minor of $M$ is $l$-coherent. Hence $M$ is an 
$l$-skeleton. As $M$ is not $q$-coherent, $M$ has a swirl-like
flower of order $q$. By Corollary~\ref{skeleton-free-swirl},
$M$ has a $\Delta_q$-minor, contradicting the assumption
that $M\in E(U_{2,q+2},U_{q,q+2},\Lambda_q,
\Delta_q)$.
\end{proof}

The next corollary follows from Theorem~\ref{totally-free}
and the observation prior to it.

\begin{corollary}
\label{3-con-win}
Let $q\geq 3$ be an integer and $\mathbb F$ be a finite field.
Then there is a function $f_{\ref{3-con-win}}(q,{\mathbb F})$
such that a $3$-connected matroid in 
$E(U_{2,q+2},U_{q,q+2},\Lambda_q,
\Delta_q)$
has at most
$f_{\ref{3-con-win}}(q,{\mathbb F})$ inequivalent 
${\mathbb F}$-representations. 
\end{corollary}

Excluding both a free swirl and a free spike is a significant
constraint, but it is not so severe that we lose all interesting
classes. We give one illustration.
Let $p$ be a prime. Recall that $p$ is a {\em Mersenne prime} if
$p=2^m-1$ for some integer $m$. It is a well-known and
widely believed conjecture that the number of Mersenne primes is
infinite. We first note that if $p$ is a Mersenne prime,
then not all free swirls are representable over $GF(p+1)$.
While this is widely known, there does not appear to be a proof
in the literature so we give one here. Readers familiar
with bias matroids of group-labelled graphs, see for example
Zaslavsky \cite{zas}, will find the proof particularly obvious.
The bound that the proof provides is certainly not tight.

\begin{lemma}
\label{exclude-swirls}
Let $\mathbb F$ be a finite field such that $|\mathbb F|-1$ is 
prime. Then there is an integer $f_{\ref{exclude-swirls}}(|\mathbb F|)$,
such that, if $n\geq f_{\ref{exclude-swirls}}(|\mathbb F|)$,
then $\Delta_n$ is not $\mathbb F$-representable.
\end{lemma}

\begin{proof}
Let $M_n$ denote a matroid whose ground set consists of a basis
$B=\{b_1,b_2,\ldots,b_n\}$ together with a set
$\{e_1,f_1,e_2,f_2,\ldots,e_n,f_n\}$ such that,
for $i\in\{1,2,\ldots,n-1\}$, the elements $e_i$ and 
$f_i$ are placed freely on the line spanned by $\{b_i,b_{i+1}\}$
and $e_n,f_n$ are placed freely on the line spanned by
$\{b_n,b_1\}$.

Recall that $\Delta_n\cong M_n\ba B$. We will refer to the members
of $\{\{e_1,f_1\},\{e_2,f_2\},\ldots,\{e_n,f_n\}\}$ as the 
{\em legs} of $M_n$. Note that a representation
of $\Delta_n$ induces a representation of $M_n$ by adding the
the points of intersection of the
legs of $\Delta_n$. 
Therefore $M_n$ is representable over a 
field $\mathbb F$ if and only if $\Delta_n$ is.
Thus we lose no generality in focussing on 
representations of $M_n$.

Let $n=|{\mathbb F}|^2$. We show that
$M_n$ is not $\mathbb F$-representable. Note that
$\Delta_3\cong U_{3,6}$ and that $U_{3,6}$ is not 
$GF(3)$-representable, so the lemma holds in this case. Thus we
may assume that $\mathbb F$ has even order and hence that $-a=a$
for all $a\in{\mathbb F}$.

Assume that $M_n$ is $\mathbb F$-representable. Consider a
standard representation of $M_n$ relative to the basis $B$.
For convenience we identify elements of the legs of $M_n$ with the
column vectors that they label. Up to scaling we may assume,
for $i\in\{1,2,\ldots,n-1\}$, that 
$e_i=(0,\ldots,0,1,1,0,\ldots,0)^T$ and that
$f_i=(0,\ldots,0,1,\alpha_i,0,\ldots,0)^T$, where 
$\alpha_i\neq 0$ and the nonzero entries are in the
$i$th and $(i+1)$th coordinates. We may also assume that
$f_n=(\alpha_n,0,\ldots,0,1)^T$. Note that $e_n$ plays no
role in this argument.

Consider a transversal $T$ of the legs of $M_n$ that contains
$f_n$. It is easily verified that the matrix labelled by
$T$ has determinant $\Pi_{f_i\in T}(\alpha_i)$.
As $n=|{\mathbb F}|^2$, there is an element $\alpha\in{\mathbb F}^*$
such that $\alpha_i=\alpha$ for at least $|{\mathbb F}^*|$
members of $\{1,2,\ldots,n-1\}$. As $M_n$ is simple, $\alpha\neq 1$.
As ${\mathbb F}^*$ has prime order, $\alpha$ generates 
${\mathbb F}^*$. Thus, for some 
$l\in\{1,2,\ldots,|{\mathbb F}|^*\}$, we have 
$\alpha^l=\alpha_n^{-1}$. It now follows that there is a transversal
of the legs that label a square matrix whose determinant is
0, contradicting the assumption that we have a representation of 
$M_n$.
\end{proof}

For a prime power $q$, let ${\mathcal R}(q)$ denote the class of
matroids representable over all fields of size at least $q$.
Thus, ${\mathcal R}(2)$ and ${\mathcal R}(3)$ are the classes of regular
and near-regular matroids respectively.

\begin{theorem}
\label{mersenne1}
Let $\mathbb F$ be a finite field and $q$ be a prime power.
If there is a Mersenne prime greater than or equal to $q-1$,
then there is a function $f_{\ref{mersenne1}}({\mathbb F},q)$
such that a $3$-connected matroid in ${\mathcal R}(q)$
has at most $f_{\ref{mersenne1}}({\mathbb F},q)$ inequivalent
$\mathbb F$-representations.
\end{theorem}

\begin{proof}
Assume that there is a Mersenne prime
greater than $q$. Then, by Lemma~\ref{exclude-swirls},
there is a free swirl that is not in ${\mathcal R}(q)$.
Certainly there is a prime $p$ greater than $q$, so there is a free
spike that is not in ${\mathcal R}(q)$. 
It follows that there is an integer $q'$ such that the 
members of ${\mathcal R}(q)$ are in
$E(U_{2,q'+2},U_{q',q'+2},\Lambda_{q'},\Delta_{q'})$.
The theorem now follows from Corollary~\ref{3-con-win}.
\end{proof}

\begin{corollary}
\label{mersenne2}
There are infinitely many Mersenne primes if and only if,
for each prime power $q$, there is a number $\rho(q)$
such that a $3$-connected member of ${\mathcal R}(q)$
has at most $\rho(q)$ inequivalent $GF(7)$-representations.
\end{corollary}

\begin{proof}
If there are infinitely many Mersenne primes, then the
corollary follows from Theorem~\ref{mersenne1}. Assume that
there are a finite number of Mersenne primes. Let 
$q$ be a prime power larger than the largest Mersenne prime.
Then ${\mathcal R}(q)$ contains all free swirls. But it is shown
in \cite{oxvewh96} that $\Delta_n$ has at least $2^n$ inequivalent
$GF(7)$ representations, so that no bound can be places on the
number of inequivalent $GF(7)$-representations of members
of ${\mathcal R}(q)$.
\end{proof}

\section{Certifying Non-representability}

As noted in the introduction, one of the main applications of the
material in this paper is in \cite{gewh1} where the following
theorem is proved.

\begin{theorem}[\cite{gewh1}~Theorem 1.1]
\label{final-task}
For any prime $p$, proving that an $n$-element matroid is not
representable over $GF(p)$ requires at most $O(n^2)$ rank evaluations.
\end{theorem}

The purpose of this section is to adumbrate a straightforward 
consequence of earlier material from this paper that is 
optimised for the application in \cite{gewh1}.
We consider a slight augmentation of the class of
$k$-coherent matroids. A matroid $M$ is {\em near} $k$-coherent
\index{near $k$-coherent matroid}
if it is connected and either $\si(M)$ or $\co(M)$
is $k$-coherent. 

\begin{corollary}
\label{k-coherent}
Let $p$ be a prime number and $k\geq 5$ be an integer. Then the
following hold.
\begin{itemize}
\item[(i)] If $M$ is a nonempty near $k$-coherent matroid,
then there is an element $e\in E(M)$ such that either 
$M\ba e$ or $M/e$ is near $k$-coherent.
\item[(ii)] Let $M$ be a near $k$-coherent and $e\in E(M)$.
Then the following hold.
\begin{itemize}
\item[(a)] If $M/e$ is not near $k$-coherent, then there is a 
$4$-separation $(A,B)$ in $M$ such that 
$e\in\cl_M(A-\{e\})\cap\cl_M(B-\{e\})$.
\item[(b)] If $M\ba e$ is not near $k$-coherent, then there
is a $4$-separation $(A,B)$ in $M$ such that
$e\in\cl^*_M(A-\{e\})\cap \cl^*_M(B-\{e\})$.
\end{itemize}
\item[(iii)] There is an integer $\mu_p$ such that a 
near $k$-coherent matroid has at most $\mu_p$ inequivalent
representations over $GF(p)$.
\end{itemize}
\end{corollary}

\begin{proof}
Consider (i). If $M$ is not $k$-coherent, then we may
either delete an element in a non-trivial parallel class
or contract an element in a non-trivial series class to 
preserve near $k$-coherence. If $M$ is $k$-coherent 
then it follows that from Corollary~\ref{unexpose2}
that $M$ has an element $e$ such that either $M\ba e$
or $M/e$ is $k$-coherent unless $M$ is a wheel or a whirl.
Assume that  $M$ is a wheel or a whirl. Let $e$ be a rim
element. Observe that $M/e$ is near $k$-coherent.

Consider condition (ii). 
Let $M$ be a near $k$-coherent matroid. Assume that $e\in E(M)$
and that $M/e$ is not near $k$-coherent. Then $e$ is not in a 
non-trivial series class. Assume that $e$ is in a non-trivial
parallel class. Then $e$ is in the guts of a 2-separation in 
$M$ and (ii) holds.
If $M$ is $3$-connected, then it follows from Lemma~\ref{x-blocks}
that $e$ is in the guts of either a 3-separation or a 4-separation
in $M$. The upgrade to the case when $M$ has non-trivial series
classes or parallel classes which do not contain $e$
is elementary and is omitted. The case when $M\ba e$ is not near
$k$-coherent follows by duality.

Consider (iii). Note that we can either delete
fixed elements or contract cofixed elements from a 
near $k$-coherent matroid to obtain a $k$-coherent matroid.
Part (iii) follows from this observation,
Lemma~\ref{extend-fixed} and Theorem~\ref{biggy}.
\end{proof}

\printindex


\begin{thebibliography}{77}

\bibitem{ai-ox} J. Aikin and J. Oxley,
The structure of crossing separations in matroids
{\em Advances in Applied Mathematics}, 
{\bf 41} (2008), 10--26.

\bibitem{be06} B. D. Beavers, {\em Circuits and structure in
matroids and graphs}, PhD Thesis, Louisiana State University, (2006). 

\bibitem{bi1} R. E. Bixby, A simple theorem on $3$-connectivity,
{\em Linear Algebra Appl.} {\bf 45} (1982), 123--126.


\bibitem{doov1} G. Ding, B. Oporowski, J. Oxley, and
D. Vertigan, Unavoidable minors of large $3$-connected 
binary matroids,
{\em J. Combin.\ Theory Ser.\ B} {\bf 66} (1996), 334--360.

\bibitem{doov} G. Ding, B. Oporowski, J. Oxley, and
D. Vertigan, Unavoidable minors of large $3$-connected matroids,
{\em J. Combin.\ Theory Ser.\ B} {\bf 71} (1997), 244--293.


\bibitem{duke} R. Duke, Freedom in matroids, 
{\em Ars Combin. } {\bf 26B} (1988), 191Ð-216.

\bibitem{ggk} J. Geelen, B. Gerards, and A. Kapoor,
The excluded minors for $GF(4)$-representable matroids,
{\em J. Combin.\ Theory Ser.\ B}, {\bf 79} (2000), 247--299.

\bibitem{gegewh06} J. Geelen, B. Gerards, and G. Whittle,
Matroid $T$-connectivity, 
{\em SIAM J. Discrete Math.}, {\bf 20} (2006), 588Ð-596.

\bibitem{gegewh09} J. Geelen, B. Gerards, and G. Whittle,
Tangles, tree-decompositions and grids in matroids 
{\em J. Combin.\ Theory Ser.\ B} {\bf 99} (2009), 657Ð-667.

\bibitem{gegewh10} J. Geelen, B. Gerards, and G. Whittle,
On inequivalent representations of matroids over non-prime fields,
{\em J. Combin.\ Theory Ser.\ B} {\bf 100} (2010), 740Ð-743.


\bibitem{gehlwh} J. Geelen, P. Hlin\v en\'y, and G. Whittle,
Bridging Separations in Matroids, 
{\em SIAM J. Discrete Math.} {\bf 18} (2004), 638Ð-646.

\bibitem{totally-free} J. Geelen, J. G. Oxley, D. Vertigan,
and G. Whittle, Totally free expansions of matroids,
{\em J. Combin.\ Theory Ser.\ B}, 
{\bf 84} (2002), 130Ð-179. 

\bibitem{short-gf5} J. Geelen, J. G. Oxley, D. Vertigan,
and G. Whittle, A short proof of non-$GF(5)$-representability
of matroids, {\em J. Combin.\ Theory Ser.\ B}, 
{\bf 91} (2004), 105Ð-121.

\bibitem{gewh1} J. Geelen and G. Whittle,
Certifying non-representability of matroids over prime fields,
submitted.

\bibitem{ka88}  
J. Kahn, On the uniqueness of matroid representations over 
$GF(4)$, {\em Bull.\ London Math.\ Soc.} {\bf 20} (1988), 5--10.

\bibitem{ku93}
J. P. S. Kung, Extremal matroid theory, in 
{\em Graph structure theory}, N. Robertson and
P. D. Seymour, eds., Amer.\ Math.\ Soc.,
Providence, Rhode Island, 1993, 21--62.

\bibitem{oxlem01} M. Lemos, and J. Oxley,
A sharp bound on the size of a connected matroid, 
{\em Trans.\ Amer.\ Math.\ Soc.} {\bf 353} (2001), 4039--4056.

\bibitem{leox03} M. Lemos, and J. Oxley,
On the minor-minimal 3-connected matroids 
having a fixed minor, {\em Europ.\ J. Combin.},
{\bf 24} (2003), 1097--1123.

\bibitem{ox92} J.  Oxley,  {\em Matroid Theory,} 
Oxford University Press, New York, 1992.

\bibitem{flower} J.  Oxley, C.  Semple, and G. Whittle,  
The structure of the 3-separations of 3-connected matroids, 
{\em J. Combin.\ Theory Ser.\ B}, {\bf 92} (2004), 257--293.

\bibitem{flower2}  J.  Oxley, C.  Semple, and G. Whittle,  
The structure of the 3-separations of 
3-connected matroids II, {\em European Journal of 
Combinatorics}, {\bf 28} (2007), 1239--1261.

\bibitem{wild} J.  Oxley, C.  Semple and G. Whittle,
Wild triangles in 3-connected matroids,
{\em J. Comb.\ Theory Ser.\ B}, {\bf 98} (2008), 291Ð-323.


\bibitem{expose} J.  Oxley, C.  Semple and G. Whittle,
Exposing 3-separations in 3-connected matroids. To appear in
{\em Advances in Applied Math.} 

\bibitem{upgrade} J.  Oxley, C.  Semple and G. Whittle,
An upgraded wheels and whirls theorem for 
3-connected matroids, submitted.

\bibitem{oxvewh96} J.  Oxley, D. Vertigan, and 
G. Whittle, On inequivalent representations of matroids
over finite fields, {\em J. Combin.\ Theory Ser.\ B},
{\bf 67} (1996), 325--343.

\bibitem{se1} P. D. Seymour, Decomposition of regular matroids,
{\em J. Combin.\ Theory Ser.\ B}, {\bf 28} (1980), 305--359. 

\bibitem{seymour1} P. D. Seymour, Recognising graphic matroids,
{\em Combinatorica} {\bf 1} (1981), 75--78.

\bibitem{tu2} W. T. Tutte, MengerÕs theorem for matroids, 
{\em J. Res.\ Nat.\ Bur.\ Standards Section B}, 
{\bf 69} (1965), 49Ð-53.

\bibitem{tu1}
W. T. Tutte, 
Connectivity in matroids,
{\em Canad.\ J. Math.} {\bf 18} (1966), 1301--1324.

\bibitem{zas} T. Zaslavsky, Biased graphs. II. The three matroids,
{\em J. Combin.\ Theory Ser.\ B}, {\bf 51} (1991), 46--72. 

\end{thebibliography}
\end{document}